\documentclass[a4paper,11pt,openany]{book}
\usepackage[margincaption,outercaption,ragged,wide]{sidecap} 

\usepackage[utf8]{inputenc} 

\usepackage{appendix}
\usepackage{fancyhdr} 

\usepackage[a4paper, left=4cm, right=3.5cm, top=3.5cm, bottom=3.5cm]{geometry}

\usepackage[T1]{fontenc} 
\usepackage{anyfontsize}
\usepackage{mathtools}
\usepackage{amsmath}
\usepackage{amssymb}
\usepackage{amsthm}
\usepackage{amsfonts} 
\usepackage{stmaryrd}
\usepackage{upgreek}
\usepackage{enumitem}
\usepackage[dvipsnames]{xcolor}
\usepackage{tikzit}
\usetikzlibrary{decorations.markings}
\usetikzlibrary{arrows.meta, arrows}
\usetikzlibrary{decorations.pathreplacing,angles,quotes}
\usetikzlibrary{patterns}
\usepackage{quiver}
\usepackage{adjustbox}
\usepackage{graphicx}
\usepackage{textcomp}
\usepackage[normalem]{ulem}
\usepackage[all]{xy}
\usepackage{cancel}
\usepackage{scalerel}
\usepackage[backend=biber, 
style=alphabetic,
minalphanames=4,
maxalphanames=5, 
maxbibnames=5]{biblatex}

\usepackage{mathpazo}
\usepackage{tgpagella}

\usepackage{hyperref}
\hypersetup{
   colorlinks,
   linkcolor={black},
   citecolor={black},
   urlcolor={blue}
}

\newlength{\pageoffset}
\newcommand{\clearemptydoublepage}{\clearpage{\pagestyle{empty}\cleardoublepage}}

\fancypagestyle{plain}{
    \fancyhf{}
    \fancyhead[LE]{\makebox[\pageoffset][l]{\thepage}}
    \fancyhead[RO]{\makebox[\pageoffset][r]{\thepage}}
}
\fancypagestyle{pretext}{
    \fancyhf{}
    \fancyhead[LE]{\makebox[\pageoffset][l]{\thepage}\leftmark}
    \fancyhead[RO]{\leftmark\makebox[\pageoffset][r]{\thepage}}
}
\fancypagestyle{maintext}{
    \fancyhf{}
    \fancyhead[LE]{\makebox[\pageoffset][l]{\thepage}\leftmark}
    \fancyhead[RO]{\rightmark\makebox[\pageoffset][r]{\thepage}}
}
\appto\frontmatter{\pagestyle{pretext}}
\appto\mainmatter{\pagestyle{maintext}}

\theoremstyle{plain}
\newtheorem{theorem}{Theorem}[section]
\newtheorem{lemma}[theorem]{Lemma}
\newtheorem{proposition}[theorem]{Proposition}
\newtheorem{corollary}[theorem]{Corollary}
\newtheorem{introthm}{Theorem}

\newtheorem{introcor}[introthm]{Corollary}

\theoremstyle{definition}

\newtheorem{prop-def}[theorem]{Proposition-Definition}
\newtheorem{definition}[theorem]{Definition} 
\newtheorem{example}[theorem]{Example}

\newtheorem{remark}[theorem]{Remark}

\usepackage{macros}

\makeatletter
\newcommand*\ifStrInTF[2]{%
  \edef\tikz@temp{{#1}{#2}}%
  \expandafter\pgfutil@in@\tikz@temp
  \ifpgfutil@in@\expandafter\pgfutil@firstoftwo\else\expandafter\pgfutil@secondoftwo\fi}
\newcommand*\ifStrEmptyTF[1]{%
  \def\tikz@temp{#1}\ifx\tikz@temp\pgfutil@empty
    \expandafter\pgfutil@firstoftwo\else\expandafter\pgfutil@secondoftwo\fi}

\def\tikz@swapanchor#1.#2\tikz@stop#3#4{#1#4#3}
\newcommand*\tikzAddAnchor[2]{%
  \ifStrInTF{.}{#1}{%
    \ifStrEmptyTF{#2}
      {\edef#1{\expandafter\tikz@swapanchor#1\tikz@stop{}{}}}
      {\edef#1{\expandafter\tikz@swapanchor#1\tikz@stop{#2}{.}}}%
  }{%
    \ifStrEmptyTF{#2}{}
      {\edef#1{#1.#2}}%
  }}
\tikzset{
  pin anchor/.style={tikz@pin@post/.append style={anchor=#1}},
  label anchor/.style={tikz@label@post/.append style={anchor=#1}},
  pin edge pin anchor/.style={
    append after command={\pgfextra\tikzAddAnchor{\tikzlastnode}{#1}\endpgfextra}%
  },
  pin edge parent anchor/.style={
    append after command={\pgfextra\tikzAddAnchor{\tikz@save@last@node}{#1}\endpgfextra}%
  }
}
\makeatother

\tikzstyle{dot}=[fill=black, draw=black, shape=circle, scale=0.5]
\tikzstyle{white-dot}=[fill=white, draw=black, shape=circle, scale=0.5]
\tikzstyle{blue-dot}=[fill=blue, draw=blue, shape=circle, scale=0.5]
\tikzstyle{red-dot}=[fill=red, draw=red, shape=circle, scale=0.5]
\tikzstyle{light-green-rounded-rectangle}=[fill={rgb,255: red,189; green,255; blue,197}, draw=black, shape=rectangle, rounded corners=3pt]
\tikzstyle{light-pink-rounded-rectangle}=[fill={rgb,255: red,255; green,225; blue,255}, draw=black, shape=rectangle, rounded corners=3pt]
\tikzstyle{light-yellow-rounded-rectangle}=[fill={rgb,255: red,255; green,254; blue,215}, draw=black, shape=rectangle, rounded corners=3pt]
\tikzstyle{thick-blue-implies}=[{=>}, draw=blue, thick]

\tikzstyle{dashed-black}=[-, dashed]
\tikzstyle{dotted-gray}=[-, dotted, draw={rgb,255: red,191; green,191; blue,191}]
\tikzstyle{red}=[-, draw=red]
\tikzstyle{thick}=[-, line width=1.5pt]
\tikzstyle{to}=[->]
\tikzstyle{red-to}=[->, draw=red]
\tikzstyle{blue-to}=[draw=blue, ->]
\tikzstyle{dash-to}=[->, dashed]
\tikzstyle{squig-to}=[->, decorate, decoration={{zigzag, segment length=4, amplitude=.9, post=lineto,postlength=2pt}}, line join=round]
\tikzstyle{red-dash-to}=[->, draw=red, dashed]
\tikzstyle{thick-blue-to}=[->, double distance=1pt, line width=1pt, draw=blue]
\tikzstyle{mapsto}=[{|->}]
\tikzstyle{directed}=[-, decoration={markings,mark=at position 0.7 with {\arrow{To[scale=1.25]}}}, postaction=decorate]
\tikzstyle{directed-dash}=[-, dashed, decoration={markings,mark=at position 0.7 with {\arrow{To[scale=1.25]}}}, postaction=decorate]
\tikzstyle{directed-dotted}=[-, dotted, decoration={markings,mark=at position 0.7 with {\arrow{To[scale=1.25]}}}, postaction=decorate]
\tikzstyle{directed-red}=[-, draw=red, decoration={markings,mark=at position 0.7 with {\arrow{To[scale=1.25]}}}, postaction=decorate]
\tikzstyle{directed-red-dash}=[-, draw=red, dashed, decoration={markings,mark=at position 0.7 with {\arrow{To[scale=1.25]}}}, postaction=decorate]
\tikzstyle{directed-blue}=[-, draw=blue, decoration={markings,mark=at position 0.7 with {\arrow{To[scale=1.25]}}}, postaction=decorate]
\tikzstyle{directed-blue-dash}=[-, draw=blue, dashed, decoration={markings,mark=at position 0.7 with {\arrow{To[scale=1.25]}}}, postaction=decorate]
\tikzstyle{light-gray-05-fill}=[-, fill={rgb,255: red,220; green,220; blue,220}, draw=none, fill opacity=0.5]
\tikzstyle{dark-gray-05-fill}=[-, fill={rgb,255: red,128; green,128; blue,128}, fill opacity=0.50, draw=none]
\tikzstyle{solid-green-fill}=[-, fill={rgb,255: red,66; green,192; blue,139}, draw=none]
\tikzstyle{light-yellow-fill-black-border}=[-, draw=black, fill={rgb,255: red,255; green,254; blue,215}]
\tikzstyle{light-yellow-05-fill}=[-, fill={rgb,255: red,255; green,254; blue,215}, draw=none, fill opacity=0.5]
\tikzstyle{light-magenta-05-fill}=[-, fill={rgb,255: red,255; green,225; blue,255}, draw=none, fill opacity=0.5]
\tikzstyle{light-cyan-05-fill}=[-, fill={rgb,255: red,196; green,255; blue,249}, draw=none, fill opacity=0.5]
\tikzstyle{dotted-fill}=[-, pattern=dots, pattern color=gray, draw=none]

\newcommand{\titlename}{Duals of higher vector bundles and cotangents of Lie $2$-groupoids}
\newcommand{\authorname}{Stefano Ronchi}
\newcommand{\authorfullname}{Stefano Ronchi}
\newcommand{\authorcity}{Vimercate, Italy}
\newcommand{\doctoralprogram}{Mathematical Sciences}
\newcommand{\dateOralExamination}{May 19th, 2025}

\newcommand{\submissionYear}{2025}

\newcommand{\keywordsOfThesis}{Duality, higher vector bundles, Lie $n$-groupoids, shifted symplectic structures, higher vector spaces, Eilenberg-Zilber theorem}
\newcommand{\abstractOfThesis}{In this thesis we define $n$-duals of $\VB$ $n$-groupoids over Lie $n$-groupoids and study their properties.  
For $n=0$ this returns the dual vector bundle construction, while for $n=1$ this returns Pradines's construction of the dual of a $\VB$ groupoid over a Lie groupoid, which includes the cotangent symplectic groupoid of Coste, Dazord and Weinstein. 
For $n=2$, we propose a new construction that shows that $\VB$ 2-duals exist for $\VB$ 2-groupoids and they are $\VB$ 2-groupoids themselves. Their canonical dual pairings are nondegenerate up to homotopy in the same sense as shifted symplectic structures. 

In particular, we can apply this construction to the tangent of a Lie 2-groupoid and obtain a cotangent $\VB$ 2-groupoid (the 2-cotangent) which is canonically 2-shifted symplectic. 
We apply this in two ways: First, to characterize 2-shifted symplectic structures on a Lie 2-groupoid as Morita equivalences between its tangent and 2-cotangent groupoid. 
Second, to compute the 2-cotangent of a Lie 1-groupoid and show it is symplectic Morita equivalent to the bar construction of the 1-cotangent.

Along the way, we develop the theory of $n$-duals for simplicial vector spaces, which covers the case where the base is a point. 
In this case, $n$-duals always exist, as they are defined by a mapping space construction. By a reformulation of the Eilenberg-Zilber theorem in terms of mapping spaces, we obtain that the canonical $n$-dual pairing is nondegenerate up to homotopy for all $n$-types.}
\author{\authorname}

\begin{document}
\frontmatter

\thispagestyle{empty}

\begin{center}
\begin{minipage}{0.8\linewidth}
    \centering
    
    \vspace{7mm}
    \includegraphics[width=\linewidth]{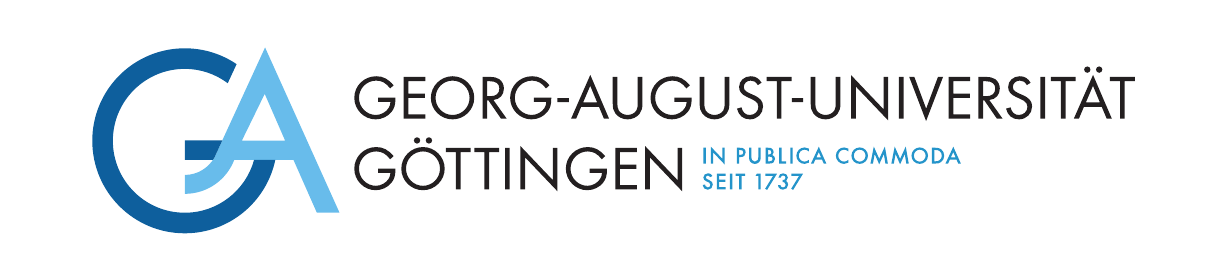}
    \par
    \vspace{21mm}
    
    \hrule
    \vspace{7mm}
    {{\huge \titlename\par}}
    \vspace{7mm}
    \hrule

    \vspace{14mm}
    {\huge Dissertation \par}
    \vspace{14mm}
    {\small for the award of the degree\par}
    {\Large ``Doctor rerum naturalium'' (Dr.\ rer.\ nat.)\par}
    \vspace{2mm}
    {\small of the Georg-August-Universität Göttingen\par}
    \vspace{14mm}
    {\small within the doctoral program\par}
    {\Large \doctoralprogram\par}
    \vspace{2mm}
    {\small of the Georg-August University School of Science (GAUSS)\par}

    \vspace{14mm}
    {\small submitted by\par}
    {\Large \authorfullname\par}
    \vspace{2mm}
    {\small from \authorcity\par}
    
    \vspace{14mm}
    {\normalsize Göttingen, \submissionYear}
    
\end{minipage}
\end{center}

\newpage


{\setlength\parindent{0pt}
\thispagestyle{empty}

\large\textbf{Thesis Committee}\\[-2mm]

\normalsize Prof.\ Dr.\ Chenchang \textsc{Zhu}\\
Mathematisches Institut, Georg-August-Universität Göttingen\\[-2mm]

Prof.\ Dr.\ Madeleine \textsc{Jotz}\\
Institut für Mathematik, Julius-Maximilians-Universität Würzburg\\[-2mm]

Prof.\ Dr.\ Ralf \textsc{Meyer}\\
Mathematisches Institut, Georg-August-Universität Göttingen\\[8mm]

\large\textbf{Examination Board}\\[-2mm]

\normalsize \textbf{Reviewer}\\
Prof.\ Dr.\ Chenchang \textsc{Zhu}\\
Mathematisches Institut, Georg-August-Universität Göttingen\\[-2mm]

\normalsize \textbf{Second Reviewer}\\
Prof.\ Dr.\ Madeleine \textsc{Jotz}\\
Institut für Mathematik, Julius-Maximilians-Universität Würzburg\\[-2mm]

\normalsize\textbf{Further members of the Examination Board}\\[-2mm]

Prof.\ Dr.\ Ralf \textsc{Meyer}\\
Mathematisches Institut, Georg-August-Universität Göttingen\\[-2mm]

Dr.\ Leonid \textsc{Ryvkin}\\
Institut Camille Jordan, Université Claude Bernard Lyon\\[-2mm]

Prof.\ Dr.\ Damien \textsc{Calaque}\\
Institut Montpelliérain Alexander Grothendieck, Université de Montpellier\\[-2mm]

Prof.\ Dr.\ Christoph \textsc{Lehrenfeld}\\
Institut für Numerische und Angewandte Mathematik, Georg-August-Universität Göttingen\\[-2mm]

\vfill
\normalsize\textbf{Date of the oral examination}\\
\dateOralExamination


\clearpage
\vspace*{3.5cm}
\thispagestyle{empty}
\begin{center}
    \bfseries\Large Abstract
\end{center}

\abstractOfThesis\\

\textbf{Keywords}:  \keywordsOfThesis\\

\textbf{DOI}: doi:10.53846/goediss-11918

\clearemptydoublepage}

\chapter*{Preface}
\markboth{Preface}{}
\addcontentsline{toc}{chapter}{Preface}
\setcounter{page}{1}

Let $E$ and $F$ be vector bundles over the same base manifold $M$. A \textit{pairing} of vector bundles between $E$ and $F$ is a bundle map
\begin{equation*}
    E \otimes F \longrightarrow \R \times M,
\end{equation*}
over the identity on $M$.\footnote{We write this as a bundle map with target $\R \times M$ to emphasize the fact that it is defined fiberwise, even though, after evaluating, the basepoint is usually omitted.}
Every vector bundle $E$ over a smooth manifold $M$ has a dual vector bundle $E^*$, which can be defined by the following universal property:
\begin{equation}\label{eq:vb-adjunction-intro}
    \VB_M(E \otimes F, \R \times M)\cong \VB_M(F, E^*), \quad \forall F\in \VB_M,
\end{equation}
where $\VB_M$ is the category of vector bundles over $M$ with bundle maps as morphisms. This property can be stated as the fact that $E^*$ classifies the space of vector bundle pairings with $E$. The dual itself has a canonical evaluation pairing 
\begin{equation*}
    \langle \cdot, \cdot \rangle: E^* \otimes E \longrightarrow \R \times M,
\end{equation*}
which corresponds to the identity under the isomorphism of hom spaces above.

In the specific example of the tangent bundle of a manifold $TM\to M$, this construction returns the cotangent bundle $T^*M$. The canonical pairing induces the \textit{tautological 1-form} $\theta \in \Omega^1(T^*M)$, which is defined by
\begin{equation}\label{eq:intro-taut-form}
    \theta_\xi (V) := \langle \xi, TqV \rangle, \quad \forall \xi \in T^*M,\; V \in T_{\xi}T^*M,
\end{equation}
where $q: T^*M \to M$ is the bundle projection. The de Rham differential of $\theta$ defines the canonical \textit{symplectic} form $\omega_{can} := -d\theta$ on $T^*M$. This means it is a closed 2-form and it defines a nondegenerate antisymmetric pairing 
\begin{equation*}
    \omega_{can}: TT^*M \wedge TT^*M \longrightarrow \R \times T^*M. 
\end{equation*}
The nondegeneracy condition is equivalent to the fact that the associated map $\omega_{can}^\flat: TT^*M \to T^*T^*M$ induced by applying the adjunction \eqref{eq:vb-adjunction-intro} is an isomorphism.
This structure is of fundamental importance in the theory of hamiltonian mechanics and is often considered to be the most fundamental example of symplectic geometry: by the Darboux theorem, any $k$-dimensional symplectic manifold $(M, \omega)$ is locally isomorphic to $(T^*\R^k, \omega_{can})$.

A Lie groupoid is a generalization of a smooth manifold (see e.g. \cite{CannasdaSilvaWeinstein2000, Mackenzie2005}). This is a (small) category where the set of object and the set of arrows are smooth manifolds, all arrows are invertible, and where the source and target maps are surjective submersions. This is commonly pictured as the diagram 
\begin{equation*}
    \huaG_1 \rightrightarrows \huaG_0,
\end{equation*}
where $\huaG_1$ is the manifold of arrows, $\huaG_0$ is the manifold of points, and the horizontal arrows are the source and target maps, which we denote by $d_0$ and $d_1$, respectively. Some examples of Lie groupoids are Lie groups, which have $\huaG_0 = pt$, manifolds, which have $\huaG_1=\huaG_0$, action groupoids $G\times M \rightrightarrows M$, where $M$ is a manifold acted on by a Lie group $G$, principal bundles, and vector bundles, among others.

A vector bundle over a Lie groupoid is known as a $\VB$-groupoid (e.g. \cite{Mackenzie2005,GraciaSazMehta2017}). This is a Lie groupoid in the category of vector bundles. It can be pictured as a diagram 
\begin{equation*}
\begin{array}{ccc}
    \huaV_1 &\rightrightarrows &\huaV_0\\
    \downarrow & &\downarrow\\
    \huaG_1 &\rightrightarrows &\huaG_0,
\end{array}
\end{equation*}
where the vertical maps are vector bundle projections and Lie groupoid morphisms, and the horizontal maps are vector bundle maps. The multiplication of the total space $\huaV_1$ is also a vector bundle map over the multiplication in $\huaG_1$, leading to the so-called \textit{interchange law},
\begin{equation*}
    (w_h \cdot v_g) + (w'_h \cdot v'_g) = (w_h + w'_h) \cdot (v_g + v'_g),
\end{equation*}
for any pairs of composable arrows $(v_g, w_h)$ and $(v'_g, w'_h)$ in $\huaV_1$ over the same pair of composable arrows $(g,h)$ in the base $\huaG_1$.

A duality for $\VB$-groupoids was identified by Pradines in \cite{Pradines1988} (see also \cite{Mackenzie2005}). 
To rephrase it in similar terms as we did for vector bundles we immediately run into the problem that the tensor product of $\VB$-groupoids is not generally a $\VB$-groupoid itself.\footnote{Pradines uses the levelwise Whitney sum of $\huaV$ and its dual as the domain of the canonical pairing and requires that it is a bilinear $\VB$-groupoid morphism, to avoid precisely this problem. With this formulation the universal property interpretation in terms of a tensor-hom adjunction is lost.} 
This is a natural point to introduce higher Lie groupoids and higher vector bundles over them, in order to consider tensor products of $\VB$-groupoids. 

We follow the simplicial point of view that started with Duskin's \textit{hypergroupoids} (see \cite{Duskin1979,Glenn1982,Duskin2001/02}) and was extended to the smooth setting in \cite{Getzler2009, Henriques2008, Zhu2009} and studied in the later works \cite{Pridham2013, Li2014, Wolfson2016,BehrendGetzler2017, RogersZhu2020}, among others. 

By applying the \textit{nerve} construction to a Lie groupoid, (see example \ref{ex:nerve-cat} for details) this can be seen as a simplicial manifold satisfying certain Kan conditions:
\begin{itemize}
    \item $\Kan(1,j)$, for $j=0,1$: the source and target map are surjective submersions,
    \item $\Kan!(m,j)$, for $j=0,\dots, m$ and any $m>1$: any sequence of $m$ arrows forming a zig-zag identifies a unique $m$-simplex. 
\end{itemize}
By abstracting these conditions one observes that the difference between them is that certain projections are surjective submersions in the case of the Kan condition $\Kan(m,j)$, and diffeomorphisms in the case of the strict Kan condition $\Kan!(m,j)$. A Lie $n$-groupoid can be defined as a simplicial manifold with all Kan conditions for $m\ge 1$, $0 \le j \le m$ and strict Kan conditions for $m > n$, $0 \le j \le m$, as we write in Definition \ref{def:Lie-n-gpd}. If only the non-strict Kan conditions are required to hold, this defines a Lie $\infty$-groupoid. When forgetting the smooth structure, these are Kan complexes (see e.g. \cite{GoerssJardine2009}), though requiring a simplicial manifold to be a set-theoretical Kan complex is not enough to obtain a Lie $\infty$-groupoid, as discussed in detail in \cite{RogersZhu2020}. 

In this picture, smooth manifolds are Lie 0-groupoids, and Lie groupoids are Lie 1-groupoids by identifying them with their nerves. 
Lie groupoid morphisms are then simplicial maps between the nerves. It is worth noting that any simplicial map with target a Lie $n$-groupoid is determined by its first $n+1$ levels, i.e. its action on the 0 to $(n+1)$-simplices of the source, as we discuss in Section \ref{sec:finite-data}. 

Following this line of thought, a $\VB$ $n$-groupoid is a simplicial vector bundle over a Lie $n$-groupoid which is itself a Lie $n$-groupoid, as we write in Definition \ref{def:vb-n-gpd}. These were defined in \cite{HoyoTrentinaglia2024}. More concretely, each level is a vector bundle, the face and degeneracy maps are bundle maps over the corresponding maps of the base, and the bundle projections commute with the face and degeneracy maps. A \textit{higher vector bundle} over a Lie $n$-groupoid $\huaG$ is a $\VB$ $\infty$-groupoid over $\huaG$. In this category, the tensor product of two higher vector bundles is simply the levelwise tensor product of vector bundles applied as a functor to all the face and degeneracy maps, as well. Then, even if the tensor product of two $\VB$-groupoids is not a $\VB$-groupoid, it is still a higher vector bundle of generally higher order. We discuss this in more detail in Remark \ref{rem:order-of-tensor-vb}. 

Going back to $\VB$-groupoids, we can now phrase their duality in terms of \textit{1-shifted simplicial pairings}, as defined in Definition \ref{def:VB-n-shifted-pairing}. Consider the additive group $B^1\R = (\R \rightrightarrows 0)$ which is topologically the Eilenberg-Maclane space $K(\R,1)$. Given a $\VB$-groupoid $\huaV \to \huaG$, Pradines's dual $\huaV^{1*} \to \huaG$ is uniquely determined\footnote{We discuss this in more detail in Remark \ref{rem:VB-1dual-def-by-ev-extending-to-simp-pairing}.} by the property that the canonical evaluation pairing 
\begin{equation*}
    ev: \huaV_1^* \otimes \huaV_1 \longrightarrow \R \times \huaG_1
\end{equation*}
is normalized and multiplicative, i.e. it extends to a 1-shifted simplicial $\VB$-pairing 
\begin{equation*}
    \langle\cdot, \cdot \rangle: \huaV^{1*} \otimes \huaV \longrightarrow B^1\R \times \huaG.
\end{equation*}
It can then be shown that $\huaV^{1*}$ satisfies the universal property
\begin{equation}\label{eq:intro-tensor-1-dual-adjunction}
    \VB^{\infty}_\huaG(\huaV \otimes \huaW, B^1\R \times \huaG) \cong \VB^{\infty}_\huaG(\huaW, \huaV^{1*}), \quad \forall \huaW\in \VB^{\infty}_\huaG,
\end{equation}
where $\VB^{\infty}_\huaG$ is the category of higher vector bundles over $\huaG$ with simplicial bundle maps over the identity of $\huaG$. We do this in Lemma \ref{lem:VB1Dual-computation}. This proves that the \textit{1-dual} $\VB$-groupoid is the object classifying 1-shifted pairings of $\VB$-groupoids. With this definition, we recover the canonical pairing by considering the pairing of $\huaV^{1*}$ and $\huaV$ corresponding to the identity of $\huaV^{1*}$ under the adjunction \eqref{eq:intro-tensor-1-dual-adjunction}.

The first immediate example of a $\VB$-groupoid is the tangent groupoid $T\huaG$ of a Lie groupoid $\huaG$. 
\begin{equation*}
\begin{array}{ccc}
    T\huaG_1 &\rightrightarrows &T\huaG_0\\
    \downarrow & &\downarrow\\
    \huaG_1 &\rightrightarrows &\huaG_0.
\end{array}
\end{equation*}
This is obtained by simply applying the tangent functor to the spaces and maps defining $\huaG$.
By dualizing the tangent groupoid $T\huaG$ of a Lie groupoid $\huaG$ with the aforementioned procedure, one obtains the cotangent groupoid $T^{1*}\huaG$. The canonical 1-shifted simplicial pairing $\langle\cdot, \cdot \rangle$ defines again a 1-shifted tautological form $\theta$ as in \eqref{eq:intro-taut-form}, which, when differentiated, gives the canonical symplectic form $\omega_{can}$. This form is closed, because it is exact, but it is also multiplicative, i.e. it induces a 1-shifted simplicial (antisymmetric) pairing
\begin{equation*}
    \omega_{can}: TT^{1*}\huaG \wedge TT^{1*}\huaG \longrightarrow B^1\R \times T^{1*}\huaG,
\end{equation*}
which is equivalent to the fact that $\omega_{can}\in \Omega^2(T^{1*}\huaG_1)$ is a closed form in the total complex of the Bott-Shulman-Stasheff double complex from \cite{BottShulmanStasheff1976}, see also \cite{AriasAbadCrainic2013,Lesdiablerets,CuecaZhu2023}. Finally $\omega_{can}$ is nondegenerate in the usual sense, so $(T^{1*}\huaG, \omega_{can})$ is a symplectic groupoid, as shown in \cite{CosteDazordWeinstein1987}. 

So far, we only re-elaborated pre-existing results. In this thesis, we unify and generalize the duality theory of vector bundles and $\VB$-groupoids to one where the canonical pairing is $n$-shifted. We define the $n$-dual $\huaV^{n*}$ of a higher vector bundle $\huaV$ over a Lie $n$-groupoid $\huaG$ as the universal object classifying $n$-shifted simplicial pairings with $\huaV$. That is, $\huaV^{n*}$ represents the functor
\begin{equation*}
\begin{split}
    \VB^{\infty}_\huaG &\longrightarrow \Vect\\
    \huaW &\longmapsto \VB_\huaG^\infty(\huaV \otimes \huaW, B^n\R\times \huaG),
\end{split}
\end{equation*}
where $B^n\R$ is now the Eilenberg-MacLane space $K(\R,n)$, or equivalently, the abelian group $\R$ seen as an $n$-group by delooping $n$ times.
If this object exists, then there is an isomorphism
\begin{equation*}
    \VB_\huaG^\infty(\huaV \otimes \huaW, B^n\R\times \huaG) \newrightleftarrows{\rho}{\tau} \VB_\huaG^\infty(\huaW, \huaV^{n*}),
\end{equation*}
which is natural in $\huaW$. With this, we can define the canonical evaluation pairing as for 1-duals:
\begin{equation*}
    \langle \cdot, \cdot \rangle = \tau(id_{\huaV^{n*}}): \huaV^{n*}\otimes \huaV \longrightarrow B^n\R \times \huaG.
\end{equation*}

When applying this duality to the tangent groupoid of a Lie $n$-groupoid one obtains the $n$-cotangent. 
Under the correspondence between Morita classes of $n$-groupoids and $n$-stacks, which was studied in different settings in \cite{Pridham2013,BehrendGetzler2017, RogersZhu2020}, the $n$-cotangent of a Lie $n$-groupoid should correspond to the cotangent stack studied in \cite{Calaque2019}. This correspondence is not clear to us at the moment, but we plan to study it in future work.

About the problem of existence of the $n$-dual, we find that the existence of the $n$-dual for a given higher vector bundle is subject to existence and regularity of the solutions of certain equations. I.e. the $n$-dual can be expressed levelwise as a limit in $\VB$, and the structure of the equations shows that it is entirely determined by solving these only for levels 0 to $n+1$. It is however not clear a priori if this is automatically a $\VB$ $n$-groupoid. This is discussed in Section \ref{sec:n-dual-construction}.

\begin{introthm}
    The $n$-dual of a higher vector bundle is given levelwise by a limit in $\VB$, and it is determined by its levels 0 to $n+1$.
\end{introthm}

Consider for the moment the case of a $\VB$ $n$-groupoid $\huaV \to \huaG$. Assuming that the $n$-dual of $\huaV$ exists, we find that the canonical evaluation pairing is nondegenerate in the homology of the higher vector bundle. The homology of a higher vector bundle is discussed after Definition \ref{def:VB-norm-cplx}, and the homological nondegeneracy condition in Definition \ref{def:VB-hndg-pairing}. This is the same nondegeneracy condition that \textit{$n$-shifted symplectic structures} are required to satisfy in \cite{Lesdiablerets}, where this notion was translated to the smooth setting from the work of \cite{PantevToenVaquieVezzosi2013}. This is the content of Theorem \ref{thm:VBndual-pairing-hndg}, which consists of a slightly more general statement than the one reported here.

\begin{introthm}
    Let $\huaV \to \huaG$ be a $\VB$ $n$-groupoid over a Lie $n$-groupoid $\huaG$, for which the $n$-dual $\huaV^{n*}$ exists. Then the $n$-dual pairing $\langle\cdot , \cdot \rangle: \huaV^{n*}_n \otimes \huaV_n \to \R \times \huaG_n$ is homologically $n$-shifted nondegenerate.
\end{introthm}

At this point we focus on the $n=2$ case and provide a concrete model of the $\VB$ 2-dual of a $\VB$ 2-groupoid. This is one of the main accomplishments of this thesis. We show the following properties in Section \ref{sec:2-dual-VB-gpd}.

\begin{introthm}
    Let $\huaV\to \huaG$ be a $\VB$ 2-groupoid. Its $\VB$ 2-dual $\huaV^{2*}$ exists and it is a $\VB$ 2-groupoid. 
\end{introthm}

\begin{introthm}
    A 2-shifted simplicial pairing $\alpha: \huaV\otimes \huaW \to B^2\R\times \huaG$ between two $\VB$ 2-groupoids $\huaV$ and $\huaW$ over the same base $\huaG$ is homologically nondegenerate if and only if it induces weak equivalences
    \begin{equation*}
        \alpha^l:\huaW \to \huaV^{2*}, \quad \alpha^r:\huaV \to \huaW^{2*}.
    \end{equation*}
\end{introthm}

This provides as a characterization of the nondegeneracy of 2-shifted symplectic structures over Lie 2-groupoids. This result is implied by the following one, which we prove in Theorem \ref{thm:we-of-VB2gpd-is-qi}. 

\begin{introthm}
    A simplicial bundle map $f:\huaV \to \huaW$ between two $\VB$ 2-groupoids over $\huaG$ over the identity of $\huaG$ is a weak equivalence if and only if it induces a quasi-isomorphism between their normalized complexes. 
\end{introthm}

This is an extension to $n=2$ of \cite[Thm. 3.5]{HoyoOrtiz2020} in the case of a map over the identity on the base. 

In the last chapter, we consider the 2-cotangent of a Lie 2-groupoid, that is the 2-dual of its tangent groupoid. We show that the definition of the tautological form given in \eqref{eq:intro-taut-form} goes through and we obtain a canonical 2-shifted symplectic structure in Theorem \ref{thm:2Cota-symp-form}.

\begin{introthm}
    Let $\huaG$ be a Lie $2$-groupoid, and $T^{2*}\huaG$ its 2-cotangent. The canonical form $\omega := -d\theta$, where $\theta \in \Omega^1(T^{2*}\huaG_2)$ is given as in \eqref{eq:intro-taut-form} with respect to the 2-dual pairing, is a 2-shifted symplectic form. 
\end{introthm}

For the 2-cotangent of a Lie 1-groupoid $\huaG$, we find a symplectic morita equivalence between $T^{2*}\huaG$ and  $\widebar{T^{1*}\huaG}$, that is, the 2-shifted symplectic $\VB$ 2-groupoid obtained by applying the Artin-Mazur bar to the 1-cotangent $T^{1*}\huaG$. 

\begin{introthm}
    Let $\huaG$ be a Lie groupoid and $(\kappa, \sigma)$ an Ehresmann connection on $\huaG$. Then, the map 
    \begin{equation*}
        \Xi: T^{2*}\huaG \to \overline{T^{1*}\huaG} 
    \end{equation*}
    defined by \eqref{eq:2Cota-1Cota-ME-def} is a symplectic Morita equivalence.
\end{introthm}

The $\VB$ 2-groupoid $\widebar{T^{1*}\huaG}$ appeared in \cite{MehtaTang2011} as the 2-groupoid integrating the Courant algebroid $A\oplus A^*$, that is the double of the Lie algebroid $A$ of $\huaG$ (See \cite{LiuWeinsteinXu1997}). So Theorem 7 frames $T^{2*}\huaG$ as an \textit{integration up to homotopy} of the Courant algebroid $A\oplus A^*$. We discuss this point of view further in Chapter 4. It is worth noting that the explicit Morita equivalence we are able to construct depends on a choice of Ehresmann connection on the groupoid, while both $T^{2*}\huaG$ and $\widebar{T^{1*}\huaG}$ are independent of choices.

In proving many of these results we use techniques that we import from the simpler case of simplicial vector spaces, which are higher vector bundles over the point.\footnote{All simplicial vector spaces are $\infty$-groupoids by a theorem of Moore that we recall in Theorem \ref{thm:MooreHornFillers}.} 
The link that allows us to import these methods from the theory of simplicial vector spaces to that of higher vector bundles over $\huaG$ lies in the fact that, by using the unique embedding of $\huaG_0$ into any other level of $\huaG$ given by the total degeneracy
\begin{equation*}
    1=(s_0)^k: \huaG_0 \to \huaG_k, 
\end{equation*} 
a higher vector bundle $\huaV$ over $\huaG$ can be pulled back to a higher vector bundle $1^*\huaV$ over the 0-groupoid $\huaG_0$. This is a ``bundle of simplicial vector spaces'': its simplicial structure is defined pointwise over $\huaG_0$.

The $\VB$ duality theory for $n>1$ is much more complex than for $n=0$ and $n=1$, starting from the fact that existence of $n$-duals is, a priori, unclear. Many of these complexities already manifest themselves in the case of simplicial vector spaces. However, this particular case is simpler, as existence of $n$-duals is guaranteed for all simplicial vector spaces by the existence of a right adjoint to the tensor product: the \textit{mapping space} functor, see e.g. \cite{GoerssJardine2009}. 
Thus, the universal property defining the $n$-dual over the point is satisfied by the mapping space $\IHom(\huaV, B^n\R)$. As we show in Theorem \ref{thm:hom-ngpd}, this is an $n$-groupoid. From now on, we call $n$-groupoids in the category of vector spaces $\VS$ $n$-groupoids.

The case of simplicial vector spaces gives a clear picture of the structure of $n$-duals, a model of which can be obtained by solving the equations computing $\IHom(\huaV, B^n\R)$. 
We discuss the structure of these equations in Section \ref{sec:computations} and we solve them explicitly in order to compute $n$-duals of $\VS$ $n$-groupoids for $n=1,2$. For $n=1$ we recover the $\VB$ 1-dual of \cite{Pradines1988} and for $n=2$ we describe a new object, which we call the 2-dual.

The aforementioned homology theory of simplicial vector spaces comes from the Dold-Kan correspondence (see e.g. \cite{Weibel1994, GoerssJardine2009}), which gives an equivalence of categories
\begin{equation*}
    \SVect \newrightleftarrows{N}{DK} \Ch_{\ge 0}(\Vect).
\end{equation*}
The functor $N$ is the \textit{normalized complex} functor, which associates to a simplicial vector space $\huaV$ the non-negative chain complex
\begin{equation*}
    N(\huaV)_k = \bigcap\limits_{i=0}^{k-1}\ker d_i^k, \qquad \partial_k = (-1)^k d_k^k.
\end{equation*}
The homology of $\huaV$ is defined as the homology of $N(\huaV)$. 

Additionally, the Dold-Kan correspondence is an equivalence between the homotopy categories, and in particular, a simplicial linear map $f:\huaV \to \huaW$ is a weak equivalence, if and only if $N(f): N(\huaV)\to N(\huaW)$ is a quasi-isomorphism. 
In other words, all the information of a $\VS$ $n$-groupoid, including its homotopy type, is encoded in its normalized complex. Thus, if the homology of $\huaV$ vanishes above degree $n$, we say that $\huaV$ is a $\VS$ $n$-type. Naturally any $\VS$ $n$-groupoid is an $n$-type.

There is a problem though, in that the Dold-Kan correspondence is not monoidal on the nose, as it does not preserve the tensor product on each side, but it does so only up to homotopy equivalence. This is the content of the Eilenberg-Zilber theorem, in its original version appearing in \cite{EilenbergZilber1953,EilenbergMacLane1954}:

Let $\huaV$ and $\huaW$ be two simplicial vector spaces. There is a natural chain homotopy equivalence given by the Eilenberg-Zilber map $EZ$ and the Alexander-Whitney map $AW$,
\begin{equation*}
    N(\huaV) \otimes N(\huaW) \newrightleftarrows{EZ}{AW} N(\huaV \otimes \huaW),
\end{equation*}
such that 
\begin{equation*}
    AW \circ EZ = id_{N(\huaV) \otimes N(\huaW)}, 
    \qquad 
    EZ \circ AW \sim id_{N(\huaV \otimes \huaW)}.
\end{equation*}

With this homotopy equivalence, we can translate an $n$-shifted pairing of two simplicial vector spaces into an $n$-shifted \textit{IM-pairing}\footnote{The name stands for ``infinitesimally multiplicative''. This terminology appears for example in \cite{BursztynCabrera2012}, \cite{BursztynDrummond2019} referring to objects previously studied in \cite{BursztynCrainicWeinsteinZhu2004}. See \cite{CuecaZhu2023} for an explanation of its use in this context.} between the normalized complexes:
\begin{equation*}
    \alpha: \huaV \otimes \huaW \to B^n\R \quad \rightsquigarrow \quad \lambda_\alpha = EZ^*\circ N(\alpha): N(\huaV) \otimes N(\huaW) \to \R[-n].
\end{equation*}
This is the same operation used to obtain the IM-pairing on the tangent complex associated to an \textit{$n$-shifted symplectic} form in \cite{Lesdiablerets}.
This allows us to define the homological nondegeneracy condition for any pairing of simplicial vector spaces, as in Definition \ref{def:homological-non-deg}.

Of course since the Eilenberg-Zilber theorem affects the tensor product, it also affects its adjoint, the mapping space. In one of our main results (Theorem \ref{thm:EZ-thm-hom}), we make this precise, by reformulating the theorem in terms of mapping spaces. 

\begin{introthm}
    There is a natural chain homotopy equivalence
    \begin{equation*}
        N(\IHom(\huaV, \huaW)) \newrightleftarrows{EZ^H}{AW^H} \IHom_{\ge 0}(N(\huaV), N(\huaW)),
    \end{equation*}
    such that
    \begin{equation*}
        EZ^H \circ AW^H = id_{\IHom_{\ge 0}(N(\huaV), N(\huaW))}, 
        \qquad 
        AW^H \circ EZ^H \sim id_{N(\IHom(\huaV, \huaW))}.
    \end{equation*}
\end{introthm}

In the specific case of the $n$-dual, this gives a chain homotopy equivalence
\begin{equation*}
    N(\huaV^{n*}) \newrightleftarrows{EZ^H}{AW^H} N
    _{\ge 0}(\huaV)^*[-n].
\end{equation*}
Because of the truncation that happens at 0 in the right side, if $N(\huaV)^*[-n]$ has non-zero negative degrees, there is some loss of information when taking the $n$-dual. 
When $\huaV$ is a $\VS$ $n$-groupoid, its normalized complex is concentrated in degrees 0 to $n$, which means that $N(\huaV)^*[-n]$ is concentrated in degrees 0 to $n$ as well, and it coincides with its zero truncation. 
Even when $\huaV$ is not a $\VS$ $n$-groupoid on the nose but it is homotopy equivalent to one (i.e. it is an $n$-type), there is no loss of homotopy information, as its homology is still concentrated in degrees $0$ to $n$. This is expressed by Theorem \ref{thm:ndual-pairing-hom-nondeg}.

\begin{introthm}
Let $\huaV$ be a simplicial vector space. The $n$-dual pairing $\langle \cdot, \cdot \rangle: \huaV^{n*} \otimes \huaV \to B^n\R$ is homologically $n$-shifted nondegenerate if and only if $\huaV$ is at most an $n$-type.
\end{introthm}

We observe that this is the right nondegeneracy condition, since the canonical 2-dual pairing is degenerate on the nose, but nondegenerate in homology. 

As an immediate consequence, we characterize nondegeneracy of pairings in terms of homotopy instead of homology (Theorem \ref{thm:hom-nondeg-homotopy-equiv}).

\begin{introthm}
An $n$-shifted simplicial pairing $\alpha$ between $\VS$ $n$-types $\huaV$ and $\huaW$ is homologically nondegenerate if and only if it induces a weak equivalence
\begin{equation*}
    \alpha^r: \huaV \to \huaW^{n*},
\end{equation*}
or equivalently, a weak equivalence
\begin{equation*}
    \alpha^l: \huaW \to \huaV^{n*}.
\end{equation*}
\end{introthm}

A main application of the above result is that $n$-duality is reflexive up to homotopy (Theorem \ref{thm:ndual-reflexive-uth}). 

\begin{introcor}
    Let $\huaV$ be at most an $n$-type. The double $n$-dual $(\huaV^{n*})^{n*}$ is weak equivalent to $\huaV$ itself. That is
    \begin{equation}
        \langle \cdot, \cdot \rangle^l: \huaV \to (\huaV^{n*})^{n*}
    \end{equation}
    is a weak equivalence. 
\end{introcor}

\subsection*{Notation and conventions}

We use the homological convention, in which the Dold-Kan correspondence is an equivalence between simplicial vector spaces and non-negative chain complexes. 
We orient edges of simplices from higher to lower number, ``as the water flows'', thus the source of an arrow is its $d_0$ face, while its target is its $d_1$ face. We also write boundaries and horns as tuples ordered by increasing face index. This has the effect that our version of the nerve of a category (Example \ref{ex:nerve-cat}) might appear to be flipped by the front-to-back symmetry (Remark \ref{rem:front-to-back}) with respect to other perhaps more common conventions. We sometimes use the terms weak equivalence and Morita equivalence interchangeably, when talking about weak equivalences of higher groupoids. Other conventions will be explained as they appear.

\subsection*{Disclaimers}

This work incorporates an edited and corrected version of the preprint \cite{RonchiZhu2024}, coauthored with my supervisor Chenchang Zhu. Her input on the preprint version of this work consisted mainly of supplying ideas and techniques during discussions and reviewing and correcting my writing, while I completed the main body of work. As required by the doctoral degree regulations, each such contribution is explicitly stated wherever it appears. This material is mostly concentrated in Theorem \ref{thm:hom-ngpd}, Section \ref{sec:simp-vect} and Chapter 2.

\subsubsection*{Declaration on the use of ChatGPT and comparable tools in the context of examinations}

In this thesis, I have used Deepseek to proofread and suggest alternative formulations for some sentences in the introduction to Chapter 1 (<1\% of the full manuscript). 
I hereby declare that I have stated all uses completely.
Missing or incorrect information will be considered as an attempt to cheat.


\chapter*{Acknowledgments}
\markboth{Acknowledgments}{}
\addcontentsline{toc}{chapter}{Acknowledgments}

There is a long list of people this work would have not been possible without. 

I am profoundly grateful to Chenchang Zhu and Madeleine Jotz for granting me the privilege of conducting this research and for their continued guidance and support through this long journey. 
I have learned a lot from their expertise and insights, and greatly appreciated their kindness through difficult times. 

I am also grateful to Ralf Meyer and Thomas Schick for helping me navigate the PhD and for their support when I needed some more time to bring this project to its conclusion. 

I cannot understate the importance of Miquel Cueca's contribution to this research.
His dedication and constant willingness to discuss math and life has been invaluable in shaping my work and perspective. 
I greatly appreciated our countless discussions together.
Most of the material in this thesis would not have been written, had he not approached me with an idea for a quick ``two-week'' exercise. 
I am hugely indebted to him and his patience for my lack of hard-drive space. 

Special thanks to Hao Xu for pointing us to the mapping space construction for simplicial vector spaces and to Matias del Hoyo for insisting that this is indeed the right way. 

There are so many people I have met during my years in Göttingen, and at conferences in this time, that have encouraged me and generously provided me with extremely valuable ideas, perspective, and fun times. 
I would like to start by thanking my colleagues and office mates and especially Camilo Angulo, Florian Dorsch, Ilias Ermeidis, Zhicheng Han, Rok Havlas, Arne Hofmann, David Kern, Kalin Krishna, Rosa Marchesini, Leonid Ryvkin, Hao Xu and Fabrizio Zanello for having made this, to quote my friend Ilias, an \textit{amaaaazing} experience. 
Cheers to all the new friends I have made in Göttingen and all the old ones that I have found again. 
(And more special thanks to Kalin for printing and handing in this manuscript for me.)

I am very grateful to Christian Blohmann, Damien Calaque, Matias del Hoyo, Ezra Getzler, Domenico Fiorenza, Marco Gualtieri, Camille Laurent-Gengoux, John Pridham, Michele Schiavina, Mathieu Stiénon, Giorgio Trentinaglia, Luca Vitagliano, and Ping Xu, for their interest, encouragement, constructive feedback, and many helpful discussions.

I thank Insitut Mittag-Leffler, where the writing of this thesis was completed, and the organizers of the wonderful program ``Cohomological Aspects of Quantum Field Theory'' Francesco Bonechi, Alberto Cattaneo, Nikita Nekrasov, and Maxim Zabzine, for giving me the opportunity to participate as a Junior Fellow. 
I have met some amazing people at this program and I thank them for their kind support in the final stages of this writing. Special thanks particularly to Leon Menger, for helping with last minute proofreading, and Eugenia Boffo and Giovanni Canepa, for making sure there was always food and good mood in the yellow house. 

I would also like to thank my teachers and professors through the years, and especially Anna Coen for her continued support and encouragement. 

I am forever grateful to my mother and my brother for their absolute support and always being up for a call, and to my father, to whom I owe much of what I am as a person and who is always with me, despite leaving this world much too soon.

I send my heartfelt thanks to the Mansilla Urban family, as well, for their warm hospitality and kindness. Lastly, but perhaps most importantly, my deepest thanks go to Anna. I am so unbelievably lucky to have met you and be able to enjoy your light every day.

\tableofcontents
\addcontentsline{toc}{chapter}{Contents}
\clearemptydoublepage

\mainmatter

\chapter{Background Notions}

In this chapter, we discuss the foundations of the theory of higher groupoids, higher vector spaces and higher vector bundles. Most of the material is a re-elaboration of more or less established literature which we reference as it occurs, except for a couple of original results which we highlight as they arise. For the sake of brevity we omit most specific references from this summary, as they are provided in the corresponding sections.

Throughout the chapter, we mostly consider two kinds of categories: categories such as $\Set$ and $\Vect$ which have all categorical limits and colimits (i.e. they are complete and cocomplete), and categories such as those of smooth manifolds $\Mfd$ and vector bundles $\VB$, which are not complete nor cocomplete. Naturally, the theory for the second kind of category is complicated by the (a priori) absence of certain limits, and we review a series of techniques to address this. 

In Section \ref{sec:simp-obj} we begin with a basic discussion of simplicial objects in a general category. We focus on the combinatorial structure of the standard simplex, which is necessary for computing $n$-duals. 
We then review some key structures on the category of simplicial sets: the monoidal structure given by the levelwise product, and the mapping set (or function complex), which is a simplicial set with simplicial maps as points, simplicial homotopies between maps as 1-simplices and higher homotopies as higher simplices. These two functors form a tensor-hom adjunction, making the category of simplicial sets a category enriched over itself. We frame this in terms of a general method to construct categories enriched over simplicial sets, which we later apply to simplicial vector spaces.

In the lengthy Section \ref{sec:higher-gpds} we introduce the simplicial model of higher groupoids, beginning with $n$-groupoids in the category of simplicial sets. A fundamental idea in this theory is the concept of an \textit{object of simplicial diagrams}, which collects arrangements of simplices of a certain shape within a simplicial object. 

In the first subsection we consider the setting of simplicial sets, and define $n$-\hspace{0pt}groupoids by imposing the Kan conditions. We then recall that Kan simplicial sets ($\infty$-groupoids) admit homotopy groups. These homotopy invariants lead to the notion of weak homotopy equivalences, which allow one to compare different $n$-groupoids up to homotopy. In this context, an $n$-groupoid represents a \textit{homotopy $n$-type}, i.e. a space with trivial homotopy groups for $m > n$. We conclude the discussion of $n$-groupoids in $\Set$ by recalling the fact that any mapping space with target a Kan simplicial set is also Kan and upgrading this statement (by assembling various results in the literature) to the fact that any mapping space with $n$-groupoid target is an $n$-groupoid. This is the content of Theorem \ref{thm:hom-ngpd}. 

In the second subsection, we discuss objects of simplicial diagrams for simplicial manifolds. As these are obtained by computing categorical limits, they might not exist as manifolds, but we show how they can be studied by considering them first as sheaves on $\Mfd$ with the pretopology of surjective submersions. With this, we state the Kan conditions that define Lie $n$-groupoids. 

In the third subsection we recall a framework to study representability of objects of simplicial diagrams and how this can be used to show that the horn spaces appearing in the Kan conditions are representable. 
We then use this to produce a new result computing the dimension of the horn spaces by writing them as fiber products over surjective submersions. This is Theorem \ref{thm:dim-gen-horns} and Corollary \ref{cor:dim-horn-spaces}. 

In the fourth subsection we discuss weak equivalences of Lie $n$-groupoids from the perspective of $n$-groupoids in a general category with a pretopology and a jointly conservative collection of point functors. This is a device that allows one to show two objects are isomorphic if and only if there is a bijection between them after applying each point functor. This is useful in concrete categories where the forgetful functor does not reflect isomorphisms, such as in $\Mfd$. We employ this to define weak equivalences \textit{stalkwise}, i.e. two objects are weakly equivalent if and only if they are so in $\SSet$ after applying each functor of points. We then describe a more practical characterization of weak equivalences which recovers the notion of Morita maps and Morita equivalences present in the literature on Lie groupoids. 

In the final subsection, we discuss how to recover a Lie $n$-groupoid from a finite set of data, interpreting it as a space with certain multiplication and inversion maps. To do this on a general category with a pretopology we introduce the coskeleton of a simplicial object, which is also a levelwise limit (so a priori not representable), and use a collection of points to prove results about it. 

In Section \ref{sec:simp-vect} we study simplicial vector spaces in more detail and review the theorem of Moore which states that they are automatically Kan complexes by virtue of being simplicial groups. In preparation to discuss the classical Dold-Kan correspondence between simplicial vector spaces $\SVect$ and non-negative chain complexes $\Ch_{\ge 0}(\Vect)$, we recall the structures of $\SVect$, $\Ch(\Vect)$, and its subcategory $\Ch_{\ge 0}(\Vect)$ as monoidal categories enriched over themselves. This gives mapping spaces for each of these categories. By considering maps up to homotopy, each of these categories has a homotopy category. 
We then recall that the Dold-Kan correspondence is not only an equivalence of categories through the normalized complex functor, but it extends to an equivalence of the homotopy categories. As such, weak equivalences of simplicial vector spaces are the same as quasi-isomorphisms between the respective normalized complexes, and homotopy groups of a simplicial vector space are the same as the homology groups of its normalized complexes. As such $\VS$ $n$-groupoids can be seen as \textit{linear $n$-types}. 
We conclude this section by showing that $\VS$ $n$-groupoids have a more restrictive finite data than general $n$-groupoids, in the sense that the Moore theorem prescribes canonical multiplication maps for each $\VS$ $n$-groupoid. This is compatible with the fact that a simplicial vector space is an $n$-groupoid if and only if its normalized complex is concentrated in degrees 0 to $n$, as a consequence of the Dold-Kan correspondence. 

In Section \ref{sec:higher-vb} we close the chapter by introducing higher vector bundles over Lie $n$-groupoids. These are simplicial objects in the category of vector bundles which are themselves Lie $n$-groupoids for some $n$ (possibly $\infty$). Their theory is simplified by the fact that surjective submersions of vector bundles are the same as fiberwise linear maps over surjective submersions that are also fiberwise surjective. The important fact which links them to the theory of simplicial vector spaces is that any higher vector bundle over the Lie $n$-groupoid $\huaG$ can be pulled back by the total degeneracy map to a higher vector bundle over the 0-groupoid $\huaG_0$. This can then be seen as a ``bundle of simplicial vector spaces'', as its simplicial structure is defined entirely pointwise over $\huaG_0$.
After recalling the construction of the dual $\VB$ 1-groupoid and the cotangent of a Lie 1-groupoid, we define the normalized complex of a simplicial vector bundle, as the fiberwise normalized complex of the pullback by the total units. This allows us to consider $\VB$ $n$-groupoids as \textit{pointwise $n$-types} and will allow us to prove many results in the $\VB$ case by using the same theory as for simplicial vector spaces applied pointwise. 
We conclude the section by taking this correspondence further and showing our main result of this chapter. This is a partial extension of a theorem relating weak equivalences of $\VB$ 1-groupoids and quasi-isomorphisms of their normalized complexes from \cite{HoyoOrtiz2020}, to one relating weak equivalences of $\VB$ 2-groupoids (over the identity on the base) to quasi-isomorphisms of their normalized complexes. 
This is the content of Theorem \ref{thm:we-of-VB2gpd-is-qi}.

\section{Simplicial objects}\label{sec:simp-obj}

A groupoid of local transformations is commonly pictured as a set of 1-\hspace{0pt}dimensional oriented edges (arrows) between points in a set. Higher categorical information comes in the form of transformations between transformations, which geometrically translates to the existence of non-trivial higher-dimensional cells.
In the simplicial formalism, higher-dimensional cells take the shape of $n$-simplices.
The geometry of simplices is encoded in the combinatorics of the category $\Delta$ of isomorphism classes of finite ordered sets  
\begin{equation*}
    [0]:= \{0\}, \quad [1]:=\{0 \le 1\}, \quad \dots, \quad [n]:= \{0 \le 1 \le \dots \le n \}, \quad \dots
\end{equation*}
and weakly order-preserving maps between them. This leads to the theory of simplicial sets and simplicial objects.
Standard references on the subject are \cite{May1967,GoerssJardine2009}. For a pedagogical introduction to simplicial sets see \cite{Friedman2012}. 

\begin{definition}
A \textbf{simplicial object} in a category $\Cat$ is a functor $\Delta^{op} \to \Cat$. The category of simplicial objects and simplicial maps in $\Cat$ is denoted by $\Simp\Cat = \Fun(\Delta^{op}, \Cat)$. 

A \textbf{simplicial set} is a functor $\Delta^{op} \to \Set$. 
We call the natural transformations between such functors \textbf{simplicial maps}. 
The category of simplicial sets and simplicial maps is denoted by $\SSet := \Fun(\Delta^{op}, \Set)$. 

Analogously, \textbf{simplicial vector spaces}, \textbf{simplicial manifolds} and \textbf{simplicial vector bundles} are simplicial objects in the categories $\Vect$ of vector spaces with linear maps, $\Mfd$ of smooth manifolds with smooth maps, and $\VB$ of smooth vector bundles with smooth bundle maps, respectively. 
\end{definition}

\begin{definition}
A \textbf{concrete category} $\Cat$ is a category equipped with a faithful functor $\Forget: \Cat \to \Set$ called the \textbf{forgetful functor} which returns the so-called \textbf{underlying set} of each object.

When the meaning is clear from context, we omit $\Forget$ and denote the underlying set of an object by the same name as that object. 
\end{definition}

\begin{remark}\label{rem:repbly-concrete-cat-limits}
The categories $\Set$, $\Vect$, $\Mfd$, and $\VB$ are concrete categories in which the forgetful functor is representable, and therefore preserves limits. (See e.g. \cite[Thm 3.4.6]{Riehl2016}). 
In fact, when $\Cat$ is $\Set$ or  $\Mfd$, then $\Forget$ is represented by the point $pt$. The forgetful functor for vector spaces $\Vect$ is $\Vect(\R, \_)$, while that for vector bundles with bundle maps $\VB$ is $\VB(\R \to pt, \_)$. 
\end{remark}

\begin{remark}
The forgetful functor of a concrete category extends immediately to a forgetful functor $\Forget: \Simp\Cat \to \SSet$, by postcomposition of every simplicial object in $\Cat$ with the forgetful functor of $\Cat$. We call this the underlying simplicial set of a simplicial object. 
\end{remark}

\begin{example}
The most immediate examples of simplicial sets are given by the hom-functors of the category $\Delta$. We call these \textbf{standard simplices} and denote them by 
\begin{equation*}
    \Delta^n := \Delta(\_, [n]).
\end{equation*}
\end{example}

The structure of simplicial objects is made more concrete by the fact that the category $\Delta$ is generated by just two classes of maps. 

\begin{definition}
For any $n \ge 0$ and $0\le j \le n$, the \textbf{coface map} $\delta^n_j: [n-1] \to [n]$ is the unique order-preserving injection whose image does not include $j$, i.e. $\delta^n_j$ ``skips $j$''. That is 
\begin{equation}\label{eq:Delta-coface-def}
    \delta^n_j(k) = \begin{cases}
        k &\text{if } k < j,\\
        k + 1 &\text{if } k \ge j.
    \end{cases}
\end{equation}

For any $n \ge 0$ and $0 \le i \le n$, the \textbf{codegeneracy map} $\sigma^n_i: [n+1] \to [n]$ is the unique order-preserving surjection where $i$ has two preimages, i.e. $\sigma^n_i$ ``repeats $i$''. That is,
\begin{equation}\label{eq:Delta-codeg-def}
    \sigma^n_i(k) = \begin{cases}
        k &\text{if } k < i,\\
        i &\text{if } k = i, i + 1,\\
        k-1 &\text{if } k > i + 1.
    \end{cases}
\end{equation}
\end{definition}

\begin{lemma}\label{lem:Delta-gen-by-cf-cd}
The morphisms in $\Delta$ are generated by the coface and codegeneracy maps \eqref{eq:Delta-coface-def}, \eqref{eq:Delta-codeg-def} subject to the \textbf{cosimplicial identities}
\begin{equation}\label{eq:cosimp-id}
    \begin{array}{lll}
        \delta^{n}_j \delta^{n-1}_i  = & \delta^n_i \delta^{n-1}_{j-1} &\text{if}\; i<j,  \\
        \sigma^{n-1}_j \sigma^{n}_i =& \sigma^{n-1}_i \sigma^{n}_{j+1} & \text{if}\; i\leq j,
        \end{array}\qquad  \sigma^{n-1}_j \delta^m_i =\left\{\begin{array}{ll}
        \delta^{n-1}_i \sigma^{n-2}_{j-1}  & \text{if}\; i<j, \\
        \id  & \text{if}\; i=j, j+1,\\
        \delta^{n-1}_{i-1} \sigma^{n-2}_j & \text{if}\; i> j+1.
        \end{array}\right.
\end{equation}

In particular, for any order-preserving map $f: [m] \to [n]$ there are unique multi-indices $I$ and $J$ of lengths $|I|$ and $|J|$ with $n = m - |I| + |J|$
such that 
\begin{equation*}
    0 \le i_1 < i_2 < \dots < i_{|I|} \le m-1, \qquad
    0 \le j_1 < j_2 < \dots < j_{|J|} \le n,
\end{equation*}
and $f$ can be written as a composition 
\begin{equation}
    f = \delta_{\widebar{J}} \sigma_I := \delta_{j_{|J|}} \circ \dots \delta_{j_1} \circ \sigma_{i_1} \circ \sigma_{i_{|I|}},
\end{equation}
where $\widebar{J}$ denotes the opposite multi-index with $\widebar{j}_k = j_{|J|-k}$.
\end{lemma}

A proof of this result can be found in \cite[Prop. VII.5.2]{MacLane1978}.\footnote{The notation in \cite{MacLane1978} is slightly different from ours, in that the empty set is included in the category $\Delta$ and called $0$.}

This immediately implies the following characterization of simplicial objects, which is often used as an alternative definition.

\begin{prop-def}\label{prop-def:simp-obj}
The data of a simplicial object $\huaX \in \Cat$ consists of: 
\begin{enumerate}[label=(\roman*)]
    \item A family of objects $\huaX_n:= \huaX([n])$ indexed by $n \ge 0$, which we call the \textbf{$n$-th level of $\huaX$}. Elements of $\huaX_n$ are called \textbf{$n$-simplices of $\huaX$}. 
    \item A family of morphisms $d_i^n := \huaX(\delta_i^n): \huaX_n \to \huaX_{n-1}$ indexed by $n\ge 0$, $0\le i\le n$, called \textbf{face maps}.
    \item A family of morphisms $s_i^n:= \huaX(\sigma_i^n): \huaX_n \to \huaX_{n+1}$ indexed by $n \ge 0$, $0 \le i \le n$ called \textbf{degeneracy maps}. Any $n$-simplex in the image of one of the degeneracy maps $s_i$ is called \textbf{degenerate}.
\end{enumerate}
This data is subject to the \textbf{simplicial identities}:
\begin{equation}\label{eq:simp-id}
    \begin{array}{lll}
    d^{n-1}_i d^{n}_j = & d^{n-1}_{j-1} d^m_i &\text{if}\; i<j,  \\
    s^{n}_i s^{n-1}_j =& s^{n}_{j+1} s^{n-1}_i & \text{if}\; i\leq j,
    \end{array}\qquad  d^m_i s^{n-1}_j =\left\{\begin{array}{ll}
    s^{n-2}_{j-1} d^{n-1}_i  & \text{if}\; i<j, \\
    \id  & \text{if}\; i=j, j+1,\\
    s^{n-2}_j d^{n-1}_{i-1} & \text{if}\; i> j+1.
    \end{array}\right.
\end{equation}

Furthermore, a simplicial map $f:\huaX \to \huaY$ between simplicial objects in $\Cat$ consists of a family of morphisms $f_n: \huaX_n \to \huaY_n$ in $\Cat$ indexed by $n\ge 0$, which commute with face and degeneracy maps at each level. That is\footnote{As usual, we omit the upper indices from the face and degeneracy maps, as they are mostly clear from context.}
\begin{equation*}
\begin{split}
    d_i f_n &= f_{n-1} d_i, \quad \text{for any } n\ge 0, \, 0\le i \le n,\\
    s_i f_n &= f_{n+1} s_i, \quad \text{for any } n\ge 0, \, 0\le i \le n.
\end{split}
\end{equation*}
\end{prop-def}

\begin{example}[The nerve of a category]\label{ex:nerve-cat}
Let $\huaC$ be a small category with set of objects $\huaC_0$ and set of morphisms $\huaC_1$. 
We denote its source and target maps by $d_0:\huaC_1 \to \huaC_0$ and $d_1:\huaC_1 \to \huaC_0$, respectively, and the map associating to every object in $\huaC_0$ the identity at that object by $1: \huaC_0 \to \huaC_1$.  
The \textbf{nerve} of $\huaC$ is the simplicial set denoted by the same name and defined levelwise by having $\huaC_0$ as level 0, and
\begin{equation*}
    \huaC_m = \{ (f_1, \dots, f_m) \in \huaC_1^{\times m} \mid d_1f_{i} = d_0 f_{i+1}, \forall i\ge 1\}
\end{equation*}
as level $m$ for any $m \ge 0$. 
That is $\huaC_m$ is the set of $m$-tuples of composable arrows. 
In our convention, the face maps of $\huaC$ are
\begin{equation*}
    d_i(f_1, \dots, f_m) = \begin{cases}
        (f_1, \dots, f_{m-1}) &\text{ if } i=0,\\
        (f_2, \dots, f_{i-1}, f_{i+1} \circ f_i, f_{i+2}, \dots, f_m) &\text{ if } 0 < i < m,\\
        (f_2, \dots, f_m) &\text{ if } i=m,
    \end{cases}
\end{equation*}
and the degeneracy maps are 
\begin{equation*}
    s_i(f_1, \dots, f_{n-i}, 1_{d_1f_{n-i}}, f_{n-{i+1}}, \dots, f_n). 
\end{equation*}

The same construction applies in the same way to categories internal to a general category $\Cat$, in which case, it yields a simplicial object in $\Cat$.
\end{example}

\begin{remark}\label{rem:canonical-form-s-d-comp}
It follows from Lemma \ref{lem:Delta-gen-by-cf-cd} that
any composition of faces and degeneracies $f: \huaX_n \to \huaX_m$ between different levels of a simplicial object $\huaX$ can be written uniquely as a composition 
\begin{equation}\label{eq:sd-comp-normal-form}
    f = s_{\widebar{I}} d_J 
    := s_{i_{|I|}} \circ \dots \circ s_{i_1} \circ d_{j_1} \circ \dots \circ d_{j_{|J|}},
\end{equation}
for unique multi-indices $I$ and $J$ of lengths $|I|$ and $|J|$ with $m = n + |I| - |J|$ such that 
\begin{equation*}
    0 \le i_1 < i_2 < \dots < i_{|I|} \le m-1, \qquad
    0 \le j_1 < j_2 < \dots < j_{|J|} \le n.
\end{equation*}
This result is obtained in practice by repeatedly applying the simplicial identities \eqref{eq:simp-id} to rewrite any given composition of faces and degeneracy maps in the normal form \eqref{eq:sd-comp-normal-form}.
\end{remark}

\begin{remark}[Front-to-back duality {\cite[\S 8.2.10]{Weibel1994}}{\cite[\href{https://kerodon.net/tag/003L}{Tag 003L}]{kerodon}}]\label{rem:front-to-back}
Every simplicial object $\huaX$ in a category $\Cat$ has an opposite $\huaX^{op}$, which consists of the same spaces but has simplicial maps $(d^{op})^n_i = d^n_{n-i}$, $(s^{op})^n_i = s^n_{n-i}$. 
This can be understood as reversing the order of the elements in each object of the category $\Delta$, or equivalently reversing the direction of all the edges and orienting them from lower to higher vertex number. 
When the simplicial object in question is the nerve of a category, its opposite is precisely the nerve of the opposite category \cite[\href{https://kerodon.net/tag/003Q}{Tag 003Q}]{kerodon}. 
The existence of opposite simplicial objects extends the well-known principle of categorical duality (see e.g. \cite[\S 1.2]{Riehl2016}) to this context.
In particular, if a property holds for any simplicial object under hypotheses symmetric with respect to the opposite construction, it must hold for its opposite as well, which can be helpful in certain proofs. 
We will make use of this property in the computation of 1- and 2-duals in Section \ref{sec:computations}. 
\end{remark}

\subsection{Standard simplices and generated subsets}

Each level $\Delta^n_m$ of the standard $n$-simplex $\Delta^n$ is the set of order-preserving maps $\Delta([m], [n])$ between $[m]$ and $[n]$. 
Each of these can be identified by the sequence of its images in $[n]$, which has length $m+1$: any $u \in \Delta^n_m$ corresponding to $u:[m] \to [n]$, can be written as 
\begin{equation*}
    u = u_0 u_1 \dots u_{m-1} u_{m} 
\end{equation*}
with $u_i =u(i)$ for all $i \in[m]$, and $0 \le u_0 \le u_1 \le \dots \le u_m \le n$.
For example, for $m=n$ we identify $id:[n] \to [n]$ with $E_n:= 012\dots n$, the sequence of all distinct numbers from 0 to $n$.

Because each face map $d_i^m: \Delta^n_m \to \Delta^n_{m-1}$ is the precomposition with the coface map $\delta_i^m: [m-1] \to [m]$, in this notation $d_i$ deletes the number in the $i$-th position:
\begin{equation}\label{eq:face-seq-notation}
    d_i u = u_0 \dots u_{i-1} u_{i+1} \dots u_{m} \in \Delta^n_{m-1}.
\end{equation}
In the same way, each degeneracy map $s_i^m: \Delta^n_m \to \Delta^n_{m+1}$ is the precomposition with the codegeneracy map $\sigma_i^m: [m+1] \to [m]$, so in this notation $s_i$ repeats the number in the $i$-th position:
\begin{equation}\label{eq:deg-splx-seq-notation}
    s_i u = u_0 \dots u_{i-1} u_i u_i u_{i+1} \dots u_{m} \in \Delta^n_{m+1}.
\end{equation}
Therefore, in this notation, a simplex is degenerate if and only if it contains repeated indices. Hence we call $E_n$ the \textbf{unique nondegenerate $n$-simplex in $\Delta^n$}. 

\begin{definition}\label{def:simp-subset-generated}
Let $\huaX$ be a simplicial object in a concrete category $\Cat$, and $x \in \huaX_n$ an $n$-simplex. The \textbf{simplicial subset generated by $x$} is the smallest simplicial subset of $\huaX$ containing $x$. We denote this by $\langle\langle x \rangle\rangle$. Each of its levels can be written as
\begin{equation*}
    \langle\langle x \rangle\rangle_m = \{ \huaX(u)(x) \in \huaX_m \mid u \in \Delta([m],[n])\}.
\end{equation*}
\end{definition}

By Lemma \eqref{lem:Delta-gen-by-cf-cd} any morphism in $\Delta$ is generated by cofaces and codegeneracies as above, so any morphism with target $[n]$ is generated by precomposing the identity with cofaces and codegeneracies. This gives the following classification of simplices in $\Delta$. 

\begin{lemma}\label{lem:Delta-n-simplex-class}
    Every $m$-simplex $u$ in $\Delta^n_m$ can be written as
    \begin{equation*}
        u = s_{\widebar{I}}d_JE_n = s_{i_{|I|}} \circ \dots \circ s_{i_1} \circ d_{j_1} \circ \dots \circ d_{j_{|J|}} (E_n),
    \end{equation*}
    for a unique pair of multi-indices $I$ and $J$ of lengths $|I|$ and $|J|$ with $n - |I| + |J| = m$ and such that 
    \begin{equation*}
        0 \le i_1 < i_2 < \dots < i_{|I|} \le m-1, \qquad
        0 \le j_1 < j_2 < \dots < j_{|J|} \le n.
    \end{equation*}
    Therefore, $\Delta^n$ is the simplicial set generated by $E_n$.
\end{lemma}

\begin{figure}[!ht]
    \centering
    \begin{adjustbox}{width=\textwidth}
    \begin{tikzpicture}[scale=0.6]
       	\begin{pgfonlayer}{nodelayer}
		\node [style=dot, label={below left:\small{0}}] (0) at (0, 0) {};
		\node [style=dot, label={below left:\small{0}}] (1) at (3, 0) {};
		\node [style=dot, label={below right:\small{1}}] (2) at (6, 0) {};
		\node [style=dot, label={below left:\small{0}}] (3) at (9, 0) {};
		\node [style=dot, label={below right:\small{2}}] (4) at (12, 0) {};
		\node [style=dot, label={above right:\small{1}}] (5) at (10.5, 2.5) {};
		\node [style=dot, label={below left:\small{0}}] (6) at (15, 0) {};
		\node [style=dot, label={above left:\small{1}}] (7) at (16.5, 2.5) {};
		\node [style=dot, label={below right:\small{2}}] (8) at (18, 0) {};
		\node [style=dot, label={above right:\small{3}}] (9) at (19.5, 2.5) {};
		\node [style=dot, label={below left:\small{0}}] (10) at (23.5, 0) {};
		\node [style=dot, label={above left:\small{1}}] (11) at (22.5, 2.5) {};
		\node [style=dot, label={above right:\small{2}}] (12) at (25, 4) {};
		\node [style=dot, label={above right:\small{3}}] (13) at (27.5, 2.5) {};
		\node [style=dot, label={below right:\small{4}}] (14) at (26.5, 0) {};
		\node [style=none] (15) at (0, -2) {{\large{$\Delta^0$}}};
		\node [style=none] (16) at (4.5, -2) {{\large{$\Delta^1$}}};
		\node [style=none] (17) at (10.5, -2) {{\large{$\Delta^2$}}};
		\node [style=none] (18) at (17.25, -2) {{\large{$\Delta^3$}}};
		\node [style=none] (19) at (25, -2) {{\large{$\Delta^4$}}};
	\end{pgfonlayer}
	\begin{pgfonlayer}{edgelayer}
		\draw [style=directed] (5) to (3);
		\draw [style=directed] (4) to (5);
		\draw [style=directed] (4) to (3);
		\draw [style=directed] (2) to (1);
		\draw [style=directed] (8) to (6);
		\draw [style=directed] (9) to (8);
		\draw [style=directed] (9) to (7);
		\draw [style=directed] (7) to (6);
		\draw [style=directed-dash] (8) to (7);
		\draw [style=directed] (9) to (6);
		\draw [style=directed] (14) to (10);
		\draw [style=directed-dash] (13) to (10);
		\draw [style=directed] (14) to (12);
		\draw [style=directed] (13) to (12);
		\draw [style=directed-dash] (13) to (11);
		\draw [style=directed] (12) to (11);
		\draw [style=directed-dotted] (12) to (10);
		\draw [style=directed] (11) to (10);
		\draw [style=directed] (14) to (11);
		\draw [style=directed] (14) to (13);
	\end{pgfonlayer}
    \end{tikzpicture}
    \end{adjustbox}
    \caption{The standard simplices in dimension up to 4. We follow the convention of orienting the edges from the higher to the lower vertex index. The dashed lines are just for ease of reading and carry no special meaning in this diagram.}
    \label{fig:standard-simplices}
\end{figure}
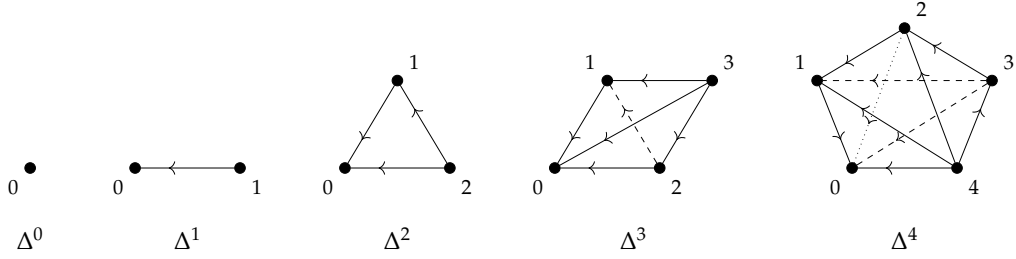

This is reflected in the geometric interpretation of simplicial sets. The standard simplex $\Delta^n$ describes the geometry of the oriented $n$-simplex. See Figure \ref{fig:standard-simplices}.

When writing a simplex in $\Delta^n$ as $u= u_0\dots u_m$ we are listing its vertices. 
If any of these are repeated, the simplex $u$ is degenerate and this represents the fact that a lower dimensional simplex can be seen as a collapsed or ``very thin'' higher dimensional simplex. 
On the other hand, taking $d_iu$ means considering the face of $u$ opposite the vertex $u_i$.

\begin{remark}\label{rem:simplices-as-simplicial-maps}
The simplicial subset generated by an $n$-simplex $x$ can also be described by observing the fundamental fact that any $n$-simplex $x$ in a concrete simplicial object $\huaX$ can be identified with the unique simplicial set map $x:\Delta^n \to \huaX$ defined by the property that $x(E_n) = x \in \huaX_n$. This follows by Lemma \ref{lem:Delta-n-simplex-class}. The image of this map is exactly the simplicial subset $\langle\langle x \rangle\rangle$ of $\huaX$ generated by $x$. 
More precisely, for any $u = s_{\widebar{I}}d_JE_n \in \Delta^n_m$, we have
\begin{equation}
    x(u) = x(s_{\widebar{I}}d_JE_n)
    =  s_{\widebar{I}}d_J x(E_n) = s_{\widebar{I}}d_J x. 
\end{equation}
This observation also shows that there is a bijection
\begin{equation*}
    \SSet(\Delta^n, \Forget(\huaX)) \cong \Forget(\huaX_n).
\end{equation*}
This can also be inferred from the Yoneda lemma, as it is the same as writing
\begin{equation*}
    \Fun(\Delta^{op}, \Set )(\Delta(\_, [n]), \Forget(\huaX)) \cong \Forget(\huaX([n])).
\end{equation*}
\end{remark}

In the same way the face and degeneracy maps of $\Delta^n$ are defined by precomposition with the cofaces and codegeneracies, \textit{postcompositions} with coface and codegeneracy maps give certain simplicial morphisms between the different standard simplices, which obey the cosimplicial identities.\footnote{Indeed, the hom functor $\Delta: \Delta^{op} \times \Delta \to \Set$ is a cosimplicial-simplicial set.}
Postcomposition with the coface maps gives the \textbf{face inclusions} $\delta_i: \Delta^{n-1} \to \Delta^n$, while postcomposition with the codegeneracy maps gives the \textbf{collapse} maps $\sigma_i: \Delta^{n+1} \to \Delta^n$.
In sequence notation, these maps act as \eqref{eq:Delta-coface-def} and \eqref{eq:Delta-codeg-def} on each number of the sequence
\begin{equation*}
    \delta_i u = (\delta_i u_1) \dots (\delta_i u_m) \in \Delta^{n-1}_m,
    \qquad
    \sigma_i u = (\sigma_i u_1) \dots (\sigma_i u_m) \in \Delta^{n+1}_m.
\end{equation*}

There are three more implicit facts related to the proof of Lemma \ref{lem:Delta-n-simplex-class} and the identification of $E_n$ with $id_{[n]}$:
\begin{enumerate}
\item The fact that both the face inclusions and the collapses are simplicial maps, which follows by associativity of composition inside $\Delta$: By Lemma \ref{lem:Delta-n-simplex-class}, it is enough to check this on the unique nondegenerate simplex which represents the identity. For example,
\begin{equation}\label{eq:Delta-simp-cosimp-commutativity}
    d_i (\delta_j E_{n}) = (\delta_j \circ id_{[n]}) \circ \delta_i = \delta_j \circ (id_{[n]} \circ \delta_i) = \delta_j (d_i E_n). 
\end{equation}
The other compatibility conditions follow in the same way. 

\item The two useful identities
\begin{equation}\label{eq:d-delta-En-identity}
\begin{split}
    d_i E_n &= id_{[n]} \circ \delta_i = \delta_i \circ id_{[n-1]} = \delta_i E_{n-1}, \text{ and }\\
    s_i E_n &= id_{[n]} \circ \sigma_i = \sigma_i \circ id_{[n+1]} = \delta_i E_{n+1}.
\end{split}
\end{equation}

\item And finally, the fact that by these identities any simplex $u$ in $\Delta^n_m$ can be written as 
\begin{equation}\label{eq:simp-cosimp-description-u}
    u = s_{\widebar{I}}d_J E_n = \delta_{\widebar{J}} \sigma_I E_m,
\end{equation}
for the same unique multi-indices $I$ and $J$. 
\end{enumerate}

We conclude this discussion of the standard simplices by observing that increasing multi-indices such as the $I$ and $J$ of Lemmas \ref{lem:Delta-gen-by-cf-cd} and \ref{lem:Delta-n-simplex-class} are in 1:1 correspondence with shuffle permutations. 

\begin{definition}[{\cite[Def. I.6.5]{May1967}}]\label{def:shuffles}
    For any integers $p,q \ge 0$, a \textbf{$(p, q)$-shuffle} is a permutation $\pi \in \mathbb{S}_{p+q}$ such that 
    \begin{equation*}
        \pi(i) < \pi(j), \quad \text{if } 0 \le i < j \le p-1 \text{ or if } p \le i < j \le p + q - 1. 
    \end{equation*}
    By writing $\mu_i = \pi(i - 1)$ for any $1 \le i \le p$ and $\nu_i = \pi(p + i - 1)$ for any $1 \le i \le q$, each shuffle $\pi$ is identified with a pair of increasing multi-indices $(\mu, \nu)$, where one item of the pair determines the other. That is 
    \begin{equation*}
        \pi = (\mu,\nu) = (\mu_1\dots\mu_p,\nu_1\dots\nu_q) = (\pi(0) \dots \pi(p-1), \pi(p) \dots \pi(p+q -1)).
    \end{equation*}
    The sign of a shuffle $\sign(\mu,\nu)$ is defined as the sign of the corresponding permutation $\pi$. 

    The set of $(p,q)$-shuffles is denoted by $\Shuf(p,q)$, and it has cardinality $\binom{p+q}{p}$. 
    As a convention, we define $\Shuf(p,q) := \emptyset$ for any negative integer $p$ or $q$. 
\end{definition}

\begin{example}
    Some sets of shuffles are
    \begin{equation*}
        \begin{split}
            & \forall n\ge 1, \, \Shuf(n,0) = \{id=(\mu=0 \dots n)\} \cong \Shuf(0,n) = \{id=(\nu=0 \dots n)\};\\
            &\Shuf(1,1) = \{id=(0,1), (1,0)\}, \quad \Shuf(2,1) = \{id=(01,2), (02,1), (12,0)\},\\ 
            &\Shuf(1,2) = \{id=(0,12), (1,02), (2,01)\},\\
            &\Shuf(3,1) = \{id=(012,3), (013, 2), (023, 1), (123, 0)\}, \dots, \\
            &\Shuf(2,2) = \{id=(01,23), (02, 13), (03, 12), (12, 03), (13, 02), (23, 01)\}, \dots
        \end{split}
    \end{equation*}
\end{example}

\begin{proposition}\label{prop:shuffles-standard-simpl-decomp}
The set of $m$-simplices of $\Delta^n$ can be decomposed as:
\begin{align}
    \text{if } m=0, \quad \Delta^n_0 &\cong pr_1\Shuf(1 , n), \notag\\
    \text{if } m \ge 1, \quad \Delta^n_m &\cong pr_1\Shuf(m+1, n-m) \label{eq:decomp-m-level-n-splx}\\
    &\qquad \sqcup \coprod_{i=1}^m \left(\coprod_{(\mu,\nu) \in \Shuf(m-i, i)} s_{\widebar{\nu}}(pr_1\Shuf(m-i+1, n-m+i))\right), \notag
\end{align}
where $pr_1$ is the projection to the first component $\mu$ of a shuffle $(\mu,\nu)$.
The first term in \eqref{eq:decomp-m-level-n-splx} is in bijection with the set of nondegenerate simplices, while the remainder is the set of degenerate simplices. 

Therefore, the cardinality of $\Delta^n_m$ is 
\begin{equation}\label{eq:cardinality-m-simplices-in-Delta-n}
    |\Delta^n_m| = \binom{n+m+1}{m+1} = \binom{n+m+1}{n}.
\end{equation} 
\end{proposition}

\begin{proof}
By Lemma \ref{lem:Delta-n-simplex-class} and the preceding discussion, the nondegenerate $m$-simplices in $\Delta^n$ are identified with strictly increasing multi-indices of length $m+1$ with values in $[n]$. 
Let $\mu \in \Delta^n_m$ be the nondegenerate $m$-simplex $\mu_1\dots \mu_m$. 
Then such a multi-index identifies an $(m+1, n-m)$-shuffle $(\mu, \nu)$, and 
\begin{equation*}
    \mu = d_\nu E_n.
\end{equation*}
Indeed, $E_n$ itself is just the $\mu$ for the unique $(n,0)$-shuffle. 
In particular, there are $\binom{n+1}{m+1}$ nondegenerate $m$-simplices in $\Delta^n$.

By the same Lemma, each nondegenerate simplex generates a number of degenerate simplices classified by another increasing multi-index, which also identifies a shuffle. 
More precisely, any $l$-simplex $\mu\in \Delta^n_l$, generates the degenerate $k$-simplices $s_{\widebar{\nu'}}\mu$ for any $k>l$, and $\nu'$ such that $(\mu',\nu')\in\Shuf(l, k-l)$ for the appropriate $\mu'$ determined as before.

Therefore, the set of $m$-simplices of $\Delta^n$ can be decomposed as in \eqref{eq:decomp-m-level-n-splx}.

With this decomposition we can compute the number of total $m$-simplices in $\Delta^n$ by using the Chu-Vandermonde identity for binomial coefficients 
\begin{equation}\label{eq:chu-vandermonde}
    \binom{p+q}{r} = \sum^r_{i=0}\binom{p}{i} \binom{q}{r-i}, \qquad \forall p,q,r\in \N. 
\end{equation}
In fact by \eqref{eq:decomp-m-level-n-splx}, we have
\begin{equation*}
\begin{split}
    |\Delta^n_m| &= \binom{n+1}{m+1} + \sum_{i=1}^{m}\binom{m}{i}\binom{n+1}{m+1-i}\\
    &= \sum_{i=0}^{m}\binom{m}{i}\binom{n+1}{m+1-i} + \binom{m}{m+1}\binom{n+1}{0},
\end{split}
\end{equation*}
where the last term is $0$. Then \eqref{eq:cardinality-m-simplices-in-Delta-n} follows from \eqref{eq:chu-vandermonde} with $p = m$, $q=n+1$ and $r = m+1$. 
\end{proof}

\begin{remark}\label{rem:total-unit}
    Because strictly increasing multi-indices such as those in Remark \ref{rem:canonical-form-s-d-comp} identify shuffles, this can be used to classify all possible compositions of face and degeneracy maps between any two levels of a simplicial object $\huaX$ by using the normal form \eqref{eq:sd-comp-normal-form}. 
    In particular, for any $m\ge 1$ there is a unique such map $1: \huaX_0 \to \huaX_m$, which is
    \begin{equation}\label{eq:total-unit}
        1 = s_{m-1} \circ \dots \circ s_1 \circ s_0,
    \end{equation}
    and it coincides with any other possible composition of degeneracy maps from $\huaX_0$ to $\huaX_m$. 
    We call this map the \textbf{total unit}. Its associated multi-index corresponds to the unique $(m,0)$-shuffle given by the identity of $\{0, \dots, m-1\}$. 
\end{remark}

\subsection{Products and mapping sets}\label{sec:SSet-products-and-mapping-spaces}

We now describe products and mapping spaces of simplicial sets. 
The material in this section is fundamental in simplicial homotopy theory and can be found for example in \cite{May1967,GoerssJardine2009, Friedman2012, Riehl2014, kerodon,stacks-project}. 

To begin with, we define the product of two simplicial objects in any category with binary products. When applied to simplicial sets, this gives a symmetric monoidal structure on $\SSet$.

\begin{definition}
    Let $\huaX$ and $\huaY$ be simplicial objects in $\Cat$. The \textbf{product} $\huaX \times \huaY$ is the simplicial set defined at any level $n \ge 0$ by 
    \begin{equation*}
        (\huaX \times \huaY)_n := \huaX_n \times \huaY_n,
    \end{equation*}
    with simplicial maps
    \begin{equation*}
        d_i(x,y) = (d_i x, d_i y), \qquad 
        s_i(x,y) = (s_i x, s_i y),
    \end{equation*}
    for all $(x,y) \in (\huaX \times \huaY)_n$ and any appropriate $i$.
\end{definition}

\begin{definition}\label{def:copowering}
    Let $\Cat$ be a category with all (finite) coproducts of an object with itself. Let $\huaS$ be a (levelwise finite) simplicial set and $\huaX$ a simplicial object in $\Cat$. The \textbf{copowering} $\huaX \otimes \huaS$ is the simplicial set defined at each level $n \ge 0$ by 
    \begin{equation}\label{eq:copowering-1}
    (\huaX \otimes \huaS)_n := \coprod_{u\in\huaS_n} \huaX_n^u,
    \end{equation}
    with simplicial maps defined by
    \begin{equation*}
        d_i(x^u) = (d_i x)^{d_i u}, \qquad s_i(x^u) = (s_i x)^{s_i u},
    \end{equation*}
    for all $x^u\in (\huaX \otimes \huaS)_n$ and any appropriate $i$.
\end{definition}

In the situation of the above definition we say that $\Simp\Cat$ is \textbf{copowered} over (levelwise finite) simplicial sets. The categories of levelwise finite-dimensional $\Simp\Vect$ and $\Simp\Mfd$ are copowered over levelwise finite simplicial sets.

\begin{remark}
    It is clear that if $\Cat = \Set$, the product and the copowering of two simplicial sets with each other coincide. 
\end{remark}

The following lemma is an immediate consequence of the definitions above and it can be proven by constructing the obvious natural isomorphisms. It is part of \cite[Thm. II.2.5]{GoerssJardine2009}.

\begin{lemma}[Associativity of the copowering]\label{lem:copow-assoc}
    Let $\huaX\in \Simp\Cat$, and $\huaS, \huaT$ be levelwise finite simplicial sets. Then 
    \begin{equation*}
        (\huaX \otimes \huaS) \otimes \huaT \cong \huaX \otimes (\huaS \times \huaT). 
    \end{equation*}
\end{lemma}

\begin{definition}\label{def:simp-homotopy-of-maps}
    Let $\huaX$, $\huaY$ be simplicial objects in $\Cat$ and $f,g:\huaX \to \huaY$ two simplicial maps between them. A \textbf{simplicial homotopy} $h$ between $f$ and $g$ is a simplicial map
    \begin{equation*}
        h: \huaX \otimes \Delta^1 \to \huaY
    \end{equation*}
    such that 
    \begin{equation*}
    \begin{split}
        &d_0 h = h \circ (id \otimes \delta_0) = f: \huaX \otimes \Delta^0 \cong \huaX \to \huaY, \text{ and }\\
        &d_1 h = h \circ (id \otimes \delta_1) = g: \huaX \otimes \Delta^0 \cong \huaX \to \huaY.
    \end{split}
    \end{equation*}
\end{definition}

A similar definition can be used to recursively define higher homotopies between homotopies. 

\begin{definition}
    Fix $n \ge 1$ and let $\huaX$, $\huaY$ be simplicial objects in $\Cat$ and $h,k:\huaX \otimes \Delta^{n-1} \to \huaY$ two simplicial $(n-1)$-homotopies. A \textbf{simplicial $n$-homotopy} $H$ between $h$ and $k$ is a simplicial map
    \begin{equation*}
        H: \huaX \otimes \Delta^n \to \huaY
    \end{equation*}
    such that 
    \begin{equation*}
    \begin{split}
        &d_0 H = H \circ (id \otimes \delta_0) = h: \huaX \otimes \Delta^{n-1} \to \huaY, \text{ and }\\
        &d_1 H = H \circ (id \otimes \delta_1) = k: \huaX \otimes \Delta^{n-1} \to \huaY.
    \end{split}
    \end{equation*}
\end{definition}

The simplicial maps between two simplicial objects are the 0-simplices of the following simplicial set, consisting of increasingly higher homotopies at each level.

\begin{definition}\label{def:mapping-set}
    The \textbf{mapping set} $\IHom(\huaX, \huaY)$ between two simplicial objects $\huaX$ and $\huaY$ in $\Cat$ is the simplicial set defined at each level $n \ge 0$ by 
    \begin{equation*}
        \IHom(\huaX, \huaY)_n := \Simp\Cat(\huaX \otimes \Delta^n, \huaY),
    \end{equation*}
    with simplicial maps defined in components by
    \begin{equation}\label{eq:mapping-set-maps}
    \begin{split}
        (d_if)_l (x^u) &:= f_l (x^{\delta^i(u)}), \quad \forall x \in \huaX_{l}, u \in \Delta^{m-1}_{l},\\
        (s_if)_l (x^u) &:= f_l (x^{\sigma^i(u)}), \quad \forall x \in \huaX_{l}, u \in \Delta^{m+1}_{l}.
   \end{split}
    \end{equation}
    This is sometimes called the \textbf{function complex}. When this is also an element in $\Simp\Cat$, we call it the \textbf{mapping space}.
\end{definition}

Simplicial $n$-homotopies can be composed across different mapping sets according to the composition defined at each level $n$ by
\begin{equation}\label{eq:composition-of-homotopies}
    \begin{split}
        \cdot \circ \cdot : &\IHom(\huaV, \huaW)_n \otimes \IHom(\huaU, \huaV)_n \longrightarrow \IHom(\huaU, \huaW)_n\\
        &ik \circ h: \huaU \otimes \Delta^n \xrightarrow{id \times diag} \huaU \otimes (\Delta^n \times \Delta^n) \cong (\huaU \otimes \Delta^n) \otimes \Delta^n\\
        &\qquad \xrightarrow{h \otimes id} \huaV \otimes \Delta^n \xrightarrow{k} \huaW,
    \end{split}
\end{equation}
where $diag$ is the diagonal inclusion, and we used the associativity property from Lemma \ref{lem:copow-assoc} in the second step. 
If $h$ is a homotopy between $f$ and $f'$ and $k$ is a homotopy between $g$ and $g'$, then their composition $h \circ k$ is a homotopy between $g\circ f$ and $g'\circ f'$.
It is straightforward to check that with this composition and the obvious identity element, $\Simp\Cat$ is an $\SSet$-category, in the sense of \cite[Def. 3.3.1]{Riehl2014}. Categories enriched over simplicial sets are also known as \textbf{simplicial categories}, see e.g. \cite[\S 3.6]{Riehl2014}, \cite[\href{https://kerodon.net/tag/00JQ}{Tag 00JQ}]{kerodon}.

If $\Cat$ is the category of sets, we denote the mapping set by $\Maps$. This is precisely the internal hom of simplicial sets with respect to the monoidal structure given by the product.
We report this adjunction in the following proposition. 
This makes $\SSet$ a \textbf{cartesian closed category}: it is closed monoidal with respect to the cartesian product.

\begin{proposition}[Enriched tensor-hom adjunction]\label{prop:enriched-tensor-hom-sset}
    Let $\huaX$, $\huaY$ and $\huaZ$ be simplicial sets. Then there are natural isomorphisms
    \begin{equation*}
        \Maps(\huaX \times \huaY, \huaZ) \cong \Maps(\huaX, \Maps(\huaY, \huaZ)) \cong \Maps(\huaY, \Maps(\huaX, \huaZ)).
    \end{equation*}
\end{proposition}

\begin{proof}
    This proof is an extension of \cite[Prop. I.5.1]{GoerssJardine2009}. We show a version of this result explicitly for simplicial vector spaces in Proposition \ref{prop:tensor-hom-svect}, since we will use it in Chapter \ref{chap:vs-duals}. 
\end{proof}

\begin{remark}\label{rem:lit-rev-simp-cat}
As mentioned above, we are using the term ``simplicial category'' to indicate a simplicially enriched category or an $\SSet$ category, following \cite{Riehl2014,kerodon}. 
There are other notions of simplicial category in the literature, such as the one in \cite[Def. II.2.1]{GoerssJardine2009}. In this definition, for $\Simp\Cat$ to be a simplicial category, in addition to the simplicial enrichment, they require that: 
\begin{enumerate}
\item The mapping set functor has a left adjoint, called the \textbf{tensoring} or \textbf{copowering} functor
\begin{equation*}
    \_ \otimes \_ : \Simp\Cat \times \SSet \to \Simp\Cat
\end{equation*}
with the associativity property in Lemma \ref{lem:copow-assoc},
\item The mapping set functor has a right adjoint, called the \textbf{cotensoring} or \textbf{powering} functor
\begin{equation*}
    \powering(\_, \_): \SSet^{op} \times \Simp\Cat \to \Simp\Cat.
\end{equation*}
In other words, these three functors must fit into a two-variable adjunction 
\begin{equation}\label{eq:pow-copow-adjunction}
    \Simp\Cat(\huaX \otimes \huaS, \huaY)
    \cong \SSet(\huaS, \IHom(\huaX, \huaY))
    \cong \Simp\Cat(\huaY, \powering(\huaS, \huaY)),
\end{equation}
for any $\huaX, \huaY \in \Simp\Cat$ and $\huaS \in \SSet$. 
\end{enumerate}

A simplicially enriched category where (1) holds is said to be \textbf{copowered over simplicial sets}. In \cite[Def. 3.7.2]{Riehl2014} the term \textbf{tensored} is used instead. 
On the other hand, a simplicially enriched category where (2) holds is said to be \textbf{powered over simplicial sets}, or analogously \textbf{cotensored}, as in \cite[Def. 3.7.3]{Riehl2014}. See also \cite[\S 14]{Shulman2006}.
\end{remark}

\begin{remark}
The category of simplicial sets satisfies all the properties in Remark \ref{rem:lit-rev-simp-cat}, and the powering functor coincides with the mapping set functor. In fact, Proposition \ref{prop:enriched-tensor-hom-sset} is a $\SSet$-enriched version of adjunction \eqref{eq:pow-copow-adjunction}, which holds for all enriched, copowered and powered categories over simplicial sets and is reported without proof in \cite[Lemma II.2.3]{GoerssJardine2009} and \cite[Rem. 3.7.4]{Riehl2014}. As we will see in Section \ref{sec:Svect-int-hom}, the category of simplicial vector spaces is enriched, copowered and powered over $\SSet$ but also over itself, and an analogous $\SVect$-enriched adjunction holds. 
\end{remark}

\begin{example}\label{ex:pow-copow-over-sets}
    Any category with self products and coproducts is powered and copowered over sets, respectively. We use this implicitly in Definition \ref{def:copowering}, where \eqref{eq:copowering-1} describes level $n$ of $\huaX \otimes \huaS$ as the copowering over sets of the level $n$ of $\huaX$ over the respective level of $\huaS$. 
\end{example}

\section{Higher Groupoids}\label{sec:higher-gpds}

The definition of higher groupoids in terms of simplicial sets or simplicial objects in a more general topos was first given in \cite{Duskin1979,Glenn1982}. This was then adapted to higher Lie groupoids in the work of \cite{Getzler2009, Henriques2008,Zhu2009}. 

We begin by discussing the definitions in the category of simplicial sets. All of these constructions generalize without any problem to higher groupoid objects in $\Vect$, as we recall in Section \ref{sec:simp-vect}. For higher Lie groupoids more care is needed and this is the bulk of the discussion in the later subsections.

\subsection{Simplicial diagrams and groupoids}\label{sec:simp-diag-sets}

The correspondence between $n$-simplices in a simplicial object $\huaX$ and simplicial maps $\Delta^n \to \huaX$ extends to simplicial diagrams in $\huaX$ of a certain shape given by a simplicial set $\huaS$ and simplicial maps $\huaS \to \huaX$. 
This is analogous to the definition of a diagram in a category $\Cat$ as a functor $\huaJ \to \Cat$ with $\huaJ$ a small indexing category representing the shape of the diagram.

\begin{definition}\label{def:simpl-diags}
    Let $\huaS$ be a levelwise finite simplicial set, and $\huaX$ a simplicial set. A \textbf{simplicial diagram} in $\huaX$ of shape $\huaS$ is a simplicial map $\huaS \to \huaX$. 
    The \textbf{set of simplicial diagrams} of shape $\huaS$ in $\huaX$ is $\SSet(\huaS, \huaX)$ and we denote it by $\hom(\huaS, \huaX)$.
\end{definition}

\begin{remark}\label{rem:diagram-space-and-powering}
The notation $\hom(\huaS, \huaX)$ comes from the fact that this is level 0 of the powering of simplicial sets with itself, which is defined by $\powering(\huaS, \huaX)_m = \hom(\huaS \otimes \Delta^m, \huaX)$. 
That is, higher levels of the powering functor are diagrams in the shape of prisms with base $\huaS$ inside the second argument.
This works for simplicial vector spaces too: as we see in Section \ref{sec:simp-vect}, the space of simplicial diagrams in a simplicial vector space can be used to define the powering of simplicial vector spaces over simplicial sets. 
\end{remark}

The most important shapes of simplicial diagrams we will consider are horns and boundaries. 

\begin{definition}
    For any integer $m \ge 0$ the \textbf{$m$-boundary} $\partial\Delta^m$ is the simplicial subset of the $m$-simplex $\Delta^m$ obtained by removing the $m$-simplex $E_m$ and all of its degenerate simplices. 
    In other words, at each level $l\ge 0$, 
    \begin{equation*}
        (\partial\Delta^m)_l =  \{ f\in (\Delta^m)_l \mid \Img f \nsupseteq \{0,\dots, m\} \} \subseteq \Delta^m_l.
    \end{equation*}
    In particular $\partial\Delta^0 = \emptyset$.

    For any integers $m\ge 0$ and $0 \le k \le m$, the \textbf{$(m,k)$-horn $\Lambda^m_k$} is the simplicial subset of the $m$-boundary $\partial\Delta^m$ obtained by removing the $k$-th face $d_kE_n$ and all its degenerate simplices. 
    In other words, at each level $l\ge 0$, 
    \begin{equation*}
        (\Lambda^m_k)_l = \{ f \in (\Delta^m)_l \mid \Img f \nsupseteq \{0,\dots,\widehat{k},\dots,m\}
        \} \subseteq (\partial\Delta^m)_l \subseteq \Delta^m_l.
    \end{equation*}
    Again, clearly $\Lambda^0_0 = \emptyset$. 
\end{definition}

The sets of $m$-boundaries and of $(m,k)$-horns in a simplicial set are denoted, respectively, by 
\begin{equation*}
    \begin{array}{cc}
    \partial^m(\huaX) := \hom(\partial\Delta^m, \huaX),
    &\Lambda^m_k(\huaX) := \hom(\Lambda^m_k,\huaX).
    \end{array}
\end{equation*}
Both of these sets can be written as fiber products (i.e. limits) of $\huaX_{m-1}$ over lower dimensional horns and boundaries,\footnote{For a more precise statement on how to write any object of simplicial diagrams as a limit see Proposition \ref{prop:simp-diag-sheaf-as-limit} in the next section.} as their elements can be written as tuples of $(m-1)$-simplices in $\huaX$ satisfying certain conditions:
\begin{equation}\label{eq:DefHornBdrySpaces}
    \begin{split}
    \partial^m(\huaX) &= \{(x_0, \dots, x_m) \in (\huaX_{m-1})^{\times m+1} \mid d_ix_j = d_{j-1}x_i, \text{ for all } i < j\},\\
    \Lambda^m_k(\huaX) &= \{(x_0, \dots, \widehat{x_k}, \dots, x_m) \in (\huaX_{m-1})^{\times m} \mid \\
    &\qquad \qquad d_ix_j = d_{j-1}x_i, \text{ for all } i < j, \text{ such that } k\neq i,j\}.
    \end{split}
\end{equation}
Here the notation $(x_0, \dots, \widehat{x_k}, \dots, x_m)$ denotes an $m$-tuple with components $x_i$ for $0\le i \le m$ with $i \neq k$.

\begin{figure}[!ht]
\centering
\begin{adjustbox}{width=\textwidth}
\pgfdeclarelayer{nodelayer}
\pgfdeclarelayer{edgelayer}
\pgfsetlayers{nodelayer,main,edgelayer}
\begin{tikzpicture}[scale=0.6,
    none/.style={}, 
    dot/.style={fill=black, draw=black, shape=circle, scale=0.40}, 
    fill-grey/.style={-, fill={rgb,255: red,220; green,220; blue,220}, draw=none, fill opacity=0.5},
    directed/.style={decoration={markings,mark=at position 0.7 with {\arrow{To[scale=1.25]}}}, postaction={decorate}},
    dashdir/.style={dashed, decoration={markings,mark=at position 0.7 with {\arrow{To[scale=1.25]}}}, postaction={decorate}},
    ]
    \begin{pgfonlayer}{nodelayer}
        \node [style=none] (25) at (20.5, 1) {};
        \node [style=none] (26) at (19.5, 3) {};
        \node [style=none] (27) at (21.5, 3) {};
        \node [style=none] (0) at (0, -1) {$\Lambda^1_0$};
        \node [style=none] (1) at (2.5, -1) {$\Lambda^1_1$};
        \node [style=none] (2) at (6, -1) {$\Lambda^2_0$};
        \node [style=none] (3) at (10.5, -1) {$\Lambda^2_1$};
        \node [style=none] (4) at (15, -1) {$\Lambda^2_2$};
        \node [style=none] (5) at (19.5, -1) {$\Lambda^3_0$};
        \node [style=none] (6) at (24, -1) {$\Lambda^3_1$};
        \node [style=dot, label={below:\small{0}}] (7) at (0, 1) {};
        \node [style=dot, label={below:\small{1}}] (8) at (2.5, 1) {};
        \node [style=dot, label={below:\small{0}}] (9) at (5, 1) {};
        \node [style=dot, label={below:\small{1}}] (10) at (7, 1) {};
        \node [style=dot, label={above:\small{2}}] (11) at (6, 3) {};
        \node [style=dot, label={below:\small{0}}] (12) at (9.5, 1) {};
        \node [style=dot, label={below:\small{1}}] (13) at (11.5, 1) {};
        \node [style=dot, label={above:\small{2}}] (14) at (10.5, 3) {};
        \node [style=dot, label={below:\small{0}}] (15) at (14, 1) {};
        \node [style=dot, label={below:\small{1}}] (16) at (16, 1) {};
        \node [style=dot, label={above:\small{2}}] (17) at (15, 3) {};
        \node [style=dot, label={below:\small{0}}] (18) at (18.5, 1) {};
        \node [style=dot, label={below:\small{1}}] (19) at (20.5, 1) {};
        \node [style=dot, label={above:\small{2}}] (20) at (19.5, 3) {};
        \node [style=dot, label={above:\small{3}}] (21) at (21.5, 3) {};
        \node [style=none] (22) at (23, 1) {};
        \node [style=none] (23) at (24, 3) {};
        \node [style=none] (24) at (26, 3) {};
        \node [style=dot, label={below:\small{0}}] (28) at (23, 1) {};
        \node [style=dot, label={below:\small{1}}] (29) at (25, 1) {};
        \node [style=dot, label={above:\small{2}}] (30) at (24, 3) {};
        \node [style=dot, label={above:\small{3}}] (31) at (26, 3) {};
    \end{pgfonlayer}
    \begin{pgfonlayer}{edgelayer}
        \draw [style=dotted-fill] (27.center)
                to (25.center)
                to (26.center)
                to cycle;
        \draw [style=directed] (10) to (9);
        \draw [style=directed] (11) to (9);
        \draw [style=directed] (13) to (12);
        \draw [style=directed] (14) to (13);
        \draw [style=directed] (17) to (15);
        \draw [style=directed] (17) to (16);
        \draw [style=directed] (20) to (18);
        \draw [style=directed] (19) to (18);
        \draw [style=directed] (21) to (20);
        \draw [style=directed] (21) to (19);
        \draw [style=dashdir] (20) to (19);
        \draw [style=dotted-fill] (23.center)
                to (24.center)
                to (22.center)
                to cycle;
        \draw [style=directed] (30) to (28);
        \draw [style=directed] (29) to (28);
        \draw [style=directed] (31) to (30);
        \draw [style=directed] (31) to (29);
        \draw [style=dashdir] (30) to (29);
        \draw [style=directed] (31) to (28);
        \draw [style=directed] (21) to (18);
    \end{pgfonlayer}
\end{tikzpicture}
\end{adjustbox}
\caption{A few low-dimensional horns. The dotted faces denote boundaries of triangles without interior.}
\label{fig:horns}
\end{figure}
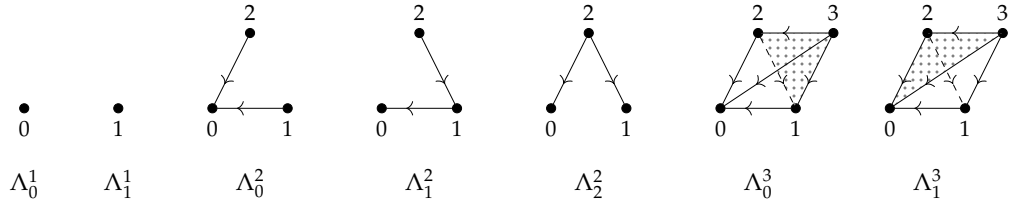

More intuitively, $\partial^m(\huaX)$ is the space of all possible boundaries of $m$-simplices in $\huaX$. In the same way, $\Lambda^m_k(\huaX)$ is the space of all possible configurations of $(m-1)$-simplices in the shape of an $(m,k)$-horn that exist in $\huaX$. This is pictured in Figure \ref{fig:horns}.

A central concept in simplicial methods and homotopy theory is that of a lifting problem. (See e.g. \cite{GoerssJardine2009}). In one of its forms, this is the question of whether for a map of simplicial sets $i: \huaS \to \huaT$ and a simplicial map $f:\huaX \to \huaY$, any diagram of the form 
\begin{equation}\label{diag:lifting-problem}
\begin{tikzcd}[ampersand replacement=\&]
	\huaS \& \huaX \\
	\huaT \& \huaY
	\arrow[from=1-1, to=1-2]
	\arrow["i"', from=1-1, to=2-1]
	\arrow["f", from=1-2, to=2-2]
	\arrow[from=2-1, to=2-2]
\end{tikzcd}
\end{equation}
admits a lift, i.e. a map $\huaT \to \huaX$ such that
\[\begin{tikzcd}[ampersand replacement=\&]
	\huaS \& \huaX \\
	\huaT \& \huaY
	\arrow[from=1-1, to=1-2]
	\arrow["i"', from=1-1, to=2-1]
	\arrow["f", from=1-2, to=2-2]
	\arrow[dashed, from=2-1, to=1-2]
	\arrow[from=2-1, to=2-2]
\end{tikzcd}\]
commutes.

To understand this, we start by considering the case where $\huaY$ is $\Delta^0 = pt$, the terminal object of $\Simp\Set$ and $f:X \to pt$ is the canonical map sending everything to the point. We will discuss the general case in Section \ref{sec:weak-eq}.

When $\huaY$ is the point, the existence of a lift --- or more properly an extension, if $i$ is an inclusion --- of any map $\huaS \to \huaX$ to a map $\huaT \to \huaX$ is equivalent to the fact that the precomposition
\begin{equation*}
    i^*: \hom(\huaT, \huaX) \to \hom(\huaS, \huaX)
\end{equation*}
is surjective. 
Moreover, the uniqueness of this lift is equivalent to $i^*$ being injective.
In terms of simplicial diagrams this can be translated as
\begin{center}
    \textit{``Any simplicial diagram of shape $\huaS$ in $\huaX$ determines a (unique)\\ simplicial diagram of shape $\huaT$ in $\huaX$.''}
\end{center}

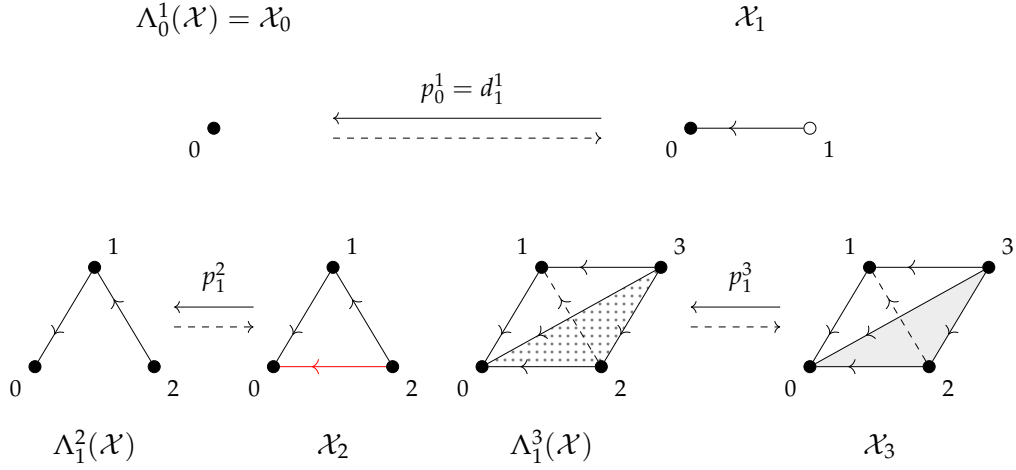
\begin{figure}[!ht]
    \centering
    \begin{adjustbox}{width=\textwidth}
    \begin{tikzpicture}[scale=0.6]
	\begin{pgfonlayer}{nodelayer}
		\node [style=none] (40) at (23, 0) {};
		\node [style=none] (41) at (27.5, 2.5) {};
		\node [style=none] (42) at (26, 0) {};
		\node [style=none] (37) at (14.75, 0) {};
		\node [style=none] (38) at (19.25, 2.5) {};
		\node [style=none] (39) at (17.75, 0) {};
		\node [style=dot, label={below left:\small{0}}] (0) at (8, 6) {};
		\node [style=dot, label={below left:\small{0}}] (1) at (3.5, 0) {};
		\node [style=dot, label={above right:\small{1}}] (2) at (5, 2.5) {};
		\node [style=dot, label={below left:\small{0}}] (3) at (9.5, 0) {};
		\node [style=dot, label={below right:\small{2}}] (4) at (12.5, 0) {};
		\node [style=dot, label={above right:\small{1}}] (5) at (11, 2.5) {};
		\node [style=dot, label={below left:\small{0}}] (6) at (14.75, 0) {};
		\node [style=dot, label={above left:\small{1}}] (7) at (16.25, 2.5) {};
		\node [style=dot, label={below right:\small{2}}] (8) at (17.75, 0) {};
		\node [style=dot, label={above right:\small{3}}] (9) at (19.25, 2.5) {};
		\node [style=none] (15) at (21.5, 8.75) {{\large{$\huaX_1$}}};
		\node [style=none] (16) at (5, -2) {{\large{$\Lambda^2_1(\huaX)$}}};
		\node [style=none] (17) at (11, -2) {{\large{$\huaX_2$}}};
		\node [style=none] (18) at (16.5, -2) {{\large{$\Lambda^3_1(\huaX)$}}};
		\node [style=none] (19) at (24.75, -2) {{\large{$\huaX_3$}}};
		\node [style=dot, label={below right:\small{2}}] (20) at (6.5, 0) {};
		\node [style=none] (21) at (8, 2.25) {{$p^2_1$}};
		\node [style=none] (22) at (9, 1.5) {};
		\node [style=none] (23) at (7, 1.5) {};
		\node [style=none] (24) at (7, 1) {};
		\node [style=none] (25) at (9, 1) {};
		\node [style=dot, label={below left:\small{0}}] (26) at (20, 6) {};
		\node [style=white-dot, label={below right:\small{1}}] (27) at (23, 6) {};
		\node [style=none] (28) at (14.25, 7) {{$p^1_0 = d^1_1$}};
		\node [style=none] (29) at (17.75, 6.25) {};
		\node [style=none] (30) at (11, 6.25) {};
		\node [style=none] (31) at (11, 5.75) {};
		\node [style=none] (32) at (17.75, 5.75) {};
		\node [style=dot, label={below left:\small{0}}] (33) at (23, 0) {};
		\node [style=dot, label={above left:\small{1}}] (34) at (24.5, 2.5) {};
		\node [style=dot, label={below right:\small{2}}] (35) at (26, 0) {};
		\node [style=dot, label={above right:\small{3}}] (36) at (27.5, 2.5) {};
		\node [style=none] (43) at (21.25, 2.25) {{$p^3_1$}};
		\node [style=none] (44) at (22.25, 1.5) {};
		\node [style=none] (45) at (20, 1.5) {};
		\node [style=none] (46) at (20, 1) {};
		\node [style=none] (47) at (22.25, 1) {};
		\node [style=none] (48) at (8, 8.75) {{\large{$\Lambda^1_0(\huaX)= \huaX_0$}}};
	\end{pgfonlayer}
	\begin{pgfonlayer}{edgelayer}
		\draw [style=light-gray-05-fill] (41.center)
			 to (40.center)
			 to (42.center)
			 to cycle;
		\draw [style=dotted-fill] (38.center)
			 to (37.center)
			 to (39.center)
			 to cycle;
		\draw [style=directed] (5) to (3);
		\draw [style=directed] (4) to (5);
		\draw [style=directed-red] (4) to (3);
		\draw [style=directed] (2) to (1);
		\draw [style=directed] (8) to (6);
		\draw [style=directed] (9) to (8);
		\draw [style=directed] (9) to (7);
		\draw [style=directed] (7) to (6);
		\draw [style=directed-dash] (8) to (7);
		\draw [style=directed] (9) to (6);
		\draw [style=directed] (20) to (2);
		\draw [style=to] (22.center) to (23.center);
		\draw [style=dash-to] (24.center) to (25.center);
		\draw [style=directed] (27) to (26);
		\draw [style=to] (29.center) to (30.center);
		\draw [style=dash-to] (31.center) to (32.center);
		\draw [style=directed] (35) to (33);
		\draw [style=directed] (36) to (35);
		\draw [style=directed] (36) to (34);
		\draw [style=directed] (34) to (33);
		\draw [style=directed-dash] (35) to (34);
		\draw [style=directed] (36) to (33);
		\draw [style=to] (44.center) to (45.center);
		\draw [style=dash-to] (46.center) to (47.center);
	\end{pgfonlayer}
    \end{tikzpicture}
    \end{adjustbox}
    \caption{Some horn fillers. The white vertex, the red edge and the shaded triangle all denote the $(n-1)$-face added as a result of filling a horn of the respective dimension $n$. Note that filling a $(1,j)$-horn for $j=0$ (respectively $j=1$) is the same as finding an edge starting from (respectively ending to) an existing point to a new one.}
    \label{fig:horn-fillers}
\end{figure}

The lifting problems with respect to the horn inclusions $\Lambda^m_k \into \Delta^m$ lead to the Kan conditions and the definition of higher groupoids. Solutions of such lifting problems are called \textbf{horn fillers}. See Figure \ref{fig:horn-fillers} for a pictorial representation of these. 

\begin{definition}\label{def:n-gpd}
    Let $n$ be a non-negative integer or $\infty$. An \textbf{$n$-groupoid} $\huaX$ is a simplicial set satisfying the Kan conditions $\Kan(m,k)$ for any $m\ge 1$, $0\le k\le m$ and the strict Kan conditions $\Kan!(m,k)$ for any $m \ge n+1$, $0 \le k \le m$. 
    These conditions are
    \begin{itemize}[leftmargin=*]
        \item $\Kan(m,k)$: The horn projection $p^m_k: \hom(\Delta^m,\huaX) \to \hom(\Lambda^m_k,\huaX)$ is surjective.
        \item $\Kan!(m,k)$: The horn projection $p^m_k: \hom(\Delta^m,\huaX) \to \hom(\Lambda^m_k,\huaX)$ is bijective.
    \end{itemize}
    An \textbf{$\infty$-groupoid} is also known as a \textbf{Kan complex}.  
    An \textbf{$n$-group} is an $n$-groupoid $\huaX$ with $\huaX_0 = pt$. 
    We say an $\infty$-groupoid has \textbf{order} $n$ if it is an $n$-groupoid but not an $(n-1)$-groupoid. 
\end{definition}

Intuitively, the Kan condition $\Kan(m,k)$ states the following: 
For any configuration of $(m-1)$-simplices in $\huaX$ forming an $(m,k)$-horn, that is, an $m$-tuple as in \eqref{eq:DefHornBdrySpaces}, there exists a horn filler, which is an $m$-simplex $x \in \huaX_m$ such that $d_i x = x_i$ for any $i\neq k$. 
For $\Kan!(m,k)$, this horn filler is unique. 

\begin{example}\label{ex:0gpd}
    Every set $X$ has an associated identity groupoid $X$ with $X$ at each level and all simplicial maps equal to the identity of $X$. This is a 0-groupoid. 
    This defines an embedding $\Set \into \SSet$ which consists in associating each set to the constant functor at that set.

    Moreover, any 0-groupoid $\huaX$ has all levels isomorphic to $\huaX_0$ and all its simplicial maps are isomorphisms. In fact, by using the simplicial identities, it can be shown inductively that, at every level $m$, $d_i^m=\phi^m$ for any $0\le i \le m$ and $s_i^m=(\phi^m)^{-1}$ for any $0 \le i \le m-1$ for some isomorphism $\phi^m: \huaX_m \to \huaX_{m-1}$.
    Hence 0-groupoids are identified with sets.

    Intuitively, this corresponds to the fact that if each 1-horn (i.e. a point) identifies a unique arrow ending in it and a unique arrow starting at it, these two must be the same, and this construction is the same as that of the nerve of the discrete category of a set. 
\end{example}

\begin{example}\label{ex:1gpd}
    A groupoid $\huaX_1 \rightrightarrows \huaX_0$ in the usual sense of small category where all arrows are invertible, is precisely a 1-groupoid in the sense of definition \ref{def:n-gpd}, when identified with its nerve as defined in Example \ref{ex:nerve-cat}. 
    Even more, a category is a groupoid if and only if its nerve is a Kan complex, as shown for example in \cite[\href{https://kerodon.net/tag/0037}{Tag 0037}]{kerodon}.

    In our convention, the source and target maps are $\bs=d^1_0$ and $\bt=d^1_1$, respectively. The unit map is $1=s^0_0$ and the multiplication is $\mathsf{m} = d_1(p^2_1)^{-1}: \Lambda^2_1(\huaX) = \huaX^{(2)} \to \huaX_1$,  where we write $(g,h) \in \huaX^{(2)}$ if $\bt(g) = \bs(h)$ and $\mathsf{m}(g,h) = h\cdot g$.  

    In this example, because of the strict Kan conditions starting at level 2, every 2-horn can be uniquely filled. This means that for any two arrows with a common endpoint, there exist a third one, which is identified with the third face of the unique filler. 
    Thus, the $(2,1)$-horn fillers are defining the multiplication and the $(2,0)$- or $(2,2)$-horn fillers are defining the inversion maps, and viceversa. 
    Going higher, the 3-horn fillers correspond to associativity. 
    This finiteness of the data needed to describe a 1-groupoid is a general feature of all $n$-groupoids that we discuss in Section \ref{sec:finite-data}.
\end{example}

Kan complexes have two important properties:
\begin{itemize}
    \item Homotopy of simplices is an equivalence relation. See e.g. \cite[Prop. 3.2]{May1967}. This allows to define \textbf{homotopy groups} for a Kan complex $\huaX$, relative to a vertex $p \in \huaX_0$ by setting
    \begin{equation*}
        \pi_n(\huaX, p) = \{ X \in \huaX_n \mid d_i(X) = 1p, \quad \forall 0 \le i \le n \}/\sim, 
    \end{equation*}
    for any $n \ge 1$, where $X \sim Y$ if there exists an $(n+1)$-simplex $W$ such that $d_0 W = X$, $d_1 W = Y$, and $d_i W = 1p$ for all $2 \le i \le n+1$. 
    Similarly, the \textbf{set of connected components} of $\huaX$ is 
    \begin{equation*}
        \pi_0(\huaX)=\{p \in \huaX_0\}/\sim,
    \end{equation*}
    with the same equivalence relation above. 
    These $\pi_n(\huaX)$ for $n \ge 1$ have group structures as defined for example in \cite[\S I.4]{May1967}.

    \item Any mapping space $\Maps(\huaX, \huaY)$ where $\huaY$ is a Kan complex is a Kan complex. See e.g. \cite[Thm. 6.9]{May1967}, \cite[Cor. I.5.3]{GoerssJardine2009}, \cite[\href{https://kerodon.net/tag/00TJ}{Tag 00TJ}]{kerodon}. 
\end{itemize}

Specializing the first property to $n$-\hspace{0pt}groupoids, we have that if $\huaX$ is an $n$-\hspace{0pt}groupoid, then its homotopy groups for $m>n$ are trivial: $\pi_m(\huaX, p) = \{1p\}$ for all $m > n$. This follows immediately by the strict Kan conditions at level $m > n$, as because of them, the only $m$-simplex in $\huaX_m$ with all faces equal to $1p \in \huaX_{m-1}$ is exactly $1p \in \huaX_m$. We make the following definitions.

\begin{definition}\label{def:htpy-equiv-n-type-set}
    A map $f:\huaX \to \huaY$ between two simplicial sets is called a \textbf{homotopy equivalence} if there exist a map $g: \huaY \to \huaX$ such that $g\circ f \simeq id_{\huaX}$ and $f \circ g \simeq id_{\huaY}$. That is, there exist simplicial homotopies as defined in Definition \ref{def:simp-homotopy-of-maps} between each composition and the respective identity. 

    A map $f:\huaX \to \huaY$ between two Kan complexes $\huaX$ and $\huaY$ is called a \textbf{weak equivalence} if it induces isomorphisms on all homotopy groups, that is 
    \begin{equation*}
        \begin{cases}
            \pi_0(\huaX) \overset{\cong}{\underset{\pi_0(f)}{\to}} \pi_0(\huaY),\\
            \pi_m(\huaX, p) \overset{\cong}{\underset{\pi_m(f)}{\to}} \pi_m(\huaY, f(p)), \quad \forall m \ge 1, \forall p \in \huaX_0.
        \end{cases}
    \end{equation*}

    We say a Kan complex $\huaX$ is an \textbf{$n$-type} if the homotopy groups $\pi_m(\huaX, p)$ are trivial for any $m > n$ and any $p \in \huaX_0$, and non-trivial for $m=n$. If $\huaX$ is a 0-type and $\pi_0(\huaX) = pt$, we say $\huaX$ is \textbf{contractible}. 
\end{definition}

\begin{example}
    Any $n$-groupoid is at most an $n$-type. In particular, the identity groupoid of any set is a 0-type. The pair groupoid of a set is contractible.
\end{example}

\begin{remark}
Homotopy equivalences of Kan complexes are the same as weak equivalences. This is a simplicial version of the Whitehead theorem that can be found for example in \cite[\href{https://kerodon.net/tag/00V2}{Tag 00V2}]{kerodon}. 
With this, $n$-types should be thought of as ``$n$-groupoids up to homotopy''. 
\end{remark}

As for the second property, we could not find an improved version of this result for $n$-groupoids in the literature, so we write one here, by assembling various other results. This result previously appeared in \cite{RonchiZhu2024}. 

\begin{theorem}\label{thm:hom-ngpd}
    Let $\huaX$ be a simplicial set, and $\huaY$ be an $n$-groupoid. Then $\Maps(\huaX, \huaY)$ is an $n$-groupoid.
\end{theorem}
\begin{proof}
    First of all, by the aforementioned \cite[Thm. 6.9]{May1967}, $\Maps(\huaX, \huaY)$ is an $\infty$-groupoid. Thus we only need to prove the strict Kan condition $\Kan(m, k)!$ for $m\ge n+1$. 
    To prove this, we use the theory of \textit{anodyne extensions} and the standard simplicial model structure on the category of simplicial sets, reviewed in e.g. \cite[Ch. I]{GoerssJardine2009}.
    A fibration (or Kan fibration) of simplicial sets is a map with the right lifting property with respect to trivial cofibrations. That is, in the lifting problem \eqref{diag:lifting-problem}, fibrations are the maps $f:\huaX \to \huaY$ such that if $i$ is a trivial cofibration, a lift always exists. In simplicial sets, trivial cofibrations are precisely the anodyne extensions, which are the saturated class of morphisms generated by the horn inclusions $i_{m,k} : \Lambda^m_k \to \Delta^m$.\footnote{A particular kind of anodyne extensions is given by the collapsible extensions we consider in Section \ref{sec:coll-ext}. We also give a definition of Kan fibrations of simplicial manifolds in \ref{def:Kan-fibration}.}  
    Then \cite[Cor. 4.6]{GoerssJardine2009} states that if $i:\huaK \to \huaL$ is an anodyne extension and $\huaC\to \huaD$ is an inclusion, the map between the pushout $(\huaK \times \huaD) \cup (\huaL \times \huaC)$ and $\huaL \times \huaD$ induced by the universal property is an anodyne extension. 
    If we take $\huaK = \Lambda^m_k$, $\huaL= \Delta^m$, $i = i_{m,k}$, and $\huaC=\huaD=\huaX$, then the map $id \times i_{m,k}: \huaX \times \Lambda^m_k \to \huaX \times \Delta^m$ is an anodyne extension, hence a trivial cofibration. 
    Note that for $m \ge n+1$ and any $0\le k\le m$, this is an isomorphism between the $(n-1)$-truncations $(\Lambda^m_k)_{\le n-1} \to (\Delta^m)_{\le n-1}$. Hence we can apply \cite[Lemma 2.14]{Pridham2013} to the $n$-groupoid $\huaY$, and $id \times i_{m,k}: \huaX \times \Lambda^m_k \to \huaX \times \Delta^m$, which is a trivial cofibration and an isomorphism on the $(n-1)$-truncations. With this, the induced maps 
    \begin{equation*}
    \hom(\huaX \times \Delta^m, \huaY) \to \hom(\huaX \times \Lambda^m_k, \huaY)    
    \end{equation*}
    are isomorphisms. But these are precisely the horn projections of the mapping set $\Maps(\huaX, \huaY)$ for $m\ge n+1$.  
\end{proof}

\subsection{Higher Lie groupoids}

Let $\Mfd$ be the category of (finite-dimensional) smooth manifolds.\footnote{This discussion can be adapted to simplicial Banach manifolds as in the work of \cite{Henriques2008,Zhu2009,RogersZhu2020}, but since we will be discussing finite-dimensional vector bundles and their duals, we have no need for this higher generality.}
To define Lie $n$-groupoids, we first need to discuss spaces of horn-shaped simplicial diagrams inside simplicial manifolds. 
A problem arises here because the sets \eqref{eq:DefHornBdrySpaces} need not have a canonical smooth manifold structure, as they are constructed through fiber products (i.e. pullbacks), which are particular limits. 
Not all fiber products exist in the category of smooth manifolds $\Mfd$, because of transversality issues.
However, fiber products (and more general limits) always exist in the category of sheaves over manifolds with respect to the Grothendieck pretopology of surjective submersions. 
In this section we describe how to leverage this to define certain diagram spaces before knowing whether they exist as manifolds. We follow \cite[\S 2]{Henriques2008},\cite[\S 1, 2]{Zhu2009},\cite[\S 3.2]{Li2014}, \cite[\S 2]{Wolfson2016}, \cite[\S 3.4, \S 4.2]{RogersZhu2020} as references for this discussion. A similar approach was also used in \cite{BehrendGetzler2017} to discuss $n$-groupoids in descent categories, which unlike $\Mfd$, are required to have all finite limits. The discussion that follows is phrased in terms of the category of smooth manifolds for the sake of concreteness, but everything works for any category with a Grothendieck pretopology satisfying the assumptions in the above references.

\begin{definition}
Consider the Yoneda embedding 
\begin{equation*}
\begin{split}
    \Yoneda: \Mfd &\into \PSh(\Mfd)\\
    M &\mapsto \Yoneda M := \Mfd(\_, M),
\end{split}
\end{equation*}
of the category of manifolds into the category of presheaves over manifolds (i.e. functors $\Mfd^{op} \to \Set$). 
We say a presheaf $\huaF$ is \textbf{representable} if $\huaF = \Yoneda M$ for some $M$.
\end{definition}

The category of smooth manifolds has a specific structure that makes it so that every \textit{representable} presheaf on $\Mfd$ is actually a sheaf: the subcanonical Grothendieck pretopology of surjective submersions. Recall that to talk about sheaves on a category $\Cat$ one needs to introduce the notion of a \textbf{Grothendieck pretopology}. This is a collection $\covers$ of morphisms in $\Cat$ such that it contains isomorphisms and it is closed under compositions and pullbacks.\footnote{I.e. the map opposite a cover in a pullback diagram is a cover itself. This operation is also known as base change.} Additionally we assume that it contains the canonical maps from any object to the terminal object of the category (\cite[Def. 2.1]{Henriques2008}, \cite[Def. 3.1]{RogersZhu2020}). Morphisms in the class $\covers$ are called \textbf{covers}.

The class of surjective submersions in $\Mfd$ is a Grothendieck pretopology. 
In fact, pullbacks where one of the maps in the starting diagram is a surjective submersion exist, 
and the pullback of a surjective submersion is a surjective submersion, as pictured in the diagram
\[\begin{tikzcd}[ampersand replacement=\&,column sep=small,row sep=scriptsize]
	X \& Y \& {\exists X\times_YZ \in \Mfd} \& Z \\
	\&\& X \& Y
	\arrow[two heads, from=1-1, to=1-2]
	\arrow[Rightarrow, from=1-2, to=1-3]
	\arrow[dotted, two heads, from=1-3, to=1-4]
	\arrow[dotted, from=1-3, to=2-3]
	\arrow["\lrcorner"{anchor=center, pos=0.125}, draw=none, from=1-3, to=2-4]
	\arrow[from=1-4, to=2-4]
	\arrow[two heads, from=2-3, to=2-4]
\end{tikzcd}\]
where we denote surjective submersions by $\onto$. The other properties can be trivially checked.

Additionally, surjective submersions form a \textbf{subcanonical} pretopology, which means that all representable presheaves over $\Mfd$ are actually \textbf{sheaves}. That is, every presheaf $\huaF \in \PSh(\Mfd)$ such that $\huaF = \Yoneda M$ for some $M\in \Mfd$, satisfies the \textbf{descent condition}: for every surjective submersion $U \onto X$, $\huaF(X)$ is the equalizer of the diagram
\[\begin{tikzcd}[ampersand replacement=\&, column sep=small]
	{\huaF(U)} \& {\huaF(U\times_X U).}
	\arrow[shift right, from=1-1, to=1-2]
	\arrow[shift left, from=1-1, to=1-2]
\end{tikzcd}\]

We will henceforth consider the Yoneda embedding as a map $\Yoneda: \Mfd \into \Sh(\Mfd)$, and consider limits in $\Mfd$ a priori as sheaves.
As anticipated, this works firstly because the category of sheaves is complete. So even if the limit of a diagram does not exist in $\Mfd$, the limit of the corresponding embedded diagram still exists in $\Sh(\Mfd)$. 
Secondly, the Yoneda embedding preserves and reflects limits (\cite[Thm. 3.4.6]{Riehl2016}). Hence, a limit of a diagram of representable sheaves in $\Sh(\Mfd)$ is representable if and only if the original diagram had a limit in $\Mfd$, which is the representing object. 
For more on pretopologies and other related technicalities we refer to \cite[\S 3.2]{RogersZhu2020}. 

With this we can define objects of simplicial diagrams for simplicial manifolds.

\begin{definition}\label{def:simp-diag-sheaf}
    Let $\huaS$ be a levelwise finite simplicial set and $\huaX$ be a simplicial manifold. The \textbf{object of $\huaS$-shaped simplicial diagrams in $\huaX$} is the sheaf 
    \begin{equation}
        \hom(\huaS, \huaX)(U) := \Simp\Mfd(U \otimes \huaS, \huaX) \in \Sh(\Mfd), 
    \end{equation}
    where $U \otimes \huaS$ is the copowering (Definition \ref{def:copowering}) of the identity groupoid of $U$ with $\huaS$. When this sheaf is representable we denote the representing object by the same name $\hom(\huaS, \huaX)$, by a slight abuse of notation, 
    and call this the \textit{space} of $\huaS$-shaped simplicial diagrams in $\huaX$. 
\end{definition}

To characterize $\hom(\huaS, \huaX)$, we write it as a limit of representable sheaves. For this we need a variant of the Yoneda embedding:
for any simplicial manifold $\huaX$, we define the simplicial-set-valued sheaf $\Yoneda\huaX \in \Simp\Sh(\Mfd)$ as the sheaf associating to each $U\in\Mfd$ the simplicial set defined at level $n$ by 
\begin{equation}\label{eq:Yoneda-simp-sheaves}
    \Yoneda\huaX_n(U) := \Mfd(U, \huaX_n), 
\end{equation}
with simplicial maps given by postcomposition of the respective simplicial maps of $\huaX$.

Additionally, for any sheaf $F \in \Sh(\Mfd)$ and any  countable set $S$ we define the sheaf $\hom(S, F)$ as
\begin{equation*}
    \hom(S, F)(U) := \prod_{r\in S} F(U) = \Mfd (S, F(U)),
\end{equation*}
where $S$ is seen as a discrete manifold. 
This is also the set of functions between $S$ and the underlying set of $F(U)$, or alternatively the powering of $F$ as a sheaf over the underlying set of $S$. \footnote{The latter is the general way of defining the powering over sets in categories with all self-products, as mentioned in Example \ref{ex:pow-copow-over-sets}.}

\begin{proposition}\label{prop:simp-diag-sheaf-as-limit}
Let $\huaS$ be a levelwise countable simplicial set and $\huaX$ be a simplicial manifold. Then, the sheaf $\hom(\huaS, \huaX)$ fits in the equalizer diagram 
    \begin{equation}\label{diag:simp-diag-sheaf-equalizer}
        \hom(\huaS,  \huaX)
        \to \prod\limits_{l\ge 0}\hom(\huaS_l, (\Yoneda\huaX)_l)
        \underset{\beta}{\overset{\alpha}{\rightrightarrows}}
        \prod\limits_{{\begin{array}{c} {\scriptstyle m,n \ge 0} \\ {\scriptstyle f \in \Delta([m],[n])}\end{array}}}\hom(\huaS_n, (\Yoneda\huaX)_m),
    \end{equation}
    where the second product is over all morphisms in $\Delta$, and the components of $\alpha$ and $\beta$ relative to any $f\in \Delta([m],[n])$ are
    \begin{equation}\label{eq:simp-diag-sheaf-equalizer-maps}
        \begin{split}
            &\alpha(g) = \huaX(f) \circ g_n : \huaS_n \to \huaX_m,\\
            &\beta(g) = g_m \circ \huaS(f) : \huaS_n \to \huaX_m.
        \end{split}
    \end{equation}
\end{proposition}

\begin{proof}
    By Definition \ref{def:simp-diag-sheaf}, $\hom(\huaS, \huaX)(U)$ is the set of maps $\phi: U \otimes \huaS \to \huaX$ such that 
    \begin{equation*}
        d_i\phi_n(u^r) = \phi_{n-1}(u^{d_ir}),
        \quad s_i\phi_n(u^r) = \phi_{n+1}(u^{s_ir}),
    \end{equation*}
    for any $u \in U$, $r\in \huaS_n$ and $n \ge 0$. 
    On the other hand, by evaluating \eqref{diag:simp-diag-sheaf-equalizer} on a manifold $U$, we have
    \begin{equation*}
        \left(\prod\limits_{l\ge 0}\hom(\huaS_l, (\Yoneda\huaX)_l)\right) (U)
        = \prod\limits_{l\ge 0}\Set(\huaS_l, \Mfd(U,\huaX_l)).
    \end{equation*}
    By Lemma \ref{lem:Delta-gen-by-cf-cd}, the equalizer \eqref{diag:simp-diag-sheaf-equalizer} evaluated on $U$ is the set of tuples of maps $(g_l)_{l\ge 0}$ such that 
    \begin{equation*}
        g_l: \huaS_l \to \Mfd(U,\huaX_l)
    \end{equation*}
    and
    \begin{equation*}
        g_{n-1}(d_ir)(u) = d_i(g_n(r)(u)), \quad 
        g_{n+1}(s_ir)(u) = s_i(g_n(r)(u)),
    \end{equation*}
    for any $u \in U$ and $r\in \huaS_n$ for any $n \ge 0$, i.e. $\hom(\huaS, \huaX)(U)$.
\end{proof}

\begin{remark}\label{rem:underlying-set-of-hom-sheaf}
    The maps \eqref{eq:simp-diag-sheaf-equalizer-maps} in the equalizer diagram \eqref{diag:simp-diag-sheaf-equalizer} are encoding the condition of $g$ being a natural transformation between the functors $\huaS$ and $\huaX$ from $\Delta^{op}$ to manifolds, where $\huaS$ is seen as a simplicial manifold with the discrete topology. 
    By Lemma \ref{lem:Delta-gen-by-cf-cd} this is the same as a reformulation of the conditions for $g$ to be a simplicial map in Proposition-Definition \ref{prop-def:simp-obj}.

    The discrete simplicial manifold $\huaS$ can more formally be written as $pt \otimes \huaS$, which is exactly describing the construction of adding a point for each simplex of $\huaS$. This is the left adjoint functor to the forgetful functor $\Forget: \Simp\Mfd \to \SSet$. Therefore, we have that if $\hom(\huaS, \huaX)$ is representable, then its underlying set is 
    \begin{equation*}
        \hom(\huaS, \huaX)(pt) = \Simp\Mfd(pt \otimes \huaS, \huaX) \cong \SSet (\huaS, \Forget(\huaX)). 
    \end{equation*}
    In particular, if $\huaX$ is just a simplicial set, $\hom(\huaS, \huaX)$ is representable in $\Set$ and it coincides with $\Maps(\huaS, \huaX)_0 = \SSet(\huaS, \huaX)$, consistently with Definition \ref{def:simpl-diags}. 
\end{remark}

\begin{example}
By considering $\Delta^k$ as a simplicial manifold with the discrete topology, and using the same argument as in the beginning of Section \ref{sec:simp-diag-sets}, we have that by the Yoneda Lemma, 
\begin{equation*}
    \hom(\Delta^n, \huaX) = \Fun(\Delta^{op}, \Mfd)(\Delta(\_, [n]), \huaX) \cong \huaX([n]) = \huaX_n.
\end{equation*}
\end{example}

We can now describe horn filling problems for simplicial manifolds, leading to the definition of Lie $n$-groupoids.

\begin{definition}\label{def:Lie-n-gpd}
Let $n$ be a non-negative integer or $\infty$. A \textbf{Lie $n$-groupoid} $\huaX$ is a simplicial manifold satisfying the Kan conditions $\Kan(m,k)$ for any $m\ge 1$, $0\le k\le m$ and the strict Kan conditions $\Kan!(m,k)$ for any $m \ge n+1$, $0 \le k \le m$. 
These conditions are
\begin{itemize}[leftmargin=*]
    \item $\Kan(m,k)$: The horn space $\hom(\Lambda^m_k,\huaX)$ is representable and the horn projection $p^m_k: \hom(\Delta^m,\huaX) \to \hom(\Lambda^m_k,\huaX)$ is a surjective submersion.
    \item $\Kan!(m,k)$: The horn space $\hom(\Lambda^m_k,\huaX)$ is representable and the projection $p^m_k: \hom(\Delta^m,\huaX) \to \hom(\Lambda^m_k,\huaX)$ is a diffeomorphism.
\end{itemize}
A Lie \textbf{$n$-group} is a Lie $n$-groupoid $\huaX$ with $\huaX_0 = pt$.

We say that a Lie $\infty$-groupoid has \textbf{order} $n$ if it is a Lie $n$-groupoid but not an $(n-1)$-groupoid.
\end{definition}

Contrary to the case of the standard simplex, the fact that $\hom(\Lambda_j^k,X)$ is a manifold is not trivial. As such it is required to hold in the definition of Kan conditions. A posteriori, each of these sheaves for a fixed level $k$ is representable if the Kan conditions hold at lower levels, by an inductive argument. We discuss this in the next section together with a formula to compute the dimension of the horn spaces.

\begin{definition}
    More generally, an \textbf{$n$-groupoid object} in a general category $(\Cat, \covers)$ with a terminal object $*$ and a subcanonical Grothendieck pretopology $\covers$ such that all canonical maps $X \to *$ are covers (see \cite[Assumptions 2.1]{Zhu2009}), is a simplicial object $\huaX$ satisfying the Kan conditions $\Kan(m,k)$ for any $m\ge 1$, $0\le k\le m$ and $\Kan!(m,k)$ for any $m \ge n+1$, $0 \le k \le m$, where 
    \begin{itemize}[leftmargin=*]
    \item $\Kan(m,k)$: The horn space $\hom(\Lambda^m_k,\huaX)$ is representable and the horn projection $p^m_k: \hom(\Delta^m,\huaX) \to \hom(\Lambda^m_k,\huaX)$ is a \textit{cover}.
    \item $\Kan!(m,k)$: The horn space $\hom(\Lambda^m_k,\huaX)$ is representable and the projection $p^m_k: \hom(\Delta^m,\huaX) \to \hom(\Lambda^m_k,\huaX)$ is an \textit{isomorphism}.
\end{itemize}
\end{definition}

\subsection{Representability and dimension of the horn spaces}\label{sec:horn-dim}

In this section we compute the dimension of the horn spaces of a Lie $\infty$-groupoid. 
More generally, we compute the dimension of any object $\hom(\huaS, \huaX)$ where $\huaS$ is a generalized horn or union of faces, as we now define. The definition of generalized horns is taken from \cite[\S 2.2.1]{Joyal2008}, while that of union of faces is inspired by the approach of \cite{BehrendGetzler2017}. 
The difference is only notational, as they both denote the same object. Nevertheless, we find it helpful to keep both perspectives in mind.

Recall that, by \eqref{eq:simp-cosimp-description-u}, the $j$-th $k$-face of the $(k+1)$ simplex can be recovered as either the image of the face inclusion $\delta_j: \Delta^k \to \Delta^{k+1}$, denoted by $\delta_j\Delta^k$, or the simplicial subset of $\Delta^{k+1}$ generated by $d_jE_{k+1}$, which we denote by $d_j\Delta^{k+1} \subseteq \Delta^{k+1}$.
The latter can also be defined as the simplicial subset of all simplices not containing the vertex $j$, that is
\begin{align*}
    (d_j\Delta^{k+1})_m := \{ f: [m] \to [k+1] \mid \Img f \nsupseteq \{j\}\}. 
\end{align*}
In geometric terms, this is the $k$-dimensional face opposite to the vertex $j$.
With this notation, $\Delta^{k} \cong d_j\Delta^{k+1} \cong \delta_j\Delta^k$ for every $j \in \{0,\dots,k\}$.

\begin{definition}\label{def:gen-horns}
    Let $A$ be a subset of $[k]$. The \textbf{generalized horn} $\Lambda_A^{k}$ is the simplicial subset of $\Delta^{k}$ defined by
    \begin{equation*}
    \Lambda_A^{k} := \bigcup_{i\notin A} d_i\Delta^{k} 
    = \bigcup_{i\notin A} \delta_i\Delta^{k-1},
    \end{equation*}
    that is, the simplicial subset given by removing the $(k-1)$-faces with index in $A$ from the standard $k$-simplex.
    Symmetrically, the \textbf{union of $k$-faces} $V_A^k$ is the simplicial subset of $\Delta^k$ defined by 
    \begin{equation*}
        V_A^k = \bigcup_{i\in A} d_i\Delta^k
        = \bigcup_{i\in A} \delta_i\Delta^{k-1}.
    \end{equation*}
    In this notation, $\Lambda^k_A = V^k_{[k]\backslash A}$.
\end{definition}

\begin{remark}
    If $A=\{j\}$, the corresponding generalized horn is the usual horn space $\Lambda^k_j = V^k_{[k]\backslash\{j\}}$. On the other hand, if the cardinality of $A$ is $|A|=k-1$, $\Lambda_A^{k} = V^k_{[k]\backslash A} \cong \Delta^{k-1}$. Moreover, if $A \subseteq B \subseteq [k]$, then $\Lambda_B^{k}\subseteq \Lambda_A^{k}$, while $V^k_A \subseteq V^k_B$.
\end{remark}

\subsubsection{Collapsible extensions}\label{sec:coll-ext}

In this section, we introduce collapsible extensions as defined in \cite[\S 2.6]{Li2014}, \cite[Definition 3.5]{RogersZhu2020}.\footnote{Collapsible extensions are also called expansions in \cite{Wolfson2016} and \cite{BehrendGetzler2017}.}
For a general category with a Grothendieck pretopology $(\Cat, \covers)$, this is a class of maps between simplicial sets that correspond to covers under the functor $\hom(\_, \huaX)$ when $\huaX$ is an $\infty$-groupoid object.
We recall this in Lemma \ref{lem:collext_covers}, which appears in \cite[Lemma 3.20]{Li2014}, \cite[Lemma 3.7]{RogersZhu2020}, and which is the fundamental result needed to prove that horn spaces are representable. 
In the context of groupoids in descent categories Lemma \ref{lem:collext_covers} appears as Lemma 3.9 in \cite{BehrendGetzler2017}. 
Collapsible extensions are a subclass of the class of \textit{anodyne extensions}, from \cite[\S IV.2]{GabrielZisman1967}, \cite[\S I.4]{GoerssJardine2009}, \cite[\S 2.2]{Joyal2008}. 
For example, the latter class additionally contains retracts of collapsible extensions. 
Anodyne extensions are trivial cofibrations in the standard model category of simplicial sets (see \cite{GoerssJardine2009}), so \cite[Corollary 2.11]{Henriques2008} is yet another (more general) version of \ref{lem:collext_covers}. 

There are two reasons why we need collapsible extensions. 
Together with the notion of collapsible simplicial sets, which they extend, they are helpful to solve the issue of representability of horn spaces,
as we report in Lemma \ref{lem:horn-representability}. 
More importantly for computing dimensions, as we recall in Lemma \ref{lem:dim_fib_prod}, there is a formula computing the dimension of fiber products of manifolds where both maps in the starting diagram are surjective submersions. 
Collapsible extensions describe exactly the dual situation of this, as a sequence of pushouts over collapsible extensions corresponds to a sequence of fiber products over surjective submersions through the functor $\hom(\_, \huaX)$.

\begin{definition}\label{def:collext}
Let $\huaS\into \huaT$ be an inclusion of (levelwise finite) simplicial sets. We call $\huaS\into \huaT$ a \textbf{collapsible extension} if there is a finite filtration of simplicial sets
\begin{equation*}
    \huaS \cong \huaS_0 \subseteq \huaS_1 \subseteq \dots \subseteq \huaS_l \cong \huaT
\end{equation*}
such that each $\huaS_i$ is obtained from $\huaS_{i-1}$ by filling a horn. That is, for every $0\le i \le l$, there is a morphism $\Lambda^{k_i}_{j_i} \to \huaS_{i-1}$ for some $m_i\ge 1$, $0\le k_i \le m_i$, such that $\huaS_i = \huaS_{i-1} \cup_{\Lambda^{k_i}_{j_i}} \Delta^{k_i}$, the pushout of the diagram $\Delta^{k_i} \hookleftarrow \Lambda^{k_i}_{j_i} \to \huaS_{i-1}$. 

If there is a collapsible extension $\Delta^0 \into \huaT$, then we simply call $\huaT$ \textbf{collapsible}.
\end{definition}

\begin{example}
    By definition, isomorphisms are collapsible extensions with filtration of length 0 and horn inclusions $\Lambda^k_j\into \Delta^k$ are collapsible extensions with filtration of length 1.
\end{example}

\begin{example}
Both the horn $\Lambda_j^k$ and the simplex $\Delta^k$ are collapsible for all $k\ge 0$, $0\le j \le k$.
\end{example}

\begin{remark}
    The composition of two collapsible extensions is a collapsible extension with respect to the filtration obtained by concatenation of the respective filtrations. Moreover, collapsible extensions are stable under pushouts, as we show in the following lemma. 
\end{remark}

\begin{lemma}\label{lem:CE-pushout}
    Given the pushout diagram of simplicial sets
\[\begin{tikzcd}
	\huaU & \huaT \\
	\huaS & {\huaS\cup_{\huaU} \huaT,}
	\arrow[from=1-1, to=2-1]
	\arrow[from=1-1, to=1-2]
	\arrow[from=2-1, to=2-2]
	\arrow[from=1-2, to=2-2]
	\arrow["\lrcorner"{anchor=center, pos=0.125, rotate=180}, draw=none, from=2-2, to=1-1]
\end{tikzcd}\]
if $\huaU\to \huaT$ is a collapsible extension, then so is $\huaS\to \huaS\cup_\huaU \huaT$. Analogously, if $\huaU \to \huaS$ is a collapsible extension, so is $\huaT\to \huaS\cup_\huaU \huaT$.
\end{lemma}
\begin{proof}
Assume $\huaU \to \huaT$ is a collapsible extension (the other case follows by the same argument, exchanging $\huaS$ and $\huaT$). Then there is a filtration by horn fillings 
\begin{align*}
    \huaU \cong \huaS_0 \subseteq \huaS_1 \subseteq \huaS_2 \subseteq \dots \subseteq \huaS_l \cong \huaT.
\end{align*}
To construct a filtration by horn fillings from $\huaS$ to the pushout $\huaS\cup_\huaU \huaT$ we simply take the pushout with $\huaS$ at each level of the given filtration. That is, we define 
\begin{align*}
    \huaS'_i = \huaS \cup_\huaU \huaS_i,
\end{align*}
and we get $\huaS'_0 \cong \huaS \cup_\huaU \huaU = \huaS$ and $\huaS'_l \cong \huaS\cup_\huaU \huaT$. Furthermore, we have at each level
\begin{align*}
    \huaS'_i = \huaS \cup_\huaU \huaS_i &= \huaS \cup_\huaU (\huaS_{i-1} \cup_{\Lambda^k_j}\Delta^k)\\ 
    &\cong \huaS'_{i-1} \cup_{\Lambda^k_j}\Delta^k\\
\end{align*}
for some $k,j$. 
Thus $\{\huaS'_i\}_i$ is a filtration by horn fillings, and $\huaS\to \huaS\cup_\huaU \huaT$ is a collapsible extension. 
\end{proof}

Another important class of collapsible extensions is given by the face inclusions. To prove this, it is necessary to use the join operation of simplicial sets. Because this is outside the scope of this thesis, we refer to \cite{Li2014} for the details. 

\begin{lemma}[{\cite[Lemma 2.44]{Li2014}}]\label{lem:ce_face_incl}
The face inclusions $\delta_I: \Delta^l \into \Delta^m$ are collapsible extensions, for any $0 \le l \le m$ and multi-index $I$ with $i_j\le l+j$ and length $|I| = m-l$. 
\end{lemma}

\begin{example}
    By Lemma \ref{lem:ce_face_incl}, the inclusion of a horn into a higher dimensional simplex $\Lambda^k_j\into \Delta^m$ is also a collapsible extension, for all $0\le k\le m$ and $0\le j\le k$. 
\end{example}

As shown in \cite[Proposition 2.12]{Joyal2008}, \cite[Lemma 3.8]{BehrendGetzler2017}, the previous example can be generalized to show that the inclusion of any generalized horn (or union of faces) is a collapsible extension. We write out a proof in our notation.

\begin{lemma}\label{lem:ghorn_incl_ce}
    Let $A$ be a non-empty proper subset of $[k]$. Then, the inclusion 
    \begin{equation*}
        i_A: V_A^{k} \into \Delta^{k}
    \end{equation*}
    is a collapsible extension.
    \end{lemma}
    
    \begin{proof}
        We argue by induction on the cardinality of the set $[k]\backslash A$, which we denote $|[k]\backslash A|=r$. If $r=1$, then $V_A^{k} = \Lambda_j^{k}$ for some $j\in[k]$, so $i_A$ is a horn inclusion. 
        Suppose now $r>1$. After choosing an element $a\notin A$, and setting $B=A\cup\{a\}$, we can write
        \begin{align}\label{eq:unionVAB}
            V_B^{k}=V_A^{k} \cup d_a\Delta^{k}.
        \end{align}
        And we have that $i_A$ factors through $i_B$ as follows:
        \begin{align*}
            V_A^{k}\subseteq V_B^{k} \subseteq \Delta^{k},
        \end{align*}
        where $i_B: V_B^{k} \into \Delta^{k}$ is a collapsible extension by inductive hypothesis, since $B$ is non-empty and $|[k]\backslash B| = r-1$. It now suffices to prove that the inclusion $V_A^{k}\subseteq V_B^{k}$ is a collapsible extension. By (\ref{eq:unionVAB}) we have the pushout square
        \begin{equation}\label{diag:pushout_AB}
        \begin{tikzcd}
        {d_a\Delta^{k} \cap V_A^{k}} & {V_A^{k}} \\
        {d_a\Delta^{k}} & {V_B^{k}}
        \arrow[from=1-1, to=2-1]
        \arrow[from=1-1, to=1-2]
        \arrow[from=1-2, to=2-2]
        \arrow[from=2-1, to=2-2]
        \arrow["\lrcorner"{anchor=center, pos=0.125, rotate=180}, draw=none, from=2-2, to=1-1]
        \end{tikzcd}
        \end{equation}
        So if the inclusion $d_a\Delta^{k} \cap V_A^{k} \subseteq d_a\Delta^{k}$ is a collapsible extension, $V_A^{k}\subseteq V_B^{k}$ is a collapsible extension.
        
        We show this by observing that the intersection $d_a\Delta^{k} \cap V_A^{k}$ is isomorphic to a union of faces $V_S^{k-1}$, for some $S \subseteq [k-1]$. Consider the coface map $\delta_a: [k-1] \to [k]$ from \eqref{eq:Delta-coface-def}. 
        This map corresponds to the face inclusion $\delta_a: \Delta^{k-1}\to\Delta^{k}$ which induces the isomorphism $\Delta^{k-1}\cong d_a\Delta^{k}$.
        Take $S := \delta_a^{-1}(A)$. Then we have that $\delta_a: \Delta^{k-1}\to\Delta^{k}$ also induces an isomorphism by restricting to the faces with index in $A$, so we get
    \begin{equation}\label{diag:intersection_S-horn}
    \begin{tikzcd}
        {V_S^{k-1}} & {d_a\Delta^{k}\cap V_A^{k}} \\
        {\Delta^{k-1}} & {d_a\Delta^{k}.}
        \arrow[hook, from=1-1, to=2-1]
        \arrow["\cong", from=2-1, to=2-2]
        \arrow["\cong", from=1-1, to=1-2]
        \arrow[hook, from=1-2, to=2-2]
    \end{tikzcd}
    \end{equation}
    The inclusion on the left is a collapsible extension by the inductive hypothesis, since $S$ is a subset of $[k-1]$ with cardinality $|S|=|A|$, and thus $|[k-1]\backslash S| = r-1 < r$. Therefore the inclusion $V_S^{k-1}\subseteq \Delta^{k-1}$ is a collapsible extension, hereby concluding the proof. 
    \end{proof}

We now report the main lemma relating collapsible extensions and covers from \cite[Lemma 3.20]{Li2014} and \cite[Lemma 3.7]{RogersZhu2020}.

\begin{lemma}\label{lem:collext_covers}
Let $\huaS \to \huaT$ be a collapsible extension of simplicial sets and let $\huaX$ be an $\infty$-groupoid in $(\Cat,\covers)$ with the usual assumptions.
Assume that $\hom(\huaS,\huaX)$ is representable. Then $\hom(\huaT,\huaX)$ is representable, and the map $\hom(\huaT,\huaX) \to \hom(\huaS,\huaX)$ is a cover. 
\end{lemma}

\begin{proof}
    If $\huaS \to \huaT$ is an isomorphism, then obviously $\hom(\huaT,\huaX) \to \hom(\huaS,\huaX)$ is an isomorphism by functoriality, and both objects are representable. 
    Let $\huaS \cong \huaS_0 \subseteq \huaS_1 \subseteq \dots \subseteq \huaS_l \cong \huaT$ be a filtration as in the definition of collapsible extensions. Without loss of generality, we may assume that $l=1$, since both collapsible extensions and covers are closed under composition. Thus $\huaT\cong \huaS\cup_{\Lambda^k_j}\Delta^k$ for some $k,j$. By definition of $n$-groupoid in $(\Cat, \covers)$, we have that $\hom(\Delta^k,\huaX) \to \hom(\Lambda_{j}^k,\huaX)$ is a cover and $\hom(\Lambda_{j}^k,\huaX)$ is representable. Then
    \begin{equation*}
    	\hom(\huaT,\huaX) = \hom(\huaS,\huaX) \times_{Hom(\Lambda_{j}^k,\huaX)}\hom(\Delta^k,\huaX)
    \end{equation*}
    is representable, and the map $\hom(\huaT,\huaX) \to \hom(\huaS,\huaX)$ is a cover because it is a pullback of a cover.    
\end{proof}

\begin{example}
Let $\huaX$ be an $\infty$-groupoid in $(\Cat,\covers)$. Because face inclusions are collapsible extensions, we have that the face maps of $\huaX$, 
\begin{align*}
    d^k_i:\hom(\Delta^k,\huaX)=\huaX_k \to \hom(\Delta^{k-1},\huaX) = \huaX_{k-1} 
\end{align*}
are covers for all $i$. 
\end{example}

As a consequence of this lemma, horn spaces of increasingly higher dimension are inductively representable. For higher Lie groups this is \cite[Lemma 2.4, Corollary 2.5]{Henriques2008}. The result for general higher groupoids in $(\Cat, \covers)$ is contained in \cite[Lemma 3.9, Remark 3.10]{RogersZhu2020}. We report a specific version of this result here and refer to these references for a more general result and a complete proof.

\begin{lemma}\label{lem:horn-representability}
    Let $\huaX$ be a simplicial object in $(\Cat, \covers)$. If $\huaX$ satisfies $\Kan(l,j)$ for all $1 \le l < m$ and $0 \le j \le l$, then $\hom(\Lambda^m_k, \huaX)$ is representable for all $0 \le k \le m$ and we denote it by $\Lambda^m_k(\huaX)$.
\end{lemma}

In particular all the horn spaces of a Lie $n$-groupoid for any $n$ are representable. 

\begin{remark}\label{rem:bdry-representability}
    Lemma \ref{lem:horn-representability} uses the fact that the horn is a collapsible subset of the standard simplex, as in Definition \ref{def:collext}.
    On the contrary, the boundary $\partial\Delta^n$ is not collapsible, so the same results cannot be applied to show representability of boundary spaces.
    In fact, the sheaf $\hom(\partial\Delta^n, \huaX)$ is a priori not representable. For example consider $\hom(\partial\Delta^2, \huaX)$.
    This can be written as the fiber product
    \begin{equation*}
        \Lambda^2_2(\huaX) \times_{(d_1d_0, d_1d_1), \huaX_0 \times \huaX_0, (d_0, d_1)} \huaX_1. 
    \end{equation*}
    Even if $\huaX$ is a Lie groupoid, the fact that either of the maps in the fiber product is a surjective submersion is not guaranteed. 
    For example, if $(d_0, d_1)$ is a surjective submersion, the Lie groupoid is said to be locally trivial \cite[Def. 1.3.2]{Mackenzie2005} and it is transitive, but not all Lie groupoids are transitive. In a similar way, we cannot conclude that $(d_1d_0, d_1d_1)$ is a surjective submersion from the fact that the face maps are surjective submersions.
\end{remark}

As previously mentioned, our dimension formula for the horn spaces relies on the following lemma about the dimension of fiber products of manifolds over two surjective submersions.

\begin{lemma}\label{lem:dim_fib_prod}
Consider the following fiber product of manifolds, where $f$ and $g$ are surjective submersions.
\begin{equation}\label{diag:fib_prod_mfd}
\begin{tikzcd}
{X\times_ZY} & Y \\
X & Z
\arrow[from=1-1, to=1-2]
\arrow[from=1-1, to=2-1]
\arrow["f"', two heads, from=2-1, to=2-2]
\arrow["g", two heads, from=1-2, to=2-2]
\arrow["\lrcorner"{anchor=center, pos=0.125}, draw=none, from=1-1, to=2-2]
\end{tikzcd}
\end{equation}

Then $X\times_Z Y$ is a manifold of dimension 
\begin{equation}
    \dim X\times_Z Y = \dim X + \dim Y - \dim Z.
\end{equation}
\end{lemma}
\begin{proof}
First of all, since at least one of $f$ and $g$ is a surjective submersion, $X\times_Z Y$ is a smooth manifold. Moreover, by stability of surjective submersions under pullbacks, all the maps in the diagram are surjective submersions.

Now it is enough to compute the dimension of the tangent space at a point $T_{(x,y)}(X\times_Z Y)$. Notice that, since $f,g$ are surjective submersions, we have that, canonically, $T_{(x,y)}(X\times_Z Y) \cong T_x X \times_{T_{f(x)}Z} T_y Y$.\footnote{This can easily be shown by comparing local coordinates.}
So the following is a pullback diagram of vector spaces. 
\[\begin{tikzcd}
	{T_{(x,y)}(X\times_{Z}Y)} & {T_yY} \\
	{T_xX} & {T_{f(x)}Z}
	\arrow["T_{(x,y)}p_Y", from=1-1, to=1-2]
	\arrow["T_{(x,y)}p_X"', from=1-1, to=2-1]
	\arrow["T_xf"', from=2-1, to=2-2]
	\arrow["T_yg", from=1-2, to=2-2]
	\arrow["\lrcorner"{anchor=center, pos=0.125}, draw=none, from=1-1, to=2-2]
\end{tikzcd}\]

All maps in this diagram are surjective, as the maps in (\ref{diag:fib_prod_mfd}) are all submersions. Also, since pullbacks preserve kernels\footnote{Can be shown by straightforward diagram chasing. For a written proof see for example \cite[Prop. I.13.1]{Mitchell1965}.}, $\ker T_xf = \ker T_{(x,y)}p_Y$. Thus, by the rank-nullity theorem,
\begin{align*}
\dim T_{(x,y)}(X\times_{Z}Y) = \dim \ker T_{(x,y)}p_Y + \dim T_y Y,
\end{align*}
and
\begin{align*}
    \dim \ker T_{(x,y)}p_Y = \dim \ker T_xf = \dim T_x X - \dim T_{f(x)}Z.
\end{align*}

Hence, 
\begin{align*}
    \dim X\times_Z Y &= \dim T_{(x,y)}(X\times_{Z}Y) \\
                    &= \dim T_x X - \dim T_{f(x)}Z + \dim T_y Y \\
                    &= \dim X + \dim Y - \dim Z.
\end{align*}
\end{proof}

\begin{figure}[!ht]
    \centering
    \includegraphics[width=\textwidth]{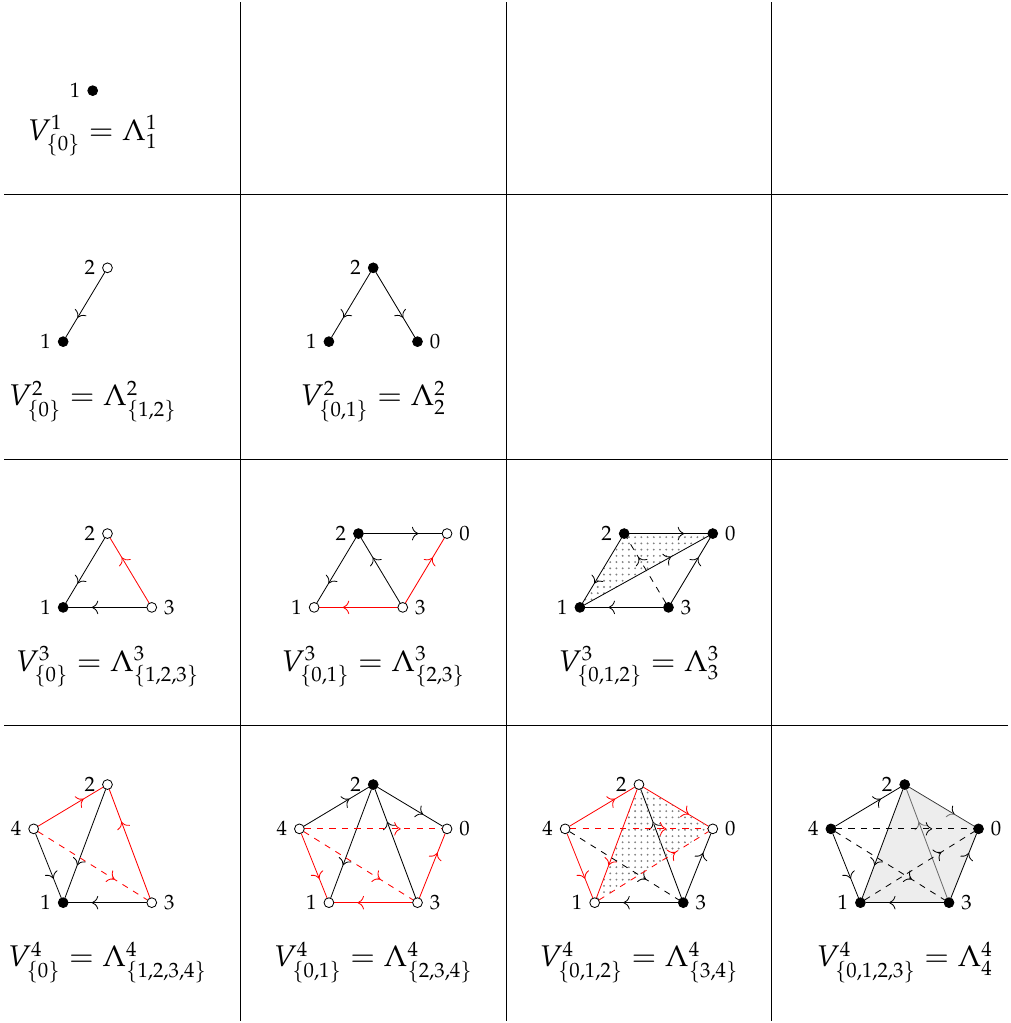}
    \caption{Some generalized horns/unions of faces. In red with white vertices the embedding of $V^{k-1}_{S_i}$ into $V^{k}_{A_i}$ used to construct $V^k_{A_{i+1}}$. All triangles are assumed to be filled in, except the dotted ones. Likewise, all tetrahedrons without missing triangles are assumed to be filled in, except the one shaded in gray. Dashed arrows are only for improving readability.}
    \label{fig:gen_horns}
\end{figure}

We now show that each union of faces can be written as a series of pushouts over a diagram where each map is a collapsible extension, so that after applying $\hom(\_,\huaX)$, the resulting generalized horn space is a series of fiber products over surjective submersions, thus falling into the hypotheses of Lemma \ref{lem:dim_fib_prod}. To do this we show exactly how horn spaces are constructed inductively by adding faces. See Figure \ref{fig:gen_horns}. 

\begin{proposition}\label{prop:pushout-filtration}
For every union of $(k-1)$-faces $V_A^k$, with $|A|=s$, we have a filtration
\begin{align*}
    \Delta^{k-1} \cong \huaT_1 \subseteq \huaT_2 \subseteq \dots \subseteq \huaT_s = V_A^k,
\end{align*}
where at each step, $\huaT_i$ is a union of faces, and
\begin{align*}
    \huaT_{i+1}\cong \huaT_i \cup_{V_{S_i}^{k-1}} \Delta^{k-1}
\end{align*}
is a pushout over the collapsible extensions $\huaT_i \leftarrow V_{S_i}^{k-1} \rightarrow \Delta^{k-1}$ where $V_{S_i}^{k-1}$ is a union of $(k-2)$-faces. 
\end{proposition}

\begin{proof}
    We enumerate the set $A$ by writing $A = \{j_1,\dots, j_s\}$, with $0\le j_1<j_2<\dots<j_s\le k$. Define
    \begin{align*}
        \huaT_1 &= d_{j_1}\Delta^k\\
        \huaT_i &= \bigcup_{l=1}^i d_{j_i}\Delta^k = V_{\{j_1,\dots,j_i\}}^k.
    \end{align*}
    For every $i=1,\dots,s-1$, we define $\huaN_i = \huaT_i \cap d_{j_{i+1}}\Delta^k$. Then clearly we can write the pushout $\huaT_{i+1}=\huaT_i \cup_{\huaN_i} d_{j_{i+1}}\Delta^k$.
    At this point all we need to prove is that for all $i$ as above, $\huaN_i \cong V_{S_i}^k$ for some $S_i$, and that the maps $d_{j_{i+1}}\Delta^k \leftarrow \huaN_i \to \huaT_i$
    are collapsible extensions. 
    
    For the first statement, we proceed as in the proof of Lemma \ref{lem:ghorn_incl_ce}. Set $S_i := \delta_{j_{i+1}}^{-1}(\{j_1,\dots,j_i\})$, for all $i=1,\dots,s-1$. Recall that the map $\delta_{j_{i+1}}: [k-1]\to[k]$ induces the face inclusion $\delta_{j_{i+1}}: \Delta^{k-1} \to \Delta^k$, which is an isomorphism on its image $d_{j_{i+1}}\Delta^k$. By restricting this map to $V_{S_i}^{k-1}$ we get an isomorphism
    \begin{equation*}
        V_{S_i}^{k-1} \overset{\cong}{\longrightarrow} \huaT_i \cap d_{j_{i+1}}\Delta^k =\huaN_i.
    \end{equation*}   
    Hence, by composing this with the original pushout diagram for $\huaT_{i+1}$, we have
    \[\begin{tikzcd}
	{V_{S_i}^{k-1}} & {\huaN_i} & {\huaT_i} \\
	{\Delta^{k-1}} & {d_{j_{i+1}}\Delta^k} & {\huaT_{i+1}}
	\arrow["\cong", from=1-1, to=1-2]
	\arrow[from=1-1, to=2-1]
	\arrow["\cong", from=2-1, to=2-2]
	\arrow[from=1-2, to=2-2]
	\arrow[from=1-2, to=1-3]
	\arrow[from=1-3, to=2-3]
	\arrow[from=2-2, to=2-3]
	\arrow["\lrcorner"{anchor=center, pos=0.125, rotate=180}, draw=none, from=2-3, to=1-2]
    \end{tikzcd}\]
    In other words $\huaT_{i+1}\cong \huaT_i \cup_{V_{S_i}^{k-1}} \Delta^{k-1}$.
    Furthermore, the inclusion $V_{S_i}^{k-1} \to \Delta^{k-1}$ is a collapsible extension by lemma \ref{lem:ghorn_incl_ce}, so $\huaN_i \to d_{j_{i+1}}\Delta^k$ is a collapsible extension as well.
    
    It remains to show that the map $\huaN_i\to \huaT_i$ is a collapsible extension for all $i = 1,\dots,{s-1}$. We have that
    \begin{align*}
        \huaN_i &= \delta_{j_{i+1}}V_{S_i}^{k-1} 
        = \delta_{j_{i+1}}\left(\bigcup_{l=1}^i d_{j_l}\Delta^{k-1}\right)
        = \bigcup_{l=1}^i d_{j_l}d_{j_{i+1}}\Delta^{k},
    \end{align*}
    by \eqref{eq:Delta-simp-cosimp-commutativity}. Because the $j_l\in A$ are ordered, by the simplicial identities \eqref{eq:simp-id}
    \begin{align*}
        \huaN_i = \bigcup_{l=1}^i d_{(j_{i+1}-1)}d_{j_{l}}\Delta^{k},
    \end{align*}
    where every element of the union is a $(k-2)$-face of a $(k-1)$-face already in $T_i$.
    Thus, we define the filtration
    \begin{align*}
        \huaN_i =: \huaM_0 \subseteq \huaM_1 \subseteq \huaM_2 \subseteq \dots \subseteq \huaM_i = T_i,
    \end{align*}
    with 
    \begin{equation}\label{eq:defM-spaces}
        \huaM_m = \bigcup_{l=1}^m d_{j_{l}}\Delta^{k} \cup \bigcup_{l=m+1}^i d_{(j_{i+1}-1)}d_{j_{l}}\Delta^{k},
    \end{equation}
    that is the union of $\huaN_i$ with the $(k-1)$-faces of $\huaT_i$ up to the $j_{m}$-th one.
    
    \begin{figure}[!ht]
    	\centering 
    	\includegraphics[width=\textwidth]{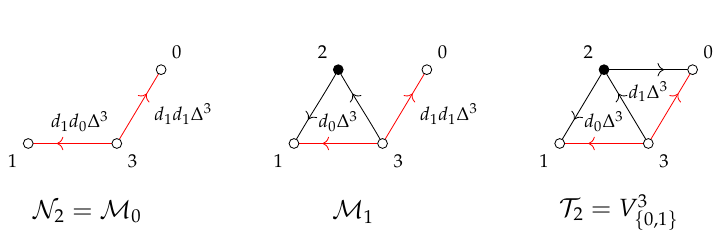}
    	\caption{Constructing $\Lambda_3^3 = V_{[3]\backslash\{3\}}^3$: The picture shows $\huaN_2$ inside $\huaT_2$ and its ``$\huaM$-filtration'' as defined above. See Figure \ref{fig:gen_horns} for the spaces $\huaT_i$ that appear in the construction of $\Lambda_3^3$.}
    \end{figure}

    We proceed by proving that each one of the inclusions in the filtration $\huaM_m \into \huaM_{m+1}$ is a collapsible extension, so that their composition, the inclusion $\huaN_i \to \huaT_i$, is a collapsible extension. 
    At each step $m$, we have the pushout diagram
    \[\begin{tikzcd}
	{\huaM_m \cap d_{j_{m+1}}\Delta^k} & {\huaM_m} \\
	{d_{j_{m+1}}\Delta^k} & {\huaM_{m+1}}
	\arrow[from=1-1, to=1-2]
	\arrow[from=1-2, to=2-2]
	\arrow[from=1-1, to=2-1]
	\arrow[from=2-1, to=2-2]
	\arrow["\lrcorner"{anchor=center, pos=0.125, rotate=180}, draw=none, from=2-2, to=1-1]
    \end{tikzcd}\]
    To show that the vertical map on the right is a collapsible extension, it suffices to show that the vertical map on the left is a collapsible extension, by stability under pushouts (Lemma \ref{lem:CE-pushout}). 
    As before, we show that the intersection $\huaM_m \cap d_{j_{m+1}}\Delta^k$ is a union of $(k-2)$-faces in $d_{j_{m+1}}\Delta^k$. Consider the order-preserving injection $\delta_{j_{m+1}}$ and call $H_{m}:= \delta_{j_{m+1}}^{-1}(\{j_1,\dots,j_m\})$. Then the face inclusion $\delta_{j_{m+1}}:\Delta^{k-1} \to d_{j_{m+1}}\Delta^k$ is an isomorphism, which restricts to an isomorphism $V_{H_m}^{k-2} \to \huaM_m \cap d_{j_{m+1}}\Delta^k$.\footnote{This can also be seen by looking at (\ref{eq:defM-spaces}), where one can see that the intersection $\huaM_m \cap d_{j_{m+1}}\Delta^k$ is composed of $d_{j_{i+1}-1}d_{j_{m+1}}\Delta^{k}$ and exactly one other $(k-2)$-face for each of the $d_{j_{l}}\Delta^{k}$ with $l=1,\dots, m$ that appear in that definition.}
	By lemma \ref{lem:ghorn_incl_ce}, the inclusion of $V_{H_m}^{k-2}$ into $\Delta^{k-1}$ is a collapsible extension. So the inclusion $\huaM_m \cap d_{j_{m+1}}\Delta^k \to d_{j_{m+1}}\Delta^k$ is a collapsible extension for all $m=0,\dots, i$, thus concluding the proof.
\end{proof}

To summarize, by writing out the filtration by pushouts we just obtained as
\begin{align*}
    \Lambda^k_{[k]\backslash A} = V^k_A 
    \cong T_{s-1} \cup_{V_{S_{s-1}}^{k-1}} \Delta^{k-1} \cong \left(T_{s-2}\cup_{V_{S_{s-2}}^{k-1}} \Delta^{k-1}\right)\cup_{V_{S_{s-1}}^{k-1}} \Delta^{k-1} \cong \dots,
\end{align*}
we can write
\begin{align*}
    \Lambda^k_{[k]\backslash A} = V^k_A 
    \cong \left(\dots\left(\left(\Delta^{k-1}\cup_{V^{k-1}_{S_1}}\Delta^{k-1}\right)\cup_{V^{k-1}_{S_2}}\Delta^{k-1}\right)\dots\right)\cup_{V^{k-1}_{S_{s-1}}}\Delta^{k-1}
\end{align*}
where all the maps in each pushout are collapsible extensions, by Lemma \ref{lem:ghorn_incl_ce} and Proposition \ref{prop:pushout-filtration}. These are the maps denoted by $\into$ in the following diagram. (Of course all the others also end up being collapsible extensions, by stability).

\[\begin{tikzcd}
	{V_{S_1}^{k-1}} & {\Delta^{k-1}} &&& {V_{S_2}^{k-1}} & {\Delta^{k-1}} \\
	{\Delta^{k-1}} & {T_1} & {T_2} && {T_{s-2}} & {T_{s-1}} & {V_{A}^k} \\
	& {V_{S_2}^{k-1}} & {\Delta^{k-1}}
	\arrow[hook, from=1-1, to=2-1]
	\arrow[hook, from=1-1, to=1-2]
	\arrow[from=1-2, to=2-2]
	\arrow[from=2-1, to=2-2]
	\arrow[from=2-2, to=2-3]
	\arrow[hook, from=3-2, to=2-2]
	\arrow[hook, from=3-2, to=3-3]
	\arrow[from=3-3, to=2-3]
	\arrow["\lrcorner"{anchor=center, pos=0.125, rotate=180}, draw=none, from=2-2, to=1-1]
	\arrow["\lrcorner"{anchor=center, pos=0.125, rotate=-90}, draw=none, from=2-3, to=3-2]
	\arrow[from=2-6, to=2-7]
	\arrow[from=1-6, to=2-6]
	\arrow["{[\dots]}"{description}, from=2-3, to=2-5]
	\arrow[hook, from=1-5, to=2-5]
	\arrow[from=2-5, to=2-6]
	\arrow[hook, from=1-5, to=1-6]
	\arrow["\lrcorner"{anchor=center, pos=0.125, rotate=180}, draw=none, from=2-6, to=1-5]
\end{tikzcd}\]

Applying the $\hom(\_,\huaX)$ functor (for $\huaX$ a Lie $\infty$-groupoid) to the expression above gives a pullback of $\huaX_{k-1}$'s over the spaces $V^{k-1}_{S_{i}}(\huaX)$ for $1\le i \le s-1$, where all maps involved in each pullback are surjective submersions. Explicitly:
\begin{align*}
    &\Lambda^k_{[k]\backslash A}(\huaX) = V^k_A(\huaX) \\
    &\cong \left(\dots\left(\left(\huaX_{k-1}\times_{V^{k-1}_{S_1}(\huaX)}\huaX_{k-1}\right)\times_{V^{k-1}_{S_2}(\huaX)}\huaX_{k-1}\right)\dots\right)\times_{V^{k-1}_{S_{s-1}(\huaX)}}\huaX_{k-1}.
\end{align*}

\begin{theorem}\label{thm:dim-gen-horns}
Let $\huaX$ be a Lie $\infty$-groupoid. Let $A$ be a subset of $[k]$ with cardinality $|A|=s$. The space $V_A^k(\huaX) = \Lambda^k_{[k] \backslash A}$ has dimension 
\begin{equation}
    \dim V^k_A(\huaX) = \sum_{i=1}^{s}(-1)^{i+1}\binom{s}{i}\dim \huaX_{k-i}.
\end{equation}
\end{theorem}

\begin{proof}
    We proceed by induction on $|A|=s$. Start by ordering $A$ as before: $A=\{j_1,\dots,j_s\}$ with $0\le j_1<j_2<\dots<j_s\le k$. Then, by the above discussion, we have that 
    \begin{align*}
        V_A^k(\huaX) \cong V_{A\backslash \{j_s\}}^k(\huaX) \times_{V^{k-1}_{S_{s-1}}(\huaX)} \huaX_{k-1},
    \end{align*}
    where $S_{s-1}=\delta_{j_s}^{-1}(\{j_1,\dots,j_s\})$ has cardinality $|S_{s-1}| = |A\backslash \{j_s\}| = s-1$, and the fiber product is over surjective submersions.
    Thus, by lemma \ref{lem:dim_fib_prod}
    \begin{align*}
        \dim V_A^k(\huaX) = \dim V_{A\backslash \{j_s\}}^k(\huaX) + \dim \huaX_{k-1} - \dim V^{k-1}_{S_{s-1}}(\huaX),
    \end{align*}
    and, by inductive hypothesis,
    \begin{align*}
        \dim V_A^k(\huaX) =& \sum_{i=1}^{s-1}(-1)^{i+1}\binom{s-1}{i}\dim \huaX_{k-i} 
        + \dim \huaX_{k-1} \\
        &\quad- \sum_{i=1}^{s-1}(-1)^{i+1}\binom{s-1}{i}\dim \huaX_{k-1-i}.
    \end{align*}
    Rewrite the last term by setting $l=i+1$ as follows:
    \begin{align*}
        - \sum_{i=1}^{s-1}(-1)^{i+1}\binom{s-1}{i}\dim \huaX_{k-1-i}
        &= - \sum_{l=2}^{s}(-1)^{l}\binom{s-1}{l-1}\dim \huaX_{k-l}\\
        &=\sum_{l=2}^{s}(-1)^{l+1}\binom{s-1}{l-1}\dim \huaX_{k-l}.
    \end{align*}
    Then
    \begin{align*}
        \dim V_A^k(\huaX) =& \binom{s-1}{1} \dim \huaX_{k-i} + \sum_{i=2}^{s-1}(-1)^{i+1}\binom{s-1}{i}\dim \huaX_{k-i} \\
        &\quad + \dim \huaX_{k-1} 
        + \sum_{l=2}^{s-1}(-1)^{l+1}\binom{s-1}{l-1}\dim \huaX_{k-l}\\
        &\quad + (-1)^{s+1}\binom{s-1}{s-1}\dim \huaX_{k-s}\\
        &= s\dim \huaX_{k-1} + (-1)^{s+1}\binom{s}{s}\dim \huaX_{k-s}\\
        &\quad + \sum_{i=2}^{s-1}(-1)^{i+1}\left(\binom{s-1}{i} + \binom{s-1}{i-1}\right)\dim \huaX_{k-i} \\
        &= \sum_{i=1}^{s}(-1)^{i+1}\binom{s}{i}\dim \huaX_{k-i}.
    \end{align*}
    Where, in the last step, we used Pascal's rule 
    \begin{align*}
        \binom{s}{i} = \binom{s-1}{i} + \binom{s-1}{i-1}
    \end{align*}
\end{proof}

\begin{corollary}[Dimension of horn spaces]\label{cor:dim-horn-spaces}
Let $\huaX$ be a Lie $\infty$-groupoid. The dimension of the horn space $\Lambda^k_j(\huaX)$ is 
\begin{equation}
\begin{split}
    \dim\Lambda^k_j(\huaX) &= \sum_{i=1}^{k-1}(-1)^{i+1}\binom{k-1}{i}\dim \huaX_{k-i}\\
    &= \binom{k}{1}\dim \huaX_{k-1}-\binom{k}{2}\dim \huaX_{k-2}\\ 
    &\quad + \binom{k}{3} \dim \huaX_{k-3} - \dots +(-1)^{k-1}\binom{k}{k} \dim \huaX_0.
\end{split}
\end{equation}
\end{corollary}

\subsection{Weak equivalences and fibrations}\label{sec:weak-eq}

To define weak equivalences and fibrations for groupoids we need to use full-generality lifting problems of the form \eqref{diag:lifting-problem}.
In the proof of Theorem \ref{thm:hom-ngpd}, we recalled a definition of fibrations in terms of a right lifting property. We now give a more explicit definition for a general category with a pretopology $(\Cat, \covers)$, which is nevertheless equivalent to the previous one for $\SSet$. 

General lifting problems define spaces of simplicial diagrams along a map $f:\huaX \to \huaY$, which are also known as \textit{objects of commutative squares} in \cite{Henriques2008}, \cite{Zhu2009}. These can also be seen as objects of simplicial diagrams in the arrow category \cite[\S 2.3.1]{Li2014}.

\begin{definition}
Let $u: \huaS \to \huaT$ be a simplicial map between simplicial sets $\huaS,\huaT$ and $f:\huaX \to \huaY$ be a simplicial map between simplicial manifolds $\huaX,\huaY$. The \textbf{object of commutative squares} $\hom(\huaS \overset{u}{\to} \huaT, \huaX \overset{f}{\to} \huaY)$ is the sheaf on $\Mfd$ whose value on $U \in \Mfd$ is the set of pairs of simplicial maps $(a,b)$ that make the following diagram commute
\[\begin{tikzcd}[ampersand replacement=\&,cramped]
{U \otimes \huaS} \& \huaX \\
{U \otimes \huaT} \& \huaY \\
\arrow["a", from=1-1, to=1-2]
\arrow["{U \otimes u}"', from=1-1, to=2-1]
\arrow["f", from=1-2, to=2-2]
\arrow["b"', from=2-1, to=2-2]
\end{tikzcd}\]
where $\otimes$ is the copowering of $\Simp\Mfd$ on $\SSet$, and $U$ is identified with its identity groupoid. 
\end{definition}

\begin{remark}
Objects of commutative squares are limits: they can be written as fiber products
\begin{equation*}
\begin{split}
    \hom(\huaS \overset{u}{\to} \huaT, \huaX \overset{f}{\to} \huaY) \cong \hom(\huaS, \huaX) \times_{\hom(\huaS, \huaY)} \hom(\huaT, \huaY),
\end{split}
\end{equation*}
where each of the factors is a limit as in Proposition \ref{prop:simp-diag-sheaf-as-limit}.\footnote{This was stated for manifolds and surjective submersions but the generalization to $(\Cat, \covers)$ is immediate.}
As we will later discuss, for some choices of $S$ and $T$, $\hom(\huaS \overset{u}{\to} \huaT, \huaX \overset{f}{\to} \huaY)$ might be representable even when some of the factors in its fiber product description are a priori not representable. This justifies writing it as an object of commutative diagrams, instead of just writing it as a fiber product. 
The fiber product description is nevertheless useful, because if the object of commutative diagrams $\hom(\huaS \overset{u}{\to} \huaT, \huaX \overset{f}{\to} \huaY)$ is representable, then its underlying set is its evaluation at $pt$, i.e. the fiber product of the underlying sets of its factors, as in Remark \ref{rem:underlying-set-of-hom-sheaf}.
\end{remark}

The lifting problem \eqref{diag:lifting-problem} with $\huaS \to \huaT$ the horn inclusion $\Lambda^m_j \into \Delta^m$ is the problem of existence of relative horn fillers --- that is, horn fillers in $\huaX$ with a prescribed image through $f$ in $\huaY$ --- and it leads to the notion of Kan fibrations. These are useful to characterize simplicial vector bundles, as we discuss in Proposition \ref{prop:SVB-is-Kan-fibration}. 

\begin{definition}\label{def:Kan-fibration}
    Let $f: \huaX \to \huaY$ be a simplicial map between simplicial objects in $(\Cat, \covers)$. We say $f$ satisfies the Kan condition $\Kan(m,j)$ if the object of commutative squares $\hom(\Lambda^m_j {\to} \Delta^m, \huaX \overset{f}{\to} \huaY)$ is representable and the canonical map 
    \begin{equation*}
        \hom(\Delta^m, \huaX) = \huaX_m \xrightarrow{(p^m_j, f)}  \hom(\Lambda^m_j {\to} \Delta^m, \huaX \overset{f}{\to} \huaY)
    \end{equation*}
    is a cover. If this map is additionally an isomorphism, we say that $f$ satisfies the strict Kan condition $\Kan!(m,j)$. 
    
    We say $f$ is a \textbf{Kan fibration} if it satisfies $\Kan(m,j)$ for all $m \ge 1$, $0 \le j \le m$. Additionally we say it is \textbf{$n$-strict} if it satisfies $\Kan!(m, j)$ for all $m \ge n + 1$. We say a Kan fibration has \textbf{order} $n$ if it is $n$-strict but not $(n-1)$-strict. 
\end{definition}

\begin{remark}
    The Kan condition $\Kan(m,j)$ for $f$ to be a Kan fibration of simplicial sets corresponds exactly to the fact that the lifting problem 
    \[\begin{tikzcd}[ampersand replacement=\&]
	\Lambda^m_j \& \huaX \\
	\Delta^m \& \huaY
	\arrow[from=1-1, to=1-2]
	\arrow["i^m_j"', hook, from=1-1, to=2-1]
	\arrow["f", from=1-2, to=2-2]
	\arrow[dashed, from=2-1, to=1-2]
	\arrow[from=2-1, to=2-2]
    \end{tikzcd}\]
    always admits a lift $\Delta^m \to \huaX$ for any pair of horizontal maps. 
    As we saw in Section \ref{sec:simp-diag-sets}, the case where $\huaY = pt$ (the identity groupoid of the point) gives the Kan conditions for $\huaX$. This is a general fact: a simplicial object $\huaX$ is an $n$-groupoid if and only if the canonical map $\huaX \to pt$ is an $n$-strict Kan fibration.
\end{remark}

\begin{remark}
    Similarly to the case of the horn space of an $\infty$-groupoid object in $(\Cat, \covers)$, which is inductively representable, a general version of Lemma \ref{lem:horn-representability} shows that if $f$ satisfies $\Kan(m,j)$ for all $m < l$, then $\hom(\Lambda^l_j {\to} \Delta^l, \huaX \overset{f}{\to} \huaY)$ is representable. We refer to \cite[Lemma 2.4]{Henriques2008}, \cite[Lemma 3.29]{Li2014} or \cite[Lemma 3.9]{RogersZhu2020} for a proof. For more results on Kan fibrations of simplicial objects and $n$-groupoid objects see for example \cite[\S 3.4]{Li2014}. 
\end{remark}

\begin{remark}\label{rem:order-of-kan-fibrations}
    As remarked in \cite[Def. 2.4]{HoyoTrentinaglia2024}, the order of a composition of Kan fibrations is the maximum between the orders of the two fibrations. Therefore, if $\huaY$ is an $n$-groupoid object of order $n$ and $f:\huaX \to \huaY$ is a Kan fibration of order $m$, then $\huaX$ is a $k$-groupoid of order $k = \max\{n,m\}$.
\end{remark}

For a map to be a weak equivalence of simplicial objects in a concrete category $(\Cat, \covers)$ it is not enough to require that the underlying map of simplicial sets is a weak equivalence. This has to do with the fact that the forgetful functor is generally not \textit{conservative}: the fact that the underlying map of $f\in \Mor(\Cat)$ is a bijection does not always imply that $f$ is an isomorphism in $\Cat$. 
This is the case in $\Mfd$. On the other hand, because the forgetful functor of $\Vect$ is monadic, it is also conservative (see e.g. \cite[Ch. 5]{Riehl2016}). In fact, as we see in Section \ref{sec:simp-vect}, weak equivalences of $\Simp\Vect$ are indeed defined as weak equivalences of the underlying simplicial sets. 
The way to resolve this problem is to choose a collection of \textit{point functors} that is large enough to detect isomorphims. We briefly recall the exposition in \cite[\S 4]{RogersZhu2020} and subsequent sections. As the authors point out, this notion originated in topos theory (see e.g. \cite[C.2.2]{Johnstone2002}, \cite[\S VII.5]{MacLaneMoerdijk1994}), and was particularly used for homotopy theory of simplicial sheaves over smooth manifolds in e.g. \cite{Dugger1998} and \cite{NikolausSchreiberStevenson2015}.

\begin{definition}[{\cite[Def. 4.1]{RogersZhu2020}}]\label{def:collection-of-points}
    Let $(\Cat, \covers)$ be a category equipped with a pretopology. A \textbf{point} of $(\Cat,\covers)$ is a functor 
    \begin{equation*}
        \mathsf{p}: \Sh(\Cat) \to \Set
    \end{equation*}
    which preserves finite limits and small colimits. 

    A collection of points $\mathsf{P}$ of $(\Cat, \covers)$ is \textbf{jointly conservative} if a morphsims of sheaves $f: F \to G$ in $\Sh(\Cat)$ is an isomorphism if and only if for all $\mathsf{p}\in \mathsf{P}$, 
    \begin{equation*}
        \mathsf{p}f: \mathsf{p}F \to \mathsf{p}G
    \end{equation*}
    is an isomorphism of $\Set$. 
    The set $\mathsf{p}F$ is also called the \textbf{$\mathsf{p}$-stalk} of the sheaf $F$. 

    If the category $\Cat$ is concrete with respect to $\Forget: \Cat \to \Set$, we say that $\mathsf{P}$ \textbf{contains the forgetful functor} if there exists a $\mathsf{p}_0 \in \mathsf{P}$ such that $\mathsf{p}_0\Yoneda = \Forget$, with $\Yoneda$ the Yoneda embedding. 

    For any simplicial object $\huaX$ in $(\Cat, \covers)$, we denote by $\mathsf{p}\huaX$ the simplicial object $\mathsf{p}\Yoneda\huaX$, where here $\Yoneda$ denotes the Yoneda embedding to simplicial sheaves over $\Cat$ defined as for manifolds in \eqref{eq:Yoneda-simp-sheaves}.
\end{definition}

\begin{example}\label{ex:points-for-mfds}
    A jointly conservative collection of points for $\Set$ or for $\Vect$ is simply given by evaluation on the point $pt$ or $\R$, respectively. Hence, by Remark \ref{rem:repbly-concrete-cat-limits}, such a collection contains the forgetful functor. 
    A jointly conservative collection of points for smooth manifolds was described in \cite[Def. 3.4.6]{Dugger1998}: Denote by $B^n_k$ the $n$-ball of radius $\frac{1}{k}$ centered at the origin in $\R^n$. Then we define the point $\mathsf{p}_n$ for any $n\ge 0$ as 
    \begin{equation*}
        \mathsf{p}_n(F) = \colim\limits_{k \to \infty} F(B^n_k), 
    \end{equation*}
    for any $F \in \Simp\Sh(\Mfd)$. 
    The collection $\mathsf{P}=\{\mathsf{p}_n\}_{n\ge 0}$ is a collection of points, where $\mathsf{p}_0\Yoneda$ is the forgetful functor. 
    An extension of this collection for the category of Banach manifolds was defined in \cite[\S 4.1]{RogersZhu2020}, where it was shown to be jointly conservative. Since this restricts to $\mathsf{P}$ over $\Mfd$, this is also jointly conservative. 
\end{example}

We summarize some results from \cite[\S 4.2]{RogersZhu2020} that we find helpful to bridge the gap between the theory of groupoids in $\Set$ and in $(\Cat, \covers)$, by understanding how the latter is generalizing the first one stalkwise. The definition of weak equivalences of simplicial objects in $(\Cat, \covers)$ is also a step in this direction.

\begin{remark}
    Let $(\Cat, \covers, \mathsf{P})$ be a category with a pretopology and a jointly conservative collection of points. By applying $\mathsf{p}\Yoneda$ to a cover $f \in \covers$, the obtained map $\mathsf{p}f_*$ is a surjection for any $\mathsf{p} \in \mathsf{P}$. In other words, any cover in $\covers$ is a \textit{stalkwise surjection} for any collection of points $\mathsf{P}$. 
    As pointed out in \cite[Rem. 6.10]{RogersZhu2020} for the case of surjective submersions in $\Mfd$, stalkwise surjections are still a weaker notion than covers in general. 
\end{remark}

\begin{remark}
    Let $(\Cat, \covers, \mathsf{P})$ be a category with a pretopology and a jointly conservative collection of points.
    In a similar way as for covers, if $f:\huaX \to \huaY$ is a Kan fibration of simplicial objects in $(\Cat, \covers)$ then for any point functor $\mathsf{p}\in \mathsf{P}$, $\mathsf{p}f$ is a Kan fibration of simplicial sets. 
    Therefore, a Kan fibration in $(\Cat, \covers)$ is in a sense a \textit{stalkwise Kan fibration} for any collection of points $\mathsf{P}$ on $(\Cat, \covers)$. 
    In particular, if $\huaX$ is an $n$-groupoid object in $(\Cat, \covers)$, $\mathsf{p}X$ is an $n$-groupoid in $\Set$ for any point functor $\mathsf{p}$.
    It follows from the same observation as in the previous remark that generally the notion of Kan fibration in $(\Cat, \covers)$ is stronger than the notion of stalkwise Kan fibration.
\end{remark}

\begin{definition}[{\cite[5.1]{RogersZhu2020}}]
    Let $(\Cat, \covers, \mathsf{P})$ be a category with a pretopology and a jointly conservative collection of points. We say a morphism $f:\huaX \to \huaY$ of simplicial objects in $\Cat$ is a \textbf{stalkwise weak equivalence} if, for any $\mathsf{p} \in \mathsf{P}$, 
    \begin{equation*}
        \mathsf{p}f: \mathsf{p}\huaX \to \mathsf{p}\huaY,
    \end{equation*}
    is a weak equivalence of simplicial sets. 
\end{definition}

\begin{remark}
    In \cite[\S 5.1]{RogersZhu2020}, the authors prove that this definition is compatible with the construction of simplicial homotopy group sheaves of $n$-groups appearing in \cite{Henriques2008}, in the sense that a morphism of $n$-groups that induces an isomorphism between the simplicial homotopy groups is a stalkwise weak equivalence, just as for $\Simp\Set$. 
\end{remark}

In \cite{RogersZhu2020}, the reason for assigning the stalkwise weak equivalences to be \textit{the} weak equivalences is that of showing that with this assignment, the category of $\infty$-groupoids in $(\Cat, \covers)$ is an incomplete category of fibrant objects (iCFO) if $\covers$ is locally stalkwise (\cite[Def. 6.5]{RogersZhu2020}) with respect to the collection of points $\mathsf{P}$. The details of this fact are outside the scope of this thesis, so we refer to their paper for more information. The main points we are interested in are that: First of all, this is the case for $\Mfd$ with surjective submersions and the collection of points defined above. Secondly, we need the following characterization of stalkwise weak equivalences in $(\Mfd, \covers_{s.sub.}, \mathsf{P})$. A version of this result for CFOs of higher geometric groupoids in a descent category appeared previously in \cite[Thm. 5.1]{BehrendGetzler2017}.

\begin{proposition}[{\cite[Cor. 7.16]{RogersZhu2020}}]\label{prop:def-we}
    A morphism $f: \huaX \to \huaY$ of Lie $\infty$-groupoids is a stalkwise weak equivalence if and only if
    $\hom(\partial\Delta^l \overset{(\delta_i\delta_l)_i}{\to} \Lambda^{l+1}_{l+1}, \huaX \overset{f}{\to} \huaY)$ is representable for all $l \ge 0$ and the map 
    \begin{equation}\label{eq:def-we}
        \begin{split}
            r^f_l:=((&d_0,\dots, d_l), p_{l+1}^{l+1}):\\
            &\hom(\Delta^l \overset{\delta_{l+1}}{\to} \Delta^{l+1}, \huaX \overset{f}{\to} \huaY)
            \to \hom(\partial\Delta^l \overset{(\delta_i\delta_l)_i}{\to} \Lambda^{l+1}_{l+1}, \huaX \overset{f}{\to} \huaY)\\
            &\qquad \cong \huaX_l \times_{f, \huaY_l, d_{l+1}} \huaY_{l+1} \to 
            \hom(\partial\Delta^l, \huaX) \times_{\hom(\partial\Delta^l, \huaY)} \Lambda^{l+1}_{l+1}(\huaY)
        \end{split}
        \end{equation}
    is a cover for all $l \ge 0$.
\end{proposition}

We refer to \cite{RogersZhu2020} for a proof of this fact. While the domain of each map $r_l^f$ is clearly representable, since $d_{l+1}$ is a cover, representability of the codomain can be shown inductively as in the case of horn spaces. This was shown in \cite[Lemma 6.28]{Li2014}, by using joins. Despite the fact that this lemma was stated there for 2-groupoid objects, the proof is immediately valid for all $n$-groupoid objects for any $n$.

\begin{lemma}[{\cite[Lemma 6.28]{Li2014}}] Let $\huaX, \huaY$ be Lie $\infty$-groupoids. 
If $\hom(\partial\Delta^l \to \Lambda^{l+1}_{l+1}, \huaX \to \huaY)$ is representable, and $r^f_l$ is a cover for all $l < m$, then $\hom(\partial\Delta^m \to \Lambda^{m+1}_{m+1}, \huaX \to \huaY)$ is representable. 
\end{lemma}

The following lemma can be shown by generalizing \cite[Lemma 6.29]{Li2014} from 2-groupoids to $n$-groupoids and using \cite[Lemma 2.5]{Zhu2009}. 

\begin{lemma}{\label{lem:we-n-gpd-n-hypercover-condition}}
Let $f:\huaX \to \huaY$ be a map of Lie $n$-groupoids such that $r^f_l$ is a cover for all $0\le l \le n-1$. If $r^f_n$ is an isomorphism, then 
\begin{equation*}
    \huaX_n \cong \hom(\partial\Delta^n \to \Delta^n, \huaX \to \huaY),
\end{equation*}
and $r^f_l$ is an isomorphism for all $l \ge n+1$. In particular $f$ is a stalkwise weak equivalence.
\end{lemma}

\begin{proof}
By using the cancellation lemma for pullbacks, since $\Lambda^{n+1}_{n+1}(\huaY) \cong \huaY_{n+1}$, we have the pullback square
\begin{equation}\label{eq:we-n-gpd-n-hypercover-condition}
\begin{tikzcd}
	{\huaX_n\cong \huaX_n\times_{\huaY_n}\huaY_n} & {\hom(\partial\Delta^n \to \Delta^n, \huaX \to \huaY)} \\
	{\huaX_n\times_{\huaY_n}\huaY_{n+1}} & {\hom(\partial\Delta^n\to \Lambda^{n+1}_{n+1}, \huaX \to \huaY)}
	\arrow[from=1-1, to=1-2]
	\arrow["{(id, s_n)}"', from=1-1, to=2-1]
	\arrow["{(id, p^{n+1}_{n+1}s_n)}", from=1-2, to=2-2]
	\arrow[from=2-1, to=2-2]
\end{tikzcd}
\end{equation}
Therefore, if the map on the bottom is an isomorphism, the map on the top is an isomorphism as well. Now, despite the fact that $f$ is not a \textit{hypercover} (see \cite[Def. 2.3]{Zhu2009}) the inductive proof of \cite[Lemma 2.5]{Zhu2009} can still be used to show that 
\begin{equation*}
    \huaX_l \cong \hom(\partial\Delta^l \to \Delta^l, \huaX \to \huaY),
\end{equation*}
for all $l>n$.
The pullback diagram \eqref{eq:we-n-gpd-n-hypercover-condition} with $l> n$ instead of $n$ still applies. Therefore  $r^f_l$ is an isomorphism for all $l \ge n+1$.
\end{proof}

If $\huaX$ and $\huaY$ are Lie $n$-groupoids, and $f:\huaX \to \huaY$ is a stalkwise weak equivalence between them, then $r^f_n$ is an isomorphism and not just a cover. We could not find an explicit proof of this known fact in the literature, but morally, this has to do with the fact that $n$-groupoids only contain homotopy information up to level $n$. For example, for $n=1$, (smooth) equivalences of categories between Lie groupoids were defined to be the weak equivalences in \cite{MoerdijkMrcun2003} (such maps are called Morita maps in \cite{Hoyo2012,HoyoOrtiz2020}). Adapted to $(\Cat, \covers)$, these are maps (functors) $f:\huaX \to \huaY$ that are:
\begin{enumerate}
    \item Essentially surjective, i.e. 
    \begin{equation*}
        r^f_1: \huaX_0 \times_{\huaY_0}\huaY_1 \to \huaY_0
    \end{equation*}
    is a cover.
    \item Fully faithful, i.e. 
    \[\begin{tikzcd}
        {\huaX_1} & {\huaY_1} \\
        {\huaX_0 \times \huaX_0} & {\huaY_0\times \huaY_0}
        \arrow[from=1-1, to=1-2]
        \arrow[from=1-1, to=2-1]
        \arrow["{(d_0,d_1)}", from=1-2, to=2-2]
        \arrow["{f_0\times f_0}"', from=2-1, to=2-2]
    \end{tikzcd}\]
    is a pullback diagram, or in other words $r^f_1$ is an isomorphism. 
\end{enumerate}
By Lemma \ref{lem:we-n-gpd-n-hypercover-condition}, $r^f_l$ is then an isomorphism for all $l \ge 1$ and in particular $f$ is a stalkwise weak equivalence of the nerves. 

Clearly, it is easier to check the finite set of conditions in Lemma \ref{lem:we-n-gpd-n-hypercover-condition} than the infinite number of them in \ref{prop:def-we}. As such, we summarize this discussion by making the following working definition. 

\begin{definition}
Let $\huaX$ and $\huaY$ be Lie $n$-groupoids.
We say a simplicial map $f: \huaX \to \huaY$ is a \textbf{weak equivalence} (or \textbf{Morita equivalence}) if the map 
    \begin{equation}
            r^f_l:
            \huaX_l \times_{f, \huaY_l, d_l} \huaY_{l+1} \to 
            \hom(\partial\Delta^l, \huaX) \times_{\hom(\partial\Delta^l, \huaY)} \Lambda^{l+1}_{l+1}(\huaY)
        \end{equation}
    is a cover for $0\leq l<n$ and an isomorphism for $l=n$.
\end{definition}
    
\begin{remark}\label{rem:w-eq-and-morita}
    As expected, the class of weak equivalences generates an equivalence relation. We say that two Lie $n$-groupoids $\huaX$ and $\huaY$ are \textbf{Morita equivalent} if they are connected by a span of Morita equivalences $\huaX \xleftarrow{\sim} \huaZ \xrightarrow{\sim} \huaY$. 
    By the factorization lemma \cite[prop. 2.5]{RogersZhu2020} in the iCFO of Lie $\infty$-groupoids, this definition is equivalent to the more common one in \cite{Zhu2009} which defines a Morita equivalence as a span of hypercovers. 
    In practice, all Morita equivalences appearing in this thesis are given by a single weak equivalence instead of a more general span.
\end{remark}

\subsection{Coskeleta and Finite Data}\label{sec:finite-data}

Because of the Kan conditions, $n$-groupoid objects in any category with a pretopology $(\Cat, \covers)$ can be reduced from an infinite simplicial object to a finite set of data. For a Lie 1-groupoid as in Example \ref{ex:1gpd}, this can be thought of as an inverse of the nerve construction. 
In this section we review the coskeleton functor and how to use its properties to obtain the finite data of an $n$-groupoid object. 
We focus particularly on 2-groupoid objects and introduce a helpful notation for their higher multiplications.
Our main references for this section are \cite{Duskin2001/02}, \cite{Zhu2009}. See also \cite[\href{https://stacks.math.columbia.edu/tag/0AMA}{Tag 0AMA}]{stacks-project} and \cite[\href{https://kerodon.net/tag/0513}{Tag 0513}]{kerodon}.
Throughout this section we consider simplicial objects in a concrete category $(\Cat, \covers)$ with a locally stalkwise pretopology with respect to a collection of points $\mathsf{P}$ which contains the forgetful functor (see Definition \ref{def:collection-of-points}). The main example is the category of smooth manifolds with surjective submersions as covers and the collection of points described in Example \ref{ex:points-for-mfds}.

Consider the \textbf{$m$-truncation} of a simplicial object $\huaX$, which is obtained by ``forgetting'' all simplicial levels higher than $\huaX_m$. This consists of $\huaX_0, \dots, \huaX_m$ and face and degeneracy maps between them. We denote it by $\tr^m\huaX$. In other words a \textbf{truncated simplicial object} is a functor $\Delta_{\le m}^{op} \to \Cat$, where $\Delta_{\le m}$ is the full subcategory of $\Delta$ generated by the elements $[0], \dots, [m]$. 
This defines a functor $\tr^m$ from the category $\Simp\Cat$ of simplicial objects in $\Cat$ to the category $\Simp\Cat_{\le m}$ of $m$-truncated simplicial objects in $\Cat$.

The truncation $tr^m$ admits a left adjoint $\sk^m$ called the \textbf{skeleton} functor:
\begin{equation*}
    \Simp\Cat(sk^m\widebar{\huaX}, \huaY) \cong \Simp\Cat_{\le m}(\widebar{\huaX}, tr^m\huaY), \qquad \forall \widebar{\huaX} \in \Simp\Cat_{\le m},\, \forall \huaY \in \Simp\Cat. 
\end{equation*}
The skeleton $sk^m\widebar{\huaX}$ is then easily described by being the simplicial object consisting only of higher degeneracies of simplices in $\widebar{\huaX}$ in all levels higher than $m$. When composed with the truncation functor, it gives an endofunctor $\Sk^m:= \sk^m\tr^m: \Simp\Cat \to \Simp\Cat$. In other words, the \textbf{skeleton} of a simplicial object $\huaX$ is its simplicial subobject generated by all simplices of dimension lower than $m$ and their higher degeneracies. It is common practice to identify each $m$-truncated simplicial object with its $m$-skeleton so that one can get away with considering only simplicial objects instead of $m$-truncated ones. In this picture, any $m$-truncated simplicial object corresponds to an $m$-skeletal simplicial object, i.e. one that is isomorphic to its $m$-skeleton.

\begin{definition}\label{def:coskeleton-general-site}
If it exists, the \textbf{$m$-coskeleton} $\Cosk^m\huaX$ of a simplicial object $\huaX \in \Simp\Cat$ is the representing object of the presheaf
\begin{equation*}
    \huaU \mapsto \Simp\Cat(\Sk^m\huaU, \huaX), \quad \forall \huaU \in \Simp\Cat. 
\end{equation*}
\end{definition}

Clearly, by this definition, if the $m$-coskeleton is defined for all simplicial objects in $\Cat$, then it is a functor (called \textbf{coskeleton functor}), which is the right adjoint of $\Sk^m$. 

\begin{remark}
A coskeleton functor can also be defined as $\cosk^m: \Simp\Cat_{\le m} \to \Simp\Cat$. This is the right adjoint of $tr^m$ as follows.
By the identification of $m$-truncated and $m$-skeletal simplicial objects, if the coskeleton of $\huaX = \Sk^m\huaX$ exists, then the representing object of the presheaf 
\begin{equation*}
    \huaU \mapsto \Simp\Cat_{\le m}(\tr^m\huaU, \tr^m\huaX)
\end{equation*}
is $\cosk^m(\tr^m(\huaX)) = \Cosk^m(\Sk^m(\huaX))$. Now if all coskeleta exist, there is an adjunction $\sk^m \dashv \tr^m \dashv \cosk^m$, between $\Simp\Cat$ and $\Simp\Cat_{\le m}$.
\end{remark}

\begin{lemma}[{\cite[\href{https://stacks.math.columbia.edu/tag/0183}{Tag 0183}]{stacks-project}}]\label{lem:cosk-as-limit}
If the category $\Cat$ has finite limits, then $\cosk^m$ and $\Cosk^m$ functors exist for all $m$. Moreover, each level of the coskeleton of any ($m$-truncated) simplicial object can be written as a limit in $\Cat$. This can be written recursively as 
\begin{equation}\label{eq:cosk-as-limit}
\Cosk^m(\huaX)_l := \begin{cases}
    \huaX_l &\text{for } l \le m,\\
    \hom(\partial\Delta^{m+1},\huaX) &\text{for } l = m+1,\\
    \hom(\partial\Delta^l, \Cosk^m(\huaX)) &\text{for } l > m+1.
\end{cases}
\end{equation}
In particular, it depends only on the $m$-skeleton of $\huaX$ (or equivalently on its $m$-truncation). 
\end{lemma}

For this proof and more details on the limit description of coskeleta we refer to \cite[\href{https://stacks.math.columbia.edu/tag/0AMA}{Tag 0AMA}]{stacks-project}. Note that for $l> m+1$, $\Cosk^m(\huaX)_l$ is an iterated boundary space. In \cite{Duskin2001/02} this is also called iterated simplicial kernel. 
We give a proof of the same result for $\Cat = \Set$ in the next lemma \ref{lem:m-coskeleton-sset}. Because $\Set$ is complete, all coskeleta of simplicial sets exist and $\Cosk^m$ defines a functor. This is also true for $\Vect$. 

In the category of simplicial manifolds, not all coskeleta exist a priori because the iterated boundary spaces in \eqref{eq:cosk-as-limit} are not always representable in $\Mfd$, as we observed in Remark \ref{rem:bdry-representability}. 
In this case, and in any other non-complete categories, we can consider $\Cosk^m(\huaX)$ to be a priori defined as a simplicial sheaf. That is, we see it as an object in $\Simp\Sh(\Cat)$. If this sheaf is representable, we identify it with its representing object and denote both by $\Cosk^m(\huaX)$.

 By definition of a collection of points $\mathsf{P}$ we have that for any $\mathsf{p} \in\mathsf{P}$,
\begin{equation*}
    \mathsf{p}\Cosk^m(\huaX) \cong \Cosk^m(\mathsf{p}\huaX),
\end{equation*}
where we recall that $\mathsf{p}F$ of a simplicial sheaf $F \in \Simp\Sh(\Cat)$ is the simplicial set defined by composing $F$ and $\mathsf{p}$ as functors, while for a simplicial object $\huaX$, $\mathsf{p}\huaX = \mathsf{p}\Yoneda\huaX$. 
If the $m$-coskeleton of a simplicial object $\huaX$ in $\Cat$ exists, and $\Cat$ is a concrete category where the forgetful functor preserves limits (such as $\Mfd$, by Remark \ref{rem:repbly-concrete-cat-limits}), then its underlying simplicial set is $\Forget(\Cosk^m(\huaX)) \cong \Cosk^m \Forget(\huaX)$, which is described as in the following lemma.

\begin{lemma}\label{lem:m-coskeleton-sset}
Let $\huaX$ be a simplicial set. Then its $m$-coskeleton $\Cosk^m\huaX$ is given at each level by 
\begin{equation*}
\cosk^m(\huaX)_l := \begin{cases}
    \huaX_l &\text{for } l \le m,\\
    \partial^{m+1}(\huaX) &\text{for } l = m+1,\\
    \partial^{l}(\Cosk^m(\huaX)) &\text{for } l > m+1,
\end{cases}
\end{equation*}
where $\partial_l$ denotes the $l$-th boundary space as in \eqref{eq:DefHornBdrySpaces}. 
The face maps of the $m$-coskeleton above level $m$ are the projections 
\begin{equation*}
    d_i: (x_0, \dots, x_l) \mapsto x_i \in \Cosk^m(\huaX)_l,
\end{equation*}
while its degeneracy maps above level $m$ are defined in components $0\le i\le m$ by 
\begin{equation}\label{eq:cosk-degeneracies}
    (s_j x)_i =\begin{cases}
    s_{j-1} d_i x  & \text{if}\; i<j, \\
    x  & \text{if}\; i=j, j+1,\\
    s_j d_{i-1} x & \text{if}\; i> j+1,
    \end{cases}
\end{equation}
for any $x \in \Cosk^m(\huaX)_l$, $l > m$.
\end{lemma}

\begin{proof}
By definition of coskeleton we have that for any $l$,
\begin{equation*}
    \Cosk^m(\huaX)_l = \SSet(\Delta^l, \Cosk^m\huaX) = \SSet(\Sk^m\Delta^l, \huaX).
\end{equation*}
If $l \le m$, then $\Sk^m\Delta^l = \Delta^l$, so this is just $\huaX_l$. 
If $l \ge m+1$, then $\Sk^m\Delta^{l} = \Sk^m\partial\Delta^{l}$, so that
\begin{equation*}
    \Cosk^m(\huaX)_l 
    = \SSet(\Sk^m\partial\Delta^l, \huaX) 
    = \SSet(\partial\Delta^l, \Cosk^m\huaX). 
\end{equation*}
If $l=m+1$ we have the special case where the latter is $\SSet(\partial\Delta^{m+1}, \huaX)$. 

The simplicial maps defined in the statement are then the obvious canonical choices that satisfy the simplicial identities. 
\end{proof}

We now recall that all coskeleta of simplicial objects satisfy certain strict Kan conditions and that conversely, simplicial sets which satisfy certain strict Kan conditions are isomorphic to their $m$-coskeleton for some $m$. Such simplicial sets are called \textbf{$m$-coskeletal}. A proof of this fact can be found for example in \cite[Prop. 2.37]{Li2014}.

\begin{lemma}[{\cite[Prop. 2.37]{Li2014}}]\label{lem:prop2.37-Li}
    Let $\huaX$ be a simplicial set. Then 
    \begin{enumerate}[label=(\roman*)]
    \item If $\huaX$ is $m$-coskeletal, then $\huaX$ satisfies $\Kan!(l,j)$ for all $l \ge m + 2$ and $0\le j \le m$. 
    \item If $\huaX$ satisfies $\Kan!(l,j)$ for all $l \ge m$, then $\huaX$ is $m$-coskeletal. In particular, if $\huaX$ is an $n$-groupoid, then it is $n+1$-coskeletal.
    \end{enumerate}
\end{lemma}

We now extend this property to general simplicial objects in $(\Cat, \covers, \mathsf{P})$, where the $m$-coskeleton might not be representable a priori. 

\begin{definition}
    Let $(\Cat,\covers,\mathsf{P})$ be a category with a pretopology and a jointly conservative collection of points.
    A simplicial object $\huaX$ in $\Cat$ is \textbf{stalkwise $m$-coskeletal} if $\mathsf{p}\huaX \cong \Cosk^m(\mathsf{p}\huaX)$ for any $\mathsf{p}\in \mathsf{P}$. 
    A simplicial object $\huaX$ in $\Cat$ is \textbf{weakly stalkwise $m$-coskeletal} if $\mathsf{p}\huaX \subseteq \Cosk^m(\mathsf{p}\huaX)$ for any $\mathsf{p}\in \mathsf{P}$. 
\end{definition}

As an immediate consequence of the collection of points $\mathsf{P}$ being jointly conservative, we have the following lemma. 

\begin{lemma}\label{lem:stalkwise-cosk-is-cosk}
    A simplicial object $\huaX$ is stalkwise $m$-coskeletal if and only if $\Cosk^m(\huaX)$ exists (i.e. it is representable as a simplicial sheaf) and it is isomorphic to $\huaX$. In this case we say that $\huaX$ is \textbf{$m$-coskeletal}. 
\end{lemma}

We also note that if $\huaX$ is weakly stalkwise $m$-coskeletal in a concrete category whose collection of points contains the forgetful functor, then $\huaX$ can be seen as a subset of its $m$-coskeleton with the structure of an object in $\Cat$. For example a weakly $m$-coskeletal simplicial manifold is a simplicial subset of its $m$-coskeleton that is a smooth simplicial manifold. 

\begin{lemma}\label{lem:n-gpd-coskeletality}
    All $n$-groupoid objects in $(\Cat, \covers, \mathsf{P})$ are $(n+1)$-coskeletal and weakly stalkwise $n$-coskeletal. That is, if $\huaX$ is an $n$-groupoid object in $(\Cat, \covers, \mathsf{P})$, then 
    \begin{equation*}
        \huaX \cong \Cosk^{n+1}(\huaX), \text{ and } 
        \mathsf{p}\huaX \subseteq \Cosk^n(\mathsf{p}\huaX), \forall \mathsf{p}\in\mathsf{P}.
    \end{equation*}
\end{lemma}

\begin{proof}
    If $\huaX$ is an $n$-groupoid in $(\Cat, \covers, \mathsf{P})$, then $\mathsf{p}\huaX$ is an $n$-groupoid in $\Set$. Therefore, by \ref{lem:prop2.37-Li}, $\mathsf{p}\huaX$ is $(n+1)$-coskeletal for any $\mathsf{p}\in \mathsf{P}$. In other words, $\huaX$ is stalkwise $n$-coskeletal, which implies it is $n$-coskeletal.

    The second part follows from the description of $\Cosk^{n}(\mathsf{p}\huaX)$ in Lemma \ref{lem:m-coskeleton-sset}, and the fact that any $(n+1)$-simplex in $\mathsf{p}\huaX \cong \Lambda^{n+1}_j (\mathsf{p}\huaX)$ can be identified with the boundary of its unique horn filler by $\Kan!(n+1,j)$, which is then an element of $\partial^{n+1}(\huaX)= \Cosk^n(\huaX)_{n+1}$. The same happens at higher levels, and this inclusion is clearly compatible with the simplicial maps. 
\end{proof}

This lemma effectively reduces the data required to describe an $n$-groupoid object $\huaX$ to just its $(n+1)$-skeleton, (or truncation), since we have 
\begin{equation*}
    \huaX \cong \Cosk^{n+1}(\Sk^{n+1}(\huaX)).
\end{equation*}
Additionally, by the strict Kan conditions at level $n+1$ and the weak $n$-\hspace{0pt}coskeletality of $\huaX$, this data can be further reduced by the following construction. 

Let $\huaX$ be an $n$-groupoid object in $(\Cat, \covers, \mathsf{P})$ with $\mathsf{P}$ containing the forgetful functor. As in the proof of the previous lemma, any $(n+1)$-simplex $X$ in $\huaX$ can be identified with its boundary.
Denote the projection $(d_0, \dots, d_{n+1}):\huaX_{n+1} \to \partial^{n+1}(\Forget\huaX)$ by $\theta_{n+1}$.
By the strict Kan conditions $\Kan!(n+1, k)$, each of the faces of $X \in \huaX_{n+1}$ is actually determined by the other ones. The operations of recovering the missing face can be seen as defining $(n+1)$-ary multiplications.

\begin{definition}\label{def:multiplication-maps}
Let $\huaX$ be an $n$-groupoid object in $(\Cat, \covers)$. Then its $k$-th \textbf{multiplication map} for $0\le k \le n+1$ is the map
\begin{equation*}
    m_k =  d_k(p^{n+1}_k)^{-1}: \Lambda^{n+1}_k(\huaX) \to \huaX_n.
\end{equation*}
By writing elements of $\Lambda^{n+1}_k(\huaX)$ as $n$-tuples 
\begin{equation*}
    (x_0, x_1, \dots, \widehat{x_k}, \dots, x_{n}, x_{n+1}), 
\end{equation*}
such that $d_i x_j = d_{j-1}x_i$ for all $i<j$ with $i, j \neq k$, we see the $k$-multiplication as an $(n+1)$-ary operation between composable $(n+1)$-tuples of $x_i \in \huaX_n$ and denote it by 
\begin{equation*}
    m_k(x_0, \dots, \widehat{x_k}, \dots, x_{n+1})
    = x_0 \dots x_{k-1} \square x_{k+1} \dots x_{n+1}, 
\end{equation*}
with the empty box in the $k$-th position. 
\end{definition}

\begin{remark}\label{rem:multiplication-defined-iso}
By identifying $(n+1)$-simplices with their boundaries, we have that the inverse of the horn projection $(p^{n+1}_k)^{-1}$ can be expressed using the multiplication map $m_k$.
Define the map
\begin{equation}\label{eq:multiplication-defined-iso}
\begin{split}
    b_k: \Lambda^{n+1}_k(\huaX) &\to \partial_{n+1}(\huaX)\\
    (x_0, \dots, \widehat{x_k}, \dots, x_{n+1}) &\mapsto (x_0, \dots, x_0 \dots x_{k-1} \square x_{k+1} \dots x_{n+1}, \dots, x_{n+1}).
\end{split}
\end{equation}
The image of this is $\huaX_{n+1} \subseteq \partial_{n+1}\huaX$, and $b_k$ is an inverse of the horn projection.
Therefore, the multiplication maps define inverses for the horn projections at level $n+1$ and the maps $p^{n+1}_jb_k: \Lambda^{n+1}_k(\huaX) \to \Lambda^{n+1}_j(\huaX)$ define isomorphisms between horn spaces relative to different indices $k$ and $j$.
\end{remark}

\begin{remark}\label{rem:multiplication-comp-with-deg}
    Let $\huaX$ be an $n$-groupoid object. By $\Kan!(n+1,k)$, for any fixed $k$, we have an isomorphism of $\Sk^{n+1}(\huaX)$ and the simplicial object 
    \begin{equation}\label{eq:sk-of-n-trunc-with-k-horn}
        \Sk^{n+1}(\Lambda^{n+1}_k(\huaX) \rightfourarrows \huaX_n \rightfourarrows \dots \rightrightarrows \huaX_0),
    \end{equation}
    The face maps between $\Lambda^{n+1}_k(\huaX)$ and $\huaX_n$ for $j\neq k$ are the projections $d_j=pr_{j+1}$ and the face map $d_k$ is $m_k$. The degeneracy maps between these two levels are defined as in \eqref{eq:cosk-degeneracies}. Then \eqref{eq:sk-of-n-trunc-with-k-horn} is a simplicial object for any $k$ and only if the multiplication maps satisfy
    \begin{equation}\label{eq:n-finite-data-comp-deg-1}
    m_k(p^{n+1}_k s_k x) = x, \text{ for all } 0\le k \le n,
    \end{equation}
    or equivalently
    \begin{equation}\label{eq:n-finite-data-comp-deg-2}
    m_k(p^{n+1}_k s_{k-1} x) = x, \text{ for all } 1\le k \le n+1.
    \end{equation}
    The two sets of conditions are equivalent due to the fact that the maps $p^{n+1}_jb_k$ are isomorphisms. If these conditions are satisfied, we say the set of multiplications $\{m_k\}_{i=0}^{n+1}$ is \textbf{compatible with the degeneracy maps}. 
\end{remark}

There is one more condition we need, in order to recover the data of a Lie $n$-groupoid from its $n$-truncations and a set of multiplications satisfying the conditions in Remarks \ref{rem:multiplication-defined-iso} and \ref{rem:multiplication-comp-with-deg}. This is an associativity condition that comes from the fact that $\huaX$ satisfies the strict Kan conditions at level $n+2$.
These involve $(n+2)$-simplices of an $n$-groupoid, which can be seen as elements of $\Cosk^{n+1}(\huaX)_{n+2}$, and as such can be represented by simplicial matrices. 
We recall this notion following \cite{Glenn1982, Duskin2001/02}.
Let $\huaX$ be a simplicial set and $m \ge 0$. By Lemma \ref{lem:m-coskeleton-sset}, the elements of $\Cosk^m\huaX_{m+2}$ are $(m+3) \times (m+2)$ matrices $(x_{ji})$ with entries $x_{ji} \in \huaX_n$ such that:
\begin{enumerate}[label=(\roman*)]
    \item the index $0\le j \le m+3$ denotes the row and $0\le i \le m+2$ the column.
    \item $d_k x_{ji} = d_{i-1} x_{jk}$ for any $0 \le j \le n+2$ and $0 \le i<k \le n+1$. That is, each row is the boundary of an $(n+1)$-simplex, an element of $\partial_{m+1}(\huaX) = \Cosk^m\huaX_{m+1}$.
    \item $x_{ji}= x_{i(j-1)}$, for any $0\le i < j \le n+2$. Denoting the matrix by $M$ this is the same as the condition that $d_id_jM = d_{j-1}d_iM$ for any $0\le i < j \le n+2$.
\end{enumerate}
We call such matrices \textbf{simplicial matrices}. By writing the elements of a simplicial matrix as $x_{ji} = x_{\mu}$, where $\mu = d_id_j E_n$, its symmetries become evident, as we illustrate in the following example.

\begin{example}\label{ex:simp-matrices}
    For $m$=1, simplicial matrices are of the form 
    \begin{equation*}
    \begin{pmatrix}
            x_{23} &x_{13} &x_{12}\\
            x_{23} &x_{03} &x_{02}\\
            x_{13} &x_{03} &x_{01}\\
            x_{12} &x_{02} &x_{01}
    \end{pmatrix}
    \end{equation*}
    where all rows form boundaries of (possibly empty) triangles in $\huaX$. This encodes the boundary of a (possibly empty) tetrahedron in $\huaX$. 

    For $m=2$, they have the form
    \begin{equation*}
        \begin{pmatrix}
            x_{234} & x_{134} & x_{124} & x_{123}\\
            x_{234} & x_{034} & x_{024} & x_{023}\\
            x_{134} & x_{034} & x_{014} & x_{013}\\
            x_{124} & x_{024} & x_{014} & x_{012}\\
            x_{123} & x_{023} & x_{013} & x_{012}
        \end{pmatrix}
    \end{equation*}
    where each row is the boundary of a (possibly empty) tetrahedron in $\huaX$, and the full matrix encodes the boundary of a (possibly empty) 4-simplex in $\huaX$.
\end{example}

\begin{definition}
We say that a multiplication map $m_k:\Lambda^{n+1}_k(\huaX) \to \huaX_n$ for $0\le k \le n+1$, is \textbf{associative} if 
given $n$-simplices 
\begin{equation*}
    \begin{split}
    &x_{ji}, \text{ for all } 0\le i \le n+1, 
    \, 0\le j \le n+2 \text{ such that } j\neq k, k+1,\text{ and } i\neq k,\\
    &x_{ki}, \text{ for all } 0\le i < k, \text{ and }
    x_{(k+1)i}, \text{ for all } k < i \le n+1, 
    \end{split}
\end{equation*}
such that $d_l x_{ji} = d_{i-1} x_{jl}$ for any $0 \le j \le n+2$ and $0 \le i<k \le n+1$, we have 
\begin{equation}\label{eq:n-gpd-associativity}
\begin{split}
    &m_k(x_{k0}, \dots, x_{k(k-1)}, \widehat{x}_{kk},
    m_k(x_{(k+2)0}, \dots, \widehat{x}_{(k+2)k}, \dots, x_{(k+2)(n+1)}),\\
    &\qquad \dots,
    m_k(x_{(n+2)0}, \dots, \widehat{x}_{(n+2)k}, \dots, x_{(n+2)(n+1)})) \\
    &=m_k(m_k(x_{00}, \dots, \widehat{x}_{0k}, \dots, x_{0(n+1)}),\\ 
    &\qquad \dots, m_k(x_{(k-1)0}, \dots, \widehat{x}_{(k-1)k}, \dots, x_{(k-1)(n+1)}),\\
    &\qquad\qquad\qquad \widehat{x}_{(k+1)k}, x_{(k+1)(k+1)}, \dots, x_{(k+1)(n+1)}).
\end{split}
\end{equation}
\end{definition}

\begin{remark}\label{rem:n-gpd-assoc-simp-mat}
In terms of simplicial matrices, associativity of $m_k$ is precisely the fact that given $n$-simplices in $\huaX$ that fit as the following components of a simplicial matrix
\begin{equation*}
    \begin{pmatrix}
            x_{00} &\dots &x_{0(k-1)} &\widehat{x}_{0k} &x_{0(k+1)} &\dots &x_{0(n+1)}\\
            \dots &\dots &\dots &\dots &\dots &\dots &\dots\\
            x_{k0} &\dots &x_{k(k-1)} &\widehat{x}_{kk} &\widehat{x}_{0(k-1)} &\dots &\widehat{x}_{k(n+1)}\\
            \widehat{x}_{(k+1)0} &\dots &\widehat{x}_{(k+1)(k-1)} &\widehat{x}_{(k+1)k} &x_{(k+1)(k+1)} &\dots &x_{(k+1)(n+1)}\\
            \dots &\dots &\dots &\dots &\dots &\dots &\dots\\
            x_{(n+2)0} &\dots &x_{(n+2)(k-1)} &\widehat{x}_{(n+2)k} &x_{(n+2)(k+1)} &\dots &x_{(n+2)(n+1)}\\
        \end{pmatrix}
\end{equation*}
where the components with a hat are undetermined, this data identifies the full simplicial matrix uniquely, which (by symmetry of the matrix) is equivalent to both $\Kan!(n+2, k)$ and $\Kan!(n+2, k+1)$ for $0\le k \le n+1$. This problem is actually overdetermined, and the consistency equation $\widehat{x}_{kk} = \widehat{x}_{(k+1)k}$ coming from the symmetry of the matrix is precisely \eqref{eq:n-gpd-associativity}. 
This is discussed in detail for $n=1,2$ in Example \ref{ex:multiplications-1gpd-assoc} and Remark \ref{rem:2-gpd-assoc-simp-mat}
\end{remark}

\begin{proposition}\label{prop:finite-data-n-gpds}
An $n$-groupoid object $\huaX$ in $(\Cat, \covers, \mathsf{P})$ is equivalent to the following finite data:
\begin{enumerate}
\item An $n$-truncated simplicial object $\widebar{\huaX}$, satisfying $\Kan(m,j)$ for all $0\le m\le n$. 
\item A set of $(n+1)$-ary multiplication maps $m_k:\Lambda^{n+1}_k(\widebar{\huaX}) \to \widebar{\huaX}_n$ for $0  \le k \le n+1$ such that:
\begin{enumerate}
\item They induce isomorphisms $p^{n+1}_jb_k: \Lambda^{n+1}_k(\widebar{\huaX}) \to \Lambda^{n+1}_j(\widebar{\huaX})$, for all $0\le j,k \le n+1$, and $b_k$ as in \eqref{eq:multiplication-defined-iso}.
\item They are compatible with the degeneracy maps $s_i:\widebar{\huaX}_n \to \Lambda^{n+1}_k(\widebar{\huaX}) \subseteq \partial_{n+1}(\widebar{\huaX})$ defined in \eqref{eq:cosk-degeneracies}, in the sense that \eqref{eq:n-finite-data-comp-deg-1}
(or equivalently \eqref{eq:n-finite-data-comp-deg-2}) holds for any $x \in \widebar{\huaX}_n$.
\item At least one of the $m_k$ is associative. (This and (2a) imply they are all associative).
\end{enumerate}
\end{enumerate}
\end{proposition}

\begin{proof}
In the previous discussion we have shown that starting with a Lie $n$-\hspace{0pt}groupoid $\huaX$ we obtain the finite data above, where the multiplications $m_k$ of a Lie $n$-\hspace{0pt}groupoid satisfy the conditions 2.a-c, and in particular they are associative because the conditions $\Kan!(n+2, k)$ hold for any $k$. 

In the other direction, by (1) we have that $\Lambda^{n+1}_k(\widebar{\huaX})$ is representable for any $k$. Fix a value of $k$. Then, as discussed in Remark \ref{rem:multiplication-comp-with-deg}, by (2b), we can construct the simplicial object
\begin{equation*}
    \widetilde{\huaX} = \Sk^{n+1}(\Lambda^{n+1}_k(\widebar{\huaX})\rightfourarrows \widebar{\huaX}_n \rightfourarrows \dots \rightrightarrows \widebar{\huaX}_0), 
\end{equation*}
as in \ref{eq:sk-of-n-trunc-with-k-horn}. This satisfies $\Kan(m, j)$ for all $0 \le m \le n$, $0 \le j \le m$, and $\Kan!(n+1, j)$ for all $0 \le j \le n+1$, by (2a). 
Therefore $\Lambda^{n+2}_j(\widetilde{\huaX})$ is representable for all $0 \le j \le n+2$. Moreover, $\widetilde{\huaX}$ is weakly stalkwise $n$-coskeletal, i.e. for any point functor $\mathsf{p} \in \mathsf{P}$, 
\begin{equation*}
    \mathsf{p}\widetilde{\huaX} \subseteq \Cosk^n(\mathsf{p}\widetilde{\huaX}). 
\end{equation*}
Then for $\mathsf{p} \in \mathsf{P}$, the set 
$\Cosk^{n+1}(\mathsf{p}\widetilde{\huaX})_{n+2}\subseteq \Cosk^{n}(\mathsf{p}\widetilde{\huaX})_{n+2}$ can be seen as the set of simplicial matrices in $\mathsf{p}\widetilde{\huaX}$ with each row in $\Lambda^{n+1}_k(\widebar{\huaX})$. Associativity in (2c) is exactly the statement that for any such matrix with a missing row there is a unique full matrix identified by it. Therefore $\Cosk^{n+1}(\mathsf{p}\widetilde{\huaX})$ satisfies $\Kan!(n+2, k)$, which means
\begin{equation*}
    \Cosk^{n+1}(\mathsf{p}\widetilde{\huaX})_{n+2} \cong \hom(\Lambda^{n+2}_k, \Cosk^{n+1}(\mathsf{p}\widetilde{\huaX})) 
    \cong \hom(\Lambda^{n+2}_k, \mathsf{p}\widetilde{\huaX}).
\end{equation*}
Because this holds for any $\mathsf{p}$ and $\mathsf{P}$ is jointly conservative, 
\begin{equation*}
    \Cosk^{n+1}(\widetilde{\huaX})_{n+2} \cong
    \hom(\Lambda^{n+2}_k, \widetilde{\huaX}),
\end{equation*}
which is representable, as previously mentioned. So $\Cosk^{n+1}(\widetilde{\huaX})_{\le n+2}$ satisfies $\Kan!(n+2, k)$ and all the Kan conditions that $\widetilde{\huaX}$ already satisfied. 
By \cite[Lemma 3.21]{Li2014}, because this satisfies all strict Kan conditions at level $n+1$, and a strict Kan condition at level $n+2$, it then satisfies all strict Kan conditions at level $n+2$. 
Finally, we prove inductively that $\Cosk^{n+1}(\widetilde{\huaX})$ is representable and it satisfies all strict Kan conditions above level $n+2$. The base case is $l =n+2$, which we just proved. For $l > n+2$ we have that
\begin{equation*}
    \Cosk^{n+1}(\mathsf{p}\widetilde{\huaX})_{l} \cong \hom(\Lambda^{l}_j, \Cosk^{n+1}(\mathsf{p}\widetilde{\huaX})),
\end{equation*}
by Lemma \ref{lem:prop2.37-Li}. Assuming that $\Cosk^{n+1}(\widetilde{\huaX})$ satisfies $\Kan!(l-1, j)$ for all $0 \le j \le l-1$, we know that $\hom(\Lambda^{l}_j, \Cosk^{n+1}(\widetilde{\huaX}))$ is representable. Thus by the jointly conservative property, $\Cosk^{n+1}(\widetilde{\huaX})_l$ is representable and $\Kan!(l, j)$ holds for all $0 \le j \le l$.
Hence $\Cosk^{n+1}(\widetilde{\huaX})$ is an $n$-groupoid object in $\Cat$.
\end{proof}

We describe explicitly the finite data of a 1-groupoid and of a 2-groupoid object in the following examples. 

\begin{example}\label{ex:multiplications-1gpd-assoc}
In the case of a 1-groupoid object $\huaG$, the three binary multiplication maps are the multiplication and the two divisions:
For any $g,h,k \in \huaG_1$ such that $d_0g = d_0h$, $d_1g = d_0k$ and $d_1h = d_1k$, we have
\begin{equation*}
    m_0(h, k) = \square h k = k^{-1} h, \quad
    m_1(g, k) = g \square k = k g, \quad 
    m_2(g, h) = g h \square = h g^{-1}. 
\end{equation*}
\begin{center}
    \begin{tikzpicture}[scale=0.6]
	\begin{pgfonlayer}{nodelayer}
		\node [style=dot, label={below left:\small{0}}] (3) at (0, 0) {};
		\node [style=dot, label={below right:\small{2}}] (4) at (3, 0) {};
		\node [style=dot, label={above right:\small{1}}] (5) at (1.5, 2.5) {};
		\node [style=dot, label={below left:\small{0}}] (18) at (6.5, 0) {};
		\node [style=dot, label={below right:\small{2}}] (19) at (9.5, 0) {};
		\node [style=dot, label={above right:\small{1}}] (20) at (8, 2.5) {};
		\node [style=dot, label={below left:\small{0}}] (22) at (13, 0) {};
		\node [style=dot, label={below right:\small{2}}] (23) at (16, 0) {};
		\node [style=dot, label={above right:\small{1}}] (24) at (14.5, 2.5) {};
		\node [style=none] (25) at (1.5, -0.75) {$h$};
		\node [style=none] (26) at (0, 1.75) {$k$};
		\node [style=none] (27) at (3.25, 1.75) {$k^{-1}h$};
		\node [style=none] (28) at (8, -0.75) {$kg$};
		\node [style=none] (29) at (6.5, 1.75) {$k$};
		\node [style=none] (30) at (9.5, 1.75) {$g$};
		\node [style=none] (31) at (14.5, -0.75) {$h$};
		\node [style=none] (32) at (13, 1.75) {$hg^{-1}$};
		\node [style=none] (33) at (16, 1.75) {$g$};
	\end{pgfonlayer}
	\begin{pgfonlayer}{edgelayer}
		\draw [style=directed] (5) to (3);
		\draw [style=directed-red] (4) to (5);
		\draw [style=directed, in=360, out=180] (4) to (3);
		\draw [style=directed] (20) to (18);
		\draw [style=directed] (19) to (20);
		\draw [style=directed-red] (19) to (18);
		\draw [style=directed-red] (24) to (22);
		\draw [style=directed] (23) to (24);
		\draw [style=directed] (23) to (22);
	\end{pgfonlayer}
\end{tikzpicture}
\end{center}

Here we keep the convention of writing tuples of faces ordered by index as in \eqref{eq:DefHornBdrySpaces} and Example \ref{ex:nerve-cat}. We use juxtaposition to indicate the actual groupoid product intended as composition of arrows and written in its usual order: $k g$ means ``$k$ after $g$''. 
The inversion map, which is often included as part of the defining data of a groupoid, is alternatively encoded in the multiplication maps $m_0$ and $m_2$. 
The very fact that $m_0$ and $m_2$ are ``divisions'' (i.e. multiplications by an inversion) follows from the fact that the multiplications define the isomorphisms $b_0, b_1, b_2$ form Remark \ref{rem:multiplication-defined-iso}. 

Another requirement on the multiplications is that $1=s_0: \huaG_0 \to \huaG_1$ acts as the unit map. This is described by compatibility with the degeneracy maps of the 2-coskeleton between levels $1$ and $2$, as defined in \eqref{eq:cosk-degeneracies}. In fact, these maps are
\begin{equation*}
\begin{array}{lcl}
    \begin{aligned}
        &s_0 : \huaG_1 \to \partial_2(\huaG)\\
        &s_0g = (g, g, 1d_1g),
    \end{aligned}
    &\qquad\qquad &
    \begin{aligned}
        &s_1 : \huaG_1 \to \partial_2(\huaG)\\
        &s_1g = (1d_0g, g, g).
    \end{aligned}
\end{array}
\end{equation*}
\begin{center}
\begin{tikzpicture}[scale=0.6]
	\begin{pgfonlayer}{nodelayer}
		\node [style=dot, label={below left:\small{0}}] (3) at (9, 0) {};
		\node [style=dot, label={below right:\small{2}}] (4) at (12, 0) {};
		\node [style=dot, label={above right:\small{1}}] (5) at (10.5, 2.5) {};
		\node [style=dot, label={below left:\small{0}}] (22) at (0, 0) {};
		\node [style=dot, label={below right:\small{2}}] (23) at (3, 0) {};
		\node [style=dot, label={above right:\small{1}}] (24) at (1.5, 2.5) {};
		\node [style=none] (25) at (10.5, -0.75) {$g$};
		\node [style=none] (26) at (9, 1.75) {$g$};
		\node [style=none] (27) at (12.25, 1.75) {$1d_0g$};
		\node [style=none] (31) at (1.5, -0.75) {$g$};
		\node [style=none] (32) at (0, 1.75) {$1d_1g$};
		\node [style=none] (33) at (3, 1.75) {$g$};
	\end{pgfonlayer}
	\begin{pgfonlayer}{edgelayer}
		\draw [style=directed] (5) to (3);
		\draw [style=directed-dash] (4) to (5);
		\draw [style=directed, in=360, out=180] (4) to (3);
		\draw [style=directed-dash] (24) to (22);
		\draw [style=directed] (23) to (24);
		\draw [style=directed] (23) to (22);
	\end{pgfonlayer}
\end{tikzpicture}
\end{center}
Compatibility of the multiplications with them is
\begin{equation*}
    m_0 p^2_0 s_0 g = (1d_1g)^{-1}g = g, \qquad 
    m_1 p^2_1 s_1 g = g(1d_0g) = g, 
\end{equation*}
or equivalently,
\begin{equation*}
    m_1 p^2_1 s_0 g = (1d_1g)g = g, \qquad 
    m_2 p^2_2 s_1 g = g(1d_0g)^{-1} = g.
\end{equation*}
These are precisely the unit equations for $1$.

The groupoid product $m_1$ is required to be associative. This means that for any $g,h,k$ such that $d_1g = d_0h$, $d_1h =d_0k$, 
\begin{equation}\label{eq:assoc-1-gpd}
    k (h g) = (k h) g.
\end{equation}
In the nerve of the groupoid, such $g,h$ and $k$ define a 3-simplex, as in \ref{ex:nerve-cat}. Because any 1-groupoid is 2-coskeletal, this 3-simplex can be written as a simplicial matrix. Each row of the matrix is the boundary of a triangle in the 1-groupoid, so it is determined by two of its three faces. This and the symmetry properties of the matrix make it so that three arrows $g,h,k$ determine the entirety of it. In addition, the associativity condition can be expressed by observing that
\begin{equation*}
    \begin{pmatrix}
            g &hg &h\\
            g &(kh)g &kh\\
            hg &k(hg) &k\\
            h &kh &k
    \end{pmatrix}
    \in \huaG_3
\end{equation*}
is the matrix determined by $g,h$ and $k$. Because it must be symmetric, we recover \eqref{eq:assoc-1-gpd}. This matrix represents the tetrahedron
\begin{center}
\begin{tikzpicture}[scale=0.6]
	\begin{pgfonlayer}{nodelayer}
		\node [style=dot, label={above right:\small{3}}] (5) at (7.75, 4.5) {};
		\node [style=dot, label={below left:\small{0}}] (22) at (0, 0) {};
		\node [style=dot, label={below right:\small{2}}] (23) at (6, 0) {};
		\node [style=dot, label={above left:\small{1}}] (24) at (2, 4.5) {};
		\node [style=none] (31) at (3, -0.5) {\small{$kh$}};
		\node [style=none] (32) at (0, 2) {\small{$k$}};
		\node [style=none] (33) at (3.75, 3.75) {\small{$h$}};
		\node [style=none] (34) at (7.75, 1.75) {\small{$g$}};
		\node [style=none] (35) at (5, 5) {\small{$hg$}};
		\node [style=none, label={[rotate=30]\small{$k(hg)=(kh)g$}}] (36) at (5, 1.5) {};
	\end{pgfonlayer}
	\begin{pgfonlayer}{edgelayer}
		\draw [style=directed] (24) to (22);
		\draw [style=directed,opacity=0.5] (23) to (24);
		\draw [style=directed] (23) to (22);
		\draw [style=directed] (5) to (23);
		\draw [style=directed] (24) to (5);
		\draw [style=directed] (22) to (5);
	\end{pgfonlayer}
\end{tikzpicture}
\end{center}
\end{example}

\begin{example}\label{ex:finite-data-2gpds}
By writing down the data of Proposition \ref{prop:finite-data-n-gpds} for a 2-groupoid object we obtain the data described in \cite[Prop-Def 2.16]{Zhu2009}.
In fact, we obtain that a 2-groupoid object $\huaX$ in $(\Cat, \covers, \mathsf{P})$ is equivalent to the following finite data:
\begin{enumerate}
\item A 2-truncated simplicial object $\huaX_2 \rightthreearrows \huaX_1 \rightrightarrows \huaX_0$, such that 
\begin{enumerate}
    \item $d^1_0, d^1_1: \huaX_1 \to \huaX_0$ are covers,
    \item $p^2_k: \huaX_2 \to  \Lambda^{2}_k(\huaX)$ are covers. 
\end{enumerate} 
\item Four morphisms called 3-multiplications $m_k:\Lambda^3_k(\huaX) \to \huaX_2$ for $0  \le k \le 3$ such that:
\begin{enumerate}
    \item They induce isomorphisms $p^{n+1}_jb_k: \Lambda^3_k(\huaX) \to \Lambda^3_j(\huaX)$, for $b_k$ as in \eqref{eq:multiplication-defined-iso}.
    \item They are compatible with the degeneracy maps $s_i:\huaX \to \Lambda^3_k(\huaX) \subseteq \partial_3(\huaX)$ defined in \eqref{eq:cosk-degeneracies}, in the sense that for any $X \in \huaX_2$,
    \begin{equation}\label{eq:finite-data-comp-deg-1}
        m_k(p^{n+1}_k s_k X) = X, \text{ for all } 0\le k \le n,
    \end{equation}
    or equivalently
    \begin{equation}\label{eq:finite-data-comp-deg-2}
        m_k(p^{n+1}_k s_{k-1} X) = X, \text{ for all } 1\le k \le n+1,
    \end{equation}
    \item they are associative, in the sense that,
    for any triangles $X_{0i4}$ and $X_{0ij}$ for $i\neq j \in \{1,2,3\}$ that fit as the faces of a 4-simplex $X_{01234}$ (with $d_id_jX_{01234} = X_{\delta_i\delta_j(01234)}$), the face $X_{123}$ can be determined equivalently by either
    \begin{enumerate}
        \item $X_{123} = m_0(X_{023}, X_{013}, X_{012})$, or
        \item $X_{123} = m_3(X_{234}, X_{134}, X_{124})$, where $X_{ij4} = m_0(X_{0j4}, X_{0i4}, X_{0ij})$.
    \end{enumerate}
    That is 
    \begin{equation}\label{eq:2gpd-associativity}
    \begin{split}
        (\square X_{034}X_{024}X_{023})(\square X_{034}&X_{014}X_{013})(\square X_{024}X_{014}X_{012})\square\\
        &= \square X_{023}X_{013}X_{012}.
    \end{split}
\end{equation}           
\end{enumerate}
\end{enumerate}

\begin{figure}[!ht]
    \centering
    \includegraphics[width=\textwidth]{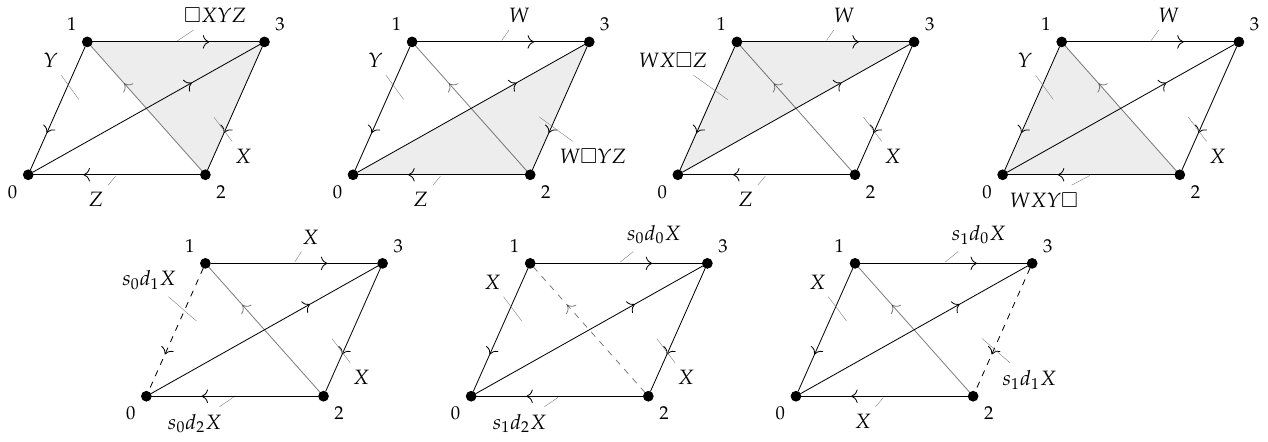}
    \caption{The multiplications $m_k$ and degeneracy maps $s^2_j$ of a 2-groupoid.}
\end{figure}

\end{example}

\begin{remark}
    More concretely, the compatibility of the multiplications of a 2-\hspace{0pt}groupoid with the degeneracies can be written as the following identities, which intuitively correspond to higher versions of the unit-defining identities for a 1-groupoid.
    \begin{equation*}
    \begin{array}{ll}
        X = \square X (s_0d_1X) (s_0d_2X),
        & X = X \square (s_0d_1X) (s_0d_2X),\\
        X = (s_0d_0X) \square X (s_1d_2X),
        & X = (s_0d_0X) X \square (s_1d_2X),\\
        X = (s_1d_0X) (s_1d_1X) \square X,
        & X = (s_1d_0X) (s_1d_1X) X \square.
    \end{array}
    \end{equation*}
\end{remark}

\begin{remark}\label{rem:2-gpd-assoc-simp-mat}
The associativity condition \eqref{eq:2gpd-associativity} for item (2c) can be rewritten in terms of existence of a simplicial matrix as in Remark \ref{rem:n-gpd-assoc-simp-mat}. 
The starting data fits in the incomplete matrix
\begin{equation*}
    \begin{pmatrix}
        ? & ? & ? & ?\\
        ? & X_{034} & X_{024} & X_{023}\\
        ? & X_{034} & X_{014} & X_{013}\\
        ? & X_{024} & X_{014} & X_{012}\\
        ? & X_{023} & X_{013} & X_{012}
    \end{pmatrix}.
\end{equation*}
Each row but the first one is a $(3,0)$-horn, so by using $m_0$ and using symmetry we have that
\begin{equation*}
    \begin{pmatrix}
        ? & \square X_{034}X_{014}X_{013} & \square X_{024}X_{014}X_{012} & \square X_{023}X_{013}X_{012}\\
        \square X_{034}X_{024}X_{023} & X_{034} & X_{024} & X_{023}\\
        \square X_{034}X_{014}X_{013} & X_{034} & X_{014} & X_{013}\\
        \square X_{024}X_{014}X_{012} & X_{024} & X_{014} & X_{012}\\
        \square X_{023}X_{013}X_{012} & X_{023} & X_{013} & X_{012}
    \end{pmatrix} 
\end{equation*}
At this point the first row can be completed in two ways: by symmetry or by using $m_0$ on the first row. These two ways must be equivalent for this to be a well-defined 4-simplex in $\huaX$. Therefore we obtain \eqref{eq:2gpd-associativity}.
\end{remark}

We conclude this section by showing that the finite data construction extends to simplicial maps with target an $n$-groupoid object, which are then determined by their $n$-truncation, subject to some conditions at level $n$. This is a fundamental property that we will use repeatedly in the rest of the thesis. 

\begin{lemma}\label{lem:finite-data-simp-map}
    Let $\huaX$ be a simplicial object and $\huaY$ be an $n$-groupoid object in $(\Cat, \covers, \mathsf{P})$. Then the data of a simplicial map $f:\huaX \to \huaY$ is equivalent to the data of an $n$-truncated simplicial map $f':\tr^n\huaX \to \tr^n\huaY$, such that for some $0\le k \le n+1$ and any $X \in \huaX_{n+1}$,
    \begin{equation}\label{eq:simp-map-multiplicative}
        f'_n(d_k x)= m_k((f'_n d_i x)_{i\neq k}). 
    \end{equation}
    We call such a map $f'$ \textbf{multiplicative}. 

    Moreover, if this holds for one value of $k$ it holds for all other values of $k$. 
\end{lemma}

\begin{proof}
    By Lemma \ref{lem:n-gpd-coskeletality}, $\huaY \cong \Cosk^{n+1}(\huaY)$. Therefore 
    \begin{equation*}
        \begin{split}
        \Simp\Cat(\huaX, \huaY) 
        &\cong \Simp\Cat (\huaX, \Cosk^{n+1}\Sk^{n+1}\huaY) 
        \cong \Simp\Cat(\Sk^{n+1}\huaX, \Sk^{n+1}\huaY)\\
        &\cong \Simp\Cat(\Sk^{n+1}\huaX, \Sk^{n+1}(\Lambda^{n+1}_k(\huaY) \rightfourarrows \tr^n{\huaY})),
        \end{split}
    \end{equation*}
    so each map $f:\huaX \to \huaY$ depends only on its components $f_0, \dots, f_{n+1}$. We now show that $f_{n+1}: \huaX_{n+1} \to \huaY_{n+1}\cong\Lambda^{n+1}_k(\huaY)$ depends only on $f_n$ and it commutes with the simplicial maps between levels $n$ and $n+1$ if and only if it satisfies the multiplicativity condition \eqref{eq:simp-map-multiplicative}. Commutativity with the face maps for $i \neq k$ implies that for any $x\in \huaX_{n+1}$, $f_{n+1}(x) = (f(d_0 x), \dots, \widehat{f(d_k x)}, \dots, f(d_n x)) = (f(d_i x))_{i \neq k}$. Since the $k$-th face map in $\Lambda^{n+1}_k(\huaY) \rightthreearrows \tr^n{\huaY}$ is precisely $m_k$, we have that commutativity with the $k$-th face map is \eqref{eq:simp-map-multiplicative}. Because of the way the degeneracy maps are defined between levels $n$ and $n+1$ (i.e. as in \eqref{eq:cosk-degeneracies}), the commutativity with these follows automatically.

\end{proof}

\section{Simplicial vector spaces}\label{sec:simp-vect}

In this section we focus on simplicial vector spaces, with a more specific attention to levelwise finite-dimensional ones over a field of characteristic zero which is assumed to be $\R$ for simplicity. Most of the material we discuss works for simplicial objects in other abelian categories such as abelian groups and modules over a more general ring, with the appropriate considerations. 
The content of this section is mostly a rewriting of that contained in the preprint \cite{RonchiZhu2024}.

The definition of space of simplicial diagrams of shape $\huaS$ inside a simplicial vector space translates directly from that of simplicial sets in Definition \ref{def:simpl-diags} because all limits exist in $\Vect$ and they are even created by the forgetful functor $\Forget$, which forms the monadic adjunction $\R[\_] \dashv \Forget$ with the free vector space functor $\R[\_]$.

\begin{definition}\label{def:simp-diag-svec}
    Let $\huaV$ be a simplicial vector space. The \textbf{space of simplicial diagrams} of shape $\huaS$ in $\huaV$ is 
    \begin{equation*}
        \hom(\huaS, \huaV) := \SSet(\huaS, \Forget(\huaV)) = \SVect(\R[\huaS], \huaV),
    \end{equation*}
    with the "pointwise" vector space structure inherited from $\huaV$. 
\end{definition}

\begin{remark}
The functor $\hom$ is level 0 of the powering functor $\powering: \SSet \times \SVect \to \SVect$ constructed as in Remark \ref{rem:diagram-space-and-powering}. 
With this, $\SVect$ is powered over $\SSet$. 
\end{remark}

With this definition we have the vector spaces of $m$-simplices $\hom(\Delta^m, \huaV) = \huaV_m$, $m$-boundaries $\hom(\partial\Delta^m, \huaV) = \partial\Delta^m(\huaV)$, and $(m,k)$-horns $\hom(\Lambda^m_k, \huaV) = \Lambda^m_k(\huaV)$, whose underlying sets are defined as in \eqref{eq:DefHornBdrySpaces}. 
It is also evident that the Kan conditions in Definition \ref{def:n-gpd} translate immediately to the category of simplicial vector spaces, where the horn projections are additionally linear. 

\begin{definition}\label{def:vs-n-gpd}
    A \textbf{$\VS$ $n$-groupoid} $\huaV$ is a simplicial vector space whose underlying set is an $n$-groupoid. A {\bf $\VS$ $n$-group}  $\huaV$ is a $\VS$ $n$-groupoid with $\huaV_0=0$. 
\end{definition}

Simplicial vector spaces are in particular simplicial groups, and it is a general fact that these are always $\infty$-groupoids, because one can construct arbitrary fillers in each dimension by using the group operations. This result is due to Moore \cite[Thm. 3 p. 18-04]{Moore1954}, see also \cite[Thm. 17.1]{May1967}, \cite[Lem. 3.1]{Curtis1971}, \cite[Lem. 8.2.8]{Weibel1994}, \cite[Lem. I.3.4]{GoerssJardine2009}. We recall it here and give an explicit proof in the case of vector spaces, adapted from \cite{nlab:simplicial_group}, as it will be useful for our later discussion. This version is also a particular case of \cite[prop. 6.6]{Wolfson2016}, which applies to the more general setting of simplicial Lie groups.

\begin{proposition}[Moore]\label{thm:MooreHornFillers}
    Simplicial groups are Kan complexes. In particular, simplicial vector spaces are $\VS$ $\infty$-groupoids. 
\end{proposition}
\begin{proof}
    Let $\huaV$ be a simplicial vector space, and $(x_0, \dots, \widehat{x_{k}}, \dots, x_n) \in \Lambda^n_k(\huaV)$ an $(n,k)$-horn in $\huaV$.
    That is, $d_i x_j = d_j x_{i+1}$ for any $i\ge j$, whenever both sides are defined. 
    We describe an algorithm to construct a horn filler $x \in \huaV_{m+1}$, i.e. an $x$ such that $d_ix = x_i$ for any $i\neq k$. 
    Let us start with the case of $0<k<n$.
    Begin with $w_0 = s_0x_0 \in \huaV_{n+1}$.
    By the simplicial identities, $d_0w_0 = x_0$. 
    Proceeding by induction on $0 \le i \le k-1$, assume that $w_{i-1}$ has been defined in a way such that $d_{j}w_{i-1} = x_{j}$ for any $0 \le j \le i-1$. 
    Then define 
    \begin{equation}
        w_i := w_{i-1} - s_id_iw_{i-1} + s_ix_i.
    \end{equation}
    Now a simple computation using the simplicial identities shows $d_iw_i = x_i$. 
    Moreover, for any $0 \le j \le i-1$, 
    \begin{equation}\label{eq:MooreThmComp1}
        \begin{split}
            d_jw_i &= d_jw_{i-1} - d_js_id_iw_{i-1} + d_js_ix_i\\
            &= d_jw_{i-1} - s_{i-1}d_jd_iw_{i-1} + s_{i-1}d_jx_i\\
            &= d_jw_{i-1} - s_{i-1}d_{i-1}d_jw_{i-1} + s_{i-1}d_{i-1}x_j = d_jw_{i-1} = x_j.
        \end{split}
    \end{equation}
    Hence $w_{k-1}$ has the property that $d_iw_{k-1}=x_i$ for $0\le i \le k-1$.
    Now, define $w_n = w_{k-1} - s_{n-1}d_nw_{k-1} + s_{n-1}x_n$. 
    As before, $d_n w_n = x_n$ and $d_iw_n = x_i$ for any $0\le i \le k-1$.
    Proceed again by induction, assuming that $w_{i+1}$ has been defined for any $k+1 \le i \le n$ and that $d_j w_{i+1} = x_j$ for any $j$ such that $0 \le j \le k-1$ or $i+1 \le j \le n$. 
    Then define 
    \begin{equation}
        w_i := w_{i+1} - s_{i-1}d_iw_{i+1} + s_{i-1}x_i.
    \end{equation}
    Again, we have that $d_iw_i = x_i$ and a computation similar to \eqref{eq:MooreThmComp1}, shows that $d_jw_{i+1}=x_j$ for any $j$ such that $0 \le j \le k-1$ or $i+1 \le j \le n$.
    Therefore, $w_{k+1}$ is the required horn filler, as we have that $d_jw_{k+1} = x_j$ for any $j\neq k$. 

    In the case of $k=n$, starting with $w_0$ as above and repeating the first induction procedure until $w_{n-1}$, we get that this is already the required horn filler. 
    
    In a similar way, for $k=0$, we start with $w_n$ as above and proceed by downward induction as in the second part above to obtain the horn filler $w_1$. 
\end{proof}

Given  a simplicial vector space $\huaV$, we define a linear map $\mu^m_k: \Lambda^m_k(\huaV) \to \huaV_m$ by the algorithm in the proof of Proposition \ref{thm:MooreHornFillers}, and we call $\mu^m_k(\lambda_k) \in \huaV_m$ the \textbf{Moore filler} of $\lambda_k \in \Lambda^m_k(\huaV)$. This is a right inverse of the respective horn projection. If $\huaV$ is further a $\VS$ $n$-groupoid, then $\mu^m_k$ is the \textit{double-sided} inverse of the horn projection $p^m_k: \huaV_m \to \Lambda^m_k(\huaV)$ for $m>n$. This makes the structure of a $\VS$ $n$-groupoid particularly rigid, as we will see in Proposition \ref{prop:FiniteDataVS}. 
Roughly speaking, in light of the discussion on finite data of $n$-groupoid objects in Section \ref{sec:finite-data}, the Moore fillers describe canonical multiplication maps for each $\VS$ $n$-groupoid. 
We will discuss this in Section \ref{sec:FiniteDataVS}, after introducing the Dold-Kan correspondence, as the latter will help in showing that the Moore fillers are compatible with the degeneracy maps in the sense of \ref{rem:multiplication-comp-with-deg}.
The following example illustrates this situation for $n=1$. 

\begin{example}\label{ex:MooreFillers1Gpd}
    Let $\huaV_1 \rightrightarrows  \huaV_0$ be a category object internal to vector spaces. This is a 2-vector space in the sense of \cite{BaezCrans2004}. A priori this data includes a multiplication map which can be seen as $m_1=d_1(p^2_1)^{-1} : \Lambda^2_1(\huaV) \to \huaV_1$ when identifying the category with its nerve, which is a simplicial vector space. By Proposition \ref{thm:MooreHornFillers}, the inverse $(p^2_1)^{-1}$ must coincide with the Moore filler $\mu^m_k$. Therefore, there is a unique possible linear multiplication on a 2-vector space with space of arrows $\huaV_1$ and space of points $\huaV_0$, that is, for any $u,w \in \huaV_1$ meeting at $x = d_1u = d_0w \in \huaV_0$,
    \begin{equation}\label{eq:VS1gpd-canon-mult-1}
        w \cdot u = m_1(u,w) = d_1\mu^2_1(u,w) = u + w - 1x.
    \end{equation}
    This was also observed in \cite[Lemma 3.2]{BaezCrans2004}.

    If $\huaV_1 \rightrightarrows  \huaV_0$ is a groupoid internal to $\Vect$, then the same argument can be applied to the division maps, which are also unique and given for any $u,v,w \in \huaV_1$, such that $d_0u = d_0v$ and $d_1v = d_1w$, by $d_i\mu^2_i$ for $i=0,2$:
    \begin{equation}\label{eq:VS1gpd-canon-mult-2}
        \begin{split}
            w^{-1} \cdot v &= m_0(v,w) = v - w + 1d_1w,\\
            v\cdot u^{-1} &= m_2(u,v) = v - u + 1d_0u.
        \end{split}
    \end{equation}
    By inserting the appropriate unit element in these formulas, we obtain that the groupoid inversion is 
    \begin{equation}\label{eq:VS1gpd-canon-mult-3}
    \begin{split}
        v^{-1} = - v + 1d_0v + 1d_1v. 
    \end{split}
    \end{equation}

    A posteriori this inversion is also well-defined if $\huaV_1 \rightrightarrows  \huaV_0$ is only a category, which means that any 2-vector space as a category internal to $\Vect$ is automatically a $\VS$ 1-groupoid.\footnote{This is an instance of a result about internal categories in $\mathsf{Grp}$ appearing in \cite{BrownSpencer1976}, which the authors credit to Duskin.} Furthermore any of these objects is completely determined by a pair of vector spaces with source, target and unit maps between them, as the multiplication and inversion can be inferred from these. We state a general version of this fact in Proposition \ref{prop:FiniteDataVS}. 

\end{example}

\subsection{Closed monoidal structures and the Dold-Kan correspondence}\label{sec:review-monoidal-DK}

The Dold-Kan correspondence establishes an equivalence of categories between the category of simplicial vector spaces $\SVect$ and that of non-negative\footnote{We use non-negative as shorthand for ``concentrated in non-negative degrees''.} chain complexes $\Ch_{\ge 0}(\Vect)$ (see for example \cite[\S 8.4]{Weibel1994}, \cite[\S 22]{May1967}, \cite[\S III.2]{GoerssJardine2009}, \cite[\href{https://kerodon.net/tag/00QQ}{Tag 00QQ}]{kerodon}).
Similarly to $\SSet$, as discussed in Section \ref{sec:SSet-products-and-mapping-spaces}, Each of these categories has a closed symmetric monoidal structure. In other words, they are equipped with a tensor product which admits a natural symmetry isomorphism and a right adjoint defining internal hom spaces.
The tensor-hom adjunction also holds in both categories in the enriched sense, by replacing the hom functor with the internal hom as in Proposition \ref{prop:enriched-tensor-hom-sset}. 
This causes it to descend to an adjunction in the respective homotopy categories as well.

\subsubsection{The simplicial category of simplicial vector spaces and its homotopy category}\label{sec:Svect-int-hom}

The monoidal product $\otimes$ in $\SVect$ is given levelwise by the tensor product of vector spaces. 
For any two simplicial vector spaces $\huaV$ and $\huaW$, the face and degeneracy maps of $\huaV \otimes \huaW$ are $d_i^\huaV \otimes d_i^\huaW$ and $s_j^\huaV \otimes s_j^\huaW$, respectively. 
This product is symmetric with respect to the symmetry isomorphism $v \otimes w \mapsto w \otimes v$.
In addition, there is a way to ``tensor'' a simplicial vector space and a simplicial set, known as \textbf{copowering}, which is defined as in Def. \ref{def:copowering}: at each level $m$ and for any $\huaV \in \Simp\Vect$ and $\huaK \in \Simp\Set$, this is
\begin{equation}\label{eq:copowering}
(\huaV \otimes \huaK)_m
:= \bigoplus_{r \in \huaK_m} \huaV_m^r = \huaV_m \otimes \R [\huaK_m]
\end{equation}
i.e. the coproduct of a copy of $\huaV_m$ for each $m$-simplex $r$ of $\huaK$, which is equivalent to the tensor product between $\huaV_m$ and the vector space generated by $\huaK_m$.
The simplicial maps of $\huaV \otimes \huaK$ are the same as for the tensor product, but we write them in the following way, using the coproduct description
\begin{equation*}
\begin{split}
    d_i(v^r) 
    &= (d_iv)^{d_i r}
    \in \huaV_{m-1}^{d_i r} 
    \subseteq (\huaV \otimes \huaK)_{m-1},\\
    s_i(v^r) 
    &= (s_iv)^{s_i r}
    \in \huaV_{m+1}^{s_i r} 
    \subseteq (\huaV \otimes \huaK)_{m+1},
\end{split}
\end{equation*}
for any $v^{r} \in \huaV_m^{r} \subseteq (\huaV \otimes \huaK)_m$. 
This is used to define $\IHom(\huaV, \huaW)$, the \textbf{mapping space} between two simplicial vector spaces as in Definition \ref{def:mapping-set}. This is a simplicial set with $m$-th level 
\begin{equation}\label{eq:HomSpaceDef}
    \IHom(\huaV, \huaW)_m
    := \SVect(\huaV \otimes \Delta^m, \huaW),
\end{equation}
where the right-hand side is the usual hom set in $\SVect$, and the face and degeneracy maps are defined as in \eqref{eq:mapping-set-maps}:
\begin{equation*}
    \begin{split}
        (d_if)_l (v^r) &:= f_l (v^{\delta^i(r)}), \quad \forall v \in \huaV_{l}, r \in \Delta^{m-1}_{l},\\
        (s_if)_l (v^r) &:= f_l (v^{\sigma^i(r)}), \quad \forall v \in \huaV_{l}, r \in \Delta^{m+1}_{l},
    \end{split}
\end{equation*}
where $\delta^i: \Delta^{m-1}_{l} \to  \Delta^m_{l}$ are the coface maps, and $\sigma^i: \Delta^{m+1}_{l} \to \Delta^m_{l}$ are the codegeneracy maps. 
Here we recall that when writing simplices in $\Delta^m_l$ as sequences of $l$ numbers in $\{0,1,2,\dots,n\}$, the coface map $\delta^i$ is the postcomposition with the unique injection of ordinals that skips $i$, and the codegeneracy map $\sigma^i$ is the postcomposition with the unique surjection of ordinals that repeats $i$.
Because each level of $\IHom(\huaV, \huaW)$ is a vector space and an immediate computation shows that the face and degeneracy maps defined in \eqref{eq:mapping-set-maps} are linear, $\IHom$ is actually a functor $\SVect \times \SVect \to \SVect$, whose level 0 coincides with the hom space $\SVect$. Its level 1 is the space of simplicial linear homotopies between simplicial linear maps, and its higher levels encode increasingly higher homotopies.
These can be composed across different mapping spaces according to the composition defined in \eqref{eq:composition-of-homotopies}.
With this composition and the obvious identity element, $\SVect$ is not just an $\SSet$-category, as previously remarked, but also an $\SVect$-category. 

We can now define the \textbf{homotopy category} $\hSVect$ as the category with objects simplicial vector spaces and with morphisms simplicial homotopy classes of simplicial linear maps. In other words, the hom spaces in $\hSVect$ are given by $\pi_0\IHom$, in the simplicial sense. For a general construction of this starting from a simplicial category see for example \cite[\href{https://kerodon.net/tag/00LW}{Tag 00LW}]{kerodon}. In this case the construction is simplified by the fact that all simplicial vector spaces (including mapping spaces) are Kan, so homotopy of maps is an equivalence relation which we denote by $\sim$. See for example \cite[Lemma I.6.1]{GoerssJardine2009}, \cite[\href{https://kerodon.net/tag/00HC}{Tag 00HC}]{kerodon} and \cite[\href{https://kerodon.net/tag/00M0}{Tag 00M0}]{kerodon} for more details. By definition, an isomorphism in the homotopy category is precisely a \textbf{homotopy equivalence}, that is a pair of maps $f: \huaV \to \huaW$ and $g: \huaW \to \huaV$ such that $fg \sim id_{\huaW}$ and $gf \sim id_{\huaV}$. In this case we write $\huaV \simeq \huaW$. The tensor product also descends in a straightforward way to the homotopy category.

Going back to the mapping space $\IHom$, this is right adjoint to the tensor product, in the enriched sense, as for simplicial sets in \ref{prop:enriched-tensor-hom-sset}. We write the natural isomorphism explicitly, as we will need it later on.

\begin{proposition}[Enriched tensor-hom adjunction for $\SVect$]\label{prop:tensor-hom-svect}
    For any $\huaU, \huaV, \huaW$ in $\SVect$, there are natural isomorphisms of simplicial vector spaces
    \begin{equation}\label{eq:tensor-hom-svect}
        \IHom(\huaU \otimes \huaV, \huaW)  
        \cong \IHom(\huaV, \IHom(\huaU, \huaW))
        \cong \IHom(\huaU, \IHom(\huaV, \huaW)).
    \end{equation}
\end{proposition}

\begin{proof}
    We write the $\SVect$-natural isomorphism between the first two spaces, the other one is obtained by precomposing with the symmetry isomorphism.

    Any $f \in \IHom(\huaU \otimes \huaV, \huaW)_m$ is written in components $f_i: \huaU_i \otimes \huaV_i \otimes \Delta^m_i \to \huaW_i$. 
    Write the $i$-simplex $t \in \Delta^l_i$ in terms of the unique nondegenerate $l$-simplex $E_l = 012\dots l \in \Delta^l_l$, as $t = s_I d_J E_l$, for some multi-indices $I,J$ such that $|I| - |J| = i - l$. 
    Then $\rho_{\huaU, \huaV, \huaW}(f) \in \IHom(\huaV, \IHom(\huaU, \huaW))_m$ is given in components by
    \begin{equation}\label{eq:tensor-hom-svect-def-rho}
        \begin{split}
            \rho_{\huaU, \huaV, \huaW}(f)_l: &\huaV_l \otimes \Delta^m_l \longrightarrow \SVect(\huaU \otimes \Delta^l, \huaW)\\
            &\rho_{\huaU, \huaV, \huaW}(f)_l(v,r)(u,t) = f_i(u, s_I d_Jv, s_I d_J r) \in \huaW_i,
        \end{split}
    \end{equation}
    for any $(v,r) \in \huaV_l \otimes \Delta^m_l$, and any $(u,t) \in \huaU_i \otimes \Delta^l_i$ with $t=s_Id_JE_l$.
    Equivalently, when seeing an element $(v,r) \in \huaV_l \otimes \Delta^m_l$ as a simplicial linear map $(v,r) : \Delta^l \to \huaV \otimes \Delta^m$, we can write $\rho_{\huaU, \huaV, \huaW}(f) \in \IHom(\huaV, \IHom(\huaU, \huaW))_m$ in components by $\rho_{\huaU, \huaV, \huaW}(f)_l(v,r) = f \circ (id_{\huaU} \otimes (v,r))$. 

    In the other direction we have maps with components 
    \begin{equation*}
        g_l: \huaV_l \otimes \Delta^m_l \to \SVect(\huaU \otimes \Delta^l, \huaW),
    \end{equation*}
    and we define 
    \begin{equation}\label{eq:tensor-hom-svect-def-tau}
        \begin{split}
            \tau_{\huaU, \huaV, \huaW}(g)_l: &\huaU_l \otimes \huaV_l \otimes \Delta^m_l \longrightarrow \huaW_l\\
            &\tau_{\huaU, \huaV, \huaW}(g)_l(u, v, r) = g_l(v,r)(u, E_l) \in \huaW_l,
        \end{split}
    \end{equation}
    where $E_l \in \Delta^l_l$ is the unique nondegenerate simplex therein, as before.
    Equivalently, we can write $\tau(g) = ev \circ (id_{\huaU} \otimes g)$, where $ev: \huaU \otimes \Delta^l \otimes \SVect(\huaU \otimes \Delta^l, \huaW) \to \huaW$ is the canonical evaluation map, which is simplicial \cite[\S I.5]{GoerssJardine2009} and linear. 
\end{proof}

From this, the usual adjunction of the underlying category can be recovered by taking the simplicial level 0 in \eqref{eq:tensor-hom-svect}. Analogously, the adjunction descends to the homotopy category by applying the functor $\pi_0$ to \eqref{eq:tensor-hom-svect}.

\begin{corollary}[Tensor-hom adjunction for $\hSVect$]\label{cor:tensor-hom-hsvect}
    For any $\huaU, \huaV, \huaW$ in $\SVect$, there are natural isomorphisms
    \begin{equation}\label{eq:tensor-hom-hsvect}
        \hSVect(\huaU \otimes \huaV, \huaW) \cong \hSVect(\huaU, \IHom(\huaV, \huaW)) \cong
        \hSVect(\huaV, \IHom(\huaU, \huaW)).
    \end{equation}
\end{corollary}

\subsubsection{The DG-category of non-negative chain complexes and its homotopy category}\label{sec:monoidal-struct-chains}

Consider the category of (bounded but not necessarily non-negative) chain complexes $\Ch(\Vect)$. Its symmetric monoidal structure is given by the tensor product defined degreewise as
\begin{equation*}
    (A \otimes B)_m = \bigoplus_{p+q=m} A_p \otimes B_q,
\end{equation*}
for any chain complexes $(A, \partial)$, $(B, \partial)$, with differential defined on homogeneous elements by $\partial(a \otimes b) = \partial a \otimes b + (-1)^{|a|} a \otimes \partial b$, where $|a|$ denotes the degree of $a$.
It is important to note that with these definitions, the symmetry isomorphism $A \otimes B \to B \otimes A$ has a sign that follows the Koszul sign rule: $a \otimes b \mapsto (-1)^{|a||b|}b \otimes a$.

The \textbf{mapping complex} $\IHom(A,B)$ is the chain complex with degree $m$ the space of degree $m$ maps of the underlying graded vector spaces 
\begin{equation*}
    \IHom(A,B)_m = \{ f: A_\bullet \to B_{\bullet + m} \},
\end{equation*}
and differential $\partial f = \partial \circ f - (-1)^{|f|} f \circ \partial$. 
Clearly this is a chain complex of vector spaces, as every level is a vector space and the differential is linear. With this, and the obvious composition rule, $\Ch(\Vect)$ becomes a category enriched over itself, also commonly known as a DG-category. See e.g. \cite[\href{https://kerodon.net/tag/00ND}{Tag 00ND}]{kerodon}.

Note that the mapping complex is generally non-zero in both positive and negative degrees and that the space of chain maps is given by the space of 0-cycles:
\begin{equation*}
    \Ch(\Vect)(A, B) = \{f: A \to B \mid \partial f = f \partial \} = \ker \partial_0.
\end{equation*}
Furthermore, two chain maps $f, g:A \to B$ are homotopic if and only if they differ by a 1-boundary $h \in \Img(\partial_1)$ of the mapping complex, i.e. $f - g = \partial h + h \partial$. This has the immediate consequence that homotopy of chain maps is an equivalence relation which is compatible with composition. Hence we define the \textbf{homotopy category} of chain complexes (see e.g. \cite[\href{https://kerodon.net/tag/00NM}{Tag 00NM}]{kerodon}) as the category $\hCh$ with objects chain complexes of vector spaces and morphisms the chain maps up to homotopy, that is $\hCh(A,B) = H_0(\IHom(A,B))$, for any two chain complexes $A$ and $B$. 
In particular, the definition of an isomorphism in the homotopy category as an invertible morphism recovers the notion of a \textbf{chain homotopy equivalence}, i.e. a pair of maps $f:A \to B$ and $g:B \to A$ such that $fg \sim id_B$ and $gf \sim id_A$, where $\sim$ denotes the homotopy relation. In this case we will write $A \simeq B$.
The monoidal structure also descends to the homotopy category, because given $f,g: A \to B$ homotopic through $h$ and $f',g':A' \to B'$ homotopic through $k$, there is a homotopy between $f \otimes f'$ and $g \otimes g'$ given by $h \otimes f' + g \otimes k$.\footnote{Observe here that by the Koszul sign rule $(g \otimes k)(a \otimes a')=(-1)^{|k||a|}g(a) \otimes k(a')=(-1)^{|a|}g(a) \otimes k(a')$.}

Consider now the full subcategory of non-negative chain complexes $\Ch_{\ge 0}(\Vect)$. This is a DG-subcategory, but we can also see it as a category enriched over itself by truncating the mapping complex at 0 and replacing its 0-chains by its 0-cycles. That is, we define, for any non-negative $A$ and $B$,
\begin{equation}\label{eq:non-neg-chain-maps}
    \IHom_{\ge 0}(A, B)_i := tr_{\ge 0}\IHom(A, B)_i = \begin{cases}
        \IHom(A,B)_i &\text{ for } i > 0\\
        \ker\partial_0 &\text{ for } i = 0\\
        0 &\text{ for } i < 0.
    \end{cases}
\end{equation}
With this, $\Ch_{\ge 0}(\Vect)$ is a category enriched over itself. It also admits a homotopy category $\hChP$, which is simply the full subcategory of $\hCh$ generated by the non-negative complexes. In other words, its hom-spaces are equivalently defined by 
\begin{equation*}
    \hChP(A,B) = H_0(\IHom_{\ge 0}(A,B)) = H_0(\IHom(A,B)) = \hCh(A,B).
\end{equation*}

We now recall that the tensor-hom adjunction in $\Ch(\Vect)$ upgrades to an adjunction of $\Ch(\Vect)$-categories. It is in fact more practical to show the enriched adjunction first. The usual adjunction then follows by taking the 0-cycles of the mapping complexes. 
This is for example the content of \cite[Exercise 10.8]{Rotman2009} and the incorrectly stated \cite[Exercise 2.7.3]{Weibel1994}, but we did not otherwise find a complete statement in the literature. The proof is quite straightforward, the main difficulty being choosing the right sign conventions so that the signs of the differentials match. 
A proof using triple complexes can be found at \cite{SE:4803304}.
We write here the natural morphisms composing the adjunction, which we will use later in the thesis and leave it to the reader to check that they form a $\Ch(\Vect)$-natural isomorphism.

\begin{proposition}\label{prop:tensor-hom-chains}
    For any $A,B,C \in \Ch(\Vect)$, there are natural isomorphisms of chain complexes
    \begin{equation}\label{eq:tensor-hom-chains}
        \IHom(A \otimes B, C) \cong \IHom(A, \IHom(B, C)) \cong
        \IHom(B, \IHom(A, C)).
    \end{equation}
\end{proposition}

\begin{proof}
The natural isomorphism $\rho_{A,B,C}: \IHom(A \otimes B, C) \to \IHom(A, \IHom(B, C))$ is given in each degree by 
\begin{equation}\label{eq:tensor-hom-chvect-def-rho}
    \begin{split}
        \rho_{A, B, C}(f): &A_{\bullet} \to \IHom_{\bullet+|f|}(B,C)\\
        &\rho_{A, B, C}(f)(a)(b) = f(a \otimes b) \in C_{|a|+|b|+|f|},
    \end{split}
\end{equation}
for any homogeneous $f \in \IHom_{|f|}(A \otimes B, C)$.

Its inverse $\tau_{A,B,C}: \IHom(A, \IHom(B, C)) \to \IHom(A \otimes B, C)$ is given by
\begin{equation}\label{eq:tensor-hom-chvect-def-tau}
    \begin{split}
        \tau_{A, B, C}(g): &(A \otimes B)_\bullet \to C_{\bullet+|g|}\\
        &\tau_{A, B, C}(g)(a \otimes b) = g(a)(b) \in C_{|a|+|b|+|g|},
    \end{split}
\end{equation}
for any homogeneous $g \in \IHom_{|g|}(A, \IHom(B, C))$.
It is worth noting that the other natural isomorphism $\rho^L_{A,B,C}: \IHom(A \otimes B, C) \to \IHom(B, \IHom(A, C))$ obtained by precomposing $\rho$ with the symmetry isomorphism has a sign. In fact, $\rho^L$ is given by
\begin{equation}\label{eq:tensor-hom-chvect-def-rho-left}
    \begin{split}
        \rho^L_{A, B, C}(f): &B_{\bullet} \to \IHom_{\bullet+|f|}(A,C)\\
        &\rho^L_{A, B, C}(f)(b)(a) = (-1)^{|a||b|} f(a \otimes b) \in C_{|a|+|b|+|f|},
    \end{split}
\end{equation}
for any homogeneous $f \in \IHom_{|f|}(A \otimes B, C)$. 
\end{proof}

This restricts to an adjunction of $\Ch_{\ge 0}(\Vect)$-categories on $\Ch_{\ge 0}(\Vect)$, and it induces a tensor-hom adjunction on the respective homotopy categories $\hCh$ and $\hChP$, where the hom-spaces are simply the homology at level 0 of the full mapping complexes. 

\begin{corollary}\label{cor:tensor-hom-positive-chains}
    For any $A,B,C \in \Ch_{\ge 0}(\Vect)$, there are natural isomorphisms of chain complexes
    \begin{equation}\label{eq:tensor-hom-non-neg-chains}
        \IHom_{\ge 0}(A \otimes B, C) \cong \IHom_{\ge 0}(A, \IHom_{\ge 0}(B, C)) \cong
        \IHom_{\ge 0}(B, \IHom_{\ge 0}(A, C)).
    \end{equation}
\end{corollary}

\begin{proof}
    It follows from the previous proposition that if the complexes in \eqref{eq:tensor-hom-chains} are isomorphic, their $0$-truncations obtained by \eqref{eq:non-neg-chain-maps} are isomorphic. That is
    \begin{equation*}
        \IHom_{\ge 0}(A \otimes B, C) \cong \IHom_{\ge 0}(A, \IHom(B, C)) \cong
        \IHom_{\ge 0}(B, \IHom(A, C)).
    \end{equation*}
    To obtain \eqref{eq:tensor-hom-non-neg-chains} it is enough to observe that for any non-negative chain complexes $A, B, C$,
    \begin{equation*}
        \IHom_{\ge 0}(A, \IHom(B, C)) = \IHom_{\ge 0}(A, \IHom_{\ge 0}(B, C)),
    \end{equation*}
    because of $A$ being non-negative. This is clear in any non-zero degree. In degree zero, the only thing to check is that for any chain map $f: A \to \IHom(B, C)$, and any $a\in A_0$, $f(a)$ is a chain map, i.e. $f(a) \in Z_0\IHom(B,C) = \IHom_{\ge 0}(B, C)_0$. But because $f$ is a chain map this follows from $A_{-1}=0$, since $\partial^{\IHom(B,C)}(f(a)) = f(\partial a) = 0$. 
\end{proof}

\begin{corollary}\label{cor:tensor-hom-hchp}
    For any $A,B,C \in \Ch_{\ge 0}(\Vect)$, there are natural isomorphisms
    \begin{equation*}
    \begin{split}
        \hCh(A \otimes B, C) &\cong \hCh(A, \IHom_{\ge 0}(B, C)) \\
        &\cong
        \hCh(B, \IHom_{\ge 0}(A, C)).
    \end{split}
    \end{equation*}
\end{corollary}

\subsubsection{The Dold-Kan correspondence and equivalences}

We recall the Dold-Kan correspondence, which holds in general for simplicial objects in abelian categories. For more details we refer to \cite[\S 8.4]{Weibel1994}, \cite[\S 22]{May1967}, \cite[\S III.2]{GoerssJardine2009}, \cite[\href{https://kerodon.net/tag/00QQ}{Tag 00QQ}]{kerodon}, for example.

\begin{definition}\label{def:Moore-normalized}
    Let $\huaV$ be a simplicial vector space. The \textbf{Moore complex} of $\huaV$ is the non-negative chain complex $C(\huaV)$ with $C(\huaV)_n = V_n$ for any $n\ge 0$ and differential given at each level $n$ by the boundary map
    \begin{equation*}
        \partial_n = \sum_{i=0}^n (-1)^{i} d_i^n.
    \end{equation*}

    The \textbf{normalized complex} of $\huaV$ is the subcomplex $N(\huaV)$ of the Moore complex with
    \begin{equation}\label{eq:normalized}
        N(\huaV)_n = \ker p^n_n = \bigcap_{i=0}^{n-1} \ker d_i^n.
    \end{equation}
    The Moore complex differential restricted to this is $\partial_n = (-1)^n d_n^n$ at each level $n$.
\end{definition}

\begin{remark}\label{rem:norm-cplx-as-ker-or-as-quot}
    Denote by $D_n\huaV$ the subspace of $\huaV_n$ generated by all degenerate $n$-simplices. 
    The composition of the canonical inclusion and the quotient projection defines an isomorphism
    \begin{equation*}
        \ker p^n_n \to \huaV_n/D_n\huaV.
    \end{equation*} 
    This is shown for example in \cite[Thm. III.2.1]{GoerssJardine2009}. Therefore, the normalized complex can alternatively be written as $N(\huaV)_n = \huaV_n/D_n\huaV$ with the differential of the Moore complex. 
\end{remark}

\begin{proposition}[Dold-Kan correspondence]\label{prop:DK-correspondence}
    The normalized complex functor $N: \Simp\Vect \to \Ch_{\ge 0}(\Vect)$ admits a right inverse functor $DK: \Ch_{\ge 0}(\Vect) \to \Simp\Vect$, and the two form an equivalence of categories between the category of simplicial vector spaces and that of non-negative chain complexes. Under this equivalence, simplicial homotopies $h: \huaV \otimes \Delta^1 \to \huaW$ between two simplicial maps $f,g:\huaV \to \huaW$, correspond to chain homotopies $N(h): N(\huaV) \to N(\huaW)[-1]$ between $N(f)$ and $N(g)$. 
\end{proposition}

Since the Dold-Kan correspondence sends homotopic maps in $\SVect$ to homotopic maps in $\Ch_{\ge 0}(\Vect)$, it descends to an equivalence of the homotopy categories as well.

\begin{corollary}\label{cor:DK-htpy-cats}
    The Dold-Kan correspondence induces an equivalence of categories between the homotopy categories $\hSVect$ and $\hChP$.
\end{corollary}

As a first consequence, the normalized complex of a simplicial vector space contains information about its order as a groupoid. 

\begin{definition}
    Let $(E,\partial)$ be a chain complex. We say $(E,\partial)$ has \textbf{amplitude} $(0,n)$ if it is concentrated in degrees 0 to $n$, that is $E_i = 0$ for $i< 0$ or $i>n$.
\end{definition}

\begin{proposition}\label{prop:VS-order-amplitude}
A simplicial vector space $\huaV$ is a $\VS$ $n$-groupoid if and only if $N(\huaV)$ has amplitude $(0,n)$. 
In particular, the order of a simplicial vector space coincides with the maximal degree for which its normalized chain complex is non-zero. 
\end{proposition}

\begin{proof}
By definition, all normalized complexes are already non-negative. 
If $\huaV$ is a $\VS$ $n$-groupoid, then clearly $N(\huaV)_i = \ker p_i^i = 0$ for any $i>n$. Hence $N(\huaV)$ has amplitude $(0,n)$.

For the converse, the theorem of Moore (Theorem \ref{thm:MooreHornFillers}) implies that all the $p^m_j$ are already surjective for all $m\ge 0$, $0 \le j \le m$. 
If $N(\huaV)$ has amplitude $(n,0)$, then $p^m_m: \huaV_m \to \Lambda^m_m(\huaV)$ is injective for all $m \ge n+1$. So they are isomorphisms with inverse the Moore filler $\mu^m_m$. 
By the dimension formula in Corollary \ref{cor:dim-horn-spaces}, all the horns have the same dimension, so the maps $p^m_j\mu^m_m: \Lambda^m_m(\huaV) \to \Lambda^m_j(\huaV)$ (which are surjective, as each of them is the composition of two surjective maps) are isomorphisms for all $m \ge n+1$ and $0\le j \le m$. 
Therefore $p^m_j$ must be an isomorphism for all $m \ge n+1$ and $0\le j \le m$ by the two-out-of-three property.
\end{proof}

\begin{definition}\label{def:order-type}
We define the \textbf{homotopy type} of a simplicial vector space as the maximal degree for which the homology of the normalized complex is non-zero. That is, we say $\huaV$ is of type $n$ --- or an \textbf{$n$-type} --- if $H_m(N(\huaV)) = 0$ for all $m > n$ and $H_n(N(\huaV)) \neq 0$. If $H_i(N(\huaV)) = 0$ for all $i$, then we say $\huaV$ is \textbf{acyclic}.
\end{definition}

\begin{remark}\label{rem:H-isom-Pi}
    It follows from the definition and Proposition \ref{prop:VS-order-amplitude} that a $\VS$ $n$-groupoid can only have homotopy type lower or equal to $n$. 
    Recall (e.g. from \cite[Thm. 22.1]{May1967}, \cite[III.2.7]{GoerssJardine2009}) that for a $\VS$ $n$-groupoid  $H_i(N(\huaV)) \cong \pi_i(\huaV, 0)$, which can then be shown to be isomorphic to the corresponding homotopy group at any other point. Hence the definition of a $\VS$ $n$-type coincides with that of a set-theoretical $n$-type in Definition \ref{def:htpy-equiv-n-type-set}.
\end{remark}

\begin{example}
    Any $\VS$ 1-groupoid of order 1 can be seen as a $\VS$ $n$-groupoid for any $n \ge 1$, but it is not an $n$-groupoid of order $n$ for any $n > 1$. 
    $\VS$ 1-groupoids can be of homotopy type 1, 0 or be acyclic. $\VS$ 2-groupoids can in addition be of type 2.
\end{example}

Homotopy equivalences and weak equivalences for simplicial vector spaces can be defined as in Definition \ref{def:htpy-equiv-n-type-set}, but they have a much simpler characterization, due to the Dold-Kan correspondence.

\begin{proposition}\label{prop:equivalences-in-svect}
    In $\SVect$, the category of simplicial vector spaces, the following are equivalent:
    \begin{enumerate}
        \item $f:\huaV \to \huaW$ is a \textbf{weak equivalence} in the sense that it induces isomorphisms between the homotopy groups at each point and between the sets of connected components. 
        \item $f:\huaV \to \huaW$ is a \textbf{homotopy equivalence} in the sense that it admits a simplicial homotopy inverse.
        \item $N(f): N(\huaV) \to N(\huaW)$ is a \textbf{quasi-isomorphism} of chain complexes, in the sense that it induces an isomorphism between the respective homologies. 
        \item $N(f): N(\huaV) \to N(\huaW)$ is a \textbf{chain homotopy equivalence}, in the sense that it admits a chain homotopy inverse. 
    \end{enumerate}
\end{proposition}

\begin{proof}
    By the Dold-Kan correspondence being an equivalence between the homotopy categories (Corollary \ref{cor:DK-htpy-cats}), (2) and (4) are equivalent.

    By the discussion in \cite[\S 1.4]{Weibel1994}, (3) and (4) are equivalent: In one direction it is clear that chain homotopy equivalences are quasi-isomorphisms. In the other, every chain complex $(A,\partial)$ of vector spaces is split in the sense that $A_n \cong H_n(A) \oplus \Img\partial_{n+1} \oplus \Img\partial_{n}$, which gives a chain homotopy equivalence between $A$ and its homology $H(A)$ that can then be composed on both sides of any quasi-isomorphism $A \to B$ to upgrade it to a chain homotopy equivalence. 

    Finally (1) and (3) are equivalent because by the Dold-Kan correspondence, for any simplicial vector space $\huaV$, $\pi_i(\huaV, 0) \cong H_i(N(\huaV))$ for any $i\ge 0$, as mentioned in Remark \ref{rem:H-isom-Pi}. 
    In addition, change of basepoint can be accounted for by using the isomorphism $\pi_i(\huaV, 0) \to \pi_i(\huaV, p)$ induced by multiplication by the unit $1p$ of each basepoint $p\in \huaV_0$.
\end{proof}

\begin{remark}\label{rem:n-type-model}
    It is important to note that all of the above equivalences preserve homotopy type but not order. 
    For example, any $n$-type $\huaV \in \SVect$ is homotopy equivalent to $DK(H(N(\huaV)))$, which is always a $\VS$ $n$-groupoid of order $n$ regardless of whether $\huaV$ is an $n$-groupoid. 
    This is one of the reasons why we do not restrict our discussion to $\VS$ $n$-groupoids of a certain fixed order, as observed and expanded upon in Remark \ref{rem:category-of-n-gpds}.
\end{remark}

\subsection{Finite Data}\label{sec:FiniteDataVS} 

In this section we improve on the results about the finite data of $n$-groupoid objects from Section \ref{sec:finite-data} in the specific situation of $\VS$ $n$-groupoids. 

To begin with, because the category of vector spaces is complete, all coskeleta exist, and the coskeleton construction in Section \ref{sec:finite-data} defines $\cosk^m$ as a functor. 
By Definition \ref{def:coskeleton-general-site}, we have the adjunctions 
\begin{equation*}
    \sk^m \dashv \tr^m \dashv \cosk^m, 
    \qquad \Sk^m \dashv \Cosk^m.
\end{equation*}

Secondly, as mentioned in Example \ref{ex:MooreFillers1Gpd} and the preceding discussion, the Moore theorem (Proposition \ref{thm:MooreHornFillers}) gives canonical multiplication maps for each $\VS$ $n$-groupoid. 
Recall that, in the case of $\VS$ $n$-groupoids, $(p^{n+1}_k)^{-1} = \mu^{n+1}_k$ for any $0 \le k \le n+1$. Therefore we can define multiplication maps by setting
\begin{equation} \label{eq:vs-ngpd-mul}
    m_k := d^{n+1}_k\mu^{n+1}_k, \quad \forall 0\le k \le n+1. 
\end{equation} 
Note that these $m_k$ depend only on the $n$-truncation of the $\VS$ $n$-groupoid. 
Because of the discussion in Section \ref{sec:finite-data}, for the $m_k$ to actually be $n$-groupoid multiplication maps, we need to check that they define isomorphisms, they are compatible with the degeneracy maps and they satisfy higher associativity. The first and the last properties follow immediately from the fact that the Moore fillers are alternating sums of simplices and degeneracies of their faces. 
We prove the second one in the following lemma.

\begin{lemma}\label{lem:MooreFillers-degeneracies}
Let $\huaV$ be a simplicial vector space. The space of degenerate $m$-simplices $D_m\huaV$ is isomorphic to the horn space $\Lambda^m_k(\huaV)$ for any $0 \le k \le m$. In particular the Moore fillers $\mu^m_k$ are compatible with the degeneracy maps in the sense that
\begin{equation*}
    \mu^m_k p^m_k s_j = s_j, \qquad \forall 0\le k \le m, 0 \le j \le m-1.
\end{equation*}
\end{lemma}

\begin{proof}
By the construction of $\mu^m_k$ in the proof of Proposition \ref{thm:MooreHornFillers}, its image is contained in $D_m\huaV$, therefore it is a right inverse of the restriction $p^m_k|_{D_m\huaV}$, which means $p^m_k|_{D_m\huaV}: D_m\huaV \to \Lambda^m_k$ is surjective.

We now show that the two spaces have the same dimension, so that $p^m_k|_{D_m\huaV}$ must be an isomorphism with inverse $\mu^m_k$. Consider the two splittings of $\huaV_m$ given by the short exact sequences
\begin{equation*}
    \begin{split}
        0\to D_m\huaV \into &\huaV_m \to \huaV_m/D_m\huaV \to 0,\\
        0\to  \ker p^m_m \to &\huaV_m  \xrightarrow{p^m_m} \Lambda^m_m(\huaV) \to 0.
    \end{split}
\end{equation*}
By Remark \ref{rem:norm-cplx-as-ker-or-as-quot}, we have an isomorphism $\huaV_m/D_m\huaV \cong \ker p^m_m$, and the dimension of these spaces is
\begin{equation*}
    \dim (\huaV_m/D_m\huaV) = \dim \ker p^m_m = \dim \huaV_m - \dim \Lambda^m_m(\huaV).
\end{equation*}
Therefore $\dim D_m\huaV = \dim \Lambda^m_m(\huaV)$, which, by Corollary \ref{cor:dim-horn-spaces} is the same as the dimension of any of the other horn spaces $\Lambda^m_k(\huaV)$. 

Hence $p^m_k|_{D_m\huaV}$ is an isomorphism with inverse $\mu^m_k$. Then for any degenerate simplex $s_j v$, $\mu^m_k p^m_k s_j v = s_j v$.
\end{proof}

Because of this, any kind of $(n+1)$-ary multiplication map between $n$-simplices that can be defined on a $\VS$ $n$-groupoid can always be written canonically as a Moore filler. Therefore to reconstruct a $\VS$ $n$-groupoid it is sufficient to have its $n$-truncation, as summarized in the following theorem.

\begin{proposition}[Finite data for VS $n$-groupoids]\label{prop:FiniteDataVS}
    The data of a $\VS$ $n$-groupoid is equivalent to the data of an $n$-truncated simplicial vector space. 
\end{proposition}

\begin{proof}
This follows from the previous discussion, by observing that the canonical multiplications defined by the Moore fillers in \ref{eq:vs-ngpd-mul} satisfy the conditions in Proposition \ref{prop:finite-data-n-gpds}, where 2b is shown in Lemma \ref{lem:MooreFillers-degeneracies}. 
\end{proof}

\begin{remark}
The fact that the multiplications of a $\VS$ $n$-groupoid can always be recovered from its $n$-truncation does not mean that there cannot be more convenient expressions for them. This is the case for the $n$-dual, as we show in Section \ref{sec:ndual-overview}. 
\end{remark}

\begin{remark}\label{ref:VSnGpd-is-n-skeletal}
Let $\huaV$ be a simplicial vector space. By Lemma \ref{lem:MooreFillers-degeneracies}, the $n$-skeleton $\Sk^n(\huaV)$ is always canonically a $\VS$ $n$-groupoid for any $n\ge 0$. 
Hence, if $\huaV$ is a $\VS$ $n$-groupoid then $\huaV \cong \Cosk^{n+1}(\Sk^{n} (\huaV)) \cong \Sk^n(\huaV)$. 
\end{remark}

\begin{remark}\label{rem:category-of-n-gpds}
One may define the category $\SVect_{\le n}^{k-mult}$ as the category of $n$-truncated simplicial vector spaces with $n$-truncated simplicial maps $f': \widebar{\huaW} \to \widebar{\huaV}$ that are multiplicative, which in this case means that $f' \circ m_k^\huaW = m_k^\huaV \circ f'$ for a choice of $0 \le k \le n+1$. As a consequence of the theorem above, this category is equivalent to the category of $\VS$ $n$-groupoids and consequently a model for $(n+1)$-vector spaces. We will however generally avoid restricting our discussion to this category, as we will often want to consider simplicial maps whose target is a $\VS$ $n$-groupoid but whose source is allowed to be a general simplicial vector space, as in the statement of the theorem. See for example the discussion in Section \ref{sec:ndual-overview}. Additionally this is not closed under the monoidal product in $\SVect$, as we show in Theorem \ref{thm:order-of-tensors}, so it is not a monoidal subcategory. 
\end{remark}

\begin{example}\label{ex:MooreFillers2Gpd}
In analogy with Example \ref{ex:MooreFillers1Gpd}, we can describe the canonical multiplications for a $\VS$ $2$-groupoid $\huaV_2 \rightthreearrows \huaV_1 \rightrightarrows \huaV_0$. For our convenience, we recall the box notation from Definition \ref{def:multiplication-maps} for the ternary multiplications of a $2$-groupoid object $\huaX$. That is, we denote each $m_i$ of three triangles fitting in a $(3,i)$-horn by writing the faces of the horn in order of increasing index with an empty box $\square$ at the $i$-th place. This allows us to see immediately which index each face corresponds to and which face is missing. Hence we write, for each $(W,X,Y,Z) \in \huaX_3 \cong \Lambda^3_k(\huaX)$, for any $0 \le k \le 3$,
\begin{equation*}
    \begin{array}{ccc}
        W = m_0(X,Y,Z) =: \square X Y Z,
        &\qquad
        &Y = m_2(W,X,Z) =: W X \square Z,\\
        X = m_1(W,Y,Z) =: W \square Y Z,
        &\qquad
        &Z = m_3(W,X,Y) =: W X Y \square.
    \end{array}
\end{equation*}
Whenever we write a 3-simplex $(W,X,Y,Z) \in \huaX_3 \cong \Lambda^3_k(\huaX)$ as a horn, we mean that the following horn conditions from \eqref{eq:DefHornBdrySpaces} hold whenever defined:
\begin{equation*}
    \begin{array}{ccc}
        d_0 W = d_0 X, &d_1 X = d_1 Y, &d_2 Y = d_2 Z \\
        d_1 W = d_0 Y, &d_2 X = d_1 Z, \\
        d_2 W = d_0 Z.
    \end{array}
\end{equation*}
In the case of a $\VS$ 2-groupoid $\huaV$, the Moore fillers give
\begin{equation}\label{eq:TriMultOverPointAll}
\begin{split}
    &W = \square X Y Z = X - Y + Z + s_0d_0 Y - s_0d_0 Z + s_1d_0 Z - s_1d_1 Z,\\
    &X = W \square Y Z = W + Y - Z - s_0d_1 W + s_0d_0 Z - s_1d_0 Z + s_1d_1 Z,\\
    &Y = W X \square Z = - W + X + Z + s_0d_1 W - s_0d_0 Z + s_1d_0 Z - s_1d_1 Z,\\
    &Z = W X Y \square = W - X + Y - s_0d_1 W + s_0d_2 W - s_1d_2 W + s_1d_2 X,
\end{split}
\end{equation}
for any $(W,X,Y,Z) \in \huaV_3 \cong \Lambda^3_k(\huaV)$.
\end{example}

\section{Higher vector bundles}\label{sec:higher-vb}

We begin with a concrete picture of simplicial vector bundles. 

\begin{definition}
A \textbf{simplicial vector bundle} $\huaV \to \huaG$ is a simplicial object in the category $\VB$ of vector bundles over $\Mfd$ with bundle maps (i.e. smooth fiberwise-linear maps). 
That is, the base $\huaG$ is a simplicial manifold and each level $q_m: \huaV_m \to \huaG_m$ is a vector bundle over $\huaG_m$, and all simplicial maps of $\huaV$ are bundle maps over the respective simplicial maps of $\huaG$. 
We can represent a simplicial vector bundle as a diagram
\[\begin{tikzcd}
	\dots & {\huaV_2} & {\huaV_1} & {\huaV_0} \\
	\\
	\dots & {\huaG_2} & {\huaG_1} & {\huaG_0}
	\arrow["{\widetilde{d}_i}", shift left=3, from=1-1, to=1-2]
	\arrow[shift right=3, from=1-1, to=1-2]
	\arrow[shift left, from=1-1, to=1-2]
	\arrow[shift right, from=1-1, to=1-2]
	\arrow[shift left=2, curve={height=-6pt}, from=1-2, to=1-1]
	\arrow["{\widetilde{s}_i}", shift left=4, curve={height=-6pt}, from=1-2, to=1-1]
	\arrow[shift left=3, curve={height=-6pt}, from=1-2, to=1-1]
	\arrow["{\widetilde{d}_i}", shift left=2, from=1-2, to=1-3]
	\arrow[from=1-2, to=1-3]
	\arrow[shift right=2, from=1-2, to=1-3]
	\arrow["{q_2}", from=1-2, to=3-2]
	\arrow["{\widetilde{s}_i}", shift left=3, curve={height=-6pt}, from=1-3, to=1-2]
	\arrow[shift left, curve={height=-6pt}, from=1-3, to=1-2]
	\arrow["{\widetilde{d}_i}", shift left, from=1-3, to=1-4]
	\arrow[shift right, from=1-3, to=1-4]
	\arrow["{q_1}", from=1-3, to=3-3]
	\arrow["{\widetilde{s}_0}", curve={height=-12pt}, from=1-4, to=1-3]
	\arrow["{q_0}", from=1-4, to=3-4]
	\arrow["{d_i}", shift left=3, from=3-1, to=3-2]
	\arrow[shift right=3, from=3-1, to=3-2]
	\arrow[shift right, from=3-1, to=3-2]
	\arrow[shift left, from=3-1, to=3-2]
	\arrow[shift left=2, curve={height=-6pt}, from=3-2, to=3-1]
	\arrow[shift left=3, curve={height=-6pt}, from=3-2, to=3-1]
	\arrow["{s_i}", shift left=4, curve={height=-6pt}, from=3-2, to=3-1]
	\arrow["{d_i}", shift left=2, from=3-2, to=3-3]
	\arrow[shift right=2, from=3-2, to=3-3]
	\arrow[from=3-2, to=3-3]
	\arrow["{s_i}", shift left=3, curve={height=-6pt}, from=3-3, to=3-2]
	\arrow[shift left, curve={height=-6pt}, from=3-3, to=3-2]
	\arrow["{d_i}", shift left, from=3-3, to=3-4]
	\arrow[shift right, from=3-3, to=3-4]
	\arrow["{s_0}", curve={height=-12pt}, from=3-4, to=3-3]
\end{tikzcd}\]
in which we denote the face maps by $(\widetilde{d}_i, d_i)$ and the degeneracy maps by $(\widetilde{s}_i, s_i)$. Equivalently the bundle projections of each level assemble into a smooth simplicial map $q:\huaV \to \huaG$.
\end{definition}

\begin{definition}\label{def:vb-n-gpd}
    Let $n$ be a positive integer or $\infty$.
    A \textbf{$\VB$ $n$-groupoid} is a simplicial vector bundle $\huaV$ over a Lie $n$-groupoid $\huaG$ that is also a Lie $n$-groupoid.
    A $\VB$ $\infty$-groupoid is also called a \textbf{higher vector bundle}.
\end{definition}

This definition is equivalent to that of defining $\VB$ $n$-groupoids as $n$-groupoid objects in $\VB$ with the Grothendieck pretopology of surjective submersions. 
Surjective submersions between vector bundles are characterized by the following lemma. 

\begin{lemma}\label{lem:surj-subm-VB}
Let $E \to M$, $F \to N$ be vector bundles. A bundle map $(\widetilde{\phi}, \phi) :E \to F$ is a surjective submersion if and only if $\phi$ is a surjective submersion and $\widetilde{\phi}^{fib}: E \to \phi^*F$ is surjective over the identity of $M$. 
\end{lemma}

\begin{proof}
   This follows from the fact that (up to a choice of splitting) for any vector bundle $E$ and any point $e_p\in E|_p$, $T_{e_p}E \cong E_p \oplus T_pM$, and the tangent map $T_{e_p}\widetilde{\phi}: T_{e_p}E \to T_{\widetilde{\phi}(e_p)} F$ splits as $(\widetilde{\phi}^{fib}|_p, T_p\phi)$ after splitting its domain and codomain, as $\widetilde{\phi}$ is linear along the fibers. 
\end{proof}

This elementary lemma characterizes surjective submersions of vector bundles as the \textit{fiberwise surjective} maps whose base map is a surjective submersion. Here we say $(\widetilde{\phi}, \phi) :E \to F$ is \textbf{fiberwise surjective} if for any $x \in M$, the map $\widetilde{\phi}:E|_x \to F|_{\phi(x)}$ is surjective. This is equivalent to the condition in the lemma that $\widetilde{\phi}^{fib}: E \to \phi^*F$ is surjective over the identity of $M$. 

The horn spaces of a higher vector bundle are representable, by the same result as in the case of Lie $\infty$-groupoids. Moreover, they are vector bundles over the corresponding horn spaces of the base. This follows from the more general result in \cite[Lemma 5.3]{HoyoTrentinaglia2024}.

This gives an additional characterization of higher vector bundles which we report from \cite[Prop. 5.4]{HoyoTrentinaglia2024}

\begin{proposition}\label{prop:SVB-is-Kan-fibration}
    Let $\huaV \to \huaG$ be a simplicial vector bundle over the Lie $\infty$-groupoid $\huaG$. 
    The following are equivalent:
    \begin{enumerate}[label=(\roman*)]
        \item $\huaV$ is a Lie $\infty$-groupoid (i.e. a higher vector bundle)
        \item $q:\huaV \to \huaG$ is a set-theoretical Kan fibration.
        \item $q:\huaV \to \huaG$ is a smooth Kan fibration.
    \end{enumerate}
\end{proposition}

Note that as a consequence of this, unlike for simplicial vector spaces, there are simplicial vector bundles which are not $\VB$ $\infty$-groupoids. From now on we will mostly focus on the latter case, that of higher vector bundles, and denote the category of higher vector bundles over a fixed Lie $n$-groupoid $\huaG$, with morphisms the simplicial maps covering the identity on the base, as $\VB^{\infty}_\huaG$. 

\begin{remark}\label{rem:order-of-vb-gpds}
Because of the considerations in Remark \ref{rem:order-of-kan-fibrations}, and Proposition \ref{prop:SVB-is-Kan-fibration}, the order of the total space $\huaV$ as a $\VB$ $n$-groupoid is always at least the order of the base $\huaG$. This excludes the existence of $\VB$ $m$-groupoids over Lie $n$-groupoids of order $n$ with $m < n$.
\end{remark}

\begin{example}[The tangent $\VB$ $n$-groupoid]\label{ex:tangent-groupoid}
Every Lie $n$-groupoid $\huaG$ has a \textbf{tangent} $\VB$ $n$-groupoid $T\huaG$ obtained by simply composing the functor $\huaG$ with the tangent functor. This is
\[\begin{tikzcd}
	\dots & {T\huaG_2} & {T\huaG_1} & {T\huaG_0} \\
	\\
	\dots & {\huaG_2} & {\huaG_1} & {\huaG_0}
	\arrow["{Td_i}", shift left=3, from=1-1, to=1-2]
	\arrow[shift right=3, from=1-1, to=1-2]
	\arrow[shift left, from=1-1, to=1-2]
	\arrow[shift right, from=1-1, to=1-2]
	\arrow[shift left=2, curve={height=-6pt}, from=1-2, to=1-1]
	\arrow["{Ts_i}", shift left=4, curve={height=-6pt}, from=1-2, to=1-1]
	\arrow[shift left=3, curve={height=-6pt}, from=1-2, to=1-1]
	\arrow["{Td_i}", shift left=2, from=1-2, to=1-3]
	\arrow[from=1-2, to=1-3]
	\arrow[shift right=2, from=1-2, to=1-3]
	\arrow["{q_2}", from=1-2, to=3-2]
	\arrow["{Ts_i}", shift left=3, curve={height=-6pt}, from=1-3, to=1-2]
	\arrow[shift left, curve={height=-6pt}, from=1-3, to=1-2]
	\arrow["{Td_i}", shift left, from=1-3, to=1-4]
	\arrow[shift right, from=1-3, to=1-4]
	\arrow["{q_1}", from=1-3, to=3-3]
	\arrow["{Ts_0}", curve={height=-12pt}, from=1-4, to=1-3]
	\arrow["{q_0}", from=1-4, to=3-4]
	\arrow["{d_i}", shift left=3, from=3-1, to=3-2]
	\arrow[shift right=3, from=3-1, to=3-2]
	\arrow[shift right, from=3-1, to=3-2]
	\arrow[shift left, from=3-1, to=3-2]
	\arrow[shift left=2, curve={height=-6pt}, from=3-2, to=3-1]
	\arrow[shift left=3, curve={height=-6pt}, from=3-2, to=3-1]
	\arrow["{s_i}", shift left=4, curve={height=-6pt}, from=3-2, to=3-1]
	\arrow["{d_i}", shift left=2, from=3-2, to=3-3]
	\arrow[shift right=2, from=3-2, to=3-3]
	\arrow[from=3-2, to=3-3]
	\arrow["{s_i}", shift left=3, curve={height=-6pt}, from=3-3, to=3-2]
	\arrow[shift left, curve={height=-6pt}, from=3-3, to=3-2]
	\arrow["{d_i}", shift left, from=3-3, to=3-4]
	\arrow[shift right, from=3-3, to=3-4]
	\arrow["{s_0}", curve={height=-12pt}, from=3-4, to=3-3]
\end{tikzcd}\]
The fact that this is a $\VB$ $n$-groupoid follows from the fact that if $p^m_k: \huaG_m \to \Lambda^m_k(\huaG)$ is a surjective submersion, or respectively a diffeomorphism, then $Tp^m_k: T\huaG_m \to \Lambda^m_k(\huaG)$ is a surjective submersion, or respectively an isomorphism of vector bundles.
\end{example}

\begin{remark}[Relation with $\VS$ $n$-groupoids] \label{rem:vb-vs-gpd}
Simplicial vector spaces can be identified with simplicial vector bundles over the identity groupoid of a point. The same goes for $\VS$ $n$-groupoids and $\VB$ $n$-groupoids over a point. 
More generally, a $\VB$ $n$-groupoid over a Lie 0-groupoid $M$ (i.e. a smooth manifold, by Example \ref{ex:0gpd}), is effectively ``a bundle of $\VS$ $n$-groupoids'' over $M$ in the sense that its simplicial structure is defined pointwise, which makes it completely described by the $\VS$ $n$-groupoid structure of its typical fiber. Observe that for any simplicial manifold $\huaG$ there is an embedding via the total unit 
\begin{equation}\label{eq:VB-total-deg}
    1: \huaG_0 \to \huaG_k
\end{equation} 
from Remark \ref{rem:total-unit}. Then any $\VB$ $n$-groupoid $\huaV\to \huaG$ may be pulled back via $1$ to $\huaG_0$ and this results in a bundle $1^*\huaV$ of $\VS$ $n$-groupoids over $\huaG_0$.\footnote{Note that taking the pullback of a vector bundle over $\huaG_m$ by $1$ is the same as restricting the vector bundle to the totally degenerate elements $\huaG_0 \overset{1}{\into}\huaG_m$.}

Because of this, the multiplication maps of any $\VB$ $n$-groupoid $\huaV \to \huaG$ restrict to the canonical ones given by the Moore fillers in \eqref{eq:vs-ngpd-mul}. In other words, the multiplications of a $\VB$ 1-groupoid over a unit on the base take the form \eqref{eq:VS1gpd-canon-mult-1} and \eqref{eq:VS1gpd-canon-mult-2}. See also \cite[Prop. (p.558)]{Pradines1988}. 
The multiplications of a $\VB$ 2-groupoid over a unit element take the form \eqref{eq:TriMultOverPointAll}. 
\end{remark}

\subsection{\texorpdfstring{$\VB$}{VB} groupoids and their duals}\label{sec:VB1-groupoids-cores-duals}

We now summarize some results in the theory of $\VB$ (1-)groupoids and their duals. References for this are \cite{Pradines1988}, \cite[\S 11.2]{Mackenzie2005}, \cite[\S 3.2]{GraciaSazMehta2017}.

\begin{definition}
    The \textbf{right core} of a $\VB$ 1-groupoid $\huaV \to \huaG$ is the kernel of the source map restricted to the units
    \begin{equation}\label{eq:right-core-def-1}
        1^*\huaV_1 \xrightarrow{\widetilde{d}_0} \huaV_0,
    \end{equation}
    which we denote by $C_R:= 1^*\ker \widetilde{d}_0 \to \huaG_0$.\footnote{Technically in the notation $1^*\ker \widetilde{d}_0$, we are taking the kernel of the map $\widetilde{d}_0: \huaV_1 \to d_0^*\huaV$, but because $1^*d_0^*\huaV = \huaV$, the pullback to the units of this kernel coincides with the kernel of \eqref{eq:right-core-def-1}.}
    Analogously, the \textbf{left core} of a $\VB$ 1-groupoid is the kernel of the target map restricted to the units $C_L := 1^*\ker \widetilde{d}_1 \to \huaG_0$.

    The \textbf{unit annihilator} of a $\VB$ 1-groupoid is the annihilator subbundle 
    \begin{equation*}
    \Ann(\widetilde{1}\huaV_0) \subseteq 1^*\huaV_1 \to \huaG_0    
    \end{equation*}
    of the subbundle $\widetilde{1}\huaV_0 \subseteq 1^*\huaV_1 \to \huaG_0$.
\end{definition}

\begin{lemma}\label{lem:VB1-core-involution-iso}
The left and right core are isomorphic through the involution of $\huaV_1$ defined on each $v\in\huaV_1$ by $-v^{-1}$.
The duals of the left and right core are isomorphic to $\Ann(\widetilde{1}\huaV_0)$. 
\end{lemma}

\begin{proof}
As discussed in Remark \ref{rem:vb-vs-gpd}, the division maps of a $\VB$ groupoid restrict to the ones given in \eqref{eq:VS1gpd-canon-mult-2} over the units. 
Hence, for any element $v \in \huaV_1|_{1p}$ for some $p\in \huaG_0$, its inverse is as described in \eqref{eq:VS1gpd-canon-mult-3}. 
Therefore this involution exchanges the cores: for any $a \in C_R$ and $b \in C_L$, 
\begin{equation*}
\begin{split}
    - a^{-1} &= a - \widetilde{1}\widetilde{d}_1 a \in C_L, \\
    - b^{-1} &= b - \widetilde{1}\widetilde{d}_0 b \in C_R.
\end{split}
\end{equation*}

Consider now the canonically split exact sequences 
\begin{equation}\label{eq:VB1DualCoreSES}
\begin{tikzcd}[ampersand replacement=\&,cramped,row sep=scriptsize]
	0 \& {1^*\ker\widetilde{d}_i} \& {1^*\huaV_1} \& {\huaV_0} \& 0 \\
	\& {\huaG_0} \& {\huaG_0} \& {\huaG_0}
	\arrow[from=1-1, to=1-2]
	\arrow["{j}", shift left, hook, from=1-2, to=1-3]
	\arrow[from=1-2, to=2-2]
	\arrow["{id-\widetilde{1}\widetilde{d}_i}", shift left, from=1-3, to=1-2]
	\arrow["{\widetilde{d}_i}", shift left, from=1-3, to=1-4]
	\arrow[from=1-3, to=2-3]
	\arrow["{\widetilde{1}}", shift left, from=1-4, to=1-3]
	\arrow[from=1-4, to=1-5]
	\arrow[from=1-4, to=2-4]
	\arrow["id", Rightarrow, no head, from=2-2, to=2-3]
	\arrow["id", Rightarrow, no head, from=2-3, to=2-4]
\end{tikzcd}
\end{equation}
for $i=0,1$, where $j$ is the inclusion of the kernel into $1^*\huaV_1$. By reading it right to left, we have that $1^*\ker\widetilde{d}_i \cong 1^*\huaV_1/\widetilde{1}\huaV_0$. Dually, for $i=0,1$, we have the isomorphism
\begin{equation*}
    \Ann(\widetilde{1}\huaV_0) \newrightleftarrows{j^*}{(id - \widetilde{1}d_i)^*} (1^*\ker\widetilde{d}_i)^*,
\end{equation*}
where one direction is given by the dual map of the involution above. 
\end{proof}

Because of this isomorphism, either choice of the left or right core is referred to as \textit{the core} of a $\VB$ groupoid. 

\begin{example}\label{ex:Lie-algebroid}
    In the case of the tangent $\VB$ groupoid of a Lie groupoid $\huaG$, the core is the underlying vector bundle of the Lie algebroid $A$, while the unit annihilator is the conormal bundle to the units $\huaG_0$. 
    The anchor map $A\to T\huaG_0$ is given by the target $Td_1$ if $A$ is the right core, and $Td_0$ if $A$ is the left core. 
    The bracket is induced from that of vector fields on $\huaG_1$ by the isomorphism of the sections of $A$ with vector fields over $\huaG_1$ invariant under right or left translations depending on the choice of core. For example the right translation $TR_g: d_1^*A^r|_g \to T_g\huaG_1$ is given in terms of the groupoid multiplication by multiplying on the right by the zero section:
    \begin{equation*}
        TR_g a = a \cdot 0_g. 
    \end{equation*}
\end{example}

The kind of translation appearing in this example can be defined for any $\VB$ groupoid and it appears in the definition of its dual, which has the dual of the core as its units. 
Because of the partially defined groupoid multiplication, only right core elements can be right translated, while only left core elements can be left translated. As such, the involution $-(\_)^{-1}$ commonly appears in the description of the dual $\VB$ groupoid in terms of right or left core, as for example in \cite[\S 11.2]{Mackenzie2005}.

We write the dual $\VB$ groupoid in terms of the unit annihilator as it appears in \cite{Pradines1988} for now, as this will be the main framework we use when computing higher duals.
Later, in Section \ref{sec:bar-of-1-cota} we will consider the cotangent $\VB$ groupoid in terms of the left core. We generalize this discussion of cores in Section \ref{sec:annihilators-and-dual-kernels} for $\VS$ $n$-groupoids and in Section \ref{sec:VB-DegAnn-Cores} for $\VB$ $n$-groupoids.

\begin{definition}[{\cite[p. 559]{Pradines1988}}]\label{def:VB1Dual}
    Let $\huaV \to \huaG$ be a $\VB$ groupoid. Its \textbf{dual} $\VB$ groupoid $\huaV^{1*}$ is 
\[\begin{tikzcd}[ampersand replacement=\&,cramped,sep=small]
	{\huaV_1^*} \& {\Ann(\widetilde{1}\huaV_0)} \\
	{\huaG_1} \& {\huaG_0}
	\arrow[shift right, from=1-1, to=1-2]
	\arrow[shift left, from=1-1, to=1-2]
	\arrow[from=1-1, to=2-1]
	\arrow[from=1-2, to=2-2]
	\arrow[shift left, from=2-1, to=2-2]
	\arrow[shift right, from=2-1, to=2-2]
\end{tikzcd}\]
with face and degeneracy maps
\begin{equation*}
\begin{array}{ll}
    (\widecheck{d}_0\eta)|_{d_0g}(v) = \eta(0_g \cdot (v - 1d_1v)),  
    &(\widecheck{1}\epsilon)|_{1p}(v) = \epsilon(v) \\
    (\widecheck{d}_1\eta)|_{d_1g}(v) = \eta((v - 1d_0v) \cdot 0_g),
\end{array}
\end{equation*}
for any $\eta \in \huaV_1|_g$, any $\epsilon \in \Ann(\widetilde{1}\huaV_0)|_{p}$ and any $v\in \huaV_1$ with the appropriate basepoint.
The multiplication of $\huaV^{1*}$ is
\begin{equation}\label{eq:VB1Dual-mult}
    \langle \xi \cdot \eta, w \cdot v \rangle = \langle \xi, w \rangle + \langle \eta, v \rangle,
\end{equation}
for any $(\eta, \xi)\in \Lambda^2_1(\huaV^{1*})|_{(g,h)}$ and $(v,w) \in\Lambda^2_1(\huaV)|_{(g,h)}$ over the (2,1)-horn $(g,h)\in \Lambda^2_1(\huaG)$.
\end{definition}

The main defining property of the multiplication of the dual $\VB$ groupoid in \cite{Pradines1988} is that it is the unique one for which the pairing 
\[\begin{tikzcd}[ampersand replacement=\&,cramped,column sep=small,row sep=scriptsize]
	{\langle\cdot, \cdot\rangle:} \& {T^*\huaG_1} \& {\huaV_1} \& {\R \times \huaG_1} \\
	{0:} \& {\Ann(\widetilde{1}\huaV_0)} \& \huaV_0 \& {0\times \huaG_0}
	\arrow["\otimes"{description}, draw=none, from=1-2, to=1-3]
	\arrow[shift left, from=1-2, to=2-2]
	\arrow[shift right, from=1-2, to=2-2]
	\arrow[from=1-3, to=1-4]
	\arrow[shift left, from=1-3, to=2-3]
	\arrow[shift right, from=1-3, to=2-3]
	\arrow[shift left, from=1-4, to=2-4]
	\arrow[shift right, from=1-4, to=2-4]
	\arrow["\otimes"{description}, draw=none, from=2-2, to=2-3]
	\arrow[from=2-3, to=2-4]
\end{tikzcd}\]
is a simplicial bundle map over the identity.\footnote{More precisely, Pradines states this by seeing the pairing as a map $\huaV^{1*} \oplus \huaV \to B^1\R \times G$, and requiring it to be a bilinear $\VB$ groupoid morphism. This is because the tensor product of two $\VB$ 1-groupoid is generally not a $\VB$ 1-groupoid, as a consequence of Theorem \ref{thm:order-of-tensors}. Nevertheless, a bilinear $\VB$-groupoid morphism from the sum is the same as a simplicial bundle map from the tensor.} This is the central example motivating our definition of duals of higher $\VB$ groupoids. In Section \ref{sec:VB1Dual-revisited} we show that our construction is a generalization of this duality for higher $\VB$ groupoids. 

\begin{remark}\label{rem:VB-1dual-def-by-ev-extending-to-simp-pairing}
We argue that the entire groupoid structure of the dual $\VB$-groupoid is defined by the fact that the canonical evaluation pairing extends to a simplicial map. In fact, normalization of the pairing (i.e. commutativity with the unit map of the tensor and that of $B^1\R$, which is the zero map) implies that $\huaV^{1*}_0 = \Ann(\widetilde{1}\huaV_0)$, and that the unit $\widecheck{1}$ is the inclusion. Moreover, the condition that the multiplication is well defined through \eqref{eq:VB1Dual-mult}, by requiring that it is invariant under a change of $(v_g, w_h)$ to a different composable pair $(v'_g, w'_h)$ such that $w_h \cdot v_g = w'_h \cdot v'_g$, is equivalent to the condition that $(\eta, \xi)$ is a $(2,1)$-horn and thus defines the source and target maps. 
\end{remark}

\begin{example}[The cotangent $\VB$ groupoid]
Let $\huaG$ be a Lie 1-groupoid. The \textbf{cotangent} $\VB$ groupoid is the $\VB$ 1-groupoid $T^{1*}\huaG$ defined by
\[\begin{tikzcd}[ampersand replacement=\&,cramped,sep=small]
	{T\huaG_1^*} \& {\Ann(T1T\huaG_0)} \\
	{\huaG_1} \& {\huaG_0}
	\arrow[shift right, from=1-1, to=1-2]
	\arrow[shift left, from=1-1, to=1-2]
	\arrow[from=1-1, to=2-1]
	\arrow[from=1-2, to=2-2]
	\arrow[shift left, from=2-1, to=2-2]
	\arrow[shift right, from=2-1, to=2-2]
\end{tikzcd}\]
with the groupoid structure in Definition \ref{def:VB1Dual}. 

The canonical symplectic structure $\omega_{can}$ of $T^*\huaG_1$ is \textbf{multiplicative} with respect to this groupoid structure. That is, 
\begin{equation*}
    \widecheck{d}_0^* \omega_{can} - \widecheck{d}_1^* \omega_{can} + \widecheck{d}_2^* \omega_{can} = \delta \omega_{can} = 0,
\end{equation*}
where each term is the pullback of $\omega_{can}$ to $\Lambda^2_1(T^{1*}\huaG)$ by one of the face map. In this case $\widecheck{d}_0$ and $\widecheck{d}_2$ are the projections to the components of the horn, while $\widecheck{d}_1$ is the multiplication map \eqref{eq:VB1Dual-mult}. 
The notation $\delta\omega_{can}$ refers to the fact that $\delta = \sum_{i=0}^n (-1)^i d_i^*$ is a differential $\Omega^p(T^{1*}\huaG_{n}) \to \Omega^p(T^{1*}\huaG_{n-1})$ that extends to $p$-forms the one for the differential cohomology of a Lie groupoid. This originated in \cite{BottShulmanStasheff1976}, and for any simplicial manifold $\huaX$, the double complex $\Omega^{\bullet}(\huaX_\bullet)$ with differentials $\delta$ and the de Rham differential $d$ is commonly referred to as the Bott-Shulman-Stasheff complex. See e.g. \cite{AriasAbadCrainic2013,Lesdiablerets,CuecaZhu2023}. 

With this, $T^{1*}\huaG$ is a \textbf{symplectic groupoid}. That is, a Lie groupoid with a multiplicative symplectic form on its space of arrows. 
The cotangent groupoid appeared in \cite{CosteDazordWeinstein1987} as the symplectic Lie groupoid integrating the cotangent Lie algebroid of the linear Poisson manifold $A^*$. In fact the Lie algebroid of $T^{1*}\huaG$ is $T^*A^*$ with the Lie algebroid structure induced by the Poisson structure of $A^*$. 
\end{example}

\subsection{Generalizations of the Dold-Kan correspondence}

The classical Dold-Kan correspondence holds for simplicial vector bundles over a manifold, so a natural question to ask is whether this can be extended to a correspondence for all simplicial vector bundles. 
This is a difficult question, which was answered for $\VB$ 1-groupoids in the work of \cite{GraciaSazMehta2017}, and for general higher $\VB$ $n$-groupoids in \cite{HoyoTrentinaglia2024}. 
We do not need to recall the entire construction here, but the upshot is that by adding extra structure (a normal weakly flat cleavage) to a $\VB$ $m$-groupoid over a Lie $n$-groupoid $\huaG$ (with $m\ge n$), this is equivalent to a representation up to homotopy of $\huaG$ concentrated in degrees $(0,m)$. 
By restricting morphisms of $\VB$ $m$-groupoids to ones that are weakly flat (i.e. preserve the cleavage in a sense), this upgrades to an equivalence of categories. 

Representations up to homotopy for Lie groupoids were introduced in \cite{AriasAbadCrainic2013}, but they also make sense for general simplicial sets as studied in \cite{AriasAbadSchaetz2012}. A more specific definition for higher Lie groupoids is given in \cite{HoyoTrentinaglia2024}. 
For our purposes it is enough to think of a representation up to homotopy of a Lie $n$-groupoid $\huaG$ as a complex of vector bundles over $\huaG_0$, with extra data given by a sequence of maps $(R_i)_{i \ge 0}$ specifying the representation data. 

Starting from a $\VB$ $n$-groupoid with a choice of cleavage, the complex of vector bundles associated to its representation up to homotopy can be obtained by a functor analogous to the normalized functor of the Dold-Kan correspondence, which actually does not depend on the cleavage at all. 
In fact, the construction of this functor only uses the fact that the pullback to the units of a $\VB$ $n$-groupoid is pointwise a $\VS$ $n$-groupoid (see Remark \ref{rem:vb-vs-gpd}) to extend the normalized functor for simplicial vector spaces.
On the other hand, the choice of cleavage is essential to determine the extra data ${R_m}$. We limit our present discussion to that of the normalized complex functor. We will also introduce cleavages in Definition \ref{def:cleavage} for other purposes, but for a more complete treatment see \cite{HoyoTrentinaglia2024}. 

\begin{definition}\label{def:VB-norm-cplx}
    Let $\huaV \to \huaG$ be a simplicial vector bundle over the simplicial manifold $\huaG$. 
    The \textbf{Moore complex} $C(\huaV)$ of $\huaV$ is the non-negative chain complex of vector bundles over $\huaG_0$ with $C(\huaV)_m = 1^*\huaV_m$, for any $m \ge 0$, and differential 
    \begin{equation*}
        \widetilde{\partial}_m = \sum_{i=0}^m (-1)^i \widetilde{d}_i^m.
    \end{equation*}
    The \textbf{normalized complex} $N(\huaV)$ of $\huaV$ is the subcomplex of the Moore complex with 
    \begin{equation*}
        N(\huaV)_m = 1^*\ker\widetilde{p}_m^m =  \bigcap_{i=0}^{m-1} 1^*\ker\widetilde{d}_i^m.
    \end{equation*}
    The normalized complex of the tangent $n$-groupoid of a Lie $n$-groupoid $\huaG$ is also called the \textbf{tangent complex} of $\huaG$ and denoted by $\huaT\huaG = N(T\huaG)$.
\end{definition}

Fiberwise, at each fixed point $p \in G_0$, these complexes follow many properties of the corresponding ones for vector spaces. In particular the homology of the normalized complex of a higher vector bundle is well-defined at any point $p \in \huaG_0$. We denote it by $H_i(\huaV)|_p := H_i(N(\huaV)|_p)$. As such, we have a natural notion of pointwise $m$-types.

\begin{definition}
Let $\huaV \to \huaG$ be a higher vector bundle over a Lie $n$-groupoid $\huaG$. We say that $\huaV$ is a \textbf{pointwise $m$-type} if the homology $H_i(\huaV)|_p$ vanishes at all $p \in \huaG_0$ for $i > m$, but not for $i=m$. 
\end{definition}

\begin{remark}
Note that the homology of a complex of vector bundles is generally not a complex of vector bundles itself, since the pointwise homology we just defined does not generally have constant rank over the whole base manifold.
\end{remark}

\begin{remark}\label{rem:order-vb-gpd-norm-cplx}
By Remark \ref{rem:order-of-vb-gpds}, the order of a $\VB$ $n$-groupoid is the maximum between the highest non-vanishing degree of its normalized complex and the order of the base. 
Clearly, a $\VB$ $n$-groupoid of order $n$ is at most a pointwise $n$-type. 
Unlike for $\VS$ $n$-groupoids, being a pointwise $m$-type does not guarantee the existence of a Morita equivalent $\VB$ $m$-groupoid: consider for example the situation of $\huaV$ being a pointwise $m$-type higher vector bundle over a Lie $n$-groupoid $\huaG$ that is an $n$-type with $n > m$.
\end{remark}

\begin{remark}
    Under the generalization of Dold-Kan in \cite{GraciaSazMehta2017}, with a choice of Ehresmann connection,
    \footnote{We recall this definition in Def. \ref{def:Ehresmann-connection}. See Remark \ref{rem:cleavage-properties} for its relation to the notion of cleavage.} the tangent $\VB$ groupoid of a Lie groupoid $\huaG$ corresponds to the adjoint representation up to homotopy of $\huaG$ relative to the connection. 
    Analogously, the cotangent $\VB$ 1-groupoid 
    corresponds to the coadjoint representation up to homotopy of $\huaG$ relative to the Ehresmann connection.

    By using the further generalization in \cite{HoyoTrentinaglia2024} and applying it to the tangent $\VB$ $n$-groupoid of a Lie $n$-groupoid, it is possible to obtain its adjoint representation. 
    Analogously, by applying it to the $n$-cotangent $\VB$ $n$-groupoids we construct in this thesis, one finds the coadjoint representation up to homotopy of a Lie $n$-groupoid (relative to a cleavage). 
    We plan to study this in future work.

    Note also that the dualization of representations up to homotopy of a Lie $n$-groupoid $\huaG$ of order $n >1$ is not as straightforward as dualizing its chain complex by a shifted dual construction. In fact, a priori, the representation data will dualize to a representation up to homotopy of the Lie $n$-groupoid $\huaG^{op}$. It is clear that this is a Morita equivalent Lie $n$-groupoid to $\huaG$, but this Morita equivalence may not be expressible as a single map by which the representation can be pulled back to $\huaG$. It is not clear to us at this point how to translate a representation up to homotopy of $\huaG^{op}$ to one of $\huaG$ in absence of the inversion map that exists for $n=1$. 
\end{remark}

\subsection{Weak equivalences of \texorpdfstring{$\VB$}{VB} 2-groupoids}

In this section we extend a particular case of \cite[Thm. 3.5]{HoyoOrtiz2020}, which states that a map of $\VB$ 1-groupoids is a Morita map if and only if it is Morita on the bases and it induces a quasi-isomorphism between the normalized complexes:
we show that a map of $\VB$ 2-groupoids over the identity on the base is a weak equivalence if and only if it induces a quasi-isomorphism between the normalized complexes.

We first recall the notion of a cleavage from \cite{HoyoTrentinaglia2021}.\footnote{We use the notion of cleavage from the first version of this preprint, which was changed in later versions. Both are equivalent.} These describe choices of $(n,k)$-horn fillers for a simplicial vector bundle. This is mostly for convenience, as we will need to repeatedly choose fillers in the proof of the theorem and we can do it ``uniformly'' by a choice of cleavage. 

\begin{definition}\label{def:cleavage}
Let $\huaV \to \huaG$ be a higher vector bundle. For any $m \ge 1$ and $0 \le k \le m$, an \textbf{$(m,k)$-cleavage} $c_{m,k}$ is a bundle map 
\begin{equation*}
    c_{m,k}: (p^m_k)^*\Lambda^m_k(\huaV) \longrightarrow \huaV_m
\end{equation*}
that is a section of the $(m,k)$-horn projection $\widetilde{p}^m_k$ over the identity of $\huaG_m$.

An $(m,k)$ cleavage is said to be \textbf{normal} if it is compatible with the degeneracy maps in the sense that for any $m \in \huaV_{m-1}|_g$, and any $0 \le i \le m-1$,
\begin{equation}\label{eq:normal-cleavage}
    c_{m,k}(\widetilde{p}^m_k \widetilde{s}_i v, s_i g) = \widetilde{s}_i v.
\end{equation}
This is equivalent to the sub-bundle of degenerate simplices $D_m\huaV$ being included in the image of $c_{m,k}$.
\end{definition}

\begin{remark}\label{rem:cleavage-properties}
We collect here some properties of cleavages:
\begin{itemize}
    \item Because of the existence of partition functions, normal cleavages for higher vector bundles always exist. This was shown in \cite[Prop. 4.5]{HoyoTrentinaglia2021}.
    \item Cleavages generalize Ehresmann connections on a $\VB$ groupoid: a normal $(1,1)$-cleavage is an Ehresmann connection in the sense of Definition \ref{def:Ehresmann-connection}. 
    \item Because of linearity of $c_{m,k}$, any horn with all faces equal to zero is always filled to the zero section. 
    \item The condition \eqref{eq:normal-cleavage} for a cleavage to be normal is similar to the compatibility condition for higher multiplication maps in Remark \ref{rem:multiplication-comp-with-deg} and item (2b) of Proposition \ref{prop:finite-data-n-gpds}. Indeed, a normal $(m,k)$-cleavage can be seen as defining a choice of smooth multiplication map for $m$-simplices by taking $d_kc_{m,k}$.
    \item Clearly, by the Kan conditions, $\VB$ $n$-groupoids have unique $(m,k)$-cleavages for all $m \ge n+1$ and $0 \le k \le m$.
    \item By the Moore Theorem (Prop. \ref{thm:MooreHornFillers}) and Lemma \ref{lem:MooreFillers-degeneracies}, any simplicial vector space has a set of unique normal cleavages given by the Moore fillers $\mu^n_k$. The images of these coincide with the degenerate subspaces $D_n\huaV$. An alternative proof of this fact is in \cite[Prop. 3.8]{HoyoTrentinaglia2024}. 
\end{itemize}
\end{remark}

\begin{theorem}\label{thm:we-of-VB2gpd-is-qi}
Let $\huaV \to \huaG$ and $\huaV' \to \huaG$ be $\VB$ 2-groupoids over the same base $\huaG$. 
Let $f:\huaV \to \huaV'$ be a simplicial bundle map over the identity of $\huaG$.
The following are equivalent:
\begin{enumerate}[label=(\roman*)]
    \item $f$ is a weak equivalence of simplicial vector bundles.
    \item At any point $p \in \huaG_0$, the map $(1^*f)|_p: 1^*\huaV|_p \to 1^*\huaV'|_p$ is a weak equivalence of simplicial vector spaces.
    \item At any point $p \in \huaG_0$, the map $N(f)|_p: N(\huaV)|_p \to N(\huaV')|_p$ is a quasi-isomorphism of chain complexes.
\end{enumerate}
\end{theorem}

\begin{proof}
The equivalence (ii) $\iff$ (iii) was shown in Proposition \ref{prop:equivalences-in-svect}. 
We now show that (ii) follows from (i). 

First of all, if $f:\huaV \to \huaV'$ is a simplicial bundle map over the identity, then 
\begin{equation*}
    r^f_l: \huaV_l \times_{f, \huaV'_l, d_{l+1}} \huaV'_{l+1} \to 
        \hom(\partial\Delta^l, \huaV) \times_{\hom(\partial\Delta^l, \huaV')} \Lambda^{l+1}_{l+1}(\huaV')
\end{equation*}
is a bundle map over the horn projection $p^{l+1}_{l+1}: \huaG_{l+1} \to \Lambda^{l+1}_{l+1}(\huaG)$, since both fiber products between the base are over the identity on one side. 
The base map $p^{l+1}_{l+1}$ is a surjective submersion for all $l \ge 0$ and a diffeomorphism for $l \ge 2$, since $\huaG$ must have order at most 2 to admit $\VB$ 2-groupoids over it, according to Remark \ref{rem:order-of-vb-gpds}. 
Therefore, by Lemma \ref{lem:surj-subm-VB}, $f$ is a weak equivalence between $\VB$ 2-groupoids if and only if $r^f_l$ is fiberwise surjective for $l=0,1$ and a fiberwise isomorphism for $l=2$. 
In particular for any point $p \in \huaG_0$, 
\begin{equation*}
    r^{1^*f}_l|_p = r^f_l|_{1p}: \huaV_l \times_{\huaV'_{l}} \huaV'_{l+1}|_{1p} \to (\hom\partial\Delta^l, \huaV) \times_{\hom(\partial\Delta^l, \huaV')} \Lambda^{l+1}_{l+1}(\huaV')|_{p^{l+1}_{l+1}1p}
\end{equation*}
is surjective for $l = 0,1$ and an isomorphism for $l=2$. 
Hence (i) implies (ii). 

To show that (ii) implies (i) we proceed level by level and choose normal cleavages $c_{n,k}$ and $c'_{n,k}$ for all $n=1,2$ and $0 \le k \le n$ for $\huaV$ and $\huaV'$, respectively.
\footnote{There are many occurrences in the following proof, where a single filler can be chosen without a specified cleavage, simply by the Kan conditions. 
We choose to fix a cleavage from the beginning for simplicity.
The only step where using a cleavage is important is when we prove injectivity for $l=2$.}

Beginning with $l=0$ we show that if 
\begin{equation*}
    r^{1^*f}_0|_p = r^f_0|_{1p}: \huaV_0 \times_{\huaV'_{0}} \huaV'_{1}|_{1p} \to \Lambda^{1}_{1}(\huaV')|_{p^1_1 1p} \cong \huaV'_0|_{p}
\end{equation*}
is surjective for all $p\in \huaG_0$, then $r^f_0$ is fiberwise surjective, that is, for any $g \in \huaG_1$, 
\begin{equation*}
    r^f_0|_{g}: \huaV_0 \times_{\huaV'_{0}} \huaV'_{1}|_{g} \to \Lambda^{1}_{1}(\huaV')|_{d_0g} \cong \huaV'_0|_{d_0g}
\end{equation*}
is surjective. Equivalently, this proves $r^f_0$ is a surjective submersion by Lemma \ref{lem:surj-subm-VB}. 

\begin{figure}[h]
\begin{adjustbox}{width=\textwidth}
\includegraphics{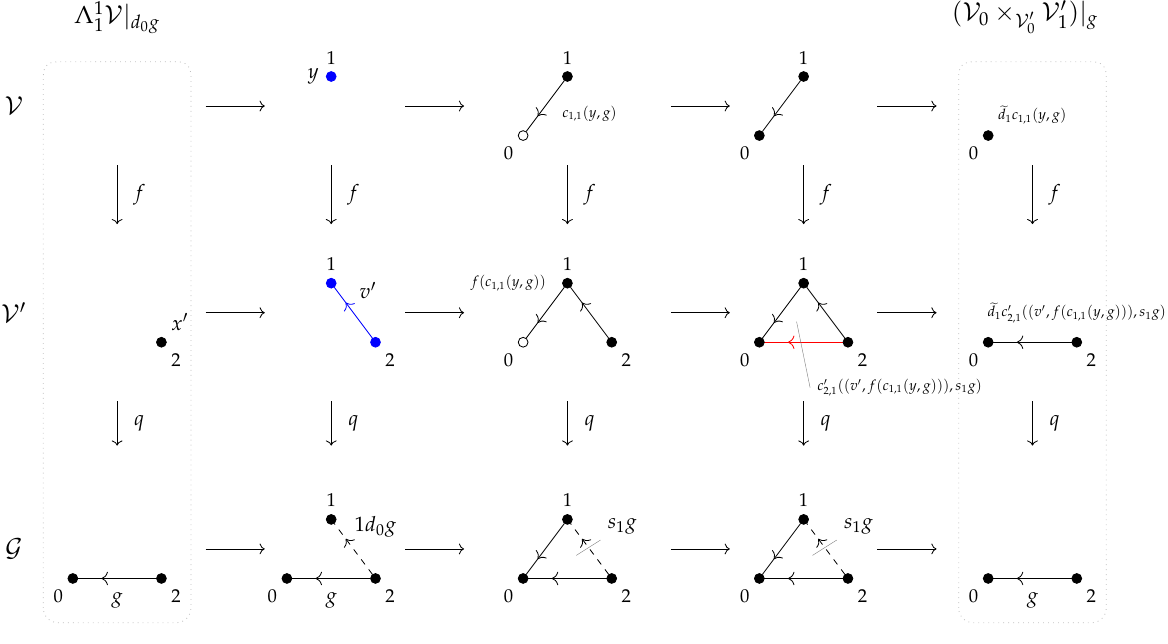}   
\end{adjustbox}
\caption{The construction of a preimage of $x' \in \Lambda^1_1(\huaV)_{d_0g}$ through $r^f_0|_g$. In blue: the existing preimage over the unit edge at $d_0g$, represented by the dashed edge. A white dot indicates a 1-dimensional horn filling, while a red edge indicates a 2-dimensional horn filling. We will use the same conventions in the other figures throughout the proof.}
\label{fig:weak-eq-vb2-gpds0}
\end{figure}

Figure \ref{fig:weak-eq-vb2-gpds0} shows how to construct a preimage of any $x'\in \huaV'_0|_{d_0g}$ through $r^f_0$ over $g$. This consists of the following four steps:
\begin{enumerate}
    \item Since $r^{1^*f}_0|_{d_0g}$ is surjective, there is a preimage $(y,v')$ of $x'$ in $(\huaV_0 \times_{\huaV'_{0}} \huaV'_{1})|_{1d_0g}$. 
    \item In $\huaV$, we fill the $(1,1)$-horn $y$ over $g$ by using the cleavage and get $c_{1,1}(y,g)$.
    \item In $\huaV'$, we get the $(2,1)$-horn filler $c'_{2,1}((v', f(c_{1,1}(y,g))),s_1g)$ over $s_1g$.
    \item Finally, $(\widetilde{d}_1c_{1,1}(y,g), \widetilde{d}_1c'_{2,1}((v', f(c_{1,1}(y,g))),s_1g))$ is a preimage of $x'$ through $r^f_0|_g$. Hence $r^f_0$ is fiberwise surjective.
\end{enumerate}
    
Moving on to $l=1$, we show that if 
\begin{equation*}
    r^{1^*f}_1|_p = r^f_1|_{1p}: \huaV_1 \times_{\huaV'_{1}} \huaV'_{2}|_{1p} \to \huaV_0 \times \huaV_0 \times_{\huaV'_0 \times \huaV'_0}\Lambda^{2}_{2}(\huaV')|_{p^2_2 1p}, 
\end{equation*}
is surjective for all $p\in \huaG_0$, then $r^f_1$ is fiberwise surjective, that is, for any $t \in \huaG_2$, 
\begin{equation*}
    r^f_1|_{t}: \huaV_1 \times_{\huaV'_{1}} \huaV'_{2}|_{t} \to \huaV_0 \times \huaV_0 \times_{\huaV'_0 \times \huaV'_0}\Lambda^{2}_{2}(\huaV')|_{p^2_2t}
\end{equation*}
is surjective. Equivalently, this proves $r^f_1$ is a surjective submersion by Lemma \ref{lem:surj-subm-VB}. 

\begin{figure}[htp]
\begin{adjustbox}{width=\textwidth}
\includegraphics{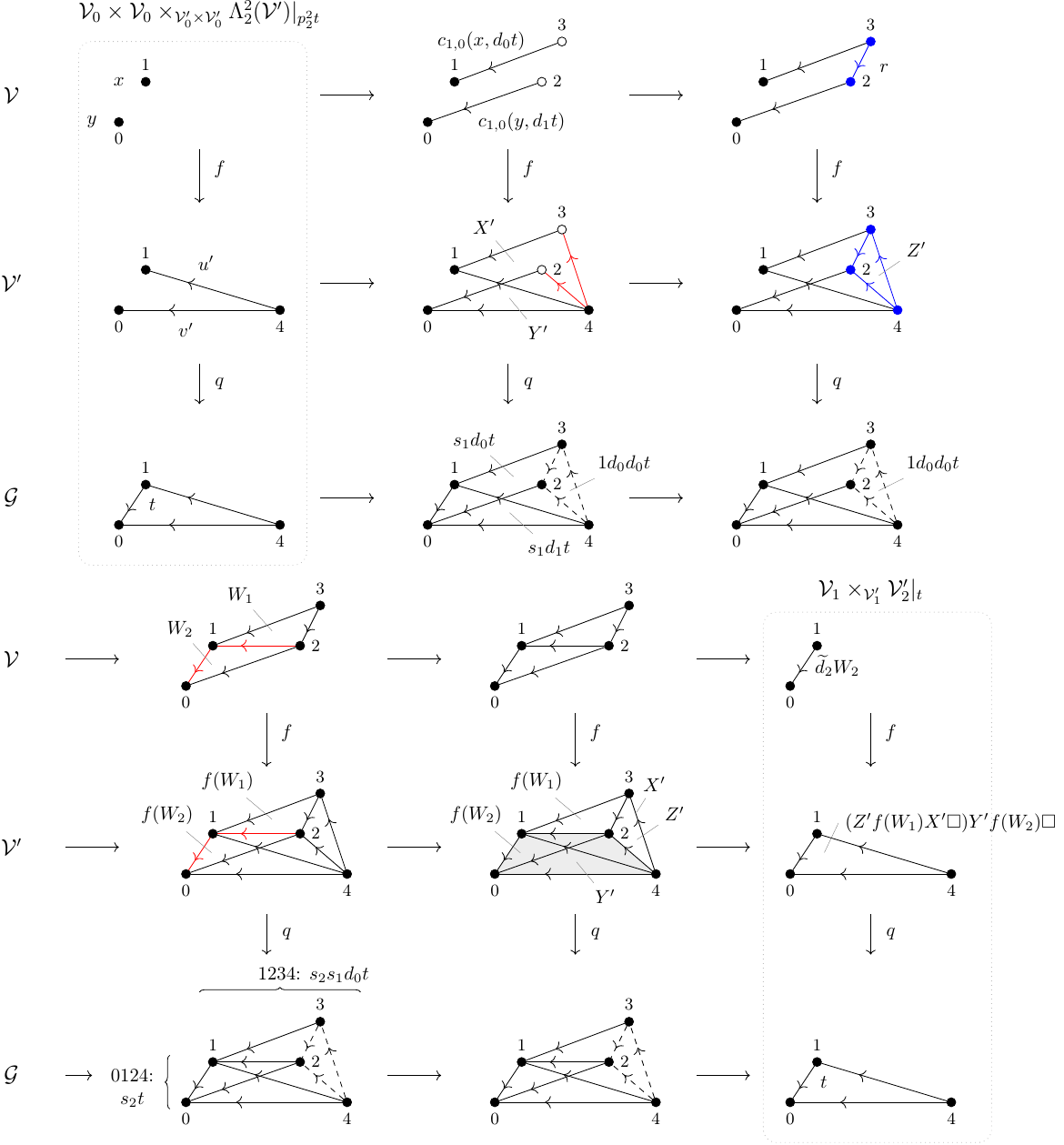}   
\end{adjustbox}
\caption{The construction of a preimage of $((x,y),(u',v')) \in \huaV_0 \times \huaV_0 \times_{\huaV'_0 \times \huaV'_0}\Lambda^{2}_{2}(\huaV')|_{p^2_2t}$ through $r^f_1|_t$. As before, a white dot indicates the added 0-face from a 1-dimensional horn filling, while a red edge indicates the added 1-face from a 2-dimensional horn filling. Analogously, the gray triangles represent 2-faces added as a result of a 3-dimensional horn filling. In blue: the existing preimage over the unit triangle at $d_0d_0t$, represented by the dashed edges.}
\label{fig:weak-eq-vb2-gpds1}
\end{figure}

Figure \ref{fig:weak-eq-vb2-gpds1} shows how to construct a preimage of any element $((x,y),(u',v'))$ in $\huaV_0 \times \huaV_0 \times_{\huaV'_0 \times \huaV'_0}\Lambda^{2}_{2}(\huaV')|_{p^2_2t}$ through $r^f_1$ over $t$. This consists of the following five steps:
\begin{enumerate}
    \item In $\huaV$ we fill the $(1,0)$-horns $x$ and $y$ over $d_1d_0t$ and $d_1d_1t$ respectively, to the edges $c_{1,0}(x,d_0t)$ and $c_{1,0}(y,d_1t)$. Then, in $\huaV'$ we use $(2,0)$-horn fillers to obtain 
    \begin{equation*}
        X' = c_{2,0}((f(c_{1,0}(x,d_0t)), u'), s_1d_0t), 
        \text{ and } Y'= c_{2,0}((f(c_{1,0}(y,d_1t)), v'), s_1d_1t).
    \end{equation*}
    \item Now the element
    \begin{equation*}
        ((\widetilde{d}_0 c_{1,0}(x,d_0t), \widetilde{d}_0 c_{1,0}(y,d_1t)), (\widetilde{d}_0 X', \widetilde{d}_0 Y')) \in \huaV_0 \times \huaV_0 \times_{\huaV'_0 \times \huaV'_0}\Lambda^{2}_{2}(\huaV')|_{1d_0d_0t}
    \end{equation*}
    has a preimage through $r^{1^*f}_1$ over $1d_0d_0t$, which we call $(r, Z')$. 
    \item In $\huaV$ we use $(2,2)$-horn fillers repeatedly to obtain 
    \begin{equation*}
        W_1 = c_{2,2}((r, c_{1,0}(x,d_0t)), s_1d_0t), 
        \quad \text{and } W_2 = c_{2,2}((\widetilde{d}_2W_1, c_{1,0}(y,d_1t)), t).
    \end{equation*}
    At this point, the construction takes place over the union of the tetrahedrons $s_2t$ and $s_2s_1d_0t$ in $\huaG$.
    \item In $\huaV'$ we use $(3,3)$-horn fillers repeatedly, i.e. we multiply triangles to obtain $Z'f(W_1)X'\square$ and $(Z'f(W_1)X'\square)Y' f(W_2) \square$. 
    \item Finally, the preimage of $((x,y),(u',v'))$ through $r^f_1|_t$ is given by
    \begin{equation*}
        (\widetilde{d}_2W_2, (Z'f(W_1)X'\square)Y' f(W_2) \square) \in \huaV_1 \times_{\huaV'_{1}} \huaV'_{2}|_{t}.
    \end{equation*}
    Hence $r^f_1$ is fiberwise surjective.
\end{enumerate}

For $l=2$, we begin as before by showing fiberwise surjectivity of $r^f_2$. That is, we show that if
\begin{equation*}
    r^{1^*f}_2|_p = r^f_2|_{1p}: \huaV_2 \times_{\huaV'_{2}} \huaV'_{3}|_{1p} \to \hom(\partial\Delta^2, \huaV) \times_{\hom(\partial\Delta^2, \huaV')} \Lambda^{3}_{3}(\huaV')|_{p^3_3 1p}, 
\end{equation*}
is bijective, then $r^f_2$ is fiberwise surjective, that is, for any tetrahedron $\tau \in \huaG_3$, 
\begin{equation*}
    r^f_2|_{\tau}: \huaV_2 \times_{\huaV'_{2}} \huaV'_{3}|_{\tau} \to \hom(\partial\Delta^2, \huaV) \times_{\hom(\partial\Delta^2, \huaV')} \Lambda^{3}_{3}(\huaV')|_{p^3_3\tau} 
\end{equation*}
is surjective. 

The intuitive principle guiding the argument is the same as before: the point is to ``translate'' the starting simplicial diagram in $\hom(\partial\Delta^2, \huaV) \times_{\hom(\partial\Delta^2, \huaV')} \Lambda^{3}_{3}(\huaV')$ over $\tau$ to a simplicial diagram over the unit tetrahedron at some point $p\in \huaG_0$. 
This then has a preimage in $(\huaV_2 \times_{\huaV'_{2}} \huaV'_{3})|_{1p}$, which needs to be translated back to a simplicial diagram of the same shape over $\tau$ instead. 
Both of these translations happen by sequences of horn fillings, which are shown in Figure \ref{fig:weak-eq-vb2-gpds2}.\footnote{
See the following Remark \ref{rem:we-of-vb2gpds-with-joins} for more comments on this construction.} 
Our naming convention in the following is that simplices in $\huaV$ are denoted by $x$ and simplices in $\huaV'$ are denoted by $x'$. All of these have subscripts indicating their set of vertices as in Figure \ref{fig:weak-eq-vb2-gpds2}. We implicitly assume that for any subscript $I$, $x'_I = f(x_I)$, when both of them exist in the diagram. Clearly, with this, $d_jx_I = x_{I\backslash \{j\}}$ and the same for $x'$.

\begin{figure}[htp]
\begin{adjustbox}{width=\textwidth}
\includegraphics{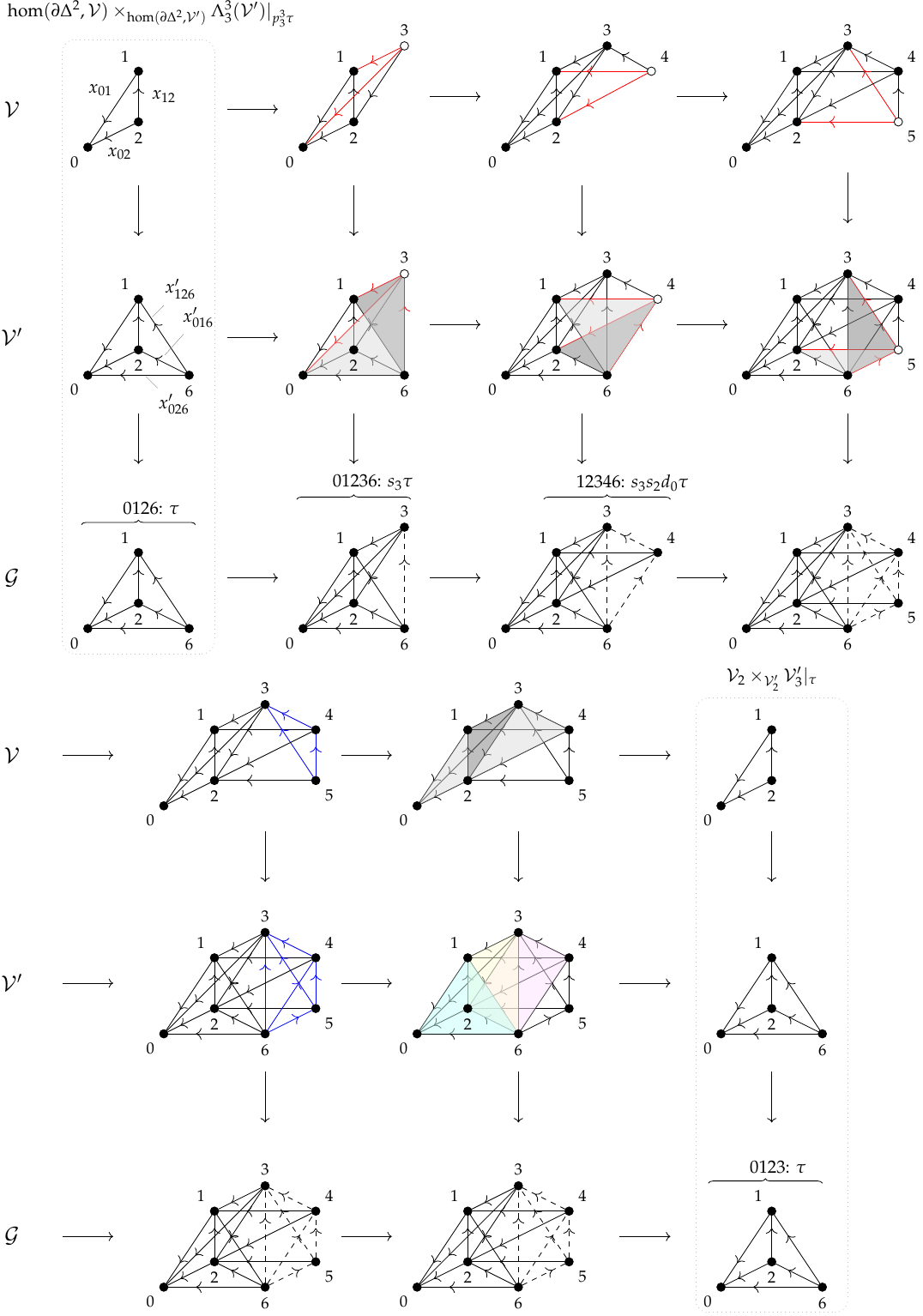}
\end{adjustbox}
\caption{Part one of the construction of a preimage of $((x_{12}, x_{02}, x_{01}), (x'_{126}, x'_{026}, x'_{016}))$ through $r^f_2|_\tau$. As before, a white dot indicates the added 0-face from a 1-dimensional horn filling, while a red edge indicates the added 1-face from a 2-dimensional horn filling. Analogously, the gray triangles represent 2-faces added as a result of a 3-dimensional horn filling. In blue: the existing preimage over the unit tetrahedron at $d_0d_1d_2\tau$, represented by the dashed edges. The tetrahedrons in color represent the ones added to fill 4-simplices.}
\label{fig:weak-eq-vb2-gpds2}
\end{figure}

The construction of a preimage consists of the following steps:
\begin{enumerate}
    \item The starting point is an element 
    \begin{equation*}
        ((x_{12}, x_{02}, x_{01}), (x'_{126}, x'_{026}, x'_{016})) \in \hom(\partial\Delta^2, \huaV) \times_{\hom(\partial\Delta^2, \huaV')} \Lambda^{3}_{3}(\huaV')|_{p^3_3\tau}.
    \end{equation*}
    Because $\huaV'$ is a $\VB$ 2-groupoid, the second component of $r^f_2$ is already an isomorphism. In other words, the $(3,3)$-horn $(x'_{126}, x'_{026}, x'_{016})$ has a unique filler $x'_{0126}$, which determines the 2-simplex $\widebar{x}'_{012} = x'_{126} x'_{026} x'_{016} \square$.\footnote{In fact, because of the isomorphism $\Lambda^3_3(\huaV')\cong \huaV'_3$, the map $r^f_2$ is a surjective submersion if and only if $((d_0, \dots, d_2), f): \huaV_2 \to \hom(\partial\Delta^2, \huaV) \times_{\hom(\partial\Delta^2, \huaV')} \huaV'_2$ is a surjective submersion. This is the content of \cite[Lemma 6.29]{Li2014} as the latter condition is the level 2 condition for $f$ to be an \textit{acyclic fibration} \cite[3.6.1]{Li2014}. These are often called \textit{hypercovers} (see e.g. \cite{Zhu2009,Wolfson2016,BehrendGetzler2017,RogersZhu2020}).}
    This implies that whatever preimage we construct, we must check that the image of its first component through $f$ is actually $\widebar{x}'_{012}$ as obtained by this multiplication. 
    \item By using the appropriate cleavages in $\huaV$ we construct: 
    \begin{enumerate}[label=(\roman*)]
    \item the $(1,0)$-filler $x_{23}$ of $x_2$ over $d_0d_1\tau$, 
    \item the $(2,1)$-filler $x_{023}$ of $(x_{23}, x_{02})$ over $d_1\tau$,
    \item and the $(2,0)$-filler $x_{013}$ of $(x_{03}, x_{01})$ over $d_2\tau$.
    \end{enumerate} 
    With the images of these through $f$ in $\huaV$ and the appropriate cleavage we construct 
    the $(2,0)$-filler $x_{236}$ of $(x'_{26}, x'_{23})$ over $s_1d_0d_1\tau$. 
    We then have: 
    \begin{enumerate}[label=(\roman*)]
    \item the unique $(3,1)$-filler $x'_{0236}$ of $(x'_{236}, x'_{026}, x'_{023})$ over $s_2d_1\tau$, which represents the multiplication $x'_{036}=x'_{236} \square x'_{026} x'_{023}$, 
    \item and the unique $(3,0)$-filler $x'_{0136}$ of $(x'_{036}, x'_{016}, x'_{013})$ over $s_2d_2\tau$, which represents the multiplication $x'_{136}=\square x'_{036}x'_{016}x'_{013}$.
    \end{enumerate}

    \item Again, by using the chosen cleavages in $\huaV$, we construct: 
    \begin{enumerate}[label=(\roman*)]
    \item the $(1,0)$-filler $x_{34}$ of $x_3$ over $1d_0d_1d_2\tau$, 
    \item the $(2,1)$-filler $x_{134}$ of $(x_{34}, x_{13})$ over $s_1d_0d_2\tau$, 
    \item and the $(2,0)$-filler $x_{124}$ of $(x_{14}, x_{12})$ over $d_0\tau$.
    \end{enumerate}
    By applying $f$ we obtain the corresponding fillers in $\huaV'$ and with the $c'_{2,0}$ cleavage we construct 
    the $(2,0)$-filler $x'_{346}$ of $(x'_{46}, x'_{34})$ over $1d_0d_1d_2\tau$. 
    We then have: 
    \begin{enumerate}[label=(\roman*)]
    \item the unique $(3,1)$-filler $x'_{1346}$ of $(x'_{346}, x'_{136}, x'_{134})$ over $s_2s_1d_0d_2\tau$, which represents the multiplication $x'_{146} = x'_{346} \square x'_{136} x'_{134}$, 
    \item and the unique $(3,0)$-filler $x'_{1246}$ of $(x'_{146}, x'_{126}, x'_{124})$ over $s_2d_0\tau$, which represents the multiplication $x'_{246} = \square x'_{146} x'_{126} x'_{124}$.
    \end{enumerate}

    \item One more time, by using cleavages, in $\huaV$ we construct:
    \begin{enumerate}[label=(\roman*)]
    \item the $(1,0)$-filler $x_{45}$ of $x_4$ over $1d_0d_1d_2\tau$, 
    \item the $(2,1)$-filler $x_{245}$ of $(x_{45}, x_{24})$ over $s_1d_0d_1\tau$,
    \item the $(2,0)$-filler $x_{235}$ of $(x_{25}, x_{23})$ over $s_1d_0d_1\tau$.
    \end{enumerate}
    By applying $f$ we obtain the corresponding fillers in $\huaV'$ and with $c'_{2,0}$ we construct
    the $(2,0)$-filler $x'_{456}$ of $(x'_{46}, x'_{45})$ over $1d_0d_1d_2\tau$.
    We then have:
    \begin{enumerate}[label=(\roman*)]
    \item the unique $(3,1)$-filler $x'_{2456}$ of $(x'_{456}, x'_{246}, x'_{245})$ over $s_2s_1d_0d_1\tau$, which represents the multiplication $x'_{256} = x'_{456} \square x'_{246} x'_{245}$, 
    \item and the unique $(3,0)$-filler $x'_{2356}$ of $(x'_{256}, x'_{236}, x'_{235})$ over $s_2s_1d_0d_1\tau$, which represents the multiplication $x'_{356} = \square x'_{256} x'_{236} x'_{235}$.
    \end{enumerate}

    \item Now we have that
    \begin{equation*}
        ((x_{45}, x_{35}, x_{34}), (x'_{456}, x'_{356}, x'_{346})) \in \hom(\partial\Delta^2, \huaV) \times_{\hom(\partial\Delta^2, \huaV')} \Lambda^{3}_{3}(\huaV')|_{1d_0d_1d_2\tau}
    \end{equation*}
    has a unique preimage $(x_{345}, x'_{3456})$ through $r^{1^*f}_2$ over $1d_0d_1d_2\tau$, by assumption.
    
    \item Finally, we can use $x_{345}$ and the three unique $(3,3)$-fillers $x_{2345}$, $x_{1234}$ and $x_{0123}$ to consecutively obtain $x_{234}$, $x_{123}$ and $x_{012}$. This corresponds to the multiplication
    \begin{equation}\label{eq:VB2Gpd-we-lv2-translation-map}
        x_{012} = ((x_{345}x_{245}x_{235}\square) x_{134} x_{124} \square ) x_{023} x_{013} \square.
    \end{equation}
    Meanwhile in $\huaV'$, because it is a $\VB$ 2-groupoid, each tetrahedron is determined by its boundary. So, once we use $f$ to map the fillers $x_{234}$, $x_{123}$ and $x_{012}$, the tetrahedrons $x'_{2346}$, $x'_{1236}$ and $x'_{0126}$ are determined uniquely and in a compatible way with the consecutive filling of the $(4,3)$-horns $x'_{23456}$, $x'_{12346}$ and $x'_{01236}$.
    These 4-simplices correspond to the associativity conditions that imply compatibility with step 1:
    \begin{equation*}
        f(x_{012}) = ((x'_{345}x'_{245}x'_{235}\square) x'_{134} x'_{124} \square ) x'_{023} x'_{013} \square = \widebar{x}'_{012} = x'_{126} x'_{026} x'_{016} \square.
    \end{equation*}
\end{enumerate}

We now show that $r^f_2$ is injective by showing its kernel vanishes at each point $\tau \in \huaG_3$. 
Any element in $\ker r^f_2|_\tau$ is of the form $(h, 0_\tau)$, with $h \in \huaV_2|_{d_3\tau}$ a triangle with zero boundary and with $f(h)=0$.
This can be easily seen by observing that the second component of $r^f_2$ is the bundle isomorphism $p^3_3:\huaV_3 \to \Lambda^3_3(\huaV)$, or alternatively, by observing that $0_{d_0\tau} 0_{d_1\tau} 0_{d_2\tau} \square = 0_{d_3\tau}$ by linearity, and that the tetrahedron corresponding to this multiplication is the unique tetrahedron with zero boundary $0_\tau$.
Then, there is an isomorphism 
\begin{equation*}
    K: \ker r^f_2|_{\tau} \to \ker r^{1^*f}_2|_{d_0d_1d_2\tau}
\end{equation*}
defined by the series of multiplications \eqref{eq:VB2Gpd-we-lv2-translation-map} in step 6 of the above construction and of Figure \ref{fig:weak-eq-vb2-gpds2}.\footnote{This is basically obtained by observing that for the zero boundary and the zero horn, the construction of a preimage above can be repeated by replacing every filler with the zero over that basepoint.} 
More explicitly,
\begin{equation*}
\begin{split}
    K((h, 0_\tau)) &= (\square 0_{s_1d_0d_1\tau}0_{s_1d_0d_1\tau}(\square 0_{s_1d_0d_2\tau} 0_{d_0\tau} (\square  0_{d_1\tau} 0_{d_2\tau} h)), 0_{1d_0d_1d_2\tau}),\\
    K^{-1}((k, 0_{1d_0d_1d_2\tau})) &= (((k 0_{s_1d_0d_1\tau}0_{s_1d_0d_1\tau}\square) 0_{s_1d_0d_2\tau} 0_{d_0\tau} \square ) 0_{d_1\tau} 0_{d_2\tau} \square, 0_\tau).
\end{split}
\end{equation*}
A straightforward computation shows that $K((h, 0_\tau))$ is indeed in the kernel of $r^{1*f}_2|_{1d_0d_1d_2}$, which is an isomorphism. So $K((h, 0_\tau))=0_{1d_0d_1d_2\tau}$ and $h = 0_{d_3\tau}$. 
Therefore $r^2_f$ is an isomorphism.

\end{proof}

\begin{remark}\label{rem:we-of-vb2gpds-with-joins}
An important guiding remark for the proof at each level is that to translate the boundary in $\huaV$ to a boundary over a unit, we are effectively constructing a simplicial diagram resembling a prism (i.e. a homotopy) in $\huaV$ each time.\footnote{In Figure \ref{fig:weak-eq-vb2-gpds2}, the ``prism'' is twisted through itself, but it should still be recognizable: its top and bottom sides are 012 and 345, respectively.} These are however not standard prisms $\Delta^l \times \Delta^1$, as their triangulations are different, so they are not isomorphic to prisms as simplicial sets. 

In \cite[Def. 6.25]{Li2014}, the simplicial diagrams appearing in \eqref{eq:def-we} for the definition of a weak equivalence are described in terms of \textit{joins} of simplicial sets, by writing $r^f_l$ for any simplicial map $f:\huaX \to \huaY$ as
\begin{equation}\label{eq:def-we-joins}
    r^f_l: \hom(\Delta^l \to \Delta^l \star \Delta^0, \huaX \to \huaY) \to \hom(\partial \Delta^l \star \Delta^0, \huaX \to \huaY).
\end{equation}
The symbol $\star$ denotes the join operation, and we have the isomorphisms
\begin{equation*}
    \Delta^{l+1} \cong \Delta^l \star \Delta^0, \text{ and } 
    \Lambda^{l+1}_{l+1} \cong \partial \Delta^l \star \Delta^0.
\end{equation*}
Roughly, a join of a simplicial set $\huaS$ with the 0-simplex $\Delta^0$ consists in adding a vertex to the simplicial set and all possible higher simplices that contain that vertex. This construction is also known as \textit{the cone over} $\huaS$. This is described in more detail in \cite[\S 2.6.1]{Li2014} and \cite{Joyal2002}. 

By observing Figures \ref{fig:weak-eq-vb2-gpds0}, \ref{fig:weak-eq-vb2-gpds1} and \ref{fig:weak-eq-vb2-gpds2} with \eqref{eq:def-we-joins} in mind, it becomes evident that the construction is compatible with joins. In fact, while translating simplicial diagrams to and from the units on the base, all the horn fillings constructing the ``prism'' in $\huaV'$ are ``the join'' with an extra point of a horn filling in $\huaV$ which construct a simplicial diagram which is the join of a ``prism'' with an extra point. 

This leads us to expect the same result of Theorem \ref{thm:we-of-VB2gpd-is-qi} to hold for all weak equivalences of higher vector bundles over Lie $\infty$-groupoids, by abstracting our argument with the help of the join operation. We plan to make this precise in future work.  
\end{remark}

\begin{remark}
The same proof can easily be adapted to the case where $f$ is a simplicial bundle map over an isomorphism of $\huaG$. 
\end{remark}

\chapter{Duals of higher vector spaces}\label{chap:vs-duals}

In this chapter we discuss $n$-duality for simplicial vector spaces. This is the basic case of a higher vector bundle over a point. 
Due to the tensor-hom adjunction in the category of simplicial vector spaces, a canonical candidate for the $n$-dual is already given to us by the mapping space construction. The content of this chapter is an updated version of part of the preprint \cite{RonchiZhu2024}. 

Recall that the theory of simplicial vector spaces is equivalent, both on the nose and up to homotopy, to that of their normalized complexes, as we discussed in Section \ref{sec:review-monoidal-DK}. 
However, by the Eilenberg-Zilber theorem, the normalized complex functor is not an equivalence of monoidal categories on the nose, but only a monoidal functor up to a canonical homotopy equivalence, given by the Eilenberg-Zilber and Alexander-Whitney maps. Therefore, to commute the functor $N$ with the tensor product, one must remember the data given by this homotopy. 
This has far-reaching consequences in our theory of duals, as the Eilenberg-Zilber theorem affects mapping spaces as much as it affects tensor products, making it a crucial point in our discussion. Therefore, we begin this chapter by reviewing this classical theorem, and proposing a new reformulation of it in terms of mapping spaces in Theorem \ref{thm:EZ-thm-hom}. We then discuss precisely how the Eilenberg-Zilber theorem causes the order of a tensor product of $\VS$-groupoids to be the sum of the orders of its factors (at least up to homotopy) in Theorem \ref{thm:order-of-tensors}.
On the other hand, the order of a mapping space depends only on (and is equal to) the order of the target, as we see in Theorem \ref{thm:VSMappingSpace-nGPD}, which is a version of Theorem \ref{thm:hom-ngpd} for simplicial vector spaces. 

In the next section we describe the core of the theory of $n$-duals and $n$-shifted simplicial pairings of simplicial vector spaces. We define the $n$-dual of a simplicial vector space $\huaV$ as the mapping space $\huaV^{n*}:=\IHom(\huaV, B^n\R)$, and review the fact that, by the tensor-hom adjunction, each $n$-shifted simplicial pairing between $\huaV$ and $\huaW$ induces two simplicial maps $\huaV \to \huaW^{n*}$ and $\huaW \to \huaV^{n*}$. With this, we define the canonical $n$-dual pairing as the pairing whose right induced map is the identity $\huaV^{n*} \to \huaV^{n*}$.

We then move to the category of chain complexes, where we define $n$-shifted IM-pairings as chain maps from the tensor product of two chain complexes to the chain complex which is zero at all degrees except for being $\R$ in degree $n$, denoted by $\R[-n]$. This terminology comes from the condition of \textit{infinitesimal multiplicativity}, which is another name for the condition that this is a chain map. We define $n$-shifted duals of chain complexes: for any chain complex $E$, this is precisely the mapping space to $\R[0]$, then shifted by $-n$, or equivalently, $\IHom(E, \R[-n])$. Importantly, this operation can also be performed internally to the category of \textit{non-negative} chain complexes, which requires a truncation at 0 of the mapping space. Then each $n$-shifted IM-pairing can be rewritten in terms of induced maps with target an $n$-shifted dual. 

A central tool in our theory is the construction of the associated IM-pairing to each simplicial pairing: By applying the normalized complex functor and precomposing with the Eilenberg-Zilber map, we can induce an IM-pairing $\lambda_\alpha: N(\huaV) \otimes N(\huaW) \to \R[-n]$ from any simplicial pairing $\alpha:\huaV \otimes \huaW$. 
This allows us to define the \textit{homological} nondegeneracy condition on a simplicial pairing: $\alpha$ is homologically nondegenerate if $\lambda_\alpha$ induces a nondegenerate pairing between the homology groups of the normalized complexes (up to a shift by $n$). 
We then prove what is possibly the main theorem of the theory of $n$-duals (Theorem \ref{thm:ndual-pairing-hom-nondeg}): The fact that the canonical $n$-dual pairing $\huaV^{n*}\otimes \huaV \to B^n\R$ is homologically nondegenerate if and only if $\huaV$ is an $n$-type. 
A perhaps unfortunate consequence of having to pass through the Eilenberg-Zilber map is that even in the case where $\huaV$ is an $n$-groupoid, the $n$-dual pairing is only nondegenerate on the nose for $n=0$ and $n=1$. 
In all other cases the $n$-dual on the simplicial vector space side contains extra information that vanishes up to homotopy with respect to the $n$-shifted dual on the chain complex side. 
This makes the theory more complex but perhaps also more interesting. 

We conclude this section by discussing two consequences of the main theorem. The first is a characterization of the homological nondegeneracy condition in terms of a \textit{homotopical} nondegeneracy condition: in Theorem \ref{thm:hom-nondeg-homotopy-equiv}, we show that an $n$-shifted simplicial pairing is homologically nondegenerate if and only if its induced maps are weak equivalences. The second one is the fact that $n$-duality is reflexive for $n$-types, up to homotopy equivalence: the left induced map of the canonical pairing is a weak equivalence. 

The third section of this chapter is devoted to the explicit computation of $n$-duals. We begin with an overview of the equations describing each level of the $n$-dual as a space of simplicial maps and introduce the necessary terminology. After that, we solve the equations for the 1-dual of a $\VS$ 1-groupoid and show that it coincides with the 1-dual $\VB$-groupoid from Section \ref{sec:VB1-groupoids-cores-duals}. The main result of this section is the computation of the 2-dual of a $\VS$ 2-groupoid, which produces a new object. We solve the equations in detail and compute the normalized complex of this object, to compare it to the 2-shifted dual of the normalized complex of the original space. 

We close the chapter with a discussion of a different possible model of a 2-dual, and that is the simplicial vector space obtained by applying the inverse $DK$ of the normalized complex functor $N$ to the 2-shifted dual of the normalized complex of the original space. 
Mainly we discuss why this route would be problematic in the case of vector bundles and why we chose not to pursue it. 
The main result is that in this case, the canonical 2-dual pairing would be the one obtained by composing the Alexander-Whitney map with the canonical pairing of chains, which is different from the canonical evaluation pairing of $\huaV_2$ as a vector space. 
In fact, the latter is not multiplicative, and as such does not extend to a simplicial pairing $DK(N(\huaV)^*[-2])\otimes \huaV \to B^2\R$.

\section{The Eilenberg-Zilber theorem revisited}

The Eilenberg-Zilber theorem describes how the Dold-Kan correspondence acts on the tensor products on each side. This originally appeared in \cite{EilenbergZilber1953}, \cite[Thm. 2.1]{EilenbergMacLane1954}. See also \cite[\S 29]{May1967}, \cite[Tags \href{https://kerodon.net/tag/00RF}{00RF},\href{https://kerodon.net/tag/00S0}{00S0}]{kerodon}.
More precisely, it describes a lax monoidal structure and an oplax monoidal structure for the normalized complex functor $N$, which assemble into a deformation retraction. 
This implies that the functor $N$ is strong monoidal between the homotopy categories, as we summarize in Corollary \ref{cor:EZ-implies-DK-homotopy-monoidal}. 

\begin{definition}\label{def:EilenbergZilberMap}
    Let $\huaV$ and $\huaW$ be two simplicial vector spaces. The (normalized) \textbf{Eilenberg-Zilber map} is the natural transformation of chain complexes 
    \begin{equation*}
        EZ : N(\huaV) \otimes N(\huaW) \longrightarrow N(\huaV \otimes \huaW)
    \end{equation*}
    defined on elementary tensors $v\otimes w \in N(\huaV)_p \otimes N(\huaW)_q$ as 
    \begin{equation}
        EZ(v\otimes w)=\sum_{(\mu,\nu)\in\Shuf(p,q)} \sign(\mu,\nu) (s_{\nu_q}\dots s_{\nu_1} v) \otimes (s_{\mu_p}\dots s_{\mu_1} w),
    \end{equation}
    where $\Shuf(p,q)$ is the set of $(p,q)$-shuffles from Definition \ref{def:shuffles}.
\end{definition}

This definition makes use of what is sometimes known as the shuffle product between chains (see e.g. \cite[Tag \href{https://kerodon.net/tag/00RF}{00RF}]{kerodon}). Additionally, this map can also be defined in the same way at the level of the Moore complexes, in what is known as the \textit{unnormalized} Eilenberg-Zilber map. It can then be observed that it preserves normalized chains by the simplicial identities. 

\begin{proposition}[Normalized Eilenberg-Zilber Theorem {\cite[Thm. 2.1a]{EilenbergMacLane1954}}]\label{prop:EZ-thm-tensor}
    Let $\huaV$ and $\huaW$ be two simplicial vector spaces. The Eilenberg-Zilber map $EZ$ admits a left inverse $AW$, called Alexander-Whitney map,
    \begin{equation*}
        N(\huaV) \otimes N(\huaW) \newrightleftarrows{EZ}{AW} N(\huaV \otimes \huaW),
    \end{equation*}
    such that 
    \begin{equation*}
        AW \circ EZ = id_{N(\huaV) \otimes N(\huaW)}, 
        \qquad 
        EZ \circ AW \sim id_{N(\huaV \otimes \huaW)}.
    \end{equation*}
    Both of these maps and the homotopy are natural in $\huaV$ and $\huaW$.
    In other words, these maps define a natural chain homotopy equivalence.
\end{proposition}

Since chain homotopy equivalences are precisely the isomorphisms in the category $\hChP$, we have the following result.

\begin{corollary}\label{cor:EZ-implies-DK-homotopy-monoidal}
    The functor $N: \hSVect \to \hChP$ is monoidal, that is, for any simplicial vector spaces $\huaV$ and $\huaW$, 
    \begin{equation*}
        N(\huaV) \otimes N(\huaW) \simeq N(\huaV \otimes \huaW). 
    \end{equation*}
    Moreover, by the Yoneda embedding, for any simplicial vector spaces $\huaU$, $\huaV$ and $\huaW$,
    \begin{equation*}
        \hChP(N(\huaU \otimes \huaV), N(\huaW)) \cong \hChP(N(\huaU) \otimes N(\huaV), N(\huaW)).
    \end{equation*}
\end{corollary}

We can now reformulate the Eilenberg-Zilber theorem in terms of the mapping complex.

\begin{theorem}[Eilenberg-Zilber theorem for mapping complexes]\label{thm:EZ-thm-hom}
    There is a natural chain homotopy equivalence
    \begin{equation*}
        N(\IHom(\huaV, \huaW)) \newrightleftarrows{EZ^H}{AW^H} \IHom_{\ge 0}(N(\huaV), N(\huaW)),
    \end{equation*}
    such that
    \begin{equation*}
        EZ^H \circ AW^H = id_{\IHom_{\ge 0}(N(\huaV), N(\huaW))}, 
        \qquad 
        AW^H \circ EZ^H \sim id_{N(\IHom(\huaV, \huaW))}.
    \end{equation*}
\end{theorem}

\begin{proof}
    In the following we write $\Ch_{\ge 0}:= \Ch_{\ge 0}(\Vect)$. 
    We define the two maps by using the inverse of the Yoneda embedding. This is 
    \begin{equation*}
        \begin{split}
            y^{-1}: &\mathsf{Fun}(\Ch_{\ge 0}, \Set)({\Ch_{\ge 0}}(\_, B), {\Ch_{\ge 0}}(\_,C)) \to \Ch_{\ge 0}(B,C)\\
            &y^{-1}(\eta) = \eta(id_{B}), \text{ where } y(f)_A = f_*.
        \end{split}
    \end{equation*}
    Hence we define natural transformations
    \begin{equation*}
        \Ch_{\ge 0}(N(\_), N(\IHom(\_,\_))) \newrightleftarrows{EZ^H_*}{AW^H_*}
        \Ch_{\ge 0}(N(\_), \IHom_{\ge 0}(N(\_), N(\_))), 
    \end{equation*}
    between these two functors between $(\Simp\Vect^{op})^2 \times \Simp\Vect$ and $\Set$.

    We obtain these in two steps by composing natural transformations in the 2-category of categories. As the first step we define
    \begin{equation}\label{eq:RT}
        \mathsf{Ch}_{\ge 0}(N(\_), N(\IHom(\_,\_)))
        \newrightleftarrows{T}{R}
        \mathsf{Ch}_{\ge 0}(N(\_\otimes\_), N(\_))
    \end{equation}
    as in the following commutative diagram:
    \begin{equation*}
    \begin{adjustbox}{width=\textwidth}
    \begin{tikzcd}[ampersand replacement=\&]
	\& {\mathsf{SVec}^{op}\times \mathsf{SVec}} \&\& {\mathsf{Ch}_{\ge 0}^{op}\times \mathsf{Ch_{\ge 0}}} \\
	\&\&\& {\mathsf{SVec}^{op}\times \mathsf{SVec}} \\
	{(\mathsf{SVec}^{op})^{\times 2}\times \mathsf{SVec}} \&\&\&\&\&\& {\mathsf{Set}} \\
	\&\&\& {\mathsf{SVec}^{op}\times \mathsf{SVec}} \\
	\& {\mathsf{SVec}^{op}\times \mathsf{SVec}} \&\& {\mathsf{Ch}_{\ge 0}^{op}\times \mathsf{Ch_{\ge 0}}}
	\arrow[""{name=0, anchor=center, inner sep=0}, "{N\times N}", from=1-2, to=1-4]
	\arrow["{DK \times DK}"{description}, from=1-4, to=2-4]
	\arrow["{\mathsf{Ch_{\ge 0}}}", curve={height=-18pt}, from=1-4, to=3-7]
	\arrow["{\mathsf{SVec}}"{description}, from=2-4, to=3-7]
	\arrow["\tau"{description}, shift right=2, between={0.2}{0.8}, Rightarrow, from=2-4, to=4-4]
	\arrow["{id \times \underline{\mathrm{Hom}}}", from=3-1, to=1-2]
	\arrow[""{name=1, anchor=center, inner sep=0}, "{id \times \underline{\mathrm{Hom}}}"{description, pos=0.8}, from=3-1, to=2-4]
	\arrow[""{name=2, anchor=center, inner sep=0}, "{\otimes \times id}"{description, pos=0.7}, from=3-1, to=4-4]
	\arrow["{\otimes \times id}"', from=3-1, to=5-2]
	\arrow["\rho"{description}, shift right=2, between={0.2}{0.8}, Rightarrow, from=4-4, to=2-4]
	\arrow["{\mathsf{SVec}}"{description}, from=4-4, to=3-7]
	\arrow[""{name=3, anchor=center, inner sep=0}, "{N\times N}"', from=5-2, to=5-4]
	\arrow["{\mathsf{Ch_{\ge 0}}}"', curve={height=18pt}, from=5-4, to=3-7]
	\arrow["{DK \times DK}"{description}, from=5-4, to=4-4]
	\arrow["\cong"{description}, shift right, between={0.2}{0.8}, Rightarrow, 2tail reversed, from=0, to=1]
	\arrow["\cong"{description}, between={0.2}{0.8}, Rightarrow, 2tail reversed, from=2, to=3]
    \end{tikzcd}
    \end{adjustbox}
    \end{equation*}
    Here the composition from the left corner to the right corner along the top path is the functor on the left hand side of \eqref{eq:RT}, while the one along the bottom path is the functor on its right hand side. Then $T$ is the top-to-bottom composition of natural isomorphisms and $R$ is the bottom-to-top composition.
    The natural isomorphism in the central diamond is given by the tensor-hom adjunction of $\SVect$ (Prop. \ref{prop:tensor-hom-svect}), while the ones in the upper and lower left quadrilaterals are given by the natural isomorphisms establishing $N \circ DK \cong id$ and $DK \circ N \cong id$ in the Dold-Kan correspondence. The triangles to the right commute because the $DK$ functor is full and faithful, by virtue of being a categorical equivalence, so naturally 
    $\SVect(DK(\_), DK(\_)) \cong \Ch_{\ge 0}(\_, \_)$.
    With this, the top-to-bottom composition $T$ and the bottom-to-top composition $R$ are inverses of each other. 

    In the second step we define
    \begin{equation*}
        \mathsf{Ch}_{\ge 0}(N(\_\otimes\_), N(\_)) \newrightleftarrows{E}{A} \mathsf{Ch}_{\ge 0}(N(\_), \IHom_{\ge 0}(N(\_), N(\_))),
    \end{equation*}
    by composing the pullbacks by the natural transformations from the Eilenberg-Zilber Theorem (Prop. \ref{prop:EZ-thm-tensor}) and the tensor-hom adjunction for non-negative chain complexes in Cor. \ref{cor:tensor-hom-positive-chains} as follows:
    \begin{equation*}
    \begin{adjustbox}{width=\textwidth}
    \begin{tikzcd}
	{\mathsf{SVec}^{op}\times \mathsf{SVec}} && {\mathsf{Ch}_{\ge 0}^{op}\times \mathsf{Ch_{\ge 0}}} && {\mathsf{Set}} \\
	\\
	{(\mathsf{SVec}^{op})^{\times 2}\times \mathsf{SVec}} && {(\mathsf{Ch}_{\ge 0}^{op})^{\times 2}\times \mathsf{Ch_{\ge 0}}} && {\mathsf{Ch}_{\ge 0}^{op}\times \mathsf{Ch_{\ge 0}}}
	\arrow["{N\times N}", from=1-1, to=1-3]
	\arrow["{EZ^*}"{description}, shift right=4, shorten <=16pt, shorten >=16pt, Rightarrow, from=1-1, to=3-3]
	\arrow["{\mathsf{Ch_{\ge 0}}}", from=1-3, to=1-5]
	\arrow["\rho_{\ge 0}"{description}, shift right=2, shorten <=15pt, shorten >=15pt, Rightarrow, from=1-3, to=3-5]
	\arrow["{\otimes \times id}", from=3-1, to=1-1]
	\arrow["{N\times N \times N}"', from=3-1, to=3-3]
	\arrow["{AW^*}"{description}, shift right=3, shorten <=16pt, shorten >=16pt, Rightarrow, from=3-3, to=1-1]
	\arrow["{\otimes \times id}"{description}, from=3-3, to=1-3]
	\arrow["{id \times \IHom_{\ge 0}}"', from=3-3, to=3-5]
	\arrow["\tau_{\ge 0}"{description}, shift right=2, shorten <=15pt, shorten >=15pt, Rightarrow, from=3-5, to=1-3]
	\arrow["{\mathsf{Ch_{\ge 0}}}"', from=3-5, to=1-5]
    \end{tikzcd}
    \end{adjustbox}
    \end{equation*}
    Clearly $E \circ A = id$, while $A \circ E \sim id$, by the Eilenberg-Zilber theorem.

    Now we define 
    \begin{equation*}
    \begin{split}
        EZ^H_* &= E \circ T:
        \\
        &\Ch_{\ge 0}(N(\_), N(\IHom(\_,\_))) \rightarrow 
        \Ch_{\ge 0}(N(\_), \IHom_{\ge 0}(N(\_), N(\_))), \\
        AW^H_* &= R \circ A:\\
        &\Ch_{\ge 0}(N(\_), \IHom_{\ge 0}(N(\_), N(\_))) 
        \rightarrow 
        \Ch_{\ge 0}(N(\_), N(\IHom(\_,\_))), \\
    \end{split}
    \end{equation*}
    and 
    \begin{equation*}
        EZ^H = EZ^H_*(id_{N(\IHom(\huaV, \huaW))}), \qquad 
        AW^H = AW^H_*(id_{\IHom_{\ge 0}(N(\huaV), N(\huaW))}). 
    \end{equation*}
    Then, by the commutative diagrams above, $EZ^H_* \circ AW^H_* = id$, which implies that $EZ^H \circ AW^H = id$. 
    By replacing all hom functors with the ones in the respective homotopy categories, and using Corollaries \ref{cor:tensor-hom-hsvect}, \ref{cor:tensor-hom-hchp}, \ref{cor:DK-htpy-cats} and \ref{cor:EZ-implies-DK-homotopy-monoidal}, the above diagrams define a natural isomorphism
    \begin{equation*}
        {\hChP}(\_, N(\IHom(\huaV, \huaW))) \cong {\hChP}(\_, \IHom_{\ge 0}(N(\huaV), N(\huaW))).
    \end{equation*} 
    Because the Yoneda embedding is fully faithful, and thus reflects isomorphisms, $N(\IHom(\huaV, \huaW))$ and $\IHom_{\ge 0}(N(\huaV), N(\huaW))$ are isomorphic through $EZ^H$ and $AW^H$ in the homotopy category. Hence these two maps form a chain homotopy equivalence.

\end{proof}

\subsubsection{Order of tensor products and mapping spaces between \texorpdfstring{$\VS$ $n-$}{VS n-}groupoids}

The fact that the chain homotopy equivalence in the Eilenberg-Zilber Theorem in Proposition \ref{prop:EZ-thm-tensor} is even a deformation retract has consequences on the order and homotopy type of tensor products of $\VS$ $n$-groupoids. The following result expands on \cite[Remark 8.5]{HoyoTrentinaglia2023}, which observes that the tensor product of two $\VS$ 1-groupoids is not generally a $\VS$ 1-groupoid. 
As the authors of this remark write, this is in contrast with what \cite{BaezCrans2004}, \cite{Roytenberg2007a} claimed, in that $\mathsf{2Vec}$ is not a monoidal category, as it is not closed under the induced monoidal product from $\SVect$.
This fact also influences tensor products of $\VB$ $n$-groupoids. In fact, the tensor product of two $\VB$ 1-groupoids is not generally a 1-groupoid again. 

\begin{theorem}\label{thm:order-of-tensors}
Let $\huaV$ be an $n$-groupoid of order $n$ and $\huaW$ be an $m$-groupoid of order $m$. 
Then $\huaV \otimes \huaW$ has order at least $n+m$ and homotopy type at most $n+m$.
\end{theorem}

\begin{proof}
This is because $(N(\huaV) \otimes N(\huaW))_{n+m} = N(\huaV)_n \otimes N(\huaW)_m$, since the other summands vanish by the hypothesis on the order of $\huaV$ and $\huaW$. For the same reason, this is the maximal non-zero degree of the tensor product of the normalized complexes.
By the Eilenberg-Zilber Theorem \ref{prop:EZ-thm-tensor}, $EZ$ admits a left inverse, hence it must be degreewise injective. Therefore  $EZ((N(\huaV) \otimes N(\huaW))_{n+m})$ is a non-zero subspace of $N(\huaV \otimes \huaW)_{m+n}$, which must then be non-zero. 
Since $N(\huaV) \otimes N(\huaW)$ has amplitude $(0,n+m)$, its homology has at most the same amplitude. This is isomorphic to the homology of $N(\huaV \otimes \huaW)$ by the Eilenberg-Zilber homotopy equivalence, hence the homotopy type of $\huaV \otimes \huaW$ is at most $n+m$. 
\end{proof}

\begin{remark}
    Theorem \ref{thm:order-of-tensors} effectively states that, a priori, a tensor product of a $\VS$ $n$-groupoid with a $\VS$ $m$-groupoid is a groupoid of order higher or equal to $n+m$, which is nevertheless homotopy equivalent to an $(n+m)$-groupoid. Models of this $(n+m)$-groupoid can be obtained by $DK(H(N(\huaV) \otimes N(\huaW)))$, as mentioned in Remark \ref{rem:n-type-model}. 
\end{remark}

A similar argument using the deformation retract in Theorem \ref{thm:EZ-thm-hom} can be used to show that a lower bound on the order of a mapping space depends on the order of the target (and on a lower bound on the amplitude of the normalized complex of the source).
We are more interested in showing that it is \textit{exactly} the order of the target, by adapting to simplicial vector spaces the result of Theorem \ref{thm:hom-ngpd}.

\begin{theorem}\label{thm:VSMappingSpace-nGPD}
    Let $\huaU$ be a simplicial vector space, and $\huaV$ be a $\VS$ $n$-groupoid. Then the mapping space $\IHom(\huaU, \huaV)$ is a $\VS$ $n$-groupoid.
\end{theorem}
\begin{proof}
A choice of bases at each level of the simplicial vector space $\huaU$ gives a simplicial set $\huaB$ such that $\huaU \cong \R [\huaB]$.
Since a linear morphism is determined by where the base vectors go, we have
$\SVect(\huaU, \huaV) = \hom(\huaB, \huaV)$.
Furthermore,  $\SVect(\huaU\otimes \Delta^m, \huaV) = \hom(\huaB\times \Delta^m, \huaV)$ because a compatible base for $(\huaU\otimes \Delta^m)_l $ is given by $\huaB_l \times \Delta^m_l$. Thus we have
\begin{equation}
    \IHom_{\SVect}(\huaU, \huaV) = \IHom_{\SSet}(\huaB, \huaV). 
\end{equation}
Then Lemma \ref{thm:hom-ngpd} implies that the underlying simplicial set of $\IHom_{\SVect}(\huaU, \huaV)$ is an $n$-groupoid. Therefore $\IHom_{\SVect}(\huaU, \huaV)$ is a $\VS$ $n$-groupoid. 
\end{proof}

\section{Dual \texorpdfstring{$n$}{n}-groupoids and pairings}\label{sec:ndual-and-pairings}

Dual spaces are commonly defined by the property of having a canonical evaluation pairing with the original space, so we begin by discussing pairings of simplicial vector spaces. 

On the simplicial side, we will consider $n$-shifted pairings, which means they target the Eilenberg-MacLane space $K(\R, n)$, as represented by the $\VS$ $n$-\hspace{0pt}groupoid $B^n\R$. This is the $\VS$ $n$-groupoid consisting of $\R$ on the $n$-th level and $0$ on the levels lower than $n$. Its full simplicial data is 
\begin{equation}\label{eq:def-BnR}
    (B^n\R)_m= \left\{\begin{array}{rl}
         0, & \text{for } m \leq n-1,\\
         \R^{(^m_n)}, & \text{for } m\ge n,  
    \end{array}\right. 
\end{equation} 
and its canonical multiplication maps given by the Moore fillers are 
\begin{equation}\label{eq:defMultiplicationBnR}
    m^{B^n\R}_k((a_i)_{0\le i \neq k \le n}) 
    = \sum_{i=0, \, i\neq k }^{n} (-1)^{i-k+1} a_i,
\end{equation} 
Recall that each of these multiplications is the $k$-th face map at level $n+1$ when writing $(B^n \R)_{n+1}=\Lambda^{n+1}_k(B^n\R)$. In this case, the other face maps $d_{i\neq k}$ are simply projections towards the $i$-th face. The normalized chain complex of $B^n\R$ is precisely the chain complex with $\R$ in degree $n$ and $0$ in all other degrees: $N(B^n\R) = \R[-n]$. 

\begin{definition}\label{def:simp-pairing}
    Let $\huaV$ and $\huaW$ be simplicial vector spaces.    
    We call a linear map $\alpha_n: \huaV_n \otimes \huaW_n \to \R$ an \textbf{$n$-shifted pairing} of $\huaV$ with $\huaW$. 
    Additionally, we say $\alpha_n$ is \textbf{multiplicative} if 
    \begin{equation}\label{eq:nShiftedPairingMultiplicative}
        \alpha_n d_0 - \alpha_n d_1 + \alpha_n d_2 - \dots + (-1)^n \alpha_n d_n = 0, 
    \end{equation}
    with $d_i = d_i^\huaV \otimes d_i^\huaW$.
    We also say $\alpha_n$ is \textbf{normalized} if
    \begin{equation}\label{eq:nShiftedPairingNormalized}
        \alpha_n s_i = 0, \quad \forall 0 \le i < n,
    \end{equation}
    with $s_i = s_i^{\huaV} \otimes s_i^{\huaW}$.

    We call a simplicial linear map $\alpha: \huaV \otimes \huaW \to B^n\R$ an \textbf{$n$-shifted simplicial pairing}.
\end{definition}

\begin{remark}
    In the language of simplicial cohomology of $\huaV \otimes \huaW$, an $n$-shifted pairing is an $n$-cochain: $\alpha_n \in C^n(\huaV \otimes \huaW)$. This is multiplicative if and only if it is closed with respect to the simplicial differential $\delta_n := \sum_{i=0}^n (-1)^i (d_i^n)^*$. In other words, \eqref{eq:nShiftedPairingMultiplicative} can be written simply as $\delta \alpha_n = 0$.
\end{remark} 

\begin{lemma}\label{lem:VS-simpl-n-pairing-mult-norm}
    There is a one-to-one correspondence between $n$-shifted simplicial pairings $\alpha:\huaV \otimes \huaW \to B^n\R$ and $n$-shifted pairings $\alpha_n:\huaV_n \otimes \huaW_n \to \R$ that are multiplicative and normalized. 
    In this case we identify $\alpha$ and $\alpha_n$.
\end{lemma}

\begin{proof}
    This follows from Lemma \ref{lem:finite-data-simp-map}.
    In fact, starting with the simplicial linear map $\alpha:\huaV \otimes \huaW \to B^n\R$, this is determined by its $n$-th level $\alpha_n$, whose compatibility with the face and degeneracy maps is equivalent to \eqref{eq:nShiftedPairingMultiplicative} and \eqref{eq:nShiftedPairingNormalized}, respectively. 
    This is because all levels $\alpha_i$ for $i<n$ must vanish, automatically implying commutativity with the face maps of the $n$-truncation of $\alpha$ and commutativity with the degeneracy maps of its $(n-1)$-truncation. 
    Meanwhile, commutativity with the face maps between levels $n+1$ and $n$ translates directly to the multiplicativity condition \eqref{eq:nShiftedPairingMultiplicative} as in Lemma \ref{lem:finite-data-simp-map}, and commutativity with the degeneracy maps between levels $n$ and $n-1$ translates directly to the normalization conditions \eqref{eq:nShiftedPairingNormalized}.
\end{proof}

We now introduce the dual $\VS$ $n$-groupoid of a simplicial vector space for each $n\ge 0$. The fact that this is a $\VS$ $n$-groupoid follows from Theorem \ref{thm:VSMappingSpace-nGPD} because $B^n\R$ is a $\VS$ $n$-groupoid.

\begin{definition}\label{def:ndual}
    Let $\huaV$ be a simplicial vector space. For each $n\ge 0$, the \textbf{dual $n$-groupoid} of $\huaV$, or \textbf{$n$-dual} for short, is
    \begin{equation}\label{eq:vs-n-dual}
        \huaV^{n*} := \IHom(\huaV, B^n\R),
    \end{equation}
    with face and degeneracy maps as defined in \eqref{eq:mapping-set-maps} and denoted by $\widecheck{d}_i$ and $\widecheck{s}_i$, respectively.

    The \textbf{$n$-dual pairing} is $\langle \cdot, \cdot \rangle = \tau(id_{\huaV^{n*}}): \huaV^{n*} \otimes \huaV \to B^n\R$, where $\tau$ is the adjunction isomorphism defined in \eqref{eq:tensor-hom-svect-def-tau}.
\end{definition}

As a first application of the notion of $n$-dual, by the tensor-hom adjunction in Prop. \ref{prop:tensor-hom-svect} we have that the data of any simplicial $n$-shifted pairing $\alpha: \huaV \otimes \huaW \to B^n\R$ is equivalent to either of two induced simplicial linear maps
\begin{equation*}
    \alpha^l: \huaW \to \huaV^{n*}, \qquad \alpha^{r}: \huaV \to \huaW^{n*}, 
\end{equation*}
which we call the \textbf{left} and \textbf{right induced map}, respectively. These are defined by applying the $\rho$ map in \eqref{eq:tensor-hom-svect-def-rho} with $m=0$. In the notation of \eqref{eq:n-dual-element-in-comp}, where we consider each $n$-simplex $q = s_I d_J E_k \in \Delta^k_n$ as indexing a component of the target element in the $n$-dual at each level $k$, we write for any $w \in \huaW_k$, $v \in \huaV_n$, and any $v' \in \huaV_k$, $w' \in \huaW_n$,
\begin{equation}\label{eq:nShiftedPairingIndMapsDef}
    (\alpha^l(w))^q (v) = \alpha(v, s^{\huaW}_I d^{\huaW}_J w), \qquad (\alpha^r(v'))^q (w') = \alpha(s^{\huaV}_I d^{\huaV}_J v', w'),
\end{equation}
with $I, J$ a pair of multi-indices such that $|I|-|J| = k - n$, so that $s^{\huaW}_I d^{\huaW}_J w\in \huaW_n$ and $s^{\huaV}_I d^{\huaV}_J v' \in \huaV_n$, as expected. 

\begin{example}\label{ex:ndual-pairing-induced-maps}
    The simplicial maps induced by the $n$-dual pairing are
    \begin{equation*}
    \langle \cdot, \cdot \rangle^l: \huaV \to (\huaV^{n*})^{n*}, \qquad \langle \cdot, \cdot \rangle^r = id: \huaV^{n*} \to \huaV^{n*}.
    \end{equation*}
    Here, the right induced map is of course the identity by definition of the $n$-dual pairing as $\tau(id)$. 
    On the other hand, if $\huaV$ is at most an $n$-type, the left induced map provides the homotopy equivalence between $\huaV$ and its double $n$-dual, as we will show in Theorem \ref{thm:ndual-reflexive-uth}. 
\end{example}

\begin{example}\label{ex:VS0dual}
    The 0-dual $\huaV^{0*}$ of any simplicial vector space $\huaV$ is exactly the identity groupoid of $\pi_0(\huaV)^* = H^0(N(\huaV))^*$. In fact, since $\huaV^{0*}$ is a 0-groupoid, it is determined entirely by its level 0. Additionally, the data of a simplicial map $f:\huaV \to B^0\R$ reduces to that of an element $f\in \huaV_0^*$ such that for any $v \in \huaV_1$, $fd_0v = fd_1v$. This is the same as saying that $f$ is constant on each connected component of $\huaV$. Thus it descends to an element $[f] \in \pi_0(\huaV)^*$.
    
    In the particular case that $\huaV$ is a 0-groupoid, i.e. the identity groupoid of the vector space $\huaV_0$, $\huaV^{0*} = \huaV_0^*$, the identity groupoid of the dual of $\huaV_0$. 
    Here the $0$-dual pairing coincides with the usual dual pairing of vector spaces and it is even nondegenerate on the nose. 
    If $\huaV$ is not a 0-groupoid on the nose, but at most a 0-type, then the 0-dual is still homotopy equivalent to $\huaV$, but the 0-dual pairing is only nondegenerate up to homotopy, as we will show in Theorem \ref{thm:ndual-pairing-hom-nondeg}. 
\end{example}

In this example we can see that, despite the fact that every simplicial vector space admits an $n$-dual for any $n\ge 0$, this might a priori come with a loss of information. Our aim in the following discussion is to make such a statement precise and explain why this happens. The main tools at our disposal are the Dold-Kan correspondence and the Eilenberg-Zilber Theorem for mapping spaces. 
Roughly speaking, the problem lies in the fact that the $n$-dual pairing is related to the homotopy equivalence $EZ^H$, and this sees only a truncated version of the $n$-shifted dual in chain complexes. 
We define the latter object in Definition \ref{def:nshifted-dual-chain-cplx} and write this result in Theorem \ref{thm:ndual-pairing-hom-nondeg}.

On the other side of the Dold-Kan correspondence from simplicial pairings, we have shifted pairings of chain complexes, which we call IM-pairings, following \cite{CuecaZhu2023}, where the associated IM-form to a shifted 2-form on a Lie $n$-groupoid was defined using a formula from \cite{Lesdiablerets}. This terminology is in reference to the theory of infinitesimally multiplicative forms and tensors appearing in \cite{BursztynCrainicWeinsteinZhu2004}, \cite{BursztynCabrera2012}, \cite{BursztynDrummond2019}.

\begin{definition}\label{def:IMPairing}
    Let $(A, \partial^A)$, $(B, \partial^B)$ be chain complexes concentrated in non-negative degrees. 
    An \textbf{$n$-shifted IM-pairing} $\lambda$ between $A$ and $B$ is a chain map $\lambda: A\otimes B \to \R[-n]$. In other words, $\lambda$ is a linear map
    \begin{equation*}
        \lambda: (A \otimes B)_n = \bigoplus_{i=0}^n A_i \otimes B_{n-i} \to \R,
    \end{equation*}
    that is \textbf{infinitesimally multiplicative}, i.e. it satisfies 
    \begin{equation}\label{eq:IMPairingDef}
        \lambda^{\alpha}(\partial u, w) + (-1)^{i+1} \lambda^{\alpha}(u, \partial w) = 0,
    \end{equation}
    for any $u\in A_{i+1}$ and $w \in B_{n-i}$.
\end{definition}

In chain complexes we have a natural notion of shifted duality, because of the existence of the $m$-shift operator $[m]$, which is defined for any chain complex $(A, \partial)$, by $A[m]_i = A_{i+m}$ with differential $\partial^{[m]} = (-1)^m \partial$. 
Because the $n$-shifted dual we are about to define is generally not concentrated in non-negative degrees, we also introduce its truncation as the $n$-shifted dual inside $\Ch_{\ge 0}(\Vect)$.

\begin{definition}\label{def:nshifted-dual-chain-cplx}
    Let $(A, \partial)$ be a non-negative chain complex. The \textbf{$n$-shifted dual} $A^*[-n]$ of $A$ is 
    \begin{equation}\label{eq:def-nshifted-dual-chain-cplx}
        A^*[-n] := \IHom(A, \R[-n]) = \IHom(A, \R[0])[-n],
    \end{equation}
    with differential $\partial^*_i = (-1)^{i+1}(\partial_{n-i+1})^t$ with $\partial_{n-i+1}: A_{n-i+1} \to A_{n-i}$ and
    \begin{equation*}
        (\partial_{n-i+1})^t: (A^*[-n])_i = A^*_{n-i}\to A^*_{n-i+1}=(A^*[-n])_{i-1}
    \end{equation*}
    its dual.  

    The \textbf{non-negative $n$-shifted dual} $A^*[-n]_{\ge 0}$ is the truncation 
    \begin{equation}\label{eq:def-nonneg-nshifted-dual-chain-cplx}
        A^*[-n]_{\ge 0} := \IHom_{\ge 0}(A, \R[-n]) = tr_{\ge 0}(\IHom(A, \R[0])[-n]),
    \end{equation}
    with the same differential. Here $tr_{\ge 0}$ is the truncation functor defined in \eqref{eq:non-neg-chain-maps}.
\end{definition}

\begin{remark}
    The sign in the differential of the $n$-shifted dual comes from the definition of the mapping complex in Section \ref{sec:monoidal-struct-chains}. 
\end{remark}

\begin{remark}\label{rem:DK-n-shift-dual-appendix}
    The existence of non-negative $n$-shifted duals and Theorem \ref{thm:EZ-thm-hom} raise the question of what the $\VS$ $n$-groupoid $DK(N(\huaV)^*[-n]_{\ge 0})$ should represent and if this is a good notion of $n$-dual of $\huaV$. We discuss this in Section \ref{sec:appendix-DK-N-nshifted-dual}.
\end{remark}

By the tensor-hom adjunction in chain complexes (Proposition \ref{prop:tensor-hom-chains}), any IM-pairing $\lambda: A \otimes B \to \R[-n]$ also induces two chain maps 
\begin{equation}\label{eq:IMPairingDef-induced-maps}
    \begin{split}
        &\lambda^l: B \to A^*[-n], \qquad \lambda^r: A \to B^*[-n]\\
        &\lambda^l(w)(v) := \rho^L\lambda(w)(v) = (-1)^{(n-j)j} \lambda(v, w), \text{ for } w\in B_{j}, v \in A_{n-j},\\
        &\lambda^r(v)(w) := \rho\lambda(w)(v) = \lambda(v, w), \text{ for } v \in A_i, w\in B_{n-i},
    \end{split}
\end{equation}
which we call \textbf{left} and \textbf{right induced map} respectively.

As for usual vector spaces, an advantage of considering the maps induced by an IM-pairing is that they are isomorphisms precisely when the pairing is nondegenerate. In chain complexes, we can additionally consider nondegeneracy of a pairing in homology, and this is equivalent to the induced maps being quasi-isomorphisms. 
By the following proposition associating an IM-pairing to each simplicial pairing, this allows to define a notion of nondegeneracy up to homotopy for simplicial pairings which we summarize in Definition \ref{def:homological-non-deg}. 
This result appeared previously in \cite{CuecaZhu2023} for $n$-shifted 2-forms, where it was interpreted in \cite[Remark 2.15]{CuecaZhu2023} as an instance of the Van Est map discussed in \cite[\S 6]{AriasAbadCrainic2011} in relation with IM-forms. The formula \eqref{eq:IMPairingAssociatedDef} originally comes from \cite{Lesdiablerets}. Here, we reformulate it for pairings in $\SVect$ in terms of the Eilenberg-Zilber map. 

\begin{proposition}\label{prop:associated-IM-pairing}
    Any multiplicative normalized $n$-shifted pairing $\alpha: \huaV_n \otimes \huaW_n \to \R$ induces an associated $n$-shifted IM-pairing $\lambda_{\alpha}: (N(\huaV) \otimes N(\huaW))_n \to \R$ between the respective normalized complexes. $\lambda_\alpha$ is defined by the composition
    \begin{equation}
        \lambda_\alpha: (N(\huaV) \otimes N(\huaW))_n \overset{EZ}{\longrightarrow}
        N(\huaV \otimes \huaW)_n \overset{N(\alpha)}{\longrightarrow} N(B^n\R)_n = \R,
    \end{equation}
    where $EZ$ is the Eilenberg-Zilber map from Definition \ref{def:EilenbergZilberMap}. More explicitly, $\lambda_\alpha$ is given, for any $v \in N(\huaV)_i$ and $w \in N(\huaW)_{n-i}$, by
    \begin{equation}\label{eq:IMPairingAssociatedDef}
        \lambda_\alpha(v,w) = \sum_{(\mu,\nu)\in\Shuf(i,n-i)} \sign(\mu,\nu) \alpha(s_{\nu_{n-i}}\dots s_{\nu_1} v, s_{\mu_{i}}\dots s_{\mu_1} w),
    \end{equation}
    where $\Shuf(i,n-i)$ is the set of $(i,n-i)$-shuffles.
\end{proposition}

\begin{proof}
   By definition of $N$ as a functor and $EZ$ as a natural transformation, $\lambda_\alpha = EZ^* N(\alpha) \in \Hom_{\Ch_{\ge 0}(\Vect)}(N(\huaV) \otimes N(\huaW), \R[-n])$, hence it is clearly a chain map and therefore an IM-pairing. For a more detailed combinatorial proof starting from \eqref{eq:IMPairingAssociatedDef} only see \cite[Lemma E.1]{CuecaZhu2023}.\footnote{In that article, this fact is proven for the IM-pairing associated to an $n$-shifted 2-form on a simplicial manifold. However, the combinatorics are the same.}
\end{proof}

\begin{definition}\label{def:homological-non-deg}
    An $n$-shifted pairing $\alpha: \huaV_n \otimes \huaW_n \to \R$ is \textbf{homologically $n$-shifted nondegenerate} if its associated IM-pairing $\lambda_\alpha$ descends to a nondegenerate pairing 
    \begin{equation*}
        \lambda_\alpha: H(N(\huaV))_i \otimes H(N(\huaW))_{n-i} \to \R,
    \end{equation*}
    for any $i \in \Z$. That is, it induces an isomorphism between the homologies of the normalized complexes up to a degree shift of $n$.
    Equivalently $\alpha$ is homologically nondegenerate if either $\lambda_\alpha^r$ or $\lambda_\alpha^l$ is a quasi-isomorphism.\footnote{If one is a quasi-isomorphism, they both are. This follows from the fact that $n$-shifted duality of unbounded chain complexes is reflexive.}
\end{definition}

Due to the fact that homological nondegeneracy of an $n$-shifted pairing requires the normalized complex of $\huaV$ to be quasi-isomorphic to the \textit{total} $n$-shifted dual, and not the truncated one, the $n$-dual pairing is only homologically nondegenerate in certain cases, which we now discuss.

\begin{theorem}\label{thm:ndual-pairing-hom-nondeg}
    Let $\huaV$ be a simplicial vector space.

    The $n$-dual pairing $\langle \cdot, \cdot \rangle: \huaV^{n*} \otimes \huaV \to B^n\R$ is homologically $n$-shifted nondegenerate if and only if $\huaV$ is at most an $n$-type.
\end{theorem}

\begin{proof}
    By definition 
    \begin{equation*}
        \lambda_{\langle \cdot, \cdot \rangle}^r = \rho \circ EZ^* \circ N (\tau (\id_{\huaV^{n*}})): N(\huaV^{n*}) \to N(\huaV)^*[-n],
    \end{equation*}
    where $\rho$ here is the natural isomorphism giving the tensor-hom adjunction in $\Ch(\Vect)$ from Proposition \ref{prop:tensor-hom-chains}. 
    By evaluating the diagrams defining $EZ^H$ in Theorem \ref{thm:EZ-thm-hom} with $\huaW=B^n\R$, on $id_{N(\huaV^{n*})}$, we have that
    \begin{equation*}
        EZ^H = \rho_{\ge 0} \circ EZ^* \circ N (\tau (\id_{\huaV^{n*}})): N(\huaV^{n*}) \to N(\huaV)^*[-n]_{\ge 0}.
    \end{equation*}
    So $EZ^H$ is the truncation at $0$ of $\lambda_{\langle \cdot, \cdot \rangle}^r$, obtained by applying the functor $tr_{\ge 0}$ defined in \eqref{eq:non-neg-chain-maps}.

    If $\huaV$ is at most an $n$-type, then the inclusion $N(\huaV)^*[-n]_{\ge 0} \into N(\huaV)^*[-n]$ is a quasi-isomorphism and restricting $\lambda_{\langle \cdot, \cdot \rangle}^r$ along this gives precisely $tr_{\ge 0}\lambda_{\langle \cdot, \cdot \rangle}^r= EZ^H$, which is a quasi-isomorphism by Theorem \ref{thm:EZ-thm-hom}. Thus $\lambda_{\langle \cdot, \cdot \rangle}^r$ is a quasi-isomorphism  by the two-out-of-three property.

    On the contrary, if $\huaV$ is an $m$-type for $m > n$, the homology of $N(\huaV)^*[-n]$ must be non-zero at the negative degree $n-m$, while the homology of $N(\huaV^{n*})$ is concentrated in non-negative degrees. Therefore $\lambda^r_{\langle \cdot, \cdot \rangle}$ cannot be a quasi-\hspace{0pt}isomorphism and $\langle\cdot , \cdot \rangle$ cannot be homologically $n$-shifted nondegenerate.
\end{proof}

\begin{remark}\label{rem:loss-of-info}
We can see immediately from the proof that if $\huaV$ is a $\VS$ $n$-\hspace{0pt}groupoid, then the right induced map $\lambda_{\langle \cdot, \cdot \rangle}^r$ is exactly $EZ^H$. 
In a sense, this is the best case scenario, where the information of $N(\huaV)^*[-n]$ does not get lost at all, not just up to homotopy as in the case of an $n$-type. 
To make this statement more precise we follow the principle that $n$-dualization corresponds to shifted dualization on the normalized complexes but the Dold-Kan correspondence only sees non-negative degrees.

First of all, if $\huaV$ has order $m$, $N(\huaV)$ has amplitude $(0, m)$, and the $n$-shifted dual $N(\huaV)^*[-n]$ has amplitude $(n-m, n)$.
If $m\le n$, $N(\huaV)^*[-n]$ is non-negative, and it coincides with the truncation $N(\huaV)^*[-n] = N(\huaV)^*[-n]_{\ge 0}$. In this case $n$-dualization causes no truncation. 
On the other hand, if $m > n$, the amplitude of $N(\huaV)^*[-n]_{\ge 0}$ is strictly smaller than that of $N(\huaV)^*[-n]$, and in particular the non-negative dual contains less information than $N(\huaV)$. If $\huaV$ is at most an $n$-type there is a truncation, but this causes no loss of information up to homotopy as the homology of $N(\huaV)$ has amplitude smaller than $(0, n)$.
Otherwise, the truncation forgets some of the homology as well. Thus, for $m > n$, there is generally a loss of information, which is only avoided up to homotopy for $m$-groupoids with homotopy type at most $n$.

This phenomenon is consistent with the situation for $n$-shifted symplectic structures explained in Remark 2.16, Example 2.20, and Example 2.23 of \cite{CuecaZhu2023}. 
\end{remark}

\begin{remark}
    When extending this construction to the category of simplicial vector bundles, with the generalized Dold-Kan correspondence in \cite{HoyoTrentinaglia2024}, extra care must be taken. Whatever the $n$-dual may be in this category, the above discussion still applies, and the truncation might additionally cause issues with smoothness and regularity. These might cause the $n$-dual to not be representable as a simplicial vector bundle and ``fall out of the category''. The simplest example of this effect is the 0-dual of a higher vector bundle over a manifold which we discuss in \ref{rem:VBnGpd-0-dual}. 
\end{remark}

\begin{example}
    As in Example \ref{ex:VS0dual}, the 0-dual pairing for a $\VS$ 0-groupoid is even nondegenerate at the level of spaces, and its IM-pairing is nondegenerate on chains. As we will show in Remark \ref{rem:VS1Dual-DK}, the 1-dual pairing for a $\VS$ 1-groupoid is also nondegenerate at the level of spaces, and its IM-pairing is nondegenerate on chains. The situation changes entirely for the 2-dual pairing, as we will see in Proposition \ref{prop:vs-2dual-norm-cplx}, since this is only homologically nondegenerate even for $\VS$ 2-groupoids. 
    This can be interpreted as the 2-dual containing additional information that is redundant up to homotopy.
\end{example}

\begin{theorem}\label{thm:hom-nondeg-homotopy-equiv}
    Let $\huaV$ and $\huaW$ be simplicial vector spaces and $\alpha : \huaV \otimes \huaW \to B^n\R$ an $n$-shifted simplicial pairing between them. The following diagrams commute:
    \begin{equation}\label{diag:nShiftedPairingsInducedMapsTriangle}
    \begin{tikzcd}[ampersand replacement=\&,cramped]
	{N(\huaW)} \& {N(\huaV^{n*})} \& {N(\huaV)} \& {N(\huaW^{n*})} \\
	\& {N(\huaV)^*[-n]} \&\& {N(\huaW)^*[-n]}
	\arrow["{N(\alpha^l)}", from=1-1, to=1-2]
	\arrow["{\lambda_\alpha^l}"', from=1-1, to=2-2]
	\arrow["{\lambda_{\langle \cdot, \cdot \rangle_{tot}}^r}", from=1-2, to=2-2]
	\arrow["{N(\alpha^r)}", from=1-3, to=1-4]
	\arrow["{\lambda_\alpha^r}"', from=1-3, to=2-4]
	\arrow["{\lambda_{\langle \cdot, \cdot \rangle_{tot}}^r}", from=1-4, to=2-4]
    \end{tikzcd}
    \end{equation}
    Additionally, if $\huaW$ is at most an $n$-type, then $\alpha$ is homologically nondegenerate if and only if $\alpha^r$ is a weak equivalence. 
    Analogously, if $\huaV$ is at most an $n$-type, then $\alpha$ is homologically nondegenerate if and only if $\alpha^l$ is a weak equivalence. 
\end{theorem}

\begin{proof}
    The commutativity of the right-hand diagram follows from the fact that $\rho EZ^* N \tau$ can be seen as a natural transformation between the functors $\SVect^{op} \times \SVect^{op} \times \SVect \to \Set$ given by
    \begin{equation*}
        \SVect(\_, \IHom(\_, \_)) \to \Ch(N(\_), \IHom(N(\_),N(\_))).
    \end{equation*}
    Evaluating the third argument at $B^n\R$ gives a natural transformation 
    \begin{equation*}
        \SVect(\_, \_^{n*}) \to \Ch(N(\_), N(\_)^*[-n]),
    \end{equation*}
    where both are functors $\SVect^{op} \times \SVect^{op} \to \Set$. 
    The naturality square at the map $((\alpha^r)^*, id)$ in $\SVect^{op} \times \SVect^{op}$ is 
    \[\begin{tikzcd}[ampersand replacement=\&,cramped,row sep=scriptsize]
	{\SVect(\huaW^{n*}, \huaW^{n*})} \& {\Ch(\Vect)(N(\huaW^{n*}), N(\huaW)^*[-n])} \\
	{\SVect(\huaV, \huaW^{n*})} \& {\Ch(\Vect)(N(\huaV), N(\huaW)^*[-n])}
	\arrow["{\rho EZ^*N\tau}", from=1-1, to=1-2]
	\arrow["{((\alpha^r)^*, id)}"', from=1-1, to=2-1]
	\arrow["{((N(\alpha^r))^*, id)}", from=1-2, to=2-2]
	\arrow["{\rho EZ^*N\tau}"', from=2-1, to=2-2]
    \end{tikzcd}\]
    which means that for $id_{\huaW^{n*}} \in \SVect(\huaW^{n*}, \huaW^{n*})$, 
    \begin{equation*}
    \begin{split}
        \lambda_{\langle \cdot, \cdot \rangle}^r \circ N(\alpha^r) &= (\rho EZ^* N \tau)(id_{\huaW^{n*}}) \circ N(\alpha^r) \\
        &= \rho EZ^* N \tau ((\alpha^r)^*, id)(id_{\huaW^{n*}}) = \rho EZ^* N \tau ( \alpha^r) = \lambda_{\alpha}^r.
    \end{split}
    \end{equation*}

    Commutativity of the left-hand diagram follows analogously.

    The other two statements follow from Theorem \ref{thm:ndual-pairing-hom-nondeg}, the two-out-of-three property for quasi-isomorphisms, and Proposition \ref{prop:equivalences-in-svect}.
\end{proof}

\begin{remark}\label{rem:hndg-pairing-admitted}
    In the setting of the theorem, if either $\huaV$ or $\huaW$ is at most an $n$-type and $\alpha$ is homologically nondegenerate, then the other is also at most an $n$-type because it is weak equivalent to a $\VS$ $n$-groupoid. This puts restrictions on which simplicial vector spaces admit homologically nondegenerate pairings with respect to a certain shift. 
    This is analogous to the fact discussed in \cite[Remark 2.16]{CuecaZhu2023} that an $m$-shifted symplectic Lie $n$-groupoid with $m < n$ must have certain vanishing homology groups, which make it interpretable as an $m$-shifted symplectic Lie $m$-groupoid with added singularities. 
\end{remark}

\begin{remark}
    A trivial but perhaps illustrative fact to observe is that in the case of the $n$-dual pairing $\alpha= \langle \cdot, \cdot \rangle: \huaV^{n*} \otimes \huaV \to B^n\R$ the right-hand diagram in \eqref{diag:nShiftedPairingsInducedMapsTriangle} becomes tautological, and the theorem only reiterates the fact that if $\huaV$ is an $n$-type to begin with, then $\lambda^r_{\langle \cdot, \cdot \rangle}$ is a quasi-isomorphism and the $n$-dual pairing is homologically nondegenerate. 
    
    On the contrary, the left-hand diagram yields an interesting property which deserves to be a theorem of its own.
\end{remark}

\begin{theorem}[$n$-duality is reflexive up to homotopy]\label{thm:ndual-reflexive-uth}
    Let $\huaV$ be at most an $n$-type. The double $n$-dual $(\huaV^{n*})^{n*}$ is weak equivalent to $\huaV$ itself. That is
    \begin{equation} \label{eq:double-W.E.}
        \langle \cdot, \cdot \rangle^l: \huaV \to (\huaV^{n*})^{n*}
    \end{equation}
    is a weak equivalence. 
\end{theorem}

\begin{proof}
    By Theorem \ref{thm:hom-nondeg-homotopy-equiv} and Theorem \ref{thm:ndual-pairing-hom-nondeg}, the left induced map in Example \ref{ex:ndual-pairing-induced-maps} is a weak equivalence.
\end{proof}

\begin{remark}
    In $\SVect$, by applying Proposition \ref{prop:equivalences-in-svect}, if $\huaV$ is an $n$-type, then $\langle \cdot, \cdot \rangle^l: \huaV \to (\huaV^{n*})^{n*}$ is even a homotopy equivalence. 
\end{remark}

\section{Computing \texorpdfstring{$n$}{n}-duals of simplicial vector spaces}\label{sec:computations}

To calculate explicitly what the $n$-dual \eqref{eq:vs-n-dual} is involves solving many linear equations. 
Determining the solution space of these linear equations for a general $n$ is not a trivial task due to the large number of equations involved. 
In this section we give an overview of the equations for a general $n$ and set out to solve them in the case of a $\VS$ $n$-groupoid for $n=1$ and $n=2$. For $n=1$, we rediscover the dual of $\VB$ groupoids from Definition \ref{def:VB1Dual} applied to the case where the base is a point. For $n=2$, we discover a new object. 

\subsection{Overview of the general computation for the \texorpdfstring{$n$}{n}-dual of a simplicial vector space}\label{sec:ndual-overview}

Let $\huaV$ be a simplicial vector space. We now write down the linear equations to compute $\huaV^{n*}$ for a general $n$. Since $\huaV^{n*}$ is a $\VS$ $n$-groupoid, its full data is contained in its $n$-truncation, by Proposition \ref{prop:FiniteDataVS}. Even though the multiplications of its $n$-simplices can be recovered by computing Moore fillers, these have a nice description in terms of the $n$-dual pairing, which we later discuss.

\begin{proposition}\label{prop:vs-n-dual-equations}
    Let $\huaV$ be a simplicial vector space. Then, any $m$-simplex $\phi \in \huaV^{n*}_m$ consists of $\binom{m+n+1}{n+1}$ components $\phi^u\in\huaV_n^*$ indexed over $n$-simplices $u\in \Delta^m_n$ defined by 
    \begin{equation}\label{eq:n-dual-element-in-comp}
        \phi^u(v):= \phi_n(v^u), \text{ for all }v^u \in \huaV_n^u \subseteq \huaV \otimes \Delta^m_n, \text{ and } u \in \Delta^m_n.
    \end{equation}
    These components satisfy the $n \cdot \binom{m+n}{n}$ \textbf{normalization equations}
    \begin{equation}\label{eq:NormalizationMostGeneral}
        \phi^{s_i t}(s_i v) = 0 , \text{ for all } 0\le i \le n-1,  \; v \in \huaV_{n-1}, \text{ and } t \in \Delta^m_{n-1},
    \end{equation}
    and the $\binom{n+m+2}{n+2}$ \textbf{multiplicativity equations}
    \begin{equation}\label{eq:MultiplicativityMostGeneral-short}
    \sum_{i=0}^n (-1)^i \phi^{d_ir}(d_i v) = 0, \text{ for all } r\in \Delta^m_{n+1}, \text{ and } v \in \huaV_{n+1}.
\end{equation}
\end{proposition}

\begin{proof}
By definition $\huaV^{n*}_m = \IHom(\huaV \otimes \Delta_m, B^n\R)$. 
At any fixed level $m$, since $B^n\R$ is a $\VS$ $n$-group, by  Lemma \ref{lem:finite-data-simp-map}, a simplicial map $\phi \in \huaV^{n*}_m = \Hom_{\Simp\Vect}(\huaV \otimes\Delta^m, B^n\R)$ is determined by its first $n$ levels, on condition they commute with the truncated simplicial structure, and they satisfy a multiplicativity condition.  More precisely, for any choice of $0\le k \le n+1$ to represent $B^n\R_{n+1}$ as the horn $\Lambda^{n+1}_k(B^n\R) \cong \R^{n+1}$,
we have
\begin{equation}\label{diag:VSnDualMultiplicativity}
    \begin{tikzcd}[ampersand replacement=\&]
	{\bigoplus\limits_{r \in \Delta^m_{n+1}} \huaV_{n+1}^r} \& {\bigoplus\limits_{u \in \Delta^m_{n}} \huaV_n^u} \& {(\huaV\otimes \Delta^m)_{\le n-1}} \\
	{\R^{n+1}} \& \R \& 0
	\arrow["{\phi_{n+1}=(\phi_n d_0,\dots, \widehat{\phi_n d_k}, \dots ,\phi_n d_n)}"', from=1-1, to=2-1]
	\arrow[shift left=2, from=1-1, to=1-2]
	\arrow["\dots"{description}, from=1-1, to=1-2]
	\arrow[shift right=2, from=1-1, to=1-2]
	\arrow["\dots"{description}, from=2-1, to=2-2]
	\arrow[shift left=2, from=2-1, to=2-2]
	\arrow[shift right=2, from=2-1, to=2-2]
	\arrow["{\phi_n}", from=1-2, to=2-2]
	\arrow["\dots"{description}, from=2-2, to=2-3]
	\arrow["\dots"{description}, from=1-2, to=1-3]
	\arrow["{\phi_{\le n-1} =0}", from=1-3, to=2-3]
	\arrow[shift left=2, from=2-2, to=2-3]
	\arrow[shift right=2, from=2-2, to=2-3]
	\arrow[shift left=2, from=1-2, to=1-3]
	\arrow[shift right=2, from=1-2, to=1-3]
    \end{tikzcd}
\end{equation}
where: the components $\phi_{\le n-1}$ vanish because $B^n\R_{\le n-1}= 0$, and, by commutativity of $\phi$ with the face maps for $i\neq k$ between level $n+1$ and $n$, 
\begin{equation*}
    \phi_{n+1} = (\phi_n d_0,\dots, \widehat{\phi_n d_k}, \dots ,\phi_n d_n).
\end{equation*}
Thus, $\phi$ is uniquely determined by $\phi_n$, subject to commutativity with the degeneracy maps between level $n$ and $n-1$ and the multiplicativity condition.
In \eqref{eq:n-dual-element-in-comp} we are just writing $\phi$ in components over the copowering $\huaV \otimes \Delta^m$. Then we have that commutativity with the degeneracy maps is equivalent to \eqref{eq:NormalizationMostGeneral}, while the multiplicativity condition gives
\begin{equation}\label{eq:MultiplicativityMostGeneral}
\begin{split}
    \phi^{d_k r}(d_k v) &= m_k^{B^n\R}((\phi^{d_ir}(d_i v))_{i\neq k}) \\
    &=\sum_{i=0, \, i\neq k }^{n} (-1)^{i-k+1} \phi^{d_i r}(d_i v), \quad \forall r\in \Delta^m_{n+1}, \; \forall v \in \huaV_{n+1},
\end{split}
\end{equation}
which can be written more synthetically (and independently of a choice of $k$) as \eqref{eq:MultiplicativityMostGeneral-short}.
The numbers of components and equations come from \eqref{eq:cardinality-m-simplices-in-Delta-n}. 
\end{proof}

\begin{remark}
In other words, equations \eqref{eq:MultiplicativityMostGeneral} and \eqref{eq:NormalizationMostGeneral} represent, respectively, the commutative diagrams for the only non-trivial face map and that for the degeneracy maps, i.e.
\[\begin{tikzcd}[ampersand replacement=\&]
	{\huaV_{n+1}} \& {\huaV_n^{d_k r}} \\
	{\R^{n+1}} \& \R
	\arrow["{(\phi^{d_0 r}, \dots, \widehat{\phi^{d_kr}},\dots, \phi^{d_n r})}"', from=1-1, to=2-1]
	\arrow["{d_k}", from=1-1, to=1-2]
	\arrow["{\phi^{d_k r}}", from=1-2, to=2-2]
	\arrow["{m_k^\R}"', from=2-1, to=2-2]
\end{tikzcd}, \quad
\begin{tikzcd}[ampersand replacement=\&]
	{\huaV_n\otimes \Delta^m_n} \& {(\huaV_{n-1}\otimes \Delta^m_{n-1})} \\
	\R \& 0
	\arrow["{\phi_n}", from=1-1, to=2-1]
	\arrow["{s_i = 0}"', from=2-2, to=2-1]
	\arrow["{\phi_{n-1} =0}", from=1-2, to=2-2]
	\arrow["{s_i}"', from=1-2, to=1-1]
\end{tikzcd}\]where $m_k^\R$ is given in \eqref{eq:defMultiplicationBnR}. 
\end{remark}

By the conditions \eqref{eq:NormalizationMostGeneral} and \eqref{eq:MultiplicativityMostGeneral-short}, each $\phi$ is determined by a subset of the components $\phi^u$ for $u\in \Delta^m_n$, which will depend on a choice. Furthermore, these components will not be completely independent of each other, but also partially related.
As such, with the computation of each level of $\huaV^{n*}$, we will have to make a choice of which of the $\phi^u$ to consider as ``independent variables'', show how they are related, and show how to determine the others from this subset without overdetermination. Since even the ``independent variables'' will have some relation among them, we will obtain each level of $\huaV^{n*}$ as a fiber product.
Once the spaces are determined, the simplicial maps of $\huaV^{n*}$ are simply projections to certain components of each element, as in \eqref{eq:mapping-set-maps}.

We now classify the components of a simplex in the $n$-dual to better interpret their meaning.

\begin{definition}\label{def:face-interior-comps}
We call \textbf{face components} of $\phi$ the components that correspond to indices $u \in \Delta^m_n$ such that $u=\delta_i q$ for some $q \in \Delta^{m-1}_{n}$ and $0\le i \le n$. These make up the faces of $\phi$ as an $m$-simplex $\phi \in \huaV^{n*}_m$. 

We call \textbf{interior components} of $\phi$ all the other components, i.e. those for which $u\neq \delta_i q$ for $q \in \Delta^{m-1}_{n}$ and any $0\le i \le n$. 
\end{definition}

\begin{remark}\label{rem:n-dual-eqs-repeat}
For each fixed $i$, the components $\phi^u$ with $u=\delta_i q$ form the $i$-th face of $\phi$, which is an $(l-1)$-simplex in $\huaV^{n*}$. Thus they are related by the equations \eqref{eq:MultiplicativityMostGeneral-short} obtained at the previous level.
This means that the equations that determine the components of an $(m-1)$-simplex must appear again in the computation of the space of $m$-simplices.
In other words, when computing each level $m$ starting from 0, the equations split into two groups: ones that are repeated from a previous level, which describe internal relations between the components of each face, and which have already been solved to compute level $m-1$; and new ones which describe relations between components of different faces and the interior of the simplex.
The ones that describe the internal relations of the components of each face of $f$ are the ones in \eqref{eq:MultiplicativityMostGeneral-short} for which $r$ is a degenerate simplex of an $(n-1)$-simplex in $\Delta^m$, while the new ones are the ones for which $r$ is nondegenerate or a degenerate simplex of a nondegenerate $n$-simplex in $\Delta^m$.
This will appear repeatedly in the computations that follow.
We discuss this in more detail in the proof of Lemma \ref{lem:n-dual-finite-equations}.
\end{remark}

\begin{remark}\label{rem:number-of-int-comp}
In this remark we use the decomposition of each level of the $m$-simplex in Proposition \ref{prop:shuffles-standard-simpl-decomp} to make some considerations about interior components at each level.

The interior components of $\phi \in \huaV^{n*}_m$, for $0 \le m <n$ can equivalently be characterized by the fact that they correspond to the degenerate simplices in $\Delta^m_{n}$ which are of the form $u=s_IE_m$ for $E_m$ the unique nondegenerate $m$-simplex in $\Delta^m_m$ and $I$ some multi-index of length $n-m$. 
Because of this, at any level $0 \le m < n$ all the interior components are always normalized in some way, by \eqref{eq:NormalizationMostGeneral}. Even more so, all the face components at these levels must also be normalized. 
Therefore these levels are described by fiber products and intersections of some annihilator spaces. 
We will classify the annihilator spaces in more detail in Section \ref{sec:annihilators-and-dual-kernels}.

On the other hand, at level $n$, $|I|=0$ and there is only one component representing the interior of an $n$-simplex in the $n$-dual, which is the one related to the generating $n$-simplex $E_n$: $\phi^{E_n}$.
This is also the unique component of an $n$-simplex that is not normalized, as it does not appear in any of the normalization equations \eqref{eq:NormalizationMostGeneral}. 
Hence $\huaV^{n*}_n$ can be written as $\huaV^{n*}_n \cong \huaV_n^* + O$, where $O$ is the space spanned by all the other components (which are necessarily normalized in some way). 
Furthermore, by recalling the definition of the $n$-dual pairing and of $\tau$ in \eqref{eq:tensor-hom-svect-def-tau}, the component $\phi^{E_n}$ is the one appearing in the $n$-dual pairing, as we have
\begin{equation}\label{eq:ndual-pairing-def}
    \langle \phi, X \rangle = \tau(id)(\phi \otimes X) = \phi^{E_n}(X) \qquad \forall \phi \in \huaV^{n*}_n, \forall X \in \huaV_n,
\end{equation}
which can be seen as an alternative definition of the pairing. 
More importantly this expresses the $n$-dual pairing as an \textit{extension} of the usual evaluation pairing of $\huaV_n$ and its dual. 
This extension is trivial in that it is set to zero on the space $O \subset \huaV^{n*}_n$ spanned by all the other components. 

Finally, $m$-simplices for $m > n$ do not admit interior components, and this is consistent with the fact that every level of the $n$-dual above $n$ contains only degenerate simplices by Lemma \ref{lem:MooreFillers-degeneracies} and Proposition \ref{prop:FiniteDataVS}.  
\end{remark}

We now present a way to describe the multiplication maps of the $n$-dual, alternative to the computation of Moore fillers. 

\begin{proposition}\label{prop:n-dual-mult-pairing}
    Let $\huaV$ be a simplicial vector space and $\huaV^{n*}$ its $n$-dual. The multiplication maps $\{\widecheck{m}_k\}$ of $\huaV^{n*}$ are uniquely determined by requiring that for any $(n+1, k)$-horn $(\phi_i)_{0\le i \neq k \le n} \in \Lambda^{n+1}_k(\huaV^{n*})$,
    \begin{equation}\label{eq:n-dual-mult-simp-id}
        \widecheck{d}_j\widecheck{m}_k((\phi_i)_{0\le i \neq k \le n}) =
        \begin{cases}
            \widecheck{d}_{k-1} \phi_j &\text{if } j<k,\\
            \widecheck{d}_k\phi_{j+1} &\text{if } j\ge k,
        \end{cases}
    \end{equation}
    and for any $(n+1)$-simplex $v \in \huaV_{n+1}$,
    \begin{equation}\label{eq:ndual-mult-pairing}
    \langle \widecheck{m}_k((\phi_i)_{0\le i \neq k \le n}), d_k v \rangle = \sum_{i=0, \, i\neq k }^{n} (-1)^{i-k+1} \langle \phi_i , d_i v \rangle.
    \end{equation}
    In other words, each multiplication of the $n$-dual can be defined by being the unique multiplication $\Lambda^{n+1}_k(\huaV^{n*}) \to \huaV^{n*}_n$ that makes the $n$-dual pairing multiplicative. 
\end{proposition}

\begin{proof}
When computing the product of any $n$-simplices forming a horn, it is clear that some of the face components of these $n$-simplices must already be the face components of the product, by the simplicial identities.\footnote{For $n=1$ this is simply saying that the endpoints of the factors in a composition of arrows determine the endpoints of the result of the composition. } Thus \eqref{eq:n-dual-mult-simp-id} determines all the face components of $\widecheck{m}_k((\phi_i)_{i\neq k})$. 
Therefore, and because each $n$-simplex has a unique interior component by Remark \ref{rem:number-of-int-comp}, to fully determine the $n$-dual multiplication it is enough to express the interior component of the product in terms of the interior components of the factors. 
There is a unique multiplicativity equation expressing exactly this relation. 
This is \eqref{eq:MultiplicativityMostGeneral} for $m=n+1$ and $r=E_{n+1}=01\dots n(n+1)$, the unique nondegenerate $(n+1)$-simplex in $\Delta^{n+1}$. 
In fact, $d_i E_{n+1} = \delta_i E_n$, so this is a relation between the $\phi^{\delta_i E_n} = (\widecheck{d}_i \phi)^{E_n}$ which are precisely the interior components of each face of an $(n+1)$-simplex $\phi \in \huaV^{n*}_{n+1}$. 
Solving for each $(\widecheck{d}_k \phi)^{E_n}$ determines the map $\widecheck{m}_k$ in terms of the interiors of the other faces of $\phi$, which form an $(n+1, k)$-horn. 
Therefore, in terms of the $n$-dual pairing, the equation \eqref{eq:MultiplicativityMostGeneral} for $r=E_{n+1}$ is exactly the multiplicativity condition of $\langle \cdot, \cdot \rangle$ as a pairing, defined in \eqref{eq:nShiftedPairingMultiplicative}. That is, \eqref{eq:ndual-mult-pairing}.
\end{proof}

\begin{remark}
The fact that $\langle \cdot, \cdot \rangle$ is normalized as a pairing as in \eqref{eq:nShiftedPairingNormalized} is equivalent to the normalization conditions \eqref{eq:NormalizationMostGeneral} for $t = E_{n-1}$. 
In fact, by the definition of the degeneracy maps in equation \eqref{eq:mapping-set-maps} and the fact that $s_i E_{n-1} = \sigma_i E_n$, for any $0 \le i \le n-1$, the $E_n$ component of $\widecheck{s}_i f$ is given by $\phi^{s_iE_{n-1}}$, and 
\begin{equation*}
    \langle \widecheck{s}_i \phi, s_i v \rangle = \phi^{\sigma_iE_{n}}(s_i v) = \phi^{s_iE_{n-1}}(s_i v) = 0.
\end{equation*}
\end{remark}

\begin{remark}\label{rem:ndual-mult-pairing-n-gpd}
If $\huaV$ is a $\VS$ $n$-groupoid, then equation \eqref{eq:MultiplicativityMostGeneral} reads
\begin{equation}\label{eq:MultiplicativityMostGeneral-n-gpd}
    \phi^{d_k r}(m_k((v_i)_{0\le i \neq k \le n})) = \sum_{i=0, \, i\neq k }^{n} (-1)^{i-k+1} \phi^{d_i r}(v_i), 
\end{equation}
for any $r\in \Delta^m_{n+1}$, and any $(v_i)_{0\le i \neq k \le n} \in \Lambda^{n+1}_k(\huaV)$, hence the name of multiplicativity condition.
Clearly, by rewriting it as \eqref{eq:MultiplicativityMostGeneral-short}, if \eqref{eq:MultiplicativityMostGeneral-n-gpd} holds for any $0\le k \le n$, then it also holds for all other $k$.

Moreover, if $\huaV$ is a $\VS$ $n$-groupoid, each multiplication map of $\huaV^{n*}$ is defined in terms of the corresponding multiplication of $\huaV$, by the property that for any $(n+1, k)$-horns $(\phi_i)_{0\le i \neq k \le n} \in \Lambda^{n+1}_k(\huaV^{n*})$ and $(v_i)_{0\le i \neq k \le n} \in \Lambda^{n+1}_k(\huaV)$,
\begin{equation}\label{eq:ndual-mult-pairing-n-gpd}
    \langle \widecheck{m}_k((\phi_i)_{0\le i \neq k \le n}), m_k((v_i)_{0\le i \neq k \le n})\rangle = \sum_{i=0, \, i\neq k }^{n} (-1)^{i-k+1} \langle \phi_i , v_i \rangle.
\end{equation}
This is how the multiplication of the dual $\VB$ groupoid was defined in \ref{def:VB1Dual}. 
\end{remark}

\begin{remark}
In the language of $n$-shifted pairings, each level $\huaV^{n*}_m$ of the $n$-dual of a VS $n$-groupoid $\huaV$ can be seen as the space of multiplicative normalized $n$-shifted pairings between $\huaV$ and $\R[\Delta^m]$, since these are equivalent to simplicial linear maps $\huaV \otimes \R[\Delta^m] = \huaV \otimes \Delta^m \to B^n\R$. 
By the adjunction in Prop. \ref{prop:tensor-hom-svect}, each of these pairings admits a left induced map, which is exactly identifying an $m$-simplex in the $n$-dual. More precisely, for any $\phi \in \huaV^{n*}_m$, this can be seen as a map $\Delta^m \to \huaV^{n*}$, and extended to a linear map $\phi^l: \R[\Delta^m] \to \huaV^{n*}$, which is precisely the left induced map of $\phi$ when seen as a pairing. 
\end{remark}

\subsection{Annihilators and generalized cores}\label{sec:annihilators-and-dual-kernels}

This section is devoted to the study of the solution spaces of the normalization equations \eqref{eq:NormalizationMostGeneral}. 
As we will see, these generalize the different isomorphic descriptions of the dual of the core of a $\VB$ 1-groupoid we discussed in Section \ref{sec:VB1-groupoids-cores-duals}.

\begin{definition}\label{def:VS-cores-deg-ann}
Let $\huaV$ be a simplicial vector space. 

For any subset $A \subseteq [m] \in \Delta$ of cardinality $|A|=j$, the associated \textbf{degree $j$ $m$-dimensional core} is $\bigcap_{i \in A} \ker d^m_i \subseteq \huaV_m$.\footnote{When defining generalized horns as in \ref{def:gen-horns}, this is the kernel of a horn projection: $\ker p^m_{[m]/A}$.}

For any subset $B \subseteq [m-1] \in \Delta$ of cardinality $|B|=j$, the associated \textbf{degree $j$ $m$-dimensional degeneracy annihilator} is $O_B := \Ann\left(\sum_{i\in B} s_i \huaV_{m-1}\right) \subseteq \huaV_{m}^*$. 
\end{definition}

We now focus on the two extreme cases needed for computing 1- and 2-duals, which we use in Section \ref{sec:computation-1dual} and Section \ref{sec:computation-2dual}. These are degree 1 $n$-dimensional cores and degree $n$ $n$-dimensional cores. As we will show, in these cases, the degeneracy annihilators are isomorphic to the duals of the cores of the corresponding degree. We comment on a possible way to investigate this for degree $1<k<n$ cores in Remark \ref{rem:general-degree-cores}.

Since different annihilators appear separately as solution spaces of the normalization equations \eqref{eq:NormalizationMostGeneral}, we refrain from identifying them as one object, as keeping track of the isomorphisms required would overcomplicate the explicit formulas in the computation. 
This will pay off further in the computation of the $\VB$ $n$-duals which we carry out in the next chapter. 

The degree $1$ $n$-dimensional degeneracy annihilators are the spaces that contain the components of an element in the $n$-dual that are normalized with respect to a single degeneracy map, $O^n_i:=\Ann(s^n_i\huaV_{n-1}^*) \subseteq \huaV_n$.\footnote{As usual with face and degeneracy maps, we will omit the upper index $n$ indicating the dimension of the ambient space, whenever this is clear from context.} We observe that for $0\le i < n$, the following exact sequences canonically split:
\begin{equation}\label{eq:VSnDualSESd_is_i}
\begin{tikzcd}[ampersand replacement=\&]
	0 \& {\ker d_i} \& {\huaV_n} \& {\huaV_{n-1}} \& 0
	\arrow[shift left=1, hook, from=1-2, to=1-3]
	\arrow[from=1-1, to=1-2]
	\arrow["{d_i}", shift left=1, from=1-3, to=1-4]
	\arrow["{s_i}", shift left=1, from=1-4, to=1-3]
	\arrow[from=1-4, to=1-5]
	\arrow["{id- s_id_i}", shift left=1, from=1-3, to=1-2],
\end{tikzcd}
\end{equation}
\begin{equation}\label{eq:VSnDualSESd_iplus1s_i}
\begin{tikzcd}[ampersand replacement=\&] 
	0 \& {\ker d_{i+1}} \& {\huaV_n} \& {\huaV_{n-1}} \& 0
	\arrow[shift left=1, hook, from=1-2, to=1-3]
	\arrow[from=1-1, to=1-2]
	\arrow["{d_{i+1}}", shift left=1, from=1-3, to=1-4]
	\arrow["{s_{i}}", shift left=1, from=1-4, to=1-3]
	\arrow[from=1-4, to=1-5]
	\arrow["{id- s_{i}d_{i+1}}", shift left=1, from=1-3, to=1-2].
\end{tikzcd}
\end{equation}
This implies that $\huaV_n \cong \ker d_i \oplus s_i(\huaV_{n-1}) \cong \ker d_{i+1} \oplus s_i(\huaV_{n-1})$.  As a result, we have the isomorphisms
\begin{equation}\label{eq:VSDual-isos-deg-ann-dual-ker}
   O_i = \Ann(s_i\huaV_{n-1}) \cong (\ker d_i)^* \cong (\ker d_{i+1})^*\cong O_{i+1}, \quad 0\le i \le n-1. 
\end{equation}
Different choices of isomorphic descriptions of the $O_i$ lead to different explicit descriptions of the $n$-dual. This generalizes the fact that the dual $\VB$ groupoid of a $\VB$ groupoid can be written in terms of the right core, the left core or the annihilator of the units of the $\VB$ groupoid.

In the other extreme case, the degree $n$ $n$-dimensional degeneracy annihilator is the space that contains fully normalized components
\begin{equation*}
    O_{01\dots (n-1)}^n := \Ann(D_n\huaV) = \bigcap_{i=0}^n O_i^n \subseteq \huaV_n^*
\end{equation*}
where $D_n\huaV = s_0\huaV_{n-1} + \dots + s_{n-1}\huaV_{n-1}$ is the space generated by all degenerate $n$-simplices in $\huaV$. As observed in Lemma \ref{lem:MooreFillers-degeneracies}, $D_n\huaV$ is isomorphic to the horn space $\Lambda^n_k(\huaV)$ via $\mu^n_k$, for any $0 \le k \le n$. The degree $n$ $n$-dimensional core $\ker p^n_k$ appears in the exact sequence
\begin{equation}\label{eq:VSnDualSESp_kmu_k}
\begin{tikzcd}[ampersand replacement=\&,cramped]
	0 \& {\ker p^n_k} \& {\huaV_n} \& {\Lambda^n_k(\huaV)} \& 0,
	\arrow[from=1-1, to=1-2]
	\arrow[shift left, hook, from=1-2, to=1-3]
	\arrow["{\ggamma^n_k = id- \mu^n_k}", shift left, from=1-3, to=1-2]
	\arrow["{p^n_k}", shift left, from=1-3, to=1-4]
	\arrow["{\mu^n_k}", shift left, from=1-4, to=1-3]
	\arrow[from=1-4, to=1-5]
\end{tikzcd}
\end{equation}
where we call the retract $\ggamma^n_k:= id-\mu^n_k$ the \textbf{$k$-th (degree $n$) core projection}. These projections also provide isomorphisms between cores for different indices $k$. From the dual sequence we have
\begin{equation}\label{eq:annihilator-core-iso}
    O_{01\dots (n-1)}^n = \Ann(D_n\huaV) \cong (\ker p^n_k)^*, \qquad \forall 0 \le k \le n,
\end{equation}
where the isomorphisms are given by the $\ggamma^*_k$. 

\begin{remark}\label{rem:general-degree-cores}
    By using the algorithm in the proof of Prop. \ref{thm:MooreHornFillers} we expect to be able to construct fillers for arbitrary generalized horns and obtain sequences such as \eqref{eq:VSnDualSESp_kmu_k}. This would allow us to classify all degree $j$ cores.
    We plan to study this in more detail in future work.
\end{remark}

\subsection{The 1-dual of a \texorpdfstring{$\VS$}{VS} 1-groupoid}\label{sec:computation-1dual}

Let $\huaV$ be a $\VS$ 1-groupoid $\huaV_1 \rightrightarrows \huaV_0$.
As remarked in Section \ref{sec:annihilators-and-dual-kernels}, the normalized components of each $m$-simplex in the $n$-dual are elements of an annihilator space, which is isomorphic to the dual of each core. For the $1$-dual, only the degree 1 1-dimensional degeneracy annihilator $O_0= \Ann(1\huaV_0)$ is relevant. Recall from Section \ref{sec:VB1-groupoids-cores-duals} that
\begin{equation*}
    O_0 = \Ann(1\huaV_0) \cong (\ker d_0)^* \cong (\ker d_1)^* \subseteq \huaV_1^*.
\end{equation*}
The core projections $\ggamma_0: \huaV_1 \to \ker d_0$ and $\ggamma_1:\huaV_1 \to \ker d_1$ appearing in \eqref{eq:VSnDualSESp_kmu_k} for $p^1_1=d_0$ and $p^1_0=d_1$ are 
\begin{equation*}
    \ggamma_0v = v - 1d_0v, \qquad \ggamma_1v = v - 1d_1v, 
\end{equation*}
and they coincide with the involution in Lemma \ref{lem:VB1-core-involution-iso} when restricted to the opposite core. 
Their dual maps $\ggamma_i^*: (\ker d_i)^* \to O_0$ give the isomorphisms between $O_0$ and $(\ker d_i)^*$. 

We summarize the computation of $\huaV^{1*}$ in the following proposition.

\begin{proposition}[1-dual of a $\VS$ 1-groupoid] \label{prop:VS1dual}
    Let $\huaV$ be a $\VS$ 1-groupoid. Then its 1-dual $\huaV^{1*}$ is
    \begin{equation*}
            \huaV_1^* 
            \rightrightarrows 
            O_0,  
    \end{equation*}
    with face maps given for any $\xi \in \huaV_1^*$ and $v \in \huaV_1$ by
    \begin{equation}
    \begin{aligned}
    \widecheck{d}_0\xi(v) = \ggamma_1^*\xi(v) = \xi(v - 1d_1v) ,\quad 
    \widecheck{d}_1\xi(v) = \ggamma_0^*\xi(v) = \xi(v - 1d_0v) ,
    \end{aligned}
    \end{equation}
    and unit map  $\widecheck{1}: O_0 \to \huaV_1^*$ the inclusion. 
    The multiplication $\widecheck{m}_1$ is
    \begin{equation}
        \langle \widecheck{m}_1(\eta,\xi), m_1(v,w) \rangle 
        = \langle \xi \cdot \eta, w \cdot v \rangle 
        = \langle \eta, v\rangle + \langle \xi, w \rangle,
    \end{equation}
    for any composable pairs $(\eta, \xi) \in \Lambda^2_1(\huaV^{1*})$, $(v,w)\in \Lambda^2_1(\huaV)$. 
\end{proposition}
\begin{proof}
According to Section \ref{sec:ndual-overview}, to compute each level of the 1-dual $\huaV^{1*}$ we only have to solve the linear equations given by the multiplicativity  \eqref{eq:MultiplicativityMostGeneral} and normalization \eqref{eq:NormalizationMostGeneral} conditions.  

Beginning with level 0, since the only 1-simplex in $\Delta^0$ is $00$, we have that an element $f \in \huaV^{1*}_0 = \Hom(\huaV, B^1\R)$ only has one component $f^{00} = f$. This has to satisfy only one normalization condition \eqref{eq:NormalizationMostGeneral} for $i=0$ and $p = 0 \in \Delta^0_0$ and one multiplicativity condition \eqref{eq:NormalizationMostGeneral} for $r=000 \in \Delta^0_2$.
The normalization condition reads $f(s_0x) = 0$ for all $x \in \huaV_0$, therefore $f \in O_0$. 
The multiplicativity condition reads $f(w\cdot v) = f(w) + f(v)$, for any $(v,w) \in \Lambda^2_1(\huaV)$. Since $w\cdot v = d_1\mu^2_1(v,w) = v + w - 1d_0w$ as in Example \ref{ex:MooreFillers1Gpd}, this is automatically satisfied by linearity and normalization of $f$. Therefore $\huaV^{1*}_0 = O_0$. 

For level 1, any $\xi\in \huaV^{1*}_1$ has three components: $\xi^{00}$, $\xi^{01}$ and $\xi^{11}$. 
There are two normalization conditions given by \eqref{eq:NormalizationMostGeneral} for each of the two 0-simplices $p = 0,1$ in $\Delta^1$, and four multiplicativity conditions given by \eqref{eq:MultiplicativityMostGeneral} for each of the four 2-simplices $r=000, 001, 011, 111$ in $\Delta^1$.
The normalization conditions read $\xi^{00}(s_0x)=0$ and $\xi^{11}(s_0x) =0$ for any $x \in \huaV_0$. Hence $\xi^{00}, \xi^{11} \in O_0$. Note that by the definition of the face maps in \eqref{eq:mapping-set-maps}, $\widecheck{d}_0\xi=\xi^{11}$ and $\widecheck{d}_1\xi= \xi^{00}$, and the normalization conditions state precisely the fact that these are in $\huaV^{1*}_0$, as should be expected. 
For $i=0,1$, the multiplicativity equations for $r=iii$ read $\xi^{ii}(w\cdot v)=\xi^{ii}(v) + \xi^{ii}(w)$ for $(v,w)\in\Lambda^2_1(\huaV)$. These are exactly the multiplicativity conditions at level 0 for each of the components $\xi^{ii}$, hence they are automatically satisfied by the same argument as before. 
The only remaining equations are the multiplicativity conditions for $001$ and $011$. We claim that they define the face maps $\widecheck{d}_0\xi=\xi^{11}$ and $\widecheck{d}_1\xi= \xi^{00}$. In fact, these equations read, for any $(v,w)\in \Lambda^2_1(\huaV)$,
\begin{equation*}
    \xi^{00}(w) = \xi^{01}(w \cdot v) - \xi^{01}(v), \qquad 
    \xi^{11}(v) = \xi^{01}(w \cdot v) - \xi^{01}(w).
\end{equation*}
Again, by Example \ref{ex:MooreFillers1Gpd}, $w \cdot v = d_1\mu^2_1(v,w)$ and we get 
\begin{equation*}
\begin{split}
        \xi^{00}(w) &= \xi^{01}(v + w - 1d_0w) - \xi^{01}(v) = \xi^{01}(w - 1d_0w),\\
        \xi^{11}(v) &= \xi^{01}(v + w - 1d_1v) - \xi^{01}(w) = \xi^{01}(v - 1d_1v).
\end{split}
\end{equation*}
Hence $\xi^{01}$ can be seen as the only independent component of $\xi$ and $\huaV^{1*}_1 = \huaV_1^*$.

By the discussion at the end of Section \ref{sec:ndual-overview}, the multiplication can now be equivalently defined without computing Moore fillers, by the multiplicativity equation \eqref{eq:MultiplicativityMostGeneral} for $r=012 \in \Delta^2_2$. Let $\phi \in \huaV^{1*}_2 \cong \Lambda^2_1(\huaV^{1*})$. Written as a horn, this is $(\phi^{12}, \phi^{01}) = (\eta, \xi)$, and $\xi \cdot \eta = m_1(\eta, \xi) = \phi^{02}$ so the equation for $012$ is 
\begin{equation*}
    (\xi \cdot \eta) ( w \cdot v) = \eta(v) + \xi (w),
\end{equation*}
as expected. The other multiplications can be defined analogously by renaming the different components of $\phi$.
\end{proof}

\begin{remark}
    As an immediate observation, this result is consistent with the dual $\VB$ groupoid construction in Definition \ref{def:VB1Dual}, as the 1-dual of a $\VS$ groupoid coincides with its 1-dual as a $\VB$ groupoid over the point.
    Moreover, if $\huaV \to M$ is a $\VB$ groupoid over the identity groupoid of a manifold $M$, then each fiber of the dual $\VB$ groupoid of $\huaV$ coincides with the 1-dual $(\huaV |_p)^{1*}$ at each point $p \in M$. 
    In particular, the restriction of any $\VB$ groupoid $\huaV \to \huaG$ to the units $\huaG_0$ is such a $\VB$ groupoid, and its fiber at any $p\in \huaG_0$ is exactly the $\VS$ 1-dual $(\huaV|_p)^{1*}$ of the fiber of $\huaV$ at $p$.
    In other words: 
    \begin{equation*}
        1^*(\huaV^{1*})|_p \cong (1^*\huaV|_p)^{1*}.
    \end{equation*}
\end{remark}

\begin{remark}\label{rem:VS1Dual-DK}
    A straightforward computation shows that the normalized complex of the 1-dual of $\huaV$ is isomorphic to 
    \begin{equation*}
        \huaV_0^* \xrightarrow{-d_1^*} (\ker d_0)^*,
    \end{equation*}
    which is precisely the 1-shifted dual $N(\huaV)^*[-1]$ of the normalized complex of $\huaV$. Therefore, in this case, $\huaV^{1*} \cong DK(N(\huaV)^*[-1])$ and the 1-dual pairing is nondegenerate ``on the nose''.
\end{remark}

\subsection{The 2-dual of a \texorpdfstring{$\VS$}{VS} 2-groupoid}\label{sec:computation-2dual}

Let $\huaV$ be a $\VS$ 2-groupoid. Following Section \ref{sec:annihilators-and-dual-kernels}, the solution spaces of the normalization conditions \eqref{eq:NormalizationMostGeneral} appearing in the 2-dual $\huaV^{2*}$ are the degree 1 2-dimensional degeneracy annihilators $O_0 = \Ann(s_0\huaV_1)$ and $O_1 = \Ann(s_1\huaV_1)$, and their intersection $O_{01} = \Ann(D_2\huaV)$, the degree 2 2-dimensional degeneracy annihilator. 

By the dual sequences of \eqref{eq:VSnDualSESd_is_i} and \eqref{eq:VSnDualSESd_iplus1s_i},  we have the isomorphisms
\begin{equation} \label{eq:O0O1}
    O_0 \cong (\ker d_0)^*\cong (\ker d_1)^*, \quad O_1 \cong (\ker d_1)^*\cong (\ker d_2)^*.
\end{equation}
Analogously, the dual of the sequence \eqref{eq:VSnDualSESp_kmu_k} for $0\le k\le 2$ gives isomorphisms
\begin{equation} \label{eq:O01}
    O_{01} \cong (\ker p^2_0)^* \cong (\ker p^2_1)^* \cong (\ker p^2_2)^*,
\end{equation}
which are the dual maps of the core projections
\begin{equation}\label{eq:VS2DualCoreProjections}
    \begin{split}
        \ggamma_0 X &= X - s_0d_1X - s_1d_2X + s_0d_2 X : \huaV_2 \to \ker p^2_0,\\
        \ggamma_1 X &= X - s_0d_0X - s_1d_2X + 1d_1d_0 X : \huaV_2 \to \ker p^2_1,\\
        \ggamma_2 X &= X - s_0d_0X - s_1d_1X + s_1d_0 X : \huaV_2 \to \ker p^2_2,
    \end{split}
\end{equation}
appearing in \eqref{eq:VSnDualSESp_kmu_k}. As in Section \ref{sec:computation-1dual}, the isomorphisms between the three different degree 2 cores are given by the restrictions of these projections. 

\begin{remark}\label{rem:VS2Dual-det-norm-comps}
    These isomorphisms translate into the useful principle that an element of $\huaV_2^*$ that satisfies a single normalization condition --- i.e. an element of $O_0$ or $O_1$ --- is determined by its evaluation on elements of any degree 1 core $\ker d^2_i$ for $i = 0,1,2$. 
    In the same way, an element of $\huaV_2^*$ that satisfies both normalization conditions --- i.e. an element of $O_{01}$ --- is determined by its evaluation on elements of any degree 2 core $\ker p^2_i$ for $i = 0,1,2$.
    We will use this repeatedly in the computation of the 2-dual.
\end{remark}

\begin{remark}
    Observe that $\ggamma_0^*|_{O_0}=\ggamma_1^*|_{O_0}$, while $\ggamma_2^*|_{O_0}$ is a different map. Similarly,  $\ggamma_1^*|_{O_1}=\ggamma_2^*|_{O_1}$, while $\ggamma_0^*|_{O_1}$ is a different map.  
\end{remark}

\begin{theorem}[2-dual of a $\VS$ 2-groupoid] \label{thm:VS2Dual}
Let $\huaV$ be a $\VS$ 2-groupoid. Then its 2-dual $\huaV^{2*}$ is
    \begin{equation} \label{eq:2-dual}
                    \huaV_2^* \times_{O_{01}} 
                    \left(O_0 \times_{O_{01}} O_1\right) \times_{O_{01}} O_0 \rightthreearrows O_0 \times_{O_{01}} O_1 \rightrightarrows O_{01},
    \end{equation}
where the fiber products at levels 1 and 2 are 
\begin{equation}\label{eq:VS2DualElementsLv12Def}
    \begin{split}
        \huaV^{2*}_1 &= \{(\eta^{001}, \eta^{011}) \in O_0 \oplus O_1 \mid  \eta^{001}(\ggamma_1X) = \eta^{011}(\ggamma_1X), \quad \forall X \in \huaV_2\},\\
        \huaV^{2*}_2 &= \{(\phi^{012}, \phi^{112}, \phi^{122}, \phi^{001}) \in \huaV_2^* \oplus 
        \left(O_0 \oplus O_1\right) \oplus O_0| \\
        &\qquad \qquad  \phi^{012}(\ggamma_0 X) = \phi^{112}(\ggamma_1 X) = \phi^{122}(\ggamma_1 X) \\
        &\qquad \qquad \text{ and } \phi^{012}(\ggamma_2 X) = \phi^{001}(\ggamma_1 X), \quad \forall X \in \huaV_2\},
    \end{split}
\end{equation}
and this is equipped with the following face and degeneracy maps\footnote{Notice that since $\eta^{011}\in O_1=\Ann(s_1\huaV_1)$, $\eta^{011}(s_1d_2(X))=0$. Thus we have a simplification (rather than a typo) in \eqref{eq:vs-2dual-face-degen-1}. This is similar for other simplifications in \eqref{eq:vs-2dual-face-degen-1} and \eqref{eq:vs-2dual-face-degen-2}. }:
\begin{equation}\label{eq:vs-2dual-face-degen-1}
    \begin{aligned}
        \widecheck{d}_0^1(\eta^{001},\eta^{011})(X)
        &= \eta^{011}(\ggamma_0 X) = \eta^{011} (X - s_0d_1 X + s_0d_2 X),
        \\
        \widecheck{d}_1^1(\eta^{001},\eta^{011})(X) 
        &= \eta^{001}(\ggamma_2 X)= \eta^{001}(X - s_1d_1 X + s_1d_0 X),
        \\
        \widecheck{s}_0^0(\epsilon) &= (\epsilon, \epsilon)
    \end{aligned}
\end{equation}
for all $(\eta^{001}, \eta^{011}) \in O_0 \times_{O_{01}} O_1$, $\epsilon \in O_{01}$, $X \in \huaV_2$, and
\begin{equation}\label{eq:vs-2dual-face-degen-2}
    \begin{aligned}
        \widecheck{d}^2_0 (\phi^{012}, \phi^{112}, \phi^{122}, \phi^{001}) 
        &= (\phi^{112}, \phi^{122}),
        \\
        \widecheck{d}^2_1 (\phi^{012}, \phi^{112}, \phi^{122}, \phi^{001}) 
        &= (\phi^{012} - \phi^{012}s_0d_0 + \phi^{001}s_1d_2,\\
        &\qquad \qquad \phi^{012} - \phi^{012}s_1d_2 + \phi^{122}s_0d_0),
        \\
        \widecheck{d}^2_2 (\phi^{012}, \phi^{112}, \phi^{122}, \phi^{001})
        &= (\phi^{001}, \phi^{012} - \phi^{012}s_1d_1 + \phi^{112}s_1d_0),
        \\
        \widecheck{s}^1_0(\eta^{001}, \eta^{011})
        &= (\eta^{001}, (\eta^{001}, \eta^{011}), \widecheck{d}_1(\eta^{001},\eta^{011})),
        \\
        \widecheck{s}^1_1(\eta^{001}, \eta^{011})
        &= (\eta^{011}, \widecheck{s}_0\widecheck{d}_0(\eta^{001}, \eta^{011}), \eta^{001}),
    \end{aligned}
\end{equation}
for all $(\phi^{012}, \phi^{112}, \phi^{122}, \phi^{001}) \in \huaV^{2*}_2$ and $(\eta^{001}, \eta^{011}) \in \huaV^{2*}_1$. 
The multiplication $\widecheck{m}_1$ is defined by the property that 
\begin{equation}\label{eq:mul-pairing}
\begin{split}
    \langle \widecheck{m}_1 (\phi, \phi', \phi''), m_1(W, Y, Z) \rangle &=
    \langle \phi \square \phi' \phi'', W \square Y Z \rangle\\
    &= \langle \phi, W \rangle 
    + \langle \phi', Y \rangle 
    - \langle \phi'', Z \rangle,
\end{split}
\end{equation}
for any $(\phi, \phi', \phi'') \in \Lambda^3_1(\huaV^{2*})$ and any $(W, Y, Z) \in \Lambda^3_1(\huaV)$.\footnote{Note that this equation relates only the interior component $(\phi\square\phi'\phi'')^{012}$ of the product to the interior components of the factors. As previously explained in Section \ref{sec:ndual-overview}, the other components can be inferred by the simplicial identities. We write these in \eqref{eq:m_1}.} 
\end{theorem}

\begin{proof}
Similarly to the computation of the 1-dual, we follow Section \ref{sec:ndual-overview} and compute $\huaV^{2*}$ level by level, by solving the linear equations given by the multiplicativity  \eqref{eq:MultiplicativityMostGeneral} and normalization \eqref{eq:NormalizationMostGeneral} conditions. 
This is however more complicated than in the 1-dual case, as there are many instances where multiple equations determine the same variable, so we will also need to check that the equations do not over-determine the solutions, that is the solution space is not empty. 
For ease of reading, we organize the proof into subsections.

\subsubsection{Level 0}\label{sec:level0} 
For $\epsilon \in \huaV^{2*}_0$, the only 2-simplex in $\Delta^0_2$ is $000$, so we have only one component $\epsilon = \epsilon^{000}$. Since $000 = s_0 00 = s_1 00$, we have two normalization conditions \eqref{eq:NormalizationMostGeneral}, and we have that $\epsilon \in \Ann(D_1\huaV_1) = O_{01}$. There is only one multiplicativity condition, which is \eqref{eq:MultiplicativityMostGeneral} for $0000$, which reads
\begin{equation*}
    \epsilon (W \square Y Z) = \epsilon(W) + \epsilon(Y) - \epsilon(Z).
\end{equation*}
But this is already implied  by  \eqref{eq:TriMultOverPointAll} and the normalization condition of $\epsilon$. Therefore $\huaV^{2*}_0=O_{01}$.

\subsubsection{Level 1}\label{sec:level1} 
Any $\eta \in \huaV^{2*}_1$ consists of four components $\eta^{000}, \eta^{001}, \eta^{011}, \eta^{111}$, which are normalized in the following way:
    \begin{equation*}
    \begin{split}
        \eta^{000}, \eta^{111} \in \Ann(D_1\huaV_1)=O_{01},\\
        \eta^{001} \in \Ann(s_0\huaV_1)=O_0,
        \quad \eta^{011} \in \Ann(s_1\huaV_1)=O_1.
    \end{split}
    \end{equation*}

The multiplicativity conditions are
    \begin{equation}\label{eq:VS2DualGenericMult}
        \eta^{d_1r}(W\square Y Z) = \eta^{d_0r}(W) + \eta^{d_2r}(Y) - \eta^{d_3r}(Z)
    \end{equation}
for all $(W,Y,Z) \in \Lambda^3_1(\huaV)$, and $r = 0000, 0001, 0011, 0111, 1111$. The two equations relative to $r = 0000$ and $r = 1111$ are automatically satisfied by the same argument as in level 0.

For the others, we make use of Remark \ref{rem:VS2Dual-det-norm-comps}, by which $\eta^{000}, \eta^{111} \in O_{01}$ are determined by their evaluation on elements in $\ker p^2_j$ for $0 \le j \le 2$.
For $r=0001$, since for any $Z\in \ker p^2_2$ we get $0\square Z Z = 0$ by \eqref{eq:TriMultOverPointAll},  \eqref{eq:VS2DualGenericMult} reads
    \begin{equation*}
            0 = \eta^{001}(Z) - \eta^{000}(Z) \quad
            \iff \quad \eta^{000}(Z) = \eta^{001}(Z), \quad \forall Z \in \ker p^2_2.
    \end{equation*}
Hence by Remark \ref{rem:VS2Dual-det-norm-comps}, $\eta^{000}$ is entirely determined by $\eta^{001}$.  To see what $\eta^{000}$ is when evaluated on a generic $X\in \huaV_2$ we use the projections in \eqref{eq:VS2DualCoreProjections} and obtain that
    \begin{equation}\label{eq:VS2Dual000}
        \eta^{000}(X) = \eta^{000}(X - s_0d_0X - s_1d_1X + s_1d_0X)= \eta^{001}(\ggamma_2 X),
    \end{equation}
by the normalization condition on $\eta^{000}$. It is easy to see that imposing \eqref{eq:VS2Dual000} is equivalent to the multiplicativity condition \eqref{eq:VS2DualGenericMult} for $r=0001$ because substituting \eqref{eq:VS2Dual000} in \eqref{eq:VS2DualGenericMult} does not impose extra conditions on $\eta^{001}$.

A symmetric argument applies to $r=0111$, by using the duality principle in Remark \ref{rem:front-to-back}. Thus, we have that 
\begin{equation}\label{eq:VS2Dual111}
    \eta^{111}(X) = \eta^{011}(\ggamma_0(X))
\end{equation}
is equivalent to the multiplicativity condition \eqref{eq:VS2DualGenericMult} for $r=0111$ and $\eta^{111}$ is completely determined by $\eta^{011}$. 

Lastly, we take care of \eqref{eq:VS2DualGenericMult} for $r=0011$. For $Y\in \ker p^2_1$, $0\square Y0 = Y$, and this equation reads $\eta^{011}(Y) = \eta^{001}(Y)$. Again, by the normalization condition,  \eqref{eq:VS2DualGenericMult} for $r=0011$ implies that
    \begin{equation}\label{eq:VS2DualLevel1FiberProductCondition}
        \eta^{011}(X - s_0d_0 X) = \eta^{001}(X - s_1d_2 X) \iff \eta^{011}(\ggamma_1 X) = \eta^{001}(\ggamma_1 X), 
    \end{equation} for any $X \in \huaV_2$. 
Conversely, with \eqref{eq:VS2DualLevel1FiberProductCondition}, using normalization conditions,  \eqref{eq:VS2DualGenericMult} for $r=0011$ becomes automatic:
    \begin{equation*}
        \begin{split}
            \eta^{011}(W + Y - Z - s_0d_1W + s_0d_0Z) 
            &= \eta^{011}(W) + \eta^{001}(Y) - \eta^{001}(Z)\\
\iff            \eta^{011}(Y - s_0d_0Y) + \eta^{011}(- Z + s_0d_0 Z) &= \eta^{001}(Y) - \eta^{001}(Z)\\
\iff            \eta^{001}(Y - \cancel{s_1d_2Y}) + \eta^{001}(- Z + \cancel{s_1d_2 Z}) &= \eta^{001}(Y) - \eta^{001}(Z)
        \end{split}
    \end{equation*}
Thus \eqref{eq:VS2DualLevel1FiberProductCondition} is equivalent to \eqref{eq:VS2DualGenericMult} for $r=0011$.

In summary, each $\eta$ is determined by the pair of components $(\eta^{001},\eta^{011})$ which satisfies \eqref{eq:VS2DualLevel1FiberProductCondition}. That is, $\huaV^{2*}_1$ is the fiber product 
    \[\begin{tikzcd}[ampersand replacement=\&]
	{\huaV^{2*}_1 \cong O_0\times_{(\ker p^2_1)^*}O_1} \& {O_1} \\
	{O_0} \& {(\ker p^2_1)^*}
	\arrow[from=1-1, to=2-1]
	\arrow["\lrcorner"{anchor=center, pos=0.125}, draw=none, from=1-1, to=2-2]
    \arrow["{\ggamma_1^*=(id-s_0d_0)^*}", from=1-2, to=2-2]
	\arrow[from=1-1, to=1-2]
	\arrow["{\ggamma_1^*=(id-s_1d_2)^*}"', from=2-1, to=2-2]
    \end{tikzcd}\]
    with face maps 
    \begin{equation*}
        \begin{split}
            \widecheck{d}_0(\eta^{001},\eta^{011}) &= \eta^{111} = \eta^{011} \circ \ggamma_{0}\\
            \widecheck{d}_1(\eta^{001},\eta^{011}) &= \eta^{000} = \eta^{001} \circ \ggamma_{2},
        \end{split}
    \end{equation*}
by \eqref{eq:VS2Dual000}, \eqref{eq:VS2Dual111}.  The degeneracy map $\widecheck{s}_0 : \huaV^{2*}_0 \to \huaV^{2*}_1$ is the obvious diagonal map.

\subsubsection{Level 2 --- Solving Equations}\label{sec:level2} An element $\phi\in \huaV^{2*}_2$ consists of 10 components, one for each 2-simplex of $\Delta^2$. Nine of them are normalized in the sense that
    \begin{equation*}
        \phi^{iii} \in \Ann(D_1\huaV)=O_{01},
        \quad \phi^{iij} \in \Ann(s_0\huaV_1)=O_0,
        \quad \phi^{ijj} \in \Ann(s_1\huaV_1)=O_1,
    \end{equation*}
for all $0\le i\le j \le 2$. The component $\phi^{012}$ is the only one with no normalization conditions. 
    
There are 15 multiplicativity conditions from \eqref{eq:MultiplicativityMostGeneral}, which we write for $m_1$ as 
\begin{equation}\label{eq:VS2DualGenericMult-c}
    \phi^{d_1r}(W\square Y Z) = \phi^{d_0r}(W) + \phi^{d_2r}(Y) - \phi^{d_3r}(Z)
\end{equation}
for all $(W,Y,Z) \in \Lambda^3_1(\huaV)$ and for all $r=ijkl$ with $0\le i\le j\le k\le l \le 2$.  

Those for $r = iiii$ with $0\le i\le 2$ are automatically satisfied, by the same argument as for level 0;    those for $r= ijjj, iijj, ijjj$ for $0\le i\le j \le 2$ can also be treated as for level 1, and they impose that each pair $(\phi^{iij},\phi^{ijj})$ forms an element of $\huaV^{2*}_1$.
Notice that each pair $(0001,0002)$, $(0111,1112)$, and $(0222,1222)$ gives two possible ways to determine $000$, $111$ and $222$, respectively. We will check later that they give consistent results by using other multiplicativity conditions in \ref{sec:consist}. Now it only remains to solve the three equations for $r = 0012, 0112, 0122$. 

We begin with $r=0012$. Then \eqref{eq:VS2DualGenericMult-c} reads
    \begin{equation}\label{eq:VS2DualMult0012}
        \phi^{012}(W\square Y Z) = \phi^{012}(W) + \phi^{002}(Y) - \phi^{001}(Z),
    \end{equation} for all $(W,Y,Z)\in \Lambda^3_1(\huaV)$.
As before, since $\phi^{002} \in O_0$, by picking $Y\in \ker d_0^2$, we can determine $\phi^{002}$ completely (see Remark \ref{rem:VS2Dual-det-norm-comps}). By plugging in $W=0$ and $Y\in \ker d_0^2$,  \eqref{eq:VS2DualMult0012} becomes
    \begin{equation}\label{eq:VS2DualMult0012with0}
        \phi^{012}(0\square YZ) = \phi^{012}(Y - Z + s_1d_1Z) = \phi^{002}(Y) - \phi^{001}(Z),
    \end{equation}
where $Z \in \ker d_0^2$ is any (2,1)-horn filler of $(0, d_2Y)$. Observe that $s_1d_2Y$ is one such horn filler, since $d_0d_2Y=d_1d_0Y = 0$. A general filler is given by $s_1d_2Y + k_1$, with $k_1\in \ker p^2_1$. In other words the space of (2,1)-fillers is an affine space modelled over $\ker p^2_1$. Since $s_1d_1(s_1d_2Y + k_1) = s_1d_2Y + s_1d_1 k_1$, by using \eqref{eq:TriMultOverPointAll}, \eqref{eq:VS2DualMult0012with0} becomes
    \begin{equation*}
        \phi^{002}(Y) = \phi^{012}(Y - k_1 + s_1d_1k_1) + \phi^{001}(s_1d_2Y + k_1).
    \end{equation*}
Thus for any $X\in \huaV_2$, since $X-s_0d_0X \in \ker d_0^2$ and $\phi^{002}\in O_0$, we get
    \begin{equation}\label{eq:VS2Dual002GenericFiller}
    \begin{split}
        \phi^{002}(X) &= \phi^{002}(X-s_0d_0X)\\ 
        &= \phi^{012}(X - s_0d_0X - k_1 + s_1d_1k_1) + \phi^{001}(s_1d_2X + k_1), 
    \end{split}
    \end{equation}
where, for the last term, we used the fact that $s_1d_2(s_0d_0X) =s_0s_0 d_1d_0X$ and $\phi^{001}\in O_0$.    
Now we need to see what conditions are imposed on $\phi^{012}$ and $\phi^{001}$ by \eqref{eq:VS2DualMult0012} after inserting \eqref{eq:VS2Dual002GenericFiller}.
With \eqref{eq:TriMultOverPointAll}, \eqref{eq:VS2DualMult0012} reads
\begin{equation*}
\begin{split}
    \phi^{012}(W + Y &- Z - s_0d_1W + s_0d_0Z - s_1d_0Z +s_1d_1Z)\\
    &= \phi^{012}(W) + \phi^{002}(Y) - \phi^{001}(Z), 
\end{split}
\end{equation*}  for all $(W, Y, Z)\in \Lambda^3_1(\huaV)$. By inserting \eqref{eq:VS2Dual002GenericFiller}
we get
    \begin{equation*} \label{eq:c012-2}
        \begin{split}
            \phi^{012}(Y - Z - s_0d_1W + s_0d_0Z - s_1d_0Z +s_1d_1Z)
            &= \phi^{012}(Y - s_0d_0Y - k_1 + s_1d_1k_1) \\
            &\qquad + \phi^{001}(s_1d_2Y + k_1 - Z),
        \end{split}
    \end{equation*}        
with $k_1$ an arbitrary element in $\ker p^2_1$. Because $d_1W = d_0Y$, $d_2Y = d_2Z$ and $\phi^{001}(\ggamma_1 X) = \phi^{001}(X - s_1d_2 X)$,this further simplifies to      
\begin{equation*}          
        \begin{split}    
            \phi^{012}(Z - s_0d_0Z + s_1d_0Z - s_1d_1Z - k_1 + s_1d_1k_1)
            &= \phi^{001}(Z - s_1d_2Z - k_1)\\
          \iff  \phi^{012}(\ggamma_2(Z - k_1)) 
            &= \phi^{001}(\ggamma_1(Z - k_1)).
        \end{split}
    \end{equation*}
This means that \eqref{eq:VS2DualMult0012} holds if and only if \eqref{eq:VS2Dual002GenericFiller} holds for all $X\in \huaV_2$ and $k_1 \in \ker p^2_1$, and 
    \begin{equation}\label{eq:VS2DualLevel2FiberProductCondition1}
        \phi^{012}(\ggamma_2 X) = \phi^{001}(\ggamma_1 X), \quad \forall X\in \huaV_2.
    \end{equation}
In fact, to use \eqref{eq:VS2Dual002GenericFiller} to determine $\phi^{002}$, we need to show the right-hand side of \eqref{eq:VS2Dual002GenericFiller} does not depend on the choice of $k_1 \in \ker p^2_1$. However, this follows precisely from \eqref{eq:VS2DualLevel2FiberProductCondition1} for $X=k_1$. In summary, \eqref{eq:VS2DualMult0012} is equivalent to \eqref{eq:VS2Dual002GenericFiller} for any $Y$ and any filler. Because this must hold for any filler, \eqref{eq:VS2DualLevel2FiberProductCondition1} must also hold.

The front-to-back symmetric case of $r=0122$ can be treated analogously by Remark \ref{rem:front-to-back}. The multiplicativity condition for $r=0122$ reads, for any $(W, X, Y) \in \Lambda^3_2(\huaV)$,
\begin{equation}\label{eq:VS2DualMult0122with0}
    \phi^{022}(X) = \phi^{012}(W X\square 0) + \phi^{122}(W) = \phi^{012}(-W + X +s_0d_1W) +\phi^{122}(W).  
\end{equation}
This is equivalent to imposing
\begin{equation}\label{eq:VS2DualLevel2FiberProductCondition2}
    \phi^{012}(\ggamma_0 X) = \phi^{122}(\ggamma_1 X), \quad \forall X \in \huaV_2
\end{equation}
and defining $\phi^{022}$ in terms of $\phi^{012}$ and $\phi^{001}$ by
\begin{equation}\label{eq:VS2Dual022GenericFiller}
    \phi^{022}(X) = \phi^{012}(X - k_1 + s_0d_1k_1) +\phi^{122}(s_0d_0X + k_1), \quad \forall X\in \huaV_2
\end{equation}    
for any choice of $k_1\in \ker p^2_1$. As before, \eqref{eq:VS2DualLevel2FiberProductCondition2} for $X=k_1$ is equivalent to \eqref{eq:VS2Dual022GenericFiller} being independent of the choice of $k_1$. 

Moving on to $r=0112$, the multiplicativity condition reads 
    \begin{equation}\label{eq:VS2DualMult0112}
        \phi^{012}(W \square Y Z) = \phi^{112}(W) + \phi^{012}(Y) - \phi^{011}(Z), \quad \forall (W,Y,Z)\in \Lambda^3_1(\huaV).
    \end{equation} 
Once more, we want to determine $\phi^{011}$ from $\phi^{012}$ and $\phi^{112}$. 
Take $Z \in \ker d^2_2$, then \eqref{eq:VS2DualMult0112} implies that
    \begin{equation*}
        \begin{split}
            \phi^{011}(Z) &= - \phi^{012}(W \square 0 Z) + \phi^{112}(W) \\
            &= \phi^{012}( - W + Z - s_0d_0Z + s_1d_0Z - s_1d_1Z) + \phi^{112}(W),
        \end{split}
    \end{equation*}
for any (2,0)-horn filler $W\in \ker d^2_1$ of $(0, d_0Z)$. Since $s_1d_0 Z - s_0d_0 Z$ is one such horn filler, an arbitrary filler is of the form $s_1d_0 Z - s_0d_0 Z + k_0$ with $k_0 \in \ker p^2_0$. Hence we have
    \begin{equation*}
        \phi^{011}(Z) = \phi^{012}(Z - s_1d_1Z - k_0) + \phi^{112}(s_1d_0Z + k_0), \quad \forall k_0 \in \ker p^2_0,
    \end{equation*} 
since $\phi^{112}\in O_0 $. 
Again for all $X\in \huaV_2$, since  $X-s_1d_2X \in \ker d^2_2$ and $\phi^{011}\in O_1$, we have 
    \begin{equation}\label{eq:VS2Dual011GenericFiller}
        \phi^{011}(X) = \phi^{011}(X - s_1d_2X) = \phi^{012}(X - s_1d_1X - k_0) +\phi^{112}(s_1d_0X + k_0),
    \end{equation}
where, for the last term, we use again the fact that $s_1d_0(s_1d_2X) =s_0s_0 d_0 d_2X$ and $\phi^{112}\in O_0$. 
With \eqref{eq:TriMultOverPointAll}, \eqref{eq:VS2DualMult0112} reads
    \begin{equation*}
    \begin{split}
        \phi^{012}(W + Y &- Z - s_0d_1W + s_0d_0Z - s_1d_0Z + s_1d_1Z)\\ 
        &= \phi^{112}(W) + \phi^{012}(Y) - \phi^{011}(Z), 
    \end{split}
    \end{equation*}  
for all $(W, Y, Z)\in \Lambda^3_1(\huaV)$.  
By inserting \eqref{eq:VS2Dual011GenericFiller} in this, we get
    \begin{equation*}
            \begin{split}
                \phi^{012}(W - Z &- s_0d_1W + s_0d_2W - s_1d_2W + s_1d_1Z)\\ 
                &= \phi^{112}(W) - \phi^{012}(Z - s_1d_1Z - k_0) - \phi^{112}(s_1d_0 Z + k_0),
            \end{split}
    \end{equation*}        
with $k_0$ an arbitrary element in $\ker p^2_0$. Because $d_0Z = d_2W$, and $\phi^{112}\in O_0$, this further simplifies to      
    \begin{equation*}       
        \begin{split} 
            \phi^{012}(W - s_0d_1W + s_0d_2W - s_1d_2W - k_0)
            &= \phi^{112}(W - s_1d_2W - k_0)\\
            \iff  \phi^{012}(\ggamma_0(W - k_0)) 
            &= \phi^{112}(\ggamma_1(W - k_0)).   
        \end{split}
    \end{equation*}
This means that \eqref{eq:VS2DualMult0112} holds if and only if \eqref{eq:VS2Dual011GenericFiller} holds for all $X\in \huaV_2$ and $k_0 \in \ker p^2_0$ and 
    \begin{equation}\label{eq:VS2DualLevel2FiberProductCondition3}
        \phi^{012}(\ggamma_0 X) = \phi^{112}(\ggamma_1 X), \quad \forall X\in \huaV_2.
    \end{equation}
The latter is precisely the fact that $\phi^{011}$ is determined by \eqref{eq:VS2Dual011GenericFiller} without depending on the choice of filler $k_0 \in \ker p^2_0$.
In summary, \eqref{eq:VS2DualMult0112} is equivalent to \eqref{eq:VS2DualLevel2FiberProductCondition3} and the fact that $\phi^{011}$ is determined by $\phi^{012}$ and $\phi^{122}$ via \eqref{eq:VS2Dual011GenericFiller}.

As a side note, \eqref{eq:VS2DualLevel2FiberProductCondition3} is also already implied by \eqref{eq:VS2DualLevel2FiberProductCondition2} and the multiplicativity condition for $r=1122$, which, by the discussion in \ref{sec:level1}, is equivalent to $\phi^{112}(\ggamma_1 X) = \phi^{122}(\ggamma_1 X)$.

\subsubsection{Level 2 --- Consistency}\label{sec:consist} Since some of the components are overdetermined, one needs to check consistency of the following conditions:
\begin{itemize}
\item     It is equivalent to determine $\phi^{000}$ by imposing \eqref{eq:VS2DualGenericMult-c} for $r = 0001$ and  for $r=0002$. This follows from
\begin{equation}\label{eq:001-002}
    \phi^{001}(\ggamma_2 X) = \phi^{002}(\ggamma_2 X), \quad \forall X \in \huaV_2.
\end{equation}
\item Analogously, it is equivalent to determine $\phi^{111}$ and $\phi^{222}$ in the two possible ways. This follows from 
\begin{equation*}
    \phi^{011}(\ggamma_0 X)= \phi^{112}(\ggamma_2 X), \quad \phi^{022}(\ggamma_0 X) = \phi^{122}(\ggamma_0 X), \quad \forall X \in \huaV_2.
\end{equation*}    
\item $\phi^{011}$ satisfies the multiplicativity condition for $r = 0011$ together with $\phi^{001}$,  when determined using the multiplicativity condition for $r=0122$. This follows from 
\begin{equation}\label{eq:001-011}
    \phi^{001} (\ggamma_1 X) = \phi^{011} (\ggamma_1 X), \quad \forall X \in \huaV_2.
\end{equation}
\item Analogously, $\phi^{002}$ and $\phi^{022}$ that we determined with \eqref{eq:VS2Dual002GenericFiller} and \eqref{eq:VS2Dual022GenericFiller} also satisfy the multiplicativity condition for $r=0022$. This follows from 
\begin{equation}\label{eq:002-022}
    \phi^{002} (\ggamma_1 X) = \phi^{022} (\ggamma_1 X), \quad \forall X \in \huaV_2.
\end{equation}
\end{itemize}
All the equations above follow from a straightforward calculation. Here we give the proof for the first one, \eqref{eq:001-002}. As the others can be obtained by similar arguments, we leave them to the interested reader.\footnote{
These consistency equations can be used to see more interesting facts, for example, by \eqref{eq:VS2Dual002GenericFiller} and \eqref{eq:VS2Dual022GenericFiller}, \eqref{eq:002-022} is $\phi^{012}(\ggamma_1 X) = \phi^{002} (\ggamma_1 X) = \phi^{022} (\ggamma_1 X)$. Together with \eqref{eq:VS2DualLevel2FiberProductCondition1}, \eqref{eq:001-011}, \eqref{eq:VS2DualLevel2FiberProductCondition2}, and \eqref{eq:VS2DualLevel2FiberProductCondition3}, this can be summarized in the identity $\ggamma_i^*\phi^{012} = \ggamma_1^*\phi^{s_jd_i012}$ for any $i=0,1,2$, $j=0,1$. This can be used to rewrite $\huaV^{2*}_2$ as a different fiber product, by changing which components of $\phi$ are taken as independent variables, as in Remark \ref{rem:VS2Dual-other-choice-vars}.}

Using \eqref{eq:VS2DualCoreProjections} and normalization conditions, \eqref{eq:VS2DualLevel2FiberProductCondition1} is equivalent to 
\begin{equation*}
    \phi^{012}(X-s_0d_0X)+\phi^{001}(s_1d_2X)=\phi^{012}(s_1d_1 X -s_1d_0X)+\phi^{001}(X), \quad \forall X \in \huaV_2.
\end{equation*}
Then by combining this with \eqref{eq:VS2Dual002GenericFiller} for $k_1=0$, we get
\begin{equation*}
    \phi^{002}(X) = \phi^{012}(s_1d_1X-s_1d_0X) + \phi^{001}(X) , \quad \forall X \in \huaV_2.
\end{equation*}
Since the image of $\ggamma_2$ is the intersection of the kernels of $d_1$ and $d_0$, by replacing $X$ with $\ggamma_2X$ we obtain \eqref{eq:001-002}. 

\subsubsection{Level 2 --- Structure Maps and Multiplication}
In the above two subsections, we have shown that an element $\phi\in \huaV^{2*}_2$ depends only on the components $\phi^{012}$, $\phi^{112}$, $\phi^{122}$, and $\phi^{001}$. Thus $\huaV^{2*}_2 \subseteq \huaV^*_2 \oplus (O_0)^2 \oplus O_1$. Furthermore, by \eqref{eq:VS2DualLevel2FiberProductCondition1}, \eqref{eq:VS2DualLevel2FiberProductCondition2}, \eqref{eq:VS2DualLevel2FiberProductCondition3}, and the previous subsection, $\huaV^{2*}_2$ is exactly the fiber product
    \begin{equation*}
        \begin{split}
            \huaV^{2*}_2 &= \huaV^*_2 \times_{(\ggamma_0^*, \ggamma_2^*), (O_{01})^2, (\ggamma_1^*\circ pr_1, \ggamma_1^*)} \times \left(\left( O_0 \times_{\ggamma_1^*, O_{01}, \ggamma_1^*} O_1 \right) \times O_0 \right)\\
            &= \{(\phi^{012}, (\phi^{112}, \phi^{122}), \phi^{001}) \mid \phi^{012}\ggamma_2 = \phi^{001}\ggamma_1, \quad 
            \phi^{012}\ggamma_0 =\phi^{122}\ggamma_1= \phi^{112}\ggamma_1
           \}.
        \end{split}
    \end{equation*}

    By the definition of the mapping space, the simplicial maps are described by \eqref{eq:mapping-set-maps}. For example, the face maps are
    \begin{equation*}
        \begin{split}
            \widecheck{d}_0(\phi^{012}, (\phi^{112}, \phi^{122}), \phi^{001}) 
            &= (\phi^{112}, \phi^{122})\\
            \widecheck{d}_1(\phi^{012}, (\phi^{112}, \phi^{122}), \phi^{001}) 
            &= (\phi^{002}, \phi^{022}) \\
            &= \left(\phi^{012} \circ (id - s_0d_0) + \phi^{001}s_1d_2 ,\right.\\ 
            &\left.\qquad\qquad\qquad\phi^{012} \circ (id - s_1d_2) + \phi^{122}s_0d_0 \right)\\
            \widecheck{d}_2(\phi^{012}, (\phi^{112}, \phi^{122}), \phi^{001}) 
            &= (\phi^{001}, \phi^{011})
            = \left(\phi^{001}, 
            \phi^{012} \circ (id - s_1d_1) + \phi^{112}s_1d_0 \right)
        \end{split}
    \end{equation*}
The degeneracy maps can be similarly computed, so we refer to the statement for their explicit description. 

Finally, a (3,1)-horn in $\huaV^{2*}$ can be written in components as $p^3_1\psi$ for a unique 3-simplex $\psi \in \huaV^{2*}_3$:
\begin{equation*}
\begin{split}
    p^3_1\psi=(&(\psi^{123},(\psi^{223},\psi^{233}),\psi^{112}),(\psi^{013},(\psi^{113},\psi^{133}),\psi^{001}),\\
    &(\psi^{012},(\psi^{112},\psi^{122}),\psi^{001})).
\end{split}
\end{equation*}
Following the discussion at the end of Section \ref{sec:ndual-overview}, the multiplication $\widecheck{m}_1$ is given by
\begin{equation} \label{eq:m_1}
    \widecheck{m}_1(p^3_1\psi) = \widecheck{d}_1 \psi  = (\psi^{023}, (\psi^{223},\psi^{233}), \psi^{002}),
\end{equation}
where $(\psi^{223},\psi^{233})$ and $\psi^{002}$ are determined as before from the horn data by the simplicial identities, while $\psi^{023}$ can be obtained from the multiplicativity condition for $r=0123$. That is, for any $(W,Y,Z) \in \Lambda^3_1(\huaV)$, we have
\begin{equation}
    \psi^{023}(W\square YZ) = \psi^{123}(W) + \psi^{013}(Y) - \psi^{012}(Z),
\end{equation}
which is the expected formula. 
\end{proof}

\begin{remark}\label{rem:VS2Dual-other-choice-vars}
    The choice of including $\phi^{112}, \phi^{122}$ and $\phi^{001}$ as ``given data'' of a certain element $\phi \in \huaV^{2*}_2$ is not unique. In fact, there are six possible combinations of choices coming from the fact that the components $001$ and $002$ are related by the multiplicativity condition relative to $r=0012$, $011$ and $112$ are related by the one for $r=0112$, while $122$ and $022$ by the one for $r=0122$. These all determine  isomorphic but different expressions of $\huaV^{2*}_2$ as a fiber product.
\end{remark} 

\subsection{The normalized complex of the 2-dual}

We compute here the normalized complex of the 2-dual for future reference. Recall that for the 1-dual of a $\VS$ 1-groupoid, as computed in Remark \ref{rem:VS1Dual-DK}, we have an isomorphism $N(\huaV^{1*}) \cong N(\huaV)^*[-1]$. 
On the contrary, in the case of a $\VS$ 2-groupoid $\huaV$, $N(\huaV^{2*})$ is only chain homotopy equivalent to $N(\huaV)^*{[-2]}$, as a consequence of Theorem \ref{thm:EZ-thm-hom} and Theorem \ref{thm:ndual-pairing-hom-nondeg}. 
The case of $n=2$ is indeed the first case in increasing order of $n$ where the $n$-dual pairing of a $\VS$ $n$-groupoid is only homologically nondegenerate and not nondegenerate on the nose.

\begin{proposition}\label{prop:vs-2dual-norm-cplx}
Let $\huaV$ be a $\VS$ 2-groupoid. Then, the normalized complex of the 2-dual of $\huaV$, $N(\huaV^{2*})$ is isomorphic to the chain complex 
\begin{equation}\label{eq:vs-2dual-norm-cplx}
\begin{tikzcd}[ampersand replacement=\&,,row sep = 0ex
    ,/tikz/column 1/.append style={anchor=base east}
    ,/tikz/column 5/.append style={anchor=base west}]
	{\huaV_0^* \oplus (\ker d^1_0)^* } \& {(\ker d^1_0)^* \oplus (\ker d^1_0)^*} \& {(\ker p^2_2)^*} \\
	{(\theta, \nu)} \& {(\nu, \nu + \theta d_1)} \\
	\& {(\xi_1, \xi_2)} \& {\xi_2 d_2 - \xi_1 d_2}
	\arrow[from=1-1, to=1-2]
	\arrow[from=1-2, to=1-3]
	\arrow[maps to, from=2-1, to=2-2]
	\arrow[maps to, from=3-2, to=3-3]
\end{tikzcd}
\end{equation}
concentrated in degrees 0,1 and 2. 
\end{proposition}

\begin{proof}
    We define the explicit isomorphism $\Psi$: at level 0, this is the identity of $N(\huaV^{2*})_0 = (\ker p^2_2)^*$. At level 1 we have
    \begin{equation*}
    \begin{split}
            \Psi_1: &N(\huaV^{2*})_1 = \ker\widecheck{d}^1_0 \longrightarrow (\ker d^1_0)^* \oplus (\ker d^1_0)^*\\
            &\Psi(\eta^{001}, \eta^{011}) = (\eta^{001}s_1, \eta^{011}s_0),
    \end{split}
    \end{equation*}
    with inverse 
    \begin{equation*}
    \begin{split}
            \Psi_1^{-1}: &(\ker d^1_0)^* \oplus (\ker d^1_0)^*\longrightarrow N(\huaV^{2*})_1 = \ker \widecheck{d}^1_0\\
            &\Psi^{-1}(\xi_1, \xi_2) = (\xi_1 \ggamma_R d_2 - \xi_2 \ggamma_R \partial, \xi_2 \ggamma_R d_1 - \xi_2 \ggamma_R d_2),
    \end{split}
    \end{equation*}
    where $\partial = d_0 - d_1 + d_2$ is the simplicial boundary map at level 2 and $\ggamma_R$ is the 1-dimensional core projection $\ggamma_R = id - 1d_0: \huaV_1 \to \ker d^1_0$. \footnote{We denote it by $\ggamma_R$ instead of $\ggamma_0$ to distinguish it from the 2-dimensional core projections $\ggamma_0, \ggamma_1$ and $\ggamma_2$ that appear later in the proof.}
    At level 2 we have 
        \begin{equation*}
    \begin{split}
            \Psi_2: &N(\huaV^{2*})_2 = \ker\widecheck{p}^2_2 \longrightarrow \huaV_0^* \oplus (\ker d^1_0)^*\\
            &\Psi(\phi^{012},0,0,\phi^{001}) = (\phi^{012}1, \phi^{001}s_1),
    \end{split}
    \end{equation*}
    with inverse 
    \begin{equation*}
    \begin{split}
            \Psi_2^{-1}: &\huaV_0^* \oplus (\ker d^1_0)^* \longrightarrow N(\huaV^{2*})_2 = \ker\widecheck{p}^2_2\\
            &\Psi^{-1}(\theta, \nu) = (\theta d_0 d_2 - \nu \ggamma_R d_2, 0, 0, \nu (d_1 - d_0)).
    \end{split}
    \end{equation*}
    The map $\Psi$ is well-defined by the normalization conditions of $\eta^{001}, \eta^{011}$ and $\phi^{001}$ and the fact that $(\ker d^1_0)^*$ is isomorphic to $\Ann(s_0\huaV_0)$. 
    We now check that $\Psi^{-1}$ is well-defined. 
    By Theorem \ref{thm:VS2Dual}, the element $\Psi^{-1}(\xi_1, \xi_2)$ is an element of $\ker \widecheck{d}^1_0$ if and only if for any $X \in \huaV_2$,
    \begin{equation}\label{eq:vs-2dual-norm-cplx1}
    \begin{cases}
        (\xi_1 \ggamma_R d_2 - \xi_2 \ggamma_R \partial)(X - s_1d_2X) = (\xi_2 \ggamma_R d_1 - \xi_2 \ggamma_R d_2)(X - s_0d_0X), \\
        (\xi_2 \ggamma_R d_1 - \xi_2 \ggamma_R d_2)(X - s_0d_1X + s_0d_2 X) = 0,
    \end{cases}
    \end{equation}
    while the element $\Psi^{-1}(\theta, \nu)$ is an element of $\ker \widecheck{p}^2_2$ if and only if for any $X \in \huaV_2$, 
    \begin{equation}\label{eq:vs-2dual-norm-cplx2}
    \begin{cases}
        (\theta d_0 d_2 - \nu \ggamma_R d_2)(\ggamma_0 X) = 0, \\
        (\theta d_0 d_2 - \nu \ggamma_R d_2)(\ggamma_2 X) = 
        (\nu (d_1 - d_0))(\ggamma_1 X), \\
        (\theta d_0 d_2 - \nu \ggamma_R d_2)(X - s_0d_0X) + (\nu (d_1 - d_0))(s_1d_2 X) = 0, \\
        (\theta d_0 d_2 - \nu \ggamma_R d_2)(X - s_1d_2X) = 0.
    \end{cases}
    \end{equation}
    Both \eqref{eq:vs-2dual-norm-cplx1} and \eqref{eq:vs-2dual-norm-cplx2} hold by a straightforward computation with the simplicial identities. 
    The fact that $\Psi^{-1}$ is actually the inverse of $\Psi$ is also a straightforward computation, where we remark that $\Psi^{-1}\Psi = id$ follows from the defining equations of $\ker \widecheck{d}_0$ and $\ker \widecheck{p}^2_2$, which are \eqref{eq:vs-2dual-norm-cplx1} and \eqref{eq:vs-2dual-norm-cplx2}, respectively, with $\Psi(\eta^{001}, \eta^{011})$ replacing $(\xi_1, \xi_2)$ in the first one and $\Psi(\phi^{012}, 0,0, \phi^{001})$ replacing $(\theta, \nu)$ in the second one.
    Finally, the differentials in \eqref{eq:vs-2dual-norm-cplx} are obtained by computing 
    \begin{equation*}
        \Psi_1 \widecheck{d}_2 \Psi^{-1}_2 (\theta, \nu) = (\nu, \nu + \theta d_1), \qquad 
        - \Psi_0 \widecheck{d}_1 \Psi^{-1}_1 (\xi_1, \xi_2) = \xi_2d_2 - \xi_1 d_2,
    \end{equation*}
    which follow once again by applying simplicial identities and evaluating the first on elements in $\ker d^1_0$ and the second on elements in $\ker p^2_2$.
\end{proof}

\begin{remark}\label{rem:vs-2dual-lambda-r}
    As expected, $N(\huaV^{2*})$ is not isomorphic to $N(\huaV)^*[-2]$, as it has one more copy of $(\ker d_0^1)^*$ in degree 1 and one in degree 2. By Definition \ref{def:nshifted-dual-chain-cplx}, the chain complex $N(\huaV)^*[-2]$ is 
    \begin{equation*}
        \huaV_0^* \xrightarrow{d_1^t} (\ker d^1_0)^* \xrightarrow{d_2^t} (\ker p^2_2)^*,
    \end{equation*}
    and the right induced map $\lambda^r_{\langle \cdot, \cdot \rangle}$ of the IM-pairing associated to the 2-dual pairing $\langle \cdot, \cdot \rangle$ is 
    \begin{equation}\label{eq:lambda-ndual-pairing-right}
    \begin{split}
        \lambda_{\langle \cdot, \cdot \rangle}^r : &N(\huaV^{2*}) \to N(\huaV)^*[-2]\\
        &(\lambda_{\langle \cdot, \cdot \rangle}^r)_2(\theta, \nu) = \theta \in \huaV_0^*,\\
        &(\lambda_{\langle \cdot, \cdot \rangle}^r)_1(\xi, \eta) = \xi - \eta \in (\ker d^1_0)^*,\\
        &(\lambda_{\langle \cdot, \cdot \rangle}^r)_0(\epsilon) = \epsilon \in (\ker p^2_2)^*.
    \end{split}
    \end{equation}
\end{remark}

\section{The model \texorpdfstring{$DK(N(\huaV)^*[-n])$}{DK(N(V)*[-n])}}\label{sec:appendix-DK-N-nshifted-dual}

As anticipated in Remark \ref{rem:DK-n-shift-dual-appendix}, we now discuss the other possible model of $n$-dual imported from chain complexes via the Dold-Kan correspondence. For simplicity, throughout this section, whenever we talk about the $n$-dual of $\huaV$ we consider $\huaV$ to be a $\VS$ $n$-groupoid.

First of all, $DK(N(\huaV)^*[-n])$ is the $n$-dual of the $\VS$ $n$-groupoid $\huaV$ for $n=0, 1$ because in these cases the deformation retract in Theorem \ref{thm:EZ-thm-hom} is actually an isomorphism, as in Example \ref{ex:VS0dual} and Remark \ref{rem:VS1Dual-DK}. 
In these cases the $n$-dual pairing is also nondegenerate on the nose and it coincides with a trivial extension of the canonical evaluation pairing of the $n$-th level $\huaV_n$ with its dual vector space to a simplicial pairing. 

In general cases (for $n\ge 2$), we would like to retain this property of the $n$-dual pairing being a trivial extension of the evaluation pairing of $\huaV_n^*$ to a canonical simplicial pairing $\huaV^{n*} \otimes \huaV \to B^n{\R}$, while requiring it to be nondegenerate only up to homotopy. 
This is precisely what we get, as we saw in Section \ref{sec:ndual-overview}. 
On the other hand, using the maps in the proof of Theorem \ref{thm:EZ-thm-hom}, one could produce a pairing on $DK(N(\huaV)^*[-n])$ by taking the evaluation pairing on chains, which is a chain map $\lambda: N(\huaV)^*[-n] \otimes N(\huaV) \to \R[-n]$, and defining $DK(AW^*\lambda)$. 
There are two reasons why this is inconvenient:

Firstly, although at level $n$, $DK(N(\huaV)^*[-n])_n \cong \huaV_n^*$, the pairing $DK(AW^*\lambda)$ will not coincide with the evaluation pairing $ev: \huaV_n^* \otimes \huaV_n \to \R$ outside the cases $n=0,1$. This is because, as we show in Proposition \ref{prop:DK2shiftedDual-counterex}, for $n\ge 2$, $ev$ does not extend trivially to a simplicial pairing $DK(N(\huaV)^*[-n]) \otimes \huaV \to B^n(\R)$ even in simple cases. 

The second reason has to do with the version of this theory for higher vector bundles. 
As shown in \cite{HoyoTrentinaglia2024}, the Dold-Kan correspondence for higher vector bundles involves the extra data given by a representation up to homotopy (RUTH) of the base simplicial manifold on the normalized complex of the higher vector bundle in question. 
While the construction of the chain complex $N$ is entirely analogous to that for simplicial vector spaces, a choice of a cleavage is needed to obtain the RUTH. 
In addition, a RUTH of a Lie $n$-groupoid $\huaG$ for $n \ge 2$ cannot be canonically dualized to a RUTH of $\huaG$, but only to a RUTH of its opposite $\huaG^{op}$ obtained by simplicial front-to-back symmetry as in Remark \ref{rem:front-to-back}.
This is due to the absence of an inversion map for arrows in Lie $n$-groupoids with $n\ge 2$ which usually gives the canonical isomorphism between a Lie groupoid and its opposite.
The RUTH is also essential to define the functor $DK$, and applying this functor to a RUTH of $\huaG^{op}$ gives a higher vector bundle on $\huaG^{op}$. 
Therefore, in this context, even assuming both $DK(N(\huaV)^*[-n])$ and $DK(AW^*\lambda)$ can be somehow defined over $\huaG$ and not its opposite, these must necessarily depend on a choice of cleavage. 
On the contrary, our approach in defining the $n$-dual and its pairing directly on the simplicial side can be extended to avoid both the choice of a cleavage and the problem with dualizing representations, thus resulting in a canonical description of the $n$-dual and its pairing. 

We now set out to show that the canonical evaluation pairing of $\huaV_2^*$ cannot be trivially extended to a simplicial pairing $DK(N(\huaV)^*[-2]) \otimes \huaV \to B^2\R$ in general.
To simplify the description of $DK(N(\huaV)^*[-2])$, we claim that given a $\VS$ 2-groupoid $\huaV$, $DK(N(\huaV)^*[-2])$ is isomorphic to the $\VS$ 2-groupoid $K^*(\huaV)$, which we define as 
\[\begin{tikzcd}
	{\huaV_2^*} & {(\ker d_0^2)^*} & {(\ker p_2^2)^*}
	\arrow[shift left=1, from=1-2, to=1-3]
	\arrow[shift right=1, from=1-2, to=1-3]
	\arrow[from=1-1, to=1-2]
	\arrow[shift left=2, from=1-1, to=1-2]
	\arrow[shift right=2, from=1-1, to=1-2]
\end{tikzcd}\]
with the following simplicial maps:
Given $\phi \in \huaV_2^*$, $\nu \in (\ker d_0^2)^*$, $\theta \in (\ker p_2^2)^*$, $v \in \huaV_2$, $h \in (\ker d_0^2)$, $k \in (\ker p_2^2)$, take
\begin{equation*}
\begin{split}
    &(\widecheck{d}_0^2\phi) (h) = \phi (h), 
    \qquad (\widecheck{d}_1^2\phi) (h) = \phi (h - s_0d_1 h),\\
    &(\widecheck{d}_2^2\phi) (h) = \phi (h - s_0d_1 h - s_1d_2 h + 1d_1d_1 h),\\
    &(\widecheck{d}_0^1\nu) (k) = \nu (k),
    \qquad (\widecheck{d}_1^1\nu) (k) = \nu (k - s_1d_2 k),
\end{split}
\end{equation*}
and
\begin{equation*}
\begin{array}{cc}
    (\widecheck{s}_0^1\nu)(v) = \nu (v - s_0d_0v),
    &(\widecheck{s}_1^1\nu)(v) = \nu (v - s_0d_0v - s_1d_1 v + 1d_0d_0 v),\\
    (\widecheck{s}_0^0\theta) (h) = \theta (h - s_1d_1 h). &
\end{array}
\end{equation*}
A straightforward computation shows these obey the simplicial identities \eqref{eq:simp-id}. This 2-truncated simplicial vector space is then equipped with the canonical $\VS$ 2-groupoid structure by using the multiplications in Example \ref{ex:MooreFillers2Gpd}. 

\begin{proposition}
    For any $\VS$ 2-groupoid $\huaV$ we have $K^*(\huaV) \cong DK(N(\huaV)^*[-2])$.
\end{proposition}

\begin{proof}
We prove this by showing that the normalized complex of $K^*(\huaV)$ is isomorphic to the 2-shifted dual of the normalized complex of $\huaV$, $N(\huaV)^*[-2]$.    

First of all, $N(\huaV)$ is the chain complex 
\begin{equation*}
    \ker p_2^2 \xrightarrow{d_2} \ker d_0^1 \xrightarrow{-d_1} \huaV_0
\end{equation*}
which is concentrated in degrees 0 to 2. Then, its $2$-shifted dual $N(\huaV)^*[-2]$ is the chain complex 
\begin{equation*}
    \huaV_0^* \xrightarrow{d_1^*} (\ker d_0^1)^* \xrightarrow{d_2^*} (\ker p_2^2)^*,
\end{equation*}
which is also concentrated in degrees 0 to 2.

The normalized chain complex of $K^*(\huaV)$ is
\begin{equation*}
    \ker\widecheck{p}^2_2
    \xrightarrow{\widecheck{d}_2} \ker\widecheck{d}_0^1 
    \xrightarrow{-\widecheck{d}_1}(\ker p_2^2)^*.
\end{equation*}

By definition of $\widecheck{d}_0^1$, $\ker\widecheck{d}_0^1 = \Ann(\ker p^2_2)\subset (\ker d_0^2)^*$. 
$\ker d_0^2$ splits as $\ker d_0^2 \cong \ker p^2_2 \oplus s_1(\ker d^1_0)$ via the splitting $h =  (h - s_1d_2h) + s_1d_2h$ for any $h \in \ker d_0^2$. So $\Ann(\ker p^2_2) \cong (\ker d^1_0)^*$, with the explicit isomorphism given by
\begin{equation*}
(\ker d^2_0)^* \supset \Ann(\ker p^2_2) \newrightleftarrows{s_1^*}{d_1^*} (\ker d^1_0)^*
\end{equation*}
This isomorphism intertwines the differentials $-\widecheck{d}^1_1: \Ann(\ker p^2_2) \to (\ker p^2_2)^*$ and $d_2^* : (\ker d_0^1)^* \to (\ker p^2_2)^*$, because
\begin{equation*}
    -(\widecheck{d}^1_1\nu)(k) = - \nu ( k - s_1d_2 k )
    = - \nu( s_1d_1(k - s_1d_2 k)) = \nu(s_1d_2 k), \quad \forall k \in \ker p^2_2.
\end{equation*} 

For $\ker\widecheck{p}^2_2 = \ker \widecheck{d}^2_0 \cap \ker \widecheck{d}^2_1$, we have $\ker\widecheck{d}_0^2 = \Ann(\ker d_0^2)\subseteq \huaV_2^*$. Moreover, for any $\phi \in \ker\widecheck{d}_0^2 \cap \ker \widecheck{d}_1^2$, we have
\begin{equation*}
    (\widecheck{d}_1^2\phi)(h) =  \phi(h) - \phi(s_0d_1 h) = -  \phi(s_0d_1 h) = 0, \quad \forall h \in \ker d_0^2.
\end{equation*}
Then, by surjectivity of $d_1: \ker d_0^2 \to \ker d_0^1$, $\phi (s_0 x) = 0$ for any $x \in \ker d_0^1$. So $\phi \in \Ann(\ker d_0^2) \cap \Ann(s_0(\ker d_0^1)) \subseteq \huaV_2^*$. The other inclusion is obvious, so $\ker\widecheck{p}^2_2 = \Ann(\ker d_0^2) \cap \Ann(s_0(\ker d_0^1))$.  Furthermore, by the splitting $\huaV_2 \cong \ker d_2^0 \oplus s_0 (\ker d_0^1) \oplus 1(\huaV_0)$ given for any $v\in \huaV_2$ by 
\begin{equation*}
    v = (v - s_0d_0 v)  + (s_0d_0 v - 1d_0d_0 v) + 1d_0d_0v, \quad \forall v\in \huaV_2,
\end{equation*}
$\Ann(\ker d_0^2) \cap \Ann(s_0(\ker d_0^1)) \cong \huaV_0^*$. 
The explicit isomorphism is 
\begin{equation*}
\huaV_2^* \supset \Ann(\ker d_0^2) \cap \Ann(s_0(\ker d_0^1)) 
\newrightleftarrows{1^*}{(d_0d_0)^*} \huaV_0^*,
\end{equation*}
which intertwines the differentials $\widecheck{d}^2_2$ and $d_1^*$. In fact,
\begin{equation*}
    (\widecheck{d}^2_2 (\theta d_0d_0))(s_1x) 
    = \theta (d_0d_0 (s_1x - s_0x - s_1x + 1d_1x))
    = \theta (d_1x), \quad \forall x \in \ker d^1_0.
\end{equation*}

By combining all of this we obtain the isomorphism of chain complexes
\[\begin{tikzcd}[ampersand replacement=\&,sep=scriptsize]
	{\ker \widecheck{p}^2_2} \& {\ker \widecheck{d}^1_0} \& {(\ker p^2_2)^*} \\
	{\huaV_0^*} \& {(\ker d_0^1)^*} \& {(\ker p^2_2)^*.}
	\arrow["{\widecheck{d}^2_2}", from=1-1, to=1-2]
	\arrow["{- \widecheck{d}^1_1}", from=1-2, to=1-3]
	\arrow["{1^*}"', shift right, from=1-1, to=2-1]
	\arrow["{s_1^*}"', shift right, from=1-2, to=2-2]
	\arrow["{d_0^*d_0^*}"', shift right, from=2-1, to=1-1]
	\arrow[Rightarrow, no head, from=1-3, to=2-3]
	\arrow["{d_2^*}"', from=2-2, to=2-3]
	\arrow["{d_1^*}"', from=2-1, to=2-2]
	\arrow["{d_1^*}"', shift right, from=2-2, to=1-2]
\end{tikzcd}\]
\end{proof}

\begin{proposition}\label{prop:DK2shiftedDual-counterex}
    Let $\huaV$ be the pair groupoid of $\R$, $\R^2 \rightrightarrows \R$. Then the 2-shifted pairing $\langle \cdot,\cdot \rangle: K^*(\huaV)_2 \otimes \huaV_2 \to \R$ on 2-simplices is not multiplicative.
\end{proposition}

\begin{proof}
Take arbitrary $(\phi,\chi,\psi)\in \Lambda^3_1(K^*(\huaV))$ (each of these is an element in $\huaV_2^*$), and $(W,Y,Z)\in \Lambda^3_1(\huaV)$. The pairing is multiplicative if and only if 
\begin{equation*}
    \langle \phi\square\chi\psi, W\square YZ\rangle
    = \langle \phi, W \rangle + \langle \chi, Y \rangle - \langle \psi, Z \rangle.
\end{equation*}
Consider the elements $\widecheck{s}_0\phi = (\phi,\widecheck{s}_0\widecheck{d}_1\phi, \widecheck{s}_0\widecheck{d}_2\phi) \in \Lambda^3_1(K^*(\huaV))$, with $\widecheck{d}_1\widecheck{s}_0\phi = \phi$, and $s_1 W = (s_0d_0W,W,s_1d_2W)$ with $d_1s_1W = W$. Then, if the pairing is multiplicative, 
\begin{equation*}
    \begin{aligned}
    \langle \phi, W \rangle &= \langle \phi\square(\widecheck{s}_0\widecheck{d}_1\phi)(\widecheck{s}_0\widecheck{d}_2\phi), (s_0d_0W)\square W(s_1d_2W) \rangle\\
    &= \langle \phi, s_0d_0W \rangle + \langle \phi - \phi s_0d_1, W\rangle - \langle \phi -\phi s_0d_1 - \phi s_1d_2 + 1d_1d_1, s_1d_2 W\rangle\\
    \end{aligned}
\end{equation*}
which is true if and only if
\begin{equation*}
    \begin{aligned}
    0  &= \langle \phi, s_0d_0W - s_0d_1 W - s_1d_2 W + s_0\stkout{d_1s_1}d_2 W + s_1\stkout{d_2s_1}d_2 W - 1d_1\stkout{d_1s_1}d_2 W\rangle\\
    &= \langle \phi, s_0\partial W - 1 d_1d_2 W\rangle,
    \end{aligned}
\end{equation*}
for $\partial= d_0 -d_1 +d_2$, the boundary map. In particular, for any $\theta \in \huaV_0^*$, we have $\theta d_0d_0 \in \huaV_2^*$, so that the multiplicativity of the pairing implies
\begin{equation*}
    \theta(d_0\partial W) = \theta(d_0d_2W) = \theta(d_1d_2W). 
\end{equation*}
Since $d_2^2$ is surjective (for any $\huaV$), this implies that $\partial^*\theta = 0$ for any $\theta \in \huaV_0^*$. 
For $\huaV = \R^2 \rightrightarrows \R$, when evaluating this on any $(x,0) \in \R^2$, we get that
\begin{equation*}
    \theta (x) = 0,
\end{equation*}
for any $\theta \in \R^*$ and $x \in \R$, yielding a contradiction.
\end{proof}

\chapter{Duals of \texorpdfstring{$\VB$}{VB} \texorpdfstring{$n$}{n}-groupoids}


In this chapter we define the $n$-dual of any higher vector bundle over a Lie $n$-groupoid by means of a universal property. 
As we will see, this is obtained at each level by a categorical limit in the category of vector bundles, which means its existence is not guaranteed, but depends on the existence and regularity of solutions for a certain set of equations. 

We begin by defining $n$-shifted pairings for higher vector bundles over the same base in Section \ref{sec:VB-pairings}. As for vector spaces, these are pairings with values in the Eilenberg-MacLane space $B^n\R$, this time seen as a trivial simplicial vector bundle over the base Lie $n$-groupoid. We then define the associated IM-pairings and the homological nondegeneracy condition. This condition previously appeared for shifted symplectic structures on a Lie groupoid. We recall their definition here, as well. 

Section \ref{sec:n-dual-construction} deals with the general construction of $n$-duals of higher vector bundles over a Lie $n$-groupoid starting with the universal property. We then develop some notation to implement the equations in Section \ref{sec:ndual-overview} pointwise over the base Lie $n$-groupoid and show that these identify each level of the $n$-dual as a limit in $\VB$ in Theorem \ref{thm:n-dual-as-limit}. At this point we show that the $\VB$ $n$-dual is consistent with the $\VS$ $n$-dual, and that it can be constructed by solving only a finite number of equations.
In the next subsection we discuss some properties of the $n$-dual pairing that carry over from the vector space version, i.e. the fact that it is homologically $n$-shifted nondegenerate in Theorem \ref{thm:VBndual-pairing-hndg} and its relation to the induced map of the IM-pairing of any other $n$-shifted pairing in Theorem \ref{thm:VBnShiftedPairingsInducedMapsTriangle}.

After that, we review the terminology necessary to discuss solution spaces of the $n$-dual equations by generalizing the notion of cores of a $\VB$ groupoid. We close this section by solving the 1-dual equations pointwise for a $\VB$ 1-groupoid and showing that this recovers the $\VB$ 1-dual from Section \ref{sec:VB1-groupoids-cores-duals}.
We then recall the analog of Theorem \ref{thm:hom-nondeg-homotopy-equiv} in this case: by using the 1-dual, one can reformulate the homological nondegeneracy of a 1-shifted pairing as the fact that its induced maps are weak equivalences. 

In Section \ref{sec:2-dual-VB-gpd} we define a model of the 2-dual of a $\VB$ 2-groupoid from the top down and show this is a $\VB$ 2-groupoid by checking the finite data conditions in Example \ref{ex:finite-data-2gpds}, and that it satisfies the 2-dual equations. We then apply Theorem \ref{thm:we-of-VB2gpd-is-qi} to Theorem \ref{thm:VBnShiftedPairingsInducedMapsTriangle} and obtain the $\VB$ 2-dual version of Theorem \ref{thm:hom-nondeg-homotopy-equiv}: a 2-shifted pairing is homologically nondegenerate if and only if its induced maps are weak equivalences. 
Applied to the 2-dual pairing, this means that 2-duality is reflexive up to homotopy for $\VB$ 2-groupoids.
In the final section of this chapter, we compute 2-duals of vector bundles and $\VB$ 1-groupoids explicitly.

\section{Pairings of higher vector bundles}\label{sec:VB-pairings}

To begin with, we establish the dualization unit that will appear in the universal property. This is the target space for simplicial pairings of higher vector bundles. 

\begin{definition}
    Let $\huaG$ be a Lie $n$-groupoid. Then $B^n\R^\huaG$ is the simplicial vector bundle over $\huaG$ defined levelwise by the trivial vector bundle
    \begin{equation*}
        B^n\R^\huaG_m := B^n\R_m \times \huaG_m \to \huaG_m,
    \end{equation*}
    for any $m\ge 0$, with the obvious face and degeneracy maps. 
\end{definition}

\begin{lemma}
    Let $\huaG$ be a Lie $n$-groupoid, then $B^n\R^\huaG$ is a $\VB$ $n$-groupoid. 
\end{lemma}

\begin{proof}
    The $n$-truncation of $B^n\R^\huaG$ is given by
    \begin{equation*}
        \R \times \huaG_n \rightfourarrows \dots \rightfourarrows 0 \times \huaG_2 \rightthreearrows 0 \times \huaG_1 \rightrightarrows 0 \times \huaG_0,
    \end{equation*}
    with the obvious simplicial maps. The $\VB$ $n$-groupoid structure is given by the multiplication maps
    \begin{equation*}
        \widetilde{m}_k(a_i|_{g_i})_{i \neq k} = \left.\left(\sum_{i=0, i\neq k}^n (-1)^{i-k+1} a_i \right)\right|_{m_k(g_i)_{i\neq k}},
    \end{equation*}
    for any $(n+1,k)$-horn $(a_i)_{i\neq k} \in \Lambda^{n+1}_k(B^n\R^\huaG)$ over $(g_i)_{i\neq k} \in \Lambda^{n+1}_k(\huaG)$.
\end{proof}

The tensor product of higher vector bundles over the same base has an obvious definition analogous to that for simplicial vector spaces.

\begin{definition}
Let $\huaV \to \huaG$ and $\huaW \to \huaG$ be higher vector bundles over the same simplicial manifold $\huaG$. Then their \textbf{tensor product} $\huaV\otimes \huaW \to \huaG$ is the higher vector bundle defined levelwise by 
\begin{equation*}
    (\huaV \otimes \huaW)_n := \huaV_n \otimes \huaW_n,
\end{equation*}
with simplicial maps 
\begin{equation*}
    \widetilde{d}^{\huaV\otimes \huaW}_i := \widetilde{d}^{\huaV}_i \otimes \widetilde{d}^{\huaW}_i,
    \qquad 
    \widetilde{s}^{\huaV \otimes \huaW}_i:= \widetilde{s}^{\huaV}_i \otimes \widetilde{s}^{\huaW}_i.
\end{equation*}
\end{definition}

\begin{remark}\label{rem:order-of-tensor-vb}
As seen from the fact that simplicial bundles over the point are simplicial vector spaces, by Remark \ref{rem:order-vb-gpd-norm-cplx} and Theorem \ref{thm:order-of-tensors} applied pointwise, the order of the tensor product of a $\VB$ $n$-groupoid of order $n$ with a $\VB$ $m$-groupoid of order $m$ is at least $n+m$. Hence the tensor product is not internal to the category of $\VB$ $n$-groupoids, but it is internal to $\VB^\infty_\huaG$. 
\end{remark}

\begin{definition}\label{def:VB-n-shifted-pairing}
Let $\huaV \to \huaG$ and $\huaW \to \huaG$ be two higher vector bundles over a Lie $n$-groupoid $\huaG$. 
Consider $\R$ as a vector bundle over a point. 
We call a bundle map $\alpha_n: \huaV_n \otimes \huaW_n \to \R \times \huaG_n$ an \textbf{$n$-shifted pairing} of $\huaV$ with $\huaW$.
Additionally, if 
\begin{equation}\label{eq:VBnShiftedPairingMultiplicative}
    \alpha_n \widetilde{d}_0 - \alpha_n \widetilde{d}_1 + \alpha_n \widetilde{d}_2 - \dots + (-1)^n \alpha_n \widetilde{d}_n = 0, 
\end{equation}
with $\widetilde{d}_i = \widetilde{d}_i^\huaV \otimes \widetilde{d}_i^\huaW$, we say $\alpha_n$ is \textbf{multiplicative}.
We also say $\alpha_n$ is \textbf{normalized} if
\begin{equation}\label{eq:VBnShiftedPairingNormalized}
    \alpha_n \widetilde{s}_i = 0, \quad \forall 0 \le i < n,
\end{equation}
with $\widetilde{s}_i = \widetilde{s}_i^{\huaV} \otimes \widetilde{s}_i^{\huaW}$.

We call a simplicial vector bundle map $\alpha: \huaV \otimes \huaW \to B^n\R \times \huaG$ an \textbf{$n$-shifted simplicial pairing}. 
\end{definition}

\begin{remark}
The fact that $\alpha_n$ is a bundle map from $\huaV \otimes \huaW$ to $\R$ is equivalent to saying that it is fiberwise linear: for any $g\in \huaG_n$, $v, v' \in \huaV|_g$ and $w, w' \in \huaW|_g$, 
\begin{equation*}
    \alpha_n(v, w) + \alpha_n(v', w') = \alpha_n(v +|_g v', w +|_g w'),
\end{equation*}
where $+|_g$ denotes addition in the fiber over $g$ of $\huaV$ or $\huaW$ depending on context.
\end{remark}

\begin{remark}
In the same way as for simplicial vector spaces, $\alpha_n$ is an $n$-cochain in the simplicial cohomology of $\huaV \otimes \huaW$. In this language, $\alpha_n \in C^n(\huaV \otimes \huaW)$ is multiplicative if and only if it is closed with respect to the simplicial differential 
\begin{equation}\label{eq:simp-differential}
\delta_n := \sum_{i=0}^n (-1)^i (\widetilde{d}_i^n)^*. 
\end{equation}
In other words, \eqref{eq:nShiftedPairingMultiplicative} can be written simply as
\begin{equation*}
    \delta \alpha_n = 0.
\end{equation*}

If $\huaV$ and $\huaW$ are $\VB$ $n$-groupoids, another equivalent way to rewrite multiplicativity is 
\begin{equation}\label{eq:VBmul-pair}
    \alpha_n(\widetilde{m}^{\huaV}_k ((X_i)_{0\le i\neq k \le n}), \widetilde{m}^{\huaW}_k((Y_i)_{0\le i\neq k \le n}) ) = \sum_{i=0, \, i\neq k }^{n} (-1)^{i-k+1} \alpha(X_i, Y_i),
\end{equation} for any $k\in \{0, \dots, n\}$ and any $(n+1, k)$-horns $(X_i)_{0\le i\neq k \le n}$ and $(Y_i)_{0\le i\neq k \le n}$ over the same basepoint in $\Lambda^{n+1}_k(\huaG)$.
\end{remark}

\begin{lemma}
    Let $\huaG$ be a $\VB$ $n$-groupoid and $\huaV, \huaW$ two higher vector bundles over it. 
    There is a one-to-one correspondence between $n$-shifted simplicial pairings $\alpha: \huaV \otimes \huaW \to B^n\R$ and $n$-shifted pairings $\alpha_n: \huaV_n \otimes \huaW_n \to \R$ that are multiplicative and normalized.
    In this case we identify $\alpha$ and $\alpha_n$.
\end{lemma}

\begin{proof}
The proof is analogous to the one for vector spaces, by using Lemma \ref{lem:finite-data-simp-map}, because $B^n\R \times \huaG$ is a $\VB$ $n$-groupoid. 
Simplicial bundle maps $\huaV \otimes \huaW \to B^n\R \times \huaG$ are in fact determined by their $n$-th level, which is required to satisfy precisely the multiplicativity and normalization conditions \eqref{eq:VBnShiftedPairingMultiplicative} and \eqref{eq:VBnShiftedPairingNormalized}. 
In the context of Lemma \ref{lem:finite-data-simp-map}, the normalization conditions ensure that the $n$-truncated map with $\alpha_n$ at level $n$ and $0$ at all lower levels is indeed an $n$-truncated simplicial map, while the multiplicativity condition \eqref{eq:VBnShiftedPairingMultiplicative} is equivalent to the condition \eqref{eq:simp-map-multiplicative} of the same name. 
\end{proof}

The definition of IM-pairing in Definition \ref{def:IMPairing} extends pointwise over $\huaG_0$ to define IM-pairings for chain complexes of vector bundles over the same base. 
After applying the normalized complex construction in Definition \ref{def:VB-norm-cplx} to the tensor product of two vector bundles, the Eilenberg-Zilber map is also defined pointwise as in Definition \ref{def:EilenbergZilberMap}.
Thus, each $n$-shifted pairing admits an associated IM-pairing, as per Proposition \ref{prop:associated-IM-pairing}. We now recall the definition of this IM-Pairing.

\begin{definition}\label{def:VBIMPairingAssociatedDef}
    Let $\huaV \to \huaG$ and $\huaW \to \huaG$ be two higher vector bundles over a $\VB$ $n$-groupoid $\huaG$, and $\alpha: \huaV \otimes \huaW \to B^n\R\times \huaG$ an $n$-shifted simplicial pairing between them. 
    Then the \textbf{associated IM-pairing} $\lambda_\alpha: (N(\huaV)\otimes N(\huaW))_n \to \R$ is defined at any point $p\in \huaG_0$ by
    \begin{equation}\label{eq:VBIMPairingAssociatedDef}
        \lambda_\alpha(v,w)_p = \sum_{(\mu,\nu)\in\Shuf(i,n-i)} \sign(\mu,\nu) \alpha(\widetilde{s}_{\nu_{n-i}}\dots \widetilde{s}_{\nu_1} v, \widetilde{s}_{\mu_{i}}\dots \widetilde{s}_{\mu_1} w),
    \end{equation}
    for any $v \in N(\huaV)_i|_p$ and any $w \in N(\huaW)_{n-i}|_p$.
\end{definition}

Once again, 
this induces a left and right map as defined pointwise at each point of $\huaG_0$ as in \eqref{eq:IMPairingDef-induced-maps}. We can now define homologically nondegenerate pairings between $\VB$ $n$-groupoids. 

\begin{definition}\label{def:VB-hndg-pairing}
    An $n$-shifted simplicial pairing $\alpha: \huaV \otimes \huaW \to B^n\R \times \huaG$ between two higher vector bundles is \textbf{homologically $n$-shifted nondegenerate} if its IM-pairing $\lambda_\alpha$ descends to a nondegenerate pairing 
    \begin{equation*}
        \lambda_\alpha|_p: H_i(N(\huaV)|_p)\otimes H_{n-i}(N(\huaW)|_p) \to \R,
    \end{equation*}
    for all $i \in \Z$ and at all points $p\in \huaG_0$.
    That is, it induces an isomorphism between the homologies of the normalized complexes, up to a degree shift of $n$.
    Equivalently, if $\huaV$ and $\huaW$ are $\VB$ $n$-groupoids, $\alpha$ is homologically $n$-shifted nondegenerate if either of the induced maps $\lambda_\alpha^r|_p$ or $\lambda_\alpha^l|_p$ is a quasi-isomorphism for each $p \in \huaG_0$.\footnote{By the same argument as in Definition \ref{def:homological-non-deg}, if one is a quasi-isomorphism at each point, the other one is as well.}
\end{definition}

An important example of homologically nondegenerate $n$-shifted pairings is given by $n$-shifted symplectic structures on Lie $m$-groupoids. These were introduced in \cite{Lesdiablerets} and studied in \cite{CuecaZhu2023}. 
We give a slightly different but equivalent definition that follows naturally from this discussion. See \cite{CuecaZhu2023} for more details. 

\begin{definition}\label{def:shifted-symplectic}
Let $\huaG$ be a Lie $m$-groupoid. An \textbf{$n$-shifted presymplectic form} on $\huaG$ is a collection of differential forms 
\begin{equation*}
    (\omega_2, \omega_3, \dots, \omega_{n+2}), \text{ with } \omega_i \in \Omega^i(\huaG_{n-i+2}),
\end{equation*}
such that 
\begin{enumerate}
\item $\omega_2 \in \Omega^2(\huaG_n)$ induces a simplicial pairing $T\huaG \otimes T\huaG \to B^2\R^\huaG$, i.e. it satisfies the multiplicativity and normalization conditions \eqref{eq:VBnShiftedPairingMultiplicative} and \eqref{eq:VBnShiftedPairingNormalized}.\footnote{Since $\omega_2$ is a 2-form, the pairing it induces is in particular antisymmetric (at level $n$).}
\item For any $i$, $d\omega_i = (-1)^{n-i+1} \delta \omega_{i+1}$, where $\delta$ is the simplicial differential \eqref{eq:simp-differential}, and $d$ is the de Rham differential.
\item For any $i$, $\omega_i$ is normalized in the sense that $s_j^*\omega_i = 0 \in \Omega^i(\huaG_{n-i+1})$. 
\end{enumerate}

If, in addition, the pairing induced by $\omega_2$ is homologically nondegenerate, then we say $\omega$ is an \textbf{$n$-shifted symplectic form}.
\end{definition}

In practice, in this thesis, we will only consider shifted symplectic forms for which the only non-zero component is $\omega_2$. In particular, for such structures, $d\omega_2 = 0$.

\begin{remark}
In (2) it is intended that $\delta\omega_2 = 0$, which is also contained in the first condition (multiplicativity of the pairing), and that $d\omega_{n+2}=0$.
In (3) normalization of $\omega_2$ is equivalent to the fact it is normalized as a pairing, so this condition is already included in (1). 

Condition (2) and the multiplicativity of condition (1) are equivalent to the fact that $\omega$ is closed with respect to the total differential $\delta + (-1)^p d$ ($p$ is the simplicial degree) in the Bott-Shulman-Stasheff double complex $\Omega^\bullet(\huaG_\bullet)$ introduced in \cite{BottShulmanStasheff1976}. The collection of forms $(\omega_i)_i$ is an element of degree $n+2$ in the total complex. It is required to have vanishing components in bidegrees $(0, n+2)$ and $(1, n+1)$ for reasons explained in \cite[p.27]{Lesdiablerets}, which originally come from \cite{PantevToenVaquieVezzosi2013}. 
\end{remark}

\section{General construction of \texorpdfstring{$\VB$}{VB} \texorpdfstring{$n$}{n}-duals}\label{sec:n-dual-construction}

We first give a definition of the $n$-dual by a universal property. 

\begin{definition}\label{def:n-dual-univ-prop}
Let $\huaV$ be a higher vector bundle over a Lie $n$-groupoid $\huaG$. 
If it exists, the \textbf{$n$-dual} of $\huaV$ is the representing object of the functor 
\begin{equation*}
\begin{split}
    \VB^{\infty}_\huaG(\huaV\otimes \_, B^n\R^\huaG): &(\VB^{\infty}_\huaG)^{op} \to \Vect\\
    &\huaW \to \VB^{\infty}_\huaG(\huaV \otimes \huaW, B^n\R^\huaG),
\end{split}
\end{equation*}
from the category of higher vector bundles over $\huaG$ with simplicial bundle maps over the identity of $\huaG$ to the category of vector spaces. 
If it exists, we denote the $n$-dual of $\huaV$ by $\huaV^{n*}$.
\end{definition}

In the following, we show that the $n$-dual of a higher vector bundle can be described as a limit, by identifying it as the kernel of a bundle map, i.e. a simplicial subbundle of a certain ambient bundle satisfying a series of equations which are a pointwise version of the ones in Proposition \ref{prop:vs-n-dual-equations}.
Because this construction is a limit, a priori this still does not solve the question of existence, which has three obstructions: 
\begin{enumerate}
    \item The equations need to admit solutions at each fiber and each level.
    \item The solution spaces must assemble into levelwise constant rank vector bundles and define a simplicial vector bundle. See Example \ref{rem:VBnGpd-0-dual} for an example of a situation in which this fails to hold.
    \item The result satisfies the appropriate Kan conditions, so that it is indeed a higher vector bundle. 
\end{enumerate}
By existence of dual vector bundles and dual $\VB$ groupoids, we know that these requirements are satisfied in the case of the $n$-dual of a $\VB$ $n$-groupoid for $n=0, 1$. We will then show these are also satisfied in the case of $n=2$.
In any case, because of the defining property being universal, if any such solution exists, it is unique up to isomorphism, and we can already use its definition to study some of its properties.

\subsection{The pointwise version of the \texorpdfstring{$n$}{n}-dual equations}

Let $\huaV \to \huaG$ be a higher vector bundle. Any $m$-simplex $g\in \huaG_m$ of the base can be seen as a simplicial map $g: \Delta^m \to \huaG$. 
Recall from Lemma \ref{lem:Delta-n-simplex-class} that any $u \in \Delta^m_l$ can be written as $u = s_Id_J E_{m}$ for some multi-indices $I,J$ such that $|I|-|J| = l - m + 1$ and $E_{m}$ the unique nondegenerate simplex in $\Delta^m_{m}$. 
With this, evaluating $g$ as a map on a $u \in \Delta^m_l$ gives 
\begin{equation*}
    g(u) = s_Id_Jg, 
\end{equation*}
which follows immediately by the fact that $g$ is a simplicial map and $g(E_m) = g$. 
Moreover, for any $u \in \Delta^{m-1}_l$, such that $u = s_Id_J E_{m-1}$ and for any $0 \le i \le m$, 
\begin{equation*}
    \begin{split}
        (d_ig)(u) 
        &= (d_i g)(s_I d_J E_{m-1}) 
        = s_Id_Jd_i g 
        = g(s_Id_Jd_iE_m) \\
        &= g(s_Id_J\delta_iE_{m-1}) 
        = g(\delta_i s_Id_J E_{m-1}) 
        = g(\delta_i u).
    \end{split}
\end{equation*}
Analogously, $(s_ig)(r) = g(\sigma_i r)$. 

We now recall the simplicial subset of $\huaG$ generated by $g$ from Definition \ref{def:simp-subset-generated}. This is the image $\langle\langle g \rangle \rangle := g(\Delta^m) \subseteq \huaG$. For example a point in $\huaG_0$ generates the simplicial subset
\begin{equation*}
	\langle\langle p\rangle\rangle = \dots \rightfourarrows \{1p\}\rightthreearrows \{1p\} \rightrightarrows \{p\},
\end{equation*}
which is isomorphic to the zero simplex $\Delta^0$. In the same way, for a nondegenerate $g\in G_1$, we have
\begin{equation*}
\begin{aligned}
\langle\langle g\rangle\rangle &= \dots \rightfourarrows \{1d_0g,s_1g,s_0g, 1d_1g\}\rightthreearrows \{s_0d_0g,g,s_0d_1g\} \rightrightarrows \{d_0g, d_1g\},
\end{aligned}
\end{equation*}
which is isomorphic to the 1-simplex $\Delta^1$. It is important to assume nondegeneracy for the isomorphism $\langle \langle g \rangle \rangle \cong \Delta^m$, because if $g \in \huaG_m$ is degenerate, the map $g:\Delta^m \to \huaG$ is not injective and $\langle\langle g \rangle\rangle$ is not isomorphic to $\Delta^m$, but to a lower-dimensional simplex. For example, if $m=1$ and $g=1p$, then $\langle\langle 1 p \rangle \rangle = \langle \langle p \rangle \rangle \cong \Delta^0$.

\begin{definition}
    Let $\huaV \to \huaG$ be a higher vector bundle. The \textbf{ simplicial fiber} of $\huaV$ at $g \in \huaG_m$ is the simplicial vector space $\textstyle\int g^*\huaV$. This is given at each level $l$ by 
    \begin{equation}\label{eq:VBSimplicialFiberDef}
        \textstyle\int g^*\huaV_l = \oplus_{u \in \Delta^m_l} (\huaV_l|_{g(u)})^u,
    \end{equation}
    the direct sum of the fibers of the pullback 
    \[\begin{tikzcd}[ampersand replacement=\&,cramped,column sep=small,row sep=scriptsize]
	{g^*\huaV} \& \huaV \\
	{\Delta^m} \& \huaG,
	\arrow[from=1-1, to=1-2]
	\arrow[from=1-1, to=2-1]
	\arrow["\lrcorner"{anchor=center, pos=0.125}, draw=none, from=1-1, to=2-2]
	\arrow[from=1-2, to=2-2]
	\arrow["g"', from=2-1, to=2-2]
    \end{tikzcd}\]
    which is a simplicial vector bundle over the discrete simplicial manifold $\Delta^m$. 
    The simplicial maps of $\textstyle\int g^*\huaV$ are induced by those of $\huaV$ and are given on each component $v^u \in (\huaV_l|_{g(u)})^u$ by 
    \begin{equation}\label{eq:VBnSimplicialFiber-def-maps}
        \begin{array}{cc}
            (\widetilde{d}_i v^u)|_{d_i g(u)} = (\widetilde{d}_i v)^{d_i u} |_{d_i g(u)}
            &(\widetilde{s}_i v^u)|_{s_i g(u)} = (\widetilde{s}_i v)^{s_i u} |_{s_i g(u)}.
        \end{array}
    \end{equation}
\end{definition}

The combinatorial structure of this simplicial fiber construction is manifestly similar to that of the copowering of simplicial objects in \ref{def:copowering}. For example, if $g\in \huaG_1$, we can picture $\textstyle\int g^*\huaV$ in a diagram such as 
\[
\begin{adjustbox}{width=\textwidth}
\begin{tikzcd}[ampersand replacement=\&,cramped,column sep=small,row sep=tiny]
	\& \dots \&\& \dots \&\& \dots \\
	{(\huaV_2|_{1d_1g})^{000}} \&\& {(\huaV_2|_{s_0g})^{001}} \&\& {(\huaV_2|_{s_1g})^{011}} \&\& {(\huaV_2|_{1d_0g})^{111}} \\
	\& {(\huaV_1|_{1d_1g})^{00}} \&\& {(\huaV_1|_{g})^{01}} \&\& {(\huaV_1|_{1d_0g})^{11}} \\
	\&\& {(\huaV_0|_{d_1g})^0} \&\& {(\huaV_0|_{d_0g})^1}
	\arrow[from=2-1, to=3-2]
	\arrow["{\widetilde{d}_0, \widetilde{d}_1, \widetilde{d}_2}", shift left=2, from=2-1, to=3-2]
	\arrow[shift right=2, from=2-1, to=3-2]
	\arrow["{\widetilde{d}_2}"', from=2-3, to=3-2]
	\arrow[shift right, from=2-3, to=3-4]
	\arrow["{\widetilde{d}_0, \widetilde{d}_1}", shift left, from=2-3, to=3-4]
	\arrow["{\widetilde{d}_1, \widetilde{d}_2}"', shift right, from=2-5, to=3-4]
	\arrow[shift left, from=2-5, to=3-4]
	\arrow["{\widetilde{d}_0}", from=2-5, to=3-6]
	\arrow[from=2-7, to=3-6]
	\arrow[shift left=2, from=2-7, to=3-6]
	\arrow["{\widetilde{d}_0, \widetilde{d}_1, \widetilde{d}_2}"', shift right=2, from=2-7, to=3-6]
	\arrow[shift left=2, curve={height=-6pt}, from=3-2, to=2-1]
	\arrow["{\widetilde{s}_0,\widetilde{s}_1}", shift left=4, curve={height=-6pt}, from=3-2, to=2-1]
	\arrow["{\widetilde{d}_0, \widetilde{d}_1}", shift left, from=3-2, to=4-3]
	\arrow[shift right, from=3-2, to=4-3]
	\arrow["{\widetilde{s}_0}", shift left=2, curve={height=-6pt}, from=3-4, to=2-3]
	\arrow["{\widetilde{s}_1}"', shift right=2, curve={height=6pt}, from=3-4, to=2-5]
	\arrow["{\widetilde{d}_1}", from=3-4, to=4-3]
	\arrow["{\widetilde{d}_0}"', from=3-4, to=4-5]
	\arrow[shift right=2, curve={height=6pt}, from=3-6, to=2-7]
	\arrow["{\widetilde{s}_0,\widetilde{s}_1}"', shift right=4, curve={height=6pt}, from=3-6, to=2-7]
	\arrow[shift left, from=3-6, to=4-5]
	\arrow["{\widetilde{d}_0, \widetilde{d}_1}"', shift right, from=3-6, to=4-5]
	\arrow["{\widetilde{s}_0}", shift left, curve={height=-6pt}, from=4-3, to=3-2]
	\arrow["{\widetilde{s}_0}"', shift right, curve={height=6pt}, from=4-5, to=3-6]
\end{tikzcd}
\end{adjustbox}
\]
As in that case, we will prefer the coproduct notation for elements of $\textstyle\int g^*\huaV$. This means a generic element is seen as a ``formal'' sum of homogeneous elements $v \in \textstyle\int g^*\huaV$ whose basepoint and index are well-defined. In practice, we will exclusively consider homogeneous elements (since everything is linear), and commit the slight abuse of notation of writing $v = v^u = v|_{g(u)}$, whenever $v \in \huaV|_{g(u)}^u$ and it is convenient to highlight the basepoint $g(u)$ or the index $u$.
As a consequence of the initial discussion, it is important to distinguish between $\textstyle\int g^*\huaV$ and the sum of fibers of the restriction of $\huaV$ to $\langle \langle g \rangle \rangle$. This is the reason why we keep the upper index $u$ in the notation: to make it clear that, in the case where $g$ is degenerate, we may still count multiple copies of the same fiber if they correspond to different indices $u$. 

\begin{remark}\label{rem:simp-fiber-at-units}
This construction actually recovers the copowering of a simplicial vector space with a standard simplex, when $g$ is a unit $1p \in \huaG_m$ for $p \in \huaG_0$. In fact, for any $p \in \huaG_0$, and the embedding $1p \in \huaG_m$ of it as a unit at any level $m\ge 0$, all the generated simplicial subsets are isomorphic to a point: $\langle\langle 1p \rangle\rangle \cong \langle\langle p \rangle\rangle \cong \Delta^0$. Therefore
\begin{equation*}
(\textstyle\int(1p)^*\huaV)_l \cong \oplus_{u \in \Delta^m_l} (\huaV_l|_{1p})^u \cong ((1^*\huaV)|_p \otimes \Delta^m)_l,
\end{equation*} 
the level $l$ of the copowering of the fiber of $1^*\huaV$ over $p$ with $\Delta^m$, as defined in Section \ref{sec:review-monoidal-DK} and Definition \ref{def:copowering}.
By comparing the simplicial maps in \eqref{eq:VBnSimplicialFiber-def-maps} with those in Definition \ref{def:copowering}, we get that $\textstyle \int (1p)^* \huaV \cong (1^*\huaV)|_p \otimes \Delta^m$, as simplicial vector spaces.
\end{remark} 

For any simplicial vector space $\huaW$ we can define the hom vector space
\begin{equation}
    \SVect(\textstyle\int g^*\huaV, \huaW)
\end{equation}
of simplicial $\R$-linear maps from $\textstyle\int g^*\huaV$ to $\huaW$ as defined in Section \ref{sec:Svect-int-hom}. 
Roughly speaking, the fact that $\textstyle\int g^*\huaV$ has an analogous simplicial structure to the copowering with $\Delta^m$, implies that $\SVect(\textstyle\int g^*\huaV, \huaW)$ has a similar form to the $m$-th level of a simplicial mapping space. 
In particular, for $\huaW = B^n\R$, $\SVect(\textstyle\int g^*\huaV, B^n\R) = (\textstyle\int g^*\huaV)^{n*}_0$ encodes a version with basepoints of the equations in \ref{sec:ndual-overview} that compute the $m$-th level of the $n$-dual of a vector space. 

\begin{lemma}\label{lem:n-dual-equations}
    Let $\huaV \to \huaG$ be a higher vector bundle over a Lie $n$-groupoid $\huaG$. Then, for any $g \in \huaG_m$, any $\phi \in \SVect(\textstyle\int g^*\huaV, B^n\R)$ consists of $\binom{m+n}{n}$ components $\phi^u \in (\huaV_n|_{g(u)}^u)^*$ indexed over $n$-simplices $u \in \Delta^m_n$, defined by 
    \begin{equation}\label{eq:VBnDualm-simplexComponents}
        \phi^u(v) = \phi_n(v|_{g(u)}^u), \quad \forall u \in \Delta^m_n
    \end{equation}
    for any homogeneous $v \in \huaV_n|_{g(u)}^u \subseteq \textstyle\int g^*\huaV$. These components satisfy the \textit{normalization} equations
    \begin{equation}\label{eq:VBnDualNormalizationMostGeneral}
        \phi^{s_it}(\widetilde{s}_i(v)) = \phi_n((\widetilde{s}_iv)^{s_it}) = 0,
    \end{equation}
    for all $0 \le i \le n-1$, $t \in \Delta^m_{n-1}$, and $v \in \huaV_{n-1}|_{g(t)}$,
    and the \textit{multiplicativity} equations 
    \begin{equation}\label{eq:VBnDualMultiplicativityMostGeneral}
        \sum_{i=0}^{n} (-1)^{i} \phi^{d_i r}(\widetilde{d}_i v) = 0,
    \end{equation}
    for all $r \in \Delta^m_{n+1}$ and all $v \in \huaV_{n+1}|_{g(r)}$. 
\end{lemma}

\begin{proof}
This follows entirely from Proposition \ref{prop:vs-n-dual-equations}. 
\end{proof}

By definition, these coincide with equations \eqref{eq:NormalizationMostGeneral} and \eqref{eq:MultiplicativityMostGeneral-short}, with the only difference being that whenever $g \neq 1p$ for some $p \in \huaG_0$, we have to remember the basepoint of each component and each input. Moreover, when $\huaV$ is an $n$-groupoid, we cannot rely on canonical multiplication maps to solve the equations as we did in the proof of Theorem \ref{thm:VS2Dual}. As a consequence, most of the discussion in Section \ref{sec:ndual-overview} carries over to this case. 
In particular, there is a choice of independent components to make, which will determine all the other ones. 
The solution spaces of \eqref{eq:VBnDualNormalizationMostGeneral} are elements of certain annihilator bundles which we study more in detail in \ref{sec:VB-DegAnn-Cores}. 
We also make the same distinction between \textit{face components} and \textit{interior components} as in Definition \ref{def:face-interior-comps}. Because they are also indexed in the same way, Remark \ref{rem:number-of-int-comp} also applies without modification. 
In fact, for any $g \in \huaG_n$, there is a unique interior component at level $n$, which is not normalized and which corresponds to the simplex $E_n \in \Delta^{n}_n$. For $\phi\in \SVect(\textstyle\int g^*\huaV, B^n\R)$, this is $\phi^{E_n} \in \huaV_n^*|_{g}$, as $g(E_n) = g$.

The simplicial vector bundle structure of the solutions space of these equations is induced by seeing the latter as a subbundle of an ambient simplicial vector bundle defined by the property that each of its $m$-simplices at $g\in\huaG_m$ consists of the components \eqref{eq:VBnDualm-simplexComponents}. 
This is one of the spaces appearing in the limit description of the $n$-dual.

\begin{definition}
    Let $\huaV \to \huaG$ be a higher vector bundle. The \textbf{ambient $n$-dual bundle} $\mathfrak{A}_n\huaV$ of $\huaV$ is the simplicial vector bundle over $\huaG$ given at each level by
    \begin{equation}
        \mathfrak{A}_n\huaV_m := \bigoplus_{\substack{u=s_Id_JE_m\\ u\in \Delta^m_n}} \left( (s_Id_J)^*\huaV_n^*\right)^u
    \end{equation} 
    with simplicial maps 
    \begin{equation}\label{eq:VBndual-ambient-face-maps}
        \begin{split}
            \widecheck{d}_i: &\mathfrak{A}_n\huaV_m|_g \to \mathfrak{A}_n\huaV_{m-1}|_{d_ig}\\
            &(\widecheck{d}_i f)^u|_{(d_i g)(u)} = f^{\delta_i u}|_{g(\delta_i u)}, \quad \forall u \in \Delta^{m-1}_n,
        \end{split}
    \end{equation}
    \begin{equation}\label{eq:VBndual-ambient-deg-maps}
        \begin{split}
            \widecheck{s}_i: &\mathfrak{A}_n\huaV_m|_g \to \mathfrak{A}_n\huaV_{m+1}|_{s_ig}\\
            &(\widecheck{s}_i f)^u|_{(s_i g)(u)} = f^{\sigma_i u}|_{g(\sigma_i u)}, \quad \forall u \in \Delta^{m+1}_n.
        \end{split}
    \end{equation}
\end{definition}

\begin{example}\label{ex:VB1dual-ambient}
    The first 3 levels of the 1-dual ambient bundle $\mathfrak{A}_1\huaV$ are 
    \begin{equation*}
        \begin{split}
            \mathfrak{A}_1\huaV_0 &= (1^*\huaV_1^*)^{00} \to \huaG_0, \\
            \mathfrak{A}_1\huaV_1 &= (\huaV_1^*)^{01}\oplus \bigoplus_{i=0}^1((1d_i)^*\huaV_1^*)^{ii} \to \huaG_1, \\
            \mathfrak{A}_1\huaV_2 &= \bigoplus_{(\mu,\nu)\in \Shuf(2,1)}\left((d_{\nu_1}^*\huaV_1^*)^{\mu_1\mu_2} \oplus ((1d_{\mu_1}d_{\mu_2})^*\huaV_1^*)^{\nu_1\nu_1}\right) \to \huaG_2.
        \end{split}
    \end{equation*}
\end{example}

\begin{example}\label{ex:VB2dual-ambient}
    The first 4 levels of the 2-dual ambient bundle $\mathfrak{A}_2\huaV$ are 
    \begin{equation*}
        \begin{split}
            \mathfrak{A}_2\huaV_0 &= (1^*\huaV_2^*)^{000} \to \huaG_0, \\
            \mathfrak{A}_2\huaV_1 &= (s_0^*\huaV_2^*)^{001} \oplus (s_1^*\huaV_2^*)^{011} \oplus \bigoplus_{i=0}^1((1d_i)^*\huaV_2^*)^{iii}  \to \huaG_1, \\
            \mathfrak{A}_2\huaV_2 &= (\huaV_2^*)^{012} \oplus \bigoplus_{(\mu,\nu)\in \Shuf(2,1)}
            \big(((s_0d_{\nu_1})^*\huaV_2^*)^{\mu_1\mu_1\mu_2} 
            \oplus ((s_1d_{\nu_1})^*\huaV_2^*)^{\mu_1\mu_2\mu_2}\big) \\
            &\qquad\qquad\qquad\qquad\qquad
            \oplus \bigoplus_{(\mu,\nu)\in \Shuf(1,2)} ((1d_{\nu_1}d_{\nu_2})^*\huaV_2^*)^{\mu_1\mu_1\mu_1} \to \huaG_2.\\
            \mathfrak{A}_2\huaV_3 &= \bigoplus_{(\mu,\nu)\in\Shuf(3,1)}(d_{\nu_1}^*\huaV_2^*)^{\mu_1\mu_2\mu_3}\\
            &\qquad
            \oplus\bigoplus_{(\mu,\nu)\in \Shuf(2,2)}
            \big(((s_0d_{\nu_1}d_{\nu_2})^*\huaV_2^*)^{\mu_1\mu_1\mu_2} 
            \oplus ((s_l1d_{\nu_1}d_{\nu_2})^*\huaV_2^*)^{\mu_1\mu_2\mu_2}\big) \\
            &\qquad
            \oplus \bigoplus_{(\mu,\nu)\in \Shuf(1,3)} ((1d_{\nu_1}d_{\nu_2}d_{\nu_3})^*\huaV_2^*)^{\mu_1\mu_1\mu_1} \to \huaG_3.\\
        \end{split}
    \end{equation*}
\end{example}

\begin{remark}
    Because each level of $\mathfrak{A}_n\huaV$ is constructed by using only the natural operations of pullbacks and Whitney sums, it has a canonical trivialization induced by that of $\huaV^*_n$. 
    Each $m$-simplex in $\mathfrak{A}_n\huaV_m$ consists of the components \eqref{eq:VBnDualm-simplexComponents}, which are indexed in the same way as the components \eqref{eq:n-dual-element-in-comp} of an $m$-simplex of the $n$-dual of a simplicial vector space as discussed in Section \ref{sec:ndual-overview}. 
    Therefore, because for any $g \in \huaG_m$ the components of an $m$-simplex of $\SVect(g^*\huaV, B^n\R)$ have the same components with the same basepoints as those of an $m$-simplex in $\mathfrak{A}_n\huaV|_g$, we have that 
    \begin{equation*}
        \coprod_{g\in\huaG_m} \SVect(\textstyle\int g^*\huaV, B^n\R) \subseteq \mathfrak{A}_n\huaV_m.
    \end{equation*}
    A priori this might not be a vector subbundle, as for example, it might not be constant rank. 

    The simplicial maps in \eqref{eq:VBndual-ambient-face-maps} and \eqref{eq:VBndual-ambient-deg-maps} have the same form as those of the simplicial mapping space as defined in \eqref{eq:mapping-set-maps}. 
    Hence they clearly satisfy the simplicial identities.
    Moreover they are well-defined bundle maps, since they are just projections to certain components of an element. 
\end{remark}

The next proposition expresses the $\VB$ $n$-dual as the solution space of the $n$-dual equations of Lemma \ref{lem:n-dual-equations} inside the ambient bundle. 

\begin{theorem}\label{thm:n-dual-as-limit}
Let $\huaV$ be a higher vector bundle over the Lie $n$-groupoid $\huaG$. If its $n$-dual exists, it is isomorphic to the simplicial subbundle of $\mathfrak{A}_n\huaV$ that satisfies the equations from Lemma \ref{lem:n-dual-equations}. That is, if $\huaV^{n*}$ exists, then for any $g\in\huaG_m$, 
\begin{equation*}
    V^{n*}_m|_g = \Simp\Vect (\textstyle\int g^*\huaV, B^n\R) \subseteq \mathfrak{A}_n\huaV_m.
\end{equation*} 
In other words, the $m$-th level of the $n$-dual of $\huaV$ is the kernel of the bundle map 
\begin{equation*}
    \begin{split}
        \mathfrak{d}_m: \quad 
        \mathfrak{A}_n\huaV_m &\longrightarrow \bigoplus_{t \in \Delta^m_{n-1}}(ev_t^*\huaV_{n-1}^*)^{n} \oplus \bigoplus_{r \in \Delta^m_{n+1}} ev_r^*\huaV_{n+1}^*\\
        \phi|_g &\mapsto \left((\widetilde{s}_i^*\phi^{s_it})_{t \in \Delta^m_{n-1}, 0\le i \le n-1}, (\sum_{i=0}^n (-1)^i\widetilde{d}_i^*\phi^{d_i r})_{r \in \Delta^m_{n+1}}\right),
    \end{split}
\end{equation*}
where $ev_t: \huaG_m \to \huaG_{m-1}$ and $ev_r: \huaG_m \to \huaG_{m+1}$ are the evaluation maps $g \mapsto g(t)$ and $g \mapsto g(r)$.
\end{theorem}

\begin{proof}
    It is enough to show that if there exist a subbundle $\huaV^{n*} \subseteq \mathfrak{A}_n\huaV$ satisfying the equations from Lemma \ref{lem:n-dual-equations}, then $\huaV^{n*}$ satisfies the universal property of the $n$-dual from Definition \ref{def:n-dual-univ-prop}, i.e. 
    \begin{equation*}
        \VB^{\infty}_\huaG(\huaV^{n*}, \huaW) \cong \VB^{\infty}_\huaG(\huaV \otimes \huaW, B^n\R^\huaG), 
    \end{equation*}
    for any higher vector bundle $\huaW$ over $\huaG$. 
    This defines the $n$-dual up to isomorphism. 

    We define isomorphisms analogous to those for the adjunction \ref{prop:tensor-hom-svect} at level 0. We rewrite them in this case for the sake of concreteness:
    For arbitrary $\alpha \in \VB^{\infty}_\huaG(\huaV \otimes \huaW, B^n\R_{\huaG})$, and any simplicial map $\beta \in \VB^{\infty}_\huaG(\huaW, \huaV^{n*})$, we define 
    \begin{equation}\label{eq:VBnDual-rho-def}
    \begin{split}
        \rho: \VB^{\infty}_\huaG(\huaV \otimes \huaW, B^n\R_{\huaG}) &\longrightarrow \VB^{\infty}_\huaG(\huaW, \huaV^{n*})\\
        (\rho(\alpha)_l (w|_g))^u (v|_{g(u)}) &= \alpha|_{g(u)}(v, \widetilde{s}_I \widetilde{d}_Jw),\\
        \forall l \ge 0, \quad g \in \huaG_l,&\quad w \in \huaW_l|_g,\quad u = s_Id_JE_l \in \Delta^{l}_n,\quad v \in \huaV_n|_{g(u)},
    \end{split}
    \end{equation}
    and its (candidate) inverse
    \begin{equation}\label{eq:VBnDual-tau-def}
    \begin{split}
        \tau: \VB^{\infty}_\huaG(\huaW, \huaV^{n*}) &\longrightarrow \VB^{\infty}_\huaG(\huaV \otimes \huaW, B^n\R_{\huaG})\\
        \tau(\beta)_n|_g(v, w) &= (\beta(w|_g))^{E_n}(v|_g)\\
        \forall g \in \huaG_n,\quad v \in \huaV_n|_g,&\quad w \in \huaW_n|_g.
    \end{split}
    \end{equation}
    In defining $\tau$ we used the fact that an $n$-shifted simplicial pairing is determined entirely by its $n$-th level, as usual. 

    We begin by showing that $\rho(\alpha)$ is well-defined. 
    That is, all of its components satisfy the normalization and multiplicativity conditions, \eqref{eq:VBnDualNormalizationMostGeneral} and \eqref{eq:VBnDualMultiplicativityMostGeneral}, necessary to make $\rho(\alpha)^u(w)$ a well-defined element of $\huaV^{n*}_l|_{g(u)}$, for any $w \in \huaW_l|_{g(u)}$ and any $g \in \huaG_l$. 

    Starting with the normalization conditions \eqref{eq:VBnDualNormalizationMostGeneral}, we have that for any $v \in \huaV_{n-1}|_{g(t)}$, $t = s_Jd_I E_l \in \Delta^{l}_{n-1}$ and $w \in \huaW_{l}|_{g}$,
    \begin{equation*}
        (\rho(\alpha)(w))^{s_i t}(\widetilde{s}_i v) = \alpha(\widetilde{s}_i v, \widetilde{s}_i \widetilde{s}_J \widetilde{d}_I w) = 0, 
    \end{equation*}
    for any $0 \le i \le n-1$, by normalization of $\alpha$ as a pairing \eqref{eq:VBnShiftedPairingNormalized}.

    For the multiplicativity conditions \eqref{eq:VBnDualMultiplicativityMostGeneral}, we have that, for any $v \in \huaV_{n+1}|_g(r)$, any $r = s_Jd_I E_l \in \Delta^{l}_{n+1}$, and any $w\in \huaW_{l}|_g$, 
    \begin{equation*}
        \sum_{i=0}^n (\rho(\alpha)(w))^{d_i r} (\widetilde{d}_i v) 
        = \sum_{i=0}^n \alpha(\widetilde{d}_i v, \widetilde{d}_i \widetilde{s}_I \widetilde{d}_J w) = 0,
    \end{equation*}
    by the multiplicativity of $\alpha$ as a pairing \eqref{eq:VBnShiftedPairingMultiplicative}. 

    Furthermore, $\rho(\alpha)$ is a simplicial morphism. Recall that the simplicial maps of $\huaV^{n*}$ are defined by \eqref{eq:VBndual-ambient-face-maps} and \eqref{eq:VBndual-ambient-deg-maps}. Then, for any $w \in \huaW_l$ and $u = s_Jd_I E_{l-1}\in \Delta^{l-1}_n$, 
    \begin{equation*}
        \begin{split}
            (d_j(\rho(\alpha)(w)))^u &= (\rho(\alpha)(w))^{\delta^j s_J d_I E_{l-1}} 
            = (\rho(\alpha)(w))^{s_J d_I d_j E_{l}} \\
            &= \alpha(\_, s_J d_I d_j w)
            = (\rho(\alpha)(d_j w))^{u},
        \end{split}
    \end{equation*}
    for any $0 \le j \le l$, by \eqref{eq:Delta-simp-cosimp-commutativity} and \eqref{eq:d-delta-En-identity}. The same holds for the degeneracy maps. 

    We now show that $\rho$ is the inverse of $\tau$.
    Let $\beta: \huaW \to \huaV^{n*}$ be a simplicial map. For any $l\ge 0$, $v \in \huaV_n|_g$, $u = s_Id_JE_l \in \Delta^{l}_n$ and $w \in \huaW_l|_{g(u)}$, we have
    \begin{equation*}
    \begin{split}
        (\rho(\tau(\beta))_l(w))^{u}(v) &= \tau(\beta)(v, s_Id_J w) = (\beta(s_Id_Jw))^{E_n}(v) \\
        &= (\widecheck{s}_I\widecheck{d}_J\beta(w))^{E_n}(w) \overset{\eqref{eq:VBndual-ambient-face-maps}}{=} (\beta(w))^{\delta_{\widebar{I}} \sigma_{\widebar{J}} E_n}(v)\\
        &\overset{\eqref{eq:simp-cosimp-description-u}}{=} (\beta(w))^{s_Id_J E_l}(v) = (\beta(w))^{u}(v),
    \end{split}
    \end{equation*}
    The other direction is even more straightforward: for any simplicial pairing $\alpha: \huaV \otimes \huaW \to B^n\R^{\huaG}$,
    \begin{equation*}
        \tau(\rho(\alpha))_n(v, w) 
        = (\rho(\alpha)_n(v))^{E_n}(w) 
        = \alpha(v, w),
    \end{equation*}
    for any $v \in \huaV_n|_g$ and $w \in \huaW_n|_g$
\end{proof}

\begin{remark}\label{rem:VB-n-pairing-ind-maps}
By symmetry of the tensor product of higher vector bundles, if the $n$-duals of both $\huaV\to \huaG$ and $\huaW \to \huaG$ exist, then
\begin{equation*}
    \VB^\infty_\huaG(\huaV \otimes \huaW, B^n\R) \cong 
    \VB^\infty_\huaG(\huaW, \huaV^{n*}) \cong
    \VB^\infty_\huaG(\huaV, \huaW^{n*}).
\end{equation*}
This is an ``adjoint-like'' property of the $n$-dual, but it is not technically an adjunction, since we did not define the $n$-dual as an actual functor. 
We plan to extend this to a well-defined adjunction in future work. 

In any case, this property allows us to define two induced maps from any $n$-shifted pairing between higher vector bundles whose $n$-duals exist. 
These are defined in the same way as in \eqref{eq:nShiftedPairingIndMapsDef} for $n$-shifted pairings of simplicial vector spaces, but using the map $\rho$ defined in \eqref{eq:VBnDual-rho-def} and its symmetric version.
To fix the notation, we define, for any $\alpha:\huaV \otimes \huaW \to B^n\R^\huaG$, the components of the \textbf{left} and \textbf{right induced maps} for $u = s_I d_J E_m \in \Delta^m_n$ at $g \in \huaG_m$ as 
\begin{equation}\label{eq:VBnShiftedPairingIndMapsDef}
    (\alpha^l(w|_g))^q (v|_{g(u)}) = \alpha(v, \widetilde{s}^{\huaW}_I \widetilde{d}^{\huaW}_J w), 
    \qquad (\alpha^r(v'|_g))^q (w'|_{g(u)}) = \alpha(\widetilde{s}^{\huaV}_I \widetilde{d}^{\huaV}_J v', w'),
\end{equation}
for any $w \in \huaW_m|_{g}$, $v \in \huaV_n|_{g(u)}$, and any $v' \in \huaV_m|_{g(u)}$, $w' \in \huaW_n|_g$.
\end{remark}

We now discuss two properties of the $n$-dual of a higher vector bundle that can be directly inferred from the $n$-dual equations from Lemma \ref{lem:n-dual-equations}. The first one is that this construction is consistent with the $n$-dual of a simplicial vector space, in the sense that the pullback to the units of the $n$-dual is pointwise the $n$-dual of the pullback to the units. 
This allows us to make some remarks about existence of 0-duals of higher vector bundles over a 0-groupoid. 

The second one allows us to show that given a candidate $\VB$ $n$-groupoid whose finite data satisfies the $n$-dual equations from Lemma \ref{lem:n-dual-equations} for levels $0 \le m \le n+1$, this satisfies all higher equations as well, and therefore it is the $n$-dual. In the discussion of the $\VS$ $n$-duals we obtained a similar result on the finiteness of equations in the opposite way: we first showed that the $\VS$ $n$-dual is always a $\VS$ $n$-groupoid and then because of that it was determined by its finite data. 
In the $\VB$ case we need to show the proposed finite data solving the equations is in fact a $\VB$ $n$-groupoid by other means, as we do for $n=2$ in Theorem \ref{thm:VB2dual-is-2-gpd}.

\begin{lemma}\label{lem:VBndual-at-units}
    Let $\huaV \to \huaG$ be a higher vector bundle over a Lie $n$-groupoid $\huaG$. If its $n$-dual $\huaV^{n*}$ exists, then
    \begin{equation*}
        (1^*\huaV^{n*})|_{p} \cong ((1^*\huaV)|_p)^{n*},
    \end{equation*}
    as simplicial vector spaces. 
    That is, the fiber of the pullback of the $n$-dual of $\huaV$ by the total unit map at any point $p \in \huaG_0$ is the $\VS$ $n$-dual of the fiber of the pullback of $\huaV$ by the unit map, at that same point. 
\end{lemma}

\begin{proof}
    Recall from Remark \ref{rem:simp-fiber-at-units} that the simplicial fiber of $\huaV$ at a unit $1p \in \huaG_m$ is $\textstyle\int (1p)^*\huaV \cong (1^*\huaV)|_p \otimes \Delta^m$. 
    Therefore, for any $m \ge 0$,
    \begin{equation*}
        (1^*\huaV^{n*})_m|_{p} \cong \SVect((1^*\huaV)|_p \otimes \Delta^m, B^n\R) = ((1^*\huaV)|_p)^{n*}_m. 
    \end{equation*}
    Moreover, the simplicial maps induced on $1^*\huaV^{n*}$ by $1^*\mathfrak{A}_n\huaV$ are the same as the simplicial maps of the $\VS$ $n$-dual (compare \eqref{eq:mapping-set-maps} with \eqref{eq:VBndual-ambient-face-maps} and \eqref{eq:VBndual-ambient-deg-maps}). So these levelwise isomorphisms upgrade to an isomorphism of simplicial vector spaces. 
\end{proof}

This result gives an immediate example of a situation where the $n$-dual of a higher vector bundle does not exist. 

\begin{remark}\label{rem:VBnGpd-0-dual}
Because of this lemma, if $\huaV$ is a higher vector bundle over a 0-groupoid $M$, its $n$-dual can be computed pointwise at each $p \in M$ by using the construction of the $\VS$ $n$-dual in Section \ref{sec:computations}. 
In general, without any assumptions on the order of $\huaV$, the pointwise solution spaces might not assemble to vector bundles over $M$ at each level. 

Consider for example the 0-dual of a $\VB$ $n$-groupoid over a 0-groupoid $M$. By Example \ref{ex:VS0dual}, we have
\begin{equation*}
    \huaV^{0*}_l|_p \cong H_0(\huaV|_p)^*,
\end{equation*}
for any $l \ge 0$ and any $p \in M$. 
This homology is generally not constant rank over $M$ and as such, $\huaV^{0*}$ is not always a well-defined $\VB$ 0-groupoid. 
On the other hand, if $\huaV$ is a $\VB$ 0-groupoid over $M$ (i.e. a classical vector bundle), then this construction coincides with the usual dual of vector bundles and it is well-defined.
\end{remark}

We conclude this section by showing that the $n$-dual can be computed by solving a finite number of equations.

\begin{lemma}\label{lem:n-dual-finite-equations}
    Let $\huaV \to \huaG$ be a higher vector bundle over a Lie $n$-groupoid.
    Then, the ambient $n$-dual bundle $\mathfrak{A}_n\huaV$ is $n$-coskeletal.
    Moreover, let $\huaV'\to \huaG$ be a simplicial subbundle of $\mathfrak{A}_n\huaV$. The following are equivalent
    \begin{enumerate}[label=(\roman*)]
        \item $\huaV'$ satisfies the $n$-dual equations of Lemma \ref{lem:n-dual-equations} at all levels $m\ge 0$.
        \item $\huaV'$ satisfies the $n$-dual equations at levels $0 \le m \le n+1$ and it is $(n+1)$-coskeletal.
    \end{enumerate}
    In particular any simplicial subbundle of $\mathfrak{A}_n\huaV$ that is a $\VB$ $n$-groupoid whose finite data satisfies the $n$-dual equations for $\huaV$ is the $n$-dual $\huaV^{n*}$.
\end{lemma}

\begin{proof}
    The simplicial subbundle $\huaV' \subseteq \mathfrak{A}_n\huaV$ satisfies the $n$-dual equations if and only if at each level $m\ge 0$, 
    \begin{equation*}
        \huaV'_m = \ker \mathfrak{d}_m, 
    \end{equation*}
    by Theorem \ref{thm:n-dual-as-limit}.
    Consider the jointly conservative collection of points $\mathsf{P}$ from Example \ref{ex:points-for-mfds}. 
    Then each point functor preserves finite limits, so 
    \begin{equation*}
        \mathsf{p}\huaV'_m = \ker \mathsf{p}\mathfrak{d}_m.
    \end{equation*}
    Which means that $\mathsf{p}\huaV'_m$ satisfies the same equations as $\huaV'_m$ but for $\mathsf{p}\huaV$ inside $\mathfrak{A}_n\mathsf{p}\huaV$. 
    We now show that $\mathfrak{A}_n\mathsf{p}\huaV$ is $n$-coskeletal and $\ker \mathsf{p}\mathfrak{d}_m$ is $(n+1)$-coskeletal for any $\mathsf{p}$, so that by Lemma \ref{lem:stalkwise-cosk-is-cosk}, $\mathfrak{A}_n\huaV$ is $n$-coskeletal and $\huaV'$ is $(n+1)$-coskeletal, both in $\Cat$. 

    We show this by proving that fiberwise, 
    \begin{equation*}
    \begin{split}
        \mathfrak{A}_n(\mathsf{p}\huaV)_m|_g = \partial_m(\mathfrak{A}_n\mathsf{p}\huaV)|_{\theta_m g}, 
        \quad \text{for any } m \ge n+1,\\
        \mathsf{p}\huaV'_m|_{g} = \partial_m(\mathsf{p}\huaV'_m)|_{\theta_m g},
        \quad \text{for any } m \ge n+2,
    \end{split}
    \end{equation*}
    for any $g \in \mathsf{p}\huaG_m$, and $\theta_m=(d_0, \dots, d_m)$ the boundary projection of $\huaG$.

    The fundamental observation is that, 
     for any simplex $u \in \Delta^m_l$ with $m \ge l+1$, we have that
    \begin{equation*}
        u = \delta_j u', \text{ with } u' \in \Delta^{m-1}_l.
    \end{equation*}
    This can be seen by using \eqref{eq:simp-cosimp-description-u} and writing $u = s_J d_I E_l = \delta_{\widebar{J}} \sigma_{\widebar{I}} E_m$, for some multi-indices $I, J$ with opposites $\widebar{I}, \widebar{J}$, respectively and such that $|J|-|I| = m - l$. 
    If $m \ge l+1$, then $|J| \ge 1$, so any $u \in \Delta^m_l$ is $\delta_j$ of some simplex $u' \in \Delta^{m-1}_{l}$. 

    Therefore, any element $\phi \in (\mathsf{p}\huaV'_m)|_g$ for $g \in \huaG_m$ with $m \ge n+1$ is entirely determined by face components, as any component of $\phi$ is a component of one of its faces:
    \begin{equation}\label{eq:lem-n-dual-finite-eq-1}
        \phi^u|_{g(u)} = \phi^{\delta_j u'}|_{g(\delta_j u)} = (\widecheck{d}_j \phi)^{u'}|_{(d_j g)(u)},
    \end{equation}
    for some $j$ and $u$. Hence, $\mathfrak{A}_n(\mathsf{p}\huaV)_m|_g = \partial_m(\mathfrak{A}_n\mathsf{p}\huaV)|_{\theta_m g}$ for any $m \ge n+1$.
    In particular this means that the ambient $n$-dual bundle is $n$-coskeletal. 

    We now show that an element $\phi\in \mathfrak{A}_n\mathsf{p}\huaV_m|_g$ for $m \ge n+2$ satisfies the $n$-dual equations if and only if its faces $d_j\phi$ for $0 \le j \le m$ satisfy them at $d_jg$. 
    
    Consider the normalization equations \eqref{eq:VBnDualNormalizationMostGeneral} at level $m$. By \eqref{eq:lem-n-dual-finite-eq-1}, these can be written as
    \begin{equation*}
        \phi^{s_it}(\widetilde{s}_iv)= \phi^{\delta_j s_i t'}(\widetilde{s}_i v) = (d_j\phi)^{s_i t'}(\widetilde{s}_i v) = 0, 
    \end{equation*}
    for all $t = \delta_j t' \in \Delta^m_{n-1}$, $0 \le i \le n-1$, and all $v \in \mathsf{p}\huaV_n|_{g(t)}$, where $g(t) = (d_jg)(t')$. These are the equations at level $m-1$ for $d_j\phi$ at $d_j g$. 

    Analogously, consider the multiplicativity equations \eqref{eq:VBnDualMultiplicativityMostGeneral} at level $m$. By \eqref{eq:lem-n-dual-finite-eq-1}, we have 
    \begin{equation*}
        \sum_{i=0}^n \phi^{d_ir}(\widetilde{d}_i v) 
        = \sum_{i=0}^n \phi^{\delta_jd_ir'}(\widetilde{d}_i v) 
        = \sum_{i=0}^n (d_j\phi)^{d_ir'}(\widetilde{d}_i v),
    \end{equation*}
    for all $r = \delta_j r' \in \Delta^m_{n+1}$, and all $v \in \mathsf{p}\huaV_n|_{g(r)}$, where $g(r) = (d_jg)(r')$. These are the equations at level $m-1$ for $d_j\phi$ at $d_jg$. 
    
    With this we have that, for any $m \ge n+2$ and $g \in \mathsf{p}\huaG_m$,
    \begin{equation*}
        \mathsf{p}\huaV'_m|_g = \left\{ (\phi^u|_{g(u)})_{u\in \Delta^m_n}  
        \mid \forall 0 \le j \le m, \; ((d_j \phi)^{u'}|_{(d_jg)(u)})_{u'\in \Delta^{m-1}_n} \in \mathsf{p}\huaV'_{m-1}|_{d_j g}\right\}
    \end{equation*}
    so, $\mathsf{p}\huaV'_m|_g = \partial_m(\mathsf{p}\huaV')|_{\theta_m g}$ for any $m \ge n+2$.
\end{proof}

\subsection{The \texorpdfstring{$n$}{n}-dual pairing}

\begin{definition}\label{def:VBnDualPairingDef}
Let $\huaV \to \huaG$ be a higher vector bundle for which $\huaV^{n*}$ exists. We define the \textbf{$n$-dual pairing} as the simplicial pairing $\langle \cdot, \cdot \rangle: \huaV^{n*} \otimes \huaV \to B^n\R^\huaG$ defined by
\begin{equation}\label{eq:VBnDualPairingDef}
    \langle \phi|_g, v|_g \rangle = \tau(id)(\phi, v) = \phi^{E_n}(v),
\end{equation}
for any $g\in \huaG_n$, $\phi\in \huaV^{n*}_n|_g$, $v \in \huaV_n|_g$, where $\tau$ is the map 
\begin{equation*}
    \tau: \VB^{\infty}_\huaG(\huaV^{n*},\huaV^{n*}) \to \VB^\infty_\huaG(\huaV^{n*} \otimes \huaV, B^n\R)
\end{equation*}
defined in \eqref{eq:VBnDual-tau-def}.
\end{definition}

By definition of $\tau$, the $n$-dual pairing is a simplicial pairing, so it satisfies the multiplicativity and normalization conditions \eqref{eq:VBnShiftedPairingMultiplicative}, \eqref{eq:VBnShiftedPairingNormalized}. 
The construction of the associated IM-pairing to a pairing of higher vector bundles happens pointwise over $\huaG_0$, and the condition of homological nondegeneracy for such a pairing is also required to hold pointwise over $\huaG_0$, by Definitions \ref{def:VBIMPairingAssociatedDef} and \ref{def:VB-hndg-pairing}. 
Therefore, the following two properties of the $n$-dual pairing of higher vector bundles follow immediately from their respective version for simplicial vector spaces discussed in Theorem \ref{thm:ndual-pairing-hom-nondeg} and \ref{thm:hom-nondeg-homotopy-equiv}.

Note that at this point we only have the latter part of Theorem \ref{thm:hom-nondeg-homotopy-equiv} for $\VB$ $n$-groupoids with $n\le 2$, as the correspondence between weak equivalences of $\VB$ $n$-groupoids over the identity map of the base and quasi-isomorphisms of their normalized complexes was only established for $n=1$ in \cite[Thm 3.5]{HoyoOrtiz2020}, and we extended it to $n=2$ in Theorem \ref{thm:we-of-VB2gpd-is-qi}. 
We will discuss this in the sections relative to the 1-dual and the 2-dual.

\begin{theorem}\label{thm:VBndual-pairing-hndg}
    Let $\huaV \to \huaG$ be a higher vector bundle over a Lie $n$-groupoid $\huaG$, for which the $n$-dual $\huaV^{n*}$ exists. If $\huaV$ is a pointwise $n$-type, then the $n$-dual pairing $\langle\cdot , \cdot \rangle: \huaV^{n*}_n \otimes \huaV_n \to \R \times \huaG_n$ is homologically $n$-shifted nondegenerate.
\end{theorem}

\begin{proof}
    This is a pointwise condition for any $p \in \huaG_0$. By Lemma \ref{lem:VBndual-at-units}, the proof of Proposition \ref{thm:ndual-pairing-hom-nondeg} can be applied at each point $p \in \huaG_0$. Hence $\lambda_{\langle \cdot, \cdot \rangle}^r: N(\huaV^{n*})|_p \to N(\huaV)|_p^*[-n]$ is a quasi-isomorphism for each point $p\in \huaG_0$.
\end{proof}

\begin{theorem}\label{thm:VBnShiftedPairingsInducedMapsTriangle}
Let $\alpha: \huaV \otimes \huaW \to B^n\R^\huaG$ be a simplicial $n$-shifted pairing between two higher vector bundles over a Lie $n$-groupoid $\huaG$. 
Assume that $\huaV$ and $\huaW$ admit $n$-duals $\huaV^{n*}$ and $\huaW^{n*}$.
Then the following diagrams commute at each point $p\in \huaG_0$.
\begin{equation}\label{diag:VBnShiftedPairingsInducedMapsTriangle}
\begin{tikzcd}[ampersand replacement=\&]
	{N(\huaW)} \& {N(\huaV^{n*})} \& {N(\huaV)} \& {N(\huaW^{n*})} \\
	\& {N(\huaV)^*[-n]} \&\& {N(\huaW)^*[-n]}
	\arrow["{N(\alpha^l)}", from=1-1, to=1-2]
	\arrow["{\lambda_\alpha^l}"', from=1-1, to=2-2]
	\arrow["{\lambda_{\langle \cdot, \cdot \rangle}^r}", from=1-2, to=2-2]
	\arrow["{N(\alpha^r)}", from=1-3, to=1-4]
	\arrow["{\lambda_\alpha^r}"', from=1-3, to=2-4]
	\arrow["{\lambda_{\langle \cdot, \cdot \rangle}^r}", from=1-4, to=2-4]
\end{tikzcd}
\end{equation}
Where $\lambda^r_{\langle \cdot, \cdot \rangle}$ is the right induced map of the IM-pairing associated to the $n$-dual pairing of $\huaV$ on the left and of $\huaW$ on the right. 
\end{theorem}

\begin{proof}
    This follows by applying the proof of commutativity of the diagrams in Theorem \ref{thm:hom-nondeg-homotopy-equiv} pointwise at each $p\in \huaG_0$.
\end{proof}

\begin{remark}
    As in Remark \ref{rem:ndual-mult-pairing-n-gpd}, because the $n$-dual of a $\VB$ $n$-groupoid $\huaV$ satisfies the same sets of equations fiberwise as the $n$-dual of a $\VS$ $n$-groupoid, assuming the $n$-dual $\huaV^{n*}$ of $\huaV$ exists, its multiplication must satisfy equation \eqref{eq:ndual-mult-pairing-n-gpd} in each fiber as well. (Recall that this corresponds to \eqref{eq:ndual-mult-pairing}, the multiplicativity condition, for $E_{n+1}\in \Delta^{n+1}_{n+1}$). Therefore, the multiplication of the $\VB$ $n$-dual can be defined fiberwise as the unique multiplication for which the $n$-dual pairing is multiplicative, just as in the vector space case. 
\end{remark}

\subsection{Degeneracy annihilators and cores}\label{sec:VB-DegAnn-Cores}

As in Section \ref{sec:annihilators-and-dual-kernels} we now discuss the solution spaces of the normalization equations and their relation to the cores of a $\VB$ groupoid from Section \ref{sec:VB1-groupoids-cores-duals}.

By equation \eqref{eq:VBnDualNormalizationMostGeneral}, any component of a map $\phi \in \SVect(\textstyle\int g^*\huaV, B^n\R)$ indexed by a degenerate simplex $s_it \in \Delta^m_n$ annihilates the image of the degeneracy map $\widetilde{s}_i: \huaV_{n-1} \to \huaV_{n}$. 
Since this is a bundle map over $s_i: \huaG_{n-1} \to \huaG_n$ and not over the identity, we need to make clear what we mean by its ``image''. Namely, we use the universal property of pullbacks 
\[\begin{tikzcd}[ampersand replacement=\&,cramped, sep=scriptsize]
	{\huaV_{n-1}} \& {s_i^*\huaV_n} \& {\huaV_n} \\
	{\huaG_{n-1}} \& {\huaG_{n-1}} \& {\huaG_n}
	\arrow["{\exists!}"', dashed, from=1-1, to=1-2]
	\arrow["{\widetilde{s}_i}", curve={height=-12pt}, from=1-1, to=1-3]
	\arrow[from=1-1, to=2-1]
	\arrow[from=1-2, to=1-3]
	\arrow[from=1-2, to=2-2]
	\arrow["\lrcorner"{anchor=center, pos=0.125}, draw=none, from=1-2, to=2-3]
	\arrow[from=1-3, to=2-3]
	\arrow["id", from=2-1, to=2-2]
	\arrow["{s_i}", from=2-2, to=2-3]
\end{tikzcd}\]
to identify $\widetilde{s}_i$ with the unique map $\huaV_{n-1} \to s_i^*\huaV_{n}$ over the identity, by an abuse of notation. 
Hence we say that the \textbf{image} of $\widetilde{s}_i$ is the subbundle $\widetilde{s}_i\huaV_{n-1} \subseteq s_i^*\huaV_n$ over $\huaG_{n-1}$. 
Now we can define the annihilator subbundle $\Ann(\widetilde{s}_i\huaV_{n-1}) \subseteq s_i^*\huaV_n^* \to \huaG_{n-1}$. More explicitly, its fiber at each $g \in \huaG_{n-1}$ is given by
\begin{equation*}
    \Ann(\widetilde{s}_i\huaV_{n-1})|_g:= \{ \phi \in \huaV_n^*|_{s_i g} \quad \mid \quad \phi(\widetilde{s}_iv) = 0, \quad \forall v \in \huaV_{n-1}|_g\}.
\end{equation*}

Now we have that for any $\phi \in \SVect(\textstyle\int g^*\huaV, B^n\R)$, \eqref{eq:VBnDualNormalizationMostGeneral} is equivalent to 
\begin{equation}\label{eq:VBnDualNormalizationAnnihilator}
    \phi^{s_i t}|_{g(s_it)} \in \Ann(\widetilde{s}_i\huaV_{n-1})|_{g(t)},
\end{equation}
for all $0 \le i \le n-1$ and $t \in \Delta^m_{n-1}$. In analogy with Definition \ref{def:VS-cores-deg-ann}, we define the following. 

\begin{definition}
    Let $\huaV$ be a higher vector bundle. For any $0\le i \le n-1$, the associated \textbf{degree 1 $n$-dimensional degeneracy annihilator} is the vector bundle $O_i:=\Ann(\widetilde{s}_i\huaV_{n-1}) \subseteq s_i^*\huaV_n$ over $\huaG_{n-1}$.
\end{definition}

Analogous sequences to \eqref{eq:VSnDualSESd_is_i} and \eqref{eq:VSnDualSESd_iplus1s_i} appear in this setting as well. These are the following canonically split short exact sequences of vector bundles over $\huaG_{n-1}$, for any $0 \le i < n$.\footnote{For $n=1$, this is the same as \eqref{eq:VB1DualCoreSES}.} 

\begin{equation}\label{eq:VBnDualCoreSES-1}
\begin{tikzcd}[ampersand replacement=\&,cramped,row sep=scriptsize]
	0 \& {s_i^*\ker\widetilde{d}_i} \& {s_i^*\huaV_n} \& {\huaV_{n-1}} \& 0 \\
	\& {\huaG_{n-1}} \& {\huaG_{n-1}} \& {\huaG_{n-1}}
	\arrow[from=1-1, to=1-2]
	\arrow["{incl}", shift left, hook, from=1-2, to=1-3]
	\arrow[from=1-2, to=2-2]
	\arrow["{id-\widetilde{s}_i\widetilde{d}_i}", shift left, from=1-3, to=1-2]
	\arrow["{\widetilde{d}_i}", shift left, from=1-3, to=1-4]
	\arrow[from=1-3, to=2-3]
	\arrow["{\widetilde{s}_i}", shift left, from=1-4, to=1-3]
	\arrow[from=1-4, to=1-5]
	\arrow[from=1-4, to=2-4]
	\arrow["id", Rightarrow, no head, from=2-2, to=2-3]
	\arrow["id", Rightarrow, no head, from=2-3, to=2-4]
\end{tikzcd}
\end{equation}

\begin{equation}\label{eq:VBnDualCoreSES-2}
\begin{tikzcd}[ampersand replacement=\&,cramped,row sep=scriptsize]
	0 \& {s_i^*\ker\widetilde{d}_{i+1}} \& {s_i^*\huaV_n} \& {\huaV_{n-1}} \& 0 \\
	\& {\huaG_{n-1}} \& {\huaG_{n-1}} \& {\huaG_{n-1}}
	\arrow[from=1-1, to=1-2]
	\arrow["{incl}", shift left, hook, from=1-2, to=1-3]
	\arrow[from=1-2, to=2-2]
	\arrow["{id-\widetilde{s}_i\widetilde{d}_{i+1}}", shift left, from=1-3, to=1-2]
	\arrow["{\widetilde{d}_{i+1}}", shift left, from=1-3, to=1-4]
	\arrow[from=1-3, to=2-3]
	\arrow["{\widetilde{s}_i}", shift left, from=1-4, to=1-3]
	\arrow[from=1-4, to=1-5]
	\arrow[from=1-4, to=2-4]
	\arrow["id", Rightarrow, no head, from=2-2, to=2-3]
	\arrow["id", Rightarrow, no head, from=2-3, to=2-4]
\end{tikzcd}
\end{equation}

The retracts in the short exact sequences define an isomorphism between the two kernels in the same way as in Lemma \ref{lem:VB1-core-involution-iso}. 
In fact, if $v \in s_i^*\ker\widetilde{d}_i$, $(id - \widetilde{s}_i \widetilde{d}_{i+1})v \in s_i^*\ker\widetilde{d}_{i+1}$, and 
\begin{equation*}
    (id - \widetilde{s}_i \widetilde{d}_{i})(id - \widetilde{s}_i \widetilde{d}_{i+1}) v = v 
    - \widetilde{s}_i \widetilde{d}_{i} v 
    - \widetilde{s}_i \widetilde{d}_{i+1} v 
    + \widetilde{s}_i \widetilde{d}_{i+1} v 
    = v.
\end{equation*}
In the same way, $(id - \widetilde{s}_i \widetilde{d}_{i+1})(id - \widetilde{s}_i \widetilde{d}_{i}) = id$. 
To employ a terminology consistent with this situation, we make the following definition. 

\begin{definition}
    Let $\huaV$ be a higher vector bundle. For any pair $(i,j)$, with $0\le i \le n-1$, $0\le j \le n$, the associated \textbf{degree 1 $n$-dimensional core} is the vector bundle $s_i^*\ker d_j \subseteq s_i^*\huaV_n$ over $\huaG_{n-1}$.
    
    For any such $i$, we call $(s_i^*\ker d_i, s_i^*\ker d_{i+1})$ the \textbf{core pair} associated to $O_i$. 
\end{definition} 

Each core pair is associated to a degeneracy annihilator because from the point of view of the dual sequences to \eqref{eq:VBnDualCoreSES-1} and \eqref{eq:VBnDualCoreSES-2}, we have the isomorphisms
\begin{equation}\label{eq:VBnDual-deg1-degAnn-isom}
    s_i^*(\ker d_i)^*
    \newrightleftarrows{(id-\widetilde{s}_i\widetilde{d}_{i})^*}{incl^*} 
    \Ann(\widetilde{s}_i\huaV_{n-1}) 
    \newrightleftarrows{incl^*}{(id-\widetilde{s}_i\widetilde{d}_{i+1})^*} 
    s_i^*(\ker d_{i+1})^*, \quad \forall 0\le i < n.
\end{equation}
Unlike for the isomorphisms in \eqref{eq:VSDual-isos-deg-ann-dual-ker}, there are no canonical isomorphisms between the pullbacks $s_i^*(\ker d_i)^*$ and $s_{i-1}^*(\ker d_i)^*$ for $0 < i < n$. Thus, a priori, the $O_i$ are not all isomorphic for all $i$. This is why we separate the cores into pairs. Diagrammatically, we have
\[
\begin{adjustbox}{width=\textwidth}
\begin{tikzcd}[ampersand replacement=\&,cramped,sep=tiny]
	\& {O_0} \&\&\& {O_1} \&\&\&\& {O_{n-1}} \\
	{s_0^*(\ker \widetilde{d}_0)^*} \&\& {s_0^*(\ker \widetilde{d}_1)^*} \& {s_1^*(\ker \widetilde{d}_1)^*} \&\& {s_1^*(\ker \widetilde{d}_2)^*} \& \dots \& {s_{n-1}^*(\ker \widetilde{d}_n)^*} \&\& {s_{n-1}^*(\ker \widetilde{d}_n)^*}
	\arrow["\cong"', tail reversed, from=1-2, to=2-1]
	\arrow["\cong", tail reversed, from=1-2, to=2-3]
	\arrow["\cong"', tail reversed, from=1-5, to=2-4]
	\arrow["\cong", tail reversed, from=1-5, to=2-6]
	\arrow["\cong"', tail reversed, from=1-9, to=2-8]
	\arrow["\cong", tail reversed, from=1-9, to=2-10]
\end{tikzcd}
\end{adjustbox}
\]
Naturally, in this terminology, the left and right cores of a $\VB$ 1-groupoid are the 1-dimensional degree 1 cores that form the core pair associated to $O_0 = \Ann(\widetilde{1}\huaV_0)$.

We now move to discussing higher degree annihilators and cores. As in Section \ref{sec:annihilators-and-dual-kernels}, after a brief general discussion, we only focus on the $n$-dimensional degree $n$ ones, since these are the only other kind of core we need to compute the 2-dual of a $\VB$ 2-groupoid. 
We leave a discussion of more general cores for future work. 

To begin with, there are only certain fibers of $\huaV \to \huaG$ which contain degenerate elements with respect to $\widetilde{s}_i$ for multiple different values of $i$. Such an element is for example $\widetilde{s}_Jv$ for some multi-index $J$ of length $|J|$ such that $v$ is a nondegenerate element of $\huaV_{n-|J|}|_{g}$ for some $g \in \huaG_{n-|J|}$. Then $\widetilde{s}_Jv \in \huaV_{\widetilde{s}_Jg}$. By the simplicial identities one can rewrite $\widetilde{s}_Jv$ to show it is in the image of $s_i$ for $|J|$ different values of $i$. If $I$ is the set of these indices $i$, then $s_J^*\huaV_n^*|_{g}$ contains the corresponding fiber of the degree 1 degeneracy annihilators $O_i|_g=\Ann(\widetilde{s}_i\huaV_{n-1})|_g$ for all $i\in I$. These can then be intersected to obtain the fiber at $g$ of the \textit{$n$-dimensional degree $|J|$ degeneracy annihilator} $O_{I}$. 

We will only limit ourselves to this heuristic definition and now consider the case of $|J|=n$. In this case, for any multi-index $J$ for which $\widetilde{s}_J$ makes sense, $\widetilde{s}_J = \widetilde{1}$, as in Remark \ref{rem:total-unit}. Hence we make the following definition.

\begin{definition}
    Let $\huaV\to \huaG$ be a higher vector bundle. 
    
    The \textbf{degree $n$ $n$-dimensional degeneracy annihilator} is the vector subbundle $O_{0\dots (n-1)}:= \Ann(D_n\huaV) \subseteq 1^*\huaV_n^*$ over $\huaG_0$, where $D_n\huaV \subseteq 1^*\huaV_n \to \huaG_0$ is the subbundle spanned by the degenerate $n$-simplices in $1^*\huaV$.

    For any $0 \le k \le n$, the associated \textbf{degree $n$ $n$-dimensional} core is the vector bundle $1^*\ker\widetilde{p}^n_k \subseteq 1^*\huaV_n$ over $\huaG_0$. 
\end{definition}

In this case we are over the total units, so the same discussion as in Section \ref{sec:annihilators-and-dual-kernels} applies: the dual of each degree $n$ $n$-dimensional core is isomorphic to $O_{0\dots (n-1)}$. 
The isomorphism is given by the dual of one of the core projections $\ggamma^n_k:= id - \mu^n_k$ defined fiberwise as in \eqref{eq:VSnDualSESp_kmu_k}, by using the fiberwise Moore fillers. 
In other words, for any $0 \le k \le n$,
\begin{equation}\label{eq:VBnDual-degn-degAnn-isom}
    O_{0\dots (n-1)} \cong 1^*(\ker p^n_k)^*. 
\end{equation}

\subsection{The dual \texorpdfstring{$\VB$}{VB}-groupoid revisited}\label{sec:VB1Dual-revisited}

We now show that the dual $\VB$ groupoid from Definition \ref{def:VB1Dual} is the 1-dual of a $\VB$ groupoid in the more general sense of Definition \ref{def:n-dual-univ-prop}.

\begin{lemma}\label{lem:VB1Dual-computation}
    Let $\huaV \to \huaG$ be a $\VB$ groupoid. Then $\huaV^{1*}$ as defined in Definition \ref{def:VB1Dual} is a simplicial subbundle of $\mathfrak{A}_1\huaV$ such that
    \begin{equation}
        \huaV^{1*}_m|_{g} := \SVect(\textstyle\int g^*\huaV, B^1\R),
    \end{equation}
    for any $g\in \huaG_m$ and $m\ge 0$.
\end{lemma}

\begin{proof}
    First of all, observe there is an embedding $\huaV^{1*} \into \mathfrak{A}_1\huaV$, which is generated by the $1$-truncated embedding
    \begin{equation*}
        \begin{split}
            \huaV^{1*}_0 = \Ann(\widetilde{1}\huaV_0) &\overset{incl}{\into} 1^*\huaV_1^* = \mathfrak{A}_1\huaV\\
            \huaV^{1*}_1 = \huaV_1 &\into (\huaV_1^*)^{01}\oplus \bigoplus_{i=0}^1((1d_i)^*\huaV_1^*)^{ii} = \mathfrak{A}_1\huaV\\
            \xi &\mapsto ((\xi)^{01}, (\widecheck{d}_1 \xi)^{00}, (\widecheck{d}_0 \xi)^{11}).
        \end{split}
    \end{equation*}
    We now solve the equations for the first three levels of the 1-dual, following Lemma \ref{lem:n-dual-finite-equations}.
    
    For level 0, let $p \in \huaG_0$. By Lemma \ref{lem:VBndual-at-units}, 
    \begin{equation*}
    \begin{split}
        \SVect(\textstyle\int p^*\huaV, B^1\R) &\cong \SVect((1^*\huaV)|_p \otimes \Delta^0, B^1\R)\\
        &= ((1^*\huaV)|_p)^{1*} = \Ann(\widetilde{1}(\huaV_0|_p)) \subseteq \huaV_1^*|_{1p},
    \end{split}
    \end{equation*}
    which is the fiber at $p$ of the vector bundle $\Ann(\widetilde{1}(\huaV_0)) \to \huaG_0$, and it coincides with $\huaV^{1*}_0$.

    For level 1, let $g \in \huaG_1$. Any $\xi \in \SVect(\textstyle\int g^*\huaV, B^1\R)$ has three components: $\xi^{00}, \xi^{01}, \xi^{11}$. By \eqref{eq:VBnDualNormalizationMostGeneral} for $t=0$ and $t=1$,
    \begin{equation*}
        \xi^{00}(\widetilde{1}x) = 0, 
        \quad \forall x \in \huaV_0|_{1d_1g},
        \quad \text{and }  \xi^{11}(\widetilde{1}x') = 0, 
        \quad \forall x' \in \huaV_0|_{1d_0g},
    \end{equation*}
    So $\xi^{00} \in \Ann(\widetilde{1}\huaV_0)|_{d_1g}$ and $\xi^{11} \in \Ann(\widetilde{1}\huaV_0)|_{d_0g}$. On the other hand, the multiplicativity equations at this level are \eqref{eq:VBnDualMultiplicativityMostGeneral} for $r=000, 001, 011, 111$. Those for $r=000,111$ are automatically satisfied by the same argument as in the proof of Proposition \ref{prop:VS1dual}, since we are over units. 
    Meanwhile, the equations for $r=001$ and $r=011$, take inputs $(v|_g,w|_{1d_1g})\in\Lambda^2_1(\huaV)|_{s_0g}$, and $(v'|_{1d_0g},w'|_g)\in\Lambda^2_1(\huaV)|_{s_1g}$, respectively, and read
    \begin{equation*}
        \begin{cases}
            \xi^{01}((w \cdot v)|_g) 
            = \xi^{01}(v|_g) + \xi^{00}(w|_{1d_1g}), \\
            \xi^{01}((w' \cdot v')|_g) 
            = \xi^{11}(v'|_{1d_0g}) + \xi^{01}(w'|_{g}). \\
        \end{cases}
    \end{equation*}
    Equivalently, 
    \begin{equation}\label{eq:VB1dualMult001011-1}
        \begin{cases}
            \xi^{01}((w \cdot v)|_g - v|_g) 
            = \xi^{00}(w|_{1d_1g}), \\
            \xi^{01}((w' \cdot v')|_g - w'|_{g}) 
            = \xi^{11}(v'|_{1d_0g}). \\
        \end{cases}
    \end{equation}
    By linearity of the multiplication and the horn conditions $\widetilde{1}\widetilde{d}_1v = \widetilde{1}\widetilde{d}_0w$ and $\widetilde{1}\widetilde{d}_1v' = \widetilde{1}\widetilde{d}_0w'$, 
    \begin{equation*}
        \begin{split}
            &(w \cdot v)|_g - v|_g 
            = (w \cdot v)|_g - (\widetilde{1}\widetilde{d}_1 v \cdot v)|_g 
            = (w - \widetilde{1}\widetilde{d}_0 w) \cdot 0_g,\\
            &(w' \cdot v')|_g - w'|_g 
            = (w' \cdot v')|_g - (w' \cdot \widetilde{1}\widetilde{d}_0 w')|_g 
            = 0_g \cdot (v' - \widetilde{1}\widetilde{d}_1 v').
        \end{split}
    \end{equation*}
    With this, \eqref{eq:VB1dualMult001011-1} is equivalent to 
    \begin{equation*}
        \begin{cases}
            \xi^{00}(w|_{1d_1g}) 
            = \xi^{01}((w - \widetilde{1}\widetilde{d}_0 w) \cdot 0_g),\\
            \xi^{11}(v'|_{1d_0g}) 
            = \xi^{01}(0_g \cdot (v' - \widetilde{1}\widetilde{d}_1 v')),
        \end{cases}
    \end{equation*}
    for any $w \in \huaV_1|_{1d_1g}$ and any $v' \in \huaV_1|_{1d_0g}$.
    This is precisely the definition of the face maps in Definition \ref{def:VB1Dual}, and it coincides fiberwise with \eqref{eq:VBndual-ambient-face-maps}. Moreover, we just showed that any $\xi \in \SVect(\textstyle\int g^*\huaV, B^1\R)$ is determined uniquely by $\xi^{01} \in \huaV_1^*|_g$, so $\SVect(\textstyle\int g^*\huaV, B^1\R) = \huaV_1^*|_g = \huaV^{1*}_1|_{g}$.
    By comparing the definition of the unit map in Definition \ref{def:VB1Dual} and \eqref{eq:VBndual-ambient-deg-maps}, it is clear that the groupoid unit maps also coincide fiberwise. 

    By Proposition \ref{prop:finite-data-n-gpds} and Lemma \ref{lem:n-dual-finite-equations} it is now enough to show that fiberwise, the multiplication defined in Definition \ref{def:VB1Dual} coincides with the one from \eqref{eq:ndual-mult-pairing-n-gpd}. 
    In this case, the 1-dual pairing is simply the pairing of $\huaV_1$ with its dual, so the two formulas coincide.
\end{proof}

We now state some consequences of the theory of $n$-shifted pairings and Theorem \ref{thm:we-of-VB2gpd-is-qi} in the case of $\VB$ 1-groupoids, which was shown in \cite[Thm. 3.5]{HoyoOrtiz2020}. We begin with recalling the simplicial maps induced by a $1$-shifted simplicial pairing. 
Let $\alpha: \huaV_1 \otimes \huaW_1 \to \R \times \huaG_1$ be a multiplicative normalized 1-shifted pairing between two $\VB$ 1-groupoids over $\huaG$. Then, by Remark \ref{rem:VB-n-pairing-ind-maps}, $\alpha$ induces the maps
\begin{equation*}
    \begin{array}{cc}
        \alpha^r : \huaV \to \huaW^{1*},
        &\alpha^l : \huaW \to \huaV^{1*}
    \end{array}
\end{equation*}
defined by
\begin{equation}\label{eq:VB-1-shifted-pairing-ind-maps}
    \begin{array}{cc}
        \alpha^r(x_p)(w_{1p}) = \alpha(\widetilde{s}_0 x_p, w_{1p}), &\alpha^l(y_p)(v_{1p}) = \alpha(v_{1p}, \widetilde{s}_0y_p),\\
        \alpha^r(v_g)(w_g) = \alpha(v_g,w_g),
        & \alpha^l(w_g)(v_g) = \alpha(v_g,w_g),
    \end{array}
\end{equation}
for all $x_p \in \huaV_0|_p$, $w_{1p} \in \huaW_1|_{1p}$, $y_p \in \huaW_0|_p$, $v_{1p} \in \huaV_1|_{1p}$, $v_{g} \in \huaV_1|_g$ and $w_{g}\in \huaW_1|_g$. 

\begin{remark}
By \eqref{eq:VB-1-shifted-pairing-ind-maps}, the 1-dual pairing $\langle \cdot, \cdot \rangle$ induces two isomorphisms:
\begin{equation*}
    \langle \cdot, \cdot \rangle^r=id: \huaV^{1*} \to \huaV^{1*},
    \quad \langle \cdot, \cdot \rangle^l: \huaV \to (\huaV^{1*})^{1*}. 
\end{equation*}
In fact, assuming $\huaV$ is a $\VB$ 1-groupoid, the 1-dual pairing is not only homologically nondegenerate, but the induced IM-pairing $\lambda_{\langle\cdot, \cdot\rangle}$ is even nondegenerate at the level of chains. 
This is a peculiarity of the $n=0$ and $n=1$ cases of $n$-duals of $\VB$ $n$-groupoids. In all other cases the $n$-dual pairing is a priori only nondegenerate in the homology, and one should not expect the $n$-duality to be reflexive on the nose. 
We will show in Corollary \ref{cor:VB-2dual-reflexive-uth} that 2-duality of $\VB$ 2-groupoids is reflexive up to homotopy, as $n$-duality of $\VS$ $n$-types is, by Theorem \ref{thm:ndual-reflexive-uth}. 
A more general reflexivity result of $n$-duality for general $\VB$ $n$-groupoids would be immediate after extending Theorem \ref{thm:we-of-VB2gpd-is-qi} to any $n$. We plan to complete this in future work. 
\end{remark}

We conclude with a characterization of homologically nondegenerate 1-\hspace{0pt}shifted simplicial pairings in terms of weak equivalences (i.e. Morita maps in the terminology of \cite{HoyoOrtiz2020}).

\begin{corollary}
    Let $\alpha: \huaV_1 \otimes \huaW_1 \to \R \times \huaG_1$ be a multiplicative normalized 1-shifted pairing between two $\VB$ 1-groupoids over $\huaG$. Then the following are equivalent. 
    \begin{enumerate}[label=(\roman*)]
        \item $\alpha$ is homologically nondegenerate.
        \item $\alpha^r$ is a weak equivalence.
        \item $\alpha^l$ is a weak equivalence.
    \end{enumerate}
    In particular, if $\omega_2$ is the leading term of a 1-shifted presymplectic form on a Lie 1-groupoid $\huaG$, it is 1-shifted symplectic if and only if $\omega_2^r = -\omega_2^l$ is a weak equivalence. 
\end{corollary}
\begin{proof}
    This follows immediately from \cite[Thm. 3.5]{HoyoOrtiz2020} and commutativity of the diagram \eqref{diag:VBnShiftedPairingsInducedMapsTriangle} from Theorem \ref{thm:VBnShiftedPairingsInducedMapsTriangle}, since $\lambda^r_{\langle,\rangle}$ is an isomorphism. In this situation, $\lambda^r_\alpha$ is a quasi-isomorphism if and only if $N(\alpha^r)$ is one as well, and this is equivalent to $\alpha^r$ being a weak equivalence. The same argument shows that (i) is equivalent to (iii).
\end{proof}

\section{The 2-dual of a \texorpdfstring{$\VB$}{VB} 2-groupoid}\label{sec:2-dual-VB-gpd}

We construct a model of the 2-dual of a $\VB$ 2-groupoid and show this is a well-defined finite data set for a $\VB$ 2-groupoid as described in Example \ref{ex:finite-data-2gpds}. We then show it is indeed a $\VB$ 2-groupoid and that it satisfies the $2$-dual equations. First, we introduce some necessary spaces and isomorphisms between them. 

Let $\huaV \to \huaG$ be a $\VB$ 2-groupoid. 
In the definition of the 2-dual $\VB$ groupoid we will use:
\begin{itemize}
    \item The 2-dimensional degree 2 degeneracy annihilator 
    \begin{equation*}
        O_{01} := \Ann(D_2\huaV) \subseteq 1^*\huaV_2^* \to \huaG_0,
    \end{equation*}
    which, by \eqref{eq:VBnDual-degn-degAnn-isom} is isomorphic to the dual $1^*(\ker \widetilde{p}^2_j)^*$ of each of the degree 2 2-dimensional cores for $j=0,1,2$.

    \item The 2-dimensional degree 1 degeneracy annihilators 
    \begin{equation*}
        O_0 := \Ann(\widetilde{s}_0 \huaV_1) \subseteq s_0^*\huaV_2^* \to \huaG_1, 
        \qquad
        O_1 := \Ann(\widetilde{s}_1 \huaV_1) \subseteq s_1^*\huaV_2^* \to \huaG_1, 
    \end{equation*}
    which are isomorphic to the dual core pairs $s_0^*(\ker\widetilde{d}_0)^* \cong s_0^*(\ker\widetilde{d}_1)^*$ and $s_1^*(\ker\widetilde{d}_1)^* \cong s_1^*(\ker\widetilde{d}_2)^*$, respectively, by \eqref{eq:VBnDual-deg1-degAnn-isom}.
    \item The pullback by $s_0^*$ of the dual of the kernel of the $(2,1)$-horn projection 
    \begin{equation*}
        s_0^*(\ker \widetilde{p}^2_1)^* \to \huaG_1,
    \end{equation*}
    with the isomorphism $\zeta^*: s_1^*(\ker\widetilde{p}^2_1)^* \to s_0^*(\ker\widetilde{p}^2_1)^*$ which is the dual of the translation
    \begin{equation}\label{eq:VB2DualTranslationZetaDef}
        \begin{split}
            \zeta:\, s_0^*\ker\widetilde{p}^2_1 &\longrightarrow s_1^*\ker\widetilde{p}^2_1\\
            k_{s_0g} &\longmapsto 0_{s_1g} \square k_{s_0g} 0_{s_0g},
        \end{split}
    \end{equation}
    over the identity of $\huaG_1$.
    Since 
    \begin{equation*}
    s_0^*(\ker\widetilde{p}^2_1)^* \subseteq s_0^*(\ker\widetilde{d}_0)^* \cong O_0, 
    \text{ and }
    s_1^*(\ker\widetilde{p}^2_1)^* \subseteq s_1^*(\ker\widetilde{d}_2)^* \cong O_1,
    \end{equation*}
    this pair of spaces is a reasonable substitute for the ``intersection'' of $O_0$ and $O_1$ over $\huaG_1$, which is not well-defined otherwise.\footnote{In the case of a $\VS$-groupoid there is no pullback and $(\ker \widetilde{p}^2_1)^* \cong O_{01}$, which \textit{is} the intersection of $O_0$ and $O_1$.}
    This will allow us to construct a fiber product analogous to that appearing at level 1 of the $\VS$ 2-dual from Theorem \ref{thm:VS2Dual}. 
    The subbundle $s_0^*\ker\widetilde{p}^2_1$ is also isomorphic to $s_0^*\ker\widetilde{p}^2_0$ and $s_1^*\ker\widetilde{p}^2_2$, through isomorphisms similar to $\zeta$: they are multiplications by two zero triangles over the degenerate tetrahedron $s_0s_1g$. By classifying all tetrahedra in $\huaV$ with two zero faces over degenerate tetrahedra in $\huaG$, one obtains the following diagram of canonical translations over $\huaG_1$. 
    
    \adjustbox{width=0.93\textwidth}{
    \begin{tikzcd}[ampersand replacement=\&,cramped]
	\&\& {s_0^*\ker\widetilde{p}_1^2} \&\&\&\& {s_1^*\ker\widetilde{p}_1^2} \&\& \\
	\\
	{s_0^*\ker\widetilde{p}_0^2} \&\&\&\&\&\&\&\& {s_1^*\ker\widetilde{p}_2^2} \\
	\&\& {s_0^*\ker\widetilde{p}_2^2} \&\&\&\& {s_1^*\ker\widetilde{p}_0^2} \\
	{d_1^*\ker\widetilde{p}^2_1} \&\&\&\&\&\&\&\& {d_0^*\ker\widetilde{p}^2_1} \\
	{d_1^*\ker\widetilde{p}^2_0} \&\& {d_1^*\ker\widetilde{p}^2_2} \&\&\&\& {d_0^*\ker\widetilde{p}^2_0} \&\& {d_0^*\ker\widetilde{p}^2_2}
	\arrow["{0_{s_1g}\square\_0_{s_0g}}", from=1-3, to=1-7]
	\arrow["{\_0_{s_0g}\square0_{1d_1g}}"', from=1-3, to=3-1]
	\arrow["{\_0_{s_1g}\square0_{s_0g}}"', from=1-7, to=3-1]
	\arrow["{\square 0_{s_1g}0_{s_0g}\_}", from=3-1, to=3-9]
	\arrow["{0_{s_1g}\_ 0_{s_0g}\square}"', from=3-9, to=1-3]
	\arrow["{0_{1d_0g}\_ 0_{s_1g}\square}"', from=3-9, to=1-7]
	\arrow["{0_{s_0g}\square 0_{s_0g}\_}", from=5-1, to=4-3]
	\arrow["{\_0_{s_1g}\square0_{s_1g}}"', from=5-9, to=4-7]
	\arrow["{\square 0_{s_0g}0_{s_0g}\_}"', from=6-1, to=4-3]
	\arrow["{0_{s_0g}0_{s_0g}\square\_}"', from=6-3, to=4-3]
	\arrow["{\_\square 0_{s_1g}0_{s_1g}}", from=6-7, to=4-7]
	\arrow["{\_0_{s_1g}0_{s_1g}\square}", from=6-9, to=4-7]
    \end{tikzcd}
    }

    In particular, it can be seen that not all pullbacks of the degree 2 2-\hspace{0pt}dimensional cores over $\huaG_1$ are canonically isomorphic to each other this way. For example, there are no canonical translations that give isomorphisms between the spaces $s_0^*\ker\widetilde{p}^2_1$ and $s_1^*\ker\widetilde{p}^2_1$ and any pullback by a face map of the degree 2 cores $1^*\ker\widetilde{p}^2_j$. 

    \item Finally, the subbundles $(\ker \widetilde{p}^2_0)^*$ and $(\ker \widetilde{p}^2_2)^*$ of $\huaV_2^* \to \huaG_2$. These will appear in the fiber products at level 2 of the $\VB$ 2-dual. To construct these we will use the isomorphisms $\nu_0^*$ and $\nu_2^*$ obtained by dualizing the translations 
    \begin{equation}\label{eq:VB2DualLevel2FiberProdIso0}
    \begin{split}
        \nu_0: \, \ker \widetilde{p}^2_0 &\longrightarrow d_0^*s_0^*\ker \widetilde{p}^2_1\\
        k_{t} &\longmapsto - \square 0_t k_t 0_{s_1d_2t},
    \end{split}
    \end{equation}
    and
    \begin{equation}\label{eq:VB2DualLevel2FiberProdIso2}
    \begin{split}
        \nu_2: \, \ker \widetilde{p}^2_2 &\longrightarrow d_2^*s_0^*\ker \widetilde{p}^2_1\\
        k_{t} &\longmapsto - 0_t k_{t} 0_{s_0d_1t} \square.
    \end{split}
    \end{equation}
    We observe that, by classifying canonical translations over $\huaG_2$ in a similar way as in the previous point, there are no canonical translations that give isomorphisms between the spaces $\ker \widetilde{p}^2_0$ and $\ker \widetilde{p}^2_2$ and any pullback by a composition of face maps of the degree 2 cores $1^*\ker\widetilde{p}^2_j$. 
\end{itemize}

With these spaces and isomorphisms, we are able to define the $\VB$ 2-dual.

\begin{definition}\label{def:VB2dual}
    Let $\huaV \to \huaG$ be a $\VB$ 2-groupoid. Its 2-dual $\huaV^{2*} \to \huaG$ is 
\[\begin{tikzcd}[ampersand replacement=\&,column sep=small,row sep=scriptsize]
	{\begin{array}{c}\huaV_2^*\times_{i^*, (\ker\widetilde{p}^2_0)^*, \nu_0^*} d_0^*\left(O_0\times_{i^*, (s_0^*\ker \widetilde{p}^2_1)^*, \zeta^*} O_1\right)\\ \times_{i^*, (\ker\widetilde{p}^2_2)^*, \nu_2^*} d^*_2 O_0\end{array}} \& {O_0\times_{i^*, (s_0^*\ker \widetilde{p}^2_1)^*, \zeta^*} O_1} \& {O_{01}} \\
	{\huaG_2} \& {\huaG_1} \& {\huaG_0}
	\arrow[shift left=2, from=1-1, to=1-2]
	\arrow[from=1-1, to=1-2]
	\arrow[shift right=2, from=1-1, to=1-2]
	\arrow[shift left=1, from=1-2, to=1-3]
	\arrow[shift right=1, from=1-2, to=1-3]
	\arrow[from=1-2, to=2-2]
	\arrow[shift left=1, from=2-2, to=2-3]
	\arrow[shift right=1, from=2-2, to=2-3]
	\arrow[shift left=2, from=2-1, to=2-2]
	\arrow[shift right=2, from=2-1, to=2-2]
	\arrow[from=2-1, to=2-2]
	\arrow[from=1-1, to=2-1]
	\arrow[from=1-3, to=2-3]
\end{tikzcd}\]
Here $O_{01} = \Ann(D_2\huaV)$, $O_0 = \Ann(\widetilde{s}_0\huaV_1) \subset s_0^*\huaV_2^*$ and $O_1 = \Ann(\widetilde{s}_1\huaV_1)\subset s_1^*\huaV_2^*$. The elements at level 1 are of the form
\begin{equation}\label{eq:VB2DualElementsLv1Def}
\begin{split}
    (\eta_{s_1g}, \xi_{s_0g}) \in &(O_0 \oplus O_1)_g \quad \text{ such that }\\ 
    &\forall k_{s_0g}\in s_0^*\ker\widetilde{p}_1^2, \quad \eta_{s_0g}(k_{s_0g})=\xi_{s_1g}(0_{s_1g} \square k_{s_0g} 0_{s_0g}),
\end{split}
\end{equation}
and those at level 2 are of the form
\begin{equation}\label{eq:VB2DualElementsLv2Def}
    \begin{split}
        &(\phi^{012}_t, (\phi^{112}_{s_0d_0t},\phi^{122}_{s_1d_0t}), \phi^{001}_{s_0d_2t}), \text{ such that }\\
        &\qquad\qquad \phi^{012}(k_t^0) = - \phi^{112}_{s_0d_0t}(\square 0_t k^0_t 0_{s_1d_2t}) = - \phi^{122}_{s_1d_0t}(\square 0_{s_1d_1t} k^0_t 0_t), \\
        &\qquad\qquad \text{ and } \phi^{012}(k_t^2) = - \phi^{001}_{s_0d_0t}(0_t k^2_{t} 0_{s_0d_1t} \square),
    \end{split}
\end{equation}
for any $k_t^0 \in \ker\widetilde{p}^2_0|_t$ and $k_t^2 \in \ker\widetilde{p}^2_2|_t$. 
The simplicial maps between levels 1 and 0 are 
\begin{equation}\label{eq:VB2DualFaceMaps1}
    \begin{split}
    \widecheck{d}_0^1(\eta_{s_0g},\xi_{s_1g})(X_{1d_0g})
    &= 
    \xi_{s_1g}
    ((\ggamma_0X_{1d_0g})\square 0_{s_1g}0_{s_1g})
    ,\quad \forall X_{1d_0g} \in (\huaV_2)_{1d_0g},\\
    \widecheck{d}_1^1(\eta_{s_0g},\xi_{s_1g})
    (X_{1d_1g})
    &=
    \eta_{s_0g}
    (0_{s_0g}0_{s_0g}\square (\ggamma_2X_{1d_1g}))
    ,\quad \forall X_{1d_1g} \in (\huaV_2)_{1d_1g},\\
    \widecheck{s}_0^0(\epsilon_{1p})&= (\epsilon_{1p}, \epsilon_{1p}).
    \end{split}
\end{equation}
The face maps between levels 1 and 2 are 
\begin{equation}\label{eq:VB2DualFaceMaps2}
    \begin{aligned}
    &\widecheck{d}_0^2(\phi^{012}_t, \phi^{112}_{s_0d_0t}, \phi^{122}_{s_1d_0t}, \phi^{001}_{s_0d_2t}) = (\phi^{112}_{s_0d_0t}, \phi^{122}_{s_1d_0t}),\\
    &\widecheck{d}_1^2(\phi^{012}_t, \phi^{112}_{s_0d_0t}, \phi^{122}_{s_1d_0t}, \phi^{001}_{s_0d_2t}) = (\phi^{002}_{s_0d_1t}, \phi^{022}_{s_1d_1t}),\\
    &\widecheck{d}_2^2(\phi^{012}_t, \phi^{112}_{s_0d_0t}, \phi^{122}_{s_1d_0t}, \phi^{001}_{s_0d_2t}) = (\phi^{001}_{s_0d_2t}, \phi^{011}_{s_1d_2t}),\\
    \end{aligned}
\end{equation}
with the following definitions:
\begin{equation}\label{eq:VB2DualFaceMapsPhi002Def}
    \phi^{002}_{s_0d_1t}(Y) = \phi^{002}_{s_0d_1t}(Y^0) = \phi^{001}_{s_0d_2t}(Z^0) + \phi^{012}_t(0_{t}\square Y^0 Z^0)
\end{equation}
for any $Y \in \huaV_2$, $Y^0 := Y - \widetilde{s}_0\widetilde{d}_0Y \in \ker\widetilde{d}^2_0|_{s_0d_1t}$, and $Z^0$ a (2,1)-horn filler of $(0_{d_2t}, (\widetilde{d}_2Y^0)_{1d_1d_2t})$ over $s_0d_2t$, that is: $\widetilde{d}_0Z^0 = 0_{d_2t}$, $\widetilde{d}_2Z^0 = \widetilde{d}_2Y^0$ (which is over $1d_1d_2t$), and $\widetilde{d}_1Z^0$ is over $d_2t$;
\begin{equation}\label{eq:VB2DualFaceMapsPhi022Def}
    \phi^{022}_{s_1d_1t}(X) = \phi^{022}_{s_1d_1t}(X^2) = \phi^{122}_{s_1d_0t}(W^2) + \phi^{012}_t(W^2 X^2 \square 0_t)
\end{equation}
for any $X \in \huaV_2$, $X^2 := X -\widetilde{s}_1\widetilde{d}_2X \in \ker\widetilde{d}^2_2|_{s_1d_1t}$, and $W^2$ a (2,1)-horn filler of $((\widetilde{d}_0X^2)_{1d_0d_0t}, 0_{d_0t})$ over $s_1d_0t$, that is: $\widetilde{d}_0W^2 = \widetilde{d}_0X^2$ (which is over ${1d_0d_0t}$), $\widetilde{d}_2W^2 = 0_{d_0t}$, and $\widetilde{d}_1W^2$ is over $d_0t$;
\begin{equation}\label{eq:VB2DualFaceMapsPhi011Def}
    \phi^{011}_{s_1d_2t}(Z) = \phi^{011}_{s_1d_2t}(Z^2) = \phi^{112}_{s_0d_0t}(W^1) - \phi^{012}_t(W^1 \square 0_t Z^2)
\end{equation}
for any $Z\in \huaV_2$, $Z^2 :=  Z -\widetilde{s}_1\widetilde{d}_2Z \in \ker\widetilde{d}^2_2|_{s_1d_2t}$, and $W^1$ a (2,0)-horn filler of $(0_{d_0t},(\widetilde{d}_0Z^2)_{1d_1d_0t})$ over $s_0d_0t$, that is: $\widetilde{d}_0W^1$ is over $d_0t$, $\widetilde{d}_1W^1 = 0_{d_0t}$, and $\widetilde{d}_2W^1 = \widetilde{d}_0Z^2$, (which is over ${1d_1d_0t}$).
The degeneracy maps between levels 1 and 2 are
\begin{equation}\label{eq:VB2DualDegeneracyMaps2}
    \begin{aligned}
    \widecheck{s}_0(\eta_{s_0g},\xi_{s_1g}) &= (\eta_{s_0g}, (\eta_{s_0g}, \xi_{s_1g}), \widecheck{d}_1(\eta_{s_0g},\xi_{s_1g})),\\
    \widecheck{s}_1(\eta_{s_0g},\xi_{s_1g}) &= (\xi_{s_1g}, \widecheck{s}_0\widecheck{d}_0(\eta_{s_0g}, \xi_{s_1g}), \eta_{s_0g}).
    \end{aligned}
\end{equation}
The 2-groupoid structure of $\huaV^{2*}$ is given by the multiplications:
\begin{equation}\label{eq:VB2DualMultiplicationsDef}
    \begin{aligned}
        &\square \chi_s \psi_t \omega_u = ((\square \chi_s \psi_t \omega_u)^{012}, (\chi^{112}, \chi^{122})_{d_0s}, \omega^{112}_{s_0d_0u}) \in (\huaV^{2*}_2)_r,\\
        &\phi_r \square \psi_t \omega_u = ((\phi_r \square \psi_t \omega_u)^{012}, (\phi^{112}, \phi^{122})_{d_0r}, \omega^{002}_{s_0d_1u}) \in (\huaV^{2*}_2)_s,\\
        &\phi_r \chi_s \square \omega_u = ((\phi_r \chi_s \square \omega_u)^{012}, (\phi^{002}, \phi^{022})_{d_1r}, \omega^{001}_{s_0d_2u}) \in (\huaV^{2*}_2)_t,\\
        &\phi_r \chi_s \psi_t \square = ((\phi_r \chi_s \psi_t \square)^{012}, (\phi^{001}, \phi^{011})_{d_2r}, \psi^{001}_{s_0d_2t}) \in (\huaV^{2*}_2)_u,
    \end{aligned}
\end{equation}
for any $(\phi_r, \chi_s, \psi_t, \omega_u)$ that form a tetrahedron in $\huaV^{2*}_3|_{(r,s,t,u)} \cong \Lambda^3_j(\huaV)|_{p^3_j(r,s,t,u)}  \subseteq (\huaV^{2*}_2)^4$, with the $012$ component of each multiplication given by
\begin{equation}\label{eq:VB2DualMultiplicationsDef012}
    \begin{aligned}
        &(\square \chi_s \psi_t \omega_u)^{012}(\square X_s Y_t Z_u) 
        = \chi^{012}_s (X_s) -  \psi^{012}_t (Y_t) +  \omega^{012}_u (Z_u)\\
        & (\phi_r \square \psi_t \omega_u)^{012} (W_r \square Y_t Z_u) 
        =  \phi^{012}_r (W_r) +  \psi^{012}_t (Y_t) -  \omega^{012}_u (Z_u)\\
        & (\phi_r \chi_s \square \omega_u)^{012} (W_r X_s \square Z_u) 
        = -  \phi^{012}_r (W_r) +  \chi^{012}_s (X_s) +  \omega^{012}_u (Z_u)\\
        & (\phi_r \chi_s \psi_t \square)^{012} (W_r X_s Y_t \square) 
        =  \phi^{012}_r (W_r) -  \chi^{012}_s (X_s) +  \psi^{012}_t (Y_t),\\
    \end{aligned}
\end{equation}
for any $(W_r, X_s, Y_t, Z_u)$ that form a tetrahedron in $\huaV_3|_{(r,s,t,u)} \cong \Lambda^3_j(\huaV)|_{p^3_j(r,s,t,u)}  \subseteq (\huaV_2)^4$.
\end{definition}

\begin{remark}
Both $\huaV^{2*}_1$ and $\huaV^{2*}_2$ are well-defined vector bundles over their respective spaces, as they are fiber products of vector bundles over fiberwise linear surjective submersions. 
\end{remark}

\begin{remark}
    In \eqref{eq:VB2DualElementsLv2Def} we use the identity
    \begin{equation*}
        \zeta\nu_0(k^0_t) = - 0_{s_1d_0t} \square (\square 0_t k^0_t 0_{s_1d_2t}) 0_{s_0d_0t} = - \square 0_{s_1d_1t} k^0_t 0_t
    \end{equation*}
    by the associativity relation given by symmetry of the simplicial matrix 
    \begin{equation*}
        \begin{pmatrix}
            0_{s_1d_0t} &0_{s_1d_0t} \square (\square 0_t k^0_t 0_{s_1d_2t}) 0_{s_0d_0t} &\square 0_t k^0_t 0_{s_1d_2t} &0_{s_0d_0t}\\
            0_{s_1d_0t} &0_{s_1d_1t} &0_t &0_t\\
            \square 0_{s_1d_1t} k^0_t 0_t &0_{s_1d_1t} &k^0_t &0_{t}\\
            \square 0_t k^0_t 0_{s_1d_2t} &0_t &k^0_t &0_{s_1d_2t}\\
            0_{s_0d_0t} &0_t &0_t &0_{s_1d_2t}
        \end{pmatrix}
    \end{equation*}
    which can be shown to be an element of $\huaV_4 = \Cosk^3(\huaV)_4$.
\end{remark}

\begin{remark}
The minus signs in \eqref{eq:VB2DualLevel2FiberProdIso0} and \eqref{eq:VB2DualLevel2FiberProdIso2} come from the fact that if $t=1p$ for $p\in \huaG_0$, then, by the canonical form of the multiplication over the units, $\square 0_{1p} k_{1p}^0 0_{1p} = - k_{1p}^0 + \widetilde{s}_0\widetilde{d}_0k_{1p}^0 = -\ggamma_1k_{1p}^0$, and $0_{1p}k_{1p}^2 0_{1p} \square = - k_{1p}^2 + \widetilde{s}_1\widetilde{d}_2k_{1p}^2 = -\ggamma_1k_{1p}^2$. So adding a minus to the definition gives $\ggamma_1$ over the units as in the vector space version. This is thus consistent with Lemma \ref{lem:VBndual-at-units}.
\end{remark}

\begin{remark}
    Recall that in the VS case we were able to express the conditions corresponding to \eqref{eq:VB2DualElementsLv1Def} and \eqref{eq:VB2DualElementsLv2Def} in terms of any 2-simplex in $\huaV$, as in \eqref{eq:VS2DualElementsLv12Def}. 
    In this case, a generic 2-simplex $X_t \in (\huaV_2)_t$ does not canonically project down to $\ker p^2_j$ as before, since the degeneracies of its faces do not live in the same fiber of $\huaV_2$ over $t$ and thus cannot be subtracted to it. 
    This, combined with the fact that only certain translations over degeneracies are allowed a priori, forces us to make some choices of specific kernels to make fiber products over, while in the VS case we could do everything only with annihilators. 
\end{remark}

Because the definitions of each component $\phi^{002}, \phi^{022}$, and $\phi^{011}$ through \eqref{eq:VB2DualFaceMapsPhi002Def}, \eqref{eq:VB2DualFaceMapsPhi022Def}, and \eqref{eq:VB2DualFaceMapsPhi011Def},  include choices, we need to check that they are well-posed and independent of choices. We show this in the following lemma, which includes some convenient identities as well. 

\begin{lemma}\label{lem:VB2DualFaceMapsPhiDefIntoMultEq}
    For any $\phi \in \huaV^{2*}_2$, $\phi^{002}, \phi^{022}$ and $\phi^{011}$ are well-defined through \eqref{eq:VB2DualFaceMapsPhi002Def}, \eqref{eq:VB2DualFaceMapsPhi022Def}, and \eqref{eq:VB2DualFaceMapsPhi011Def}, respectively. 
    Additionally, the pairs $(\phi^{002}, \phi^{022})$ and $(\phi^{001}, \phi^{011})$ are well-defined elements of $\huaV^{2*}_1$.
    
    Furthermore, for any $\phi \in \huaV^{2*}_2$, \eqref{eq:VB2DualFaceMapsPhi002Def} is equivalent to
    \begin{equation}\label{eq:VB2DualFaceMapsPhiMultCond0012}
        \phi^{012}_t(W) - \phi^{012}_t(X) + \phi^{002}_{s_0d_1t}(Y) - \phi^{001}_{s_0d_2t}(Z) = 0,
    \end{equation}
    for any $(W,X,Y,Z) \in (\huaV_3)_{s_0t}$;   
    \eqref{eq:VB2DualFaceMapsPhi022Def} is equivalent to 
    \begin{equation}\label{eq:VB2DualFaceMapsPhiMultCond0112}
        \phi^{112}_{s_0d_0t}(W) - \phi^{012}_t(X) + \phi^{012}_t(Y) - \phi^{011}_{s_1d_2t}(Z) = 0,
    \end{equation}
    for any $(W,X,Y,Z) \in (\huaV_3)_{s_1t}$;
    and \eqref{eq:VB2DualFaceMapsPhi011Def} is equivalent to 
    \begin{equation}\label{eq:VB2DualFaceMapsPhiMultCond0122}
        \phi^{122}_{s_1d_0t}(W) - \phi^{022}_{s_1d_1t}(X) + \phi^{012}_t(Y) - \phi^{012}_t(Z) = 0,
    \end{equation}
    for any $(W,X,Y,Z) \in (\huaV_3)_{s_2t}$.
\end{lemma}

\begin{proof}
    The proof of well-definedness is very similar to that for vector spaces in Theorem \ref{thm:VS2Dual}, namely we have to show that \eqref{eq:VB2DualFaceMapsPhi002Def}, \eqref{eq:VB2DualFaceMapsPhi022Def} and \eqref{eq:VB2DualFaceMapsPhi011Def} are independent of the choice of filler.
    Starting with $\phi^{002}_{s_0d_1t}$, this depends on a (2,1)-horn filler of $(0_{d_2t}, (d_2Y^0)_{1d_1d_2t})$. The space of such fillers is an affine space modelled on $k^1\in (\ker\widetilde{p}^2_1)_{s_0d_1t}$. In fact, after fixing such a filler $\widebar{Z}$, which exists by $\Kan(2,1)$, any other filler is given by $\widebar{Z} + k^1$, with $k^1\in (\ker\widetilde{p}^2_1)_{s_0d_2t}$. Then 
    \begin{equation*}
        \begin{split}
            \phi^{001}_{s_0d_2t}(\widebar{Z}) + \phi^{012}_t(0_{t}\square Y^0 \widebar{Z}) &= \phi^{001}_{s_0d_2t}(\widebar{Z} + k^1) + \phi^{012}_t(0_{t}\square Y^0 (\widebar{Z} + k^1))\\
            \phi^{001}_{s_0d_2t}(k^1) + \phi^{012}_t(0_{t}\square 0_{s_0d_1t} (k^1)) &= 0 \\
            \phi^{012}_t(0_{t}\square 0_{s_0d_1t} (k^1)) &= -\phi^{001}_{s_0d_2t}(k^1),
        \end{split}
    \end{equation*}
    by linearity of multiplication. Clearly, in the last line we have the inverse of \eqref{eq:VB2DualLevel2FiberProdIso2}, so this holds if and only if 
    \begin{equation*}
        \phi^{012}_t(k^2) = - \phi^{001}_{s_0d_0t}(0_t k^2 0_{s_0d_1t}\square),
    \end{equation*}
    for any $k^2 \in (\ker\widetilde{p}^2_2)_t$, which is one of the defining conditions of $\huaV_2^{2*}$ in \eqref{eq:VB2DualElementsLv2Def}.
    For $\phi^{022}$ the proof is analogous. Fix a filler $\widebar{W}$ of the (2,1)-horn $(\widetilde{d}_0X^2, 0_{d_0t})$ over $s_1d_0t$. Any other filler is given by $\widebar{W} + k^1$ for some $k^1 \in (\ker \widetilde{p}^2_1)_{s_1d_0t}$. Then $\phi^{022}$ is well-defined if and only if 
    \begin{equation*}
            \phi^{122}_{s_1d_0t}(k^1) = - \phi^{012}_t(k^1 0_{s_1d_2t} \square 0_t)
    \end{equation*}
    for any $k^1 \in (\ker \widetilde{p}^2_1)_{s_1d_0t}$, which is equivalent to 
    \begin{equation*}
            - \phi^{122}_{s_1d_0t}(\square 0_{s_1d_1t} k^0 0_t) = \phi^{012}_t(k^0)
    \end{equation*}
    for any $k^0 \in (\ker \widetilde{p}^2_0)_{t}$. This is again a condition in \eqref{eq:VB2DualElementsLv2Def}.
    Finally, by the same argument, $\phi^{011}$ is well-defined if and only if for any $k^0 \in (\ker \widetilde{p}^2_0)_{s_0d_0t}$, 
    \begin{equation*}
        \phi^{112}_{s_0d_0t}(k^0) = \phi^{012}_t(k^0 \square 0_t 0_{s_1d_2t}).
    \end{equation*}
    Now, for any $k^0 \in (\ker \widetilde{p}^2_0)_{s_0d_0t}$, we have that
    \begin{equation*}
            \square k^0 k^0 0_{1d_1d_0t} - \square k^0 0_{s_0d_0t} 0_{1d_1d_0t}= \square 0_{s_0d_0t} k^0 0_{1d_1d_0t} = \widetilde{s}_0\widetilde{d}_0k^0 - k^0,
    \end{equation*}
    because $\huaV$ is a $\VB$ 2-groupoid, so by Example \ref{ex:finite-data-2gpds}, its multiplications are compatible with the degeneracy maps and they are fiberwise linear over $\huaG$.
    Because of this, the simplicial matrix
    \begin{equation*}
        \begin{pmatrix}
            (\widetilde{s}_0\widetilde{d}_0k^0 - k^0) &0_{s_0d_0t} &k^0 &0_{1d_1d_0t}\\
            (\widetilde{s}_0\widetilde{d}_0k^0 - k^0) &0_{t} &((\widetilde{s}_0\widetilde{d}_0k^0 - k^0)0_{t} \square 0_{s_1d_2t}) &0_{1d_1d_0t}\\
            0_{s_0d_0t} &0_t &0_t &0_{s_1d_2t}\\
            k^0 &(k^0 \square 0_t 0_{s_1d_2t}) &0_t &0_{s_1d_2t}\\
            0_{1d_0d_2t} &0_{s_1d_2t} &0_{s_1d_2t} &0_{s_1d_2t}\\
        \end{pmatrix}
    \end{equation*}
    is an element of $\huaV_4$ over $s_2s_1t$, which implies that 
    \begin{equation*}
        k^0 \square 0_t 0_{s_1d_2t} = (\widetilde{s}_0\widetilde{d}_0k^0 - k^0)0_{t} \square 0_{s_1d_2t},
    \end{equation*}
    for any $k^0 \in (\ker \widetilde{p}^2_0)_{s_0d_0t}$. 
    Note now that $k^0 - \widetilde{s}_0\widetilde{d}_0k^0 \in \ker (\widetilde{p}^2_1)_{s_0d_0 t}$, so, by normalization of $\phi^{112}$ and \eqref{eq:VB2DualElementsLv2Def},
    \begin{equation*}
    \begin{split}
        \phi^{112}_{s_0d_0t}(k^0) = \phi^{112}_{s_0d_0t}(k^0 - \widetilde{s}_0\widetilde{d}_0k^0) &= - \phi^{012}_t((k^0 - \widetilde{s}_0\widetilde{d}_0k^0)0_{t} \square 0_{s_1d_2t})\\ 
        &= \phi^{012}_t(k^0 \square 0_t 0_{s_1d_2t}).
    \end{split}
    \end{equation*}
    Therefore, $\phi^{011}$ is well-defined.

    We now show that for any $\phi\in \huaV^{2*}|_t$, the pairs $(\phi^{002}, \phi^{022})$ and $(\phi^{001}, \phi^{011})$ satisfy \eqref{eq:VB2DualElementsLv1Def}. For the first pair, we need to show that for any $k\in \ker\widetilde{p}_1^2|_{s_0d_1t}$,
    \begin{equation*}
        \phi^{002}_{s_0d_1t}(k)=\phi^{022}_{s_1d_1t}(0_{s_1d_1t} \square k 0_{s_0d_1t}).
    \end{equation*}
    By inserting the definitions of $\phi^{002}$ and $\phi^{022}$ in \eqref{eq:VB2DualFaceMapsPhi002Def} and \eqref{eq:VB2DualFaceMapsPhi022Def} with fillers set to the respective zeroes, since both $\widetilde{d}_2k$ and $\widetilde{d}_0k$ vanish, this is equivalent to 
    \begin{equation*}
        \phi^{012}_{t}(0_t\square k 0_{s_0d_2t})=\phi^{012}_{t}(0_{s_1d_0t}(0_{s_1d_1t} \square k 0_{s_0d_1t})\square 0_t)
    \end{equation*}
    which follows from the associativity condition encoded in the symmetry of the simplicial matrix 
    \begin{equation*}
    \begin{pmatrix}
        0_{s_1d_0t}
        &0_{s_1d_1t}
        &0_{t}
        &0_{t}\\ 
        0_{s_1d_0t}
        &(0_{s_1d_1t} \square k 0_{s_0d_1t}) 
        &(0_{s_1d_0t}(0_{s_1d_1t} \square k 0_{s_0d_1t})\square 0_t)
        &0_t\\ 
        0_{s_1d_1t} &(0_{s_1d_1t} \square k 0_{s_0d_1t}) &k &0_{s_0d_1t}\\
        0_t &(0_t\square k 0_{s_0d_2t}) &k &0_{s_0d_2t}\\
        0_{t} &0_{t} &0_{s_0d_1t} &0_{s_0d_2t}\\
    \end{pmatrix}
    \end{equation*}
    over $s_3s_0t$. For the second pair we show equivalently that
    \begin{equation*}
        \phi^{001}_{s_0d_2t}(0_{s_1d_2t}k \square 0_{s_0d_2t}) = \phi^{011}_{s_1d_2t}(k),
    \end{equation*} 
    for any $k \in \ker \widetilde{p}^2_1|_{s_1d_2t}$. 
    By inserting \eqref{eq:VB2DualFaceMapsPhi011Def} with the zero filler, this is equivalent to 
    \begin{equation*}
        \phi^{001}_{s_0d_2t}(0_{s_1d_2t}k \square 0_{s_0d_2t}) = - \phi^{012}_{t}(0_{s_0d_0t}\square 0_t k).
    \end{equation*}
    We define $k':= 0_{s_0d_0t}\square 0_t k \in \ker \widetilde{p}^2_2|_{s_1d_2t}$, and by \eqref{eq:VB2DualElementsLv2Def}, the above is equivalent to 
    \begin{equation*}
        \phi^{001}_{s_0d_2t}(0_{s_1d_2t} (0_{s_0d_0t} k' 0_t \square) \square 0_{s_0d_2t}) = \phi^{001}_{t}(0_t k' 0_{s_0d_1t} \square),
    \end{equation*}
    which follows from the associativity condition encoded in the symmetry of the simplicial matrix 
    \begin{equation*}
    \begin{pmatrix}
        0_{s_0d_0t}
        &0_{t}
        &0_{t}
        &0_{s_1d_2t}\\ 
        0_{s_0d_0t}
        &k' 
        &0_t
        &(0_{s_0d_0t} k' 0_t \square)\\ 
        0_t &k' &0_{s_0d_1t} &(0_t k' 0_{s_0d_1t} \square)\\
        0_t &0_t &0_{s_0d_1t} &0_{s_0d_2t}\\
        0_{s_1d_2t} &(0_{s_0d_0t}k'0_t \square) &(0_{s_1d_2t} (0_{s_0d_0t} k' 0_t \square) \square 0_{s_0d_2t}) &0_{s_0d_2t}\\
    \end{pmatrix}
    \end{equation*}
    over $s_2s_0t$. 

    We now show that \eqref{eq:VB2DualFaceMapsPhi002Def} is equivalent to \eqref{eq:VB2DualFaceMapsPhiMultCond0012}. The same argument can then be used to show that \eqref{eq:VB2DualFaceMapsPhi011Def} is equivalent to \eqref{eq:VB2DualFaceMapsPhiMultCond0112} and that \eqref{eq:VB2DualFaceMapsPhi022Def} is equivalent to \eqref{eq:VB2DualFaceMapsPhiMultCond0122}. Starting with \eqref{eq:VB2DualFaceMapsPhi002Def}, we have that for any $\phi\in (\huaV^{2*})_t$,
    \begin{equation*}
        \phi^{002}_{s_0d_1t}(Y) = \phi^{001}_{s_0d_2t}(Z^0) + \phi^{012}_t(0_{t}\square Y^0 Z^0)
    \end{equation*}
    for any $Y \in \huaV_2$, $Y^0 := Y - \widetilde{s}_0\widetilde{d}_0Y \in (\ker\widetilde{d}_0^2)_{s_0d_1t}$, and any $Z^0$ a (2,1)-horn filler of $(0_{d_2t}, (\widetilde{d}_2Y^0)_{1d_1d_2t})$ over $s_0d_2t$. 
    Consider now any $W \in (\huaV_2)_t$ such that $\widetilde{d}_1W = \widetilde{d}_0Y$, and define $Z:= Z^0 + \widetilde{s}_0\widetilde{d}_2W$. Clearly $\widetilde{d}_0Z = \widetilde{d}_2W$, but also 
    \begin{equation*}
        \widetilde{d}_2Z = \widetilde{d}_2Z^0 + 1\widetilde{d}_1\widetilde{d}_2W = \widetilde{d}_2Y^0 + 1\widetilde{d}_1\widetilde{d}_0Y = \widetilde{d}_2Y.
    \end{equation*}
    So the data of $(W, Y, Z^0)$ as above is equivalent to the data of $(W,Y,Z)\in \Lambda^3_1(\huaV)$. Moreover $\phi^{001}_{s_0d_2t}(Z^0) = \phi^{001}_{s_0d_2t}(Z)$, by normalization. Hence \eqref{eq:VB2DualFaceMapsPhi002Def} is equivalent to
    \begin{equation*}
        \begin{split}
            \phi^{002}_{s_0d_1t}(Y) 
            &= \phi^{001}_{s_0d_2t}(Z) + \phi^{012}_t(W_{t}\square Y_{s_0d_1t} Z_{s_0d_2t} - W_t \square (\widetilde{s}_0\widetilde{d}_0Y)(\widetilde{s}_0\widetilde{d}_0Z))\\
            &= \phi^{001}_{s_0d_2t}(Z) + \phi^{012}_t(W_{t}\square Y_{s_0d_1t} Z_{s_0d_2t}) - \phi^{012}_t(W_t \square (\widetilde{s}_0\widetilde{d}_1W)(\widetilde{s}_0\widetilde{d}_2W))\\
            &= \phi^{001}_{s_0d_2t}(Z) + \phi^{012}_t(W_{t}\square Y_{s_0d_1t} Z_{s_0d_2t}) - \phi^{012}_t(W_t),
        \end{split}
    \end{equation*}
    for any $(W,Y,Z)\in \Lambda^3_1(\huaV)$, which is then exactly \eqref{eq:VB2DualFaceMapsPhiMultCond0012}. 
\end{proof}

The following two results show that the data above is actually the 2-dual of the $\VB$ 2-groupoid $\huaV$.

\begin{theorem}\label{thm:VB2dual-is-2-gpd}
    The finite data in Definition \ref{def:VB2dual} defines a $\VB$ 2-groupoid. 
\end{theorem}

As the proof of this theorem is a series of lengthy computations to check that the data in Definition \ref{def:VB2dual} satisfies the conditions of Example \ref{ex:finite-data-2gpds}, we leave it to Section \ref{sec:VB2dual-2-gpd-proof}. 
We now proceed to showing that $\huaV^{2*}$ as defined above  satisfies the 2-dual equations from Lemma \ref{lem:n-dual-equations}. 

\begin{theorem}\label{thm:V2star-is-VB-2-dual}
    Let $\huaV \to \huaG$ be a $\VB$ 2-groupoid. Then $\huaV^{2*}$ is the 2-dual of $\huaV$. 
\end{theorem}

\begin{proof}
    By Lemma \ref{lem:n-dual-finite-equations}, and Theorem \ref{thm:VB2dual-is-2-gpd}, we only need to show that $\huaV^{2*}$ is a simplicial subbundle of $\mathfrak{A}_2\huaV$ and that it satisfies the 2-dual equations from Lemma \ref{lem:n-dual-equations} up to level 3.
    To begin with, there is an embedding $\mathfrak{a}: \huaV^{2*} \into \mathfrak{A}_2\huaV$ which is defined by matching the obvious components at level 0 and 2, and at level 1 by setting 
    \begin{equation}\label{eq:VB2dual-into-ambient-lv1}
        (\eta_{s_0g}, \xi_{s_1g}) \mapsto ((\widecheck{d}_1(\eta, \xi))^{000}, (\eta)^{001}, (\xi)^{011}, (\widecheck{d}_0(\eta, \xi))^{111}),
    \end{equation}
    where level 1 of $\mathfrak{A}_1\huaV$ is described in Example \ref{ex:VB2dual-ambient}.

    For level 0, by Lemma \ref{lem:VBndual-at-units} and Theorem \ref{thm:VS2Dual}, we get, for any $p \in \huaG_0$,
    \begin{equation*}
        \SVect(\textstyle\int p^*\huaV, B^2\R) \cong ((1^*\huaV)|_p)^{2*} = \Ann(D_2\huaV)|_p = O_{01}|_p \subseteq \huaV_2^*|_{1p},
    \end{equation*}
    which is exactly $\huaV^{2*}$ at $p$.

    For level 1 and $g\in \huaG_1$, we get that any $\eta \in \SVect(\textstyle\int g^*\huaV, B^2\R)$ has four components: $\eta^{000}, \eta^{001}, \eta^{011}, \eta^{111}$. 
    These are normalized according to \eqref{eq:VBnDualNormalizationMostGeneral} so that
    \begin{equation*}
        \begin{array}{ll}
            \eta^{000} \in O_{01}|_{d_1g}, 
            &\eta^{111} \in O_{01}|_{d_0g}, \\
            \eta^{001} \in O_0|_g, 
            &\eta^{011} \in O_1|_g.
        \end{array}
    \end{equation*}
    Which are the same normalization conditions expected from \eqref{eq:VB2DualElementsLv1Def} and \eqref{eq:VB2dual-into-ambient-lv1}. 
    There are multiplicativity equations \eqref{eq:VBnDualMultiplicativityMostGeneral} for each $r=0000$, 0001, 0011, 0111, 1111. 
    Those for $0000$ and $1111$ are equations over the total units so they are automatically satisfied the usual argument as in the $\VS$ case. 

    For $0001$, we have
    \begin{equation}\label{eq:VB2DualMult0001-1}
        \eta^{001}|_{s_0g}(W) - \eta^{001}|_{s_0g}(X) + \eta^{001}|_{s_0g}(Y) - \eta^{000}|_{1d_1g}(Z) = 0,
    \end{equation}
    for any $(W,X,Y,Z) \in \huaV_3|_{s_1s_0g}$.
    By normalization of $\eta^{000}$, for any $Z\in\huaV_2$, we have $\eta^{000}(Z) = \eta^{000}(\ggamma_2 Z)$, with $\ggamma_2 Z \in \ker\widetilde{p}^2_2|_{1d_1g}$.
    Consider now the translation over $s_1s_0g$ given by 
    \begin{equation*}
        \begin{split}
        d_1^*1^*\ker \widetilde{p}^2_2 &\to s_0^*\ker\widetilde{p}^2_2\\
        k|_{1d_1g} &\mapsto 0_{s_0g} 0_{s_0g} \square k|_{1d_1g}.
        \end{split}
    \end{equation*}
    Evaluating \eqref{eq:VB2DualMult0001-1} on the tetrahedron corresponding to the translation of $\ggamma_2Z$ gives
    \begin{equation}\label{eq:VB2DualMult0001-2}
        \eta^{000}|_{1d_1g}(Z) = \eta^{000}|_{1d_1g}(\ggamma_2 Z) = \eta^{001}|_{s_0g}( 0_{s_0g} 0_{s_0g} \square (\ggamma_2 Z)), 
    \end{equation}
    for any $Z\in \huaV_2|_{1d_1g}$. 
    Reinserting this into \eqref{eq:VB2DualMult0001-1} (rewritten with $Y=WX\square Z$) does not yield any extra conditions: 
    \begin{equation*}
        \begin{split}
            \eta^{001}|_{s_0g}(W) - \eta^{001}|_{s_0g}(X) &= -\eta^{001}|_{s_0g}(WX\square Z) + \eta^{001}|_{s_0g}(0_{s_0g} 0_{s_0g} \square (\ggamma_2 Z)), \\
            \eta^{001}|_{s_0g}(W - X) &= - \eta^{001}|_{s_0g}((W - \widetilde{s}_0 \widetilde{d}_0 W + \widetilde{s}_0\widetilde{d}_0W)(W - W + X)\\
            &\qquad\qquad\qquad\square (Z - Z + \widetilde{s}_0 \widetilde{d}_2 W + \widetilde{s}_1 \widetilde{d}_2 X - \widetilde{s}_1 \widetilde{d}_2 W))\\
            &= \eta^{001}|_{s_0g}(\widetilde{s}_0\widetilde{d}_1W - W + X),
        \end{split}
    \end{equation*}
    by linearity of the multiplication of $\huaV$, its compatibility with the degeneracies as in Example \ref{ex:finite-data-2gpds} and the normalization of $\eta^{001}$ in the last step. 
    So we get that \eqref{eq:VB2DualMult0001-1} is equivalent to \eqref{eq:VB2DualMult0001-2}, which, with \eqref{eq:VBndual-ambient-face-maps}, is precisely the definition of $\widecheck{d}_0$ we gave in \eqref{eq:VB2DualFaceMaps1}. 

    For $0111$ the multiplicativity equation \eqref{eq:VBnDualMultiplicativityMostGeneral} is front-to-back symmetric to the previous one, as defined in Remark \ref{rem:front-to-back}, so we apply a front-to-back symmetric argument with the translation 
    \begin{equation*}
        \begin{split}
        d_0^*1^*\ker \widetilde{p}^2_0 &\to s_1^*\ker\widetilde{p}^2_0\\
        k|_{1d_0g} &\mapsto k|_{1d_0g}\square 0_{s_1g} 0_{s_1g}
        \end{split}
    \end{equation*}
    and obtain that \eqref{eq:VBnDualMultiplicativityMostGeneral} for $0111$ is equivalent to the fact that for any $W \in \huaV_2|_{1d_0g}$,
    \begin{equation*}
        \eta^{111}|_{1d_0g}(W) = 
        \eta^{111}|_{1d_0g}(\ggamma_0W) = 
        \eta^{011}|_{s_1g}((\ggamma_0W)\square 0_{s_1g} 0_{s_1g}),
    \end{equation*}
    which is precisely the definition of $\widecheck{d}_1$ in \eqref{eq:VB2DualFaceMaps1}. 

    Finally, for $0011$ we have
    \begin{equation}\label{eq:VB2DualMult0011-1}
        \eta^{011}|_{s_1g}(W) - \eta^{011}|_{s_1g}(X) + \eta^{001}|_{s_0g}(Y) - \eta^{001}|_{s_0g}(Z) = 0, 
    \end{equation}
    for any $(W,X,Y,Z) \in \huaV_3|_{s_2s_0g}$.
    By normalization of $\eta^{001}$ and $\eta^{011}$, this is equivalent to 
    \begin{equation*} 
        \eta^{001}|_{s_0g}(Y - Z - \widetilde{s}_0\widetilde{d}_0(Y-Z)) = \eta^{011}|_{s_1g}(X - W - \widetilde{s}_1\widetilde{d}_2(X - W) ),  
    \end{equation*}
    for any $(W,X,Y,Z) \in \huaV_3|_{s_2s_0g}$.
    Denote $k:=Y - Z - \widetilde{s}_0\widetilde{d}_0(Y-Z) \in (s_0^*\ker \widetilde{p}^2_1)|_g$. 
    By observing that $\widetilde{s}_1\widetilde{d}_2(X - W) = \widetilde{s}_1 \widetilde{d}_1(Z - \widetilde{s}_0\widetilde{d}_0 Z)$ by the simplicial identities, and by using compatibility of $\square$ with degeneracies, we have that\footnote{This whole computation takes place over the tetrahedron $s_0s_1g = s_2s_0g$, so all basepoints in each sum match.}
    \begin{equation*}
    \begin{split}
    X - W - \widetilde{s}_1\widetilde{d}_2(X - W) 
    &= W \square Y Z 
    - W \square (\widetilde{s}_0\widetilde{d}_1W) (\widetilde{s}_0 \widetilde{d}_2W) \\
    &\qquad - (\widetilde{s}_0\widetilde{d}_0(Z - \widetilde{s}_0\widetilde{d}_0Z))\square Z Z = 0_{s_0g} \square k 0_{s_1g}.
    \end{split}
    \end{equation*}
    In the last equality we additionally used linearity of the multiplication and the simplicial identities $\widetilde{d}_1W = \widetilde{d}_0Y$.
    Because $(W,X,Y,Z)$ was an arbitrary tetrahedron we have now shown that \eqref{eq:VB2DualMult0011-1} is equivalent to 
    \begin{equation*}
        \eta^{001}|_{s_0g}(k) = \eta^{011}|_{s_1g}(\zeta k) 
        = \eta^{011}|_{s_1g}(0_{s_1g} \square k 0_{s_0g}), 
    \end{equation*}
    for any $k \in (s_0^*\ker \widetilde{p}^2_1)_g$, which is precisely the fiber product condition for  $\huaV^{2*}_1$ at $g \in \huaG_1$ we gave in \eqref{eq:VB2DualElementsLv1Def}.
    In other words, \eqref{eq:VB2DualMult0011-1} for any tetrahedron $(W,X,Y,Z)$ over $s_2s_0g$ is equivalent to \eqref{eq:VB2DualMult0011-1} holding for the tetrahedron related to the translation $\zeta$ over $s_2s_0g$ defined in \eqref{eq:VB2DualTranslationZetaDef}, for any $k \in \ker\widetilde{p}^2_1|_{s_0g}$.
    In conclusion, we have shown that $\mathfrak{a}:\huaV^{2*}_1 \into \mathfrak{A}_2\huaV$ satisfies the 2-dual equations at level 1. 

    At level 2, for any $t\in \huaG_2$, each element $\phi \in \SVect(\textstyle\int t^*\huaV, B^2\R)$ has 10 components, which are normalized according to \eqref{eq:VBnDualNormalizationMostGeneral} in the following way:
    \begin{equation}\label{eq:VB2DualLevel2Normalization}
        \begin{array}{lll}
            \phi^{222} \in O_{01}|_{d_0d_1t}, 
            &\phi^{112} \in O_0|_{d_0t}, 
            &\phi^{122} \in O_1|_{d_0t}.\\
            \phi^{111} \in O_{01}|_{d_0d_2t},
            &\phi^{002} \in O_0|_{d_1t}, 
            &\phi^{022} \in O_1|_{d_1t}.\\
            \phi^{000} \in O_{01}|_{d_1d_2t}, 
            &\phi^{001} \in O_0|_{d_2t}, 
            &\phi^{011} \in O_1|_{d_2t},\\
            \phi^{012} \in \huaV_2^*|_t. & &
        \end{array}
    \end{equation}
    The latter being the only unnormalized component and the unique interior one. The components of each $\phi\in \huaV^{2*}|_t$ through the embedding $\mathfrak{a}$ satisfy the same normalization conditions. 
    There are 15 multiplicativity equations at this level, divided into three types:
    \begin{itemize}
        \item The ones indexed by $iiii \in \Delta^2_3$, for $0 \le i \le 2$, which are solved in the same way as for level 0 and automatically hold. 
        \item The ones indexed by $iiij$, $iijj$ and $ijjj$ for $0\le i\neq j \le 2$, which relate internally the components of each face of the 2-simplex $\phi$. They can be treated in the same way as for level 1. By the simplicial identities and the fact that, as we show in Lemma \ref{lem:VB2DualFaceMapsPhiDefIntoMultEq}, $(\phi^{002}, \phi^{022})$ and $(\phi^{001}, \phi^{011})$ are elements of $\huaV^{2*}_1$, we have that for any $\phi \in \huaV^{2*}_t$, $\mathfrak{a}(\phi)$ satisfies these equations. 
        \item The ones indexed by $0012$, $0112$ and $0122$ are precisely the equations \eqref{eq:VB2DualFaceMapsPhiMultCond0012}, \eqref{eq:VB2DualFaceMapsPhiMultCond0112} and \eqref{eq:VB2DualFaceMapsPhiMultCond0122}. As we show in Lemma \ref{lem:VB2DualFaceMapsPhiDefIntoMultEq}, these are already equivalent to the equations used to define the components $\phi^{002}, \phi^{022}$, and $\phi^{011}$ of an element in $\huaV^{2*}$. 
    \end{itemize}
    Altogether this means that the last thing to check for $\huaV^{2*}$ to satisfy the 2-dual equations at level 2 is the consistency of the different sets of equations, as we checked it for the $\VS$ 2-dual in the proof of Theorem \ref{thm:VS2Dual}. This is basically taken care of by the simplicial identities and the fact that $(\phi^{002}, \phi^{022})$ and $(\phi^{001}, \phi^{011})$ are elements of $\huaV^{2*}_1$ by Lemma \ref{lem:VB2DualFaceMapsPhiDefIntoMultEq}.

    Finally, we discuss the 2-dual equations at level 3. As we have shown in the proof of Lemma \ref{lem:n-dual-finite-equations}, at level $3$, all components are face components, as they describe the boundary of a 3-simplex. 
    In particular, they can be seen as describing four (overlapping) sets of components normalized as in \eqref{eq:VB2DualLevel2Normalization}, which describe the four faces of the boundary of an empty tetrahedron. 
    This solves the normalization equations \eqref{eq:VBnDualNormalizationMostGeneral}. 
    Most of the multiplicativity equations can also be traced back to the ones for level two, these are the ones corresponding to simplices $iiii, iiij, iijj, ijjj, iijk, ijjk, ijkk$ for $0 \le i\neq j \neq k \le 3$ in $\Delta^3_3$. 
    The only remaining equation relates the four interior components of the faces of the 3-simplex $\phi \in \SVect(\textstyle\int \tau^*\huaV, B^2\R)$ for $\tau\in \huaG_3$ and reads: for any $(W, X, Y, Z) \in \huaV|_\tau$, 
    \begin{equation*}
        \phi^{123}(W) - \phi^{023}(X) + \phi^{013}(Y) - \phi^{012}(Z) = 0.
    \end{equation*}
    This is precisely the equation defining the interior component of the multiplication maps in \eqref{eq:VB2DualMultiplicationsDef012}. Therefore $\huaV^{2*}$ satisfies the 2-dual equations for $\huaV$ at level 3.
\end{proof}

We conclude with a characterization of homologically nondegenerate pairings that follows from Theorem \ref{thm:we-of-VB2gpd-is-qi}. 
This can be applied to show that a 2-shifted presymplectic form on a Lie 2-groupoid is 2-shifted symplectic if it induces a Morita equivalence between the tangent and the 2-cotangent of the Lie 2-groupoid. 
With this we can also show that 2-duality of $\VB$ 2-groupoids is reflexive up to homotopy.

We begin by writing explicitly the maps induced by a 2-shifted simplicial pairing as defined in Remark \ref{rem:VB-n-pairing-ind-maps}. 
Let $\alpha: \huaV_2 \otimes \huaW_2 \to \R$ be a multiplicative normalized 2-shifted pairing between two $\VB$ 2-groupoids over $\huaG$. Then $\alpha$ induces the maps 
\begin{equation*}
    \alpha^r : \huaV \to \huaW^{2*}, \qquad
    \alpha^l : \huaW \to \huaV^{2*},
\end{equation*}
where $\alpha^r$ is defined as
\begin{equation}
    \begin{split}
        \alpha^r_2(X_t)&= (\alpha(X_t, \_), (\alpha(\widetilde{s}_0\widetilde{d}_0X_t, \_), \alpha(\widetilde{s}_1\widetilde{d}_0X_t, \_)), \alpha(\widetilde{s}_0\widetilde{d}_2X_t, \_)),\\
        \alpha^r_1(v_g)&= (\alpha(\widetilde{s}_0v_g, \_), \alpha(\widetilde{s}_1v_g, \_)),\\
        \alpha^r_0(x_p)&= \alpha(\widetilde{1}x_p, \_), 
    \end{split}
\end{equation}
for any $X_t\in \huaV_2|_t$, $v_g\in \huaV_1|_g$ and $x_p\in \huaV_0|_p$.
The other map $\alpha^l$ is defined by inserting elements of $\huaW$ in $\alpha$ in a similar way. 

\begin{corollary}
    Let $\alpha: \huaV_2 \otimes \huaW_2 \to \R \times \huaG_2$ be a multiplicative normalized 2-shifted pairing between two $\VB$ 2-groupoids over $\huaG$. Then the following are equivalent:
    \begin{enumerate}[label=(\roman*)]
        \item $\alpha$ is homologically 2-shifted nondegenerate.
        \item $\alpha^r$ is a weak equivalence.
        \item $\alpha^l$ is a weak equivalence. 
    \end{enumerate}

    In particular, if $\omega_2$ is the leading term of a 2-shifted presymplectic form on a Lie 2-groupoid $\huaG$, it is 2-shifted symplectic if and only if 
    \begin{equation*}
    \omega_2^r = - \omega_2^l: T\huaG \to T^{2*}\huaG   
    \end{equation*}
    is a weak equivalence. 
\end{corollary}

\begin{proof}
    This follows immediately from Theorem \ref{thm:we-of-VB2gpd-is-qi} and commutativity of the diagram \eqref{diag:VBnShiftedPairingsInducedMapsTriangle} from Theorem \ref{thm:VBnShiftedPairingsInducedMapsTriangle}, where the 2-dual pairing is homologically $2$-shifted nondegenerate by Theorem \ref{thm:VBndual-pairing-hndg}. In this situation the induced chain map $\lambda_\alpha^r$ is a quasi-isomorphism if and only if $N(\alpha^r)$ is a quasi-isomorphism, and the latter is equivalent to $\alpha^r$ being a weak equivalence. 
\end{proof}

With this result and the fact that the 2-dual pairing of a $\VB$ 2-groupoid is homologically nondegenerate by Theorem \ref{thm:VBndual-pairing-hndg}, we have the $\VB$ version of reflexivity up to homotopy from Theorem \ref{thm:ndual-reflexive-uth} for $n=2$. 

\begin{corollary}\label{cor:VB-2dual-reflexive-uth}
    The double 2-dual $(\huaV^{2*})^{2*}$ of a $\VB$ 2-groupoid $\huaV \to \huaG$ is weak equivalent to $\huaV$, via the left induced map of the 2-dual pairing. 
\end{corollary}

\subsection{Proof of Theorem \ref{thm:VB2dual-is-2-gpd}}\label{sec:VB2dual-2-gpd-proof}

The following Lemma breaks down the proof of Theorem \ref{thm:VB2dual-is-2-gpd} by recalling the conditions of Example \ref{ex:finite-data-2gpds}.

\begin{lemma}
    The face and degeneracy maps in Definition \ref{def:VB2dual} define a 2-truncated simplicial vector bundle. The multiplications are well-defined maps $\widecheck{m}_j: \Lambda^3_j(\huaV^{2*}) \to \huaV^{2*}_2$, for $j=0,1,2,3$, through \eqref{eq:VB2DualMultiplicationsDef} and \eqref{eq:VB2DualMultiplicationsDef012}. The data above defines a $\VB$ 2-groupoid over $\huaG$, meaning it satisfies the conditions of Example \ref{ex:finite-data-2gpds} in the site of vector bundles with surjective bundle maps over surjective submersions as covers. In summary, $\huaV^{2*}$ is a $\VB$ 2-groupoid.
\end{lemma}

\begin{proof}
    We begin by checking the simplicial identities between levels 0, 1, and 2. 
    This amounts to the fact that the embedding $\mathfrak{a}$ defined in the proof of Theorem \ref{thm:V2star-is-VB-2-dual} is a well-defined simplicial map. We nevertheless include this explicit computation to provide more examples of computations with elements of a $\VB$ 2-groupoid.
    The following simplicial identities can be shown immediately from the definitions:
    \begin{equation*}
        \begin{array}{lll}
            \widecheck{s}_1\widecheck{s}_0 = \widecheck{s}_0\widecheck{s}_0,    \qquad 
            &\widecheck{d}_0^1\widecheck{s}_0^0 = id, \qquad 
            &\widecheck{d}_0^2\widecheck{s}_0^1 = id, \\
            &\widecheck{d}_1^1\widecheck{s}_0^0 = id,
            &\widecheck{d}_0^2\widecheck{s}_1^1 = \widecheck{s}_0^0\widecheck{d}_0^1.
        \end{array}
    \end{equation*}
    The missing ones are
    \begin{equation*}
        \begin{array}{ll}
            \widecheck{d}_1^2\widecheck{s}_0^1 = id,    \qquad\qquad 
            &\widecheck{d}_0\widecheck{d}_1 = \widecheck{d}_0\widecheck{d}_0, \\
            \widecheck{d}_1^2\widecheck{s}_1^1 = id, 
            &\widecheck{d}_0\widecheck{d}_2 = \widecheck{d}_1\widecheck{d}_0, \\
            \widecheck{d}_2^2\widecheck{s}_0^1 = \widecheck{s}_0^0\widecheck{d}_1^1,
            &\widecheck{d}_1\widecheck{d}_2 = \widecheck{d}_1\widecheck{d}_1,  \\
            \widecheck{d}_2^2\widecheck{s}_1^1 = id.
        \end{array}
    \end{equation*}
    \begin{itemize}
    \item For any $(\eta, \xi) \in (\huaV^{2*}_1)_g$, $\widecheck{d}_1\widecheck{s}_0(\eta,\xi)_g = \widecheck{d}_1(\eta, (\eta,\xi), \widecheck{d}_1(\eta,\xi))_{s_0g}$. 
    By \eqref{eq:VB2DualFaceMapsPhi002Def}, the first component of this is
    \begin{equation*}
        \left(\widecheck{s}_0(\eta,\xi)\right)^{002}_{s_0g} (Y) 
        = (\widecheck{d}_1(\eta, \xi))_{1d_1g}(Z^0) + \eta_{s_0g}(0_{s_0g}\square Y^0_{s_0g}Z^0_{1d_1g})
    \end{equation*}
    for any $Y \in (\huaV_2)_{s_0g}$, $Y^0 := Y - \widetilde{s}_0 \widetilde{d}_0 Y$, and any (2,1)-filler $Z^0$ of the horn $(0_{1d_1g}, (\widetilde{d}_2Y^0)_{1d_1g})$ over $1d_1g$. Since this is a filler over a unit, we can use the same one as in solving equation \eqref{eq:VS2DualMult0012with0} in the computation of the $\VS$ case, and take $Z^0 = (\widetilde{s}^1\widetilde{d}_2Y^0)_{1d_1g}$. Hence we get $(\widecheck{d}_1(\eta, \xi))_{1d_1g}(\widetilde{s}^1\widetilde{d}_2Y^0) =0$ by normalization, and
    \begin{equation*}
    \begin{split}
        \left(\widecheck{s}_0(\eta,\xi)\right)^{002}_{s_0g} (Y)
        &=\eta_{s_0g}(0_{s_0g}\square Y^0_{s_0g}(\widecheck{d}_1(\eta, \xi))_{1d_1g}) 
        \\
        &= \eta(Y^0) = \eta(Y),
    \end{split}        
    \end{equation*} 
    since $0_{s_0g}=\widetilde{s}_0\widetilde{d}_0Y^0$, and $(\widetilde{s}_0\widetilde{d}_0Y^0, Y^0, Y^0, \widetilde{s}_1\widetilde{d}_2Y^0)_{s_1s_0g} = \widetilde{s}_1Y^0$ is a degenerate tetrahedron in $\huaV_3$ and $\huaV$ is a 2-groupoid (so its multiplications are compatible with degeneracies). By \eqref{eq:VB2DualFaceMapsPhi022Def}, the second component of $\widecheck{d}_1\widecheck{s}_0(\eta,\xi)_g$ is 
    \begin{equation*}
        \left(\widecheck{s}_0(\eta,\xi)\right)^{022}_{s_1g} (X) 
        = \xi_{s_1g}(W^2) + \eta_{s_0g}(W^2 X^2 \square 0_{s_0g})
    \end{equation*}
    for any $X \in (\huaV_2)_{s_1g}$, $X^2 := X - \widetilde{s}_1\widetilde{d}_2 X$, and $W^2$ a (2,1)-filler of $((\widetilde{d}_0X^2)_{1d_0g}, 0_g)$ over $s_1g$. But $W^2 = X^2$ is such a filler, so, since $\widetilde{s}_0 X^2 = (X^2, X^2, \widetilde{s}_0\widetilde{d}_1X^2, 0_{s_0g})$, we have 
    \begin{equation*}
    \begin{split}
        \left(\widecheck{s}_0(\eta,\xi)\right)^{022}_{s_1g} (X) 
        &= \xi_{s_1g}(X^2) + \eta_{s_0g}(\widetilde{s}_0\widetilde{d}_1X^2)\\ 
        &= \xi_{s_1g}(X^2)  = \xi_{s_1g}(X),
    \end{split}
    \end{equation*}
    by normalization of $\eta$. Therefore $\widecheck{d}_1\widecheck{s}_0(\eta,\xi)_g = (\eta,\xi)_g$.
    
    \item We similarly compute $\widecheck{d}_1\widecheck{s}_1(\eta,\xi)_g = \widecheck{d}_1(\xi, \widecheck{s}_0\widecheck{d}_0(\eta,\xi), \eta)$, for any $(\eta, \xi) \in (\huaV^{2*}_1)_g$. By \eqref{eq:VB2DualFaceMapsPhi002Def}, with the filler $Z^0 = Y^0$ (which is a valid (2,1)-filler of $(0_{g}, (\widetilde{d}_2Y^0)_{1d_1g})$ over $s_0g$), and the same arguments as above, the first component of this is
    \begin{equation*}
    \begin{split}
        \left(\widecheck{s}_1(\eta,\xi)\right)^{002}_{s_0g} (Y) 
        &= \eta_{s_0g}(Y^0) + \xi_{s_1g}(0_{s_1g}\square Y^0_{s_0g}Y^0_{s_0g}) \\
        &= \eta_{s_0g}(Y^0) + \xi_{s_1g}(\widetilde{s}_1\widetilde{d}_1 Y^0) =  \eta_{s_0g}(Y) 
    \end{split}
    \end{equation*}
    for any $Y \in (\huaV_2)_{s_0g}$, and $Y^0 := Y - \widetilde{s}_0 \widetilde{d}_0 Y$.
    Similarly, by \eqref{eq:VB2DualFaceMapsPhi022Def}, with the filler $W^2 = \widetilde{s}_0\widetilde{d}_0X^2$ (a valid (2,1)-filler of $((\widetilde{d}_0X^2)_{1d_0g}, 0_{1d_0g})$ over $1d_0g$), the second component is  
    \begin{equation*}
    \begin{split}
        \left(\widecheck{s}_1(\eta,\xi)\right)^{022}_{s_1g} (X) 
        &= (\widecheck{d}_0(\eta,\xi))_{1d_0g}(\widetilde{s}_0\widetilde{d}_0X^2)\\
        &\qquad \qquad + \xi_{s_1g}((\widetilde{s}_0\widetilde{d}_0X^2) X^2 \square 0_{s_1g})\\
        &= \xi_{s_1g}(X^2) =  \xi_{s_1g}(X), 
    \end{split}
    \end{equation*}
    for any $X \in (\huaV_2)_{s_1g}$, $X^2 := X - \widetilde{s}_1\widetilde{d}_2 X$. Therefore $\widecheck{d}_1\widecheck{s}_1(\eta,\xi)_g = (\eta, \xi)_g$.

    \item For the same arbitrary element, $\widecheck{d}_2\widecheck{s}_1(\eta,\xi)_g = \widecheck{d}_2(\xi, \widecheck{s}_0\widecheck{d}_0(\eta,\xi), \eta)$, has first component $\eta$ and second component given by \eqref{eq:VB2DualFaceMapsPhi011Def}. 
    By choosing the filler $W^1 = (\widetilde{s}_1\widetilde{d}_0 Z^2 - \widetilde{s}_0\widetilde{d}_0 Z^2)_{1d_0g}$ of the (2,0)-horn $(0_{1d_0g}, (\widetilde{d}_0 Z^2)_{1d_0g})$ over $1d_0g$, we get 
    \begin{equation*}
    \begin{split}
        &\left(\widecheck{s}_1(\eta,\xi)\right)^{011}_{s_1g} (Z)\\ 
        &= (\widecheck{d}_0(\eta,\xi))_{1d_0g}(\widetilde{s}_1\widetilde{d}_0 Z^2 - \widetilde{s}_0\widetilde{d}_0 Z^2) - \xi_{s_1g}((\widetilde{s}_1\widetilde{d}_0 Z^2 - \widetilde{s}_0\widetilde{d}_0 Z^2) \square 0_{s_1g} Z^2)\\
        &=  - \xi_{s_1g}((\widetilde{s}_1\widetilde{d}_0 Z^2)\square Z^2 Z^2 - (\widetilde{s}_0\widetilde{d}_0 Z^2) \square Z^2 0_{s_1g})\\
        &= \xi_{s_1g}(Z^2 - \widetilde{s}_1\widetilde{d}_1Z^2) = \xi_{s_1g}(Z),
    \end{split}
    \end{equation*}
    for any $Z \in (\huaV_2)_{s_1g}$ and $Z^2 := Z - \widetilde{s}_1\widetilde{d}_2 Z$, where in the last line we used linearity of the multiplication of $\huaV$. So $\widecheck{d}_2\widecheck{s}_1(\eta,\xi)_g = (\eta,\xi)_g$.

    \item Similarly, $\widecheck{d}_2\widecheck{s}_0(\eta,\xi)_g = \widecheck{d}_2(\eta, (\eta,\xi), \widecheck{d}_1(\eta,\xi))_{s_0g}$, has first component 
    \begin{equation*}
        (\widecheck{d}_1(\eta,\xi))_{1d_1g} (\_)= \eta_{s_0g}(0_{s_0g} 0_{s_0g} \square (\ggamma_2\_)_{1d_1g}).    
    \end{equation*}
    We need to show the second component is the same, by using \eqref{eq:VB2DualFaceMapsPhi011Def}. In this case this gives
    \begin{equation*}
        \left(\widecheck{s}_0(\eta,\xi)\right)^{011}_{1d_1g} (Z) = \eta_{s_0g}(W^1) - (W^1 \square 0_{s_0g} Z^2)
    \end{equation*}
    for any $Z \in (\huaV_2)_{1d_1g}$, with $Z^2 := Z - \widetilde{s}_1\widetilde{d}_2 Z$ and $W^1$ a (2,0)-horn filler of $(0_{g}, (\widetilde{d}_0 Z^2)_{1d_1g})$ over $s_0g$. 
    Note that we can write $W^1 = W^1 \square 0_{s_0g} (\widetilde{s}_0\widetilde{d}_0Z^2)$, since $\widetilde{s}_0W^1 = (W^1, W^1, \widetilde{s}_0\widetilde{d}_1 W^1, \widetilde{s}_0\widetilde{d}_2 W^1)$, $\widetilde{s}_0\widetilde{d}_1 W^1 = 0$, and $\widetilde{d}_2 W^1 = \widetilde{d}_0 Z^2$. So, by linearity of the multiplication,
    \begin{equation*}
        \begin{split}
            \left(\widecheck{s}_0(\eta,\xi)\right)^{011}_{1d_1g} (Z) 
            &= \eta_{s_0g}(W^1 \square 0_{s_0g} (\widetilde{s}_0\widetilde{d}_0Z^2) - W^1 \square 0_{s_0g} Z^2)\\
            &= \eta_{s_0g}(0_{s_0g}\square 0_{s_0g} (- Z^2 + \widetilde{s}_0\widetilde{d}_0 Z^2)).
        \end{split}
    \end{equation*}
    Observe now that $- Z^2 + \widetilde{s}_0\widetilde{d}_0 Z^2 = - \ggamma_1Z$, with $\ggamma_1$ the core projection defined pointwise as in \eqref{eq:VS2DualCoreProjections}. Moreover, by using the simplicial identities, it can readily be shown that $-\widetilde{d}_1\ggamma_1Z = d_2\ggamma_2 Z $ for any $Z \in (\huaV_2)_{1d_1g}$. So, the simplicial matrix 
    \begin{equation*}
        \begin{pmatrix}
            0_{s_0g} &0_{s_0g} &0_{s_0g} &0_{1d_1g}\\
            0_{s_0g} &0_{s_0g} &(0_{s_0g} 0_{s_0g} \square (\ggamma_2 Z)_{1d_1g}) &(\ggamma_2 Z)_{1d_1g}\\
            0_{s_0g} &0_{s_0g} &0_{s_0g} &0_{1d_1g}\\
            0_{s_0g} &(0_{s_0g} \square 0_{s_0g} (-\ggamma_1 Z)_{1d_1g}) &0_{s_0g} &(-\ggamma_1 Z)_{1d_1g}\\
            0_{1d_1g} &(\ggamma_2 Z)_{1d_1g} &0_{1d_1g} &(-\ggamma_1 Z)_{1d_1g}
        \end{pmatrix}
    \end{equation*}
    is a valid element in $(\huaV_4)_{s_0s_0s_0g} \cong (\cosk^3(\huaV)_4)_{s_0s_0s_0g}$, and it encodes the associativity relation 
    \begin{equation*}
        0_{s_0g} \square 0_{s_0g} (-\ggamma_1 Z)_{1d_1g} = 0_{s_0g} 0_{s_0g} \square (\ggamma_2 Z)_{1d_1g}.
    \end{equation*}
    Therefore 
    \begin{equation*}
        \left(\widecheck{s}_0(\eta,\xi)\right)^{011}_{1d_1g} (Z) = \eta_{s_0g}(0_{s_0g} 0_{s_0g} \square (\ggamma_2 Z)_{1d_1g}) = \widecheck{d}_1(\eta,\xi)_g,
    \end{equation*}
    and $\widecheck{d}_2\widecheck{s}_0(\eta,\xi)_g = \widecheck{s}_0\widecheck{d}_1(\eta,\xi)_g$.

    \item Consider an arbitrary $\phi\in (\huaV^{2*}_2)_t$, then, for any $X \in (\huaV_2)_{1d_0d_0t} = (\huaV_2)_{1d_0d_1t}$,
    \begin{equation*}
        \widecheck{d}_0 \widecheck{d}_0 \phi (X)
        = \phi^{122}_{s_1d_0t} ((\ggamma_0 X)_{1d_0d_0t} \square 0_{s_1d_0t} 0_{s_1d_0t}).
    \end{equation*}
    Meanwhile,
    \begin{equation*}
        \begin{split}
        \widecheck{d}_0 \widecheck{d}_1 \phi (X)
        &= \phi^{022}_{s_1d_1t} ((\ggamma_0 X)_{1d_0d_1t} \square 0_{s_1d_1t} 0_{s_1d_1t})\\
        &= \phi^{122}_{s_1d_0t} ((\ggamma_0 X)_{1d_0d_0t} \square 0_{s_1d_0t} 0_{s_1d_0t})\\
        &\quad + \phi^{012}_t(((\ggamma_0 X)_{1d_0d_0t} \square 0_{s_1d_0t} 0_{s_1d_0t})((\ggamma_0 X)_{1d_0d_1t} \square 0_{s_1d_1t} 0_{s_1d_1t}) \square 0_t)\\
        &= \phi^{122}_{s_1d_0t} ((\ggamma_0 X)_{1d_0d_0t} \square 0_{s_1d_0t} 0_{s_1d_0t}),
        \end{split}
    \end{equation*}
    where in the second to last line we used \eqref{eq:VB2DualFaceMapsPhi022Def} evaluated on the 2-simplex $((\ggamma_0 X)_{1d_0d_0t} \square 0_{s_1d_0t} 0_{s_1d_0t})$, with (2,1)-filler of $(\widetilde{d}_0(\ggamma_0 X), 0_{d_0t})$ over $s_1d_0t$ chosen to be $((\ggamma_0 X)_{1d_0d_1t} \square 0_{s_1d_1t} 0_{s_1d_1t})$. In the last line we use the associativity relation
    \begin{equation}\label{eq:VB2DualProofExampleSimpMat}
        ((\ggamma_0 X)_{1d_0d_0t} \square 0_{s_1d_0t} 0_{s_1d_0t})((\ggamma_0 X)_{1d_0d_1t} \square 0_{s_1d_1t} 0_{s_1d_1t}) \square 0_t = 0_t,
    \end{equation} 
    encoded by the simplicial matrix
    \begin{equation*}
        \begin{pmatrix}
            (\ggamma_0X) &(\ggamma_0 X)_{1d_0d_0t} \square 0_{s_1d_0t} 0_{s_1d_0t} &0_{s_1d_0t} &0_{s_1d_0t}\\
            (\ggamma_0X) &(\ggamma_0 X)_{1d_0d_1t} \square 0_{s_1d_1t} 0_{s_1d_1t} &0_{s_1d_1t} &0_{s_1d_1t}\\
            (\ggamma_0 X)_{1d_0d_0t} \square 0_{s_1d_0t} 0_{s_1d_0t} &(\ggamma_0 X)_{1d_0d_1t} \square 0_{s_1d_1t} 0_{s_1d_1t} &(LS) &0_{t}\\
            0_{s_1d_0t} &0_{s_1d_1t} &0_{t} &0_{t}\\
            0_{s_1d_0t} &0_{s_1d_1t} &0_{t} &0_{t}
        \end{pmatrix}
    \end{equation*}
    in $(\huaV_4)_{s_3s_2t} \cong (\cosk^3(\huaV)_4)_{s_3s_2t}$, where $(LS)$ is shorthand for the left side of \eqref{eq:VB2DualProofExampleSimpMat}. Therefore $\widecheck{d}_0 \widecheck{d}_0 = \widecheck{d}_0 \widecheck{d}_1$.

    \item By the evident fact that all formulas involved in the proof above are symmetric under the simplicial front-to-back duality $i \leftrightarrow n-i$ as discussed in Remark \ref{rem:front-to-back}, the proof of $\widecheck{d}_1\widecheck{d}_1 = \widecheck{d}_1\widecheck{d}_2$ goes in the same way as the previous one.
    
    \item Finally, we have that, for any $\phi\in (\huaV^{2*}_2)_t$ and any $X \in (\huaV_2)_{1d_0d_2t} = (\huaV_2)_{1d_1d_0t}$,
    \begin{equation*}
        \widecheck{d}_1 \widecheck{d}_0 \phi (X)
        = \phi^{112}_{s_0d_0t} (0_{s_0d_0t} 0_{s_0d_0t} \square (\ggamma_2 X)_{1d_1d_0t}).
    \end{equation*}
    On the other hand, 
    \begin{equation*}
        \begin{split}
        \widecheck{d}_0 \widecheck{d}_2 \phi (X)
        &= \phi^{011}_{s_1d_2t} ((\ggamma_0 X)_{1d_0d_2t} \square 0_{s_1d_2t} 0_{s_1d_2t})\\
        &= \phi^{112}_{s_0d_0t} (0_{s_0d_0t} 0_{s_0d_0t} \square (\ggamma_2 X)_{1d_1d_0t})\\
        &\quad - \phi^{012}_t((0_{s_0d_0t} 0_{s_0d_0t} \square (\ggamma_2 X)_{1d_1d_0t}) \square 0_t ((\ggamma_0 X)_{1d_0d_2t} \square 0_{s_1d_2t} 0_{s_1d_2t}))\\
        &= \phi^{122}_{s_1d_0t} (0_{s_0d_0t} 0_{s_0d_0t} \square (\ggamma_2 X)_{1d_1d_0t}),
        \end{split}
    \end{equation*}
    where, as before, in the second to last line we used \eqref{eq:VB2DualFaceMapsPhi011Def} evaluated on the triangle $((\ggamma_0 X)_{1d_0d_2t} \square 0_{s_1d_2t} 0_{s_1d_2t})$, with (2,0)-filler of $(0_{d_0t}, (\widetilde{d}_0(\ggamma_0 X)))$ over $s_0d_0t$ chosen to be $(0_{s_0d_0t} 0_{s_0d_0t} \square (\ggamma_2 X)_{1d_1d_0t})$. In the last line we use the fact that
    \begin{equation}\label{eq:VB2DualProofExampleSimpMat2}
        (0_{s_0d_0t} 0_{s_0d_0t} \square (\ggamma_2 X)_{1d_1d_0t}) \square 0_t ((\ggamma_0 X)_{1d_0d_2t} \square 0_{s_1d_2t} 0_{s_1d_2t}) = 0_t,
    \end{equation} 
    which we now show. To begin with, we have the associativity relation 
    \begin{equation}\label{eq:VB2DualProofExampleSimpMat3}
        \ggamma_0 X = (0_{s_0d_0t} 0_{s_0d_0t} \square \ggamma_2 X) 0_{s_0d_0t} 0_{s_0d_0t} \square
    \end{equation}
    encoded by the simplicial matrix
    \begin{equation*}
        \begin{pmatrix}
            0_{s_0d_0t} &0_{s_0d_0t} &(0_{s_0d_0t} 0_{s_0d_0t} \square \ggamma_2 X) &(\ggamma_2 X)_{1d_1d_0t}\\
            0_{s_0d_0t} &0_{s_0d_0t} &0_{s_0d_0t} &0_{1d_1d_0t}\\
            0_{s_0d_0t} &0_{s_0d_0t} &0_{s_0d_0t} &0_{1d_1d_0t}\\
            (0_{s_0d_0t} 0_{s_0d_0t} \square \ggamma_2 X) &0_{s_0d_0t} &0_{s_0d_0t} &(RS)\\
            (\ggamma_2 X)_{1d_1d_0t} &0_{1d_1d_0t} &0_{1d_1d_0t} &(\ggamma_0 X)_{1d_1d_0t}
        \end{pmatrix}
    \end{equation*}
    in $(\huaV_4)_{s_0s_0s_0t}$, with $(RS)$ the right side of \eqref{eq:VB2DualProofExampleSimpMat3}. Note that the last row is a valid tetrahedron in $(\huaV_3)_{1d_1d_0t}$ because 
    \begin{equation*}
        (\ggamma_2 X)_{1d_1d_0t} 0_{1d_1d_0t} 0_{1d_1d_0t} \square = (\ggamma_0 X)_{1d_1d_0t},
    \end{equation*}
    by the canonical form of the multiplication on the VS $2$-groupoid fiber $(\huaV)_{1d_1d_0t}$, being as in \eqref{eq:TriMultOverPointAll}. Now call $(0_{s_0d_0t} 0_{s_0d_0t} \square \ggamma_2 X) = \widebar{X} \in (\huaV_2)_{s_0d_0t}$. Consider the simplicial matrix
    \begin{equation*}
        \begin{pmatrix}
            \widebar{X}_{s_0d_0t} &0_{s_0d_0t} &0_{s_0d_0t} &(\widebar{X} 0_{s_0d_0t} 0_{s_0d_0t} \square)\\ 
            \widebar{X}_{s_0d_0t} &(LS) &0_t &(\blacksquare)\\
            0_{s_0d_0t} &0_{t} &0_{t} &0_{s_1d_2t}\\
            0_{s_0d_0t} &0_{t} &0_{t} &0_{s_1d_2t}\\
            (\widebar{X} 0_{s_0d_0t} 0_{s_0d_0t} \square) &(\blacksquare) &0_{s_1d_2t} &0_{s_1d_2t}
        \end{pmatrix}
    \end{equation*}
    in $(\huaV_4)_{s_2s_1t}$, where $(LS)$ is the left side of \eqref{eq:VB2DualProofExampleSimpMat2}, and 
    \begin{equation*}
        (\blacksquare) = ((\widebar{X} 0_{s_0d_0t} 0_{s_0d_0t} \square) \square 0_{s_1d_2 t} 0_{s_1d_2 t}).
    \end{equation*}
    This encodes the identity $(LS) = 0_t$. Therefore $\widecheck{d}_0 \widecheck{d}_2 = \widecheck{d}_1 \widecheck{d}_0$.
\end{itemize}
    This concludes the proof that the data above defines a 2-truncated simplicial vector bundle over $\huaG$. 

    Moving to the multiplications, the proof of well-definedness contains several steps:
    \begin{enumerate}
        \item Any element of $\huaV_2^*$ can be entirely defined by their evaluation on multiplications of elements in $\huaV_2$ as in \eqref{eq:VB2DualMultiplicationsDef012}.
        \item Each element of $\huaV_2^*$ defined in \eqref{eq:VB2DualMultiplicationsDef012} is independent of the choice of decomposition of its argument into a multiplication. 
        \item Each item in \eqref{eq:VB2DualMultiplicationsDef} is an element of $\huaV^{2*}_2$, i.e. it satisfies the conditions in \eqref{eq:VB2DualElementsLv2Def}.
        \item Each multiplication of elements in \eqref{eq:VB2DualMultiplicationsDef} forms a tetrahedron with its factors (with the expected ordering of faces). 
    \end{enumerate}

    We will show each of these statements for $\widecheck{m}_1$, as the proofs for the other multiplications are entirely analogous, and can be readily reproduced by following the same arguments.

    As for the first step, as in the usual definition of the 1-dual of a $\VB$ 1-groupoid in Definition \ref{def:VB1Dual}, it is enough to define the $\huaV_2^*$ component of $\widecheck{m}_1$ when evaluated on a multiplication $\widetilde{m}_1$ of elements in $\huaV_2$, because any element in $\huaV_2$ can be written as a $1$-multiplication of some 2-simplices over prescribed 2-simplices on the base. 
    More precisely, suppose we need to define $\widecheck{m}_1(\phi_r, \psi_t, \omega_u) = \phi_r \square \psi_t \omega_u$, for any $(\phi_r, \psi_t, \omega_u) \in \Lambda^3_1(\huaV^{2*})_{(r,t,u)}$ with $(r,t,u) \in \Lambda^3_1(\huaG)$ such that $r\square t u = s$. 
    Then any $X_s \in (\huaV_2)_s$ can be written as $X_s = W_r\square Y_t Z_u$, for $(W_r, Y_t, Z_u) \in \Lambda^3_1(\huaV)_{(r,t,u)}$ if and only if the bundle map $(\widetilde{m}_1, m_1)$ is fiberwise surjective. 
    Note that since both $\huaV$ and $\huaG$ are 2-groupoids, $\Lambda^3_1(\huaV) \cong \huaV_3$, $\Lambda^3_1(\huaG) \cong \huaG_3$, and the multiplications correspond to the face maps under this isomorphism. 
    Thus, since $\huaV$ is a $\VB$ 2-groupoid over the Lie 2-groupoid $\huaG$, $(\widetilde{m}_1, m_1)$ is a fiberwise surjective bundle map over a surjective submersion (i.e. a cover in $\VB$). 

    For step 2, we want to prove that the definition of $(\widecheck{m}_1)^{012}$ we give is independent choice of a decomposition of $X_s$. In other words, for any $(3,1)$-horn $(\phi_r, \psi_t, \omega_u) \in \Lambda^3_1(\huaV^{2*})_{(r,t,u)}$, and any $X_s$ which decomposes as 
    \begin{equation*}
        X_s = W_r\square Y_t Z_u = W'_r \square Y'_t Z'_u
    \end{equation*}  
    for some $(W_r, Y_t, Z_u), (W'_r, Y'_t, Z'_u) \in \Lambda^3_1(\huaV)_{(r,t,u)}$, we have
    \begin{equation}\label{eq:VB2DualProof012OfMult1IsWellDef}
        \begin{split}
            &(\phi_r \square \psi_t \omega_u)^{012}(W_r \square Y_t Z_u) = \phi^{012}_r (W) + \psi^{012}_t (Y) - \omega^{012}_u (Z) \\
            &= (\phi_r \square \psi_t \omega_u)^{012}(W'_r \square Y'_t Z'_u) = \phi^{012}_r (W') + \psi^{012}_t (Y') - \omega^{012}_u (Z') \\
            \iff \qquad & \phi^{012}_r (W - W') + \psi^{012}_t (Y - Y') - \omega^{012}_u (Z - Z') = 0.
        \end{split}
    \end{equation}
    Observe now that 
    \begin{equation*}
        \begin{array}{ccc}
            \widetilde{d}_0 W = \widetilde{d}_0 X =\widetilde{d}_0 W',
            & \widetilde{d}_1 Y = \widetilde{d}_1 X =\widetilde{d}_1 Y',
            & \widetilde{d}_1 Z = \widetilde{d}_2 X =\widetilde{d}_1 Z'.
        \end{array}
    \end{equation*}
    So that 
    \begin{equation*}
        \begin{array}{ccc}
            W - W' \in (\ker \widetilde{d}_0)_r,
            & Y - Y' \in (\ker \widetilde{d}_1)_t,
            & Z - Z' \in (\ker \widetilde{d}_1)_u.
        \end{array}
    \end{equation*}
    Choose any (2,0)-filler $U$ of $(\widetilde{d}_2(Y -Y')_{d_2t}, 0_{d_2t}) = (\widetilde{d}_2(Z - Z')_{d_2u}, 0_{d_2u})$ over $s_1d_2t = s_1d_2u$, and write 
    \begin{equation*}
        \begin{array}{cc}
            R^Y := \square (Y - Y')_t 0_t U_{s_1d_2t} \in (\huaV_2)_{s_0d_0t}, 
            & R^Z := \square (Z - Z')_u 0_u U_{s_1d_2u} \in (\huaV_2)_{s_0d_0u}.
        \end{array}
    \end{equation*}
    Now consider the simplicial matrix
    \begin{equation*}
        \begin{pmatrix}
            (W - W')_r &0_{r} &R^Y_{s_0d_1r} &R^Z_{s_0d_2r}\\ 
            (W - W')_r &0_{s} &(Y - Y')_t &(Z - Z')_u\\ 
            0_{r} &0_{s} &0_{t} &0_{u}\\
            R^Y_{s_0d_0t} &(Y - Y')_t &0_{t} &U_{s_1d_2t}\\
            R^Z_{s_0d_0u} &(Z - Z')_u &0_u &U_{s_1d_2u}
        \end{pmatrix}
    \end{equation*}
    in $(\huaV_4)_{s_1(r,s,t,u)}$. It can easily be shown that this is a well-defined simplicial matrix, by checking that all rows but the first one are boundaries of tetrahedra in $\huaV_3$. In particular the second row identifies a tetrahedron in $\huaV_3$ because 
    \begin{equation*}
        W_r\square Y_t Z_u = W'_r \square Y'_t Z'_u \implies (W - W')_r \square (Y - Y')_t (Z - Z')_u = 0,
    \end{equation*}
    by linearity of the multiplication of $\huaV$. 
    Then this implies the first line is also a boundary of a tetrahedron in $\huaV_3$, which gives 
    \begin{equation*}
        W - W' = \square 0_s R^Y_t R^Z_u.
    \end{equation*}
    Hence, by the multiplicativity equations \eqref{eq:VB2DualFaceMapsPhiMultCond0012}, \eqref{eq:VB2DualFaceMapsPhiMultCond0112}, from Lemma \ref{lem:VB2DualFaceMapsPhiDefIntoMultEq}, 
    \begin{equation*}
        \begin{split}
            \phi^{012}_r (W - W') &=  \phi^{012}_r (\square 0_s R^Y_t R^Z_u) 
            = \phi^{001}_{s_0 d_2 r}(R^Z) - \phi^{002}_{s_0d_1r}(R^Y),\\
            \psi^{012}_t (Y - Y') &= \psi^{012}_t (R^Y \square 0_t U_{s_1d_2t}) = \psi^{112}_{s_0 d_0 t}(R^Y) - \psi^{011}_{s_1d_2t}(U), \\
            \omega^{012}_u (Z - Z') &= \omega^{012}_t (R^Z \square 0_u U_{s_1d_2u}) = \omega^{112}_{s_0 d_0 u}(R^Z) - \omega^{011}_{s_1d_2t}(U).
        \end{split}
    \end{equation*}
    (Note the evaluations in the first one are well-defined because $d_2r = d_0u$, $d_1r = d_0 t$, by the simplicial identities on $(r,s,t,u)$.) 
    Finally, by the horn conditions on $(\phi_r, \psi_t, \omega_u)$, 
    \begin{equation*}
        \begin{split}
            \widecheck{d}_2\phi &= (\phi^{001}, \phi^{011}) = (\omega^{112}, \omega^{122}) = \widecheck{d}_0 \omega,\\
            \widecheck{d}_1\phi &= (\phi^{002}, \phi^{022}) = (\psi^{112}, \psi^{122}) = \widecheck{d}_0 \psi,\\
            \widecheck{d}_2\psi &= (\psi^{001}, \psi^{011}) = (\omega^{001}, \omega^{011}) = \widecheck{d}_2 \omega,
        \end{split}
    \end{equation*}
    so that \eqref{eq:VB2DualProof012OfMult1IsWellDef} is satisfied:
    \begin{equation*}
        \begin{split}
            &\phi^{012}_r (W - W') + \psi^{012}_t (Y - Y') - \omega^{012}_u (Z - Z') \\
            &= \phi^{001}_{s_0 d_2 r}(R^Z) - \phi^{002}_{s_0d_1r}(R^Y) + \psi^{112}_{s_0 d_0 t}(R^Y)\\
            &\qquad - \psi^{011}_{s_1d_2t}(U) - \omega^{112}_{s_0 d_0 u}(R^Z) + \omega^{011}_{s_1d_2t}(U) = 0.
        \end{split}
    \end{equation*}

    For step 3, for any $(\phi_r, \psi_t, \omega_u) \in \Lambda^3_1(\huaV^{2*})_{(r,t,u)}$ with $s=r\square tu$, we need to show that 
    \begin{equation*}
        \begin{split}
            &(\phi \square \psi \omega)^{012}_s(k^0_s) = - \phi^{112}_{s_0d_0r}(\square 0_s k^0_s 0_{s_1d_2s})\\
            &\text{and } (\phi \square \psi \omega)^{012}_s(k^2_s) = - \omega^{002}_{s_0d_1u}(0_s k^2_s 0_{s_0d_1s} \square).
        \end{split}
    \end{equation*}

    Define $h^0_r := \square k^0_s 0_t 0_u \in (\ker\widetilde{p}^2_0)_r$, so that 
    \begin{equation*}
        (\phi \square \psi \omega)^{012}_s(k^0_s) = (\phi \square \psi \omega)^{012}_s(h^0_r \square 0_t 0_u) = \phi^{012}_r(h^0_r) = - \phi^{112}_{s_0d_0r}(\square 0_r h^0_r 0_{s_1d_2r}),
    \end{equation*}
    by \eqref{eq:VB2DualElementsLv2Def} for $\phi$.
    By existence of the simplicial matrix 
    \begin{equation*}
        \begin{pmatrix}
            \square 0_r h^0_r 0_{s_1d_2r} &0_r &h^0_r &0_{s_1d_2r}\\
            \square 0_s k^0_s 0_{s_1d_2s} &0_s &k^0_s &0_{s_1d_2s}\\
            0_r &0_s &0_t &0_u\\
            h^0_r &k^0_s &0_t &0_u\\
            0_{s_1d_0u} &0_{s_1d_1u} &0_u &0_u
        \end{pmatrix}
    \end{equation*}
    in $\huaV_4$ over $s_2(r,s,t,u)$, 
    \begin{equation*}
        (\phi \square \psi \omega)^{012}_s(k^0_s)= - \phi^{112}_{s_0d_0r}(\square 0_r h^0_r 0_{s_1d_2r}) = - \phi^{112}_{s_0d_0r}(\square 0_s k^0_s 0_{s_1d_2s}).
    \end{equation*}

    Define $h^1_u := 0_r k^2_s 0_t \square \in (\ker\widetilde{p}^2_1)_u$, so that 
    \begin{equation*}
        (\phi \square \psi \omega)^{012}_s(k^2_s) = (\phi \square \psi \omega)^{012}_s(0_r \square 0_t h^1_u) = \omega^{012}(h^1_u) = - \omega^{002}_{s_0d_1u}(0_u h^1_u \square 0_{s_0d_2u}),
    \end{equation*}
    where the last equality follows from \eqref{eq:VB2DualFaceMapsPhiMultCond0012} for $\omega$. 
    By existence of the simplicial matrix
    \begin{equation*}
        \begin{pmatrix}
            0_r &0_s &0_t &0_u\\
            0_r &k^2_s &0_t &h^2_u\\
            0_s &k^2_s &0_{s_0d_1s} &0_s k^2_s 0_{s_0d_1s} \square \\
            0_t &0_t &0_{s_0d_1t} &0_{s_0d_2t}\\
            0_u &h^1_u &0_u h^1_u \square 0_{s_0d_2u} &0_{s_0d_2u}
        \end{pmatrix}
    \end{equation*}
    in $\huaV_4$ over $s_0(r,s,t,u)$, 
    \begin{equation*}
        (\phi \square \psi \omega)^{012}_s(k^2_s) = - \omega^{002}_{s_0d_1u}(0_u h^1_u \square 0_{s_0d_2u}) = - \omega^{002}_{s_0d_1u}(0_s k^2_s 0_{s_0d_1s} \square).
    \end{equation*}

    We now check step 4 for $\widecheck{m}_1$. For any $(\phi_r, \psi_t, \omega_u) \in \Lambda^3_1(\huaV^{2*})_{(r,t,u)}$ with $s=r\square tu$, we need to show that
    \begin{equation*}
        \widecheck{d}_0(\phi \square \psi \omega)=\widecheck{d}_0\phi,
        \qquad \widecheck{d}_1(\phi \square \psi \omega)=\widecheck{d}_1\psi,
        \qquad \widecheck{d}_2(\phi \square \psi \omega)=\widecheck{d}_1\omega.
    \end{equation*}
    The conditions on $\widecheck{d}_0$ and on the first component of $\widecheck{d}_2$ are true by definition \eqref{eq:VB2DualMultiplicationsDef}. Hence we need to check the ones for the $002$, $022$ and $011$ components of $(\phi \square \psi \omega)$. 
    That is, we need to show that 
    \begin{equation*}
        (\phi \square \psi \omega)^{002}_{s_0d_1s} = \psi^{002}_{s_0d_1t},
        \qquad (\phi \square \psi \omega)^{022}_{s_1d_1s} = \psi^{022}_{s_1d_1t},
        \qquad (\phi \square \psi \omega)^{011}_{s_1d_2s} = \omega^{022}_{s_1d_1u}.
    \end{equation*}

    Starting from the $002$ component, consider an arbitrary $Y\in (\huaV_2)_{s_0d_1s}$, and define $Y^0 : = Y - \widetilde{s}_0\widetilde{d}_0Y$. 
    Then choose any (2,1)-filler $Z^0$ of $(0_{d_2s}, \widetilde{d}_2Y^0)$ over $s_0d_2s = s_0d_1u$, and (2,1)-filler $U$ of $(0_{d_2t}, \widetilde{d}_2Y^0) = (0_{d_2u}, \widetilde{d}_2Z^0)$ over $s_0d_2t = s_0d_2u$. 
    Now, by \eqref{eq:VB2DualMultiplicationsDef}, \eqref{eq:VB2DualFaceMapsPhi002Def} reads
    \begin{equation*}
        (\phi \square \psi \omega)^{002}_{s_0d_1s}(Y) = \omega^{002}_{s_0d_1u}(Z^0) + (\phi\square \psi \omega)^{012}_s(0_s \square Y^0 Z^0).
    \end{equation*}
    Observe that the simplicial matrix
    \begin{equation*}
        \begin{pmatrix}
            0_r &0_r &0_t &0_u\\
            0_r &0_r \square (0_t\square Y^0 U)(0_u\square Z^0 U) &(0_t\square Y^0 U) &(0_u\square Z^0 U)\\
            0_s &0_s \square Y^0 Z^0 &Y^0_{s_0d_1s} &Z^0_{s_0d_2s}\\
            0_t &(0_t\square Y^0 U) &Y^0_{s_0d_1t} &U_{s_0d_2t}\\
            0_u &(0_u\square Z^0 U) &Z^0_{s_0d_1u} &U_{s_0d_2u}
        \end{pmatrix}
    \end{equation*}
    is a valid element of $\huaV_4$ over $s_0(r,s,t,u)$, so that 
    \begin{equation*}
        \begin{split}
            (\phi \square \psi \omega)^{002}_{s_0d_1s}(Y) &= \omega^{002}_{s_0d_1u}(Z^0) + (\phi\square \psi \omega)^{012}_s(0_r \square (0_t\square Y^0 U)(0_u\square Z^0 U))\\
            &= \omega^{002}_{s_0d_1u}(Z^0) + \psi^{012}_t(0_t\square Y^0 U) -\omega^{012}_u(0_u\square Z^0 U)\\
            &= \omega^{001}_{s_0d_2u}(U) + \psi^{012}_t(0_t\square Y^0 U) 
            = \psi^{002}_{s_0d_1t}(Y^0) = \psi^{002}_{s_0d_1t}(Y).
        \end{split}
    \end{equation*}
    Here we used Lemma \ref{lem:VB2DualFaceMapsPhiDefIntoMultEq} repeatedly, the horn condition $\omega^{001}_{s_0d_2u} = \psi^{001}_{s_0d_2t}$, and the normalization of $\psi^{002}$.

    For the $022$ component, consider an arbitrary $X\in (\huaV_2)_{s_1d_1s}$ and define $X^2:= X - \widetilde{s}_1\widetilde{d}_2X$. Then choose any (2,1)-horn filler $W^2$ of $(\widetilde{d}_0X^2, 0_{d_0s})$ over $s_1d_0s = s_1d_0r$, and any (2,1)-horn filler $U$ of $(\widetilde{d}_0W^2, 0_{d_1r})=(\widetilde{d}_0X^2, 0_{d_0t})$ over $s_1d_1r = s_1d_0t$.
    Now, by \eqref{eq:VB2DualMultiplicationsDef}, \eqref{eq:VB2DualFaceMapsPhi022Def} reads
    \begin{equation*}
        (\phi \square \psi \omega)^{022}_{s_1d_1s}(X) = \phi^{122}_{s_1d_0r}(W^2) + (\phi\square \psi \omega)^{012}_s(W^2 X^2 \square 0_s).
    \end{equation*}
    Observe that the simplicial matrix 
    \begin{equation*}
        \begin{pmatrix}
            W^2_{s_1d_0r} &U_{s_1d_1r} &(W^2 U\square 0_r) &0_r\\
            W^2_{s_1d_0s} &X^2_{s_1d_1s} &W^2 X^2 \square 0_s &0_s\\
            U_{s_1d_0t} &X^2_{s_1d_1t} &(U X^2 \square 0_t) &0_t\\
            (W^2 U \square 0_r) &(W^2 U \square 0_r)\square (U X^2 \square 0_t) 0_u &(U X^2 \square 0_t) &0_u\\
            0_r &0_s &0_t &0_u
        \end{pmatrix}
    \end{equation*}
    is a valid element of $\huaV_4$ over $s_3(r,s,t,u)$, so that 
    \begin{equation*}
        \begin{split}
            (\phi \square \psi \omega)^{022}_{s_1d_1s}(X) &= \phi^{122}_{s_1d_0r}(W^2) + (\phi\square \psi \omega)^{012}_s((W^2 U \square 0_r)\square (U X^2 \square 0_t) 0_u)\\
            &= \phi^{122}_{s_1d_0r}(W^2) + \phi^{012}_r(W^2 U \square 0_r) + \psi^{012}_t(U X^2 \square 0_t)\\
            &= \phi^{022}_{s_1d_1r} + \psi^{012}_t(U X^2 \square 0_t)
            = \psi^{022}_{s_1d_1t}(X^2) = \psi^{022}_{s_1d_1t}(X).
        \end{split}
    \end{equation*}
    Here we used again Lemma \ref{lem:VB2DualFaceMapsPhiDefIntoMultEq} repeatedly, the horn condition $\phi^{022}_{s_1d_1r} = \psi^{122}_{s_1d_0t}$, and the normalization of $\psi^{022}$. 

    Finally, for the $011$ component, consider an arbitrary $Z \in (\huaV_2)_{s_1d_2s}$ and define $Z^2:= Z - \widetilde{s}_1\widetilde{d}_2Z$. Then choose a (2,0)-filler $W^1$ of $(0_{d_0s}, \widetilde{d}_0Z^2)$ over $s_0d_0s = s_0d_0r$ and a (2,1)-filler $U$ of $(\widetilde{d}_2W^1, 0_{d_2r}) = (\widetilde{d}_0Z^2, 0_{d_0u})$ over $s_1d_2r = s_1d_0u$. 
    Now, by \eqref{eq:VB2DualMultiplicationsDef}, \eqref{eq:VB2DualFaceMapsPhi011Def} reads
    \begin{equation*}
        (\phi \square \psi \omega)^{011}_{s_1d_2s}(Z) = \phi^{112}_{s_0d_0r}(W^1) + (\phi\square \psi \omega)^{012}_s(W^2 X^2 \square 0_s).
    \end{equation*}
    Observe that the simplicial matrix 
    \begin{equation*}
        \begin{pmatrix}
            W^1_{s_0d_0r} &(W^1\square 0_r U) &0_r &U_{s_1d_2r}\\
            W^1_{s_0d_0s} &W^1\square 0_s Z^2 &0_s &Z^2_{s_1d_2s}\\
            (W^1\square 0_r U) &(W^1\square 0_r U)\square 0_t (U Z^2 \square 0_u) &0_t &(U Z^2 \square 0_u)\\
            0_r &0_s &0_t &0_u\\
            U_{s_1d_0u} &Z^2_{s_1d_1u} &(U Z^2 \square 0_u) &0_u
        \end{pmatrix}
    \end{equation*}
    is a valid element of $\huaV_4$ over $s_2(r,s,t,u)$, so that 
    \begin{equation*}
        \begin{split}
            (\phi \square \psi \omega)^{011}_{s_1d_2s}(Z) &= \phi^{112}_{s_0d_0r}(W^1) + (\phi\square \psi \omega)^{012}_s((W^1\square 0_r U)\square 0_t (U Z^2 \square 0_u))\\
            &= \phi^{112}_{s_0d_0s}(W^1) - \phi^{012}_r(W^1 \square 0_r U) + \omega^{012}_u(U Z^2 \square 0_u)\\
            &= \phi^{011}_{s_1d_2r}(U) + \omega^{012}_u(U Z^2 \square 0_u) 
            = \omega^{022}_{s_1d_1u}(Z^2) = \omega^{022}_{s_1d_1u}(Z).
        \end{split}
    \end{equation*}
    Here we used again Lemma \ref{lem:VB2DualFaceMapsPhiDefIntoMultEq} repeatedly, the horn condition $\phi^{011}_{s_1d_2r} = \omega^{122}_{s_1d_0u}$, and the normalization of $\omega^{022}$. 

    At this point we have some finite data as in Example \ref{ex:finite-data-2gpds}, so we need to check that this satisfies the conditions therein. 
    Beginning with the first condition, we need to show $\widecheck{d}_0^1, \widecheck{d}_1^1$ and $\widecheck{p}^2_i$ for $i=0,1,2$ are fiberwise surjective bundle maps over surjective submersions. 
    Since $\huaG$ is a Lie 2-groupoid, its base maps are all surjective submersions. 

    Starting with $\widecheck{d}_0^1$, we need to show that for any $g \in G_1$ and $\epsilon \in (\huaV^{2*}_0)_{d_0g}$, there exists a $(\eta, \xi) \in (\huaV^{2*}_1)_g$ such that $\widecheck{d}_0(\eta,\xi) = \epsilon$.
    Consider $X^0 \in (\ker \widetilde{d}^2_0)_{s_0g}$ and $Y^2 \in (\ker \widetilde{d}^2_2)_{s_1g}$. Then, by choosing splittings
\begin{equation}\label{eq:VB2DualSplittingLevel1}
    \begin{tikzcd}[ampersand replacement=\&,cramped,sep=scriptsize]
        0 \& {(\ker\widetilde{p}^2_1)_{s_0g}} \& {(\ker\widetilde{d}_0^2)_{s_0g}} \& {(\ker\widetilde{d}^1_1)_{1 d_1 g}} \& 0, \\
        0 \& {(\ker\widetilde{p}^2_1)_{s_1g}} \& {(\ker\widetilde{d}^2_2)_{s_1g}} \& {(\ker\widetilde{d}^1_0)_{1 d_0 g}} \& 0,
        \arrow[from=1-1, to=1-2]
        \arrow[hook, from=1-2, to=1-3]
        \arrow["{\widetilde{d}_2}", shift left, from=1-3, to=1-4]
        \arrow[from=1-4, to=1-5]
        \arrow[from=2-1, to=2-2]
        \arrow[hook, from=2-2, to=2-3]
        \arrow["{\widetilde{d}_0}", shift left, from=2-3, to=2-4]
        \arrow[from=2-4, to=2-5]
        \arrow["{\sigma_0}", shift left, dashed, from=1-4, to=1-3]
        \arrow["{\sigma_2}", shift left, dashed, from=2-4, to=2-3]
    \end{tikzcd}
\end{equation}
which always exist for each $g\in G_1$, we can write
\begin{equation*}
    X^0 = k + \sigma_0(c^1) \qquad Y^2 = h + \sigma_2(c^0),
\end{equation*}
for some $k \in (\ker\widetilde{p}^2_1)_{s_0g}$, $c^1 \in (\ker\widetilde{d}^1_1)_{1 d_1 g}$, $h \in (\ker\widetilde{p}^2_1)_{s_1g}$ and $c^0 \in (\ker\widetilde{d}^1_1)_{1 d_1 g}$.
Thus we can define
\begin{equation}
    \begin{split}
        \eta_{s_0g}(X) = \eta_{s_0g}(X^0) 
        &:= \epsilon_{1d_0g}(\square (0_{s_1g} \square k 0_{s_0g}) 0_{s_1g} 0_{s_1g})\\
        \xi_{s_1g}(Y) = \xi_{s_1g}(Y^2)
        &:= \epsilon_{1d_0g}(\square h_{s_1g} 0_{s_1g} 0_{s_1g}),
    \end{split}
\end{equation}
for any $X \in (\huaV_2)_{s_0g}$, $Y \in (\huaV_2)_{s_1g}$, with $X^0 := X - \widetilde{s}_0\widetilde{d}_0 X$, $Y^2 := Y - \widetilde{s}_1\widetilde{d}_2$, by using the usual isomorphism of annihilators and duals of kernels. 
Then, clearly $(\eta, \xi)$ satisfies \eqref{eq:VB2DualElementsLv1Def}, so it is a well-defined element of $(\huaV^{2*}_1)_g$, and $\widecheck{d}_0(\eta,\xi) = \epsilon$. 
Hence $\widecheck{d}_0^1$ is fiberwise surjective. 
An analogous construction can be used to show that $\widecheck{d}^1_1$ is fiberwise surjective: for any $\epsilon \in (\huaV^{2*}_0)_{1d_1g}$, setting 
\begin{equation}
    \begin{split}
        \eta_{s_0g}(X) = \eta_{s_0g}(X^0) 
        &:= \epsilon_{1d_0g}(0_{s_0g} 0_{s_0g} k_{s_0g} \square)\\
        \xi_{s_1g}(Y) = \xi_{s_1g}(Y^2)
        &:= \epsilon_{1d_0g}(0_{s_0g} 0_{s_0g} (0_{s_1g} h_{s_1g} \square 0_{s_0g}) \square),
    \end{split}
\end{equation}
gives an element $(\eta, \xi) \in (\huaV^{2*}_1)_g$ such that $\widecheck{d}_1(\eta,\xi) = \epsilon$.

Moving on to $\widecheck{p}^2_1$, we need to show that for any $t\in G_2$, and $((\eta_0, \xi_0), (\eta_2, \xi_2)) \in \Lambda^2_1(\huaV^{2*})_{d_0t,d_2t}$, there exists an element $\phi \in (\huaV^{2*}_2)_t$ such that 
\begin{equation*}
    \widecheck{p}^2_1(\phi) = ((\eta_0, \xi_0), (\eta_2, \xi_2)).
\end{equation*}
Since $t$ is fixed, we can always choose a splitting
\[\begin{tikzcd}[ampersand replacement=\&,cramped,sep=scriptsize]
	0 \& {(\ker\widetilde{d}_1)_{t}} \& {(\huaV_2)_t} \& {(d_1^*\huaV_1)_{t}} \& 0,
	\arrow[from=1-1, to=1-2]
	\arrow[hook, from=1-2, to=1-3]
	\arrow["{\widetilde{d}_1}", shift left, from=1-3, to=1-4]
	\arrow[from=1-4, to=1-5]
	\arrow["{\sigma_1}", shift left, dashed, from=1-4, to=1-3]
\end{tikzcd}\]
so that any $X \in (\huaV_2)_t$ can be written as $X = X^1 + \sigma_1\widetilde{d}_1X$, for $X^1 \in (\ker\widetilde{d}_1)_t$. 
With this, we can define 
\begin{equation*}
    \phi^{012}(X) := \eta_0(\square X^1  0_t Z) - \xi_2(Z),
\end{equation*}
for $X^1 = X - \sigma_1\widetilde{d}_1X$, and $Z$ any (2,0)-filler of $((\widetilde{d}_2X^1)_{d_2t}, 0_{d_2t})$ over $s_1d_2t$. 
We claim this is a well-defined element of $(\huaV_2^*)_t$ and that $\phi:=(\phi^{012}, (\eta_0, \xi_0), \eta_2)$ is the required filler. Hence we have to show that this definition does not depend on the filler $Z$ and that $\phi^{011} = \xi_2$. Since any other (2,0)-filler of $((\widetilde{d}_2X^1)_{d_2t}, 0_{d_2t})$ over $s_1d_2t$ can be written as $Z + k$ for $k \in (\ker \widetilde{p}^2_0)_{s_1d_2t}$, the definition of $\phi^{012}$ is filler independent if and only if 
\begin{equation}\label{eq:VB2DualProofp21CoverFillerIndep}
    \eta_0(\square 0_t  0_t k) - \xi_2(k) = 0.
\end{equation}
Note that the horn condition $\widecheck{d}_1(\eta_0, \xi_0) = \widecheck{d}_0(\eta_2, \xi_0)$ can be written as
\begin{equation*}
        \xi_2(k) = \eta_0(0_{s_0d_0t} 0_{s_0d_0t} \square (\ggamma_2(\square k 0_{s_1d_2 t} 0_{s_1d_2 t})))
        = \eta_0(0_{s_0d_0t} 0_{s_0d_0t} \square (\square 0_{s_1d_2 t} 0_{s_1d_2 t} k)),
\end{equation*}
for any $k \in (\ker \widetilde{p}^2_0)_{s_1d_2t}$, where we used the fact that, for such $k$,
\begin{equation*}
    \ggamma_2(\square k 0_{s_1d_2 t} 0_{s_1d_2 t}) = \square 0_{1d_0d_2t} 0_{1d_0d_2t} (\square k 0_{s_1d_2 t} 0_{s_1d_2 t}) = \square 0_{s_1d_2 t} 0_{s_1d_2 t} k.
\end{equation*}
Here the first equality follows from the canonical form of the multiplication on the VS 2-groupoid fiber over $1d_0d_2t$, while the second one follows from the associativity condition encoded by the simplicial matrix
\begin{equation*}
    \begin{pmatrix}
        \square 0_{1d_0d_2t} 0_{1d_0d_2t} (\square k 0_{s_1d_2 t} 0_{s_1d_2 t})
        &0_{1d_0d_2t} 
        &0_{1d_0d_2t}  
        &(\square k 0_{s_1d_2 t} 0_{s_1d_2 t})\\ 
        \square 0_{s_1d_2 t} 0_{s_1d_2 t} k 
        &0_{s_1d_2 t} 
        &0_{s_1d_2 t} 
        &k\\ 
        0_{1d_0d_2t} &0_{s_1d_2 t} &0_{s_1d_2 t} &0_{s_1d_2 t}\\
        0_{1d_0d_2t} &0_{s_1d_2 t} &0_{s_1d_2 t} &0_{s_1d_2 t}\\
        \square k 0_{s_1d_2 t} 0_{s_1d_2 t}
        &k
        &0_{s_1d_2 t}
        &0_{s_1d_2 t}\\
    \end{pmatrix}
\end{equation*}
which is an element of $\huaV_4$ over $s_3s_2s_1d_2t$. Now \eqref{eq:VB2DualProofp21CoverFillerIndep} can be rewritten as
\begin{equation*}
    \eta_0(\square 0_t 0_t k) = \eta_0(\eta_0(0_{s_0d_0t} 0_{s_0d_0t} \square (\square 0_{s_1d_2 t} 0_{s_1d_2 t} k))),
\end{equation*}
which holds, since 
\begin{equation*}
    \square 0_t 0_t k = 0_{s_0d_0t} 0_{s_0d_0t} \square (\square 0_{s_1d_2 t} 0_{s_1d_2 t} k),
\end{equation*}
by the associativity condition encoded by the simplicial matrix
\begin{equation*}
    \begin{pmatrix}
        0_{s_0d_0 t}
        &0_{s_0d_0 t}
        &0_{s_0d_0t} 0_{s_0d_0t} \square (\square 0_{s_1d_2 t} 0_{s_1d_2 t} k)  
        &(\square 0_{s_1d_2 t} 0_{s_1d_2 t} k)\\ 
        0_{s_0d_0 t} &0_{t} &0_{t} &0_{s_1d_2 t}\\ 
        0_{s_0d_0 t} &0_{t} &0_{t} &0_{s_1d_2 t}\\
        \square 0_t 0_t k &0_{t} &0_{t} &k\\
        \square 0_{s_1d_2 t} 0_{s_1d_2 t} k
        &0_{s_1d_2 t}
        &0_{s_1d_2 t}
        &k\\
    \end{pmatrix}
\end{equation*}
which is an element of $\huaV_4$ over $s_2s_1 t$. So $\phi^{012}$ defined above is well-defined. Moreover, $\widecheck{d}_0\phi = (\eta_0, \xi_0)$, and the first component of $\widecheck{d}_2\phi$ is $\eta_2$. The second component of $\widecheck{d}_2\phi$ is defined by \eqref{eq:VB2DualFaceMapsPhi011Def}, and by choosing the (2,0)-filler of $(0_{d_0t}, (\widetilde{d}_0Z)_{1d_1d_0t})$ over $s_0d_0t$ to be $(\square X^1  0_t Z)$ in that equation, we get 
\begin{equation*}
    \begin{split}
        \phi^{011}_{s_1d_2t}(Z) &= \eta_0(\square X^1  0_t Z) - \phi^{012}_t((\square X^1  0_t Z)\square 0_t Z)\\
        &= \eta_0(\square X^1  0_t Z) - \eta_0(\square X^1  0_t Z) + \xi_2(Z) 
        = \xi_2(Z),
    \end{split}
\end{equation*}
since $(\square X^1  0_t Z)\square 0_t Z = X^1$. This proves that $\widecheck{p}^2_1$ is fiberwise surjective, i.e. a cover. 

We now show $\widecheck{p}^2_2$ is fiberwise surjective. That is we show that for any $t\in G_2$, and $((\eta_0, \xi_0), (\eta_1, \xi_1)) \in \Lambda^2_2(\huaV^{2*})_{d_0t,d_1t}$, there exists an element $\phi \in (\huaV^{2*}_2)_t$ such that $\widecheck{p}^2_0(\phi) = ((\eta_0, \xi_0), (\eta_1, \xi_1))$.
Once again, we choose a splitting
\[\begin{tikzcd}[ampersand replacement=\&,cramped,sep=scriptsize]
	0 \& {(\ker\widetilde{d}_2)_{t}} \& {(\huaV_2)_t} \& {(d_2^*\huaV_1)_{t}} \& 0,
	\arrow[from=1-1, to=1-2]
	\arrow[hook, from=1-2, to=1-3]
	\arrow["{\widetilde{d}_2}", shift left, from=1-3, to=1-4]
	\arrow[from=1-4, to=1-5]
	\arrow["{\sigma_2}", shift left, dashed, from=1-4, to=1-3]
\end{tikzcd}\]
so that any $Y \in (\huaV_2)_t$ can be written as $Y = Y^2 + \sigma_2\widetilde{d}_2Y$, for $Y^2 \in (\ker\widetilde{d}_2)_t$. 
With this we define
\begin{equation*}
    \phi^{012}_t(Y) := \xi_1(X) - \xi_0(\square X Y^2 0_t),
\end{equation*}
for $Y^2 = Y - \sigma_2d_2Y$, and $X$ any (2,0)-filler of $((\widetilde{d}_1Y^2)_{d_1t}, 0_{d_1t})$ over $s_1d_1t$. Since we want $\eta_1 = \phi^{002}$, we can reformulate \eqref{eq:VB2DualFaceMapsPhi002Def} to define $\phi^{001}$ for any $Z \in (\huaV_2)_{s_0d_2t}$ as
\begin{equation*}
    \phi^{001}(Z) := \eta_1(\widebar{Y}) - \phi^{012}(0_t \square \widebar{Y} Z^0),
\end{equation*}
for $Z^0 = Z - s_0d_0Z$ and $\widebar{Y}$ any (2,1)-filler of $(0_{d_1t}, (\widetilde{d}_2 Z)_{1d_1d_1t})$ over $s_0d_1t$. Then we claim that $\phi:=(\phi^{012}, (\eta_0, \xi_0), \phi^{001})$ is a (2,2)-filler of $((\eta_0, \xi_0), (\eta_1, \xi_1))$.
First of all, $\phi^{012}$ is well-defined: it is invariant under change of filler if and only if 
\begin{equation*}
    \xi_1(k) = \xi_0(\square k 0_t 0_t),
\end{equation*}
for any $k \in \ker \widetilde{p}^2_0$. By using the horn condition $\widecheck{d}_0(\eta_0, \xi_0) = \widecheck{d}_1(\eta_1, \xi_1)$, this is equivalent to 
\begin{equation*}
    \xi_0((\square k 0_{s_1d_1t} 0_{s_1d_1t}) \square 0_{s_1d_0t} 0_{s_1d_0t}) = \xi_0(\square k 0_t 0_t).
\end{equation*}
This holds since 
\begin{equation*}
    (\square k 0_{s_1d_1t} 0_{s_1d_1t}) \square 0_{s_1d_0t} 0_{s_1d_0t} = \xi_0(\square k 0_t 0_t),
\end{equation*}
for any $k \in \ker \widetilde{p}^2_0$, follows from the associativity condition encoded by the simplicial matrix 
\begin{equation*}
    \begin{pmatrix}
        \square k 0_{s_1d_1t} 0_{s_1d_1t}
        &(\square k 0_{s_1d_1t} 0_{s_1d_1t}) \square 0_{s_1d_0t} 0_{s_1d_0t}
        &0_{s_1d_0t}
        &0_{s_1d_0t}\\ 
        \square k 0_{s_1d_1t} 0_{s_1d_1t} &k &0_{s_1d_1t} &0_{s_1d_1t}\\ 
        (\square k 0_t 0_t) &k &0_{t} &0_{t}\\
        \square 0_t 0_t k &0_{t} &0_{t} &k\\
        0_{s_1d_0 t} &0_{s_1d_1t} &0_{t} &0_{t}\\
        0_{s_1d_0 t} &0_{s_1d_1t} &0_{t} &0_{t}\\
    \end{pmatrix}
\end{equation*}
over $s_3s_2t$. 
$\phi^{001}$ is also well-defined: it is filler-independent if and only if
\begin{equation*}
\begin{split}
\eta_1(h) &= \phi^{012}(0_t \square h 0_{s_0d_2t})\\
&= \xi_1((0_{s_1d_1t}\square h 0_{s_0d_1t})) - \xi_0(\square (0_{s_1d_1t}\square h 0_{s_0d_1t})(0_t \square h 0_{s_0d_2t}) 0_t)    
\end{split}
\end{equation*}
for any $h \in \ker (\widetilde{p}^2_1)_{s_0d_1t}$, where we chose $(0_{s_1d_1t}\square h 0_{s_0d_1t})$ as the (2,0)-filler of $((\widetilde{d}_1h)_{d_1t}, 0_{d_1t})$ over $s_1d_1t$, present in the definition of $\phi^{012}$. But since $(\eta_1, \xi_1) \in (\huaV^{2*})_1$, $\eta_1(h) = \xi_1((0_{s_1d_1t}\square h 0_{s_0d_1t}))$, so $\phi^{001}$ is filler-independent if and only if
\begin{equation*}
    \xi_0(\square (0_{s_1d_1t}\square h 0_{s_0d_1t})(0_t \square h 0_{s_0d_2t}) 0_t) = 0.
\end{equation*}
But this is true, since the simplicial matrix
\begin{equation*}
    \begin{pmatrix}
        0_{s_1d_0t}
        &0_{s_1d_1t}
        &0_{t}
        &0_{t}\\ 
        \square (0_{s_1d_1t}\square h 0_{s_0d_1t})(0_t \square h 0_{s_0d_2t}) 0_t 
        &(0_{s_1d_1t}\square h 0_{s_0d_1t}) &(0_t \square h 0_{s_0d_2t}) &0_t\\ 
        0_{s_1d_1t} &(0_{s_1d_1t}\square h 0_{s_0d_1t}) &h &0_{s_0d_1t}\\
        0_t &(0_t \square h 0_{s_0d_2t}) &h &0_{s_0d_2t}\\
        0_{t} &0_{t} &0_{s_0d_1t} &0_{s_0d_2t}\\
    \end{pmatrix}
\end{equation*}
over $s_3s_0t$ encodes the associativity condition 
\begin{equation*}
    \square (0_{s_1d_1t}\square h 0_{s_0d_1t})(0_t \square h 0_{s_0d_2t}) 0_t = 0_{s_1d_0t}.
\end{equation*}
Clearly $\widecheck{d}_0 \phi = (\eta_0, \xi_0)$. We now check that $\widecheck{d}_1 \phi = (\eta_1, \xi_1)$: By \eqref{eq:VB2DualFaceMapsPhi002Def}, we have
\begin{equation*}
    \begin{split}
        \phi^{002}_{s_0d_1t}(Y) &= \phi^{001}_{s_0d_2t}(Z) + \phi^{012}_t(0_{t}\square Y^0 Z) \\
        &= \eta_1(Y) - \phi^{012}_t(0_{t}\square Y^0 Z) + \phi^{012}_t(0_{t}\square Y^0 Z) = \eta_1(Y),
    \end{split}
\end{equation*}
where $Z$ is any (2,1)-filler of $(0_{d_2t}, \widetilde{d}_2 Y^0)$ over $s_0d_2t$, and in the definition of $\phi^{001}$ we can choose as filler $Y^0 := Y - \widetilde{s}_0\widetilde{d}_0Y$. 
By \eqref{eq:VB2DualFaceMapsPhi022Def}, we have
\begin{equation*}
    \begin{split}
        \phi^{022}_{s_1d_1t}(X) &= \xi_0(W) + \phi^{012}_t(W X^2 \square 0_t)\\
        & = \xi_0(W) + \xi_1(X) - \xi_0(\square X (W X^2 \square 0_t) 0_t) = \xi_1(X),
    \end{split}
\end{equation*}
where $W$ is any (2,1)-filler of $(\widetilde{d}_0X, 0_{d_0t})$ over $s_1d_0t$, and in the definition of $\phi^{012}$, we chose $X^2 = X - \widetilde{s}_1\widetilde{d}_2X$ as filler. We also used that $W X^2 \square 0_t$ is in the kernel of $\widetilde{d}_2$ and that $\square X (W X^2 \square 0_t) 0_t = W$. Hence $\widecheck{p}^2_2$ is fiberwise surjective, i.e. a cover. 

Finally, a $(2,0)$-filler of any $((\eta_1, \xi_1),(\eta_2,\xi_2)) \in \Lambda^2_0(\huaV^{2*})_{(d_1t,d_2t)}$ over $t$ can be constructed by entirely analogous methods, first by using \eqref{eq:VB2DualFaceMapsPhi002Def} to deduce the restriction of $\phi^{012}$ on $\ker\widetilde{d}_0$, and then using the other identities in Lemma \ref{lem:VB2DualFaceMapsPhiDefIntoMultEq} to obtain $\phi^{112}$ and $\phi^{122}$.

We now check the conditions of Example \ref{ex:finite-data-2gpds} on the multiplications of $\huaV^{2*}$ are satisfied. First, we need to check that the multiplications induce isomorphisms between the horn spaces $\Lambda^3_i(\huaV^{2*})$ and $\Lambda^3_j(\huaV^{2*})$ for any pair of $i,j$, with $i=0,1,2,3$ and $j=0,1,2,3$. We check this for a single fixed pair of indices and argue that it can be checked in exactly the same way for any other pair of indices. Take for example $i=0$ and $j=1$. We show the map 
\begin{equation*}
    \begin{split}
        &\Lambda^3_1(\huaV^{2*}) \to \Lambda^3_0(\huaV^{2*})\\
        &(\phi, \psi, \omega) \to (\phi\square\psi\omega, \psi, \omega),
    \end{split}
\end{equation*}
is an isomorphism with inverse
\begin{equation*}
    \begin{split}
        &\Lambda^3_0(\huaV^{2*}) \to \Lambda^3_1(\huaV^{2*})\\
        &(\chi, \psi, \omega) \to (\square\chi\psi\omega, \psi, \omega).
    \end{split}
\end{equation*}
First of all, both maps are well-defined, by step 4 above. 
Moreover, by definition, for any $(\phi, \psi, \omega)\in \Lambda^3_1(\huaV^{2*})_{(r,t,u)}$ and any $(\chi, \psi, \omega)\in \Lambda^3_1(\huaV^{2*})_{(s,t,u)}$,
\begin{equation*}
    \begin{split}
        \square(\phi\square\psi\omega)\psi\omega = ((\square(\phi\square\psi\omega)\psi\omega)^{012}, (\phi^{112}, \phi^{122}), \omega^{112}),\\
        (\square\chi\psi\omega)\square\psi\omega = (((\square\chi\psi\omega)\square\psi\omega)^{012}, (\chi^{112}, \chi^{122}), \omega^{002}).\\
    \end{split}
\end{equation*}
By the horn gluing conditions $\omega^{112} = \phi^{001}$ and $\omega^{002} = \chi^{001}$. Finally, we need to show that 
\begin{equation*}
    (\square(\phi\square\psi\omega)\psi\omega)^{012} = \phi^{012}, \text{ and } ((\square\chi\psi\omega)\square\psi\omega)^{012} = \chi^{012}.
\end{equation*}
For any tetrahedron $(W,X,Y,Z) in (\huaV_3)_{(r,s,t,u)}$, we have $X= W\square YZ$ and $W = \square XYZ$, so that 
\begin{equation*}
    \begin{split}
        (\square(\phi\square\psi\omega)\psi\omega&)^{012}(\square X Y Z) = (\phi\square\psi\omega)^{012}_s(X) -\psi^{012}_t(Y) + \omega^{012}_u (Z) \\
        &= \phi^{012}_r(W) + \psi^{012}_t(Y) - \omega^{012}_u(Z) -\psi^{012}_t(Y) + \omega^{012}_u (Z)\\ 
        &= \phi^{012}_r(W),
    \end{split}
\end{equation*}
and 
\begin{equation*}
    \begin{split}
        ((\square\chi\psi\omega)\square\psi\omega&)^{012}(W \square Y Z) = (\square\chi\psi\omega)^{012}_r(W) + \psi^{012}_t(Y) - \omega^{012}_u (Z) \\
        &= \chi^{012}_s(X) - \psi^{012}_t(Y) + \omega^{012}_u(Z) + \psi^{012}_t(Y) - \omega^{012}_u (Z)\\
        &= \chi^{012}_s(X).
    \end{split}
\end{equation*}
Hence the multiplications induce isomorphisms between the horn spaces.

The second condition is compatibility of the multiplications with the degeneracy maps of the 2-coskeleton of $\huaV^{2*}_{\le 2}$. 
We show this in the form \eqref{eq:finite-data-comp-deg-2}, which in this case amounts to the three equations
\begin{equation*}
    \phi = \phi \square (\widecheck{s}_0 \widecheck{d}_1 \phi) (\widecheck{s}_0 \widecheck{d}_2 \phi),
    \quad \phi = (\widecheck{s}_0 \widecheck{d}_0 \phi) \phi \square (\widecheck{s}_1 \widecheck{d}_2 \phi), 
    \quad \text{and } \phi = (\widecheck{s}_1 \widecheck{d}_0 \phi) (\widecheck{s}_1 \widecheck{d}_1 \phi) \phi \square,
\end{equation*} 
for any $\phi \in \huaV^{2*}_2$.

We will now show these equations are satisfied for any $\phi \in \huaV^{2*}$.
Starting with the first one, we have
\begin{equation*}
    \phi \square (\widecheck{s}_0 \widecheck{d}_1 \phi) (\widecheck{s}_0 \widecheck{d}_2 \phi) = ((\phi \square (\widecheck{s}_0 \widecheck{d}_1 \phi) (\widecheck{s}_0 \widecheck{d}_2 \phi))^{012}, (\phi^{112}, \phi^{122}), \phi^{001}),
\end{equation*}
where $(\widecheck{s}_0\widecheck{d}_2 \phi)^{002} = \phi^{001}$, by the simplicial identities of the 2-truncation of $\huaV^{2*}$ we previously showed.
Also, for any $(W_t, Y_{s_0d_1t}, Z_{s_0d_2t}) \in \Lambda^3_1(\huaV)$,
\begin{equation*}
\begin{split}
    (\phi \square (\widecheck{s}_0 \widecheck{d}_1 \phi) (\widecheck{s}_0 \widecheck{d}_2 \phi))^{012} (W_t\square Y_{s_0d_1t} Z_{s_0d_2t})
    &= \phi^{012}_t(W) + \phi^{002}_{s_0d_1t}(Y) - \phi^{001}_{s_0d_2t}(Z) \\
    &= \phi^{012}_t (W_t\square Y_{s_0d_1t} Z_{s_0d_2t}),
\end{split}
\end{equation*}
by \eqref{eq:VB2DualFaceMapsPhiMultCond0012}.
For the second one, 
\begin{equation*}
    (\widecheck{s}_0 \widecheck{d}_0 \phi) \phi \square (\widecheck{s}_1 \widecheck{d}_2 \phi) = (((\widecheck{s}_0 \widecheck{d}_0 \phi) \phi \square (\widecheck{s}_1 \widecheck{d}_2 \phi))^{012}, (\phi^{112}, \phi^{122}), \phi^{001}),
\end{equation*}
since $\widecheck{d}_1 \widecheck{s}_0 \widecheck{d}_1 \phi = \widecheck{d}_0 \phi$ and $\widecheck{d}_2 \widecheck{s}_1 \widecheck{d}_2\phi = \widecheck{d}_2\phi$, and, for any $(W_{s_0d_0t}, X_t, Z_{s_0d_2t}) \in \Lambda^3_2(\huaV)$, 
\begin{equation*}
\begin{split}
    ((\widecheck{s}_0 \widecheck{d}_0 \phi) \phi \square (\widecheck{s}_1 \widecheck{d}_2 \phi))^{012}(W_{s_0d_0t} X_t \square Z_{s_0d_2t}) 
    &= - \phi^{112}_{s_0d_0t}(W) + \phi^{012}_t(X) + \phi^{011}_{s_1d_2}(Z) \\
    &= \phi^{012}_t(W_{s_0d_0t} X_t \square Z_{s_0d_2t}),
\end{split}
\end{equation*}
by \eqref{eq:VB2DualFaceMapsPhiMultCond0112}.
Finally for the third one, 
\begin{equation*}
    (\widecheck{s}_1 \widecheck{d}_0 \phi) (\widecheck{s}_1 \widecheck{d}_1 \phi) \phi \square = (((\widecheck{s}_1 \widecheck{d}_0 \phi) (\widecheck{s}_1 \widecheck{d}_1 \phi) \phi \square)^{012}, (\phi^{112}, \phi^{122}), \phi^{001}),
\end{equation*}
since $\widecheck{d}_2\widecheck{s}_1\widecheck{d}_0 \phi = \widecheck{d}_0 \phi$, and, for any $(W_{s_1d_0t}, X_{s_1d_1t}, Y_{t}) \in \Lambda^3_3(\huaV)$, 
\begin{equation*}
\begin{split}
    ((\widecheck{s}_1 \widecheck{d}_0 \phi) (\widecheck{s}_1 \widecheck{d}_1 \phi) \phi \square)^{012}(W_{s_1d_0t} X_{s_1d_1t} Y_{t} \square) 
    &= \phi^{122}_{s_1d_0t}(W) - \phi^{022}_{s_1d_1t}(X) + \phi^{012}_t(Y)\\
    &= \phi^{012}_t(W_{s_1d_0t} X_{s_1d_1t} Y_{t} \square),
\end{split}
\end{equation*}
by \eqref{eq:VB2DualFaceMapsPhiMultCond0122}.

Finally, we check that the multiplications $\widecheck{m}_j$ are associative. Consider any 2-simplices $\phi_{0 i 4}$ and $\phi_{0ij}$ for $i \neq j \in \{1,2,3\}$ in $\huaV^{2*}_2$ that fit into the boundary of a 4-simplex over $t_{01234} \in \huaG_4$ as in Example \ref{ex:finite-data-2gpds}.
We show that 
\begin{equation*}
    (\square \phi_{034} \phi_{024} \phi_{023}) 
    (\square \phi_{034} \phi_{014} \phi_{013})
    (\square \phi_{024} \phi_{014} \phi_{012})
    \square 
    = \square \phi_{023} \phi_{013} \phi_{012},
\end{equation*}
where both sides have basepoint $t_{123}$ by associativity on the base. 
First of all, all the components of each side that are not the $012$ one are equal in the expected pairs, by repeatedly applying the simplicial identities.
Secondly, for any $X_{123} \in (\huaV_2)_{t_{123}}$, this can be decomposed as 
\begin{equation*}
    X_{123} = (\square X_{034} X_{024} X_{023}) 
    (\square X_{034} X_{014} X_{013})
    (\square X_{024} X_{014} X_{012})
    \square 
    = \square X_{023} X_{013} X_{012},
\end{equation*} 
for a 4-simplex $X_{01234}$ in $(\huaV_2)_{t_{01234}}$, by the Kan conditions of $\huaV_2$ (including the ones giving associativity in $\huaV_2$).
Then 
\begin{equation*}
\begin{split}
    ((\square &\phi_{034} \phi_{024} \phi_{023}) 
    (\square \phi_{034} \phi_{014} \phi_{013})
    (\square \phi_{024} \phi_{014} \phi_{012})
    \square)^{012} (X_{123}) \\
    &= (\square \phi_{034} \phi_{024} \phi_{023})^{012}(\square X_{034} X_{024} X_{023} )
    - (\square \phi_{034} \phi_{014} \phi_{013})^{012}(\square X_{034} X_{014} X_{013})\\
    &\quad + (\square \phi_{024} \phi_{014} \phi_{012})^{012}(\square X_{024} X_{014} X_{012})\\
    &= (\phi_{023})^{012}(X_{023}) - (\phi_{013})^{012}(X_{013}) + (\phi_{012})^{012}(X_{012})\\
    &= (\square \phi_{023} \phi_{013} \phi_{012})^{012} (X_{123}),
\end{split}
\end{equation*}
so the $\widecheck{m}_j$ are associative.
\end{proof}

\section{Examples of \texorpdfstring{$\VB$}{VB} \texorpdfstring{$n$}{n}-duals}\label{sec:vb-n-dual-examples}

We conclude this chapter by presenting some cases of $n$-duals for $n \le 2$ for $\VB$ groupoids of order less than 2, where the general formulas simplify. The case of the 0-dual of a $\VB$ $n$-groupoid over a manifold was considered in Remark \ref{rem:VBnGpd-0-dual}. 

\begin{example}[1- and 2- dual of a $\VB$ 0-groupoid]
    The most trivial example of a higher vector bundle is that of a $\VB$ 0-groupoid, which is just the unit groupoid of a vector bundle $E \to M$, for $M$ a smooth manifold. By the discussion in \ref{rem:VBnGpd-0-dual} this is the only order that always admits a 0-dual, which is just the unit $\VB$-groupoid of the dual vector bundle $E^*\to M$.
    For higher duals, observe that, in this case, all the $\widetilde{p}^m_j$ are isomorphisms for $m\ge 1$, $0\le j\le m$, where $\widetilde{p}^1_j = \widetilde{d}^1_j$. Therefore, all the various kernels present in higher duals vanish, and we have that the 1-dual is given by the finite data 
    \[\begin{tikzcd}[ampersand replacement=\&]
        {E^*} \& {0\times M} \\
        M \& M
        \arrow["{(0,q)}", shift left=1, from=1-1, to=1-2]
        \arrow["id", shift left=1, from=2-1, to=2-2]
        \arrow["id"', shift right=1, from=2-1, to=2-2]
        \arrow["{(0,q)}"', shift right=1, from=1-1, to=1-2]
        \arrow[from=1-2, to=2-2]
        \arrow["q"', from=1-1, to=2-1]
    \end{tikzcd}\]
    where the $(2,j)$-horn spaces are all isomorphic to the Whitney sum of two copies of $E^*$ (since source and target coincide with the bundle projection), and the multiplications are given by addition and subtraction of each fiber. That is, $\Lambda^2_j(E^{1*}) = E^* \oplus E^*$, and $\widecheck{m}_0(e_p,f_p) = e_p-f_p$, $\widecheck{m}_1(e_p,f_p) = e_p + f_p$, $\widecheck{m}_2(e_p,f_p) = - e_p + f_p$. 
    Meanwhile, the 2-dual is given by
    \[\begin{tikzcd}[ampersand replacement=\&, sep=scriptsize]
        {E^*} \& {0\times M} \& {0\times M} \\
        M \& M \& M
        \arrow[shift left=1, from=1-2, to=1-3]
        \arrow[shift left=1, from=2-2, to=2-3]
        \arrow[shift right=1, from=2-2, to=2-3]
        \arrow[from=1-3, to=2-3]
        \arrow[from=1-2, to=2-2]
        \arrow[from=2-1, to=2-2]
        \arrow["id", shift left=2, from=2-1, to=2-2]
        \arrow[shift right=2, from=2-1, to=2-2]
        \arrow["{(0,q)}", shift left=2, from=1-1, to=1-2]
        \arrow[shift right=2, from=1-1, to=1-2]
        \arrow[from=1-1, to=1-2]
        \arrow["q"', from=1-1, to=2-1]
        \arrow[shift right=1, from=1-2, to=1-3]
    \end{tikzcd}\]
    where the $(3,j)$-horn spaces are all isomorphic to the Whitney sum of three copies of $E^*$, and the multiplications are given by alternating sums in each fiber. That is, $\Lambda^3_j(E^{2*}) = E^* \oplus E^*\oplus E^*$, and
    \begin{equation*}
    \begin{split}
        &\widecheck{m}_0(d_p,e_p,f_p) = d_p - e_p + f_p = \widecheck{m}_3(d_p,e_p,f_p), \\
        &\widecheck{m}_1(d_p,e_p,f_p) = d_p + e_p - f_p, \quad \widecheck{m}_2(d_p,e_p,f_p) = - d_p + e_p + f_p
    \end{split}
    \end{equation*} 
    
    These results are consistent with the pointwise interpretation from \ref{lem:VBndual-at-units}, where a $\VB$ $n$-groupoid over a 0-groupoid is defined pointwise for each fiber, as each point $p$ of $M$ has an associated VS n-groupoid $(1^*E)_p$. Then the fiber of the $n$-dual at each point is given by the $n$-dual of this vector space groupoid. 
\end{example}
    
\begin{example}[2-dual of a $\VB$ 1-groupoid]\label{ex:VB2DualVB1Gpd}
    Consider the $\VB$ 1-groupoid $\huaV \to \huaG$, that is, up to level 2,
\[\begin{tikzcd}[ampersand replacement=\&, sep=scriptsize]
	{\Lambda^2_1(\huaV)} \& {\huaV_1} \& {\huaV_0} \\
	{\Lambda^2_1(\huaG)} \& {\huaG_1} \& {\huaG_0}
	\arrow[shift left=2, from=1-1, to=1-2]
	\arrow[shift right=2, from=1-1, to=1-2]
	\arrow[from=1-1, to=1-2]
	\arrow[from=1-1, to=2-1]
	\arrow[shift left, from=1-2, to=1-3]
	\arrow[shift right, from=1-2, to=1-3]
	\arrow[from=1-2, to=2-2]
	\arrow[from=1-3, to=2-3]
	\arrow[shift left=2, from=2-1, to=2-2]
	\arrow[shift right=2, from=2-1, to=2-2]
	\arrow[from=2-1, to=2-2]
	\arrow[shift left, from=2-2, to=2-3]
	\arrow[shift right, from=2-2, to=2-3]
\end{tikzcd}\]
    where we are identifying level 2 with the (2,1)-horn $\Lambda^2_1(\huaV) \cong \huaV_1\times_{\widetilde{d}_1,\huaV_0, \widetilde{d}_0}\huaV_1$, since $\huaV$ is a $\VB$ 1-groupoid. Hence the horn maps $\widetilde{p}^2_j$ are isomorphisms for any $j=0,1,2$. 
    For the same reason, the 2-dimensional degree 2 cores of $\huaV$ vanish. We now show how to write the 2-dual of $\huaV$ in terms of its degree 1 1-dimensional cores $C_R= 1^*\ker\widetilde{d}_0$ and $C_L= 1^*\ker\widetilde{d}_1$. As we showed in Lemma \ref{lem:VB1-core-involution-iso} and recalled in Section \ref{sec:VB-DegAnn-Cores}, these two form a core pair and they are isomorphic through the involution given by the core projections $\ggamma_R, \ggamma_L$:
    \begin{equation}\label{eq:VB2Dual1Gpd-core-proj}
        \begin{array}{cc}
        \begin{aligned}
            \ggamma_L: C_R &\longrightarrow C_L\\ 
            a &\mapsto -a^{-1} =  a - \widetilde{1}\widetilde{d}_0a,
        \end{aligned}
            &
        \begin{aligned}
            \ggamma_R: C_L &\longrightarrow C_R\\
            b &\mapsto -b^{-1} = b - \widetilde{1}\widetilde{d}_1b.
        \end{aligned}
        \end{array}
    \end{equation}
    First of all, by the isomorphism \eqref{eq:VBnDual-degn-degAnn-isom}, and the vanishing of the 2-dimensional degree 2 cores, $\huaV^{2*}_0 = 0 \times \huaG_0$.
    By using the isomorphisms in \eqref{eq:VBnDual-deg1-degAnn-isom}, the degree 1 degeneracy annihilators $O_0$ and $O_1$ are isomorphic to $s_0^*(\ker\widetilde{d}^2_0)^*$ and $s_1^*(\ker\widetilde{d}^2_2)^*$, respectively. But also 
    \begin{equation}\label{eq:VB2Dual1Gpd-deg1-2dim-cores}
        \begin{aligned}
            &s_0^*\ker\widetilde{d}^2_0 = \{(0_g, v_{1d_1g}) \in (V_1 \times V_1)_{(g,1d_1g)} \mid g \in \huaG_1, \widetilde{d}_0v_{1d_1g} = 0_{d_1g}\} \cong d_1^*C_R,\\
            &s_1^*\ker\widetilde{d}^2_2 = \{(v_{1d_0g}, 0_g) \in (V_1 \times V_1)_{(1d_0g,g)} \mid g \in \huaG_1, \widetilde{d}_1v_{1d_0g} = 0_{d_1g}\} \cong d_0^*C_L.
        \end{aligned}
    \end{equation} 
    So, as a vector bundle over $\huaG_1$, $\huaV^{2*}_1 \cong d_1^*C_R^* \oplus d_0^*C_L^*$. For the same reasons, $\huaV^{2*}_2 \cong (\huaV_1 \times_{\huaV_0} \huaV_1)^* \oplus d_0^*d_1^*C_R^* \oplus d_0^*d_0^*C_L^* \oplus d_2^*d_1^*C_R^*$. Therefore the 2-dual of $\huaV$ is
    \[\begin{tikzcd}[ampersand replacement=\&]
        {\begin{array}{c}\Lambda^2_1(\huaV)^* \oplus d_0^*d_1^*C_R^* \\ \quad \oplus d_0^*d_0^*C_L^*\oplus d_2^*d_1^*C_R^*\end{array}} \& {d_1^*C_R^*\oplus d_0^*C_L^*} \& {0\times \huaG_0} \\
        {\Lambda^2_1(\huaG)} \& {\huaG_1} \& {\huaG_0}
        \arrow[shift left=1, from=1-2, to=1-3]
        \arrow[shift left=1, from=2-2, to=2-3]
        \arrow[shift right=1, from=2-2, to=2-3]
        \arrow[from=1-3, to=2-3]
        \arrow[from=1-2, to=2-2]
        \arrow[shift right=1, from=1-2, to=1-3]
        \arrow[from=1-1, to=2-1]
        \arrow[shift left=2, from=2-1, to=2-2]
        \arrow[shift right=2, from=2-1, to=2-2]
        \arrow[from=2-1, to=2-2]
        \arrow[shift left=2, from=1-1, to=1-2]
        \arrow[shift right=2, from=1-1, to=1-2]
        \arrow[from=1-1, to=1-2]
    \end{tikzcd}\]
    Here, the face maps between levels 1 and 0 are the obvious 0 maps and the degeneracy map is also just the inclusion of the 0 section at the unit of each $p\in G_0$, $\widecheck{s}_0^0(0,p)=(0_{1p},0_{1p})$. To write the face maps between levels 2 and 1 we consider a generic element $\Phi \in \huaV^{2*}_2|_{(g,h)}$ for $(g,h)\in \Lambda^2_1(\huaG)$. Then $\Phi^{012} := \phi \in \Lambda^2_1(\huaV)^*|_{(g,h)}$, and by composing the isomorphisms \eqref{eq:VBnDual-deg1-degAnn-isom} and \eqref{eq:VB2Dual1Gpd-deg1-2dim-cores}, 
    \begin{equation}\label{eq:VB2Dual1Gpd-indep-face-comps}
        \begin{aligned}
            &\Phi^{112}|_{(g,1d_1g)}(0_g,a) = \alpha|_{d_1g}(a),\\
            &\Phi^{122}|_{(1d_0g,g)}(b,0_g) = \beta|_{d_0g}(b),\\
            &\Phi^{001}|_{(h,1d_1h)}(0_h,a') = \gamma|_{d_1h}(a'),
        \end{aligned}
    \end{equation}
    for $\alpha \in C_R^*|_{d_1g}$, $\beta \in C_L^*|_{d_0g}$ and $\gamma \in C_L^*|_{d_1h}$, and any arbitrary $a \in C_R|_{d_1g}$, $a' \in C_R|_{d_1h}$ and $b \in C_L|_{d_0g}$. 
    Then, with this and \eqref{eq:VB2DualFaceMapsPhi002Def} - \eqref{eq:VB2DualFaceMapsPhi011Def}, we compute, for any $a' \in C_R|_{d_1h}$, $b\in C_L|_{d_0g}$ and $b' \in C_L|_{d_0h = d_1g}$,
    \begin{equation}\label{eq:VB2Dual1Gpd-dep-face-comps}
        \begin{aligned}
            &\Phi^{002}|_{(hg,1d_1h)}(0_{hg}, a) = \gamma|_{d_1h}(a) + \phi|_{(g,h)}(0_g, 0_h\cdot a),\\
            &\Phi^{022}|_{(1d_0g, hg)}(b,0_{hg}) = \beta|_{d_0g}(b) + \phi|_{(g,h)}(b \cdot 0_g, 0_h),\\
            &\Phi^{011}|_{(1d_0h, h)}(b',0_{h}) = -\alpha|_{d_1g}((b')^{-1}) - \phi|_{(g,h)}((b')^{-1}\cdot 0_g, 0_h\cdot b'),\\
        \end{aligned}
    \end{equation}
    where we used the fact that the horn fillers appearing in \eqref{eq:VB2DualFaceMapsPhi002Def} - \eqref{eq:VB2DualFaceMapsPhi011Def} are unique and these are identified with the horns themselves as elements of $\huaV_2 \cong \Lambda^2_1(\huaV)$. Note also that in the last equation we can write $\alpha|_{d_1g}((b')^{-1})$ by the horn condition $d_1g=d_0h$. With this, the face maps are given as in \eqref{eq:VB2DualFaceMaps2}. We write them explicitly for the 2-cotangent of a 1-groupoid in \eqref{eq:2Cota-1Gpd-face-maps}.
    On the other hand, the degeneracy maps are given by \eqref{eq:VB2DualDegeneracyMaps2}: for any $(\alpha, \beta) \in d_1^*C_R^*\oplus d_0^*C_L^*|_g$, 
    \begin{equation}\label{eq:VB2Dual1Gpd-deg-maps}
        \begin{split}
            &\widecheck{s}_0(\alpha, \beta) = ((\widetilde{d}_2\ggamma_R)^*\alpha|_{(g, 1d_1g)}, (\alpha, \beta), 0),\\
            &\widecheck{s}_1(\alpha, \beta) = ((\widetilde{d}_0\ggamma_L)^*\beta_{(1d_0g, g)}, (0, 0), \alpha),
        \end{split}
    \end{equation}
    where $(\widetilde{d}_2\ggamma_R)^*: d_1^*C_R^* \to O_0 \subseteq s_0^*\Lambda^2_1(\huaV)^*$ and $(\widetilde{d}_0\ggamma_R)^*: d_0^*C_R^* \to O_1 \subseteq s_1^*\Lambda^2_1(\huaV)^*$ are the compositions of the isomorphisms described in \eqref{eq:VBnDual-deg1-degAnn-isom} and \eqref{eq:VB2Dual1Gpd-deg1-2dim-cores}. Here $\widetilde{d}_0$ and $\widetilde{d}_2$ are just the projections $\Lambda^2_1(\huaV) \to \huaV_1$.  

    We now describe the multiplications. For ease of reading we give a name to the following three maps appearing in the definitions of the face components \eqref{eq:VB2Dual1Gpd-dep-face-comps}, which are maps over the identity of $\Lambda^2_1(\huaG)$.
    \begin{equation}\label{eq:VB2Dual1Gpd-lambda-horn-lifts}
    \begin{array}{ccc}
    \begin{aligned}
        \llambda_0: (d_1d_2)^*C_R|_{(g,h)} &\to \Lambda^2_1(\huaV)|_{(g,h)}\\
        a'|_{d_1h} &\mapsto (0_g, a' \cdot 0_h),
    \end{aligned}
    & &
    \begin{aligned}
        \llambda_1: (d_0d_2)^*C_L&|_{(g,h)} \to \Lambda^2_1(\huaV)|_{(g,h)}\\
        b'|_{d_1g = d_0h} &\mapsto ((b')^{-1} \cdot 0_g, 0_h \cdot b'),
    \end{aligned}
    \\
    & & \\
    \begin{aligned}
        \llambda_2: (d_0d_1)^*C_L|_{(g,h)} &\to \Lambda^2_1(\huaV)|_{(g,h)}\\
        b|_{d_0g} &\mapsto (0_g \cdot b, 0_h).
    \end{aligned}
    & &
    \end{array}
    \end{equation}
    
    \begin{center}
    \begin{adjustbox}{width=\textwidth}
    \begin{tikzpicture}[scale=0.6]
	\begin{pgfonlayer}{nodelayer}
		\node [style=dot, label={below left:\small{0}}] (3) at (0, 0) {};
		\node [style=white-dot, label={below right:\small{2}}] (4) at (3, 0) {};
		\node [style=white-dot, label={above right:\small{1}}] (5) at (1.5, 2.5) {};
		\node [style=dot, label={below left:\small{0}}] (18) at (7.5, 0) {};
		\node [style=dot, label={below right:\small{2}}] (19) at (10.5, 0) {};
		\node [style=dot, label={above right:\small{1}}] (20) at (9, 2.5) {};
		\node [style=dot, label={below left:\small{0}}] (22) at (15, 0) {};
		\node [style=dot, label={below right:\small{2}}] (23) at (18, 0) {};
		\node [style=dot, label={above right:\small{1}}] (24) at (16.5, 2.5) {};
		\node [style=none] (25) at (1.5, -1) {\small{$a'\cdot 0_{hg}$}};
		\node [style=none] (26) at (0, 1.75) {\small{$a'$}};
		\node [style=none] (27) at (3.25, 1.5) {\small{$0_g$}};
		\node [style=none] (28) at (9, -1) {\small{$0_{hg}$}};
		\node [style=none] (29) at (7, 1.5) {\small{$0_h\cdot b'$}};
		\node [style=none] (30) at (11.5, 1.5) {\small{$(b')^{-1}\cdot 0_g$}};
		\node [style=none] (31) at (16.5, -1) {\small{$0_{hg}\cdot b$}};
		\node [style=none] (32) at (14.75, 1.5) {\small{$0_h$}};
		\node [style=none] (33) at (18.5, 1.5) {\small{$0_g\cdot b$}};
		\node [style=none] (34) at (1.5, -2.5) {\small{$\llambda_0(a')$}};
		\node [style=none] (35) at (9, -2.5) {\small{$\llambda_1(b')$}};
		\node [style=none] (36) at (16.5, -2.5) {\small{$\llambda_2(b)$}};
	\end{pgfonlayer}
	\begin{pgfonlayer}{edgelayer}
		\draw [style=directed] (5) to (3);
		\draw [style=directed-dash] (4) to (5);
		\draw [style=directed, in=360, out=180] (4) to (3);
		\draw [style=directed] (20) to (18);
		\draw [style=directed] (19) to (20);
		\draw [style=directed-dash] (19) to (18);
		\draw [style=directed-dash] (24) to (22);
		\draw [style=directed] (23) to (24);
		\draw [style=directed] (23) to (22);
	\end{pgfonlayer}
    \end{tikzpicture}
    \end{adjustbox}
    \end{center}
    Now consider an arbitrary tetrahedron in $\huaV^{2*}$.
    We write this as a boundary and denote it by
    \begin{equation*}
        ((\phi_0, \alpha_0, \beta_0, \gamma_0), (\phi_1, \alpha_1, \beta_1, \gamma_1), (\phi_2, \alpha_2, \beta_2, \gamma_2), (\phi_3, \alpha_3, \beta_3, \gamma_3)) \in \partial^3(\huaV^{2*}) \cong \huaV_3,
    \end{equation*}
    where the components obey the boundary conditions\footnote{An unfortunate notational overlap becomes evident at this point, where $\gamma$ is used to denote both the fourth component of a generic element of $\huaV^{2*}$, and in denoting the core projections $\boldsymbol{\upgamma}_R, \boldsymbol{\upgamma}_L$. We trust that the meaning of each $\gamma$ will be clear from context and the slightly bolder font used for the core projections. }
    \begin{equation}\label{eq:VB2Dual1Gpd-3-horn-conditions}
        \begin{array}{lll}
            \begin{aligned}
                &\alpha_0 = \alpha_1,\\
                &\beta_0 = \beta_1,
            \end{aligned}
            &
            \begin{aligned}
                &\alpha_2 = \gamma_0 + \llambda_0^*\phi_0,\\
                &\beta_2 = \beta_0 + \llambda_2^*\phi_0,
            \end{aligned}
            &
            \begin{aligned}
                &\gamma_2 + \llambda_0^*\phi_2 = \gamma_1 + \llambda_0^*\phi_1,\\
                &\beta_2 + \llambda_2^*\phi_2 = \beta_1 + \llambda_2^*\phi_1,
            \end{aligned}
            \\
            & & \\
            \begin{aligned}
                &\alpha_3 = \gamma_0,\\
                &\beta_3 = \ggamma_R^*\eta_0 - \llambda_1^*\phi_0,
            \end{aligned}
            &
            \begin{aligned}
                &\gamma_3 + \llambda_0^*\phi_3 = \gamma_1,\\
                &\beta_3 + \llambda_2^*\phi_3 = \ggamma_R^*\eta_3 - \llambda_1^*\phi_1,
            \end{aligned}
            &
            \begin{aligned}
                &\gamma_3 = \gamma_2,\\
                &\boldsymbol{\upgamma}_R^*\eta_3 - \llambda_1^*\phi_3 = \boldsymbol{\upgamma}_R^*\eta_2 - \boldsymbol{\uplambda}_1^*\phi_2.
            \end{aligned}
            \\
        \end{array}
    \end{equation}
    For any given $(3,j)$-horn in $\huaV$, its faces obey a subset of these boundary conditions. The $j$-th face is determined by the multiplications 
    \begin{align}\label{eq:VB2Dual1Gpd-def-multiplications}
        \square (\phi_1, \alpha_1, \beta_1, \gamma_1)(\phi_2, \alpha_2, \beta_2, \gamma_2)(\phi_3, \alpha_3, \beta_3, \gamma_3) &= (\phi_0, \alpha_1, \beta_1, \alpha_3),\\
        (\phi_0, \alpha_0, \beta_0, \gamma_0)\square(\phi_2, \alpha_2, \beta_2, \gamma_2)(\phi_3, \alpha_3, \beta_3, \gamma_3) &= (\phi_1, \alpha_0, \beta_0, \gamma_3 + \llambda_0^* \phi_3),\notag\\
        (\phi_0, \alpha_0, \beta_0, \gamma_0)(\phi_1, \alpha_1, \beta_1, \gamma_1)\square(\phi_3, \alpha_3, \beta_3, \gamma_3) &= (\phi_2, \gamma_0 + \llambda_0^*\phi_0, \beta_0 + \llambda_2^*\phi_0, \gamma_3),\notag\\
        (\phi_0, \alpha_0, \beta_0, \gamma_0)(\phi_1, \alpha_1, \beta_1, \gamma_1)(\phi_2, \alpha_2, \beta_2, \gamma_2)\square &= (\phi_3, \gamma_0, \ggamma_R^*\eta_0 - \llambda_1^*\phi_0, \gamma_2),\notag
    \end{align}
    where the $012$ components satisfy the equation 
    \begin{equation}\label{eq:VB2Dual1Gpd-def-multiplications-intcomp}
        \phi_0|_{(g,h)} (v, w) - \phi_1|_{(g,kh)}(v, u \cdot w) + \phi_2|_{(hg,k)}(w \cdot v,u) - \phi_3|_{(k,h)}(w,u) = 0,
    \end{equation}
    for any tetrahedron $((v,w),(v, u \cdot w),(w \cdot v, u),(w,u)) \in \partial^3(\huaV)\cong\huaV_3$ over the tetrahedron $((g,h),(g, k \cdot h),(h \cdot g, k),(h,k)) \in \partial^3(\huaG)\cong\huaG_3$, corresponding to the triple of composable arrows $(v,w,u)$ in $\huaV$ over $(g,h,k)$ in $\huaG$, as in Example \ref{ex:nerve-cat}.

    It is possible to obtain an expression of the 2-dual in terms of just $C_R$ or just $C_L$, just by composing with the isomorphisms above where appropriate, but we preferred to avoid having to keep track of them. 
\end{example}

\chapter{2-cotangents and shifted symplectic structures}

In this short chapter we study the 2-cotangent of a Lie 2-groupoid. Most importantly, we show that it carries a natural 2-shifted symplectic structure. 

When the Lie 2-groupoid is given by the nerve of a Lie groupoid $\huaG$, we show how the 2-cotangent $T^{2*}\huaG$ is Morita equivalent to the 2-groupoid $\widebar{T^{1*}\huaG}$ obtained by using the Artin-Mazur bar construction on the 1-cotangent of $\huaG$, as we show in Theorem \ref{thm:2Cota-1Cota-ME}. The Morita equivalence we construct is also compatible with the 2-shifted symplectic structures, as we show in Theorem \ref{thm:2Cota-1Cota-SME}.

The symplectic 2-groupoid $\widebar{T^{1*}\huaG}$ was proposed in \cite{MehtaTang2011} as the global object integrating the Courant algebroid $A \oplus A^*$ (see \cite{LiuWeinsteinXu1997}). This is meant in the sense that the tangent complex of $\widebar{T^{1*}\huaG}$ is isomorphic (up to a splitting) to the chain complex 
\begin{equation}\label{eq:double-of-A-courant-complex}
    T^*\huaG_0 \to A \oplus A^* \to T\huaG_0,
\end{equation}
which is the chain complex obtained from $A\oplus A^*$ by taking the anchor and its dual as differentials. The integration problem for Courant algebroids has been approached in many ways since Ševera put forward the idea that Courant algebroids should integrate to some version of symplectic Lie 2-groupoids in \cite{Severa2005}. Since the integration proposed there produces infinite-dimensional groupoids, the search for finite-dimensional models of integrations has produced various results in different contexts. We point the interested reader to the introductions of \cite{MehtaTang2018a} and the more recent \cite{CuecaZhu2023}, which put the approaches in \cite{LiBlandSevera2012, MehtaTang2011, MehtaTang2018,Severa2005, SeveraSiran2019, ShengZhu2017} in context. 

Under the correspondence in \cite{Roytenberg2002}, the Courant algebroid $A \oplus A^*$ corresponds to the 2-shifted cotangent degree 2 dg-manifold $T^*[2]A[1]$ of the Lie algebroid $A$ considered as a degree 1 dg-manifold. By the symplectic Morita equivalence in \ref{thm:2Cota-1Cota-SME}, the tangent complex of the 2-cotangent $T^{2*}\huaG$ is quasi-isomorphic to the tangent complex of $\widebar{T^{1*}\huaG}$, which is \eqref{eq:double-of-A-courant-complex}. In this sense we can think of the 2-cotangent of a Lie groupoid $\huaG$ to be an \textit{integration up to homotopy} of the 2-shifted cotangent of its Lie algebroid. This is to say that the operation of taking the 2-shifted cotangent of a Lie algebroid corresponds up to homotopy to the operation of taking the 2-cotangent of a Lie groupoid, under a (partial) analog of the Lie correspondence.

\section{The 2-cotangent of a Lie 2-groupoid}   

Let $\huaG$ be a Lie $2$-groupoid. It was shown in Example \ref{ex:tangent-groupoid} that $T\huaG\to \huaG$ defines a $\VB$ 2-groupoid with all the structural maps given by the tangent lifts of the maps on $\huaG$. Using the duality introduced in the previous sections we can construct a $\VB$ 2-groupoid dual to this.

\begin{definition}
Let $\huaG$ be a Lie $2$-groupoid. The \textbf{2-cotangent} $T^{2*}\huaG\to \huaG$ of $\huaG$ is the $\VB$ 2-dual of the tangent $T\huaG\to \huaG$. That is 
\begin{equation*}
\begin{adjustbox}{width=\textwidth}
\begin{tikzcd}[ampersand replacement=\&,column sep=small,row sep=scriptsize]
	{\begin{array}{c}T^*\huaG_2\times_{i^*, (\ker Tp^2_0)^*, \nu_0^*} d_0^*\left(O_0\times_{i^*, (s_0^*\ker Tp^2_1)^*, \zeta^*} O_1\right)\\ \times_{i^*, (\ker Tp^2_2)^*, \nu_2^*} d^*_2 O_0\end{array}} \& {O_0\times_{i^*, (s_0^*\ker Tp^2_1)^*, \zeta^*} O_1} \& {O_{01}} \\
	{\huaG_2} \& {\huaG_1} \& {\huaG_0}
	\arrow[shift left=2, from=1-1, to=1-2]
	\arrow[from=1-1, to=1-2]
	\arrow[shift right=2, from=1-1, to=1-2]
	\arrow[shift left=1, from=1-2, to=1-3]
	\arrow[shift right=1, from=1-2, to=1-3]
	\arrow[from=1-2, to=2-2]
	\arrow[shift left=1, from=2-2, to=2-3]
	\arrow[shift right=1, from=2-2, to=2-3]
	\arrow[shift left=2, from=2-1, to=2-2]
	\arrow[shift right=2, from=2-1, to=2-2]
	\arrow[from=2-1, to=2-2]
	\arrow[from=1-1, to=2-1]
	\arrow[from=1-3, to=2-3]
\end{tikzcd}
\end{adjustbox}
\end{equation*}
where 
\begin{equation*}
\begin{split}
    O_0 = \Ann(Ts_0T\huaG_1&) \subseteq s_0^*T^*\huaG_2, \quad O_1 = \Ann(Ts_1T\huaG_1)\subseteq s_1^*T^*\huaG_2,\\
    &O_{01} = \Ann(D_2T\huaG) \subseteq 1^*T^*\huaG_2,
\end{split}   
\end{equation*}
and the fiber products and simplicial structure maps are as in Definition \ref{def:VB2dual}.
\end{definition}

We propose the 2-cotangent as the canonical model for a 2-shifted symplectic 2-groupoid. To show that our candidate canonical symplectic form is homologically 2-shifted nondegenerate we need to study the tangent complex of $T^{2*}\huaG$. For this, we will make use of the structure of simplicial double vector bundle of the tangent $TT^{2*}\huaG$. 

\subsection{Simplicial double vector bundles}

\begin{definition}
    A \textbf{simplicial double vector bundle} (simplicial DVB for short), is a double object in the category of simplicial vector bundles. That is, a simplicial DVB over a simplicial manifold $\huaG$ is a triple $(\huaD, \huaA, \huaB)$ of simplicial manifolds, where $\huaA$ and $\huaB$ are simplicial vector bundles over $\huaG$, and $\huaD$ has two compatible simplicial vector bundle structures over $\huaA$ and $\huaB$, respectively. The compatibility condition is that for each $m$, $(\huaD_m, \huaA_m, \huaB_m)$ is a DVB over $\huaG_m$ and the simplicial maps of $\huaD$ are double vector bundle maps. We represent $(\huaD, \huaA, \huaB)$ by 
    \[\begin{tikzcd}[ampersand replacement=\&,cramped,sep=scriptsize]
	\huaD \& \huaB \\
	\huaA \& \huaG
	\arrow["{q_\huaB^\huaD}", from=1-1, to=1-2]
	\arrow["{q_\huaA^\huaD}"', from=1-1, to=2-1]
	\arrow["{q_\huaB}", from=1-2, to=2-2]
	\arrow["{q_\huaA}"', from=2-1, to=2-2]
    \end{tikzcd}\]
    The \textbf{core} of a simplicial DVB $(\huaD, \huaA, \huaB)$ is the simplicial vector bundle $\huaC \to \huaG$ defined at each level $m \ge 0$ as the core of $(\huaD_m, \huaA_m, \huaB_m)$, i.e.
    \begin{equation*}
        \huaC_m := \ker q_{\huaA_m}^{\huaD_m} \cap \ker q_{\huaB_m}^{\huaD_m},
    \end{equation*}
    with the simplicial maps induced from $\huaD$.
\end{definition}

We will now focus on the structure of $\huaD$ as a simplicial vector bundle over $\huaA$ and drop the subscripts on the relevant projections, to denote $q:= q_\huaA$ and $\widetilde{q}:= q_\huaB^\huaD$. 
Recall (e.g. from \cite[Thm. 9.1.6]{Mackenzie2005}) that for any DVB $(\huaD_m, \huaA_m, \huaB_m)$, there is an exact sequence
\[\begin{tikzcd}[ampersand replacement=\&,cramped, sep=scriptsize]
	0 \& {q_{m}^*\huaC_m} \& {\huaD_m} \& {q_{m}^*\huaB_m} \& 0
	\arrow[from=1-2, to=1-3]
	\arrow["{\widetilde{q}}", from=1-3, to=1-4]
	\arrow[from=1-1, to=1-2]
	\arrow[from=1-4, to=1-5]
\end{tikzcd}\]
of vector bundles over $\huaA_m$. By taking this sequence at every level $m$ of a simplicial DVB, since both maps in the sequence are simplicial morphisms, we get an exact sequence 
\[\begin{tikzcd}[ampersand replacement=\&,cramped, sep=scriptsize]
	0 \& {q^*\huaC} \& \huaD \& {q^*\huaB} \& 0
	\arrow[from=1-2, to=1-3]
	\arrow["{\widetilde{q}}", from=1-3, to=1-4]
	\arrow[from=1-1, to=1-2]
	\arrow[from=1-4, to=1-5]
\end{tikzcd}\]
of simplicial vector bundles over $\huaA$.
The pullback of this sequence by the total units $(1_\huaA, 1)$ of $\huaA \to \huaG$ is an exact sequence 
\[\begin{tikzcd}[ampersand replacement=\&,cramped,sep=scriptsize]
	0 \& {q^*1^*\huaC} \& {1_\huaA^*\huaD} \& {q^*1^*\huaB} \& 0
	\arrow["{\widetilde{q}}", from=1-3, to=1-4]
	\arrow[from=1-1, to=1-2]
	\arrow[from=1-4, to=1-5]
	\arrow[from=1-2, to=1-3]
\end{tikzcd}\]
of simplicial vector bundles over $\huaA_0$, as it is the core sequence of the simplicial DVB $(1_\huaA^*\huaD, \huaA_0, 1^*\huaB)$ over $\huaG_0$ (where $\huaA_0$ and $\huaG_0$ denote the respective identity 0-groupoids). This immediately translates to an exact sequence of chain complexes of vector bundles over $\huaA_0$ between the Moore complexes of each simplicial vector bundle. By applying the snake lemma repeatedly, and noting that the face maps and horn projections are surjective submersions (so their kernels are subbundles and their cokernels vanish), we obtain an exact sequence 
\begin{equation}\label{eq:SDVB-coreseq-norm-cplx}
\begin{tikzcd}[ampersand replacement=\&,cramped,sep=scriptsize]
0 \& {q^*N(\huaC)} \& {N(\huaD)} \& {q^*N(\huaB)} \& 0
\arrow[from=1-2, to=1-3]
\arrow["{N(\widetilde{q})}",from=1-3, to=1-4]
\arrow[from=1-1, to=1-2]
\arrow[from=1-4, to=1-5]
\end{tikzcd}
\end{equation}
of chain complexes of vector bundles over $\huaA_0$ as well. 

\begin{example}\label{ex:TT2starG-as-DVB}
    An immediate example of a simplicial DVB is the tangent bundle of a simplicial vector bundle. 
    We will focus on the tangent 2-groupoid $TT^{2*}\huaG$ of the 2-cotangent $T^{2*}\huaG$, which fits in the diagram 
    \[\begin{tikzcd}[ampersand replacement=\&,cramped,sep=scriptsize]
	{TT^{2*}\huaG} \& T\huaG \\
	{T^{2*}\huaG} \& \huaG
	\arrow["Tq", from=1-1, to=1-2]
	\arrow[from=1-1, to=2-1]
	\arrow[from=1-2, to=2-2]
	\arrow["q"', from=2-1, to=2-2]
    \end{tikzcd}\]
    It is a well-known fact in the theory of double vector bundles (see e.g. \cite[Ex. 9.1.7, \S 3.4]{Mackenzie2005}) that the core of the tangent DVB of a vector bundle is isomorphic to the vector bundle itself, via the vertical lift map. We recall this isomorphism explicitly at each simplicial level and show that it induces a simplicial isomorphism between $q^*C$ and $q^*T^{2*}\huaG$.
    At each fixed level $m$ and each point $\phi_g \in T^{2*}\huaG_m|_g$, we have that $q^*C_m|_{\phi_g} = \ker Tq|_{\phi_g} \subseteq T_{\phi_g}T^{2*}\huaG$,  the vertical bundle at $\phi_g$, which is the tangent space to the fiber $T^{2*}\huaG_m|_g$ at $\phi_g$. Hence, the isomorphism $q^*C_m|_{\phi_g}\cong T^{2*}G_m|_g$ is given by the vertical lift at $\phi_g$. 
    For any $\psi_g \in T^{2*}G_m|_g$, this is defined by the vector $\psi_g^\uparrow|_{\phi_g} \in T_{\phi_g}T^{2*}\huaG_m$ such that, for any $f\in C^{\infty}(T^{2*}\huaG_m)$,
    \begin{equation*}
        \psi_g^\uparrow|_{\phi_g} \cdot f := \left.\frac{d}{dt}\right|_{t=0}f(\phi_g + t\psi_g).
    \end{equation*}
    Notice that vertical lifts commute with the simplicial structure maps: for example, for any face map $\widecheck{d}_i$ of $T^{2*}\huaG$, we have that
    \begin{equation*}
        T_{\phi_g}\widecheck{d}_i(\psi_g^\uparrow|_{\phi_g}) \cdot f = \left.\frac{d}{dt}\right|_{t=0}f(\widecheck{d}_i\phi_g + t\widecheck{d}_i\psi_g) = (\widecheck{d}_i(\psi_g))^{\uparrow}|_{\widecheck{d}_i\phi_g}.
    \end{equation*}
    The same computation works in the case of the degeneracy maps. So the vertical lift is a simplicial bundle map, and we can write the core sequence of $TT^{2*}\huaG$ over $T^{2*}\huaG$ as the exact sequence
    \begin{equation}\label{diag:VB2CotangentTangentCoreSeq}
    \begin{tikzcd}[ampersand replacement=\&,cramped,sep=scriptsize]
	0 \& {q^*T^{2*}\huaG} \& {TT^{2*}\huaG} \& {q^*T\huaG} \& 0
	\arrow["Tq", from=1-3, to=1-4]
	\arrow[from=1-1, to=1-2]
	\arrow[from=1-4, to=1-5]
	\arrow["{(\_)^{\uparrow}}", from=1-2, to=1-3]
    \end{tikzcd}
    \end{equation}
    of vector bundles over $T^{2*}\huaG$.
\end{example}

\begin{lemma}\label{lem:NormCplxSDVB}
    Let $(\huaD, \huaA, \huaB)$ be a simplicial double vector bundle over a simplicial manifold $\huaG$ with core $\huaC \to \huaG$. Then the normalized complex $N(\huaD) \to \huaA_0$ of $\huaD$ as a simplicial vector bundle over $\huaA \overset{q}{\to} \huaG$ fits into the following short exact sequence of chain complexes of vector bundles over $A_0$.
    \[\begin{tikzcd}[ampersand replacement=\&,cramped,sep=scriptsize]
	0 \& {q^*N(\huaC)} \& {N(\huaD)} \& {q^*N(\huaB)} \& 0
	\arrow[from=1-2, to=1-3]
	\arrow[from=1-3, to=1-4]
	\arrow[from=1-1, to=1-2]
	\arrow[from=1-4, to=1-5]
    \end{tikzcd}\]

    Consequently, consider two simplicial DVBs $(\huaD, \huaA, \huaB)$ and $(\huaD', \huaA, \huaB')$, with the same side $\huaA \overset{q}{\to} \huaG$. Let $f: N(\huaD) \to N(\huaD')$ be a chain map between their normalized complexes over $A$ which covers the identity on $\huaA_0$. If the following diagram commutes, and if $f^{\huaC}$ and $f^{\huaB}$ are quasi-isomorphisms at each point $a_p \in \huaA_0$, then $f$ is a quasi-isomorphism at each point $a_p \in \huaA_0$. 
    
    \begin{equation}\label{diag:NormCplxSDVBMap}    
    \begin{tikzcd}[ampersand replacement=\&,cramped,sep=scriptsize]
	0 \& {q^*N(\huaC)|_{a_p}} \& {N(\huaD)|_{a_p}} \& {q^*N(\huaB)|_{a_p}} \& 0 \\
	0 \& {q'^*N(\huaC')|_{a_p}} \& {N(\huaD')|_{a_p}} \& {q'^*N(\huaB')|_{a_p}} \& 0
	\arrow[from=1-2, to=1-3]
	\arrow[from=1-3, to=1-4]
	\arrow[from=1-1, to=1-2]
	\arrow[from=1-4, to=1-5]
	\arrow["{f^{\huaC}}"', from=1-2, to=2-2]
	\arrow["{f^{\huaB}}", from=1-4, to=2-4]
	\arrow[from=2-2, to=2-3]
	\arrow[from=2-3, to=2-4]
	\arrow["f"', from=1-3, to=2-3]
	\arrow[from=2-1, to=2-2]
	\arrow[from=2-4, to=2-5]
    \end{tikzcd}
    \end{equation}
\end{lemma}

\begin{proof}

    By the previous discussion, the core sequence for simplicial DVBs gives a short exact sequence \eqref{eq:SDVB-coreseq-norm-cplx} at the level of normalized complexes over $\huaA_0$. If \eqref{diag:NormCplxSDVBMap} commutes, then, by naturality of the long exact sequence in homology, we get, at each point $a_p\in \huaA_0$, and at each $i$, the commutative diagram
    \[
    \begin{adjustbox}{width=\textwidth}
    \begin{tikzcd}[ampersand replacement=\&,cramped, column sep=small, row sep=scriptsize]
	\dots \& {q^*H_{i+1}(\huaB)|_{a_p}} \& {q^*H_i(\huaC)|_{a_p}} \& {H_i(\huaD)|_{a_p}} \& {q^*H_i(\huaB)|_{a_p}} \& {q^*H_{i-1}(\huaC)|_{a_p}} \& \dots \\
	\dots \& {q'^*H_{i+1}(\huaB')|_{a_p}} \& {q'^*H_i(\huaC')|_{a_p}} \& {H_i(\huaD')|_{a_p}} \& {q'^*H_i(\huaB')|_{a_p}} \& {q'^*H_{i-1}(\huaC')|_{a_p}} \& \dots
	\arrow[from=1-3, to=1-4]
	\arrow[from=1-4, to=1-5]
	\arrow[from=1-2, to=1-3]
	\arrow[from=1-5, to=1-6]
	\arrow["{f^{\huaC}}"', from=1-3, to=2-3]
	\arrow["{f^{\huaB}}"', from=1-5, to=2-5]
	\arrow[from=2-3, to=2-4]
	\arrow[from=2-4, to=2-5]
	\arrow["f"', from=1-4, to=2-4]
	\arrow[from=2-2, to=2-3]
	\arrow[from=2-5, to=2-6]
	\arrow["{f^{\huaB}}"', from=1-2, to=2-2]
	\arrow[from=1-1, to=1-2]
	\arrow[from=2-1, to=2-2]
	\arrow["{f^{\huaC}}"', from=1-6, to=2-6]
	\arrow[from=2-6, to=2-7]
	\arrow[from=1-6, to=1-7]
    \end{tikzcd}
    \end{adjustbox}\]
    where the $f^{\huaC}$ and $f^{\huaB}$ are isomorphisms on the homology. Hence, by the five lemma, $f$ is an isomorphism on the homology at each level $i$.
\end{proof}

\subsection{The canonical shifted symplectic form}

We now come to the canonical 2-shifted symplectic form on $T^{2*}\huaG$.

\begin{definition}\label{def:2Cota-taut-form-def}
    Let $\huaG$ be a Lie $2$-groupoid, and $T^{2*}\huaG$ its 2-cotangent. The \textbf{tautological 1-form} $\Theta \in \Omega^1(T^{2*}\huaG_2)$ is defined at each point $\phi_t \in T^{2*}\huaG_2$, by 
    \begin{equation}\label{eq:2Cota-taut-form-def}
        \Theta(V_{\phi_t}) = \langle \phi_t, Tq V_{\phi_t} \rangle = \phi^{012}_t(Tq V_{\phi_t}),
    \end{equation}
    for any $V_{\phi_t}\in TT^{2*}\huaG_2$, with $\langle \cdot, \cdot \rangle$ the 2-dual pairing, and $q: T^{2*}\huaG_2 \to \huaG_2$, the bundle projection.
\end{definition}

\begin{remark}
    Let $\mathbf{p}: T^{2*}\huaG_2 \to T^*\huaG_2$ be the projection to the $012$ component of each element $\phi$. Then $\Theta = \mathbf{p}^*\theta_{can}$, where $\theta_{can}$ is the tautological form of $T^*\huaG_2$. 

    Analogously, the canonical shifted symplectic form we introduce in the next proposition is exactly $\mathbf{p}^*\omega_{can} = -d \mathbf{p}^*\theta_{can}$, where $\omega_{can}$ is the canonical symplectic form on the cotangent bundle $T^*\huaG_2$.
\end{remark}

\begin{theorem}\label{thm:2Cota-symp-form}
    Let $\huaG$ be a Lie $2$-groupoid, and $T^{2*}\huaG$ its 2-cotangent. The 2-form $\omega := -d\Theta \in \Omega^2(T^{2*}\huaG)$ is a 2-shifted symplectic form. We will refer to $\omega$ as the \textbf{canonical symplectic form} on the 2-cotangent of $\huaG$. 
\end{theorem}

\begin{proof}
    Since $\omega$ is even exact, it is closed with respect to the de Rham differential. 
    
    We now show it is multiplicative, i.e. $\delta\omega = 0$. 
    First of all, $\delta\omega = - \delta d \Theta = - d \delta \Theta$. Furthermore, for any tetrahedron $(V_0, V_1, V_2, V_3) \in TT^{2*}\huaG_3|_{\lambda_{(r,s,t,u)}}$, over $\lambda_t = (\phi_r, \chi_s, \psi_t, \zeta_u) \in T^{2*}\huaG_3$, we have that, since $q$ is a simplicial map, $Tq$ is a simplicial vector bundle map, and 
    \begin{equation*}
        Tq V_1 = Tq(V_0 \square V_2 V_3) = (TqV_0) \square (TqV_2)(TqV_3).
    \end{equation*}
    Therefore, 
    \begin{equation*}
        \begin{split}
        (\delta\Theta)_{\lambda_{(r,s,t,u)}} (V_0, V_1, V_2, V_3) &= 
        \langle \phi_r, Tq V_0 \rangle 
        - \langle \chi_s, Tq V_1 \rangle 
        + \langle \psi_t, Tq V_2 \rangle 
        - \langle \zeta_u, Tq V_3 \rangle\\
        &= \langle \phi_r, Tq V_0 \rangle 
        - \langle \phi_r \square \psi_s \zeta_u, (TqV_0)\square (TqV_2)(TqV_3) \rangle \\
        &\qquad + \langle \psi_t, Tq V_2 \rangle 
        - \langle \zeta_u, Tq V_3 \rangle =0,
        \end{split} 
    \end{equation*}
    by multiplicativity of the 2-dual pairing. Hence $\delta\Theta =0$, which implies $\delta\omega =0$, as well.

    Since $(T\widecheck{s}_i)^*\omega = - (T\widecheck{s}_i)^*d \Theta = - d(T\widecheck{s}_i)^*\Theta$, for $i =0,1$, if $\Theta$ is normalized, then $\omega$ is also normalized. But for any $V \in TT^{2*}\huaG_1|_{(\eta, \xi)_g}$, we have
    \begin{equation*}
        \begin{split}
            ((T\widecheck{s}_0)^*\Theta)_{(\eta, \xi)_g}(V)
            &= \langle \widecheck{s}_0(\eta, \xi)_g, Tq T\widecheck{s}_0 V \rangle \\
            &= \langle \widecheck{s}_0(\eta, \xi)_g, Ts_0 Tq V \rangle 
            = \eta(Ts_0 Tq V) = 0,
        \end{split}
    \end{equation*}
    since $Tq$ is a simplicial vector bundle map $TT^{2*}\huaG \to T\huaG$ and $\eta \in \Ann(Ts_0 TG_1)$. Similarly, $(T\widecheck{s}_1)^*\Theta = 0$, since $\xi \in \Ann(Ts_1 TG_1)$, and $\langle \widecheck{s}_1 (\eta, \xi), \_ \rangle = \xi$. 

    Finally, we show that $\omega$ is homologically 2-shifted nondegenerate as a pairing $TT^{2*}\huaG \wedge TT^{2*}\huaG \to B^2\R^{\huaG}$. As shown in Example \ref{ex:TT2starG-as-DVB}, $TT^{2*}\huaG$ is a simplicial double vector bundle, with core sequence as in \eqref{diag:VB2CotangentTangentCoreSeq}. 
    Therefore, by Lemma \ref{lem:NormCplxSDVB}, the tangent complex $\huaT(T^{2*}\huaG)$ of the 2-cotangent $\VB$ 2-groupoid fits in the short exact sequence 
    \[\begin{tikzcd}[ampersand replacement=\&,cramped,sep=small]
        0 \& {q^*N(T^{2*}\huaG)} \& {\huaT(T^{2*}\huaG)} \& {q^*\huaT(\huaG)} \& 0 \\
        \arrow[from=1-2, to=1-3]
        \arrow[from=1-3, to=1-4]
        \arrow[from=1-1, to=1-2]
        \arrow[from=1-4, to=1-5]
    \end{tikzcd}\]   
    of complexes of vector bundles over $T^{2*}\huaG_0$.

    We now claim that the map $\lambda_\omega^r$ induced by the associated $IM$-pairing $\lambda_\omega$ fits in the commutative diagram 
    \begin{equation}\label{diag:VB2CotangentOmegaTangentCommDiag}
    \begin{tikzcd}[ampersand replacement=\&,cramped, sep=scriptsize]
	0 \& {q^*N(T^{2*}\huaG)} \& {\huaT(T^{2*}\huaG)} \& {q^*\huaT(\huaG)} \& 0 \\
	0 \& {q^*\huaT(\huaG)^*[-2]} \& {\huaT(T^{2*}\huaG)^*[-2]} \& {q^*N(T^{2*}\huaG)^*[-2]} \& 0
	\arrow["{(\_)^{\uparrow}}", from=1-2, to=1-3]
	\arrow["Tq", from=1-3, to=1-4]
	\arrow[from=1-1, to=1-2]
	\arrow[from=1-4, to=1-5]
	\arrow["{-q^*\lambda_{\langle,\rangle}^r}", from=1-2, to=2-2]
	\arrow["{\lambda_\omega^r}", from=1-3, to=2-3]
	\arrow["{q^*\lambda_{\langle,\rangle}^l}", from=1-4, to=2-4]
	\arrow["{}"', from=2-2, to=2-3]
	\arrow[from=2-1, to=2-2]
	\arrow[from=2-4, to=2-5]
	\arrow["{}"', from=2-3, to=2-4]
    \end{tikzcd}
    \end{equation}
    of chain complexes over $T^{2*}\huaG_0$, where $\lambda_{\langle,\rangle}^r$ and $\lambda_{\langle,\rangle}^l$ are the two maps induced by the IM-pairing associated to the 2-dual pairing, and the bottom row is the 2-shifted dual of the top row. With this, since $\lambda_{\langle,\rangle}^r$ and $\lambda_{\langle,\rangle}^l$ are (pointwise) quasi-isomorphisms by Theorem \ref{thm:VBndual-pairing-hndg}, we are in the situation of Lemma \ref{lem:NormCplxSDVB}. Therefore $\lambda_\omega^r$ is a quasi-isomorphism at each point and $\omega$ is homologically 2-shifted nondegenerate. 

    We now show commutativity of \eqref{diag:VB2CotangentOmegaTangentCommDiag}. 
    To show this we need to understand how $\omega$ acts on vertical lifts. 
    A similar computation to the one below can also be found in \cite[Appendix A.2]{JotzLean2018}. 
    Consider a splitting $\sigma$ of the sequence 
    \eqref{diag:VB2CotangentTangentCoreSeq} at level 2, that is 
    \[\begin{tikzcd}[ampersand replacement=\&,cramped]
	0 \& {q^*T^{2*}\huaG_2} \& {TT^{2*}\huaG_2} \& {q^*T\huaG_2} \& 0
	\arrow["Tq", from=1-3, to=1-4]
	\arrow[from=1-1, to=1-2]
	\arrow[from=1-4, to=1-5]
	\arrow["{(\_)^{\uparrow}}", from=1-2, to=1-3]
	\arrow["{\sigma}", curve={height=-12pt}, from=1-4, to=1-3]
    \end{tikzcd}\]
    where $Tq\sigma = id$. This is equivalent to the choice of a linear (Ehresmann) connection on $T^{2*}\huaG_2$. Since horizontal lifts $\sigma(X)$ are linear vector fields for any $X\in \VF(\huaG_2)$, this allows one to see vector fields over $T^{2*}\huaG_2$ as generated over $C^{\infty}(T^{2*}\huaG_2)$ by core vector fields (vertical lifts) and linear vector fields (horizontal lifts). 
    Hence we are able to characterize $\Theta$, $\omega$ and the Lie bracket by their action on vertical and horizontal lifts. 
    Consider $V \in \VF(T^{2*}\huaG_2)$ such that $TqV = X \in \VF(\huaG_2)$ is a well-defined vector field. 
    This can be written as $V = \psi^\uparrow + \sigma(X)$, for some $\psi \in \Gamma(T^{2*}\huaG_2)$, and we have, at any point $\phi_t \in T^{2*}\huaG_2$
    \begin{equation}\label{eq:VB2CotangentTautFormCoreLinear}
        \begin{split}
            \Theta(\sigma(X))|_{\phi_t} &= \langle \phi_t , Tq \sigma (X) \rangle = \langle \phi_t, X \rangle,\\
            \Theta(\psi^\uparrow)|_{\phi_t} &= \langle \phi_t , Tq \psi^\uparrow \rangle = 0.
        \end{split}
    \end{equation}
    In terms of functions in $C^\infty(T^{2*}\huaG_2)$, then, $\Theta(\psi^\uparrow)$ is identically $0$, while $\Theta(\sigma(X)) = ev_{X}$, the linear function given by the 2-dual evaluation on $X$. Now, by definition of $\omega$ we have
    \begin{equation*}
        \omega(V,W) = -(d\Theta)(V,W) = - \huaL_{V}\Theta(W) + \huaL_{W}\Theta(V) + \Theta([V,W]),
    \end{equation*}
    for any $V,W \in \VF(T^{2*}\huaG_2)$. Hence, for $X=TqV$,
    \begin{equation}\label{eq:VB2CotangentOmegaOnVertLift}
        \begin{split}
            \omega(V, \chi^\uparrow) = \omega(\psi^\uparrow, \chi^\uparrow) + \omega(\sigma(X), \chi^\uparrow) 
            &= -\huaL_{\psi^\uparrow} 0 + \huaL_{\chi^\uparrow} 0 + \Theta([\psi^\uparrow, \chi^\uparrow])\\
            &\qquad -\huaL_{\sigma(X)}0 + \huaL_{\chi^\uparrow} ev_{X} + \Theta([\sigma(X), \chi^\uparrow])\\
            &= {\chi^\uparrow} \cdot ev_{X} = q^*\langle \chi, X\rangle.
        \end{split}
    \end{equation}
    Here we used the fact that vertical lifts commute, and ${\chi^\uparrow} \cdot ev_{X} = q^*\langle \chi, X\rangle$ from \cite[\S 3.4, Eq. 21]{Mackenzie2005}, and the fact that the bracket $[\sigma(X), \chi^\uparrow]$ is a vertical lift from \cite[Thm. 3.4.5]{Mackenzie2005}, so that term also vanishes, by \eqref{eq:VB2CotangentTautFormCoreLinear}.               
    
    Write $\huaT(\huaG) = N(T\huaG)$ as
    \[\begin{tikzcd}[ampersand replacement=\&,row sep = 0ex
    ,/tikz/column 1/.append style={anchor=base east}
    ,/tikz/column 5/.append style={anchor=base west}]
	{A_2} \& {A_1} \& {T\huaG_0}
	\arrow["{-Td_1}", from=1-2, to=1-3]
	\arrow["{Td_2}", from=1-1, to=1-2]
    \end{tikzcd}\]
    where $A_2 := 1^*\ker Tp^2_2$ and $A_1 := 1^*\ker Td^1_0$. 
    Recall the chain isomorphism $\Psi$ from Proposition \ref{prop:vs-2dual-norm-cplx}, which allows us to write $N(T^{2*}\huaG)$ as
    \[\begin{tikzcd}[ampersand replacement=\&,cramped,column sep=scriptsize,row sep=tiny]
	{T^*\huaG_0 \oplus A_1^*} \& {A_1^*\oplus A_1^*} \& {A_2^*} \\
	{(\theta, \nu)} \& {(\nu, \nu +\theta \circ Td_1)} \\
	\& {(\eta, \xi)} \& {(\xi -\eta)\circ Td_2}
	\arrow[from=1-2, to=1-3]
	\arrow[from=1-1, to=1-2]
	\arrow[maps to, from=2-1, to=2-2]
	\arrow[maps to, from=3-2, to=3-3]
    \end{tikzcd}\]
    Then the following equations hold by a straightforward computation:
    \begin{equation}\label{eq:VB2CotangentNormIsoPsiEquations}
        \begin{array}{ll}
        \langle \Psi^{-1}(\theta, \nu), T1 x \rangle = \theta(x),
        &\langle \widecheck{s}_0 \Psi^{-1}(\eta, \xi), Ts_1 c \rangle = \eta(c),\\
        \langle \widecheck{s}_1 \Psi^{-1}(\eta, \xi), Ts_0 c \rangle = \xi(c),
        &\langle \widecheck{1} \Psi^{-1}(\beta), k \rangle = \beta(k),
        \end{array}
    \end{equation}
    for any $x \in T\huaG_0$, $c\in A_1$, $k \in A_2$ and $\beta \in A_2^*$, by definition of $\Psi$. 
    With this notation we also have that the maps induced by the pairing are, by \eqref{eq:lambda-ndual-pairing-right},
    \begin{equation*}
        \begin{array}{ll}
            \lambda_{\langle \cdot, \cdot \rangle}^r : N(T^{2*}\huaG) \to N(T\huaG)^*[-2] 
            &\lambda_{\langle \cdot, \cdot \rangle}^l : N(T\huaG) \to N(T^{2*}\huaG)^*[-2]\\
            (\lambda_{\langle \cdot, \cdot \rangle}^r)_2(\theta, \nu) =  \theta \in T^*\huaG_0, 
            &(\lambda_{\langle \cdot, \cdot \rangle}^l)_2(k) = k \in A_2,\\
            (\lambda_{\langle \cdot, \cdot \rangle}^r)_1(\xi, \eta) = \xi - \eta \in A_1^*, 
            &(\lambda_{\langle \cdot, \cdot \rangle}^l)_1(c) = (c, -c) \in A_1 \oplus A_1,\\
            (\lambda_{\langle \cdot, \cdot \rangle}^r)_0(\epsilon) = \epsilon \in A_2^*, 
            &(\lambda_{\langle \cdot, \cdot \rangle}^l)_0(x) = (x,0) \in T\huaG_0 \oplus A_1,
        \end{array}
    \end{equation*}
    where we used that $A_2\cong A_2^{**}$, $A_1 \oplus A_1 \cong (A_1^* \oplus A_1^*)^*$, and $T\huaG_0 \oplus A_1 \cong (T^*\huaG_0 \oplus A_1^*)^*$.
    With this, clearly, $(-1)^i(\lambda^r_i)^t = \lambda^l_{2-i}$. 
    We now consider diagram \eqref{diag:VB2CotangentOmegaTangentCommDiag} at a point $\epsilon$ of $T^{2*}\huaG_0|_p$. We get that commutativity of the left square at $\epsilon$ is equivalent to the equations
    \begin{equation}\label{eq:VB2CotangentOmegaTangentCommDiag-Left}
        \begin{split}
            &(\lambda_\omega^r)((\Psi^{-1}(\beta))^\uparrow|
            _{\epsilon})(H) = - \lambda_{\langle,\rangle}^r(\beta)(T_\epsilon q H),\\
            &(\lambda_\omega^r)((\Psi^{-1}(\eta, \xi))^\uparrow|
            _{\epsilon})(C) = - \lambda_{\langle,\rangle}^r(\eta, \xi)(T_\epsilon q C),\\
            &(\lambda_\omega^r)((\Psi^{-1}(\theta,\nu))^\uparrow|
            _{\epsilon})(U) = - \lambda_{\langle,\rangle}^r((\theta,\nu))(T_\epsilon q U),\\
        \end{split}
    \end{equation}
    for $\beta \in A_2^*|_p$, $H \in \huaT_2(T^{2*}\huaG)|_\epsilon$, $(\eta, \xi)  \in A_1^* \oplus A_1^*|_p$, $C \in \huaT_1(T^{2*}\huaG)|_\epsilon$, $(\theta, \nu) \in T^*\huaG_0 \oplus A_1^*|_p$, and $U \in TT^{2*}\huaG_0|_\epsilon$.
    On the other hand commutativity with the right square is equivalent to the equations
    \begin{equation}\label{eq:VB2CotangentOmegaTangentCommDiag-Right}
        \begin{split}
            &(\lambda_\omega^r)(U)((\Psi^{-1}(\theta,\nu))^\uparrow|
            _{\epsilon}) = \lambda_{\langle \cdot, \cdot \rangle}^l(T_\epsilon q U)(\theta, \nu),\\
            &(\lambda_\omega^r)(C)((\Psi^{-1}(\eta,\xi))^\uparrow|
            _{\epsilon}) = \lambda_{\langle \cdot, \cdot \rangle}^l(T_\epsilon q C)(\eta, \xi),\\
            &(\lambda_\omega^r)(H)((\Psi^{-1}(\beta))^\uparrow|
            _{\epsilon})= \lambda_{\langle \cdot, \cdot \rangle}^l(T_\epsilon q H)(\beta).\\
        \end{split}
    \end{equation}
    By antisymmetry of $\omega$ we have that for any $H \in \huaT_2(T^{2*}\huaG)|_\epsilon$, $C, C' \in \huaT_1(T^{2*}\huaG)|_\epsilon$, and $U \in TT^{2*}\huaG_0|_\epsilon$,
    \begin{equation*}
        \begin{split}
            &\lambda_\omega(H, U) 
            = \omega(H, T\widecheck{1}U) 
            = - \omega(T\widecheck{1}U, H) 
            = - \lambda_\omega(U, H)\\
            &\lambda_\omega(C, C') 
            = \omega(T\widecheck{s}_1C, T\widecheck{s}_0C') - \omega(T\widecheck{s}_0C, T\widecheck{s}_1C')
            = \lambda_\omega(C', C) 
        \end{split}
    \end{equation*}
    Therefore, because $(-1)^i(\lambda^r_i)^t = \lambda^l_{2-i}$, the first and the last equations in \eqref{eq:VB2CotangentOmegaTangentCommDiag-Left} are equivalent to the last and the first equations in \eqref{eq:VB2CotangentOmegaTangentCommDiag-Right}, respectively, while the middle equations on either side are equivalent.  
    These remaining equations are now straightforward to check, by using \eqref{eq:VB2CotangentOmegaOnVertLift}, \eqref{eq:VB2CotangentNormIsoPsiEquations} and the fact that vertical lifts commute with simplicial structure maps as explained in Example \ref{ex:TT2starG-as-DVB}, in fact:
    \begin{equation*}
        \begin{split}
            (\lambda_\omega^r)(H)((\Psi^{-1}(\beta))^\uparrow|_\epsilon)
            &= \omega(H, T\widecheck{1}(\Psi^{-1}(\beta))^\uparrow|_\epsilon) 
            = \langle \widecheck{1}\Psi^{-1}(\beta), T_\epsilon q H \rangle \\
            &= \beta (H) 
            = \lambda_{\langle,\rangle}^l(\beta)(T_\epsilon q H),\\
            (\lambda_\omega^r)(U)((\Psi^{-1}(\theta,\nu))^\uparrow|_\epsilon) 
            &= \omega(T\widecheck{1}U, (\Psi^{-1}(\theta,\nu))^\uparrow|_\epsilon)
            = \langle \Psi^{-1}(\theta,\nu),  \widecheck{1}T_\epsilon qU \rangle\\
            &= \theta(T_\epsilon qU)
            = \lambda_{\langle \cdot, \cdot \rangle}^l(T_\epsilon q U)(\theta, \nu),\\
            (\lambda_\omega^r)((\Psi^{-1}(\eta, \xi))^\uparrow|_\epsilon)(C) 
            &= \omega (T\widecheck{s}_1(\Psi^{-1}(\eta, \xi))^\uparrow|_\epsilon, T\widecheck{s}_0 C)\\
            &\qquad\qquad - \omega (T\widecheck{s}_0(\Psi^{-1}(\eta, \xi))^\uparrow|_\epsilon, T\widecheck{s}_1 C)\\
            &= \omega ((\widecheck{s}_1\Psi^{-1}(\eta, \xi))^\uparrow|_\epsilon, T\widecheck{s}_0 C)\\
            &\qquad\qquad - \omega ((\widecheck{s}_0\Psi^{-1}(\eta, \xi))^\uparrow|_\epsilon, T\widecheck{s}_1 C)\\
            &= - \langle \widecheck{s}_1\Psi^{-1}(\eta, \xi), Ts_0 C \rangle
            +  \langle \widecheck{s}_0\Psi^{-1}(\eta, \xi), Ts_1 C \rangle\\
            &= \eta(T_\epsilon q C) - \xi(T_\epsilon q C)
            = - \lambda_{\langle,\rangle}^r(\eta, \xi)(T_\epsilon q C).
        \end{split}
    \end{equation*}
\end{proof}

\section{The 2-cotangent of a 1-groupoid}

In this section we consider the case of a Lie 1-groupoid $\huaG$. In this case there are two versions of a cotangent 2-groupoid we can associate to $\huaG$. One is the 2-cotangent $T^{2*}\huaG$ and the other is the Artin-Mazur bar $\overline{T^{1*}\huaG}$ of the 1-cotangent $T^{1*}\huaG$ considered as a double groupoid. (Any simplicial vector bundle can be seen as a bisimplicial object where one direction is the nerve of its vector bundle structure.) We refer to \cite[\S 2.2]{MehtaTang2011} for the bar construction. We will show that these two ``versions'' of a 2-cotangent space for $\huaG$ are Morita equivalent. To begin with, we write down both of these objects explicitly. As usual, in order to avoid having to keep track of too many inversions, we differentiate explicitly between the two versions of the Lie algebroid of $\huaG$: $A_R$, defined as the right core of $T\huaG$ and $A_L$, defined as the left core thereof.

Another reason why we prefer to use the left Lie algebroid $A_L$ in defining $\overline{T^{1*}\huaG}$ is that this allows the Morita equivalence $\Xi$ in \eqref{eq:2Cota-1Cota-ME-def} to be expressed in terms of an Ehresmann connection splitting the core sequence for $A_R$ as in the more usual convention used in \cite{AriasAbadCrainic2013}.

Following Example \ref{ex:VB2DualVB1Gpd}, the 2-cotangent $T^{2*}\huaG$ is given by 
\[\begin{tikzcd}[ampersand replacement=\&,cramped,sep=scriptsize]
	\begin{array}{c} (T\huaG_1 \times_{T\huaG_0} T\huaG_1)^* \oplus d_0^*d_1^*A_R^* \\ \quad\oplus d_0^*d_0^*A_L^*\oplus d_2^*d_1^*A_R^* \end{array} \& {d_1^*A_R^*\oplus d_0^*A_L^*} \& {0\times \huaG_0} \\
	{\huaG_1 \times_{\huaG_0} \huaG_1} \& {\huaG_1} \& {\huaG_0}
	\arrow[shift left=2, from=1-1, to=1-2]
	\arrow[shift right=2, from=1-1, to=1-2]
	\arrow[from=1-1, to=1-2]
	\arrow["{q^{2*}_2}"', from=1-1, to=2-1]
	\arrow[shift left, from=1-2, to=1-3]
	\arrow[shift right, from=1-2, to=1-3]
	\arrow["{q^{2*}_1}"', from=1-2, to=2-2]
	\arrow["{q^{2*}_0}"', from=1-3, to=2-3]
	\arrow[shift left=2, from=2-1, to=2-2]
	\arrow[shift right=2, from=2-1, to=2-2]
	\arrow[from=2-1, to=2-2]
	\arrow[shift left, from=2-2, to=2-3]
	\arrow[shift right, from=2-2, to=2-3]
\end{tikzcd}\]
This has simplicial maps defined as in Example \eqref{ex:VB2DualVB1Gpd}, just by replacing $\widetilde{d}_i$ and $\widetilde{s}_i$ with $Td_i$ and $Ts_i$, respectively, for any $i$. In particular, to fix notation, a 2-simplex $(\phi, \alpha, \beta, \gamma) \in T^{2*}\huaG_2|_{(g,h)}$ has faces
\begin{equation}\label{eq:2Cota-1Gpd-face-maps}
\begin{split}
    \widecheck{d}_0(\phi, \alpha, \beta, \gamma) 
    &= (\alpha, \beta) \in T^{2*}\huaG_1|_g,\\
    \widecheck{d}_1(\phi, \alpha, \beta, \gamma) 
    &= (\gamma + \phi(0_g, \_ \cdot 0_h), \beta + \phi(0_g \cdot \_, 0_h)) \in T^{2*}\huaG_1|_{hg},\\
    \widecheck{d}_2(\phi, \alpha, \beta, \gamma)
    &= (\gamma, -\alpha(\_^{-1}) - \phi(\_^{-1} \cdot 0_g, 0_h \cdot \_)) \in T^{2*}\huaG_1|_{h}.
\end{split}
\end{equation}
By the content of the previous section this is 2-shifted symplectic with respect to the canonical form $\omega^{2*}$ as in Theorem \ref{thm:2Cota-symp-form}. In this particular case the tautological form $\Theta^{2*} \in \Omega^1(T^{2*}\huaG_2)$ defined in \eqref{eq:2Cota-taut-form-def} can be written as  
\begin{equation}\label{eq:2Cota-taut-form-1gpd}
\begin{split}
\Theta^{2*}|_{(\phi, \alpha, \beta, \gamma)_{(g,h)}}(V) &= \langle (\phi, \alpha, \beta, \gamma),Tq^{2*} V\rangle\\
&= \phi|_{(g,h)}(Td_0Tq^{2*}V, Td_2Tq^{2*}V),        
\end{split}
\end{equation}
when evaluated on any $V \in TT^{2*}\huaG_2|_{(\phi, \alpha, \beta, \gamma)_{(g,h)}}$. This is because, in this case, $Tq^{2*}V \in T\Lambda^2_1(\huaG)|_{(g,h)}$ can be written as the pair $(Td_0Tq^{2*}V, Td_2Tq^{2*}V)$. 

\subsection{A different 2-cotangent}\label{sec:bar-of-1-cota}

Because any $\VB$-groupoid is in particular a double groupoid, we can obtain a different ``cotangent 2-groupoid'' of a Lie 1-groupoid $\huaG$, by applying the Artin-Mazur bar construction from \cite[\S 2.2]{MehtaTang2011} to the 1-cotangent $T^{1*}\huaG$ in the form
\[\begin{tikzcd}[ampersand replacement=\&,cramped,sep=small]
	{T^*\huaG_1} \& {A_L^*} \\
	{\huaG_1} \& {\huaG_0}
	\arrow[shift right, from=1-1, to=1-2]
	\arrow[shift left, from=1-1, to=1-2]
	\arrow["{q^{1*}_1}"', from=1-1, to=2-1]
	\arrow["{q^{1*}_0}"', from=1-2, to=2-2]
	\arrow[shift left, from=2-1, to=2-2]
	\arrow[shift right, from=2-1, to=2-2]
\end{tikzcd}\]
with face maps
\begin{equation*}
\begin{array}{ll}
    (\widecheck{d}_0\eta)|_{d_0g}(b) = \eta|_g(0_g \cdot b),
    &(\widecheck{d}_1\eta)|_{d_1g}(b') = - \eta|_g((b')^{-1} \cdot 0_g).
\end{array}
\end{equation*}
With this, we get that $\overline{T^{1*}\huaG}$ is a Lie 2-groupoid\footnote{The hypotheses of this theorem demand that to get a Lie 2-groupoid (and not just a local one), the starting double groupoid must be \textit{full}. I.e. in this case, the map $(q, \widecheck{d}_0): T^*\huaG_1 \to d_0^*A^*$ needs to be surjective. This is true because this is the horn projection $\widecheck{p}^1_1$. In other words, every $\VB$ 1-groupoid is a full double groupoid.} by \cite[Thm 4.5]{MehtaTang2011}, and it is given at each level by 
\begin{equation*}
    \begin{split}
        &\overline{T^{1*}\huaG}_2 :=  (T^*\huaG_1 \times_{d_0q, \huaG_0, q} A_L^*) \times_{d_1 q\circ pr_1, \huaG_0, d_0} \huaG_1 = d_0^*T^*\huaG_1 \oplus d_0^*d_0^*A_L^*,\\
        &\overline{T^{1*}\huaG}_1 := d_0^*A_L^* \times_{q,\huaG_0,d_0} \huaG_1  = d_0^*A_L^*,\\
        &\overline{T^{1*}\huaG}_0 := \huaG_0 \cong 0 \times \huaG_0,
    \end{split}
\end{equation*}
where we identify $\overline{T^{1*}\huaG}_0$ with $0 \times \huaG_0$ to see $\overline{T^{1*}\huaG}$ as a $\VB$-groupoid over $\huaG$. In a diagram, this is the $\VB$ 2-groupoid
\[\begin{tikzcd}[ampersand replacement=\&,cramped]
	{d_0^*T^*\huaG_1 \oplus d_0^*d_0^*A_L^*} \& {d_0^*A_L^*} \& {0\times \huaG_0} \\
	{\Lambda^2_1(\huaG)} \& {\huaG_1} \& {\huaG_0}
	\arrow[shift left, from=1-2, to=1-3]
	\arrow[shift left, from=2-2, to=2-3]
	\arrow[shift right, from=2-2, to=2-3]
	\arrow[from=1-3, to=2-3]
	\arrow[from=1-2, to=2-2]
	\arrow[shift right, from=1-2, to=1-3]
	\arrow[from=1-1, to=2-1]
	\arrow[shift left=2, from=2-1, to=2-2]
	\arrow[shift right=2, from=2-1, to=2-2]
	\arrow[from=2-1, to=2-2]
	\arrow[shift left=2, from=1-1, to=1-2]
	\arrow[shift right=2, from=1-1, to=1-2]
	\arrow[from=1-1, to=1-2]
\end{tikzcd}\]
It is easy to see that the other two levels are also vector bundles over the corresponding levels of $\huaG$, and the face and degeneracy maps are bundle maps. In fact these are, for any\footnote{To be consistent with our convention of writing any horn $(g,h) \in \Lambda^2_1(\huaG)$ by increasing index of its faces, we write elements $(\eta, \mu, h) \in \overline{T^{1*}\huaG}_2|_{(g,h)}$ in this order, so that $\eta|_g \in d_0^*T^*\huaG_1|_{(g,h)}$, $\mu|_{d_0g} \in d_0^*d_0^*A_L^*|_{(g,h)}$ and $h \in \huaG_1$ such that $d_0h = d_1g$. For consistency, we then write elements $(\mu, g) \in \overline{T^{1*}\huaG}_1|_{g}$, with the basepoint on the right as well.} $(\eta, \mu, h) \in \overline{T^{1*}\huaG}_2|_{(g,h)}$, and any $\mu \in d_0^*A_L^*|_{g},$
\begin{equation*}
    \begin{split}
        &\overline{d}_0(\eta, \mu, h) = (\mu, g) \in d_0^*A_L^*|_{g},\\
        &\overline{d}_1(\eta, \mu, h) = (\widecheck{d}_0\eta + \mu, hg) =  (\eta(0_g \cdot \_) + \mu, hg) \in d_0^*A_L^*|_{hg},\\
        &\overline{d}_2(\eta, \mu, h) = (\widecheck{d}_1\eta, h) = (- \eta(\_^{-1} \cdot 0_g), h) \in d_0^*A_L^*|_{h},\\
        &\overline{s}_0(\mu) = (0_g, \mu, 1d_1 g) \in \overline{T^{1*}\huaG}_2|_{s_0g},\\
        &\overline{s}_1(\mu) = (\widehat{\mu}_{1d_0g}, 0_{1d_0g}, g) =  (\mu(\ggamma_L\_), 0_{1d_0g}, g)  \in \overline{T^{1*}\huaG}_2|_{s_1g},\\
    \end{split}
\end{equation*}
where $\widehat{\mu}_{1d_0g}$ is the element of $\Ann(1T\huaG_0) \subset T^*\huaG_1$ corresponding to $\mu_{d_0g} \in A_L^*$ by the usual isomorphism.

We now describe the multiplications of this $\VB$ 2-groupoid. Consider an arbitrary tetrahedron in $\widebar{T^{1*}\huaG}$. This is an object over the tetrahedron corresponding to the composable triple $(g,h,k)$ in level 3 of the nerve of $\huaG$. We write it as a boundary and denote it by
\begin{equation}\label{eq:2Cota-1CotaBar-generic-tetra}
    \left((\eta_0, \mu_0)|_{(g,h)}, (\eta_1, \mu_1)|_{(g,kh)}, (\eta_2, \mu_2)|_{(hg,k)}, (\eta_3, \mu_3)|_{(h,k)}\right) \in \widebar{T^{1*}\huaG}_3 \cong \partial^3(\widebar{T^{1*}\huaG}), 
\end{equation}
where the components obey the boundary conditions
\begin{equation*}
\begin{array}{lll}
    \mu_0 = \mu_1, 
    & \widecheck{d}_0\eta_1 + \mu_1 = \widecheck{d}_0\eta_2 + \mu_2,
    & \widecheck{d}_1\eta_2 = \widecheck{d}_1\eta_3,\\
    & & \\
    \widecheck{d}_0\eta_0 + \mu_0 = \mu_2, 
    & \widecheck{d}_1\eta_1 = \widecheck{d}_0\eta_3 + \mu_3,
    & \widecheck{d}_1\eta_0 = \mu_3.
\end{array}
\end{equation*}
As usual, the faces of any given $(3,j)$-horn obey the subset of these boundary conditions where $\eta_j$ and $\mu_j$ do not appear. These are determined by the multiplications
\begin{equation}\label{eq:2Cota-1Cota-bar-multiplications}
\begin{split}
    \square(\eta_1, \mu_1)(\eta_2, \mu_2)(\eta_3, \mu_3) &= (\eta_1 - (\eta_3)^{-1}\cdot \eta_2, \mu_1)|_{(g,h)},\\
    (\eta_0, \mu_0)\square(\eta_2, \mu_2)(\eta_3, \mu_3) &= (\eta_0 + (\eta_3)^{-1}\cdot \eta_2, \mu_0)|_{(g,h)},\\
    (\eta_0, \mu_0)(\eta_1, \mu_1)\square(\eta_3, \mu_3) &= (\eta_3 \cdot (\eta_1 - \eta_0), \widecheck{d}_0\eta_0 + \mu_0)|_{(hg,k)},\\
    (\eta_0, \mu_0)(\eta_1, \mu_1)(\eta_2, \mu_2)\square&= (\eta_2 \cdot (\eta_1 - \eta_0)^{-1}, \widecheck{d}_1\eta_0)|_{(h,k)}.
\end{split}
\end{equation}
In particular, the four $\eta$ components of the boundary of a tetrahedron such as \eqref{eq:2Cota-1CotaBar-generic-tetra} satisfy the equation
\begin{equation*}
    \eta_1|_g - \eta_0|_g = (\eta_3|_h)^{-1} \cdot \eta_2|_{hg}.
\end{equation*}

Because $T^{1*}\huaG$ is a symplectic groupoid, by \cite[Thm 4.19]{Alvarez2023}, $\widebar{T^{1*}\huaG}$ is 2-shifted symplectic with respect to the form $\widebar{\omega} := P^*\omega_{can}$, with $P$ the composition of the canonical projections
\begin{equation}\label{diag:2Cota-Bar-Proj-Def}
\begin{tikzcd}[ampersand replacement=\&,cramped,sep=small]
	{P:} \& {d_0^*T^*\huaG_1 \oplus d_0^*d_0^*A_L^*} \& {T^*\huaG_1 \oplus d_0^*A_L^*} \& {T^*\huaG_1} \\
	\& {\Lambda^2_1(\huaG)} \& {\huaG_1} \& {\huaG_1}
	\arrow[from=1-2, to=1-3]
	\arrow["{\widebar{q}_2}"', from=1-2, to=2-2]
	\arrow[from=1-3, to=1-4]
	\arrow["{\widebar{q}_1}"', from=1-3, to=2-3]
	\arrow["{q^{1*}_1}"', from=1-4, to=2-4]
	\arrow[""{name=0, anchor=center, inner sep=0}, "{d_0}", from=2-2, to=2-3]
	\arrow[Rightarrow, no head, from=2-3, to=2-4]
	\arrow["\lrcorner"{anchor=center, pos=0.125}, draw=none, from=1-2, to=0]
\end{tikzcd}
\end{equation}
and $\omega_{can}$ the canonical symplectic form on $T^*\huaG$. With respect to the tautological 1-form
\begin{equation*}
    \Theta^{1*} = \theta_{can} \in \Omega^1(T^{1*}\huaG_1) = \Omega^1(T^*\huaG_1),
\end{equation*}
we define $\widebar{\Theta}:= P^*\theta_{can} \in \Omega^1(\widebar{T^{1*}\huaG}_2)$ and we have 
\begin{equation}\label{eq:2Cota-Bar-SymplForm-Def}
    \widebar{\omega} = -dP^*\theta_{can} = -d \widebar{\Theta}.
\end{equation} 
Because $P$ is a bundle map over $d_0$, so it commutes with the projections as in diagram \eqref{diag:2Cota-Bar-Proj-Def}, the identity
\begin{equation}\label{eq:2Cota-Bar-taut-form}
    \begin{split}
        \widebar{\Theta}|_{(\eta, \mu, h)_{(g,h)}}(V) = \theta_{can}|_{\eta_g}(TP(V)) &= \langle \eta_g, (Tq^{1*}_1 \circ TP)(V) \rangle\\
        &= \langle \eta_g, (Td_0 \circ T\widebar{q}_2)(V) \rangle\\
    \end{split}
\end{equation}
holds for any $V \in T\widebar{T^{1*}\huaG}_2|_{(\eta, \mu, h)_{(g,h)}}$.

The normalized complex of $\widebar{T^{1*}\huaG}$ is readily computed: it is $N(T^{1*}\huaG)[-1]$. 

\begin{lemma}\label{lem:1CotaBar-norm-cplx}
    The normalized complex of $\widebar{T^{1*}\huaG}$ is concentrated in degrees 2 and 1, and it is isomorphic to the complex $N(T^{1*}\huaG)[-1] = N(T\huaG)^*[-2]$:
    \begin{equation}\label{eq:1CotaBar-norm-cplx}
    \begin{array}{ccccc}
        T^*\huaG_0 &\overset{\widebar{\partial}}{\longrightarrow} &A^*_L &\longrightarrow &0 \\
        \theta &\longmapsto &- \theta \circ Td_0 & &
    \end{array}
    \end{equation}
\end{lemma}

\begin{proof}
    A quick computation shows
    \begin{equation*}
            N(\overline{T^{1*}\huaG})_2|_{p} = \ker \overline{p}^2_2 = \{(\phi_{1p}, 0)\in \overline{T^{1*}\huaG}_2 \mid \widecheck{d}_0\phi_{1p} = 0\} \cong \ker\widecheck{d}_0|_{1p} \cong T^*\huaG_0,
    \end{equation*}
    with the last isomorphism given by $\phi \mapsto \theta:= \phi\circ T1$, $\phi:=\theta Td_1 \mapsfrom \theta$, since $\widecheck{d}_0 \phi = \phi\circ (id - T1Td_1) =0$. The differential is then $\overline{d}_2(\phi, 0) = \widecheck{d}_1\phi$, which when paired with $b \in A_L$ gives $\widecheck{d}_1\phi(b) = \phi(b - T1Td_0b) = \phi(T1Td_1b - T1Td_0b) = -\theta(Td_0b)$. So we have that $N(\overline{T^{1*}\huaG})\cong N(T\huaG)^*[-2]$, as pictured in \eqref{eq:1CotaBar-norm-cplx}.
\end{proof}

\subsection{Equivalence of the two models}

The Morita equivalence $\Xi: T^{2*}\huaG \to \overline{T^{1*}\huaG}$ that we construct in this section depends on the choice of an Ehresmann connection on $\huaG$. We recall some necessary properties of such objects from \cite{AriasAbadCrainic2013}. These are the same as normal $(1,1)$-cleavages for $T\huaG$ from Definition \ref{def:cleavage}.

\begin{definition}\label{def:Ehresmann-connection}
An \textbf{Ehresmann connection} is a splitting $(\kappa, \sigma)$ of the sequence 
\[\begin{tikzcd}[ampersand replacement=\&,cramped,column sep=scriptsize,row sep=tiny]
	0 \& {d_1^*A_R} \& {T\huaG_1} \& {d_0^*T\huaG_0} \& 0
	\arrow["{Td_0}", from=1-3, to=1-4]
	\arrow["r", from=1-2, to=1-3]
	\arrow[from=1-1, to=1-2]
	\arrow[from=1-4, to=1-5]
	\arrow["\kappa", curve={height=-12pt}, from=1-3, to=1-2]
	\arrow["\sigma", curve={height=-12pt}, from=1-4, to=1-3]
\end{tikzcd}\]
which coincides with the canonical splitting $(id - T1Td_0, T1)$ of the same sequence over the units of $\huaG$, given in \eqref{eq:VBnDualCoreSES-1}. Here $r$ is the right translation map, which is given at any $g\in \huaG_1$ and for any $a \in A_R|_{d_1g}$ by $r_g(a) = a \cdot 0_g$.
Recall that a choice of either $\sigma$ or $\kappa$ determines the other one, as we have that any $v_g \in T\huaG_1|_g$ splits as 
\begin{equation}\label{eq:2Cota-Ehresmann-conn-splitting}
    v_g = \kappa_g(v_g) \cdot 0_g + \sigma_g(Td_0 v_g).
\end{equation}    
\end{definition}

Such a choice of connection allows us to define the quasi-actions of the adjoint representation up to homotopy of $\huaG$ on $\rho = Td_1: A \to TM$ as 
\begin{equation}\label{eq:Ehresmann-conn-quasi-actions}
\begin{split}
    &\lambda_g(x_{d_0g}) = Td_1 \sigma_g(x_{d_0g})\\
    &\lambda_g(a_{d_0g}) = - \kappa_g(0_g \cdot (a_{d_0g}))
\end{split}
\end{equation}
for all $x_{d_0g} \in T_{d_0g}\huaG_0$ and all $a_{d_0g} \in A_R|_{d_0g}$. 

Recall furthermore that the failure of $\sigma$ to be multiplicative is encoded in the basic curvature $R^{bas}\in \Gamma(\Lambda^2_1(\huaG), \Hom(d_0^*T\huaG_0, d_1^*A_R))$ defined equivalently by
\begin{equation}
   \begin{split}
    R^{bas}_{(g,h)}(x_{d_0g}) 
    &= \left(\sigma_{hg}(x_{d_0g}) - \sigma_h(Td_1\sigma_g(x_{d_0g})) \cdot \sigma_g(x_{d_0g})\right) \cdot 0_{hg}^{-1}\\
    &= - \kappa_{hg}(\sigma_h(Td_1\sigma_g(x_{d_0g})) \cdot \sigma_g(x_{d_0g})).
   \end{split}
\end{equation}
This also encodes the failure of $\lambda$ to be groupoid representations, as we recall in the following lemma. We refer to \cite[Prop. 2.15]{AriasAbadCrainic2013} for a proof.

\begin{lemma}
For any $v \in T_g\huaG_1$ and $(g,h) \in \Lambda^2_1(\huaG)$, 
\begin{equation}\label{eq:Cota2-Ehresmann-conn-prop1}
    R^{bas}_{(g,h)}(Td_0v) = - \kappa_{hg}((\sigma_hTd_1v) \cdot v) - \kappa_h(0_h \cdot (\kappa_g(v))^{-1}).
\end{equation}
In particular, for $a \in A_R|_{d_0g}$, 
\begin{equation*}
    R^{bas}_{(g,h)}(Td_1a) = - \kappa_{hg}(0_{hg} \cdot a^{-1}) - \kappa_h(0_h \cdot (\kappa_g(0_g \cdot a^{-1}))^{-1})
    = \lambda_{hg} a - \lambda_h \lambda_g a. 
\end{equation*}
\end{lemma}

We now define $\Xi: T^{2*}\huaG \to \overline{T^{1*}\huaG}$ as 
\begin{equation}\label{eq:2Cota-1Cota-ME-def}
    \begin{split}
        &\Xi_0 = 0 \\
        &\Xi_1(\alpha, \beta)_g (b) = (\beta(b) - \alpha(\kappa_g(0_g \cdot b)), g)\\
        &\Xi_2(\phi, \alpha, \beta, \gamma)_{(g,h)}(v, b) 
        = \big(\alpha(\kappa_g (v)) 
        - \gamma\left(\kappa_{hg}((\sigma_hTd_1v) \cdot v)\right)\\ 
        &\qquad\qquad + \phi_{(g,h)}\left(v, - \kappa_{hg}((\sigma_hTd_1v) \cdot v) \cdot 0_h + \sigma_h Td_1 v \right), \Xi_1(\alpha, \beta)(b), h \big)
    \end{split}
\end{equation}
for any $(\alpha, \beta)_g \in T^{2*}\huaG_1|_g$, $b \in A_L|_{d_0g}$, $(\phi, \alpha, \beta, \gamma) \in T^{2*}\huaG_2|_{(g,h)}$, and $v \in T_g\huaG_1$.

\begin{remark}
    The map $\Xi$ can be rewritten in an alternative form including the basic curvature by using \eqref{eq:Cota2-Ehresmann-conn-prop1} to rewrite the arguments of $\gamma$ and $\phi$ in the first component of $\Xi_2$ in \eqref{eq:2Cota-1Cota-ME-def}. By the splittings \eqref{eq:2Cota-Ehresmann-conn-splitting} of $0_h \cdot (\kappa_g(v))^{-1}$ and of $v$ we have
    \begin{equation*}
    \begin{split}
        0_h \cdot (\kappa_g(v))^{-1} &= \kappa_h(0_h\cdot(\kappa_g(v))^{-1}) \cdot 0_h + \sigma_h Td_0 (0_h \cdot (\kappa_g(v))^{-1}) \\
        &= \kappa_h(0_h\cdot(\kappa_g(v))^{-1}) \cdot 0_h + \sigma_h Td_1 \kappa_g(v)  \\
        &= \kappa_h(0_h\cdot(\kappa_g(v))^{-1}) \cdot 0_h + \sigma_h Td_1 v - \sigma_h Td_1 \sigma_g Td_0 v.
    \end{split}
    \end{equation*}
    With this and \eqref{eq:Cota2-Ehresmann-conn-prop1}, the second component of the argument of $\phi$ becomes
    \begin{equation*}
    \begin{split}
        - \kappa_{hg}((\sigma Td_1 v ) \cdot v) \cdot 0_h + \sigma_h Td_1 v 
        &=  \kappa_h(0_h\cdot(\kappa_g(v))^{-1}) \cdot 0_h \\
        &\qquad + R^{bas}_{(g,h)}(Td_0v) \cdot 0_h + \sigma_h Td_1 v \\
        &= 0_h \cdot (\kappa_g(v))^{-1} + \sigma_h Td_1 \sigma_g Td_0 v\\
        &\qquad + R^{bas}_{(g,h)}(Td_0v) \cdot 0_h 
    \end{split}
    \end{equation*}
    Hence $\Xi$ can be rewritten as
    \begin{equation}\label{eq:2Cota-1Cota-ME-def-alt}
    \begin{split}
        &\Xi_0 = 0 \\
        &\Xi_1(\alpha, \beta)_g (b) = (\beta(b) - \alpha(\kappa_g(0_g \cdot b)), g)\\
        &\Xi_2(\phi, \alpha, \beta, \gamma)_{(g,h)}(v, b) \\
        &\quad = \left(\alpha(\kappa_g (v)) 
        + \gamma\left(\kappa_h(0_h \cdot (\kappa_g(v))^{-1}) + R^{bas}_{(g,h)}(Td_0 v)\right)\right.\\ 
        &\qquad \left. + \Phi_{(g,h)}\left(v, 0_h \cdot (\kappa_g(v))^{-1} + \sigma_h Td_1 \sigma_g Td_0 v + R^{bas}_{(g,h)}(Td_0v) \cdot 0_h \right), \Xi_1(\alpha, \beta)(b) \right)
    \end{split}
    \end{equation}
    This is helpful in some computations such as the ones in the proof of the following lemma.

    Notice that we could also rewrite this in terms of the quasi-actions \eqref{eq:Ehresmann-conn-quasi-actions}, to obtain yet another form of $\Xi$ in terms of the adjoint representation. This notably shows that $\Xi_1(\alpha, \beta)_g (b) = \beta(b) + \alpha(\lambda_g(b^{-1}))$, but other than perhaps helping with intuition, this form of $\Xi$ is not particularly useful in the following discussion.
\end{remark}

\begin{lemma}
    The map $\Xi: T^{2*}\huaG \to \overline{T^{1*}\huaG}$ defined by \eqref{eq:2Cota-1Cota-ME-def} is a well-defined simplicial vector bundle map over the identity of $\huaG$. 
\end{lemma}

\begin{proof}
    By definition, each level of $\Xi$ is clearly a bundle map over the identity of the respective level of $\huaG$. Since both $T^{2*}\huaG$ and $\overline{T^{1*}\huaG}$ are $\VB$ 2-groupoids with level 0 of vanishing rank, it is enough to show that $\Xi$ commutes with face and degeneracy maps between levels 1 and 2, and that it is multiplicative with respect to $\widebar{m}_2$ from \eqref{eq:2Cota-1Cota-bar-multiplications}.
    
    By definition, $\widebar{d}_0\Xi_2 = \Xi_1\widecheck{d}_0$. 
    The identity $\widebar{d}_1\Xi_2 = \Xi_1\widecheck{d}_1$ follows by a straightforward computation:
    For any arbitrary $(\phi, \alpha, \beta, \gamma) \in T^{2*}\huaG_2|_{(g,h)}$ and $b \in A_L|_{d_0g}$,
    \begin{equation*}
    \begin{split}
        (\Xi_1\widecheck{d}_1(\phi, \alpha, \beta, \gamma))(b) &= \beta(b) + \phi(0_g \cdot b, 0_h) - \gamma(\kappa_{hg}(0_{hg}\cdot b))\\
        &\qquad\qquad - \phi(0_g, \kappa_{hg}(0_{hg}\cdot b)\cdot 0_h)\\ 
        &= \beta(b) + \phi(0_g \cdot b, - \kappa_{hg}(0_{hg}\cdot b)\cdot 0_h) - \gamma(\kappa_{hg}(0_{hg}\cdot b)).
    \end{split}
    \end{equation*}
    Because $\sigma_h Td_1(0_g \cdot b) = 0_h$,
    \begin{equation*}
    \begin{split}
        (\widebar{d}_1\Xi_2(\phi, \alpha, \beta, \gamma))(b) &= \cancel{\alpha(\kappa_g(0_g \cdot b))} - \gamma(\kappa_{hg}(0_h \cdot 0_g \cdot b))\\
        &\qquad + \phi(0_g \cdot b, - \kappa_{hg}(0_h \cdot 0_g \cdot b) \cdot 0_h + 0_h)\\
        &\qquad + \beta(b) - \cancel{\alpha(\kappa_g(0_g \cdot b))} = (\Xi_1\widecheck{d}_1(\phi, \alpha, \beta, \gamma))(b).
    \end{split}
    \end{equation*}
    The identity $\widebar{d}_2\Xi_2 = \Xi_1\widecheck{d}_2$ is also straightforward, by using the alternative form of $\Xi_2$ in \eqref{eq:2Cota-1Cota-ME-def}. This is because, for any $b \in A_L|_{d_1g = d_0h}$,  $Td_0(b^{-1}\cdot 0_g)=0$, so the curvature and $\sigma$ terms vanish in the following:
    \begin{equation*}
    \begin{split}
        (\widebar{d}_2 \Xi_2 (\phi,\alpha,\beta, \gamma))(b) &= - \alpha (\kappa_g(b^{-1} \cdot 0_g)) - \gamma(\kappa_h(0_h\cdot(\kappa_g(b^{-1}\cdot 0_g))^{-1}))\\
        &\qquad - \phi(b^{-1} \cdot 0_g, 0_h \cdot (\kappa_g(b^{-1}\cdot 0_g))^{-1})\\
        &= - \alpha(b^{-1}) - \gamma(\kappa_h(0_h \cdot b)) - \phi(b^{-1} \cdot 0_g, 0_h \cdot b)\\
        &= (\Xi_1\widecheck{d}_2 (\phi,\alpha,\beta, \gamma))(b).
    \end{split}
    \end{equation*}

    We now show the map $\Xi$ preserves the multiplication $m_2$. That is
    \begin{equation}\label{eq:2Cota-1Cota-Xi-multiplicative}
    \begin{split}
        (\Xi_2&(\phi_0, \alpha_0, \beta_0, \gamma_0)|_{(g,h)})
        (\Xi_2(\phi_1, \alpha_1, \beta_1, \gamma_1)|_{(g,kh)})
        \square
        (\Xi_2(\phi_3, \alpha_3, \beta_3, \gamma_3)|_{(h,k)})
        =\\
        &\Xi_2((\phi_0, \alpha_0, \beta_0, \gamma_0)|_{(g,h)}(\phi_1, \alpha_1, \beta_1, \gamma_1)|_{(g,kh)}\square (\phi_3, \alpha_3, \beta_3, \gamma_3)|_{(h,k)}),
    \end{split}
    \end{equation}
    for any $(\phi_i, \alpha_i, \beta_i, \gamma_i)_{i \neq 2} \in \Lambda^3_2(\huaV)|_{(g,h,k)}$. For any such element, the following horn conditions from \eqref{eq:VB2Dual1Gpd-3-horn-conditions} hold:
    \begin{equation*}
    \begin{array}{lll}
        \begin{cases}
            \alpha_0 = \alpha_1, & \\
            \beta_0 = \beta_1,
        \end{cases}
        &
        \begin{cases}
            \alpha_3 = \gamma_0, & \\
            \beta_3 = \ggamma_R^* \alpha_0 - \llambda_1^*\phi_0,
        \end{cases}
        &
        \begin{cases}
            \gamma_3 + \llambda_0^*\phi_3 = \gamma_1, & \\
            \beta_3 + \llambda_2^*\phi_3 = \ggamma_R^* \alpha_1 - \llambda_1^*\phi_1.
        \end{cases}
    \end{array}
    \end{equation*}
    Recall that any element in the fiber of $T\huaG_1$ at $hg$ can be written as a multiplication $w\cdot v$ with $w\in T\huaG_1|_h$ and $v \in T\huaG_1|_g$, and consider an arbitrary element $(w \cdot v, b) \in (d_0^*T\huaG_1 \oplus d_0^*d_0^*A_L)|_{(hg,k)}$.
    By definition of the multiplications in \eqref{eq:2Cota-1Cota-bar-multiplications}, the left-hand side of \eqref{eq:2Cota-1Cota-Xi-multiplicative} evaluated on this element is (dropping the basepoints for convenience)
    \begin{equation*}
    \begin{split}
        \big((\Xi_2&(\phi_0, \alpha_0, \beta_0, \gamma_0))
        (\Xi_2(\phi_1, \alpha_1, \beta_1, \gamma_1))
        \square
        (\Xi_2(\phi_3, \alpha_3, \beta_3, \gamma_3))
        \big)(w \cdot v,b)=\\
        \bigg(
        &\alpha_3(\kappa_h(w)) 
        - \gamma_3(\kappa_{kh}((\sigma_{kh}Td_1w)\cdot w))\\
        &\qquad 
        + \phi_3|_{(h,k)}(w, -\kappa_{kh}((\sigma_{kh}Td_1w)\cdot w) \cdot 0_k + \sigma_kTd_1w) \\
        &\qquad 
        + \alpha_1(\kappa_g(v)) 
        - \gamma_1(\kappa_{khg}((\sigma_{kh}Td_1v)\cdot v))\\
        &\qquad 
        + \phi_1|_{(g,kh)}(v, -\kappa_{khg}((\sigma_{kh}Td_1v)\cdot v) \cdot 0_{kh} + \sigma_{kh}Td_1v)\\
        &\qquad
        - \alpha_0(\kappa_g(v)) 
        + \gamma_0(\kappa_{hg}((\sigma_hTd_1v)\cdot v))\\
        &\qquad
        - \phi_0|_{(g,h)}(v, - \kappa_{hg}((\sigma_hTd_1v)\cdot v) \cdot 0_h + \sigma_h Td_1 v),\\
        &\quad\begin{multlined}
        \alpha_0(\kappa_g(0_g \cdot b)) 
        - \gamma_0(\kappa_{hg}(0_{hg} \cdot b)) 
        + \phi_0|_{(g,h)}(v, -\kappa_{hg}(0_{hg} \cdot b) \cdot 0_h)\\
        \qquad + \beta_0(b) - \alpha_0(\kappa_g(0_g \cdot b)),\quad h \bigg),
        \end{multlined}
    \end{split}
    \end{equation*}
    where in the second component we used that $\sigma_h Td_1(0_g \cdot b) = 0_h$ by linearity of the splitting.
    By the definition of the multiplications of the 2-cotangent in \eqref{eq:VB2Dual1Gpd-def-multiplications}, the right-hand side of \eqref{eq:2Cota-1Cota-Xi-multiplicative} when evaluated on the same element $(w\cdot v, b)$ is 
    \begin{equation*}
    \begin{split}
        (\Xi_2(\phi_2, &\gamma_0 + \llambda_0^*\phi_0, \beta_0 + \llambda_2^*\phi_0, \gamma_3))(w \cdot v,b) = \\
        &\bigg( 
        \gamma_0(\kappa_{hg}(w\cdot v))
        + \phi_0(\llambda_0(\kappa_{hg}(w\cdot v)))
        - \gamma_3(\kappa_{khg}((\sigma_k Td_1 w) \cdot w \cdot v))\\
        &\qquad 
        + \phi_2(w\cdot v, - \kappa_{khg}((\sigma_k Td_1 w) \cdot w \cdot v) \cdot 0_k + \sigma_k Td_1 w), \\
        &\quad
        \beta_0(b) + \phi_0(\llambda_2 b) - \gamma_0(\kappa_{hg}(0_{hg} \cdot b)) - \phi_0(\llambda_0(\kappa_{hg}(0_{hg} \cdot b))), \quad h \bigg),
    \end{split}
    \end{equation*}
    where $\phi_2$ is determined by equation \eqref{eq:VB2Dual1Gpd-def-multiplications-intcomp}. 
    Immediately, the second components of each side are the same, simply by writing out $\llambda_0$ and $\llambda_2$ as in \eqref{eq:VB2Dual1Gpd-lambda-horn-lifts}. 
    Equality of the first components of each side of \eqref{eq:2Cota-1Cota-Xi-multiplicative} is equivalent to
    \begin{equation*}
    \begin{split}
        \alpha_3(&\kappa_h(w)) 
        + \alpha_3(\kappa_{hg}((\sigma_hTd_1v)\cdot v))\\
        &- \gamma_1(\kappa_{khg}((\sigma_{kh}Td_1v)\cdot v))
        - \gamma_3(\kappa_{kh}((\sigma_{kh}Td_1w)\cdot w))\\
        &- \phi_0|_{(g,h)}(v, - \kappa_{hg}((\sigma_hTd_1v)\cdot v) \cdot 0_h + \sigma_h Td_1 v)\\
        &+ \phi_1|_{(g,kh)}(v, -\kappa_{khg}((\sigma_{kh}Td_1v)\cdot v) \cdot 0_{kh} + \sigma_{kh}Td_1v)\\
        &+ \phi_3|_{(h,k)}(w, -\kappa_{kh}((\sigma_{kh}Td_1w)\cdot w) \cdot 0_k + \sigma_kTd_1w)\\
        = &\alpha_3(\kappa_{hg}(w\cdot v))
        - \gamma_3(\kappa_{khg}((\sigma_k Td_1 w) \cdot w \cdot v))\\
        &\qquad + \phi_0(0_g, \kappa_{hg}(w\cdot v) \cdot 0_h) 
        - \phi_0 (v,w)\\
        &\qquad 
        + \phi_1 (v, ( - \kappa_{khg}((\sigma_k Td_1 w) \cdot w \cdot v) \cdot 0_k + \sigma_k Td_1 w)\cdot w)\\
        &\qquad
        + \phi_3(w, - \kappa_{khg}((\sigma_k Td_1 w) \cdot w \cdot v) \cdot 0_k + \sigma_k Td_1 w),
    \end{split}
    \end{equation*}
    by using the horn conditions $\alpha_0 = \alpha_1$ and $\gamma_0 = \alpha_3$, and writing out the $\phi_2$ term according to equation \eqref{eq:VB2Dual1Gpd-def-multiplications-intcomp}.
    This follows from the following four equations.  
    \begin{enumerate}[label={(\roman*)}]
        \item $\alpha_3(\kappa_h(w)) + \alpha_3(\kappa_{hg}((\sigma_hTd_1v)\cdot v)) 
        = \alpha_3(\kappa_{hg}(w\cdot v))$,
        \item $\phi_0(v, - \kappa_{hg}((\sigma_hTd_1v)\cdot v) \cdot 0_h + \sigma_h Td_1 v) 
        = \phi_0(v, w - \kappa_{hg}(w\cdot v) \cdot 0_h)$,
        \item $\begin{aligned}[t]
        \phi_1(v&, -\kappa_{khg}((\sigma_{kh}Td_1v)\cdot v) \cdot 0_{kh} + \sigma_{kh}Td_1v)\\
        =& \phi_1 (v, ( - \kappa_{khg}((\sigma_k Td_1 w) \cdot w \cdot v) \cdot 0_k + \sigma_k Td_1 w)\cdot w),
        \end{aligned}$
        \item $\begin{aligned}[t]
        - \gamma_1&(\kappa_{khg}((\sigma_{kh}Td_1v)\cdot v))
        - \gamma_3(\kappa_{kh}((\sigma_{kh}Td_1w)\cdot w)) \\
        &+ \phi_3(w, -\kappa_{kh}((\sigma_{kh}Td_1w)\cdot w) \cdot 0_k + \sigma_kTd_1w)\\
        =& - \gamma_3(\kappa_{khg}((\sigma_k Td_1 w) \cdot w \cdot v))\\
        &+ \phi_3(w, - \kappa_{khg}((\sigma_k Td_1 w) \cdot w \cdot v) \cdot 0_k + \sigma_k Td_1 w),
        \end{aligned}$
    \end{enumerate}
    To show (i) we show that 
    \begin{equation}\label{eq:2Cota-1Cota-kappa-wv}
        \kappa_{hg}(w\cdot v) = \kappa_h(w) + \kappa_{hg}((\sigma_h Td_1 v) \cdot v).
    \end{equation}
    This follows by splitting $w$ over $h$ according to \eqref{eq:2Cota-Ehresmann-conn-splitting}. With that, 
    \begin{equation*}
        w \cdot v = (\kappa_h(w) \cdot 0_h + \sigma_hTd_0w) \cdot (0_g + v) = \kappa_h(w) \cdot 0_{hg} + (\sigma_h Td_1 v) \cdot v.
    \end{equation*}
    Then by definition of $\kappa_{hg}$ as a left inverse of $\_\cdot 0_{hg}$, $\kappa_{hg}(\kappa_h(w) \cdot 0_{hg}) = \kappa_h(w)$, so \eqref{eq:2Cota-1Cota-kappa-wv} holds, which implies (i). 
    Equation (ii) also follows from \eqref{eq:2Cota-1Cota-kappa-wv} and splitting $w$ over $h$. In fact
    \begin{equation*}
    \begin{split}
        w - \kappa_{hg}(w\cdot v) \cdot 0_h 
        &= \kappa_h(w) \cdot 0_h + \sigma_hTd_0w - \kappa_h(w) \cdot 0_h - \kappa_{hg}((\sigma_h Td_1 v) \cdot v)\cdot 0_h\\
        &= - \kappa_{hg}((\sigma_h Td_1 v) \cdot v)\cdot 0_h + \sigma_hTd_0w.
    \end{split}
    \end{equation*}
    For equations (iii) and (iv) we use the fact that by splitting $(\sigma_k Td_1w) \cdot w$ over $kh$ we obtain 
    \begin{equation}\label{eq:2Cota-1Cota-kappa-uwv}
    \begin{split}
        \kappa_{khg}((\sigma_k Td_1 w) \cdot w \cdot v) &= \kappa_{khg}((\kappa_{kh}((\sigma_k Td_1 w) \cdot w) \cdot 0_{kh} + \sigma_{kh} Td_0 w) \cdot (0_g + v))\\
        &= \kappa_{kh}((\sigma_k Td_1 w) \cdot w) + \kappa_{khg}((\sigma_{kh} Td_1 v) \cdot v).
    \end{split}
    \end{equation}
    With this and the same splitting of $(\sigma_k Td_1w) \cdot w$ over $kh$, the second component of the argument of $\phi_1$ on the right-hand side of (iii) becomes
    \begin{equation*}
    \begin{split}
        ( - \kappa_{khg}&((\sigma_k Td_1 w) \cdot w \cdot v) \cdot 0_k + \sigma_k Td_1 w)\cdot w\\
        &= - \kappa_{khg}((\sigma_k Td_1 w) \cdot w \cdot v) \cdot 0_{kh} + (\sigma_k Td_1 w) \cdot w\\
        &= - \kappa_{kh}((\sigma_k Td_1 w) \cdot w) \cdot 0_{kh}  - \kappa_{khg}((\sigma_{kh} Td_1 v) \cdot v) \cdot 0_{kh}\\
        &\qquad + \kappa_{kh}((\sigma_k Td_1 w) \cdot w) \cdot 0_{kh} + \sigma_{kh} Td_0 w\\
        &= - \kappa_{khg}((\sigma_{kh} Td_1 v) \cdot v) \cdot 0_{kh}  + \sigma_{kh} Td_1 v,
    \end{split}
    \end{equation*}
    which is the second component of the argument of $\phi_1$ on the left-hand side of (iii), which then holds. 
    Finally, (iv) holds by a similar argument after rewriting $\gamma_1$ by using the horn condition $\gamma_1 = \gamma_3 + \llambda_0^*\phi_3$: this gives 
    \begin{equation*}
    \begin{split}
        -\gamma_1(\kappa_{khg}((\sigma_{kh} Td_1 v) \cdot v)) &= \gamma_3(\kappa_{khg}((\sigma_{kh} Td_1 v) \cdot v))\\ &\qquad + \phi_3(0_h, \kappa_{khg}((\sigma_{kh} Td_1 v) \cdot v)\cdot 0_k).
    \end{split}
    \end{equation*}
    So by \eqref{eq:2Cota-1Cota-kappa-uwv} the terms in $\gamma_3$ agree on both sides of the equation. By multiplying each side of \eqref{eq:2Cota-1Cota-kappa-uwv} by $0_k$ on the right, the terms in $\phi_3$ also agree. 
\end{proof}

\begin{theorem}\label{thm:2Cota-1Cota-ME}
    Let $\huaG$ be a Lie groupoid and $(\kappa, \sigma)$ an Ehresmann connection on $\huaG$. Then, the map $\Xi: T^{2*}\huaG \to \overline{T^{1*}\huaG}$ defined by \eqref{eq:2Cota-1Cota-ME-def} is a Morita equivalence. Hence $T^{2*}\huaG$ is Morita equivalent to $\overline{T^{1*}\huaG}$. 
\end{theorem}

\begin{proof}
    We show that $\Xi$ is a Morita equivalence by using Theorem \ref{thm:we-of-VB2gpd-is-qi}, and showing equivalently that it induces quasi-isomorphisms between the normalized complexes  of $T^{2*}\huaG$ and $\overline{T^{1*}\huaG}$ at each point of $\huaG_0$.

    By Lemma \ref{lem:1CotaBar-norm-cplx}, $N(\widebar{T^{1*}\huaG})=N(T\huaG)^*[-2]$. 
    We claim that $N(\Xi) = \lambda_{\langle \cdot, \cdot \rangle}^r$, the map induced by the IM-pairing associated to the 2-dual pairing, which we know to be a quasi-isomorphism between $N(T^{2*}\huaG)$ and $N(T\huaG)^*[-2]$ by \ref{thm:VBndual-pairing-hndg}. 
    In fact, $N(T^{2*}\huaG)$ can be rewritten in degrees 2,1,0 as 
    \[\begin{tikzcd}[ampersand replacement=\&,cramped, column sep=scriptsize, row sep=tiny]
	{T^*\huaG_0 \oplus A_R^*} \& {A_R^*\oplus A_L^*} \& 0 \\
	{(\theta, \nu)} \& {(\nu, -\theta \circ Td_0 + \nu \circ \ggamma_R)}
	\arrow["{\widecheck{\partial}}", from=1-1, to=1-2]
	\arrow[from=1-2, to=1-3]
	\arrow[maps to, from=2-1, to=2-2]
    \end{tikzcd}\]
    by composing $\Psi$ from Proposition \ref{prop:vs-2dual-norm-cplx} with the isomorphism which is the identity on 2-chains and $(id, \upgamma_R^*): A_R^* \oplus A_R^* \to A_R^* \oplus A_L^*$ on 1-chains, for the core projection $\ggamma_R$ as in \eqref{eq:VB2Dual1Gpd-core-proj}.
    In this picture $N(\Xi)_p$ and $\lambda_{\langle  \cdot, \cdot \rangle}^r$ coincide, as a straightforward computation shows:
    \begin{equation*}
        \begin{split}
            &N(\Xi_1)(\eta, \xi)_p = \xi - \eta \circ \ggamma_R = (\lambda_{\langle \cdot, \cdot \rangle}^r)_1(\eta, \xi)\\
            &N(\Xi_2)(\theta, \nu)_p = \theta = (\lambda_{\langle \cdot, \cdot \rangle}^r)_2(\theta, \nu).
        \end{split}
    \end{equation*}

\end{proof}

\begin{theorem}\label{thm:2Cota-1Cota-SME}
    Let $\huaG$ be a Lie groupoid and $(\kappa, \sigma)$ an Ehresmann connection on $\huaG$. Let $\Xi: T^{2*}\huaG \to \overline{T^{1*}\huaG}$ be defined by \eqref{eq:2Cota-1Cota-ME-def}. 
    Then there exists a form $B \in \Omega^2(T^{2*}\huaG_1)$ such that
    \begin{equation}\label{eq:2Cota-1Cota-SME}
        \omega^{2*} - \Xi_2^*\widebar{\omega} = DB.
    \end{equation}
    That is, the map $\Xi$ is a symplectic Morita equivalence between $T^{2*}\huaG$ and $\overline{T^{1*}\huaG}$. 
\end{theorem}

\begin{proof}
    We show that this holds for $B = - d (K \theta_{can})$, where $\theta_{can}$ is the tautological form on $T^*\huaG = T^{1*}\huaG_1$, which induces a 1-form on  $d_1^*A_R^* \oplus d_0^*A_L^* = T^{2*}\huaG_1$ at each point $g \in \huaG_1$ by applying 
    \begin{equation*}
        K|_g = \kappa_g^* \circ pr_1 :  T^{2*}\huaG_1|_g = d_1^*A_R^* \oplus d_0^*A_L^*|_g \overset{pr_1}{\to} d_1^*A_R^*|_g   \overset{\kappa_g^*}{\to} T^*_g\huaG_1,
    \end{equation*}
    where $pr_1$ is the canonical projection. In particular $K$ is a bundle map over the identity of $\huaG_1$, which means $q^{1*}_1 K = q^{2*}_1$.
    With this, for any $W \in TT^{2*}\huaG_1|_{(\alpha, \beta)_g}$
    \begin{equation}\label{eq:2Cota-1Cota-SME-KTheta1}
    \begin{split}
        (K^*\theta_{can})|_{(\alpha,\beta)_g}(W) 
        &= \theta_{can}|_{K|_g(\alpha,\beta)}(TK W) \\
        &= (\kappa^*_g \alpha)|_g(Tq^{1*}_1TK W)
        = \alpha(\kappa_g(Tq^{2*}_1W)).        
    \end{split}
    \end{equation}
    For this $B$ we have that \eqref{eq:2Cota-1Cota-SME} becomes 
    \begin{equation*}
    \begin{split}
        -d \Theta^{2*} + \Xi_2^*d \widebar{\Theta} 
        &=  \delta (d (K \theta_{can}))\\
        -d \Theta^{2*} + d \Xi_2^* \widebar{\Theta} 
        &=  d \delta (K \theta_{can}),
    \end{split}
    \end{equation*}
    which holds if and only if 
    \begin{equation}\label{eq:2Cota-1Cota-SME-2}
        \Theta^{2*} - \Xi_2^* \widebar{\Theta} = - \delta (K \theta_{can}).
    \end{equation}
    We evaluate each side of this equation on an arbitrary $V \in TT^{2*}\huaG_2|_{(\phi, \alpha, \beta, \gamma)_{(g,h)}}$. For brevity's sake, define $V^{\downarrow}_i := Td_iTq^{2*}_2V$ for $0 \le i \le 2$, so that $Tq^{2*}V = (V^{\downarrow}_0, V^{\downarrow}_2) \in T\Lambda^2_1(\huaG)|_{(g,h)} \cong \Lambda^2_1(T\huaG)|_{g,h}$.
    By \eqref{eq:2Cota-taut-form-1gpd}, \eqref{eq:2Cota-Bar-taut-form}, and \eqref{eq:2Cota-1Cota-ME-def}, evaluating the right-hand side of \eqref{eq:2Cota-1Cota-SME-2} on $V$ gives
    \begin{equation*}
    \begin{split}
        (\Theta^{2*} - &\Xi_2^* \widebar{\Theta})|_{(\phi, \alpha, \beta, \gamma)_{(g,h)}}(V)\\ 
        &= \phi|_{(g,h)}(V^{\downarrow}_0, V^{\downarrow}_2) 
        - \alpha|_{d_1g}(\kappa_g(V^{\downarrow}_0)) 
        + \gamma|_{d_1h}(\kappa_{hg}((\sigma_hTd_1V^{\downarrow}_0)\cdot V^{\downarrow}_0)) \\
        &\qquad \qquad - \phi|_{(g,h)}(V^{\downarrow}_0, - \kappa_{hg}((\sigma_hTd_1V^{\downarrow}_0)\cdot V^{\downarrow}_0) \cdot 0_h + \sigma_hTd_1V^{\downarrow}_0)\\
        &= - \alpha|_{d_1g}(\kappa_g(V^{\downarrow}_0)) 
        + \gamma|_{d_1h}(\kappa_{hg}((\sigma_hTd_1V^{\downarrow}_0)\cdot V^{\downarrow}_0)) \\
        &\qquad \qquad + \phi|_{(g,h)}(0_g, V^{\downarrow}_2 + \kappa_{hg}((\sigma_hTd_1V^{\downarrow}_0)\cdot V^{\downarrow}_0) \cdot 0_h - \sigma_hTd_1V^{\downarrow}_0).
    \end{split}
    \end{equation*}    
    By using \eqref{eq:2Cota-1Cota-SME-KTheta1} and the fact that $Tq^{2*}_1T\widecheck{d}_i V = Td_i Tq^{2*}_2V = V^{\downarrow}_i$, evaluating the left-hand side of \eqref{eq:2Cota-1Cota-SME-2} on $V$ gives
    \begin{equation*}
    \begin{split}
        - (\delta(K &\theta_{can}))|_{(\phi, \alpha, \beta, \gamma)_{(g,h)}}(V) \\
        &= - (\widecheck{d}_0^*(K \theta_{can}) 
        - \widecheck{d}_1^*(K \theta_{can}) 
        + \widecheck{d}_2^*(K \theta_{can}))|_{(\phi, \alpha, \beta, \gamma)_{(g,h)}}(V)\\
        &= - (K \theta_{can})|_{\widecheck{d}_0(\phi, \alpha, \beta, \gamma)}(T\widecheck{d}_0V) 
        + (K \theta_{can})|_{\widecheck{d}_1(\phi, \alpha, \beta, \gamma)}(T\widecheck{d}_1V)\\ 
        &\qquad\qquad\qquad - (K \theta_{can})|_{\widecheck{d}_2(\phi, \alpha, \beta, \gamma)}(T\widecheck{d}_2V) \\
        & = - \alpha|_{d_1g}(\kappa_g(Tq^{2*}_1T\widecheck{d}_0V)) 
        + \gamma|_{d_1h}(\kappa_{hg}(Tq^{2*}_1T\widecheck{d}_1V))\\
        & \qquad \qquad + \phi|_{(g,h)}(0_g, \kappa_{hg}(Tq^{2*}_1T\widecheck{d}_1V)\cdot 0_h ) 
        - \gamma|_{d_1h}(\kappa_h(Tq^{2*}_1T\widecheck{d}_2V))\\
        & = - \alpha|_{d_1g}(\kappa_g(V^{\downarrow}_0)) 
        + \gamma|_{d_1h}(\kappa_{hg}(V^{\downarrow}_1) - \kappa_h(V^{\downarrow}_2)) 
        + \phi|_{(g,h)}(0_g, \kappa_{hg}(V^{\downarrow}_1)\cdot 0_h ).
    \end{split}
    \end{equation*}
    Therefore \eqref{eq:2Cota-1Cota-SME-2} is equivalent to 
    \begin{equation}\label{eq:2Cota-1Cota-SME-3}
    \begin{cases}
        \kappa_{hg}((\sigma_hTd_1V^{\downarrow}_0)\cdot V^{\downarrow}_0) = \kappa_{hg}(V^{\downarrow}_1) - \kappa_h(V^{\downarrow}_2), & \\
        V^{\downarrow}_2 + \kappa_{hg}((\sigma_hTd_1V^{\downarrow}_0)\cdot V^{\downarrow}_0) \cdot 0_h - \sigma_hTd_1V^{\downarrow}_0 = \kappa_{hg}(V^{\downarrow}_1)\cdot 0_h. &
    \end{cases}
    \end{equation}
    By using the splitting \eqref{eq:2Cota-Ehresmann-conn-splitting}, we have
    \begin{equation*}
        V^\downarrow_2 = \kappa_h(V^{\downarrow}_2) \cdot 0_h + \sigma_h Td_0 V^{\downarrow}_2 = \kappa_h(V^{\downarrow}_2) \cdot 0_h + \sigma_h Td_1 V^{\downarrow}_0.
    \end{equation*}
    Because of this, the first equation in \eqref{eq:2Cota-1Cota-SME-3} implies the second: 
    \begin{equation*}
        V^{\downarrow}_2 + \kappa_{hg}(V^{\downarrow}_1) \cdot 0_h - \kappa_h(V^{\downarrow}_2) \cdot 0_h - \sigma_hTd_1V^{\downarrow}_0 = \kappa_{hg}(V^{\downarrow}_1)\cdot 0_h
    \end{equation*}
    To show the first equation holds we use the fact that $V^{\downarrow}_1 = V^{\downarrow}_2\cdot V^{\downarrow}_0$ as this is the face indexed by 1 of a (2,1)-horn, and the splitting of $V^{\downarrow}_2$ above:
    \begin{equation*}
    \begin{split}
        \kappa_{hg}(V^{\downarrow}_1) - \kappa_h(V^{\downarrow}_2) 
        &= \kappa_{hg}(V^{\downarrow}_2\cdot V^{\downarrow}_0) - \kappa_h(V^{\downarrow}_2)\\
        &= \kappa_{hg}((\kappa_h(V^{\downarrow}_2) \cdot 0_h + \sigma_h Td_1 V^{\downarrow}_0)\cdot (0_g + V^{\downarrow}_0)) - \kappa_h(V^{\downarrow}_2)\\
        &= \kappa_{hg}((\kappa_h(V^{\downarrow}_2) \cdot 0_{hg}) + (\sigma_h Td_1 V^{\downarrow}_0 \cdot V^{\downarrow}_0)) - \kappa_h(V^{\downarrow}_2)\\
        &= \kappa_{hg}((\sigma_hTd_1V^{\downarrow}_0)\cdot V^{\downarrow}_0).
    \end{split}
    \end{equation*}
\end{proof}



\begingroup
\setlength{\emergencystretch}{3em}
\printbibliography[heading=bibintoc]
\endgroup

\end{document}